\documentclass[11pt,reqno]{amsart}

\usepackage{Cubic-NLW-Diagrams}
\usepackage{Cubic-NLW-packages}

\begin{document}

\author[Bjoern Bringmann]{Bjoern Bringmann$^\dagger$}
\thanks{$\dagger$ Corresponding author}
\address{Bjoern Bringmann, School of Mathematics, Institute for Advanced Study, Princeton, NJ 08540 \& Department of Mathematics, Princeton University, Princeton, NJ 08540}
\email{bjoern@ias.edu}

\author[Yu Deng]{Yu Deng}
\address{Yu Deng, Department of Mathematics, University of Southern California, Los Angeles,  CA 90089, USA }
\email{yudeng@usc.edu}
\author[Andrea R. Nahmod]{Andrea R. Nahmod}
\address{Andrea R. Nahmod,  
Department of Mathematics and Statistics,  University of Massachusetts Amherst,  Amherst, MA 01003}
\email{nahmod@math.umass.edu}
\author[Haitian Yue]{Haitian Yue}
\address{Haitian Yue,
Institute of Mathematical Sciences, ShanghaiTech University, Shanghai, 201210, China}
\email{yuehaitian@shanghaitech.edu.cn}
\date{}

\title{Invariant Gibbs measures for the three dimensional cubic nonlinear wave equation}

\linespread{1.14}

\begin{abstract}
We prove the invariance of the Gibbs measure under the dynamics of the three-dimensional cubic wave equation, which is also known as the hyperbolic $\Phi^4_3$-model. This result is the hyperbolic counterpart to seminal works on the parabolic $\Phi^4_3$-model by Hairer \cite{H14} and Hairer-Matetski \cite{HM18}. 

The heart of the matter lies in establishing local in time existence and uniqueness of solutions on the statistical ensemble, which is achieved by using a para-controlled Ansatz for the solution, the analytical framework of the random tensor theory, and the combinatorial molecule estimates.  

The singularity of the Gibbs measure with respect to the Gaussian free field brings out a new caloric representation of the Gibbs measure and a synergy between the parabolic and hyperbolic theories embodied in the analysis of heat-wave stochastic objects. Furthermore from a purely hyperbolic standpoint our argument relies on key new ingredients that include a hidden cancellation between sextic stochastic objects and a new bilinear random tensor estimate.

 \end{abstract}
\maketitle

\tableofcontents

\section{Introduction}\label{section:introduction}

In this article, we study the (renormalized) cubic nonlinear wave equation in three dimensions. We begin our discussion with the corresponding Hamiltonian, which is formally given by 
\begin{equation}\label{intro:eq-Hamiltonian}
H[u,\partial_t u] = \int_{\T^3} \dx \bigg( \frac{u^2}{2} + \frac{|\nabla u|^2}{2} + \frac{(\partial_t u )^2}{2} \bigg) +  \int_{\T^3} \dx \bigg( \frac{\lcol u^4 \rcol}{4} + \infty \cdot \frac{u^2}{2} \bigg).
\end{equation}
Here, $\mathbb{T}^3 := (\R/ 2\pi \mathbb{Z})^3$ is the three-dimensional torus, $\lcol u^4 \rcol$ is the Wick-ordered quartic power, and $\infty \cdot u^2$ represents an additional mass-renormalization. We further define the canonical symplectic form
\begin{equation}\label{intro:eq-symplectic}
\omega \Big( (u_0,u_1) , (v_0,v_1) \Big) := \int_{\T^3} \dx \Big( u_0 v_1 - u_1 v_0 \Big). 
\end{equation}
The nonlinear evolution equation corresponding to the Hamiltonian $H$ and symplectic form $\omega$ is given by the (renormalized) cubic nonlinear wave equation. It can be written as 
\begin{equation}\label{intro:eq-NLW}
\begin{cases}
\big( \partial_t^2 + 1 - \Delta \big) u = - \lcol u^3 \rcol -\, \infty \cdot u \qquad (t,x) \in \mathbb{R} \times \mathbb{T}^3, \\
(u, \langle \nabla \rangle^{-1} \partial_t u)\big|_{t=0} = (\phi^{\cos},\phi^{\sin}),
\end{cases}
\end{equation}
where $\lcol u^3 \rcol$ denotes the Wick-ordered cubic power and $\phi^{\cos},\phi^{\sin}\colon \T^3 \rightarrow \R$. We note that the initial condition is written in terms of $\langle \nabla \rangle^{-1}\partial_t u$, which is purely notational but convenient (see Section \ref{section:ansatz}). Furthermore, the Gibbs measure corresponding to the Hamiltonian $H$ and symplectic form $\omega$ is formally given by 
\begin{equation}\label{intro:eq-Gibbs}
`` \mathrm{d}\mu(\phi^{\cos},\phi^{\sin}) = \mathcal{Z}^{-1} \exp\Big( - H[\phi^{\cos}, \langle \nabla \rangle \phi^{\sin}]\Big) \mathrm{d}\phi^{\cos} \mathrm{d}\phi^{\sin}. "
\end{equation}
We now state a formal version of our main theorem, but postpone a rigorous version to Theorem \ref{intro:thm-rigorous} below.

\begin{theorem}[Global well-posedness and invariance, formal version]\label{intro:thm-formal}
The cubic nonlinear wave equation \eqref{intro:eq-NLW} is almost surely globally well-posed with respect to the Gibbs measure $\mu$. Furthermore, the global solution leaves the Gibbs measure invariant. 
\end{theorem}

We now briefly discuss both the physical and mathematical significance of Theorem \ref{intro:thm-formal}, but postpone more detailed discussions to the following subsections. From a physical perspective, the interest in our main theorem stems from its connections to quantum field theory. In fact, the Gibbs measure from \eqref{intro:eq-Gibbs} is essentially identical to the so-called $\Phi^4_3$-measure, whose construction was one of the crowning achievements of the constructive quantum field theory program (see e.g. \cite{GJ73,GJ87}). 

From a mathematical perspective, the interest in our main theorem and its proof stem from connections to several areas of mathematical research. At its core, Theorem \ref{intro:thm-formal} concerns a refined understanding of how randomness gets transported by the flow of a nonlinear equation
which involves probability theory, and partial differential equations. Such study of propagation of randomness under nonlinear evolution equations is in itself of independent interest and also plays an important role in mathematical wave turbulence \cite{DNY20,DH21} and singular parabolic SPDEs \cite{H14,GIP15,K16}. In our proof, we further utilize methods from combinatorics (see Section \ref{section:analytic2}), elementary geometry (see Section \ref{section:counting}), harmonic analysis (see Section \ref{section:para}), random matrix theory (see Section \ref{section:bilinear} and Section \ref{section:linear}), and nonlinear Fourier analysis which permeates much of the paper.

\subsection{Earlier results}\label{section:intro-earlier-results}

In this subsection, we give a broad overview of the relevant previous literature. A more detailed discussion of several individual articles is postponed until the next subsection (Subsection \ref{section:intro-earlier-methods}). 

\subsubsection{\protect{The $\Phi^4_d$-models}} 
Many of the most significant contributions in constructive quantum field theory, singular SPDEs, and random dispersive equations concern the so-called $\Phi^4_d$-models. The starting point of our discussion is the energy 
\begin{equation}\label{intro:eq-V}
H(\phi) = \int_{\T^d} \dx \bigg( \frac{|\phi|^2}{2} + \frac{|\nabla \phi|^2}{2} + \frac{|\phi|^4}{4} \bigg). 
\end{equation}
In \eqref{intro:eq-V}, we omit the possible renormalization of the potential energy term, since its precise form depends on the dimension $d$. Equipped with \eqref{intro:eq-V}, we formally define the $\Phi^4_d$-measure as  
\begin{equation}\label{intro:eq-Phi4d}
``\mathrm{d}\Phi^4_d(\phi) = \mathcal{Z}^{-1} \exp\bigg( - \int_{\T^d} \dx  \bigg( \frac{|\phi|^2}{2} + \frac{|\nabla \phi|^2}{2} + \frac{|\phi|^4}{4} \bigg) \bigg) \mathrm{d}\phi. " 
\end{equation}
In addition to the $\Phi^4_d$-measure, the energy  in \eqref{intro:eq-V} induces three different evolution equations, which are called dynamical $\Phi^4_d$-models:
\begin{enumerate}
    \item A Langevin equation, which is given by the cubic heat equation with space-time white noise,
    \item a real-valued Hamiltonian equation, which is given by the cubic wave equation,
    \item and a complex-valued Hamiltonian equation, which is given by the cubic Schr\"{o}dinger equation. 
\end{enumerate}
The main problems then concern the construction of the $\Phi^4_d$-measure, the probabilistic well-posedness of the three evolution equations, as well as the invariance\footnote{The Gibbs measure of the wave and Schr\"{o}dinger dynamics only slightly differs from the $\Phi^4_d$-measure. In the wave setting, the Gibbs measure is given by a product of the $\Phi^4_d$-measure and a Gaussian measure, see e.g. \eqref{intro:eq-Gibbs}. In the Schr\"{o}dinger setting, the Gibbs measure is given by the complex-valued instead of real-valued $\Phi^4_d$-measure.} of the $\Phi^4_d$-measure under the three evolution equations. We now discuss parts of the literature on $\Phi^4_d$-models, which is also illustrated in Figure \ref{figure:literature-cubic}.   

The construction and properties of the $\Phi^4_d$-measure were extensively studied by constructive quantum field theorists. The construction of the $\Phi^4_1$-measure\footnote{The reason why the $\Phi^4_1$-measure was not studied by many constructive quantum field theorists is that it corresponds to a field theory in $(1+0)$-dimensional Minkowski space, which is just the time-axis.} is rather elementary and, for example, can be found in \cite{B94}. The $\Phi^4_2$ and $\Phi^4_3$-measure were first constructed in \cite{Nelson66} and \cite{GJ73}, respectively. For a more detailed overview of the relevant literature on this construction, we refer the reader to the introduction of \cite{GH21} and the monograph \cite{GJ87}. Aizenman \cite{A81} and Fr\"{o}hlich \cite{F82} proved the marginal triviality of the $\Phi^4_d$-measure in dimension $d\geq 5$. In recent seminal work, Aizenman and Duminil-Copin \cite{ADC21} also proved the marginal triviality of the (real-valued) $\Phi^4_4$-measure. In the complex-valued setting, which is most important for the Schrödinger equation, marginal triviality (as in \cite{ADC21}) is expected to hold but has not yet been proven. For partial progress, which addresses the small coupling regime, we refer to \cite{CCF83,BBS14}. Loosely speaking, marginal triviality implies that the $\Phi^4_d$-measure in dimension $d\geq 4$  essentially yields a Gaussian measure, for any renormalization of the potential energy term in \eqref{intro:eq-V} giving rise to a well defined measure.
In particular, this invalidates the question regarding the invariance of the $\Phi^4_d$-measure under the dynamical $\Phi^4_d$-models. We note that it is still possible to study the triviality of the dynamics (see e.g. the discussion in \cite{HairerTalk}), but we view this as a separate problem from the invariance of the Gibbs measure. Recently, there has also been tremendous progress on the derivation of the $\Phi^4_d$ and related measures from many-body quantum systems \cite{FKSS17,FKSS20,FKSS22,LNT15,LNT21}, which goes beyond their construction.

The parabolic $\Phi^4_d$-models, i.e., the cubic heat equation with space-time white noise, was first studied in one and two dimensions in \cite{I87} and \cite{DPD03}, respectively. In three dimensions, the probabilistic well-posedness was first shown by Hairer \cite{H14} using the theory of regularity structures. The probabilistic well-posedness has now also been shown using para-controlled calculus \cite{CC18,GIP15} and renormalization group methods \cite{K16}. The invariance of the $\Phi^4_3$-measure under the parabolic $\Phi^4_3$-dynamics was first obtained by Hairer and Matetski \cite{HM18}. More recent articles, such as \cite{AK20,GH21,MW17,MW20}, address the stochastic quantization in the spirit of Parisi and Wu \cite{PW81}, i.e., the construction of the $\Phi^4_3$-measure using the parabolic $\Phi^4_3$-model as explained in Subsection \ref{section:ansatz-caloric} below.

The hyperbolic $\Phi^4_d$-model, i.e., the cubic wave equation with initial data drawn from its associated Gibbs measure (and/or space-time white noise), was first studied in one dimension by Friedlander \cite{F85} and Zhidkov \cite{Z94}. In two dimension, the invariance of the Gibbs measure was stated as a theorem in \cite[Theorem 111]{B99}. For a complete proof, we also refer the reader to \cite{OT20}. In three dimensions, the invariance of the Gibbs measure is the main theorem of this article. 

In the case of the cubic Schr\"{o}dinger equation, the invariance of its associated Gibbs measure was proven in one and two dimensional tori in seminal works of Bourgain \cite{B94,B96}. In three dimensions, this question is a famous open problem, which will be briefly discussed in Subsection \ref{section:intro-open} below.

\begin{figure}
\begin{tabular}{
!{\vrule width 1pt}>{\centering\arraybackslash}P{2cm}
!{\vrule width 1pt}>{\centering\arraybackslash}P{2.5cm}
!{\vrule width 1pt}>{\centering\arraybackslash}P{2.5cm}
!{\vrule width 1pt}>{\centering\arraybackslash}P{2.5cm}
!{\vrule width 1pt}>{\centering\arraybackslash}P{2.5cm}
!{\vrule width 1pt}} 
\noalign{\hrule height 1pt} & & & & 
 \\[-1.5ex] 
Dimension &  Measure & Heat & Wave  & Schr\"{o}dinger 
\\[4pt] \noalign{\hrule height 1pt} \rule{0pt}{14pt}
$d=1$	
&  \cellcolor{Green!10} 
& \cite{I87} \cellcolor{Green!10} 
& \cellcolor{Green!10} \cite{F85,Z94}
& \cellcolor{Green!10} 	\cite{B94} 
\\[4pt] \noalign{\hrule height 1pt} \rule{0pt}{14pt}
$d=2$ 
& \cite{Nelson66} \cellcolor{Green!10} 
& \cite{DPD03} \cellcolor{Green!10} 
&  \cite{B99}	\cellcolor{Green!10} 
& \cellcolor{Green!10} \cite{B96}		
\\[4pt] \noalign{\hrule height 1pt} \rule{0pt}{14pt}
$d=3$
& \cite{GJ73}	 \cellcolor{Green!10} 
&  \cite{H14} \cellcolor{Green!10} 
& 		\cellcolor{Green!10} \tikzmarkin{a}(0.34,-0.24)(-0.34,0.46) \emph{This article}\tikzmarkend{a}
&  \emph{Open}
\\[4pt] \noalign{\hrule height 1pt} \rule{0pt}{14pt}
$d=4$  
&	\multicolumn{4}{c!{\vrule width 1pt}}{\cellcolor{red!10}\cite{ADC21}}   
\\[4pt] \noalign{\hrule height 1pt} \rule{0pt}{14pt}
$d\geq 5$ 
&	\multicolumn{4}{c!{\vrule width 1pt}}{\cellcolor{red!10}\cite{A81,F82}}   
\\[3pt] \noalign{\hrule height 1pt} 
\end{tabular}
\caption{\small{Existence and invariance of the Gibbs measure for the cubic stochastic heat, wave, and Schrödinger equations.}}
\label{figure:literature-cubic}
\end{figure}

\subsubsection{Invariant Gibbs measures in other dispersive settings}\label{section:intro-other-settings}
In one and two spatial dimensions Gibbs measures can be constructed on tori not just for $\Phi^4_d$ but also for  $\Phi^{p+1}_d$ with $p\geq 5$ odd, which contain (defocusing) Wick ordered interaction polynomial potentials $\lcol u^{p}\rcol$ (or $\lcol |u|^{p-1}u\rcol$) in \eqref{intro:eq-V}.  The construction of such measures also stems from seminal works of Nelson \cite{Nelson66, Nelson73} 
in the context of quantum field theory (see also \cite{Simon, GJ87}). A concise exposition of the construction can be found in \cite{OT18, OT20}. 

Regarding the dynamics, see Figure \ref{figure:literature-dispersive}, in one dimension, the invariance of such Gibbs measures with respect to the dynamics of the corresponding Wick ordered nonlinear Schr\"odinger equation with $p$-th power nonlinearity was proved by Bourgain in \cite{B94}.  In two dimensions, the invariance of such Gibbs measures under the flow of the corresponding Wick ordered nonlinear \emph{wave} equation with $p$-th power nonlinearity was proved in \cite{B99, OT20}  via a linear-nonlinear decomposition, known nowadays as Bourgain's trick (or as DaPrato-Debussche's trick \cite{DaPD} in the stochastic community). However for the Wick ordered nonlinear \emph{Schr\"odinger} equation with $p$-th power nonlinearity, the invariance of its associated Gibbs measures in two dimensions (in the strong sense as in \cite{B96}) does not follow from the same type of arguments. Indeed, in \cite{DNY19}, the second, third and fourth authors of this paper introduced the method of random averaging operators to obtain this result which will be discussed in Subsection \ref{section:intro-tensor} below.

Gibbs measures on tori can also be constructed for (suitably renormalized) Hartree nonlinearities $(V_{\beta} \ast |u|^2) u$, where the interaction potential $V_{\beta}$ behaves like $|x|^{-(d-\beta)}$ (see e.g. \cite{B97,B20I,OOT20}).  Their invariance under the flow of the corresponding Hartree nonlinear wave and Schr\"odinger equations was also studied. This was done in the Schr\"odinger case in three dimensions for $\beta >2$ by Bourgain{\footnote{In \cite{B97}, Bourgain treats also the two dimensional case where the results hold for any $\beta >0$.}} in \cite{B97} and for (say) $\beta>0.99$ by the second, third and fourth authors of this paper in \cite{DNY21}. In the nonlinear wave case this was achieved  by Oh, Okamoto and Tolomeo in \cite{OOT20} for $\beta >1$ and by the first author of this paper in \cite{B20II} for any $\beta >0$. The case of $0< \beta \leq 1/2$ is of special significance since the Gibbs measure is then singular with respect to the Gaussian free field (see \cite{B20I,OOT20}) just as the case in this paper.

In addition to the nonlinear wave and Schrödinger equations with odd power-type and Hartree nonlinearities discussed above, invariant Gibbs measures have also been studied in several other settings. For instance, there has been research on invariant Gibbs measures for derivative nonlinearities \cite{BLS21, Deng_BO,Tz_BO}, quadratic nonlinearities \cite{GKO18, OOT21}, exponential nonlinearities \cite{ORTW20,ORTW21,SunT20b}, radially-symmetric settings \cite{BB14a,BB14b,BB14c,D12,T06}, KdV and generalized KdV \cite{B94,CK21,ORT16,R16}, fractional dispersion relations \cite{ST20a,SunT21}, and lattice models \cite{AKV20}.

For integrable equations, one can also construct and study the invariance of weighted Wiener measures associated to higher order conservation laws than the Hamiltonian. There has been substantial activity in this regard and we do not go into further details here, but refer the interested reader to a few references  such as \cite{
DTzV, GLV1, GLV2, NORS, NRSS, ThTz, TzV1, TzV2, Zhidkov_KdV, Zhidkov_NLS} for some examples.

\begin{figure}
\begin{tabular}{
!{\vrule width 1pt}>{\centering\arraybackslash}P{2cm}
!{\vrule width 1pt}>{\centering\arraybackslash}P{2cm}
!{\vrule width 1pt}>{\centering\arraybackslash}P{4cm}
!{\vrule width 1pt}>{\centering\arraybackslash}P{4cm}
!{\vrule width 1pt}} 
\noalign{\hrule height 1pt} & & &  
 \\[-1.5ex] 
Dimension &  Nonlinearity &  Wave  & Schr\"{o}dinger 
\\[4pt] \noalign{\hrule height 1pt} \rule{0pt}{14pt}
\multirow{3}*{$d=2$}
&  $|u|^2 u$
& \cite{B99}
& \cite{B96}
\\[4pt]  \cline{2-4} \rule{0pt}{14pt}
& $|u|^{2r} u$
& \cite{B99}
& \cite{DNY19}
\\[4pt] \noalign{\hrule height 1pt} \rule{0pt}{14pt}
\multirow{4}*{$d=3$}
& $u^2$
& \cite{GKO18,OOT21}
&  \cellcolor{Black!50} 
\\[4pt]  \cline{2-4} \rule{0pt}{14pt}
& $(V_\beta \ast |u|^2 )u $
& \begin{tabular}{rl}
& \\[-2ex]
$\beta>1:$ & \cite{OOT20} \\
$\beta>0:$ & \cite{B20II} \\[0.5ex]
\end{tabular}		
& \begin{tabular}{rl}
& \\[-2ex]
$\beta>2:$ & \cite{B97} \\
$\beta>1-\epsilon:$ & \cite{DNY21} \\
$\beta>0:$ & \emph{Open} \\[0.5ex]
\end{tabular}		
\\[4pt]  \cline{2-4} \rule{0pt}{14pt}
& $|u|^{2} u$
& \emph{This article}
& \emph{Open} 
\\[3pt] \noalign{\hrule height 1pt} 
\end{tabular}
\caption{\small{Invariance of the Gibbs measure for wave and Schrödinger equations.}}
\label{figure:literature-dispersive}
\end{figure}

\subsection{Earlier ideas and methods} \label{section:intro-earlier-methods}
In this subsection, we describe some of the main ideas from the previous literature in more detail. 

\subsubsection{Scaling heuristics}\label{section:intro-scaling} The notion of \emph{probabilistic scaling} was introduced in \cite{DNY19,DNY20} by the second, third and fourth authors of this paper. It provides a guiding principle when studying nonlinear dispersive equations in the random or stochastic setting, just as what the ordinary (deterministic) scaling does with the corresponding deterministic well-posedness problems. The idea is to compare the (expected) regularity of the nonlinear evolution with the linear evolution (or noise), under the assumption that all the frequencies are comparable, and determine the threshold where the nonlinear evolution becomes rougher than the linear one; see \cite{DNY19,DNY20} for details.

In \cite{DNY20} it is proved that any nonlinear Schr\"{o}dinger equation with odd power nonlinearity is probabilistically locally well-posed at any probabilistically subcritical regularity. We expect the same guiding principle to work for other dispersive equations, including the wave equation (albeit with some additional twists, see Remark \ref{intro:addrem}). In this regard, one can calculate that the probabilistic critical index for (\ref{intro:eq-NLW}) is $s_{pr}=-3/4$, while the regularity of the Gibbs measure initial data is at $H^{-1/2-}$ (see Figure \ref{intro:scalingfig}). Thus (\ref{intro:eq-NLW}) is probabilistically subcritical and expected to be locally well-posed.
\begin{figure}
\begin{tabular}{|P{3cm}|P{3cm}|P{3cm}|P{3cm}|P{3cm}|} \hline  \rule{0pt}{14pt} 
Equation &2D wave & 2D Schr\"{o}dinger & 3D wave & 3D Schr\"{o}dinger 
\\[4pt]  \cline{1-5} \rule{0pt}{16pt}
Gibbs measure regularity $s_G-$ & $ 0-$ & $ 0-$ & $ -1/2-$ & $ -1/2-$ \\
[4pt]  \cline{1-5} \rule{0pt}{16pt}
Prob. critical regularity $s_{pr} $ &$ -3/4$&$ -1/2$&$ -3/4$&$ -1/2$\\[4pt]\cline{1-5}
\end{tabular}
\caption{\small{Probabilistic critical and Gibbs measure regularities for 2D and 3D wave and Schr\"{o}dinger equations.}}
\label{intro:scalingfig}
\end{figure}
\begin{remark}\label{intro:addrem} Unlike the Schr\"{o}dinger case, there are extra obstacles for nonlinear wave equations, linked to high-high-to-low interactions, which may result in the actual probabilistic well-posedness threshold being \emph{strictly higher} than the scaling predictions. For example the 4D quadratic wave equation is probabilistic critical in $H^{-3/2}$, but one actually expects probabilistic ill-posedness already in $H^{-1-}$, due to high-high-to-low interactions. Note that such discrepancy also exist and is well-known in the context of stochastic heat equations and rough path theory. See \cite[Remark 1.8]{OO21} for a related discussion.
\end{remark}
\subsubsection{The para-controlled approach}\label{section:intro-para-heat}

We first discuss the para-controlled approach of \cite{GIP15} in the context of singular parabolic SPDEs. 
We focus on the parabolic $\Phi^4_3$-model which, in the sense of the parabolic scaling from \cite[Section 8]{H14}, is subcritical. It is formally given by 
\begin{equation}\label{intro:eq-heat}
\begin{cases}
\big( \partial_s + 1 - \Delta \big) \Phi = - \lcol \Phi^3 \rcol - \infty \cdot \Phi + \sqrt{2} \mathrm{d}W \qquad (s,x) \in (s_0,\infty) \times \T^3, \\
\Phi\big|_{s=s_0}=\phi, 
\end{cases}
\end{equation}
where $\mathrm{d}W$ denotes space-time white noise, $s_0 \in \R$ is the initial time, and $\phi\colon \T^3_x \rightarrow \R$. We use $s\in \R$ for the parabolic time-variable since $t\in \R$ is reserved for the time-variable of wave equations\footnote{The reason for using different variable names is that the stochastic objects in this article involve both heat and wave propagators (see Subsection \ref{section:diagram-mixed}).}. Due to the smoothing effects of the heat propagator, the precise structure of the initial data in  \eqref{intro:eq-heat} is not important, and we only assume that $\phi \in \C_x^{-1/2-\epsilon}(\T^3)$. Instead, the difficulties in \eqref{intro:eq-heat} stem from the stochastic forcing $\mathrm{d}W$. It can be written as 
\begin{equation}\label{intro:eq-dW}
\mathrm{d}W_s(x) = \sum_{n\in \Z^3} e^{i\langle n,x\rangle} \mathrm{d}W_s(n), 
\end{equation}
where $(W_s(n))_{n\in \Z^3}$ are standard complex-valued Brownian motions satisfying the constraint $\overline{W_s(n)}=W_s(-n)$. The parabolic space-time regularity of $\mathrm{d}W$ is $(2+3) (-1/2-)=-5/2-$, where the time-dimension is counted twice and the three spatial dimensions are counted once. \\

\begin{figure}
\begin{tabular}{|P{2cm}|P{1.25cm}|P{1.25cm}|P{1.25cm}|P{1.25cm}|} \hline  \rule{0pt}{14pt} 
Regularity &$-1/2-$ & $1/2-$ & $1-$ & $3/2-$ 
\\[4pt]  \cline{1-5} \rule{0pt}{16pt}
Object & $\shlinear$ & $\shcubic$ & X & Y  \\[4pt]  \cline{1-5} 
\end{tabular}
\caption{\small{Random structure of the cubic stochastic heat equation.}}
\label{figure:structure-heat}
\end{figure}

The first step in the analysis of \eqref{intro:eq-heat} is the analysis of linear and cubic stochastic objects. We define the linear stochastic object $\shlinear$ as the stationary solution to the stochastic heat equation
\begin{equation}\label{intro:eq-shlinear}
(\partial_s +1 - \Delta) \, \shlinear = \sqrt{2} \mathrm{d}W. 
\end{equation}
Since the heat equation gains two derivatives and $\mathrm{d}W$ has regularity $-5/2-$, the linear stochastic object $\shlinear$ has regularity $-1/2-$. We also define\footnote{The renormalization in \eqref{intro:eq-shcubic} only contains the Wick-ordering, but does not contain the additional renormalization in \eqref{intro:eq-heat}. This is because $``\infty \cdot \Phi"$ cancels double-resonances, which occur in terms such as $\scalebox{0.8}{$\shcubic \shquadratic$}$, but not the single resonance which occurs in the cubic nonlinearity $\big( \shlinear \big)^3$.} the cubic stochastic object $\shcubic$ as the stationary solution to 
\begin{equation}\label{intro:eq-shcubic}
(\partial_s + 1 - \Delta) \shcubic = \lcol \big( \, \shlinear\, \big)^3 \rcol. 
\end{equation}
Since the heat equation gains two derivatives and $\shlinear$ has regularity $-1/2-$, the cubic stochastic object has regularity $2+3\cdot (-1/2-)=1/2-$.\\

The next step in the analysis of \eqref{intro:eq-heat} is the analysis of the nonlinear remainder $\Psi$, which is defined as
\begin{equation}\label{intro:eq-Psi}
\Psi = \Phi - \shlinear + \shcubic. 
\end{equation}
Due to \eqref{intro:eq-heat}, \eqref{intro:eq-shlinear}, and \eqref{intro:eq-shcubic}, the nonlinear remainder $\Psi$ satisfies
\begin{equation}\label{intro:eq-heat-Psi}
\begin{cases}
\big( \partial_s + 1 - \Delta \big) \Psi = - 3 \, \shquadratic \Big( - \shcubic + \Psi \Big) -3 \, \shlinear \Big( - \shcubic + \Psi \Big)^2 - \Big( - \shcubic + \Psi \Big)^3 \\
\hspace{18ex}- \infty \cdot \Big( \shlinear - \shcubic + \Psi\Big),  \\
\Psi\big|_{s=s_0}=\phi - \shlinear (s_0) + \shcubic (s_0). 
\end{cases}
\end{equation}
where $\shquadratic = \, \lcol  \big( \, \shlinear \big)^2 \rcol$ is the Wick-ordered square. While subtracting $\shlinear$ and adding $\shcubic$ has removed the lowest regularity terms on the right-hand side of \eqref{intro:eq-heat}, \eqref{intro:eq-heat-Psi} still cannot be solved using a direct contraction argument. Since $\shquadratic$ has regularity $-1-$, the Duhamel integral of the high$\times$low-interaction
\begin{equation}\label{intro:eq-heat-high-low}
\shquadratic \paragg \Big( \Psi - \shcubic \Big) 
\end{equation}
has regularity at most $1-$. As a result, we also expect that $\Psi$ has regularity at most $1-$. However, this regularity is insufficient to even define the high$\times$high-interaction 
\begin{equation}\label{intro:eq-heat-high-high}
\shquadratic \parasim \Psi, 
\end{equation}
which prevents us from closing a direct contraction argument for $\Psi$. The main idea of the para-controlled approach of Gubinelli, Imkeller, and Perkowski \cite{GIP15}, which was first applied to the stochastic heat equation in \cite{CC18}, is to separate the treatment of \eqref{intro:eq-heat-high-low} and \eqref{intro:eq-heat-high-high}. To this end, the nonlinear remainder is written as  $\Psi=X+Y$, where $Y$ is a smooth nonlinear remainder with regularity $3/2-$ and $X$ is a para-controlled component with regularity $1-$. Namely the para-controlled component $X$ absorbs the contribution of \eqref{intro:eq-heat-high-low} and solves\footnote{In \cite{GIP15,CC18}, the equation for $X$ is further simplified using commutator estimates for the Duhamel integral, but we omit these (important) aspects from our discussion.}
\begin{equation}\label{intro:eq-heat-X}
(\partial_s + 1 - \Delta) X = - 3 \shquadratic \paragg \Big( - \shcubic + X + Y \Big). 
\end{equation}
Together with \eqref{intro:eq-heat-Psi}, \eqref{intro:eq-heat-X} also determines the nonlinear evolution equation for the smooth nonlinear remainder $Y$. In order to control the high$\times$high-interaction in \eqref{intro:eq-heat-high-high}, it is rewritten as 
\begin{equation}\label{intro:eq-heat-high-high-decomposed}
\shquadratic \parasim \Psi = \shquadratic \parasim X + \shquadratic \parasim Y.
\end{equation}
Since the sum of the regularities of $\shquadratic$ and $Y$ is positive, the $Y$-term in \eqref{intro:eq-heat-high-high-decomposed} is easy to estimate. While the regularity of $X$ is insufficient to control the $X$-term in \eqref{intro:eq-heat-high-high-decomposed}, \eqref{intro:eq-heat-X} yields detailed structural information on $X$, which is sufficient to control (a renormalization of) the  $X$-term in \eqref{intro:eq-heat-high-high-decomposed}. Ultimately, it is possible to solve the system of heat equations for $(X,Y)$, which yields the well-posedness of \eqref{intro:eq-heat} and the decomposition 
\begin{equation}\label{intro:eq-heat-description}
\Phi = \shlinear - \shcubic + X + Y. 
\end{equation}
For the convenience of the reader, we collected the terms in the treatment of \eqref{intro:eq-heat} in Figure \ref{figure:structure-heat}. \\

While the para-controlled approach was initially developed for parabolic equations, Gubinelli, Koch, and Oh \cite{GKO18} applied similar ideas to the three-dimensional quadratic wave equation with space-time white noise. On the surface, the random structure of the local dynamics in \cite{GKO18}, which consists of explicit stochastic objects and para-controlled terms, looks similar to \eqref{intro:eq-heat-description}. However, the required estimates in the dispersive setting are quite different from the estimates in the parabolic setting. For instance, the Schauder estimates in \cite{GIP15} are replaced by (multilinear) dispersive estimates in \cite{GKO18}, which are much more delicate. The para-controlled approach to wave equations was also used in \cite{B20II,OOT20,OOT21} and, as further discussed in Subsection \ref{section:intro-main-ideas} and Subsection \ref{section:ansatz-para}, will also be used in this article.

\subsubsection{Random tensors and random averaging operators}\label{section:intro-tensor}
The random tensor theory was developed in \cite{DNY20} by the second, third and fourth authors of this paper to understand at a granular level the propagation of randomness in the context of dispersive equations. This theory can be viewed as the dispersive counterpart of the existing parabolic theories  \cite{H14, GIP15, K16} mentioned above.  In \cite{DNY20} it was applied to nonlinear Schr\"odinger equations with power nonlinearities to obtain 
probabilistic local well-posedness in the full subcritical range relative to the probabilistic scaling as described in Subsection \ref{section:intro-scaling}. Furthermore, the solution constructed using this framework has an explicit expansion in the Fourier side in terms of multilinear Gaussians with random tensor coefficients.

The method of random averaging operators developed in \cite{DNY19} and mentioned in Subsection \ref{section:intro-other-settings} is a precursor of the random tensor theory. It controls to the first order the propagation of randomness beyond the linear evolution. The idea is to include the high-low interactions in the ansatz, but instead of putting them in a low regularity space in the ansatz, as it is done in the para-controlled calculus{\footnote{See Subsection \ref{section:intro-para-heat} for a more detailed explanation of \cite{GIP15} in the parabolic setting and of \cite{GKO18} in the wave setting.}},  one writes them as an operator applied to the high frequency linear evolution of the random initial data.
This operator, whose coefficients are independent with the modes of the linear evolution,  contains the randomness information of the low frequency components of the solution, which is carried by two operator norm estimates in an inductive argument. More precisely the random averaging operators take the form
\begin{equation}\label{random_matrix}
\mathcal P (u_{\mathrm{lin}}) = \sum_{N}\sum_{L\ll N} \mathcal{P}_{NL} (P_N u_{\mathrm{lin}})=\sum_{N} \mathcal{P}_{N} (P_N u_{\mathrm{lin}}) \,\text{ and }\,\bigl(\mathcal{P}_{N} (P_N u_{\mathrm{lin}})\bigr)_n = \sum_{\langle n'\rangle\sim N}\sum_{L\ll N} h^{NL}_{nn'}\, g_{n'},
\end{equation}
where $N$ and $L$ are the frequencies of the high and low inputs, the subscript $n$ means the $n$-th mode, the  $g_{n'}$ are i.i.d. Gaussian random variables and the Fourier kernel matrix $h^{NL}_{nn'}$ of $\mathcal{P}_{NL}$  is independent{\textbf{\footnote{Independence between the high frequencies of the linear evolution and the low frequencies of the nonlinear part was first used by Bringmann in \cite{B21}}}} with $P_N u_{\mathrm{lin}}$ and carries the random information of $P_{\leq L} u$.

To understand the random structures at a finer level, one is naturally lead to studying multilinear expressions of the form
\begin{equation}\label{random_tensor}
\sum_{\substack{n_1,n_2,\cdots,n_r\\\langle n_j\rangle\sim N_j, \,j=1,\cdots, r}}\,  h_{nn_1n_2\cdots n_r} \,\, \prod_{j=1}^r g_{n_j}, 
\end{equation}
and their associated random $(r,1)$-tensors $h =h_{nn_1n_2\cdots n_r}$, which depend on the low frequency  components $P_{\leq L} u$ of the solution and are independent with the product of Gaussians $\prod_{j=1}^r g_{n_j}$ that lives in the high frequencies (namely $N_j\,(1\leq j\leq r)$, which are all $\gg L$).
 The random tensor terms like \eqref{random_tensor} describe all possible mixed high-low interactions in the nonlinear evolution of the random data. The final random tensor Ansatz of the solution contains a finite expansion of such random tensor terms and a smooth remainder. The regularity of the remainder can be determined by the order of the expansion, which is possibly large but finite given any fixed subcritical equation.
 The key in propagating this explicit random structure under the flow is to control the appropriate operator and $\ell^2$ norms of the random tensor $h$ in the Ansatz by an inductive argument.  See \cite{DNY20} for details.
 
We remark that, in this paper, it is possible to use the random tensor theory in the proof of local well-posedness. However, due to the subcritical nature of (\ref{intro:eq-NLW}), plus the smoothing effect of the wave propagator,  the various terms we encounter here will have better regularity than in \cite{DNY20}, hence
the full power of the random tensor \emph{Ansatz} is not needed for the current (hyperbolic $\Phi_3^4$) problem. In fact, trying to introduce the inductive Ansatz of \cite{DNY20} will make the proof conceptually more intricate by unraveling the finer random structure of the remainder (yielding also better estimates), but such fine structure is not necessary for our main purpose, so we did not choose this option here. On the other hand, our proof does rely heavily on the tensor \emph{norms and estimates} developed in \cite{DNY20}, see Subsection \ref{section:intro-main-structure} for further discussions.

\subsubsection{Globalization}\label{section:intro-earlier-global}

While the methods discussed in Subsections \ref{section:intro-para-heat}-\ref{section:intro-tensor} yield detailed structural information on the corresponding local dynamics, we are ultimately interested in proving both global well-posedness and invariance (as in Theorem \ref{intro:thm-formal}), which requires additional ingredients. 

 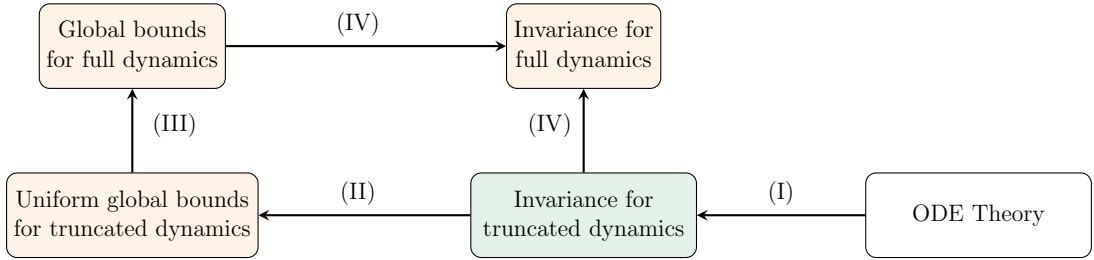
\begin{figure}
\begin{tikzpicture}[node distance=2cm,scale=0.75, every node/.style={scale=0.75}]
\node (ODEtheory)[wrect] at (7,0) 
{\large ODE Theory};
\node (InvTruncated)[grect]  at (0,0) 
{\large Invariance for \\ \large truncated dynamics};
\node (GBTruncated)[orect]  at (-8,0) 
{\large Uniform global bounds \\ \large for truncated dynamics};
\node (GBFull)[orect] at (-8,3) 
{\large Global bounds \\ \large for full dynamics};
\node (InvFull)[orect] at (0,3) 
{\large Invariance for \\ \large full dynamics};
\draw[arrow] (ODEtheory) to (InvTruncated);
\draw[arrow] (InvTruncated) to (GBTruncated);
\draw[arrow] (GBTruncated) to (GBFull);
\draw[arrow] (GBFull) to (InvFull);
\draw[arrow] (InvTruncated) to (InvFull);
\node[align=center] at (3.5,0.4) {\large \ref{intro:item-global-1}};
\node[align=center] at (-4,0.4) {\large \ref{intro:item-global-2}};
\node[align=center] at (-7.25,1.6) {\large \ref{intro:item-global-3}};
\node[align=center] at (-4,3.4) {\large \ref{intro:item-global-4}};
\node[align=center] at (-0.6,1.6) {\large \ref{intro:item-global-4}};
\end{tikzpicture}
\caption{\small{Globalization}}
\label{figure:globalization}
\end{figure}

The general strategy for proving global well-posedness and invariance, which is known as Bourgain's globalization argument \cite{B94}, is illustrated in Figure \ref{figure:globalization} and consists of the following four steps: 
\begin{enumerate}[label=(\Roman*)]
    \item\label{intro:item-global-1} First, one defines a truncated Hamiltonian, the corresponding truncated evolution equation, and the corresponding truncated Gibbs measure. Due to ODE-methods, the truncated Gibbs measure is invariant under the truncated evolution equation.
    \item\label{intro:item-global-2} Second, one uses the invariance of the truncated Gibbs measure to obtain global bounds for the truncated dynamics. In this step, it is essential that the global bounds are uniform in the truncation parameter. Furthermore, the global bounds have to control not just the solution itself, but also the individual components of the Ansatz (such as $X$ in \eqref{intro:eq-heat-X}). 
    \item\label{intro:item-global-3} Third, one proves the global well-posedness of the full evolution equation using the global bounds from Step \ref{intro:item-global-2} and stability arguments. 
    \item\label{intro:item-global-4} Finally, one proves the invariance of the Gibbs measure under the full evolution equation using the global well-posedness from Step \ref{intro:item-global-3} and the invariance for the truncated evolution equation. 
\end{enumerate}

In many implementations of Bourgain's globalization argument, the most difficult aspects concern Step \ref{intro:item-global-2}. The difficulties are most prominent when the Gibbs measure is singular with respect to the Gaussian free field, in which case one often works with a representation of the Gibbs measure which is not preserved by the dynamics. The difficulties for singular Gibbs measures were first addressed by the first author in \cite{B20II} and later by Oh, Okamoto, and Tolomeo in \cite{OOT21} and will be further discussed in Subsection \ref{section:intro-main-globalization} below.

\subsection{Setting and main result} \label{section:intro-setting}
After the discussion of the earlier literature, we now proceed with a more detailed discussion of our main theorem. In order to make rigorous statements, we first introduce frequency-truncated versions of the Hamiltonian \eqref{intro:eq-Hamiltonian}, the nonlinear wave equation \eqref{intro:eq-NLW}, and the Gibbs measure \eqref{intro:eq-Gibbs}. For any dyadic frequency-scale $N$, we define 
\begin{equation}\label{intro:eq-Hamiltonian-N}
H_{\leq N}[u,\partial_t u] = \int_{\T^3} \dx \bigg( \frac{u^2}{2} + \frac{|\nabla u|^2}{2} + \frac{(\partial_t u )^2}{2} \bigg) +  \int_{\T^3} \dx \bigg( \frac{\lcol (P_{\leq N} u)^4 \rcol}{4} + \gamma_{\leq N} \cdot \frac{(P_{\leq N} u)^2}{2} \bigg).
\end{equation}
Here, $(P_{\leq N})_N$ is a sequence of sharp Fourier-truncation operators which are defined in \eqref{intro:eq-frequency-projection}. The quartic term $\lcol (P_{\leq N} u)^4 \rcol$  is the frequency-truncated Wick-ordered quartic power which is as in \eqref{intro:eq-Wick-4}. Finally, the renormalization constant $\gamma_{\leq N}$ is defined as in Definition \ref{sec6:def:Gamma} and diverges logarithmically. Together with the symplectic form \eqref{intro:eq-symplectic}, the frequency-truncated Hamiltonian $H_{\leq N}$ induces a frequency-truncated cubic wave equation, which is given by 
\begin{equation}\label{intro:eq-NLW-N}
\begin{cases}
\big( \partial_t^2 + 1 - \Delta \big) u_{\leq N} = - P_{\leq N} \Big[ \lcol (P_{\leq N} u_{\leq N})^3 \rcol +\, \gamma_{\leq N} \cdot P_{\leq N} u_{\leq N} \Big] \qquad (t,x) \in \mathbb{R} \times \mathbb{T}^3, \\
(u_{\leq N}, \langle \nabla \rangle^{-1} \partial_t u_{\leq N})\big|_{t=0} = (\phi^{\cos},\phi^{\sin}). 
\end{cases}
\end{equation}
In order to define the frequency-truncated Gibbs measure, we first define the three-dimensional Gaussian free field $\mathscr{g}$. At a formal level, it can be written as 
\begin{equation}\label{intro:eq-GFF-formal}
``\mathrm{d}\mathscr{g}(\phi)= \mathcal{Z}^{-1} \exp\bigg( - \int_{\T^3} \dx  \bigg( \frac{|\phi|^2}{2} + \frac{|\nabla \phi|^2}{2}  \bigg) \bigg) \mathrm{d}\phi",  
\end{equation}
which  is strongly related to \eqref{intro:eq-Phi4d}. At a rigorously level, it is defined by 
\begin{equation}\label{intro:eq-GFF}
\mathscr{g} = \operatorname{Law}\Big( \sum_{n\in \Z^3} \frac{g_n}{\langle n \rangle} e^{inx} \Big),
\end{equation}
where $g_0$ is a real-valued standard Gaussian and $(g_n)_{n\in \Z^3 \backslash \{0 \}}$ is  a sequence of complex-valued standard Gaussians such that $g_{n_1}$ and $g_{n_2}$ are independent if $n_1 \neq \pm n_2$ and $g_{-n}=\overline{g_n}$. Equipped with \eqref{intro:eq-GFF}, we define the frequency-truncated Gibbs measure\footnote{Due to the $\langle \nabla \rangle^{-1}$-multiplier in the initial condition for $\partial_t u_{\leq N}$ in \eqref{intro:eq-NLW-N}, the Gibbs measure in \eqref{intro:eq-Gibbs-N} contains the product $\mathscr{g}\otimes \mathscr{g}$ and not $\mathscr{g}\otimes (\langle \nabla \rangle_{\#} \mathscr{g})$, where $\langle \nabla \rangle_{\#}$ denotes the push-forward under $\langle \nabla \rangle$.} by 
\begin{equation}\label{intro:eq-Gibbs-N}
\mathrm{d}\mu_{\leq N}(\phi^{\cos},\phi^{\sin}) = \mathcal{Z}^{-1} \exp\bigg( - \int_{\T^3} \dx \bigg( \frac{\lcol (P_{\leq N} \phi^{\cos})^4 \rcol}{4} + \gamma_{\leq N} \cdot \frac{(P_{\leq N} \phi^{\cos})^2}{2} \bigg) \bigg) \mathrm{d}\big( \mathscr{g} \otimes  \mathscr{g} \big)(\phi^{\cos},\phi^{\sin}). 
\end{equation}
In  \cite{GJ73,BG18}, it has been proven\footnote{While \cite{BG18} does not consider sharp cutoffs, the argument only requires minor modifications (see also  Appendix \ref{section:heat-appendix}).} that $(\mu_{\leq N})_{N}$ weakly converges to a unique limit, which we define as
\begin{equation}\label{intro:eq-Gibbs-limit}
\mu = \operatorname{w-lim}\displaylimits_{N\rightarrow \infty} \mu_{\leq N}.
\end{equation}
Equipped with \eqref{intro:eq-NLW-N}, \eqref{intro:eq-Gibbs-N}, and \eqref{intro:eq-Gibbs-limit}, we can now state a rigorous version of our main result. In the statement, $\mathscr{H}_x^{-1/2-\epsilon}$ refers to the Cartesian product of Sobolev spaces at regularity $-1/2-\epsilon$, see \eqref{prelim:eq-cartesian-sobolev} below.

\begin{theorem}[Global well-posedness and invariance, rigorous version]\label{intro:thm-rigorous}
For any frequency-scale $N\geq 1$ and $(\phi^{\cos},\phi^{\sin}) \in \mathscr{H}_x^{-1/2-\epsilon}(\T^3)$, let $u_{\leq N}$ be the solution of the frequency-truncated cubic wave equation \eqref{intro:eq-NLW-N} with initial data $u_{\leq N}[0]=(\phi^{\cos},\phi^{\sin})$. In addition, let $\mu$ be the Gibbs measure from \eqref{intro:eq-Gibbs-limit}. Then, for $\mu$-almost every $(\phi^{\cos},\phi^{\sin})$ and all $T\geq 1$, the limiting dynamics 
\begin{equation}\label{intro:eq-limiting-dynamics}
u[t] = \lim_{N\rightarrow \infty} u_{\leq N}[t]
\end{equation}
exists in $C_t^0 \mathscr{H}_x^{-1/2-\epsilon}([-T,T]\times \T^3)$. Furthermore, the Gibbs measure is invariant under the limiting dynamics, i.e., 
\begin{equation}\label{intro:eq-invariance}
\Law_\mu \big( u[t] \big) = \mu
\end{equation}
for all $t\in \R$. 
\end{theorem}

In addition to the invariance of the Gibbs measure, our argument further yields detailed information on the random structure of the limiting dynamics $u$ in space-time. For a detailed discussion, we refer to Subsection \ref{section:intro-main-structure} below. 
We now end this section with three comments on our main theorem. 

\begin{enumerate}
\item Since the infinite-volume limit of the Gibbs measure has been constructed in \cite{FO76} (see also \cite{AK21,GH21,MW20}) and the wave equation exhibits finite speed of propagation, Theorem \ref{intro:thm-rigorous} can likely be extended from the periodic to the Euclidean setting.
\item In Theorem \ref{intro:thm-rigorous}, we do not obtain the flow property for the limiting dynamics. We nevertheless believe that, after a modification of the stability theory in Section \ref{section:global}, our estimates should be strong enough to yield the flow property. However, due to the length of this article, we did not pursue this question in detail. 
\item While Theorem \ref{intro:thm-rigorous} concerns only the Gibbs measure, our argument also yields the probabilistic local well-posedness (see Section \ref{section:ansatz}) for the Gaussian free field. Therefore, this article further improves earlier results on the probabilistic local well-posedness of the cubic nonlinear wave equation with Gaussian initial data \cite{OPT19,OWZ21}.
\end{enumerate}

\subsection{Main ideas} \label{section:intro-main-ideas}

We now discuss the main ingredients and novelties in the proof of Theorem \ref{intro:thm-rigorous}. For expository purposes, we split the discussion into three separate parts: 
\begin{enumerate}[label=(\Roman*)]
\item representation of the Gibbs measure,
\item random structure of the dynamics, 
\item and globalization. 
\end{enumerate}

\subsubsection{Representation of the Gibbs measure}\label{section:intro-main-gibbs}

We now describe the properties of the Gibbs measures $(\mu_{\leq N})_N$ from \eqref{intro:eq-Gibbs-N} which are most relevant for the proof of Theorem \ref{intro:thm-rigorous}. There are two different properties of $(\mu_{\leq N})_N$ which create challenges in the well-posedness theory of the cubic wave equation: 
\begin{enumerate}[label=(\alph*)]
    \item  \label{intro:item-Gibbs-challenge-1} The spatial regularity of typical samples is $-1/2-$. 
    \item \label{intro:item-Gibbs-challenge-2} The Fourier coefficients of the samples are probabilistically dependent. 
\end{enumerate}
In lower-dimensional models such as in \cite{B94,B96,DNY19}, the Gibbs measure is absolutely continuous with respect to the Gaussian free field, and the probabilistic dependencies in \ref{intro:item-Gibbs-challenge-2} can mostly be avoided. However, just as in \cite{B20II,OOT21}, the limiting Gibbs measure in our setting is singular with respect to the Gaussian free field, and hence \ref{intro:item-Gibbs-challenge-2} is an important aspect of our problem. 

\begin{figure}
\begin{tikzpicture}[node distance=2cm,scale=0.75, every node/.style={scale=0.75}]

\node (sample)[rect]  at (-5,0) {\Large $\phi^{\cos}$};
\node (heat)[rect]  at (0,0) {\Large \, $\Phi_{\leq N}^{\cos}\colon$ \\[1ex] \Large $\Phi_{\leq N}^{\cos}(s_0) = \phi^{\cos}$};
\node (wave)[rect]  at (6.5,0) {\\[-1ex] \Large \, $u_{\leq N}\colon$ \\[1ex] \Large $  u_{\leq N}(0) = \Phi_{\leq N}^{\cos}(0)$ };

\draw[arrow] (sample) to (heat);
\node[align=center] at (-2.75,0.4) {\eqref{intro:eq-heat-N}};
\draw[arrow] (heat) to (wave);
\node[align=center] at (2.75,0.4) {\eqref{intro:eq-NLW-N}};
\end{tikzpicture}
\caption{\small{The caloric initial data.}}
\label{figure:caloric}
\end{figure}
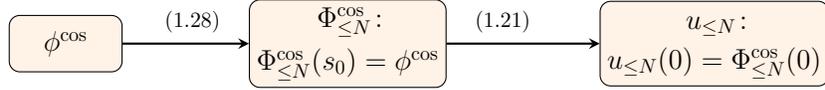

While the challenges resulting from \ref{intro:item-Gibbs-challenge-1} and \ref{intro:item-Gibbs-challenge-2} cannot be avoided, it is possible to separate them. To this end, we do not directly solve the cubic wave equation \eqref{intro:eq-NLW-N}, but instead proceed as illustrated in Figure \ref{figure:caloric}. 
 We start with samples $(\phi^{\cos},\phi^{\sin})$ from the truncated Gibbs measure $\mu_{\leq N}$. Then, we solve the frequency-truncated cubic stochastic heat equation 
\begin{equation}\label{intro:eq-heat-N}
\begin{cases}
\big( \partial_s + 1 - \Delta \big) \Phi_{\leq N}^{\cos}= - P_{\leq N} \Big[ \lcol (P_{\leq N} \Phi_{\leq N}^{\cos})^3 \rcol + \gamma_{\leq N} P_{\leq N} \Phi_{\leq N}^{\cos} \Big] + \sqrt{2} \mathrm{d}W^{\cos} \qquad (s,x) \in (s_0,\infty) \times \T^3, \\
\Phi_{\leq N}^{\cos}\big|_{s=s_0}=\phi^{\cos}. 
\end{cases}
\end{equation}
The choice of the initial time $s_0 <0$ is discussed in Subsection \ref{section:diagram-parabolic}, but we encourage the reader to simply think of $s_0=-1$. Then, we solve the frequency-truncated cubic wave equation \eqref{intro:eq-NLW-N}, but with initial data 
\begin{equation}\label{intro:eq-initial}
\big(u_{\leq N}, \langle \nabla \rangle^{-1} \partial_t u_{\leq N}\big)\Big|_{t=0} = \big( \Phi_{\leq N}^{\cos}(0), \phi^{\sin} \big). 
\end{equation}
Due to the invariance of the Gibbs measure under the Langevin dynamics, i.e., \eqref{intro:eq-heat-N}, the law of $(\Phi_{\leq N}^{\cos}(0),\phi^{\sin})$ is still given by $\mu_{\leq N}$. At a heuristic level, one can think of the replacement of $(\phi^{\cos},\phi^{\sin})$ by $(\Phi_{\leq N}^{\cos}(0),\phi^{\sin})$ as a (highly non-trivial) change of coordinates. The advantage of the initial data in \eqref{intro:eq-initial} is that, using the para-controlled approach discussed in Subsection \ref{section:intro-para-heat}, we have a detailed description of the random structure of $\Phi(0)$. More precisely, using similar notation as in Subsection \ref{section:intro-para-heat}, we write
\begin{equation}\label{intro:eq-initial-caloric}
    (\Phi_{\leq N}^{\cos}(0),\phi^{\sin}) = \initial{blue} - \initial{green}[][\leqN] + \initial{red}[][\leqN],
\end{equation}

where\footnote{The cosine-superscripts in the stochastic diagrams \scalebox{0.9}{$\shlinear[][\cos]$} and \scalebox{0.8}{$\shcubic[\leqN][\cos]$} emphasize that the stochastic forcing in their definitions is different from the stochastic forcing in the definition of the initial velocity, see Subsection \ref{section:diagram-parabolic}.}

\begin{itemize}
    \item $\initial{blue}=\Big(\shlinear[][\cos](0),\phi^{\sin}\Big)$ is Gaussian random data with spatial regularity $-1/2-$,
    \item $\initial{green}[][\leqN]=\Big( \shcubic[\leqN][\cos](0),0 \Big)$ is a cubic Gaussian chaos with spatial regularity $1/2-$, 
    \item and $\initial{red}[][\leqN]=\Big( \Psi_{\leq N}^{\cos}(0), 0 \Big)$ is a nonlinear remainder with spatial regularity $1-$. 
\end{itemize}

As will be explained below the remark, we refer to $\initial{blue}$, $\initial{green}[][\leqN]$, and $\initial{red}[][\leqN]$ as the caloric initial data.

\begin{remark}[Colors and shapes of the caloric initial data]
Due to the listed spatial regularities,  the reader may wish to think of the blue caloric initial data $\initial{blue}$, the green caloric initial data $\initial{green}[][\leqN]$, and the red caloric initial data $\initial{red}[][\leqN]$ as high, medium, and low-frequency functions, respectively. The choice of colors then matches the visible spectrum of light, since blue, green, and red light also corresponds to high, medium, and low-frequencies.  

The different shapes of the caloric initial data, i.e., the circle of the blue caloric initial data $\initial{blue}$, the diamond of the green caloric initial data $\initial{green}[][\leqN]$, and the pentagon of the red caloric initial data $\initial{red}[][\leqN]$ have no special meaning. The reason for using different shapes is to make the paper accessible to colorblind researchers and readable in black and white print. 

We note that, just as the blue caloric initial data $\initial{blue}$, the stochastic objects of the heat equation, such as $\shlinear[]$ and $\scalebox{0.8}{\shcubic[\leqN][\cos]}$, also contain a circle. However, since the blue caloric initial data only appears next to full arrows (see e.g. \eqref{intro:eq-NLW-Ansatz} below), whereas the stochastic diagrams of the heat equation contain dotted arrows, this should not cause any confusion. 
\end{remark}

The motivation behind \eqref{intro:eq-initial-caloric} is similar to the motivation behind the caloric gauge. The caloric gauge, which was introduced by Tao in \cite{Tao04}, is used in the analysis of geometric wave equations and constructed using geometric heat equations. In comparison, \eqref{intro:eq-initial-caloric} yields a representation of the initial data for a random wave equation using a stochastic heat equation. Since gauge transformations and representations of measures can be viewed as coordinates, there is a strong connection between both ideas. This is our reason for referring to $\initial{blue}$, $\initial{green}[][\leqN]$, and $\initial{red}[][\leqN]$ as the caloric initial data. 

We also remark that a similar representation as in \eqref{intro:eq-initial-caloric} was used in a note of Hairer \cite{HairerNote} for a proof of the singularity of the $\Phi^4_3$-measure. 

\begin{remark} Instead of directly using the cubic stochastic heat equation (as in Figure \ref{figure:caloric}), the earlier works \cite{B20II} and \cite{OOT21} relied on the drift measures from \cite{BG20} and the variational method from \cite{BG20}, respectively. The present approach offers advantages in the globalization argument, which will be discussed in Subsection \ref{section:intro-main-globalization} and Section \ref{section:global}.
\end{remark}

\subsubsection{Random structure of solutions to the cubic wave equation} \label{section:intro-main-structure}

Motivated by the representation of the Gibbs measure from \eqref{intro:eq-initial-caloric}, we now consider the caloric initial value problem 
\begin{equation}\label{intro:eq-NLW-caloric}
\begin{cases}
\big( \partial_t^2 + 1 - \Delta \big) u_{\leq N} = - P_{\leq N} \Big[ \lcol (P_{\leq N} u_{\leq N})^3 \rcol +\, \gamma_{\leq N} \cdot P_{\leq N} u_{\leq N} \Big] \qquad (t,x) \in \mathbb{R} \times \mathbb{T}^3, \\
(u_{\leq N}, \langle \nabla \rangle^{-1} \partial_t u_{\leq N})\big|_{t=0} = \initial{blue}-\initial{green}[][\leqN] + \initial{red}[][\leqN]. 
\end{cases}
\end{equation}
In order to prove Theorem \ref{intro:thm-rigorous}, we require detailed information on the random structure of the solution $u_{\leq N}$ of \eqref{intro:eq-NLW-caloric}. This random structure is expressed through the Ansatz 
\begin{equation}\label{intro:eq-NLW-Ansatz}
u_{\leq N}= \slinear - \scubic[\leqN] - \slinear[green][\leqN] + 3 \squintic[\leqN] + X^{(1)}_{\leq N} + X^{(2)}_{\leq N} + Y. 
\end{equation}
In this introduction, we only give a brief description of the terms in \eqref{intro:eq-NLW-Ansatz} and postpone a detailed (and self-contained) description until Section \ref{section:ansatz} below. The first four summands in \eqref{intro:eq-NLW-Ansatz} are explicit stochastic objects, which are defined as follows:
\begin{itemize}
    \item $\slinear$ is the linear evolution of $\initial{blue}$ and has spatial regularity $-1/2-$, 
    \item $\scubic[\leqN]$ is the cubic Picard iterate originating from $\initial{blue}$ and has spatial regularity $0-$, 
    \item $\slinear[green][\leqN]$ is the linear evolution of $\initial{green}[][\leqN]$ and has spatial regularity $1/2-$, 
    \item and $\squintic[\leqN]$ is the quintic Picard iterate originating from $\initial{blue}$ and has spatial regularity $1/2-$. 
\end{itemize}
The fifth and sixth terms in \eqref{intro:eq-NLW-Ansatz} are para-controlled components, which are morally\footnote{For the precise definition of $X^{(1)}_{\leq N}$ and $X^{(2)}_{\leq N}$, which contain additional terms, we refer to \eqref{ansatz:eq-X11-a}-\eqref{ansatz:eq-X11-d} and \eqref{ansatz:eq-X21} below.} defined as follows: 
\begin{itemize}
    \item $X^{(1)}_{\leq N}$ contains high$\times$low$\times$low-interactions between $\slinear$ and two other components of $\eqref{intro:eq-NLW-Ansatz}$ and has spatial regularity $1/2-$, 
    \item and $X^{(2)}_{\leq N}$ contains high$\times$high$\times$low-interactions between $\slinear$, $\slinear$, and one other component of \eqref{intro:eq-NLW-Ansatz} and has spatial regularity $1/2-$. 
\end{itemize}
Finally, $Y$ is a nonlinear remainder, which contains all remaining interactions and has spatial regularity greater than $1/2$. In Figure \ref{figure:structure-NLW}, we list the terms from \eqref{intro:eq-NLW-Ansatz} together with their spatial regularities. For comparison, Figure \ref{figure:structure-NLW} also contains lists of the corresponding terms in \cite{GKO18} and \cite{B20II}. \\

\begin{figure}
\begin{subfigure}{\textwidth}
\begin{center}
\begin{tabular}{|P{2cm}|P{1.25cm}|P{1.25cm}|P{1.25cm}|P{1.25cm}|P{1.25cm}|P{1.25cm}|P{1.25cm}|} \hline \rule{0pt}{14pt} 
Regularity &$-1/2-$ & $0-$ & \multicolumn{4}{c|}{$1/2-$}  & $>1/2$ \\[4pt]  \cline{1-8} \rule{0pt}{14pt} 
Object & $\slinear[blue]$ & $\scubic$ & $\slinear[green]$ & $\squintic$ & $X^{(1)}$ & $X^{(2)}$ & Y \\[4pt]  \cline{1-8} 
\end{tabular}
\end{center}
\subcaption{\small{The random structure of solutions to the cubic wave equation in this article.}}
\label{figure:structure-cubic-NLW}
\end{subfigure}
\begin{subfigure}{0.8\textwidth}~\\[-1ex]
\begin{center}
\begin{tabular}{|P{2cm}|P{1.25cm}|P{1.25cm}|P{1.25cm}|P{1.25cm}|} \hline \rule{0pt}{14pt}
Regularity &$-1/2-$ & $\beta-$  & $1/2-$  & $>1/2$ \\[4pt]  \cline{1-5} \rule{0pt}{14pt} 
Object & $\slinear$ & $\Bricubic$ & $X^{(1)}$ & $Y$ \\[4pt]  \cline{1-5} 
\end{tabular}
\end{center}
\subcaption{\small{The random structure of solutions to the Hartree-nonlinear wave equation  with interaction potential $V(x)\simeq |x|^{-(3-\beta)}$, where $\beta>0$, in \cite{B20II}. }}
\end{subfigure}
\begin{subfigure}{0.8\textwidth}~\\[-1ex]
\begin{center}
\begin{tabular}{|P{2cm}|P{1.25cm}|P{1.25cm}|P{1.25cm}|P{1.25cm}|} \hline \rule{0pt}{14pt}
Regularity &$-1/2-$ & \multicolumn{2}{c|}{$1/2-$}  & $>1/2$ \\[4pt]  \cline{1-5} \rule{0pt}{14pt} 
Object & $\GKOlinear$ & $\GKOquadratic$ & $X^{(1)}$ & $Y$ \\[4pt]  \cline{1-5} 
\end{tabular}
\end{center}
\subcaption{\small{The random structure of solutions to the quadratic wave equation in \cite{GKO18}.}}
\end{subfigure}
\caption{\small{The random structure of solutions to the cubic wave equation in this article, the Hartree-nonlinear wave equation with interaction potential $V(x)\simeq |x|^{-(3-\beta)}$ in \cite{B20II}, and the quadratic wave equation in \cite{GKO18}.}}
\label{figure:structure-NLW}
\end{figure}

In the following, we discuss the main novelties of our analysis of the caloric initial value problem \eqref{intro:eq-NLW-caloric} and the random structure \eqref{intro:eq-NLW-Ansatz}. \\

\emph{Lattice counting, tensors, and molecule estimates for wave equations:}
As discussed in Subsection \ref{section:intro-tensor}, our proof relies heavily on the framework of tensor norms and estimates developed in \cite{DNY20}. Most of the tensors used here are explicitly defined, and can be constructed starting from the \emph{base tensor} $h^b$, which is basically
\[(h^b)_{nn_1n_2n_3}=\mathbf{1}\big\{n=n_1+n_2+n_3,\langle n\rangle=\langle n_1\rangle\pm\langle n_2\rangle\pm\langle n_3\rangle+O(1)\big\},\] see (\ref{base-tensor}) for the precise definition. A large part of the proof consists of estimating various operator and $\ell^2$ norms of these tensors, which will rely on the merging and moment estimates from \cite{DNY20}, see Lemma \ref{counting:lem-merging} and Proposition \ref{counting:prop-moment}, as well as corresponding estimates for the base tensor $h^b$. By Schur's lemma, these estimates for $h^b$ then reduce to certain integer vector counting estimates, which are proved by suitable geometric arguments. See Section \ref{section:counting} for details.

In later parts of the proof, we will also encounter higher order (septic and nonic) stochastic objects and need to estimate their regularity. In principle this could be done using the tensor estimates as described above, but due to the high order (and rich structure) of these objects, such an argument would lead to an enumeration of  about $40$ different cases, which is way too tedious. Instead, we will adopt the notion of \emph{molecules} introduced by the second author and Hani in \cite{DH21}, which wraps up the key information and omits the irrelevant properties of these objects. By doing so we can reduce the number of cases to less than $10$, and then proceed with an algorithmic-enumeration hybrid method. See Section \ref{section:analytic2} for details.

We remark that, while \cite{DNY20} provides uniform tensor and counting estimates for arbitrarily high order terms in the case of \emph{Schr\"{o}dinger} equations, things are not quite the same for \emph{wave} equation. On the one hand, the wave case is easier with extra smoothing and larger subcritical range; on the other hand, we currently do not have the uniform estimates as in \cite{DNY20}. This may be related to the different behavior of dispersion relations (such as Galilean vs. Lorentzian symmetry). Another possible reason is that, due to the smoothing effect for wave, the high-high-to-high and high-high-to-low interactions behave differently, which is not the case of Schr\"{o}dinger. As a result, there are certain terms in the expansion that need special treatments (i.e. those leading to renormalizations) while in the Schr\"{o}dinger case all terms are treated equally and satisfy the same uniform bound. Nevertheless, we believe that it is possible to establish a wave (or heat) version of the theory in \cite{DNY20} that also takes into account renormalizations, which would be needed in the further study of random data wave equations.
\\

\emph{Bilinear operator estimates:} As said above, most parts of our proof can be done by utilizing the (linear) operator norms of tensors in the same way as in \cite{DNY20}. However, there is one place where we need to estimate the \emph{bilinear operator norm} of an explicitly defined tensor, say of form
\[\|G\|_{k_1\times k_2\to k}:=\sup\bigg\{\bigg\|\sum_{k_1,k_2}G_{kk_1k_2}u_{k_1}v_{k_2}\bigg\|_{\ell_k^2}:\|u\|_{\ell_{k_1}^2}=\|v\|_{\ell_{k_2}^2}=1\bigg\}\,\,\,\,\mathrm{where}\,\,\,\, G_{kk_1k_2}=\sum_{k_3}(h^b)_{kk_1k_2k_3}g_{k_3}\] with $\{g_{k_3}\}$ being independent Gaussians, see Lemma \ref{bilinear:lem-abstract}.

Such bilinear norms are not studied in \cite{DNY20}. One may control them by the corresponding linear operator norms by ``concatenating" some of the variables (for example we have $\|G\|_{k_1\times k_2\to k}\leq\|G\|_{k_1k_2\to k}$ if one views $(k_1,k_2)$ as one big variable and $u\otimes v$ a function of that variable), but the resulting estimates from \cite{DNY20}, despite being sharp for the \emph{linear} norm, will not suffice for the original \emph{bilinear} norm. 

Instead, we will exploit the bilinear nature of the norm by estimating it directly, using a variation of the ideas from \cite{DNY20}, to obtain improvements upon the linear estimates of \cite{DNY20}, which are enough for our purposes. Such probabilistic estimates go \emph{beyond} the regime of random matrices, and have not appeared in any earlier works to the best of our knowledge. See Section \ref{section:bilinear} for details.
\\

\emph{The $1533$-cancellation:} In addition to the counting, tensor, and molecule estimates, which control individual stochastic objects and random operators, we also rely on a delicate hidden cancellation involving multiple stochastic objects. In our analysis of \eqref{intro:eq-NLW-caloric}, we naturally encounter the square of the cubic stochastic object, i.e.,
\begin{equation}
\label{intro:eq-NLW-33}
\Big( \scubic[\leqN] \Big)^2(t,x). 
\end{equation}
In Lemma \ref{diagram:lem-C33} below, we prove that \eqref{intro:eq-NLW-33} logarithmically diverges in the sense of space-time distributions as $N\rightarrow \infty$. Since \eqref{intro:eq-NLW-33} contains two cubic stochastic objects, we refer to the divergence of \eqref{intro:eq-NLW-33} as the $33$-divergence. The $33$-divergence cannot be removed by a further renormalization, since the renormalization is already dictated by the construction of the Gibbs measure, which does not encounter a similar problem\footnote{A quick way to see this is by looking at the cubic stochastic object $\shcubic$ from the stochastic heat equation, which has spatial regularity $1/2-$ and whose square is therefore well-defined.}. Fortunately, in spite of the $33$-divergence, the linear combination
\begin{equation}
\label{intro:eq-NLW-1533} 
6 \, \slinear[blue][\leqN](t,x) \squintic[\leqN](t,x) + \Big( \scubic[\leqN] \Big)^2(t,x)
\end{equation}
has a well-defined limit as $N\rightarrow \infty$. Since \eqref{intro:eq-NLW-1533} contains one linear, one quintic, and two cubic stochastic objects, we refer to the cancellation in \eqref{intro:eq-NLW-1533} as the \oftt-cancellation. In the nonlinearity of \eqref{intro:eq-NLW-caloric}, it is possible to always group the stochastic objects together as in \eqref{intro:eq-NLW-1533}, and therefore the $33$-divergence does not prevent the well-posedness of \eqref{intro:eq-NLW-caloric}. \\

\emph{The heat-wave stochastic objects:} In Subsection \ref{section:intro-main-gibbs}, we discussed the caloric initial data $\initial{blue}$, $\initial{green}[][\leqN]$, and $\initial{red}[][\leqN]$, which is constructed using the cubic stochastic heat equation \eqref{intro:eq-heat-N}. In interactions involving only the blue caloric initial data $\initial{blue}$, the origin of the caloric initial data is not important and it is possible to work with any representation of the Gaussian free field (such as in \eqref{intro:eq-GFF}). However, in interactions involving both the blue and green caloric initial data $\initial{blue}$ and $\initial{green}[][\leqN]$,  the origin of the caloric initial data is essential. In our analysis of such interactions, we utilize long-hand stochastic diagrams, which are more detailed than the short-hand diagrams in  \eqref{intro:eq-NLW-Ansatz}. For example, we write 
\begin{equation}\label{intro:eq-NLW-longhand}
P_{\leq N} \Big[ \lcol \big( P_{\leq N} \slinear \big)^2 \rcol  \, P_{\leq N} \slinear[green][\leqN] \Big] = \scalebox{0.875}{\quinticmixednl[\leqN][][0]} \hspace{-2ex}+ 6 \hspace{-4ex} \scalebox{0.875}{\quinticmixednl[\leqN][][1]} + 6  \hspace{-5ex} \scalebox{0.875}{\quinticmixednl[\leqN][][2]}.
\end{equation}
The long-hand diagrams on the right-hand side of \eqref{intro:eq-NLW-longhand} will be defined in Section \ref{section:diagrams}. At this point, we only emphasize that the long-hand diagrams contain two different types of arrows. While the dotted arrows represent the heat propagator, the full arrows represent the wave propagator.

We note that this aspect of our argument was not present in the earlier works \cite{B20II,OOT21}, since the corresponding random shifts of the Gaussian random data have spatial regularity greater than $1/2$.

\subsubsection{Globalization}\label{section:intro-main-globalization}
In Subsection \ref{section:intro-earlier-global}, we discussed Bourgain's globalization argument \cite{B94}, which has been a successful paradigm for proving both probabilistic global well-posedness and the invariance of Gibbs measures. In this subsection, we focus only on aspects specific to our article and the earlier works \cite{B20II,OOT21}. \\

We first focus our discussion on our proof of global bounds, which correspond to 
Step \ref{intro:item-global-2} in Figure \ref{figure:globalization}. As discussed in Subsection \ref{section:intro-earlier-global}, our global bounds have to control the individual components in our Ansatz \eqref{intro:eq-NLW-Ansatz}. As a result, it is natural to try to simply iterate\footnote{In each step of this iteration, one has to utilize the invariance of the Gibbs measure $\mu_{\leq N}$ under the dynamics of \eqref{intro:eq-NLW-N}, which serves as a substitute for a conservation law. However, this aspect is not the main subject of our discussion.}
our local theory. 

The implementation of this iteration, however, presents a challenge. To illustrate this, let $\tau>0$ be a small time. Using our local theory, it follows that $u_{\leq N}$ exhibits the random structure \eqref{intro:eq-NLW-Ansatz} on the time-interval $[0,\tau]$. In particular, it follows that
\begin{equation}\label{intro:eq-global-tau}
u_{\leq N}[\tau]= \bigg( \slinear - \scubic[\leqN] - \slinear[green][\leqN] + 3 \squintic[\leqN] + X^{(1)}_{\leq N} + X^{(2)}_{\leq N} + Y \bigg)[\tau]. 
\end{equation}
Since the right-hand side of \eqref{intro:eq-global-tau} looks very different from the caloric representation of the initial data in \eqref{intro:eq-NLW-caloric}, it is difficult to directly establish the desired random structure on the time-interval $[\tau,2\tau]$. \\

In \cite{B20II}, the first author addressed a similar difficulty using the following three steps:
\begin{enumerate}[label=(\alph*),leftmargin=6ex]
    \item \label{intro:eq-item-B20-a} First, rewrite the random structure of $u_{\leq N}$   from \eqref{intro:eq-NLW-Ansatz} only in terms of the initial data $u_{\leq N}[0]$.  In our context, the random structure would have to be written as a function of only the sum $u_{\leq N}[0]=\initial{blue}-\initial{green}[][\leqN]+\initial{red}[][\leqN]$, and not as a function of the tuple $\big(\, \initial{blue},\initial{green}[][\leqN],\initial{red}[][\leqN]\big)$.  
    \item\label{intro:eq-item-B20-b} Second, show that new random structure of $u_{\leq N}$ from Step \ref{intro:eq-item-B20-a} satisfies a semigroup property. 
    \item Third, use the invariance of the Gibbs measure and the semigroup property from Step \ref{intro:eq-item-B20-b} to obtain global control of the new random structure.
\end{enumerate}
This strategy is rather general and should, at least in principle, be applicable to many dispersive equations. However, its main drawback is its computational and notational complexity, which stems from both Step \ref{intro:eq-item-B20-a} and Step \ref{intro:eq-item-B20-b}. For more details, we refer to \cite[Section 3 and 9]{B20II}. \\

Instead of the argument of \cite{B20II}, we introduce a new version of Bourgain's globalization argument \cite{B94}, which is also partially inspired by \cite{OOT21}. Our approach consists of the following three steps:
\begin{enumerate}[label=(\alph*)]
\item\label{intro:item-global-BNDY-a} First, we obtain local bounds on the solution $u_{\leq N}$ in a certain weaker norm $\|\cdot\|_W$. It is important that bounds in $\| \cdot \|_W$ can easily be iterated in time. 
\item\label{intro:item-global-BNDY-b} Second, we use the invariance of the Gibbs measure $\mu_{\leq N}$ and Step \ref{intro:item-global-BNDY-a} to obtain global bounds in the weaker norm $\| \cdot \|_W$  uniformly in the truncation parameter $N$. 
\item Third, we use the global bounds in the weaker norm $\| \cdot \|_W$ from Step \ref{intro:item-global-BNDY-b} and our main estimates (from Subsection \ref{section:main-estimates}) to obtain global control of the random structure. 
\end{enumerate}
In our article, the role of the weaker norm $\| \cdot \|_W$ from Step \ref{intro:item-global-BNDY-a} is played by the nonlinear smoothing norm  (see Definition \ref{ansatz:def-nonlinear-smoothing}). In contrast to the argument in \cite[(6.95)]{OOT21}, our nonlinear smoothing estimates (Proposition \ref{ansatz:prop-nonlinear-smoothing}) relies on multilinear dispersive effects. 

This strategy is seemingly less general than the strategy of \cite{B20II}, since it is unclear when a suitable weaker norm $\|\cdot\|_W$ can be found. If applicable, however, it is much simpler.  \\

We now briefly discuss our stability theory, which corresponds to 
Step \ref{intro:item-global-3} in Figure \ref{figure:globalization}. In the earlier works \cite{B20II,OOT21}, the stability theory is based on the full initial data $u_{\leq N}[0]$. As can be seen from \cite[Section 9]{B20II} and \cite[Proposition 6.5 and Remark 6.6]{OOT21}, earlier stability estimates require a tremendous computational effort, which would be practically infeasible in our setting. In our approach, the stability theory is instead stated in terms of the three individual components $\initial{blue}$, $\initial{green}[][\leqN]$, and $\initial{red}[][\leqN]$ (see Proposition \ref{global:prop-stability}), which simplifies earlier methods. For a more detailed discussion, we refer to Subsection \ref{section:global-stability}.

\subsection{Open problems}\label{section:intro-open}
Near the end of this introduction, we discuss some open problems in random dispersive equations. While the methods of this paper do not provide complete answers to either of these open problems, we believe that they address certain aspects.

\subsubsection{Geometric wave equations} 
In our previous discussion of parabolic singular SPDEs, we primarily focused on the parabolic $\Phi^4_d$-models. Due to the previously mentioned literature, such as \cite{CC18,H14,GIP15,K16}, the parabolic $\Phi^4_d$-models are now rather well-understood. Very recently, several works on singular parabolic SPDEs focused on geometric equations, such as the geometric stochastic heat equation \cite{BGHZ19} and the stochastic Yang-Mills equation \cite{CCHS20,CCHS22}. 

From both a mathematical and physical perspective, it would be interesting to obtain similar results for geometric wave and Schr\"{o}dinger equation with random initial data or stochastic forcing. At the time of writing, some initial progress in this direction has been made in \cite{BLS21,KLS20}, but almost all questions in this area remain wide open. We hope that the different counting and tensor estimates of this article, which were briefly discussed in  Subsection \ref{section:intro-main-structure}, will also be useful in the geometric setting.  

\subsubsection{The cubic nonlinear Schr\"{o}dinger equation in three dimensions} As shown in Figure \ref{figure:literature-cubic}, the 3D cubic nonlinear Schr\"{o}dinger equation is now the only classical $\Phi_d^4$ measure-invariance problem that still remains open. In fact, we expect it to be much harder than all positive results in Figure \ref{figure:literature-cubic} combined, for the sole reason that it is probabilistically \emph{critical}, namely $s_{pr}=s_G$ in Figure \ref{intro:scalingfig}. This means that the nonlinear evolution is only marginally (if at all) smoother than the linear evolution, which makes a fundamental challenge in the local theory.

Recently, the second author and Hani have resolved for the first time a closely related probabilistically critical problem, in the context of wave turbulence, see \cite{DH21}. We expect the full Gibbs measure problem to utilize all these state-of-art techniques, including (i) the full random tensor Ansatz in \cite{DNY20}, (ii) the deep molecular analysis in \cite{DH21}, and (iii) the caloric data and possible cancellations coming from mixed heat-Schr\"{o}dinger objects already exploited in this paper.
\subsubsection{Long time dynamics in regimes without an invariant measure} Another challenging problem is to control the propagation of randomness (in various forms) beyond the short time regime, in the absence of any invariant quantity (such as a Gibbs measure). For the corresponding parabolic problem there has now been some significant progress, see e.g. \cite{MW17,MW20}, but in the dispersive and wave settings the picture is relatively less clear. Below we just list a few directions along these lines:
\begin{enumerate}
\item Wave turbulence: this basically concerns (large) finite-time propagation of statistical quantities, according to some effective (wave kinetic) equations. Recently some fundamental progress has been made \cite{CG19,CG20,DH21,DH21II,DNY20,ST21}, and more developments are expected in the near future.
\item Quasi-invariance problems: this concerns infinite-time propagation of non-invariant measures in terms of absolute continuity with respect to the initial ensemble. We refer the reader to \cite{GLT22, GOTW18, NRSS, OT20II, PTV20, Tz15} and references therein for some examples, but getting to lower regularities is still a challenge. 
\item Global well-posedness problems: this concerns infinite-time propagation of regularity or of random structure (such as para-controlled terms, random tensors, etc.) of solutions. Here the prevailing idea has been to apply some form of energy or modified energy estimate, and some progress in this direction has been made in the Euclidean setting in,  for example,  \cite{BriScatter20, BriScatter21, DLM19, DLM20, KMV19, OP16, Poc17} while for the periodic settings we have \cite{CO12}.
\end{enumerate}
\subsection{Overview of the paper}

\tikzset{overviewrect/.style={rectangle, rounded corners, minimum width=4cm, minimum
height=1cm,text centered,align=center, draw=black, fill=white!10},
arrow/.style={thick,->,>=stealth}}
\tikzset{title/.style={rectangle,  minimum width=3cm, minimum
height=1cm,text centered,align=center, draw=black, fill=orange!10},
arrow/.style={thick,->,>=stealth}}
\tikzset{theorem/.style={rectangle, minimum width=4cm, minimum
height=1cm,text centered,align=center, draw=black, fill=red!10},
arrow/.style={thick,->,>=stealth}}

\tikzstyle{innerWhite} = [semithick, white,line width=1.4pt, shorten >= 4.5pt]

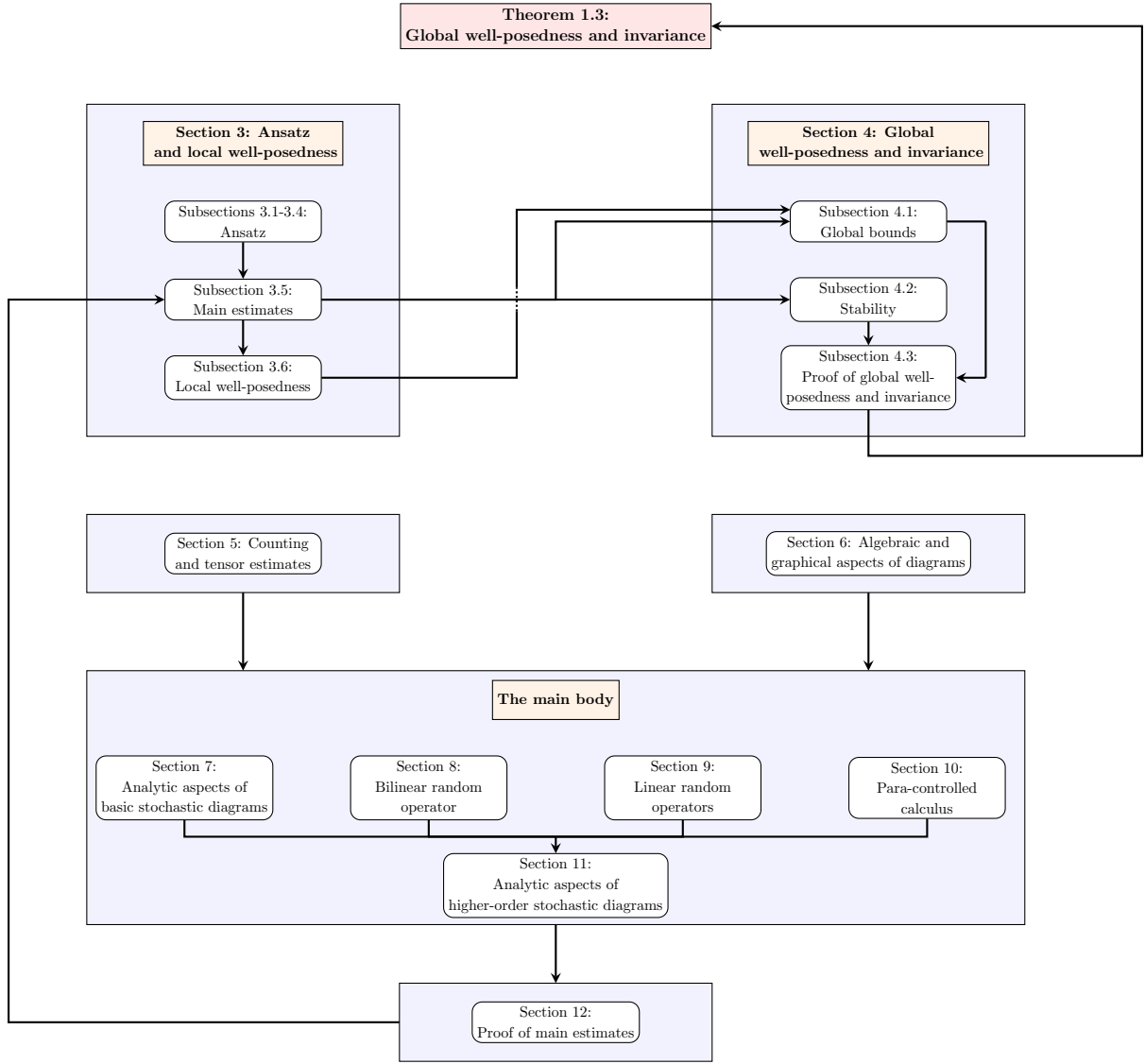
\begin{figure}[t]
\begin{center}
\begin{tikzpicture}[node distance=2cm,scale=0.55, every node/.style={scale=0.55}]

\node (thm)[theorem] at (0,2.5) {\textbf{\large \begin{NoHyper}Theorem \ref{intro:thm-rigorous}:\end{NoHyper}} \\ \textbf{\large Global well-posedness  and invariance}};

\node (local)[title]  at (-8,-0.5) 
{\textbf{Section \begin{NoHyper}\ref{section:ansatz}\end{NoHyper}: Ansatz} \\  \textbf{ and local well-posedness}};
\node (ansatz)[overviewrect]  at (-8,-2.5) 
{Subsections \begin{NoHyper}\ref{section:ansatz-caloric}-\ref{section:ansatz-para}\end{NoHyper}:\\ Ansatz};
\node (mestimate)[overviewrect]  at (-8,-4.5) {Subsection \begin{NoHyper}\ref{section:main-estimates}\end{NoHyper}:\\ Main estimates};
\node (lwp)[overviewrect]  at (-8,-6.5) 
{Subsection \begin{NoHyper}\ref{section:ansatz-proof-lwp}\end{NoHyper}:\\ Local well-posedness};
\begin{scope}[on background layer]
\draw [fill=blue!5] (-12,0.5) rectangle (-4,-8); 
\end{scope}

\node (global)[title]  at (8,-0.5) 
{\textbf{Section \begin{NoHyper}\ref{section:global}\end{NoHyper}: Global} \\  \textbf{ well-posedness and invariance}};
\node (gbounds)[overviewrect]  at (8,-2.5) 
{Subsection \begin{NoHyper}\ref{section:global-bounds}\end{NoHyper}:\\ Global bounds};
\node (stability)[overviewrect]  at (8,-4.5) {Subsection \begin{NoHyper}\ref{section:global-stability}\end{NoHyper}:\\ Stability};
\node (pgwp)[overviewrect]  at (8,-6.5) 
{Subsection \begin{NoHyper}\ref{section:global-proof}\end{NoHyper}:\\ Proof of global well-\\posedness and invariance};
\begin{scope}[on background layer]
\draw [fill=blue!5] (12,0.5) rectangle (4,-8); 
\end{scope}

\begin{scope}[on background layer]
\draw [fill=blue!5] (-12,-10) rectangle (-4,-12); 
\end{scope}
\node (prep)[overviewrect]  at (-8,-11) {Section \begin{NoHyper}\ref{section:counting}\end{NoHyper}:  Counting \\ and tensor estimates};

\begin{scope}[on background layer]
\draw [fill=blue!5] (12,-10) rectangle (4,-12); 
\end{scope}
\node (prep)[overviewrect]  at (8,-11) {Section \begin{NoHyper}\ref{section:diagrams}\end{NoHyper}:  Algebraic and \\ graphical aspects of diagrams};

\begin{scope}[on background layer]
\draw [fill=blue!5] (12,-20.5) rectangle (-12,-14); 
\end{scope}
\node (body)[title]  at (0,-14.75) 
{\textbf{The main body}};
\node (basic)[overviewrect] at (-9.5,-17) 
{\begin{NoHyper}Section \ref{section:analytic}\end{NoHyper}: \\
Analytic aspects of \\ basic stochastic diagrams};
\node (bilinear)[overviewrect] at (-3.25,-17) 
{\begin{NoHyper}Section \ref{section:bilinear}\end{NoHyper}: \\
Bilinear random \\ operator};
\node (linear)[overviewrect] at (3.25,-17) 
{\begin{NoHyper}Section \ref{section:linear}\end{NoHyper}: \\
Linear random \\ operators};
\node (para)[overviewrect] at (9.5,-17) 
{\begin{NoHyper}Section \ref{section:para}\end{NoHyper}: \\
Para-controlled \\calculus};
\node (high)[overviewrect] at (0,-19.5) 
{\begin{NoHyper}Section \ref{section:analytic2}\end{NoHyper}: \\
Analytic aspects of \\ higher-order stochastic diagrams};

\begin{scope}[on background layer]
\draw [fill=blue!5] (4,-22) rectangle (-4,-24); 
\end{scope}
\node (pmain)[overviewrect] at (0,-23) 
{\begin{NoHyper}Section \ref{section:proof-main-estimates}\end{NoHyper}: \\
Proof of main estimates};

\draw[arrow] (ansatz) -- (mestimate);
\draw[arrow] (mestimate) -- (lwp);

\draw[thick] ($(gbounds)+(3,0)$) -- (gbounds);
\draw[thick] ($(gbounds)+(3,0)$) -- ($(pgwp)+(3,0)$); 
\draw[arrow] ($(pgwp)+(3,0)$) -- (pgwp);
\draw[arrow] (stability) -- (pgwp);


\draw[arrow] (-8,-12) -- (-8,-14); 
\draw[arrow] (8,-12) -- (8,-14);

\draw[thick] (basic) |- (0,-18.25);
\draw[thick] (bilinear) |- (0,-18.25);
\draw[thick] (linear) |- (0,-18.25);
\draw[thick] (para) |- (0,-18.25);
\draw[arrow]  (0,-18.25) -- (high);

\draw[arrow] (0,-20.5) -- (0,-22); 

\draw[arrow] (mestimate) -- (stability);
\draw[arrow] (mestimate) -- (0,-4.5) -- (0,-2.5) -- (gbounds);
\draw[thick] (lwp) -- (-1,-6.5) -- (-1,-4.8);
\draw[thick,densely dotted] (-1,-4.8) -- (-1,-4.2);
\draw[arrow] (-1,-4.2) -- (-1,-2.2) -- ($(gbounds)+(-2,0.3)$);

\draw[arrow] (-4,-23) -- (-14,-23) -- (-14,-4.5) -- (mestimate);

\draw[arrow] (pgwp) -- (8,-8.5) -- (15,-8.5) -- (15,2.5) -- (thm); 

\end{tikzpicture}
\end{center}
\caption{\small{In this figure, we display the dependencies between different sections of this article. The heart of our article lies in Section \ref{section:ansatz} and Section \ref{section:global}, which are self-contained.}}\label{figure:organization}
\end{figure}

We now give a quick overview of the structure of this paper, which is also illustrated in Figure \ref{figure:organization}. In Section \ref{section:prelim}, we define parameters which are used throughout this article and recall basic definitions and lemmas from the literature. The heart of our article lies in the Ansatz and local theory (Section \ref{section:ansatz}) and global theory (Section \ref{section:global}). Despite utilizing estimates from the body of this article, both sections can be read without reading the rest of this article. In Section \ref{section:counting}, we prove several counting and tensor estimates, which will be used extensively in Section \ref{section:analytic}-Section \ref{section:analytic2}. In Section \ref{section:diagrams}, which is essentially self-contained, we address algebraic aspects of the stochastic diagrams. In Section \ref{section:analytic}-Section \ref{section:para}, we estimate basic stochastic objects, bilinear operators, linear operators, and para-controlled operators. These four sections depend heavily on the tensor estimates from Section \ref{section:counting}, but can be read independently of each other. In Section \ref{section:analytic2}, we estimate higher-order stochastic objects, and the arguments partially rely on earlier estimates in Section \ref{section:analytic}-Section \ref{section:para}. Finally, in Section \ref{section:proof-main-estimates}, we prove our main estimates (from Subsection \ref{section:main-estimates}), which collect all of our estimates in a few main propositions. \\

\textbf{Acknowledgements:} B.B. thanks Jonathan Mattingly, Felix Otto, Igor Rodnianski, and Terence Tao for helpful conversations. B.B. was partially supported by the NSF under Grant No. DMS-1926686. Y.D. was partially supported by the NSF under grant No. DMS-1900251 and the Sloan Research Fellowship. A.N. was partially supported by the NSF under Grants No. DMS-2101381 and DMS-205274,  and by the Simons Foundation Collaboration Grant on Wave Turbulence (Nahmod's Award ID  651469).  H. Y. was supported by a start-up funding of ShanghaiTech University.

This material is based upon work supported by both the National Science Foundation under Grant No. DMS-1929284 and the Simons Foundation Institute Grant Award ID 507536 while the authors were in residence at the Institute for Computational and Experimental Research in Mathematics in Providence, RI, during the fall 2021 semester.

\section{Prelimininaries}\label{section:prelim}

In this section, we introduce the parameters which will be used throughout this article. Furthermore, we recall a few basic definitions and lemmas from the literature. The expert may wish to skip this section and directly continue with Section \ref{section:ansatz} and Section \ref{section:global}, which are at the heart of this article. 

\subsection{Parameters} The most important parameters of our argument are $\delta_1,\delta_2,\epsilon,\eta \in (0,1)$, which are chosen such that 
\begin{equation}
0 <\epsilon \ll \delta_2 \ll \eta \ll \delta_1 \ll 1. 
\end{equation}
These parameters will be used to describe spatial regularities of terms in our Ansatz and frequency-restrictions in our para-products. We also choose additional parameters such that
\begin{equation}
0<c_4 \ll c_3 \ll c_2 \ll c_1 \ll \theta \ll \big(b_--\frac{1}{2}\big)\ll \big(b-\frac{1}{2}\big)\ll \big(b_+-\frac{1}{2}\big) \ll \epsilon \qquad \text{and} \qquad \delta_1 \ll \nu \ll 1. 
\end{equation}
For notational convenience, we also define 
\begin{equation}
\Theta := \theta^{-1} \qquad \text{and} \qquad C_j := c_j^{-1}, \quad \text{where } 1\leq j \leq 4. 
\end{equation}
While the above parameters are fixed throughout this article, we also introduce positive constants $\alpha,c$, and $C$, which are allowed to take different values in each estimate. Despite not having been fixed, the constant $\alpha$ never depends on any of the parameters. The constants $C$ and $c$ are permitted to be chosen depending on $b_-,b,b_+,\delta_1,\delta_2,\eta$, and $\nu$, i.e., 
\begin{equation}\label{intro:eq-c}
C=C(b_-,b,b_+,\delta_1,\delta_2,\eta,\nu) \qquad \text{and} \qquad c=c(b_-,b,b_+,\delta_1,\delta_2,\eta,\nu). 
\end{equation}
However, $C$ and $c$ are not allowed to depend on $c_1,c_2,c_3$, or $c_4$. Since $C$ and $c$ will only be used a finite number of times, it is therefore possible to guarantee that 
\begin{equation}
c_1 \ll c \qquad \text{and} \qquad  C \ll C_1. 
\end{equation}

\subsection{Basic definitions}\label{section:prelim-basic} In this subsection, we make the following basic definitions, which are primarily meant to fix our notation.  \\ 

\emph{\basicdef Fourier transform.} For any smooth function $f\colon \T^3 \rightarrow \mathbb{C}$, we define the spatial Fourier transform by 
\begin{equation}
    \widehat{f}(n) =  \big(\mathcal{F}_x f \big) (n):= \frac{1}{(2\pi)^3} \int_{\T^3}\dx f(x) e^{-i \langle n ,x \rangle}
\end{equation}
for all $n\in \Z^3$. Furthermore, for any $u\colon \R \times \T^3 \rightarrow \mathbb{C}$, we define the space-time Fourier transform by 
\begin{equation}
\big( \mathcal{F}_{t,x} u \big)(n,\rho) = \frac{1}{(2\pi)^6}
\int_{\R} \dt \int_{\T^3} \dx \, u(t,x) e^{-it\rho} e^{-i\langle n,x\rangle}
\end{equation}
for all $\rho\in \R$ and $n\in \Z^3$. In the following, we will often write the time-frequency variable $\rho$ as $\pm \langle n \rangle + \lambda$. \\

\emph{\basicdef Sums of frequency variables.} We represent sums of frequency variables by concatenating the indices of the individual frequency variables. For example, if $n_1,n_2,n_3 \in \Z^3$, we write
\begin{equation*}
n_{12} := n_1 + n_2 \quad \text{and} \quad n_{123} := n_1 + n_2 + n_3. 
\end{equation*}

\emph{\basicdef Norms on $\Z^3$.} For any $n=(n_1,n_2,n_3) \in \Z^3$, we define
\begin{equation}
|n| := \| n\|_2 = \big( n_1^2 + n_2^2 + n_3^2\big)^{1/2} \qquad \text{and} \qquad \langle n \rangle := \big( 1+ \| n\|_2^2 \big)^{1/2}.  
\end{equation}
We also define the $\ell^\infty$-norm by 
\begin{equation}
|n|_\infty := \max\big( |n_1|, |n_2|, |n_3|\big). 
\end{equation}
\emph{\basicdef Indicator functions.} For any statement $\mathcal{S}$, we write the corresponding indicator function as $\mathbf{1}\big\{ \mathcal{S}\big\}$. For example, if $n_0,n_1,n_2,n_3 \in \Z^3$, we write 
\begin{equation}
\mathbf{1}\big\{ n_0 =n_{123} \big\} =
\begin{cases}
\begin{tabular}{ll}
$1$ & if $n_0=n_{123}$,\\
$0$ & else. 
\end{tabular}
\end{cases}
.
\end{equation}

\emph{\basicdef Frequency-projections.} We now define the sharp cut-off function $1_N\colon \Z^3 \rightarrow \{0,1\}$, where $N\geq 1$ is a dyadic frequency-scale, as follows: If $N=1$, we define
\begin{equation}\label{intro:eq-cutoff-a}
1_N(n) = \mathbf{1}\big\{ |n|_\infty \leq 1 \big\}. 
\end{equation}
If $N\geq 2$, we define 
\begin{equation}\label{intro:eq-cutoff-b}
1_N(n) = \mathbf{1}\big\{ N/2 < |n|_\infty \leq N \big\}. 
\end{equation}
Since the zero frequency will not play a role, below we will ignore the discrepancy between (\ref{intro:eq-cutoff-a}) and (\ref{intro:eq-cutoff-b}) and only use (\ref{intro:eq-cutoff-b}). Furthermore, we define 
\begin{equation}
1_{\leq N}(n):= \sum_{K\leq N} 1_{K}(n). 
\end{equation}

Finally, we define the sequences of sharp frequency-truncation operators $(P_N)_N$ and $(P_{\leq N})_N$ as Fourier-multipliers with symbol $1_N$ or $1_{\leq N}$, respectively. To be precise, we define
\begin{equation}\label{intro:eq-frequency-projection}
\widehat{P_N f}(n) = 1_N(n) \widehat{f}(n) \qquad \text{and} \qquad \widehat{P_{\leq N} f}(n) = 1_{\leq N}(n) \widehat{f}(n). 
\end{equation}
The reason for using the $\ell^\infty$-norm instead of the $\ell^2$-norm in \eqref{intro:eq-cutoff-a} and \eqref{intro:eq-cutoff-b} is due to the mapping properties of $(P_N)_N$ and $(P_{\leq N})_N$.  \\

\emph{\basicdef Dependence on multiple parameters:} Many of the expressions in this article depend on multiple parameters, such as  frequency-scales $(N_j)_{j=1}^J \subseteq 2^{\mathbb{N}_0}$ or phase-functions $(\varphi_j)_{j=1}^{J} \subseteq \{ \cos, \sin \}$, where $J\in \mathbb{N}$. To simplify the notation, we often express this dependence by simply writing $N_\ast$ or $\varphi_\ast$, respectively. For instance, if a function $f\colon \T^3 \rightarrow \R$ depends on the three frequency-scales $N_1,N_2$, and $N_3$, we express this dependence by writing $f[N_\ast](x)$ instead of $f[N_1,N_2,N_3](x)$.\\

\emph{\basicdef Signs.} In linear combinations of dispersive symbols (see e.g. Section \ref{section:counting}), we repeatedly fix the corresponding signs. The fixed signs will be written as $(\pm_j)_{j \in \Nc} \in \{ +, - \}^{\Nc}$, where $\Nc$ is a finite index set. If $\pm_j$ occurs more than once in the same proof, it will refer to the same sign, and we also write $\mp_j$ for the opposite sign. 

In contrast, whenever we write $\pm$, the sign will not be fixed and different occurrences of $\pm$ may refer to different signs. \\

\emph{\basicdef Wick-ordering.}  For any $N\geq 1$, we first define 
\begin{equation}
    \sigma_{\leq N}^2 := \sum_{n\in \Z^3} \frac{1}{\langle n \rangle^2}, 
\end{equation}
which corresponds to the expectation of $\|\phi\|_{L^2_x}^2$ with respect to the Gaussian free field. For any $f\colon \T^3 \rightarrow \R$, we further define
\begin{align}
\lcol (P_{\leq N} f)^2 \rcol &=(P_{\leq N} f)^2 - \sigma_{\leq N}^2  \label{intro:eq-Wick-2},\\ 
\lcol (P_{\leq N} f)^3 \rcol &= (P_{\leq N} f)^3 -3  \sigma_{\leq N}^2 P_{\leq N} f \label{intro:eq-Wick-3},\\ 
\lcol (P_{\leq N} f)^4 \rcol &=  (P_{\leq N} f)^4 -6  \sigma_{\leq N}^2 (P_{\leq N} f)^2 + 3 \sigma_{\leq N}^4 \label{intro:eq-Wick-4}.
\end{align}

\emph{\basicdef Smooth cutoff-function in time.}  We fix a smooth, compactly supported function $\chi\colon \R \rightarrow [0,1]$ satisfying $\chi(t)=1$ for all $t \in [-1,1]$ and $\chi(t)=0$ for all $t\not \in [-2,2]$.\\

\emph{\basicdef Probability theory.} If $(\Omega,\mathcal{E},\mathbb{P})$ is a probability space, $E\in \mathcal{E}$ is an event, and $A\geq 1$, then the event $E$ is called $A$-certain with respect to $\mathbb{P}$ if 
\begin{equation}
\mathbb{P} \big( E \big) \geq 1 - c^{-1} \exp\big( -c A^c \big). 
\end{equation}
We recall that, since $c$ is assumed to be as in \eqref{intro:eq-c}, it does not depend on $A$. \bigskip

\subsection{Function spaces and linear estimates}

For any $\alpha \in \R$ and $f\colon \T^3 \rightarrow \mathbb{C}$, we define the Sobolev-norm
\begin{equation}
\| f \|_{H_x^\alpha}^2 := \sum_{n\in \Z^3} \langle n \rangle^{2\alpha} |\widehat{f}(n)|^2.
\end{equation}
The corresponding Sobolev space $H_x^\alpha$ is defined as the space of all distributions with finite $H_x^\alpha$-norm. Furthermore, we define
\begin{equation}\label{prelim:eq-cartesian-sobolev}
    \mathscr{H}_x^\alpha := H_x^\alpha \times H_x^\alpha. 
\end{equation}

\begin{definition}[$X^{s,b}$-spaces] 
Let $u\colon \R \times \T^3\rightarrow \mathbb{C}$. For any $s\in \R$ and $b\in \R$, we define the global $X^{s,b}$-norm by 
\begin{equation}\label{defofxsb}
\| u \|_{X^{s,b}(\R)} := \big\| \big\langle n \big\rangle^s \big\langle |\rho| - \langle n \rangle \big\rangle^b \big( \mathcal{F}_{x,t}u\big)(n,\rho) \big\|_{\ell_n^2 L_\rho^2(\T^3\times \R)}. 
\end{equation}
For any closed interval $\Jc\subseteq \R$, we also define the local $X^{s,b}$-norm by 
\begin{equation}
\| u \|_{X^{s,b}(\Jc)} := \inf \Big\{ \| v \|_{X^{s,b}(\R)} \colon 
v(t)=u(t) \text{ for all } t\in \Jc \Big\}. 
\end{equation}
The corresponding function spaces are denoted by $X^{s,b}(\R)$ and $X^{s,b}(\Jc)$, respectively. 
\end{definition}

For any $\varphi \in \{ \cos, \sin \}$, we define $\varphi\big( t \langle \nabla \rangle\big)$ as the Fourier multiplier with symbol $n \mapsto \varphi( t\langle n \rangle)$. Equipped with this definition, we can write the solution of the inhomogeneous linear wave equation 
\begin{equation*}
\big( \partial_t^2 - \Delta \big) \phi = F, \qquad \big( \phi, \langle \nabla \rangle^{-1} \partial_t \phi \big)\big|_{t=0}= \big( \phi^{\cos}, \phi^{\sin} \big),
\end{equation*}
as
\begin{equation}\label{prelim:eq-inhomogeneous-linear}
\phi = \cos(t\langle \nabla \rangle) \phi^{\cos} + \sin(t\langle \nabla \rangle) \phi^{\sin} + \int_0^t \dt^\prime \frac{\sin((t-t^\prime)\langle \nabla \rangle)}{\langle \nabla \rangle} F(t^\prime). 
\end{equation}
The integral on the right-hand side of \eqref{prelim:eq-inhomogeneous-linear} is called the Duhamel integral and, since it is used repeatedly in this article, we define
\begin{equation}\label{prelim:eq-duhamel}
\Duh F (t) := \int_0^t \dt^\prime \frac{\sin((t-t^\prime)\langle \nabla \rangle)}{\langle \nabla \rangle} F(t^\prime). 
\end{equation}
For notational convenience, we also define the localized Duhamel integral
\begin{equation}
\Duh_{\chi} F(t) := \chi(t) \Duh \big[ \chi F \big](t),
\end{equation}
where the cut-off function $\chi$ is as in Subsection \ref{section:prelim-basic}. \\

In working with $X^{s,b}$-norm, propagators, and the Duhamel integral, it is convenient to work with a twisted space-time Fourier-transform. For any $F\colon \R\times \T^3 \rightarrow \mathbb{C}$, $n\in \Z^3$, and $\lambda \in \R$, we define 
\begin{equation}
\label{tilted}\widetilde{F}^{\pm}(n,\lambda):=\mathcal{F}_{x,t}{[F]}(n,\lambda\pm \langle n\rangle),
\end{equation}
Similarly, if $h_n(t)$ is a function of $n\in \Z^3$ and $t\in \R$, we define
\begin{equation}\label{prelim:eq-twisted-FT-h}
\widetilde{h}^\pm_n(\lambda):=\mathcal{F}_{t}{[h]}(n,\lambda\pm \langle n\rangle).
\end{equation}
Equipped with this notation, we can state the following estimate for the localized Duhamel integral.

\begin{remark}\label{rmk:notation}
When the precise sign $+$ or $-$ in \eqref{tilted} or \eqref{prelim:eq-twisted-FT-h} is irrelevant, we sometimes simplify the notation by omitting the superscripts.
\end{remark}

We will use the following notation for the shifted Fourier transform of a spacetime function $F$:
 where $\mathcal{F}_{x,t}{(F)}$ is the spacetime Fourier transform of $F$. Similarly, assume $h_n(t)$ is a function of $(n,t)$, and then we can also define $\widetilde{h}^\pm_n(\lambda):=\mathcal{F}_{t}{[h]}(n,\lambda\pm \langle n\rangle)$.

We fix a smooth function $\chi(t)$ that is $1$ for $|t|\leq 1$ and $0$ for $|t|\geq 2$.
The Duhamel integral is defined as
\begin{equation}\label{eq:duhamel}
IF(t)=\int_0^t \frac{\sin {((t-s)\langle \nabla \rangle)}}{\langle \nabla \rangle} F(s)\,\mathrm{d}s,\quad \mathcal{I}F(t)=\chi(t)\cdot I(\chi(s)\cdot F(s)).
\end{equation}

\begin{lemma}[Lemma 4.1 in \cite{DNY20}]\label{lin} For all $F\colon \R \times \T^3\rightarrow \mathbb{C}$, we have the identity
\begin{equation}\label{linfor2}
\widetilde{\Duh_\chi \hspace{-0.5ex} F}^{\pm}(n,\lambda)=\int_{\mathbb{R}}K^{\pm}(\lambda,\sigma) \cdot \langle n\rangle^{-1}\cdot \widetilde{F}^{\pm}(n,\sigma)\,\mathrm{d}\sigma,
\end{equation} 
where the kernels $K^+(\lambda,\sigma)$ and $K^- (\lambda,\sigma)$ satisfy 
\begin{equation}\label{esti}|K^\pm(\lambda,\sigma)|\lesssim_B\bigg(\frac{1}{\langle\lambda\rangle^B}+\frac{1}{\langle\lambda\mp\sigma\rangle^B}\bigg)\frac{1}{\langle \sigma\rangle} \lesssim \frac{1}{\langle \lambda \rangle\langle \lambda\mp\sigma\rangle}
\end{equation}  
for all $B\geq 1$ and $\lambda,\sigma \in \R$. Furthermore, the derivatives
$\partial_\lambda K^\pm(\lambda,\sigma)$ and $\partial_\sigma K^\pm(\lambda,\sigma)$ hold the same bound as \eqref{esti}.
\end{lemma}

We now list a few basic properties of $X^{s,b}$-spaces.

\begin{lemma}[Energy estimate]\label{prep:lem-xsb}
Let $\Jc\subseteq \R$ be a compact interval containing zero and let $\phi$ be as in \eqref{prelim:eq-inhomogeneous-linear}. Then, we have for all $s\in \R$ and $b\in (1/2,1)$ that 
\begin{align*}
&\big\| \phi \big\|_{L_t^\infty H_x^s(\Jc \times \T^3)} 
+ \big\| \langle \nabla \rangle^{-1} \partial_t \phi \big\|_{L_t^\infty  H_x^s(\Jc \times \T^3)} + \big\| \phi \big\|_{X^{s,b}(\Jc)} \\
\lesssim& \,  \big(1+|\Jc|\big)^2 
\big( \big\| \phi^{\cos} \big\|_{H_x^s}
+ \big\| \phi^{\sin} \big\|_{H_x^s}
+ \big\| F\big\|_{X^{s-1,b-1}(\Jc)} \big). 
\end{align*}
\end{lemma}

\begin{lemma}[Time-localization]\label{sttime} 
Let $\Jc \subseteq \R$ be a compact time-interval of size $\tau:= |\Jc|\in (0,1)$, let $s\in \R$, and let $-1/2<b_1<b_2<1/2$. Then, we have for all $F\in X^{s,b_2}(\Jc)$ that
\begin{equation}
\big\| F \big\|_{X^{s,b_2}(\Jc)} \lesssim \tau^{b_2-b_1} \big\| F \big\|_{X^{s,b_1}(\Jc)}. 
\end{equation}
If instead $1/2<b_1<b_2 \leq 1$ and $u\in X^{s,b_2}(\Jc)$ satisfies $u(t_0)=0$ for some $t_0 \in \Jc$, then
\begin{equation}
\big\| u \big\|_{X^{s,b_2}(\Jc)} \lesssim \tau^{b_2-b_1} \big\| u \big\|_{X^{s,b_1}(\Jc)}. 
\end{equation}
\end{lemma}

\begin{lemma}[\protect{Gluing lemma, \cite[Lemma 4.5]{B20II}}]\label{prep:lem-gluing}
Let $s\in \mathbb{R}$, let $-1/2<b'<1/2$ and let $\Jc_1$ and $\Jc_2$ be compact intervals satisfying $\Jc_1\cap \Jc_2 \neq \varnothing$. Then, we have for all $F: (\Jc_1\cup \Jc_2)\times \mathbb{T}^3 \to \mathbb{R}$ that
\[
\|F\|_{X^{s,b'}(\Jc_1\cup \Jc_2)}\lesssim \|F\|_{X^{s,b'}(\Jc_1)} + \|F\|_{X^{s,b'}(\Jc_2)}.
\]
\end{lemma}

\subsection{Multiple stochastic integrals}
\label{section:prelim-multiple-stochastic}

In this subsection, we briefly recall the definition and elementary properties of multiple stochastic integrals. For more detailed discussions, we refer to the textbook \cite{N06}, the lecture notes \cite{M14}, or \cite[Subsection 4.6]{B20II}. 

We let $(\Omega,\mathcal{E},\mathbb{P})$ be a probability space and let $(W_s(n))_{n\in \Z^3}$ be a sequence of Gaussian processes on $(\Omega,\mathcal{E},\mathbb{P})$ satisfying the following properties: 

\begin{enumerate}[label={(\roman*)},leftmargin=9mm]
    \item $\Wp[][s][0]$ is a two-sided  real-valued standard Brownian motion and $(\Wp[][s][n])_{n\in \Z^3\backslash \{0\}}$ is a sequence of two-sided complex-valued standard Brownian motions. To be precise, we impose the normalization $\mathbb{E}[|\Wp[][s][n]|^2]=s$ for all $s\in \R$ and $n \in \mathbb{Z}^3$.
    \item For all $m,n\in \Z^3$ satisfying $m\neq \pm n$, the processes $\Wp[][s][n]$  and $\Wp[][s][m]$ are independent.
    \item For all $n \in \Z^3$,  $\overline{\Wp[][s][n]}=\Wp[][s][-n]$. 
\end{enumerate}
For any bounded and measurable set $A\subseteq \mathbb{R} \times \Z^3$, we define
\begin{equation*}
W[A] = \sum_{n\in \Z^3} \int_{\R} \dW[][s][n] \mathbf{1}\big\{ (s,n)\in A\big\},
\end{equation*}
where the right-hand side is understood as a sum of It\^{o}-integrals.

\begin{definition}[Multiple stochastic integrals]\label{prelim:def-multiple-stochastic}
Let $\mathcal{N}$ be a finite index set and let $f\colon (\R \times \Z^3)^{\Nc} \rightarrow \mathbb{C}$ be an elementary function of the form
\begin{equation*}
f\big( (s_j,n_j)_{j\in \Nc} \big) 
= \sum_{\ell_{\Nc}} a_{\ell_\Ac} \mathbf{1}\Big\{ (s_j,n_j)_{j\in \Nc} \in \operatorname{\raisebox{-1pt}{\scalebox{1.5}{$\otimes$}}}\displaylimits_{j \in \Nc} A_{\ell_j}\Big\},
\end{equation*}
which satisfies the same boundedness and ``off-diagonal"-assumptions as in \cite[Subsection 4.6]{B20II} (see also \cite[Section 1]{N06} and \cite[Section 4]{M14}). Then, we define the multiple stochastic integral by 
\begin{equation}\label{prelim:eq-multiple-stochastic}
\sum_{n_\Nc \in (\Z^3)^{\Nc}} \int_{\R^{\Ac}} 
\medotimes_{j\in \Nc} \hspace{-0.25ex} \dW[][s_j][n_j] f\big( (s_j,n_j)_{j\in \Nc} \big) 
:= \sum_{\ell_\Nc} a_{\ell_{\Nc}} \prod_{j\in \Nc} W[A_{\ell_{j}}]. 
\end{equation}
For a general function $f\in L^2\big((\R\times \Z^3)^{\Nc}\big)$, the multiple stochastic integral is defined by a density argument and \eqref{prelim:eq-multiple-stochastic}.
\end{definition}

In the following, we state two lemmas, which will be used repeatedly throughout this article. The first lemma controls higher-order moments of multiple stochastic integral.

\begin{lemma}[Gaussian hypercontractivity]\label{prelim:lem-hypercontractivity}
For all finite index sets $\Nc$, all $f\in L^2\big((\R\times \Z^3)^{\Nc}\big)$, and all $p\geq 2$, it holds that 
\begin{align*}
&\E \Big[ \Big| \sum_{n_\Nc \in (\Z^3)^{\Nc}} \int_{\R^{\Nc}} 
\medotimes_{j\in \Nc} \hspace{-0.5ex} \dW[][s_j][n_j] f\big( (s_j,n_j)_{j\in \Nc} \big)  \Big|^{p} \Big]^{1/p} \\\lesssim&( \sqrt{p})^{\# \Nc} \E \Big[ \Big| \sum_{n_\Nc \in (\Z^3)^{\Nc}} \int_{\R^{\Nc}} 
\medotimes_{j\in \Nc}  \hspace{-0.5ex} \dW[][s_j][n_j] f\big( (s_j,n_j)_{j\in \Nc} \big)  \Big|^{2} \Big]^{1/2}.
\end{align*}
\end{lemma}

The next lemma consists of a product formula for multiple stochastic integrals. To state it, however, we first require the following definition.

\begin{definition}[Pairings and contractions]
Let $\Nc$ and $\Kc$ be two finite disjoint index sets. A collection of two-elements sets $\Pc$ is called a pairing of $\Nc$ and $\Kc$ if all sets in $\Pc$ are disjoint and contain exactly one element of $\Nc$ and one element of $\Kc$. We denote the paired elements of $\Nc$ by $\Nc_p$, the paired elements of $\Kc$ by $\Kc_p$, and the unpaired elements of $\Nc \cup \Kc$ by $\Uc$. 

For any $f\in L^2\big( (\R \times \Z^3)^{\Nc} \big)$ and $g\in L^2\big( (\R \times \Z^3)^{\Kc}\big)$, and pairing $\Pc$, we define the contraction
\begin{equation*}
\begin{aligned}
&\big( f \otimes_{\Pc} g \big)\big( (s_u,n_u)_{u\in \Uc} \big)\\
=& \sum_{\substack{(n_j)_{j\in \Nc_p} }} 
\sum_{(n_k)_{k\in \Kc_p}} 
\int \medotimes_{j\in \Nc_p} \hspace{-1ex} \ds_j \int \medotimes_{k\in \Kc_p}  \hspace{-1ex} \ds_k \prod_{\{j,k\}\in \Pc} \Big( \mathbf{1}\big\{ n_{jk} = 0 \big\} \delta\big( s_j - s_k \big) \Big)  f \big( (s_j,n_j)_{j\in \Nc}\big) 
g\big( (s_k,n_k)_{k\in \Kc} \big). 
\end{aligned}
\end{equation*}
\end{definition}

\begin{lemma}[Product formula]\label{prelim:lem-product-formula}
For all disjoint finite index sets $\Nc$ and $\Kc$, $f\in L^2\big( (\R \times \Z^3)^{\Nc} \big)$, and $g\in L^2\big( (\R \times \Z^3)^{\Kc}\big)$, it holds that 
\begin{align*}
&\bigg(\sum_{n_\Nc \in (\Z^3)^{\Nc}} \int_{\R^{\Nc}}
\medotimes_{j\in \Nc} \hspace{-0.5ex} \dW[][s_j][n_j] f\big( (s_j,n_j)_{j\in \Nc} \big) \bigg) \cdot \bigg(
\sum_{n_\Kc \in (\Z^3)^{\Kc}} \int_{\R^{\Kc}} 
\medotimes_{k\in \Kc} \hspace{-0.5ex} \dW[][s_k][n_k] g\big( (s_k,n_k)_{k\in \Kc} \big) \bigg) \\
=& \sum_{\Pc} \sum_{n_{\Uc}} \int_{\R^{\Uc}} \medotimes_{u\in \Uc} \hspace{-0.5ex} \dW[][s_u][n_u] \big( f \otimes_{\Pc} g \big) \big( (s_u,n_u)_{u\in \Uc} \big). 
\end{align*}
\end{lemma}
The aforementioned reference \cite[Section 1]{N06} contains a proof of a symmetric version of Lemma \ref{prelim:lem-product-formula}. For a non-symmetric version similar as in Lemma \ref{prelim:lem-product-formula}, we also refer the reader to \cite[Lemma 10.3]{H14}. \\

In the following (see e.g. Section \ref{section:diagrams}), we will also need multiple stochastic integrals in two independent copies $(\Wp[\cos][s][n])_{n\in \Z^3}$ and $(\Wp[\sin][s][n])_{n\in \Z^3}$ of $(\Wp[][s][n])_{n\in \Z^3}$. By using functions $f\colon (\R \times \Z^3 \times \{ \cos,\sin\})^{\Nc} \rightarrow \mathbb{C}$, Definition \ref{prelim:def-multiple-stochastic}, Lemma \ref{prelim:lem-hypercontractivity}, and Lemma \ref{prelim:lem-product-formula} can easily be extended to this more general setting. To simplify the notation, we also make the following definition. For any $m \geq 1$, $n_1,\hdots,n_m \in \Z^3$, and $\varphi_1,\hdots,\varphi_m \in \{ \cos, \sin\}$, we define 
\begin{equation}\label{diagram:eq-multiple-stochastic}
\mathcal{S}\mathcal{I}[n_j,\varphi_j \colon 1 \leq j \leq m] := 2^{\frac{m}{2}} \int_{[-\infty,0]^m} \medotimes_{j=1}^m \hspace{-0.5ex} \dW[\varphi_j][s_j](n_j) \Big( \prod_{j=1}^m \langle n_j \rangle e^{s_j \langle n_j \rangle^2} \Big),
\end{equation}
where the right-hand side is the multiple stochastic integral of
\begin{equation*}
f\big( (s_j^\prime,n_j^\prime,\varphi_j^\prime)_{j=1}^m \big) 
= 2^{m/2} \prod_{j=1}^m \Big( \mathbf{1}\{ \varphi_j^\prime = \varphi_j \} \mathbf{1}\{ n_j^\prime = n_j \} \langle n_j \rangle \mathbf{1}\{ s_j^\prime \leq 0 \} e^{s_j^\prime \langle n_j \rangle^2} \Big). 
\end{equation*}
We emphasize that  $\mathcal{S}\mathcal{I}[n_j,\varphi_j \colon 1 \leq j \leq m]$ is (essentially) normalized, i.e., 
\begin{equation*}
\mathbb{E}\Big[ \Big|\mathcal{S}\mathcal{I}[n_j,\varphi_j \colon 1 \leq j \leq m] \Big|^2\Big]\sim_m 1.
\end{equation*}
For any $\pm_1,\hdots,\pm_m \in \{+,-\}$, we also define
\begin{equation}
\mathcal{S}\mathcal{I}[n_j,\pm_j \colon 1 \leq j \leq m] = \sum_{\substack{\varphi_1,\hdots,\varphi_m \\ \in \{ \cos,\sin \}}} \Big( \prod_{j=1}^m c_{\pm_j,\varphi_j} \Big) \mathcal{S}\mathcal{I}[n_j,\varphi_j \colon 1 \leq j \leq m],
\end{equation}
where the coefficients are dictated by  the choice of $\pm_j$ and the expressions of cosine and sine as sums of complex exponentials; i.e. $c_{\pm_j,\cos} =1/2$ and $c_{\pm_j,\sin}=(\mp_j) i/2$.

\section{Ansatz and local well-posedness}\label{section:ansatz}

We now examine the frequency-localized nonlinear wave equation
\begin{equation}\label{ansatz:eq-frequency-truncated-NLW-raw}
\begin{cases}
(\partial_t^2 + 1 - \Delta) u_{\leq N} = - P_{\leq N} \Big( \lcol  (P_{\leq N} u_{\leq N})^3 \rcol + \gamma_{\leq N} u_{\leq N} \Big) \\
u_{\leq N}[0]= (\phi^{\cos}, \phi^{\sin}).
\end{cases}
\end{equation}
We recall that $N$ is a frequency-truncation parameter, $P_{\leq N}$ is a sharp frequency-truncation, ${\lcol (P_{\leq N} u_{\leq N})^3 \rcol}$ is the Wick-ordered cubic power, and $\gamma_{\leq N}$ is an additional renormalization. Furthermore, 
\begin{equation}\label{ansatz:eq-initial-convention}
u_{\leq N}[0] = \Big( u_{\leq N}(0), \langle \nabla \rangle^{-1} \partial_t u_{\leq N}(0)\Big). 
\end{equation}
As was already mentioned in the introduction, the $\langle \nabla \rangle^{-1}$-operator acting on $\partial_t u_{\leq N}(0)$ is a notational convenience, which emphasizes symmetry\footnote{The significance of this aspect of our formulation will be explained in full detail in Section \ref{section:diagrams}.} in (the Gaussian parts of) the random initial data $\phi^{\cos}$ and $\phi^{\sin}$. 

The main result of this section is the probabilistic local well-posedness of \eqref{ansatz:eq-frequency-truncated-NLW-raw} on the support of the Gibbs measure. A qualitative version of our local well-posedness result can be stated as follows: 

\begin{proposition}[Qualitative local well-posedness]\label{ansatz:prop-lwp-qualitative}
For any  $0<\tau\ll 1$, there exists a Borel measurable event $\operatorname{Local}_\tau \subseteq \mathscr{H}_x^{-1/2-\epsilon}$ such that the following two properties hold: 
\begin{itemize}
    \item[(i)] (High probability) It holds that 
    \begin{equation*}
    \mu\big( \operatorname{Local}_\tau  \big) \geq 1 - c_1^{-1} \exp\big( - c_1 \tau^{-c_1} \big). 
    \end{equation*}
    \item[(ii)] (Convergence) For all $(\phi^{\cos},\phi^{\sin})\in \operatorname{Local}_\tau $, the solutions $u_{\leq N}$ of \eqref{ansatz:eq-frequency-truncated-NLW-raw} converge in
    \begin{equation*}
        L_t^\infty \mathscr{H}_x^{-1/2-\epsilon}([-\tau,\tau]\times \T^3)
    \end{equation*}
    as $N$ tends to infinity. 
\end{itemize}
\end{proposition}

While Proposition \ref{ansatz:prop-lwp-qualitative} is a significant result in its own right, it will not enter into our proof of the main theorem (Theorem \ref{intro:thm-rigorous}) and is stated (and proven) only for illustrative purposes. Instead, we will rely on a more quantitative version, which is the subject of Proposition \ref{ansatz:prop-lwp-quantitative} below.

In addition to the local well-posedness, our globalization argument (see Section \ref{section:global}) also relies on a nonlinear smoothing estimate. While the solution $u_{\leq N}$ lives at spatial regularities $-1/2-$, the Duhamel integral of the  nonlinearity has regularity $0-$. The nonlinear smoothing effect will be captured through the following norm, which provides better $L_t^\infty L_x^\infty$-control. 

\begin{definition}[Nonlinear smoothing norm]\label{ansatz:def-nonlinear-smoothing}
For any closed interval $\Jc \subseteq \R$, we define the nonlinear smoothing norm by 
\begin{equation*}
\| \, u \, \|_{\NSN(\Jc)} = \inf_{\substack{v,w \colon \\ u=v+w}} \Big( 
\big\|\, v \, \big\|_{(L_t^\infty \C_x^{-\epsilon}\cap X^{-\epsilon,b})(\Jc \times \T^3)}
+ \big\| \, w \, \big\|_{X^{1/2+\delta_2,b}(\Jc \times \T^3)} \Big). 
\end{equation*}
\end{definition}

The definition of the $\NSN$-norm is motivated by our Ansatz from \eqref{intro:eq-NLW-Ansatz}. The terms $\scubic$, $\squintic$, $\SXXone$, and $\SXXtwo$ will be placed in $L_t^\infty \C_x^{-\epsilon} \medcap X^{-\epsilon,b}$, whereas the nonlinear remainder $Y$ will be placed in $X^{1/2+\delta_2,b}$. 
Equipped with Definition \ref{ansatz:def-nonlinear-smoothing}, we can now state our nonlinear smoothing estimate. 

\begin{proposition}[Nonlinear smoothing]\label{ansatz:prop-nonlinear-smoothing}
Let $M\geq 1$ and let $0<\tau\ll 1$. Then, there exists a Borel measurable event $\NSE_{M,\tau} \subseteq \mathscr{H}_x^{-1/2-\epsilon}$ such that the following two properties hold:
\begin{itemize}
    \item[(i)] (High probability) It holds that 
    \begin{equation*}
        \mu_{\leq M}\Big( \NSE_{M,\tau}\Big) \geq 1 - c_1^{-1} \exp\big(-c_1 \tau^{-c_1}\big). 
    \end{equation*}
    \item[(ii)] (Nonlinear smoothing) For all $(\phi^{\cos},\phi^{\sin})\in \NSE_{M,\tau}$, all $N\leq M$, and all $T\geq 1$, the solution $u_{\leq N}$ of \eqref{ansatz:eq-frequency-truncated-NLW-raw} satisfies
    \begin{equation*}
    \Big\| P_{\leq N} \Duh \Big[ \mathbf{1}\big\{ 0\leq  t\leq \tau \big\} \Big( \lcol (P_{\leq N}u_{\leq N})^3 \rcol + \gamma_{\leq N} u_{\leq N} \Big) \Big] \Big\|_{\NSN([-T,T])} \leq T^\alpha. 
    \end{equation*}
\end{itemize}
\end{proposition}

\begin{remark}
The most significant aspect of Proposition \ref{ansatz:prop-nonlinear-smoothing}, which will be further discussed in Section \ref{section:global}, is that the formulation does not explicitly involve our Ansatz (from Subsection \ref{section:ansatz-caloric}-Subsection \ref{section:ansatz-para}). This makes it particularly useful in our globalization argument.
\end{remark}

We now briefly describe the remainder of this section.
In Subsection \ref{section:ansatz-caloric}, we briefly discuss the caloric representation of the Gibbs measure $\mu$. A more detailed discussion, however, is postponed until Subsection \ref{section:diagram-parabolic} below. In Subsection \ref{section:ansatz-explicit}, Subsection \ref{section:ansatz-remainder}, and Subsection \ref{section:ansatz-para}, we derive our Ansatz for the solution $u_{\leq N}$. In Subsection \ref{section:main-estimates}, we state the main estimates needed in our local theory. Finally, in Subsection \ref{section:ansatz-proof-lwp}, we prove the local well-posedness of the cubic nonlinear wave equation \eqref{ansatz:eq-frequency-truncated-NLW-raw}. More precisely, we prove Proposition \ref{ansatz:prop-lwp-qualitative}, its quantitative version from Proposition \ref{ansatz:prop-lwp-quantitative} below, and the nonlinear smoothing estimate from Proposition \ref{ansatz:prop-nonlinear-smoothing}.  

\subsection{The caloric initial data and initial value problem}\label{section:ansatz-caloric}

In this section, we briefly discuss the caloric initial data 
\begin{equation*}
\initial{blue}, \quad \initial{green}[][\leqM], \quad \text{and} \quad \initial{red}[][\leqM],
\end{equation*}
which provides us with a representation of the Gibbs measure $\mu_{\leq M}$. In order to promptly continue with our Ansatz (Subsection \ref{section:ansatz-explicit}-Subsection \ref{section:ansatz-para}), we postpone\footnote{Since Section \ref{section:diagrams} is essentially self-contained, the reader can also skip ahead to Section \ref{section:diagrams} and then return to Section \ref{section:ansatz}.}  a more detailed discussion until Section \ref{section:diagrams}.

The starting point of our construction is an abstract probability space $(\mathfrak{Z},\mathcal{Z},\mathbb{Q})$ and random functions $(\phi^{\cos}_{\leq M},\phi^{\sin}_{\leq M})\colon \mathfrak{Z} \rightarrow \mathscr{H}_x^{-1/2-\epsilon}$ satisfying
\begin{equation*}
\Law_{\mathbb{Q}}\Big( \big( \phi^{\cos}_{\leq M},\phi^{\sin}_{\leq M} \big) \Big) = \mu_{\leq M}. 
\end{equation*}
That is, we represent the (frequency-truncated) Gibbs measure through the law of two random functions. However, two different properties of $(\phi^{\cos}_{\leq M},\phi^{\sin}_{\leq M})$ make it rather difficult to directly solve the initial value problem \eqref{ansatz:eq-frequency-truncated-NLW-raw}:
\begin{enumerate}[label=(\roman*)]
    \item\label{ansatz:item-caloric-motivation-1} The initial data $(\phi^{\cos}_{\leq M},\phi^{\sin}_{\leq M})$ only has regularity $-1/2-$.
    \item\label{ansatz:item-caloric-motivation-2} The Fourier-coefficients of $\phi^{\cos}_{\leq M}$ are non-Gaussian and probabilistically dependent. 
\end{enumerate}
The first property \ref{ansatz:item-caloric-motivation-1} requires us to use probabilistic methods in our local theory, since the regularity of the initial data is far below the deterministic threshold for local well-posedness. Unfortunately, many of the probabilistic methods for dispersive equations (see e.g. \cite{DNY19,DNY20,GKO18}) require both Gaussian and probabilistically independent Fourier-coefficients. To separate the two issues in \ref{ansatz:item-caloric-motivation-1} and \ref{ansatz:item-caloric-motivation-2}, we represent $\mu_{\leq M}$ as the sum of two low-regularity functions, which have good probabilistic properties, and a high-regularity function, which (potentially) has bad probabilistic properties. As will be explained in Section \ref{section:diagrams}, such a representation can be obtained by solving the cubic stochastic heat equation (or parabolic $\Phi^4_3$-model) with initial data given by $(\phi^{\cos}_{\leq M},\phi^{\sin}_{\leq M})$. In the following proposition, we record the most important aspects of the construction. 

\begin{proposition}[Caloric initial data]\label{ansatz:prop-caloric}
Let $M\geq 1$, $(\mathfrak{Z},\mathcal{Z},\mathbb{Q})$, and $(\phi^{\cos}_{\leq M},\phi^{\sin}_{\leq M})$ as above. Let $A\geq 1$ and let $(\Omega,\mathcal{E},\mathbb{P})$ be the ambient probability space from Section \ref{section:diagrams}. Then, the Gibbs measure $\mu_{\leq M}$ can be represented as 
\begin{equation*}
\mu_{\leq M} = \Law_{\mathbb{P}\otimes \mathbb{Q}} \Big( \initial{blue} - \initial{green}[][\leqM] + \initial{red}[][\leqM] (A,\phi^{\cos}_{\leq M}) \Big)
\end{equation*}
and we have the following properties: 
\begin{enumerate}[label=(\roman*)]
    \item\label{ansatz:item-caloric-1} (Gaussian) The random function $\initial{blue}$ is as in Definition \ref{diagram:definition-caloric}. In particular, it is Gaussian, $\mathcal{E}$-measurable, and has spatial regularity $-1/2-\epsilon$. 
    \item\label{ansatz:item-caloric-2} (Cubic Gaussian chaos) The random function $\initial{green}[][\leqM]$ is as in Definition \ref{diagram:definition-caloric}. In particular, it is a cubic Gaussian chaos, $\mathcal{E}$-measurable, and has spatial regularity $1/2-\epsilon$. 
    \item\label{ansatz:item-caloric-3} (Remainder) There exists an $A$-certain event $E_A \in \mathcal{E} \otimes \mathcal{Z}$ such that, on this event, the remainder satisfies 
    \begin{equation*}
        \Big\| \initial{red}[][\leqM](A,\phi^{\cos}_{\leq M}) \Big\|_{\mathscr{H}_x^{1-\epsilon}} \leq A. 
    \end{equation*}
\end{enumerate}
\end{proposition}

\begin{remark}
\begin{enumerate}[label=(\alph*)]
    \item As was already mentioned in the introduction, the term ``caloric" is motivated by the connection with the caloric gauge from \cite{Tao04}. 
    \item The representation in Proposition \ref{ansatz:prop-caloric} is more detailed than in the work of the first author \cite[Theorem 1.1]{B20II}. In \cite{B20II}, the form of the cubic Gaussian chaos $\initial{green}[][\leqM]$ is not essential and it is simply hidden in the remainder $\initial{red}[][\leqM]$. In the present setting, however, the precise form of $\initial{green}[][\leqM]$ is crucial for several aspects of our argument.
    \item We note that the remainder $\initial{red}[][\leqM]=\initial{red}[][\leqM](A,\phi^{\cos}_{\leq M})$ in the caloric representation is $A$-dependent. In particular, if we want to increase the probability of the event $E_A$ by increasing $A$, we also have to change our representation. Fortunately, this does not create any difficulties in our argument. 
\end{enumerate}
\end{remark}

Equipped with the caloric initial data (from Proposition \ref{ansatz:prop-caloric}) we now turn to the cubic nonlinear wave equation. To this end, 
  we let $M\geq N \geq 1$ be frequency-scales and examine the caloric initial value problem
\begin{equation}\label{ansatz:eq-frequency-truncated-NLW-orig}
\begin{cases}
(\partial_t^2 + 1 - \Delta) u_{N} = - P_{\leq N} \Big( \lcol  (P_{\leq N} u_{\leq N})^3 \rcol + \gamma_{\leq N} u_{\leq N} \Big) \\
u_{\leq N}[0]= \initial{blue}[][] - \initial{green}[][\leqM] + \initial{red}[][\leqM].
\end{cases}
\end{equation}

\begin{remark}\label{ansatz:rem-ucal}
In the rest of this subsection, we will mostly consider solutions of the caloric initial value problem \eqref{ansatz:eq-frequency-truncated-NLW-orig}. However, in some arguments, such as in the proofs of Proposition \ref{ansatz:prop-lwp-qualitative} and Proposition \ref{ansatz:prop-nonlinear-smoothing}, we need to consider solutions of both the general initial value problem \eqref{ansatz:eq-frequency-truncated-NLW-raw} and the caloric initial value problem \eqref{ansatz:eq-frequency-truncated-NLW-orig}. Whenever necessary, we distinguish the solutions by writing either $u_{\leq N}^{\phi}$ or $u_{\leq N}^{\textup{cal}}$, respectively.
\end{remark}

In the following three subsections, we rigorously derive our Ansatz for the solution of  \eqref{ansatz:eq-frequency-truncated-NLW}, which was briefly discussed in Subsection \ref{section:intro-main-ideas}. While the derivation consists of several steps, they can be categorized into the following three groups. 
\begin{enumerate}[label=(\Roman*)]
\item \label{ansatz:item-group-1} Subtract explicit stochastic objects, which correspond to Picard iterates, from $u_{\leq N}$. 
\item \label{ansatz:item-group-2}  Derive the nonlinear wave equation for the remainder $v_{\leq N}$ and exhibit cancellations. 
\item \label{ansatz:item-group-3}  Define the para-controlled components $\Xone_{\leq N}$ and $\Xtwo_{\leq N}$ and the smooth nonlinear remainder $Y_{\leq N}$.
\end{enumerate}

From a practical perspective, it is convenient to exhibit the cancellations in Step \ref{ansatz:item-group-2} before further decomposing $v_{\leq N}$ into $\Xone_{\leq N}$, $\Xtwo_{\leq N}$, and $Y_{\leq N}$ in Step \ref{ansatz:item-group-3}. Otherwise, it is more difficult to correctly group the terms together.

\subsection{The explicit stochastic objects in $u_{\leq N}$}\label{section:ansatz-explicit}

In order to define the explicit stochastic objects, it is convenient to write the renormalization constant $\gamma_{\leq N}$ as 
\begin{equation}\label{ansatz:eq-decomposition-gamma}
\gamma_{\leq N} = \Gamma_{\leq N} + \big( \gamma_{\leq N} - \Gamma_{\leq N} \big),
\end{equation}
where  $\Gamma_{\leq N}$ is the renormalization multiplier from Definition \ref{sec6:def:Gamma}. 
 The motivation behind \eqref{ansatz:eq-decomposition-gamma} is that $\Gamma_{\leq N}$ exactly cancels certain resonant terms, so it is easier to base our stochastic objects on $\Gamma_{\leq N}$ instead of $\gamma_{\leq N}$. The remainder $\gamma_{\leq N}- \Gamma_{\leq N}$, which turns out to be a Fourier multiplier with symbol bounded by $\langle n \rangle^\epsilon$, can be treated perturbatively. Of course, it would have been possible to directly renormalize the equation \eqref{ansatz:eq-frequency-truncated-NLW} using $\Gamma_{\leq N}$ instead of $\gamma_{\leq N}$. However, both the physical and parabolic literature (see e.g. \cite{BHZ19,GJ87,H18}) prefer local counterterms, thus precluding the direct use of renormalization multipliers. Inserting 
\eqref{ansatz:eq-decomposition-gamma} into \eqref{ansatz:eq-frequency-truncated-NLW-orig}, we now write the equation as
\begin{equation}\label{ansatz:eq-frequency-truncated-NLW}
\begin{cases}
(\partial_t^2 + 1 - \Delta) u_{ \leq N} = - P_{\leq N} \Big( \lcol  (P_{\leq N} u_{\leq N})^3 \rcol + \Gamma_{\leq N} u_{\leq N} \Big) - (\gamma_{\leq N} - \Gamma_{\leq N}) P_{\leq N} u_{\leq N} \\
u_{\leq N}[0]= \initial{blue}[][] - \initial{green}[][\leqM] + \initial{red}[][\leqM].
\end{cases}
\end{equation}

Throughout the rest of this section, we use our shorthand notation for stochastic diagrams. The longhand notation, which was already mentioned in Subsection \ref{section:intro-main-ideas}, will not be used here. We first define the linear stochastic object $\slinear[blue]$ as the solution to
\begin{equation}\label{ansatz:eq-linear-blue}
(\partial_t^2 +1 - \Delta) \slinear[blue] =0, \qquad \slinear \hspace{0ex}  [0]= \initial{blue}[][]. 
\end{equation}
To simplify the notation, we also define  $\slinear[blue][\leqN] = P_{\leq N} \slinear[blue][]$. 
Similarly, we define the linear evolution of the green caloric initial data by 
\begin{align}
(\partial_t^2 +1 - \Delta) \slinear[green][\leqM] &=0, \qquad \slinear[green][\leqM] \hspace{0ex}  [0]= \initial{green}[][\leqM] \label{ansatz:eq-linear-greenM}. 
\end{align}
The two stochastic objects in \eqref{ansatz:eq-linear-blue} and \eqref{ansatz:eq-linear-greenM} correspond to the zeroth Picard iterate of \eqref{ansatz:eq-frequency-truncated-NLW} which initial data given by only $\initial{blue}[][]$ or  $\initial{green}[][\leqM]$, respectively. As we will see in Section \ref{section:analytic} below, the spatial regularity of $\slinear[blue][\leqN]$ is $-1/2-$ and the spatial regularity of $\slinear[green][\leqM]$ is $1/2-$. 
To simplify the notation, we also define
\begin{equation}\label{ansatz:eq-quadratic-cubic-nl}
\squadratic[\leqN] \overset{\textup{def}}{=} \lcol \big( \slinear[blue][\leqN] \big)^2 \rcol \qquad \text{and} 
\qquad \scubicnl[\leqN] \overset{\textup{def}}{=} \lcol \big( \slinear[blue][\leqN] \big)^3 \rcol.
\end{equation}
We emphasize that the right-hand side in \eqref{ansatz:eq-quadratic-cubic-nl} contains the Wick-ordering (\eqref{intro:eq-Wick-2}--\eqref{intro:eq-Wick-3})
but does not contain the renormalization multiplier $\Gamma_{\leq N}$. The reason is that $\Gamma_{\leq N}$ is designed to cancel a double-resonance, which cannot occur in the quadratic or cubic stochastic object. 
We now turn to higher-order Picard iterates. To this end,  we define the cubic stochastic object by 
\begin{equation}\label{ansatz:eq-cubic}
(\partial_t^2 +1 - \Delta) \scubic [\leqN] =P_{\leq N} \scubicnl[\leqN], \qquad \scubic[\leqN] \hspace{0ex} [0] =0. 
\end{equation}
In contrast to $\slinear[blue][\leqN]$, we emphasize that $\scubic[\leqN]$ is not the low-frequency projection of a limiting stochastic object. In particular, the two stochastic objects $\scubic[\leqN]$ and $P_{\leq N} \scubic[\leqK]$, where $K>N$, are not equal. As we will see in Section \ref{section:analytic} below, the spatial regularity of $\scubic[\leqN]$ is $0-$. Finally, we define the quintic stochastic object by 
\begin{equation}\label{ansatz:eq-quintic}
\begin{cases}
3 \, (\partial_t^2 +1 - \Delta) \squintic[\leqN] = P_{\leq N}  \Big( 3 \squadratic[\leqN]  \scubic[\leqN] - \Gamma_{\leq N} \slinear \Big) \qquad (t,x) \in \R \times \T^3, \\
\squintic[\leqN] \hspace{0ex}[0]=0. 
\end{cases}
\end{equation}
We note that a factor of three has been included on the left-hand side of \eqref{ansatz:eq-quintic}, which simplifies the combinatorics below. 
As we will see in Section \ref{section:analytic} below, the spatial regularity of $\squintic[\leqN]$ is $1/2-$.

\subsection{The nonlinear remainder $v_{\leq N}$ and cancellations}\label{section:ansatz-remainder}

Equipped with the explicit stochastic diagrams from Subsection \ref{section:ansatz-explicit}, we now write the solution $u_{\leq N}$ as 
\begin{equation}\label{ansatz:eq-nonlinear}
u_{\leq N} = \slinear[blue] - \scubic[\leqN] - \slinear[green][\leqM] + 3 \squintic[\leqN] + v_{\leq N}. 
\end{equation}
Using the definitions of $\slinear[blue]$,  $\scubic[\leqN]$,  $\slinear[green][\leqM]$, and  $\squintic[\leqN]$, it follows that $v_{\leq N}$ solves the nonlinear wave equation
\begin{equation}\label{ansatz:eq-NLW-vN}
\begin{aligned}
&(\partial_t^2 +1 - \Delta) v_{\leq N} \\
=& - P_{\leq N} \bigg[ 3 \squadratic[\leqN] \Big( 3 \squintic[\leqN] - P_{\leq N} \slinear[green][\leqM] + P_{\leq N} v_{\leq N} \Big) + \Gamma_{\leq N} \Big( -\scubic[\leqN] - \slinear[green][\leqM] +3 \squintic[\leqN] + P_{\leq N} v_{\leq N} \Big) \\ 
 +& 3 \, \slinear[blue][\leqN] \Big( - \scubic[\leqN] - P_{\leq N} \slinear[green][\leqM] + 3 \squintic[\leqN] + P_{\leq N} v_{\leq N} \Big)^2   \\
 +& \Big( - \scubic[\leqN] - P_{\leq N} \slinear[green][\leqM] + 3 \squintic[\leqN] + P_{\leq N} v_{\leq N} \Big)^3  \\
 +& \big(\gamma_{\leq N} - \Gamma_{\leq N}\big) \Big( \slinear[blue][\leqN] -\scubic[\leqN] -  P_{\leq N} \slinear[green][\leqM] +3 \squintic[\leqN] + P_{\leq N} v_{\leq N}  \Big) \bigg]. 
\end{aligned}
\end{equation}
and has initial data given by $v_{\leq N}[0]=\initial{red}[][\leqM]$. In the rest of this subsection, we exhibit cancellations in the right-hand side of  \eqref{ansatz:eq-NLW-vN}.

In \cite{GKO18}, Gubinelli, Koch, and Oh already observed a cancellation in products of $\slinear[blue][\leqN]$ and (the para-controlled component of) the nonlinear remainder $v_{\leq N}$, which was also used in  \cite{B20II,OOT20,OOT21}. Throughout this article, we refer to this cancellation as the \sine-cancellation, since it crucially relies on the $\sin((t-t^\prime) \langle \nabla \rangle)$-multiplier in the Duhamel integral. In Lemma \ref{counting:lem-Sine-symmetrization}, we present a new perspective on the \sine-cancellation, which links it to a certain complete derivative. Our reason for mentioning the \sine-cancellation here is that, since it is used in several earlier works, it may have already been on the reader's mind. For the derivation of our Ansatz, however, the \sine-cancellation only plays a minor role (despite its importance in our estimates). Instead, we will now focus on another cancellation. 

As shown in Lemma \ref{diagram:lem-C33} below, the expectations (or resonant parts) of the linear$\times$quintic and cubic$\times$cubic-stochastic objects are ill-defined. More precisely, both expectations
\begin{equation}\label{ansatz:eq-diverging-objects}
\mathbb{E} \Big[ \slinear[blue][\leqN] \squintic[\leqN] \Big] \qquad \text{and} \qquad \mathbb{E} \Big[  \scubic[\leqN] \scubic[\leqN] \Big]
\end{equation}
diverge logarithmically as $N\rightarrow \infty$. This divergence cannot be removed through additional renormalizations, since the renormalization is already dictated by the construction of the Gibbs measure. While the individual expectations in \eqref{ansatz:eq-diverging-objects} diverge as $N$ tends to infinity, however, the linear combination 
\begin{equation}\label{ansatz:eq-cancellation}
\mathbb{E}\Big[ 6 \, \slinear[blue][\leqN] \squintic[\leqN] + \scubic[\leqN] \scubic[\leqN]  \Big] 
\end{equation}
has a well-defined limit. Thus, there is a major cancellation between the linear$\times$quintic and cubic$\times$cubic-stochastic objects, which we now call the \oftt-cancellation. As mentioned in Subsection \ref{section:intro-main-ideas}, the  \oftt-cancellation is one of the main novelties of this article. While the proof of the \oftt-cancellation is postponed until Section \ref{section:diagrams}, we give a heuristic motivation. 

\begin{remark}[Heuristic argument for the \oftt-cancellation]\label{ansatz:remark-cancellation-oftt} 
For a typical sample $(u_0,u_1)$ of the Gibbs measure $\mup_{\leq N}$, it is well-known that the Wick-ordered square $\lcol (P_{\leq N} u_0)^2 \rcol$ is well-defined (with uniform control in $N$). This follows easily from the construction of the Gibbs measure (see e.g. \cite{BG18,MW20}) and can also be obtained from our representation of $\mup_{\leq N}$ (see Proposition \ref{ansatz:prop-caloric}). As long as we believe in the invariance of the Gibbs measure, i.e., our main theorem, this suggests that $\lcol (P_{\leq N} u_{\leq N}(t))^2 \rcol$ is well-defined for all $t \in \R$. After inserting our Ansatz from \eqref{ansatz:eq-nonlinear} and using the binomial formula for Wick-polynomials, we obtain that 
\begin{align}
    &6 \, \slinear[blue][\leqN] \squintic[\leqN] + \scubic[\leqN] \scubic[\leqN] \\
    =& (P_{\leq N} u_{\leq N}(t))^2 \notag \\
    -&  \squadratic[\leqN] + 2 \, \slinear[blue][\leqN] \Big( \scubic[\leqN] + P_{\leq N}\slinear[green][\leqM] \Big) \label{ansatz:eq-heuristic-1} \\
    -& 2 \, \slinear[blue][\leqN] P_{\leq N} v_{\leq N} \label{ansatz:eq-heuristic-2}  \\
    +& 2 \, \scubic[\leqN] \Big( - P_{\leq N}\slinear[green][\leqM] + 3 \squintic[\leqN] + P_{\leq N} v_{\leq N} \Big) 
    - \Big( - P_{\leq N}\slinear[green][\leqM] + 3 \squintic[\leqN] + P_{\leq N} v_{\leq N} \Big)^2 
    \label{ansatz:eq-heuristic-3}.  
\end{align}
Based solely on the explicit formulas for our stochastic objects and the \sine-cancellation, one can prove that \eqref{ansatz:eq-heuristic-1} is well-defined. While $v_{\leq N}$ has regularity $1/2-$, and one may therefore be worried about the product in \eqref{ansatz:eq-heuristic-2}, this term can be treated as in the previous literature \cite{GKO18,B20II}. Finally, the terms in \eqref{ansatz:eq-heuristic-3} are clearly well-defined, as the sum of the regularities is positive. Thus, we therefore expect that the left-hand side
\[
    6 \, \slinear[blue][\leqN] \squintic[\leqN] + \scubic[\leqN] \scubic[\leqN] 
\]
is well-defined. Of course, this motivation does not constitute a proof, since it already presumes the validity of our main result.

A second heuristic argument for the \oftt-cancellation can be obtained from the wave kinetic equation. In fact, (\ref{ansatz:eq-cancellation}) is just the sum in $k$ of the quantity
\begin{equation}\label{ansatz:wke}\mathbb{E}\Big[ 6 \, \Fc_x\slinear[blue][\leqN](t,k) \cdot\overline{\Fc_x\squintic[\leqN](t,k)} + \left|\Fc_x\scubic[\leqN](t,k)\right|^2 \Big]. \end{equation} Let $u_{\leq N}^{\mathrm{Gau}}$ be the solution to (\ref{intro:eq-NLW-caloric}) with initial data $\initial{blue}[][]$, then (\ref{ansatz:wke}) is just the \emph{first order term} (that converges to a constant multiple of $t$ when $N\to\infty$) if one expands $\Ec(t,k):=\Eb|\Fc_x(u_{\leq N}^{\mathrm{Gau}})(t,k)|^2$ as a power series in time $t$:
\[\Ec(t,k)=\Ec(0,k)+ (\ref{ansatz:wke})+O(t^2).\] Now, the wave turbulence theory predicts that in the limit $N\to\infty$, the energy profile $\Ec(t,k)$ should asymptotically be $\Ec(t,k)\approx \varphi(N^{-1}t,N^{-1}k)$ where $\varphi$ solves the corresponding wave kinetic equation
\begin{multline}\label{ansatz:wke3}\partial_t\varphi(t,z)=\sum_{\epsilon_j\in\{\pm\}}\int_{|z|,|z_j|\leq 1}\delta(z-z_1-z_2-z_3)\delta(|z|-\epsilon_1|z_1|-\epsilon_2|z_2|-\epsilon_3|z_3|)\frac{\varphi(t,z)\varphi(t,z_1)\varphi(t,z_2)\varphi(t,z_3)}{|k|}\\\times\bigg[\frac{1}{|z|\varphi(t,z)}-\frac{\epsilon_1}{|z_1|\varphi(t,z_1)}-\frac{\epsilon_2}{|z_2|\varphi(t,z_2)}-\frac{\epsilon_3}{|z_3|\varphi(t,z_3)}\bigg]\,\mathrm{d}z_1\mathrm{d}z_2\mathrm{d}z_3.
\end{multline} Moreover, as $N\to\infty$, we can check that the initial data $\varphi(0,z)\to |z|^{-2}$, while $\varphi(t,z)=|z|^{-2}$ is a \emph{stationary} solution to (\ref{ansatz:wke3}). Thus it is expected that $\Ec(t,k)$ should be constant in time to the top order, hence (\ref{ansatz:wke}) should vanish to the top order, and so should (\ref{ansatz:eq-cancellation}). As the top order term in (\ref{ansatz:eq-cancellation}) is at most logarithmic in $N$, the lower order terms will converge and thus (\ref{ansatz:eq-cancellation}) has a well-defined limit.
\end{remark}

In order to state the evolution equations concisely, we first introduce additional notation associated with the \oftt-cancellation. We define
\begin{equation}\label{ansatz:eq-mathfrakC}
\begin{aligned}
\mathfrak{C}^{(1,5)}_{\leq N}(t) &:= \mathbb{E} \Big[ \slinear[blue][\leqN](t,x) \squintic[\leqN](t,x) \Big], \quad 
\mathfrak{C}^{(3,3)}_{\leq N}(t) := \mathbb{E} \Big[ \scubic[\leqN](t,x) \scubic[\leqN](t,x) \Big], \\
\mathfrak{C}_{\leq N}(t)&:= 6 \mathfrak{C}^{(1,5)}_{\leq N}(t) +\mathfrak{C}^{(3,3)}_{\leq N}(t). 
\end{aligned}
\end{equation}
We note that, due to the translation-invariance of the random initial data, the functions $\mathfrak{C}^{(1,5)}_{\leq N}(t)$, $\mathfrak{C}^{(3,3)}_{\leq N}(t)$, and $\mathfrak{C}_{\leq N}(t)$ do not depend on $x\in \mathbb{T}^3$. In order to isolate the divergences in 
 $\slinear[blue][\leqN] \squintic[\leqN]$ and $\scubic[\leqN] \scubic[\leqN]$, we now want to decompose them as 
\begin{equation}\label{ansatz:eq-oftt-decomposition}
\slinear[blue][\leqN] \squintic[\leqN] 
= \bigg( \slinear[blue][\leqN] \squintic[\leqN] - \mathfrak{C}^{(1,5)}_{\leq N}  \bigg) + \mathfrak{C}^{(1,5)}_{\leq N} 
\qquad \text{and} \qquad  
\scubic[\leqN] \scubic[\leqN] 
= \bigg( \scubic[\leqN] \scubic[\leqN] - \mathfrak{C}^{(3,3)}_{\leq N} \bigg) + \mathfrak{C}^{(3,3)}_{\leq N}. 
\end{equation}
After inserting the decomposition \eqref{ansatz:eq-oftt-decomposition} into the evolution equation \eqref{ansatz:eq-NLW-vN}, we obtain new terms involving either $\mathfrak{C}^{(1,5)}_{\leq N}$ or $\mathfrak{C}^{(3,3)}_{\leq N}$. In order to utilize the \oftt-cancellation, we have to group all new terms together and ensure that $\mathfrak{C}^{(1,5)}_{\leq N}$ and $\mathfrak{C}^{(3,3)}_{\leq N}$ never occur individually, but only in combination with their corresponding stochastic objects or through their linear combination $\mathfrak{C}_{\leq N}$. In order to prove this algebraic result, it is convenient to make the following definition.

\begin{definition}[Symbols]\label{ansatz:def-symbols}
We define the set of basic symbols by 
\begin{equation}
\Symb^b := \Big \{ \slinear[blue], \, \scubic, \, \slinear[green][\leqM],\, \squintic, \, v\Big\}. 
\end{equation}
For each symbol $\zeta \in \Symb^b$ and dyadic scale $N$, we define $\zeta_{\leq N}$ as the random function obtained by adding the subscript ``$\leqN$'' to the stochastic object. For example, $\zeta=\slinear[blue]$ leads to $\zeta_{\leq N}=\slinear[blue][\leqN]$ and $\zeta= \scubic$ leads to $\zeta_{\leq N}=\scubic[\leqN]$. We also collect objects with regularity at least $0-$ in $\Symb^b_0$ and $1/2-$ in $\Symb^b_{1/2}$, i.e.,  we define
\begin{equation}
\Symb^b_0 := \Big \{  \scubic, \, \slinear[green][\leqM],\, \squintic, \, v\Big\} \qquad \text{and} \qquad 
\Symb^b_{1/2} := \Big \{  \, \slinear[green][\leqM],\, \squintic, \, v\Big\}. 
\end{equation}
\end{definition}
To avoid confusion, we remark that the set $\Symb_{1/2}^b$ will only be used in the proofs of our main estimates (see Section \ref{section:proof-main-estimates}), but not in their statements. 

In addition to the symbol notation from Definition \ref{ansatz:def-symbols}, we also introduce a modified product $\Pi^\ast_{\leq N}$, which subtracts the desired multiples of $\mathfrak{C}^{(1,5)}_{\leq N}$ and $\mathfrak{C}^{(3,3)}_{\leq N}$. 

\begin{definition}[Modified product] \label{ansatz:def-product}
We let $N\geq 1$,  $\zeta^{(1)} \in \Symb^b$, and $\zeta^{(2)}, \zeta^{(3)} \in \Symb^b_0$. Then, we define a symmetric modified product $\Pi^\ast_{\leq N}\big(\zeta^{(1)}_{\leq N},  \zeta^{(2)}_{\leq N},  \zeta^{(3)}_{\leq N}\big)$ as follows:
\begin{alignat*}{3}
\Pi^\ast_{\leq N} \Big(  \slinear[blue][\leqN] , \scubic[\leqN], \scubic[\leqN]\Big) &:=  \slinear[blue][\leqN] \Big( \Big( \scubic[\leqN] \Big)^2  -  \mathfrak{C}^{(3,3)}_{\leq N} \Big), &&&\\ 
\Pi^\ast_{\leq N}\Big(   \slinear[blue][\leqN], \squintic[\leqN], \squintic[\leqN] \Big) &:=  
\slinear[blue][\leqN]  \squintic[\leqN] \squintic[\leqN] - 2 \mathfrak{C}_{\leq N}^{(1,5)} \squintic[\leqN],  &&&\\ 
\Pi^\ast_{\leq N} \Big(   \slinear[blue][\leqN], \squintic[\leqN],  \zeta^{(3)}_{\leq N} \Big) &:=  
\slinear[blue][\leqN]  \squintic[\leqN]  \zeta^{(3)}_{\leq N} -  \mathfrak{C}_{\leq N}^{(1,5)}  \zeta^{(3)}_{\leq N}\qquad \qquad &\text{if }  \zeta^{(3)} \neq  \squintic[], &&\\ 
\Pi^\ast_{\leq N} \Big(   \slinear[blue][\leqN], \scubic[\leqN],  \zeta^{(3)}_{\leq N} \Big) &:=  
\slinear[blue][\leqN]  \scubic[\leqN]  \zeta^{(3)}_{\leq N} \qquad \qquad  &\text{if } \zeta^{(3)} \neq  \scubic[], \squintic[],  &&\\  
\Pi^\ast_{\leq N} \Big(   \slinear[blue][\leqN],  \zeta^{(2)}_{\leq N}, \zeta^{(3)}_{\leq N} \Big) &:=  
\slinear[blue][\leqN]   \zeta^{(2)}_{\leq N} \zeta^{(3)}_{\leq N}\qquad \qquad  &\text{if } \zeta^{(2)},  \zeta^{(3)} \neq  \scubic[], \squintic[],  &&\\ 
\Pi^\ast_{\leq N} \Big(   \scubic[\leqN], \scubic[\leqN], \scubic[\leqN]\Big) &:=  \Big( \scubic[\leqN] \Big)^3 - 3 \mathfrak{C}^{(3,3)}_{\leq N} \scubic[\leqN], &&&\\ 
\Pi^\ast_{\leq N} \Big(   \scubic[\leqN], \scubic[\leqN], \zeta^{(3)}_{\leq N}\Big) &:=  \Big( \Big( \scubic[\leqN] \Big)^2  - \mathfrak{C}^{(3,3)}_{\leq N} \Big) \zeta^{(3)}_{\leq N} \qquad \qquad  &\text{if }  \zeta^{(3)} \neq   \scubic[], &&\\ 
\Pi^\ast_{\leq N} \Big(   \scubic[\leqN],  \zeta^{(2)}_{\leq N},  \zeta^{(3)}_{\leq N}\Big) &:=   \scubic[\leqN]  \zeta^{(2)}_{\leq N} \zeta^{(3)}_{\leq N} \qquad \qquad  &\text{if } \zeta^{(2)}, \zeta^{(3)} \neq   \scubic[], &\\ 
\Pi^\ast_{\leq N} \Big(   \zeta^{(1)}_{\leq N},  \zeta^{(2)}_{\leq N},  \zeta^{(3)}_{\leq N}\Big) &:=  \zeta^{(1)}_{\leq N}  \zeta^{(2)}_{\leq N}  \zeta^{(3)}_{\leq N} \qquad \qquad  &\text{ if } \zeta^{(1)} \neq  \slinear[blue], \, \, \zeta^{(1)}, \zeta^{(2)},  \zeta^{(3)}\neq   \scubic[]. &&
\end{alignat*}
\end{definition}

\begin{remark}
In  Definition \ref{ansatz:def-product}, we only considered products involving at most one factor of $\slinear[blue][\leqN]$. While products with two factors of $\slinear[blue][\leqN]$ also occur in our nonlinear wave equation for $v_{\leq N}$, they will not be rewritten using the modified product notation. 
\end{remark}

We now prove the following lemma, which concerns the algebraic structure of the third, fourth, and fifth summand in \eqref{ansatz:eq-NLW-vN}.

\begin{lemma}[Grouping]\label{ansatz:lem-grouping}
There exist two maps 
\begin{equation}
\Aone \colon \Symb^b \rightarrow \Z \qquad \text{and} \qquad \Athree \colon \Symb^b \times \Symb_0^b \times \Symb_0^b \rightarrow \Z,
\end{equation}
which map combinations of symbols to absolute constants, such that the following algebraic identity holds: 
\begin{align}
& 3 \, \slinear[blue][\leqN] \Big( - \scubic[\leqN] - P_{\leq N}\slinear[green][\leqM] + 3 \squintic[\leqN] + P_{\leq N} v_{\leq N} \Big)^2  
+ \Big( - \scubic[\leqN] - P_{\leq N}\slinear[green][\leqM] + 3 \squintic[\leqN] +  P_{\leq N} v_{\leq N} \Big)^3 \label{ansatz:eq-algebraic-identity-1}\\
=& - 18 \mathfrak{C}^{(1,5)}_{\leq N}  \slinear[blue][\leqN]  
-\mathfrak{C}_{\leq N}  \sum_{\zeta \in \Symb^b} \Aone(\zeta) \zeta_{\leq N}  \label{ansatz:eq-algebraic-identity-2}\\
-& \sum_{\zeta^{(1)} \in \Symb^b} \sum_{\zeta^{(2)}, \zeta^{(3)} \in \Symb_0^b} \Athree(\zeta^{(1)}, \zeta^{(2)}, \zeta^{(3)} ) 
 \,    \Pi^\ast_{\leq N} \Big( \zeta^{(1)}_{\leq N}, \zeta^{(2)}_{\leq N}, \zeta^{(3)}_{\leq N} \Big). \label{ansatz:eq-algebraic-identity-3}
\end{align}
Furthermore, there exists a map 
\begin{equation*}
\Atone \colon \Symb^b \rightarrow \Z,
\end{equation*}
which also maps symbols to absolute constants, such that 
\begin{equation}
\big(\gamma_{\leq N} - \Gamma_{\leq N}\big) P_{\leq N} \Big( \slinear[blue][\leqN] -\scubic[\leqN] - \slinear[green][\leqM] +3 \squintic[\leqN] + v_{\leq N}  \Big)
=-  (\gamma_{\leq N} - \Gamma_{\leq N})  \sum_{\zeta \in \Symb^b} \Atone(\zeta) \zeta_{\leq N}.  \label{ansatz:eq-algebraic-identity-gamma}
\end{equation}
\end{lemma}
 
 The first term in \eqref{ansatz:eq-algebraic-identity-2}, i.e., $- 18 \mathfrak{C}_{\leq n}^{(1,5)} P_{\leq N} \slinear[blue][\leqN]$ is still divergent, but it will ultimately cancel with a term in the first summand of \eqref{ansatz:eq-NLW-vN}. 
From the proof, one can easily work out the precise values of the maps $\Aone$, $\Athree$, and $\Atone$. However, since the precise values are completely irrelevant for our argument, we do not keep track of them. 

\begin{proof} The second algebraic identity \eqref{ansatz:eq-algebraic-identity-gamma} is trivial and we can explicitly choose 
\begin{equation*}
\Big( \Atone\Big(\slinear[blue]\Big), \, \Atone\Big(\scubic \Big), \, \Atone\Big( \slinear[green][\leqM] \Big), \, \Atone\Big( \squintic \Big), \, \Atone\Big( v \Big) \Big) = \Big( -1 ,1 , 1 , -3 , -1 \Big). 
\end{equation*}
The proof of the first algebraic identity follows from a simple (but tedious) calculation. Due to the central importance of this identity to our work, we present the full details. 
To simplify the notation, we denote the linear space spanned by all summands \eqref{ansatz:eq-algebraic-identity-2} and \eqref{ansatz:eq-algebraic-identity-3} by $\mathfrak{L}$, i.e., 
\begin{equation}
\mathfrak{L} := \operatorname{span} \Big(  \Big \{ \mathfrak{C}_{\leq N}  P_{\leq N} \zeta_{\leq N} \colon \zeta \in \Symb^b \} \cup \Big\{  P_{\leq N} \Pi^\ast_{\leq N}\Big(\zeta^{(1)}_{\leq N},\zeta^{(2)}_{\leq N},\zeta^{(3)}_{\leq N} \Big) \colon \zeta^{(1)} \in \Symb^b, \zeta^{(2)}, \zeta^{(3)} \in \Symb^b_0 \Big\}\Big). 
\end{equation} 
In the following, most calculations will be performed up to elements in $\mathfrak{L}$, which we denote by writing $\operatorname{mod}\mathfrak{L}$.  For the first summand in \eqref{ansatz:eq-algebraic-identity-1}, we obtain that 
\begin{align*}
&3 \, \slinear[blue][\leqN] \Big( - \scubic[\leqN] - P_{\leq N}\slinear[green][\leqM] + 3 \squintic[\leqN] +  P_{\leq N} v_{\leq N} \Big)^2 \\
=& 3 \, \slinear[blue][\leqN] \Big( \scubic[\leqN] \Big)^2 + 27 \, \slinear[blue][\leqN] \Big( \squintic[\leqN] \Big)^2  + 18\, \slinear[blue][\leqN] \squintic[\leqN] \Big(   - \scubic[\leqN] - P_{\leq N}\slinear[green][\leqM]  + P_{\leq N} v_{\leq N} \Big)  \mod \mathfrak{L}. 
\end{align*}
For the second summand in \eqref{ansatz:eq-algebraic-identity-2}, we obtain that 
\begin{align*}
&\Big( - \scubic[\leqN] - P_{\leq N}\slinear[green][\leqM] + 3 \squintic[\leqN] + P_{\leq N} v_{\leq N} \Big)^3 \\
&= - \Big( \scubic[\leqN] \Big)^3 + 3 \Big( \scubic[\leqN] \Big)^2 \Big( - P_{\leq N}\slinear[green][\leqM] + 3 \squintic[\leqN] + P_{\leq N} v_{\leq N} \Big) \mod \mathfrak{L}. 
\end{align*}
By adding both expressions and grouping terms together, we obtain that 
\begin{equation} \label{ansatz:eq-algebraic-identity-4}
\begin{aligned}
&3 \, \slinear[blue][\leqN] \Big( - \scubic[\leqN] - P_{\leq N}\slinear[green][\leqM] + 3 \squintic[\leqN] + P_{\leq N} v_{\leq N} \Big)^2  + \Big( - \scubic[\leqN] - P_{\leq N}\slinear[green][\leqM] + 3 \squintic[\leqN] + P_{\leq N} v_{\leq N} \Big)^3 \\
=& 3 \, \slinear[blue][\leqN] \Big( \scubic[\leqN] \Big)^2 
+ \Big[ 27 \, \slinear[blue][\leqN] \Big( \squintic[\leqN] \Big)^2 + 9 \Big( \scubic[\leqN] \Big)^2   \squintic[\leqN] \Big]  
 -\Big[  18 \, \slinear[blue][\leqN] \squintic[\leqN] \scubic[\leqN] +  \Big( \scubic[\leqN] \Big)^3  \Big]  \\
 +& \Big[ \Big( 18 \,  \slinear[blue][\leqN] \squintic[\leqN] + 3 \Big( \scubic[\leqN] \Big)^2 \Big) \Big( - P_{\leq N}\slinear[green][\leqM]+ P_{\leq N} v_{\leq N} \Big) \Big]  \mod \mathfrak{L}. 
\end{aligned}
\end{equation}
For the first summand in \eqref{ansatz:eq-algebraic-identity-4}, we have that 
\begin{align*}
3  \, \slinear[blue][\leqN] \Big( \scubic[\leqN] \Big)^2 &= 3 \mathfrak{C}_{\leq N}^{(3,3)} \slinear[blue][\leqN] 
+ 3 \Pi^\ast_{\leq N}  \Big( \slinear[blue][\leqN], \scubic[\leqN], \scubic[\leqN] \Big) \\
&=-18 \mathfrak{C}_{\leq N}^{(1,5)} \slinear[blue][\leqN] + 3\mathfrak{C}_{\leq N} \slinear[blue][\leqN]  
+ 3 \Pi^\ast_{\leq N} \Big( \slinear[blue][\leqN], \scubic[\leqN], \scubic[\leqN] \Big) \\ 
&= -18 \mathfrak{C}_{\leq N}^{(1,5)} \slinear[blue][\leqN]  \mod \mathfrak{L}. 
\end{align*}
This yields the main term in \eqref{ansatz:eq-algebraic-identity-2}. The remaining three summands in \eqref{ansatz:eq-algebraic-identity-4} are all elements of $\mathfrak{L}$. Indeed, for the second summand in \eqref{ansatz:eq-algebraic-identity-4}, it holds that 
\begin{align*}
&27 \, \slinear[blue][\leqN] \Big( \squintic[\leqN] \Big)^2 + 9 \Big( \scubic[\leqN] \Big)^2   \squintic[\leqN]\\
=& 27 \bigg[\, \slinear[blue][\leqN] \Big( \squintic[\leqN] \Big)^2 
- 2 \mathfrak{C}^{(1,5)}_{\leq N} \squintic[\leqN] \bigg] + 
9  \bigg[ \bigg( \Big( \scubic[\leqN] \Big)^2  - \mathfrak{C}^{(3,3)}_{\leq N} \bigg) \squintic[\leqN] \bigg]
+ 9 \mathfrak{C}_{\leq N} \squintic[\leqN] \in \mathfrak{L}. 
\end{align*}
For the third summand in \eqref{ansatz:eq-algebraic-identity-4}, it holds that 
\begin{align*}
&18 \, \slinear[blue][\leqN] \squintic[\leqN] \scubic[\leqN] +  \Big( \scubic[\leqN] \Big)^3 \\
=& 18 \bigg[  \bigg( \, \slinear[blue][\leqN] \squintic[\leqN] - \mathfrak{C}^{(1,5)}_{\leq N} \bigg) \scubic[\leqN] \bigg] 
+\bigg[ \Big( \scubic[\leqN] \Big)^3 - 3 \mathfrak{C}^{(3,3)}_{\leq N} \scubic[\leqN] \bigg] 
+ 3 \mathfrak{C}_{\leq N} \scubic[\leqN] \in \mathfrak{L}.
\end{align*}
For the fourth summand in \eqref{ansatz:eq-algebraic-identity-4}, it holds that 
\begin{align*}
&\Big( 18 \, \slinear[blue][\leqN] \squintic[\leqN] + 3 \Big( \scubic[\leqN] \Big)^2 \Big) \Big( - P_{\leq N}\slinear[green][\leqM]+ P_{\leq N} v_{\leq N} \Big)\\
=& 18  \bigg[ \Big( \slinear[blue][\leqN] \squintic[\leqN] - \mathfrak{C}^{(1,5)}_{\leq N} \Big)
 \Big( - P_{\leq N}\slinear[green][\leqM]+P_{\leq N} v_{\leq N} \Big) \bigg] 
+ 3  \bigg[ \Big( \Big( \scubic[\leqN] \Big)^2 - \mathfrak{C}_{\leq N}^{(3,3)} \Big) \Big( - P_{\leq N}\slinear[green][\leqM]+ P_{\leq N} v_{\leq N} \Big) \bigg] \\
+& 3 \mathfrak{C}_{\leq N} \Big( - P_{\leq N}\slinear[green][\leqM]+ P_{\leq N} v_{\leq N} \Big) \in \mathfrak{L}. 
\end{align*}
As a result, \eqref{ansatz:eq-algebraic-identity-4} can be written as in \eqref{ansatz:eq-algebraic-identity-2} and \eqref{ansatz:eq-algebraic-identity-3}, which completes the proof. 
\end{proof} 

After applying Lemma \ref{ansatz:lem-grouping}, the equation for $v_{\leq N}$ takes the form
\begin{align}
&(\partial_t^2 +1 - \Delta) v_{\leq N}  \notag \\
=& - P_{\leq N} \bigg[ 9  \squadratic[\leqN] \squintic[\leqN] - \Gamma_{\leq N}  \scubic[\leqN] + 3 \Gamma_{\leq N}  \squintic[\leqN]  - 18 \mathfrak{C}_{\leq N}^{(1,5)} P_{\leq N} \slinear[blue][\leqN] \bigg]  \label{ansatz:eq-NLW-vN-a} \\
+& P_{\leq N} \bigg[ 3 \squadratic[\leqN] P_{\leq N} \slinear[green][\leqM]  + \Gamma_{\leq N}  \slinear[green][\leqM]  \bigg] \label{ansatz:eq-NLW-vN-b} \\ 
 -& P_{\leq N} \bigg[ 3 \squadratic[\leqN] P_{\leq N} v_{\leq N} + \Gamma_{\leq N} v_{\leq N} \bigg] \label{ansatz:eq-NLW-vN-c} \\
+&  \sum_{\zeta^{(1)} \in \Symb^b} \sum_{\zeta^{(2)}, \zeta^{(3)} \in \Symb_0^b} \Athree(\zeta^{(1)}, \zeta^{(2)}, \zeta^{(3)} ) 
 \,  P_{\leq N} \Pi^\ast_{\leq N} \Big( \zeta^{(1)}_{\leq N}, \zeta^{(2)}_{\leq N}, \zeta^{(3)}_{\leq N} \Big)   \label{ansatz:eq-NLW-vN-d} \\
 +&   \mathfrak{C}_{\leq N} \sum_{\zeta \in \Symb^b} \Aone(\zeta) P_{\leq N} \zeta_{\leq N}  +    (\gamma_{\leq N} - \Gamma_{\leq N}) \sum_{\zeta \in \Symb^b} \Atone(\zeta) P_{\leq N}\zeta_{\leq N} . \label{ansatz:eq-NLW-vN-e}
\end{align}
We emphasize that the first summand \eqref{ansatz:eq-NLW-vN-a} is rather complex. In fact, each of the three terms 
\begin{equation}\label{ansatz:eq-septic-resonant}
 \Gamma_{\leq N}  \scubic[\leqN] , \qquad \Gamma_{\leq N}  \squintic[\leqN], \qquad \text{and} \qquad  \mathfrak{C}_{\leq N}^{(1,5)} P_{\leq N} \slinear[blue][\leqN]
\end{equation}
will cancel a resonant term originating from $ \squadratic[\leqN] \squintic[\leqN]$. For the first two terms in \eqref{ansatz:eq-septic-resonant}, this is to be expected, since the cancellation of resonances is the whole purpose of the renormalization multiplier $\Gamma_{\leq N}$. For the third term in \eqref{ansatz:eq-septic-resonant}, it is a consequence of the \oftt-cancellation. 

\subsection{\protect{The para-controlled components $\Xone_{\leq N}$ and $\Xtwo_{\leq N}$}}\label{section:ansatz-para}

Instead of starting directly with rigorous definitions, we first motivate the para-controlled Ansatz for $X^{(1)}$ and $X^{(2)}$ with a heuristic discussion. In the previous subsection, we derived the evolution equation for the nonlinear remainder $v_{\leq N}$ and exhibited the \oftt-cancellation. Ideally, we would like to control $v_{\leq N}$ directly through a contraction argument, but this turns out not to be possible. The nonlinear term in \eqref{ansatz:eq-NLW-vN-d} contains the 
(high$\times$high$\rightarrow$low)$\times$low-interaction 
\begin{equation}\label{ansatz:eq-heuristic-hhtll}
\Duh \bigg[ P_{\lesssim 1} \Big( \slinear[blue][N] \cdot P_N v_{\leq N} \Big) P_{\lesssim 1} v_{\leq N} \bigg]
\end{equation}
The Duhamel integral in \eqref{ansatz:eq-heuristic-hhtll} does not experience any gain through multilinear dispersive effects. To see the reason for this, let us denote the frequencies of $\slinear[blue][N]$, $P_N v_{\leq N}$, and $P_{\lesssim 1} v_{\leq N}$ by $n_1$, $n_2$, and $n_3$, respectively. In the (high$\times$high$\rightarrow$low)$\times$low-interaction, the multilinear dispersive symbol satisfies
\begin{equation*}
|\langle n_1 + n_2 + n_3 \rangle - \langle n_1 \rangle + \langle n_2 \rangle - \langle n_3 \rangle| 
\lesssim \langle n_1 + n_2 \rangle + \langle n_3 \rangle \lesssim 1
\end{equation*}
and therefore cannot yield any gain. In order to obtain uniform control of \eqref{ansatz:eq-heuristic-hhtll} in $N$, we therefore need the regularities of $\slinear[blue][N]$ and $v_{\leq N}$ to add to a positive number. Since $\slinear[blue][N]$ has regularity $-1/2-$, this would require that $v_{\leq N}$ has regularity greater than $1/2$. Unfortunately, this is prevented by at least two different interactions in \eqref{ansatz:eq-NLW-vN-a}-\eqref{ansatz:eq-NLW-vN-e}, which will now be discussed. \\

\emph{High$\times$low$\times$low-interactions:} We first discuss certain high$\times$low$\times$low-interactions, which yield contributions of regularity at most $1/2-$. For illustrative purposes, we focus on 
\begin{equation}\label{ansatz:eq-heuristic-hll}
\Duh \bigg[ \slinear[blue][N] P_{\lesssim 1} v_{\leq N} P_{\lesssim 1} v_{\leq N} \bigg]. 
\end{equation}
Unfortunately, \eqref{ansatz:eq-heuristic-hll} experiences no gain through multilinear dispersive effects. Indeed, if $n_1$, $n_2$, and $n_3$ denote the frequencies of $\slinear[blue][N]$, $P_{\lesssim 1} v_{\leq N}$, and $P_{\lesssim 1} v_{\leq N}$, the multilinear dispersive symbol satisfies
\begin{equation*}
| \langle n_{123}\rangle - \langle n_1 \rangle - \langle n_2 \rangle - \langle n_3 \rangle | 
\lesssim \langle n_2 \rangle +\langle n_3 \rangle \lesssim 1.
\end{equation*}
As a result, the only gain of regularity comes from the $\langle \nabla \rangle^{-1}$-multiplier in the Duhamel integral. Since $\slinear[blue][N]$ has regularity $-1/2-$, the term \eqref{ansatz:eq-heuristic-hll} has regularity at most $1/2-$. While the regularity of \eqref{ansatz:eq-heuristic-hll} is therefore too low for a direct contraction argument, it exhibits a useful random structure. Indeed, it morally behaves like the Duhamel integral of $\slinear[blue][N]$. The para-controlled component  $\Xone_{\leq N}$, which will be introduced below, is used to capture \eqref{ansatz:eq-heuristic-hll} and similar terms.   \\

\emph{High$\times$high$\times$low-interactions:} We now discuss certain high$\times$high$\times$low-interactions, which will also only yield contributions of regularity at most $1/2-$. For illustrative purposes, we consider 
\begin{equation}\label{ansatz:eq-heuristic-hhl}
\Duh \bigg[ \lcol \slinear[blue][N] \slinear[blue][N] \rcol P_{\lesssim 1} v_{\leq N} \bigg]. 
\end{equation}
In contrast to the high$\times$low$\times$low-interactions from \eqref{ansatz:eq-heuristic-hll}, the high$\times$high$\times$low-interactions in \eqref{ansatz:eq-heuristic-hhl} benefits from multilinear dispersive effects. The reason is that, for frequencies $n_1$, $n_2$, and $n_3$ satisfying $|n_1|\sim |n_2| \sim N$ and $|n_3|\sim 1$, the multilinear dispersive symbols 
\begin{equation*}
(n_1,n_2,n_3) \mapsto \pm_{123} \langle n_{123} \rangle \pm_1 \langle n_1 \rangle \pm_2 \langle n_2 \rangle \pm_3 \langle n_3 \rangle 
\end{equation*}
are essentially equidistributed over all values $\lesssim N$. Based on this heuristic\footnote{For the precise details, we refer to Lemma \ref{para:lem-strichartz-XXtwo} below.}, we expect a gain of half of a derivative from multilinear dispersive effects and one derivative from the $\langle \nabla \rangle^{-1}$-multiplier. In total, we then expect that \eqref{ansatz:eq-heuristic-hhl} has regularity
\begin{equation*}
2 \cdot \big( -1/2- \big) + 1/2 + 1 = 1/2-. 
\end{equation*}
While the regularity of \eqref{ansatz:eq-heuristic-hhl} is therefore too low for a direct contraction argument, it exhibits a (different) useful random structure. Indeed, it morally behaves like the Duhamel integral of 
$\lcol \slinear[blue][N] \slinear[blue][N] \rcol $. The para-controlled component $\Xtwo_{\leq N}$, which will be introduced below, is used to capture \eqref{ansatz:eq-heuristic-hhl} and similar terms. \\

In addition to the high$\times$low$\times$low and high$\times$high$\times$low-interactions above, which yield contributions of regularity below $1/2$, there is a third problematic frequency-interaction. While its expected regularity is better than $1/2$, it still prevents us from directly closing a contraction argument. \\

\emph{Resonant-interaction:} For the discussion of the resonant interaction, we assume that the para-controlled components $\Xone_{\leq N}$ and $\Xtwo_{\leq N}$ have been defined and we now examine the nonlinear remainder
\begin{equation}\label{ansatz:eq-YN-heuristic} 
Y_{\leq N}:= v_{\leq N}-\Xone_{\leq N} - \Xtwo_{\leq N},
\end{equation}
 which is expected to have regularity $s>1/2$. In contrast to the para-controlled components, the nonlinear remainder $Y_{\leq N}$ exhibits no random structure and is treated using only deterministic methods. After inserting the para-controlled Ansatz into \eqref{ansatz:eq-NLW-vN-c}, we encounter the high$\times$(high$\times$high$\rightarrow$low)-interaction\footnote{For simplicity, we omitted the renormalization of $\slinear[blue][N] \slinear[blue][N]$, which is irrelevant in this term.}
\begin{equation}\label{ansatz:eq-heuristic-hhhtl}
\Duh \bigg[ \slinear[blue][N] P_{\lesssim 1} \Big( \slinear[blue][N] P_N Y_{\leq N} \Big) \bigg]. 
\end{equation}
In contrast to \eqref{ansatz:eq-heuristic-hhtll} above, the resonant product $ P_{\lesssim 1}\big( \slinear[blue][N] P_N Y_{\leq N} \big) $ in \eqref{ansatz:eq-heuristic-hhhtl} is well-defined, since $Y_{\leq N}$ has regularity greater than $1/2$. Our issues with \eqref{ansatz:eq-heuristic-hhhtl} are of a different kind, and instead concern only its low regularity. Due to the absence of multilinear dispersive effects in \eqref{ansatz:eq-heuristic-hhhtl}, its regularity is equal to 
\begin{equation*}
1 + \big( -1/2- \big) + \big( -1/2-\big) + s = s-. 
\end{equation*}
As a result, the regularity of  \eqref{ansatz:eq-heuristic-hhhtl} is (just barely) too low to close a contraction argument for $Y_{\leq N}$. Fortunately, \eqref{ansatz:eq-heuristic-hhhtl} exhibits a random structure, and will also be included in our first para-controlled component $\Xone_{\leq N}$. After this informal discussion, we now proceed with rigorous definitions. \\

In order to isolate the problematic frequency-interactions, we dyadically decompose the nonlinear terms in \eqref{ansatz:eq-NLW-vN-a}-\eqref{ansatz:eq-NLW-vN-e}. For this, we need dyadically-localized versions of the modified products from Definition \ref{ansatz:def-product}. The precise formulas for the dyadic components are rather complicated, and we encourage the reader to ignore these technicalities on first reading. 

\begin{definition}[Dyadically-localized modified product]\label{ansatz:def-dyadic-modified-product}
Let $N\geq 1$. For all $1\leq N_1, N_2 \leq N$, we first define
\begin{align}
\mathfrak{C}^{(1,5)}_{\leq N}[N_1,N_2](t) &:= \mathbb{E} \Big[ \slinear[blue][N_1] P_{N_2} \squintic[\leqN] \Big], \\
\mathfrak{C}^{(3,3)}_{\leq N}[N_1,N_2](t) &:= \mathbb{E} \Big[ P_{N_1} \scubic[\leqN] \, P_{N_2} \scubic[\leqN] \Big]. 
\end{align}
For all $\zeta^{(1)} \in \Symb^b$, $\zeta^{(2)},\zeta^{(3)} \in \Symb^b_0$, and frequency-scales $N,N_1,N_2,N_3$, we define the frequency-localized modified product as follows: 
\begin{alignat*}{3}  
&\Pi^\ast_{\leq N} \Big( P_{N_1} \slinear[blue][\leqN] , P_{N_2} \scubic[\leqN], P_{N_3} \scubic[\leqN]\Big)  
:=  
\slinear[blue][N_1] \Big( P_{N_2} \scubic[\leqN] P_{N_3} \scubic[\leqN]  -  \mathfrak{C}^{(3,3)}_{\leq N}[N_2,N_3]\Big), &&&\allowdisplaybreaks[3]\\[1ex]
&\Pi^\ast_{\leq N}\Big(  P_{N_1}  \slinear[blue][\leqN], P_{N_2} \squintic[\leqN], P_{N_3} \squintic[\leqN] \Big) \\
&:= 
\slinear[blue][N_1]  P_{N_2} \squintic[\leqN] P_{N_3} \squintic[\leqN] 
- \mathfrak{C}_{\leq N}^{(1,5)}[N_1,N_2] P_{N_3} \squintic[\leqN] 
- \mathfrak{C}_{\leq N}^{(1,5)}[N_1,N_3] P_{N_2} \squintic[\leqN],  
&&&\allowdisplaybreaks[3]\\[1ex] 
&\Pi^\ast_{\leq N} \Big( P_{N_1}   \slinear[blue][\leqN], P_{N_2} \squintic[\leqN],  P_{N_3} \zeta^{(3)}_{\leq N} \Big) 
:= 
\slinear[blue][N_1]  P_{N_2} \squintic[\leqN]  P_{N_3} \zeta^{(3)}_{\leq N} -  \mathfrak{C}_{\leq N}^{(1,5)}[N_1,N_2]  P_{N_3} \zeta^{(3)}_{\leq N}\qquad \qquad &\hspace{-15ex}\text{if }  \zeta^{(3)} \neq  \squintic[], &&\allowdisplaybreaks[3]\\[1ex]
&\Pi^\ast_{\leq N} \Big( P_{N_1}   \slinear[blue][\leqN], P_{N_2} \scubic[\leqN],  P_{N_3} \zeta^{(3)}_{\leq N} \Big) 
:= 
\slinear[blue][N_1]  P_{N_2}\scubic[\leqN]  P_{N_3}\zeta^{(3)}_{\leq N} \qquad \qquad  &\hspace{-15ex}\text{if } \zeta^{(3)} \neq  \scubic[], \squintic[],  &&\allowdisplaybreaks[3]\\[1ex]  
&\Pi^\ast_{\leq N} \Big(  P_{N_1}  \slinear[blue][\leqN], P_{N_2}  \zeta^{(2)}_{\leq N}, P_{N_3} \zeta^{(3)}_{\leq N} \Big) 
:= 
\slinear[blue][N_1]   P_{N_2} \zeta^{(2)}_{\leq N} P_{N_3} \zeta^{(3)}_{\leq N}\qquad \qquad  &\hspace{-15ex}\text{if } \zeta^{(2)},  \zeta^{(3)} \neq  \scubic[], \squintic[],  &&\allowdisplaybreaks[3]\\[1ex] 
&\Pi^\ast_{\leq N} \Big( P_{N_1}   \scubic[\leqN], P_{N_2} \scubic[\leqN], P_{N_3} \scubic[\leqN]\Big) \\
&:=  P_{N_1} \scubic[\leqN] P_{N_2} \scubic[\leqN]  P_{N_3} \scubic[\leqN]  
-  \mathfrak{C}^{(3,3)}_{\leq N}[N_2,N_3] P_{N_1} \scubic[\leqN]
-  \mathfrak{C}^{(3,3)}_{\leq N}[N_1,N_3] P_{N_2} \scubic[\leqN]
-  \mathfrak{C}^{(3,3)}_{\leq N}[N_1,N_2] P_{N_3} \scubic[\leqN], \hspace{-10cm}&&&\allowdisplaybreaks[3]\\[1ex] 
&\Pi^\ast_{\leq N} \Big( P_{N_1}   \scubic[\leqN], P_{N_2} \scubic[\leqN], P_{N_3} \zeta^{(3)}_{\leq N}\Big) 
:= \Big( P_{N_1} \scubic[\leqN] P_{N_2} \scubic[\leqN]  - \mathfrak{C}^{(3,3)}_{\leq N}[N_1,N_2] \Big) P_{N_3} \zeta^{(3)}_{\leq N} \qquad \qquad  &\hspace{-15ex}\text{if }  \zeta^{(3)} \neq   \scubic[], &&\allowdisplaybreaks[3]\\[1ex] 
&\Pi^\ast_{\leq N} \Big(  P_{N_1}  \scubic[\leqN], P_{N_2} \zeta^{(2)}_{\leq N}, P_{N_3}  \zeta^{(3)}_{\leq N}\Big) 
:=  P_{N_1} \scubic[\leqN]  P_{N_2} \zeta^{(2)}_{\leq N} P_{N_3} \zeta^{(3)}_{\leq N} \qquad \qquad  &\hspace{-15ex}\text{if } \zeta^{(2)}, \zeta^{(3)} \neq   \scubic[], &\allowdisplaybreaks[3]\\[1ex] 
&\Pi^\ast_{\leq N} \Big( P_{N_1}   \zeta^{(1)}_{\leq N}, P_{N_2}  \zeta^{(2)}_{\leq N}, P_{N_3} \zeta^{(3)}_{\leq N}\Big)
:=   P_{N_1} \zeta^{(1)}_{\leq N} P_{N_2} \zeta^{(2)}_{\leq N}  P_{N_3} \zeta^{(3)}_{\leq N} \qquad \qquad  &\hspace{-15ex}\text{if } \zeta^{(1)} \neq  \slinear[blue], \, \, \zeta^{(1)}, \zeta^{(2)},  \zeta^{(3)}\neq   \scubic[]. &&
\end{alignat*}
\end{definition}

Equipped with Definition \ref{ansatz:def-dyadic-modified-product}, we now define the following trilinear para-product operators, which capture the bad interactions from our informal discussion. In contrast to Definition \ref{ansatz:def-product} and Definition \ref{ansatz:def-dyadic-modified-product}, we allow for two linear evolutions in our arguments. The reason is that, once we restrict to specific frequency-interactions, fewer of the renormalization or cancellation-terms in \eqref{ansatz:eq-NLW-vN-a} and \eqref{ansatz:eq-NLW-vN-c} are relevant. 

\begin{definition}[Trilinear para-product operators]\label{ansatz:def-para-product}
We make the following three definitions:
\begin{enumerate}
    \item High$\times$low$\times$low-interactions: For $\zeta^{(2)},\zeta^{(3)} \in \Symb_0^b$, we define 
    \begin{equation*}
    \HLL \Big( \slinear[blue][\leqN], \zeta^{(2)}_{\leq N} , \zeta^{(3)}_{\leq N} \Big)
    := \sum_{\substack{N_1,N_2,N_3 \leq N\colon \\ N_2,N_3 \leq N_1^\eta}} \Pi^{\ast}_{\leq N}\Big( \slinear[blue][N_1], P_{N_2} \zeta^{(2)}_{\leq N}, P_{N_3} \zeta^{(3)}_{\leq N} \Big). 
    \end{equation*}
    Furthermore, we also define 
    \begin{align*}
    \HLL \Big( \slinear[blue][\leqN], \slinear[blue][\leqN], \squintic[\leqN] \Big) 
    &:= \sum_{\substack{N_1,N_2,N_3 \leq N \colon \\ N_2,N_3 \leq N_1^\eta}} 
    \slinear[blue][N_1] \Big( \slinear[blue][N_2] P_{N_3} \squintic[\leqN] - \mathfrak{C}^{(1,5)}_{\leq N}[N_2,N_3] \Big), \\ 
       \HLL \Big( \slinear[blue][\leqN], \slinear[blue][\leqN], v_{\leq N}\Big) 
    &:=\sum_{\substack{N_1,N_2,N_3 \leq N \colon \\ N_2,N_3 \leq N_1^\eta}} 
    \slinear[blue][N_1] \slinear[blue][N_2] P_{N_3}v_{\leq N}. 
    \end{align*}
    \item High$\times$high$\times$low-interactions: For $\zeta^{(3)}\in \big\{ \squintic[], v \big\}$, we define 
    \begin{equation*}
    \HHL \Big( \slinear[blue][\leqN], \slinear[blue][\leqN], \zeta^{(3)}_{\leq N} \Big) 
    = \sum_{\substack{N_1,N_2,N_3\leq N\colon \\ \min(N_1,N_2) > \max(N_1,N_2)^\eta \\ N_3 \leq \max(N_1,N_2)^\eta }} \lcol \slinear[blue][N_1] \slinear[blue][N_2] \rcol P_{N_3} \zeta^{(3)}_{N_3}.
    \end{equation*}
    \item Resonant-interaction: For any function $Y\colon \R \times \T^3 \rightarrow \R$, we define \begin{align*}
    &\RES \Big( \slinear[blue][\leqN], \slinear[blue][\leqN], P_{\leq N} Y \Big) \\
    =& \sum_{\substack{N_1,N_2,N_3 \leq N \colon \\ N_3 > \max(N_1,N_2)^\eta }}
     \sum_{N_{13},N_{23}} \sum_{n_1,n_2,n_3 \in \Z^3} \bigg[ 
     \Big( \mathbf{1}\{ N_{13}\leq N_2^\eta \} + \mathbf{1}\{ N_{23} \leq N_1^\eta \} \Big) 
     \Big( \prod_{j=1}^3 1_{N_j}(n_j) \Big) \\
     &\times 1_{N_{13}}(n_{13}) 1_{N_{23}}(n_{23}) \, \slinear[blue][N_1](n_1) \slinear[blue][N_2](n_2) P_{N_3} \widehat{Y}(n_3) e^{i \langle n_{123}, x \rangle} \bigg].
    \end{align*}
\end{enumerate}
\end{definition}

\begin{remark}[On $\RES$]
In the definition of the resonant-interaction $\RES$, we omitted the Wick-ordering of $\slinear[blue][N_1] \slinear[blue][N_2]$ but double-count the case when both $N_{13}\leq N_2^\eta$ and $N_{23}\leq N_1^\eta$. Both of these choices, which simplify the notation in our treatment of $\Xone$ below, are made primarily for convenience. As it turns out, both contributions essentially correspond to terms with at least two resonances between $n_1$, $n_2$, and $n_3$ (instead of just one resonance), which are easily estimated.
\end{remark}

Equipped with the tri-linear para-product operators from Definition \ref{ansatz:def-para-product}, we can now state the para-controlled equations. As dictated by \eqref{ansatz:eq-YN-heuristic}, we decompose
\begin{equation}
v_{\leq N} = \Xone_{\leq N} + \Xtwo_{\leq N} + Y_{\leq N}. 
\end{equation}
The notation for the para-controlled components is chosen in agreement with \cite{DNY20}, since 
$\Xone_{\leq N}$ and  $\Xtwo_{\leq N}$ are related to $(1,1)$ and $(2,1)$-tensors, respectively. We first define the evolution equation for $\Xone_{\leq N}$ by 
\begin{align}
    &(\partial_t^2 + 1 - \Delta) \Xone_{\leq N} \notag \\
    =&- 6 P_{\leq N} \HLL \Big( \slinear[blue][\leqN], \slinear[blue][\leqN], 3 \squintic[\leqN] -P_{\leq N} \slinear[green][\leqM]+ v_{\leq N} \Big) \label{ansatz:eq-X11-a}\\
    +& \sum_{\zeta^{(2)},\zeta^{(3)} \in \Symb_0^b} \Athree\big( \slinear[blue],\zeta^{(2)}, \zeta^{(3)}\big) P_{\leq N}\HLL \Big( \slinear[blue][\leqN], \zeta^{(2)}_{\leq N}, \zeta^{(3)}_{\leq N} \Big) \label{ansatz:eq-X11-b} \\
    +& \Aone\big(\, \slinear[blue] \big) \mathfrak{C}_{\leq N} P_{\leq N} \slinear[blue][\leqN] + \Atone\big(\,\slinear[blue]\big) (\gamma_{\leq N}- \Gamma_{\leq N}) P_{\leq N} \slinear[blue][\leqN] \label{ansatz:eq-X11-c} \\
    -& 3 P_{\leq N}\RES \Big( \slinear[blue][\leqN], \slinear[blue][\leqN], Y_{\leq N} \Big)
    \label{ansatz:eq-X11-d}
\end{align}
and impose the initial condition 
\begin{equation}\label{ansatz:eq-X11-initial}
\Xone_{\leq N}(0) = \partial_t \Xone_{\leq N}(0) = 0. 
\end{equation}
Secondly, we define the evolution equation for $\Xtwo_{\leq N}$ by 
\begin{equation}\label{ansatz:eq-X21}
(\partial_t^2 + 1 - \Delta) \Xtwo_{\leq N} = - 3 P_{\leq N} \HHL \Big( \, \slinear[blue][\leqN], \slinear[blue][\leqN] , 3 \squintic[\leqN] -P_{\leq N} \slinear[green][\leqM] + v_{\leq N} \Big)
\end{equation}
and impose the initial condition
\begin{equation}\label{ansatz:eq-X21-initial}
    \Xtwo_{\leq N}(0) = \partial_t \Xtwo_{\leq N}(0) = 0. 
\end{equation}
It remains to state the evolution equation for $Y_{\leq N}$, which is essentially determined by the evolution equations for $v_{\leq N}$, $\Xone_{\leq N}$, and $\Xtwo_{\leq N}$. In the evolution equation for $Y_{\leq N}$, however, we will make use of the para-controlled structure of $\Xone_{\leq N}$ and $\Xtwo_{\leq N}$ through the double Duhamel trick. That is, we will be inserting the Duhamel integral of the right-hand sides \eqref{ansatz:eq-X11-a}-\eqref{ansatz:eq-X11-d} and \eqref{ansatz:eq-X21} instead of $\Xone_{\leq N}$ and $\Xtwo_{\leq N}$, respectively. In order to emphasize the double Duhamel trick in our notation, we now introduce the following two operators.

\begin{definition}[Para-controlled operators]\label{ansatz:def-operators-para}
We define two operators $\XXone$ and $\XXtwo$ by 
\begin{align}
\XXone &= \XXone [v_{\leq N}, Y_{\leq N}] = 
\Duh \Big[ \eqref{ansatz:eq-X11-a} +  \eqref{ansatz:eq-X11-b} +  \eqref{ansatz:eq-X11-c} +  \eqref{ansatz:eq-X11-d} \Big], \\
\XXtwo &= \XXtwo [v_{\leq N}, Y_{\leq N}] = 
\Duh \Big[ \eqref{ansatz:eq-X21} \Big]. 
\end{align}
\end{definition}
From the definition, it is clear that $\XXone$ can be written as the sum of explicit stochastic objects, linear operators in $v_{\leq N}$ and $Y_{\leq N}$, and bilinear operators in $v_{\leq N}$. Similarly, $\XXtwo$ can be written as a sum of explicit stochastic objects and linear operators in $v_{\leq N}$.  Before stating the evolution equation for $Y_{\leq N}$, we also expand the symbols from Definition \ref{ansatz:def-symbols}. 

\begin{definition}[Para-controlled symbols]\label{ansatz:def-symbols-para}
We define the set of para-controlled symbols by 
\begin{equation}
\Symb^p := \Big \{ \slinear[blue], \, \scubic, \, \slinear[green][\leqM],\, \squintic, \, \SXXone, \, \SXXtwo, Y \Big\}. 
\end{equation}
Furthermore, we also define 
\begin{equation}
\Symb^p_0 := \Big \{ \scubic, \, \slinear[green][\leqM],\, \squintic, \, \SXXone, \, \SXXtwo, Y \Big\}. 
\qquad \text{and} \qquad 
\Symb^p_{1/2} := \Big \{  \slinear[green][\leqM],\, \squintic, \, \SXXone, \, \SXXtwo, Y \Big\}. 
\end{equation}
\end{definition}
Similar as for $\Symb_{1/2}^b$, the set $\Symb_{1/2}^p$ will only occur in the proofs but not the statements of our main estimates (see Section \ref{section:proof-main-estimates}).

Equipped with both the para-controlled operators (Definition \ref{ansatz:def-operators-para}) and symbols (Definition \ref{ansatz:def-symbols-para}), we can now write the evolution equation for the smooth nonlinear remainder $Y_{\leq N}$. It is given by 
\begin{align}
    &(\partial_t^2 + 1 - \Delta) Y_{\leq N} \notag\\
    =& - P_{\leq N} \bigg[ 9 \squadratic[\leqN] \squintic[\leqN] - \Gamma_{\leq N}  \scubic[\leqN]  - 18 \mathfrak{C}^{(1,5)}_{\leq N} P_{\leq N} \slinear[blue][\leqN]  
    - 9 \Big( 2 \HLL + \HHL \Big) \Big( \, \slinear[blue][\leqN], \slinear[blue][\leqN], \squintic[\leqN] \Big) \bigg] \allowdisplaybreaks[3]  \label{ansatz:eq-YN-1} \\
    +&  3 P_{\leq N}  \bigg[  \squadratic[\leqN] P_{\leq N} \slinear[green][\leqM] - \Big( 2 \HLL + \HHL \Big) \Big( \slinear[blue][\leqN], \slinear[blue][\leqN], P_{\leq N} \slinear[green][\leqM] \Big) \bigg] \label{ansatz:eq-YN-2} \\
    -&3 P_{\leq N} \bigg[ \squadratic[\leqN] \XXone - 
    \Big( 2 \HLL + \HHL \Big) \Big( \slinear[blue][\leqN] , \slinear[blue][\leqN], \XXone \Big) \bigg] \allowdisplaybreaks[3] \label{ansatz:eq-YN-3} \\
    -& P_{\leq N} \bigg[ 3 \squadratic[\leqN] \XXtwo - 
    \Big( 6 \HLL + 3 \HHL \Big) \Big( \slinear[blue][\leqN] , \slinear[blue][\leqN], \XXtwo \Big)
    + \Gamma_{\leq N} \Big( 3 \squintic[\leqN] - P_{\leq N} \slinear[green][\leqM]+ v_{\leq N} \Big) \bigg] \allowdisplaybreaks[3] \label{ansatz:eq-YN-4} \\
    -& 3 P_{\leq N} \bigg[ \squadratic[\leqN] Y_{\leq N} - 
    \Big( 2 \HLL + \HHL + \RES \Big) \Big( \slinear[blue][\leqN] , \slinear[blue][\leqN], Y_{\leq N} \Big) \bigg] \allowdisplaybreaks[3] \label{ansatz:eq-YN-5} \\
    +& \sum_{\zeta^{(2)},\zeta^{(3)} \in \Symb_0^p } \Athree\big(\, \slinear, \zeta^{(2)}, \zeta^{(3)} \big)  P_{\leq N} \Big( \Pi^\ast_{\leq N}- \HLL \Big) \Big( \slinear[blue][\leqN], \zeta^{(2)}_{\leq N}, \zeta^{(3)}_{\leq N} \Big) \allowdisplaybreaks[3] \label{ansatz:eq-YN-6} \\
    +& \sum_{\zeta^{(1)},\zeta^{(2)},\zeta^{(3)} \in \Symb_0^p} \Athree \big( \zeta^{(1)}, \zeta^{(2)}, \zeta^{(3)} \big) P_{\leq N} \Pi^\ast_{\leq N} \Big( P_{\leq N} \zeta^{(1)}, P_{\leq N} \zeta^{(2)}, P_{\leq N} \zeta^{(3)}  \Big) \allowdisplaybreaks[3] \label{ansatz:eq-YN-7} \\
    +& \mathfrak{C}_{\leq N} \sum_{\zeta \in \Symb_0^b} \Aone(\zeta) P_{\leq N} \zeta_{\leq N}
    + (\gamma_{\leq N}-\Gamma_{\leq N}) \sum_{\zeta \in \Symb_0^b} \Atone(\zeta) P_{\leq N} \zeta_{\leq N}. \label{ansatz:eq-YN-8}
    \end{align}
    
    From the derivation of the Ansatz, we obtain the following reformulation of the nonlinear wave equation \eqref{ansatz:eq-frequency-truncated-NLW}. 
    
    \begin{proposition}[Para-controlled reformulation]
    Let $M\geq N \geq 1$ and let $(v_{\leq N},\Xone_{\leq N},\Xtwo_{\leq N},Y_{\leq N})$ satisfy
    \begin{equation*}
        v_{\leq N} = \Xone_{\leq N} + \Xtwo_{\leq N} + Y_{\leq N}
    \end{equation*}
    and be a solution of the para-controlled nonlinear wave equations \eqref{ansatz:eq-X11-a}-\eqref{ansatz:eq-X21-initial} and \eqref{ansatz:eq-YN-1}-\eqref{ansatz:eq-YN-8}. Then, 
    \begin{equation*}
    u_{\leq N} = \slinear[blue] - \scubic[\leqN] - \slinear[green][\leqM] + 3 \squintic[\leqN] + \Xone_{\leq N} + \Xtwo_{\leq N} + Y_{\leq N}
    \end{equation*}
    solves the nonlinear wave equation \eqref{ansatz:eq-frequency-truncated-NLW}.
    \end{proposition}
    
    \begin{remark}\label{ansatz:rem-para-formulation}
    In the following, we often refer to $(v_{\leq N},\Xone_{\leq N},\Xtwo_{\leq N},Y_{\leq N})$ as the solutions of the para-controlled nonlinear wave equations.
    \end{remark}
    
    \subsection{Main estimates}\label{section:main-estimates} 
    In this subsection, we state the main estimates of our local well-posedness theory. For expository purposes, the estimates are split over four separate propositions. In the first proposition, we state probabilistic Strichartz and regularity estimates for  $\XXone$ and $\XXtwo$. In the remaining three propositions, we state our estimates of the nonlinear terms in the evolution equation for $Y_{\leq N}$. The three propositions treat \eqref{ansatz:eq-YN-1}-\eqref{ansatz:eq-YN-5}, \eqref{ansatz:eq-YN-6}, and \eqref{ansatz:eq-YN-7}-\eqref{ansatz:eq-YN-8}, respectively. In other words, we distinguish terms which contain two, one, or zero linear stochastic objects. 
    
    \begin{proposition}[Probabilistic Strichartz and regularity estimates for $\XXone$ and $\XXtwo$]\label{ansatz:prop-Xj} For all $A\geq 1$, there exists an $A$-certain event $\Event_A \in \mathcal{E}$ such that the following estimates hold: For all $j=1,2$, $N\geq 1$,  $T\geq 1$, and all closed intervals $0 \in \Jc \subseteq [-T,T]$, we have that
    \begin{equation}\label{ansatz:eq-Xj}
    \begin{aligned}
    &\Big\| \XXj[v_{\leq N}, Y_{\leq N}] \Big\|_{X^{1/2-\delta_1,b}(\Jc)} 
    + \Big\| \XXj[v_{\leq N}, Y_{\leq N}] \Big\|_{L_t^\infty \C_x^{1/2-\delta_1}(\Jc)} \\
    &\leq A T^\alpha \Big( 1 + \| v_{\leq N} \|_{X^{-1,b}(\Jc)}^2 + \| Y_{\leq N} \|_{X^{1/2+\delta_2,b}(\Jc)} \Big).
    \end{aligned}
    \end{equation}
    \end{proposition}
    
    \begin{remark}[On the $X^{-1,b}$-norm]
    The right-hand side of \eqref{ansatz:eq-Xj} only contains the $X^{-1,b}$-norm of $v_{\leq N}$. In our local theory, this aspect of Proposition \ref{ansatz:prop-Xj} is rather irrelevant, and \eqref{ansatz:eq-Xj} would be just as useful if the $X^{-1,b}$-norm is replaced by the much stronger $X^{1/2-\delta_1,b}$-norm. In the global theory (see Section \ref{section:global}), however, the bound by weaker norms is essential. In fact, our globalization bounds (Proposition \ref{global:prop-para-controlled}) are based on a Gronwall-type inequality, which  only works if there is at most one high-regularity norm on the right-hand side of  \eqref{ansatz:eq-Xj}. In contrast, Proposition \ref{ansatz:prop-Y-one} and Proposition \ref{ansatz:prop-Y-no} utilize more than one high-regularity norm, and will be revisited below (see Lemma \ref{global:lem-one-zero-revisited}). 
    \end{remark} 
    
    We now turn to the main estimates for the smooth nonlinear remainder $Y_{\leq N}$. After each proposition, we point the reader to the sections of our article which are most relevant for the proof. We first start with terms involving two linear stochastic objects, i.e, \eqref{ansatz:eq-YN-1}-\eqref{ansatz:eq-YN-5}. 
    
    \begin{proposition}[Terms involving two linear stochastic objects]\label{ansatz:prop-Y-two}
    For all $A\geq 1$, there exists an $A$-certain event $\Event_A \in \mathcal{E}$ on which the following estimates hold for all $M \geq N\geq 1$, $T\geq 1$, and closed intervals $0\in \Jc \subseteq [-T,T]$.  
    \begin{enumerate}[label=(\roman*)]
        \item\label{item:prop-Y-two-1} (Explicit stochastic objects) It holds that 
        \begin{align*}
            &\bigg\| P_{\leq N} \bigg[ 9 \squadratic[\leqN] \squintic[\leqN] - \Gamma_{\leq N} \Big( \scubic[\leqN] + \slinear[green][\leqN] \Big) - 18 \mathfrak{C}^{(1,5)}_{\leq N} P_{\leq N} \slinear[blue][\leqN]  \\
            &\hspace{2ex} 
   - 9 \Big( 2 \HLL + \HHL \Big) \Big( \, \slinear[blue][\leqN], \slinear[blue][\leqN], \squintic[\leqN] \Big) \bigg] \bigg\|_{X^{-1/2+\delta_2,b_+-1}(\Jc)} \\
    +&\bigg\| P_{\leq N}  \bigg[  \squadratic[\leqN] P_{\leq N} \slinear[green][\leqM] - \Big( 2 \HLL + \HHL \Big) \Big( \slinear[blue][\leqN], \slinear[blue][\leqN], P_{\leq N} \slinear[green][\leqM] \Big) \bigg] \bigg\|_{X^{-1/2+\delta_2,b_+-1}(\Jc)}  \\
    \leq& \, A T^\alpha. 
        \end{align*}
        \item \label{item:prop-Y-two-2} (Para-controlled calculus) It holds that 
        \begin{align*}
        &\bigg\| \squadratic[\leqN] \XXone - 
     \Big( 2 \HLL + \HHL \Big) \Big( \slinear[blue][\leqN] , \slinear[blue][\leqN], \XXone \Big) 
     \bigg\|_{X^{-1/2+\delta_2,b_+-1}(\Jc)} \\
    +& \bigg\|  3 \squadratic[\leqN] \XXtwo - 
    \Big( 6 \HLL + 3 \HHL \Big) \Big( \slinear[blue][\leqN] , \slinear[blue][\leqN], \XXtwo \Big)
     \\
     &\hspace{2ex}+ \Gamma_{\leq N} \Big( 3 \squintic[\leqN]- P_{\leq N} \slinear[green][\leqM] + v_{\leq N} \Big)  \bigg\|_{X^{-1/2+\delta_2,b_+-1}(\Jc)} \\
    \leq& \, A T^\alpha \Big( 1+ \| v_{\leq N} \|_{X^{-1,b}(\Jc)}^2 + \| Y_{\leq N} \|_{X^{1/2+\delta_2,b}(\Jc)} \Big).
        \end{align*}
        \item \label{item:prop-Y-two-3} (The $Y$-term) It holds that 
        \begin{align*}
            &\bigg\|  \squadratic[\leqN] Y_{\leq N} - 
    \Big( 2 \HLL + \HHL + \RES \Big) \Big( \slinear[blue][\leqN] , \slinear[blue][\leqN], Y_{\leq N} \Big) \bigg\|_{X^{-1/2+\delta_2,b_+-1}(\Jc)} \\ 
            \leq& \, A T^\alpha \| Y_{\leq N} \|_{X^{1/2+\delta_2,b}(\Jc)}. 
        \end{align*}
    \end{enumerate}
    \end{proposition}
    
    The explicit stochastic objects in  \ref{item:prop-Y-two-1} are studied from an algebraic perspective in Section \ref{section:diagrams} and estimated in Section \ref{section:analytic} and Section \ref{section:analytic2}. The para-controlled calculus in \ref{item:prop-Y-two-2}, which concerns random operators containing $\XXone$ and $\XXtwo$, is performed in Section \ref{section:para}. This estimate also heavily relies on operator-norm estimates of our tensors (Section \ref{section:counting}). Finally, the linear random operator in $Y_{\leq N}$ from \ref{item:prop-Y-two-3} is treated in Subsection \ref{section:quad}.

    \begin{proposition}[Terms involving one linear stochastic object]\label{ansatz:prop-Y-one}
    For all $A\geq 1$, there exists an $A$-certain event $\Event_A \in \mathcal{E}$ such that the following estimates hold: For all $N\geq 1$, $T\geq 1$, closed intervals $0 \in \Jc \subseteq [-T,T]$, and $\zeta^{(2)},\zeta^{(3)} \in \Symb_0^p$, it holds that 
    \begin{equation}\label{ansatz:eq-Y-one}
    \begin{aligned}
    &\bigg\| \Big( \Pi^\ast_{\leq N}- \HLL \Big) \Big( \slinear[blue][\leqN], \zeta^{(2)}_{\leq N}, \zeta^{(3)}_{\leq N} \Big) \bigg\|_{X^{-1/2+\delta_2,b_+-1}(\Jc)}  \\
    \leq& \, A  T^\alpha \Big( 1+ \| v_{\leq N}\|_{X^{-1,b}}^4 + \| Y_{\leq N} \|_{X^{1/2+\delta_2,b}(\Jc)}^2 \Big). 
    \end{aligned}
    \end{equation}
    \end{proposition}
    
    The argument leading to \eqref{ansatz:eq-Y-one} primarily depends on whether two, one, or none of the arguments $\zeta^{(2)}$ and $\zeta^{(3)}$ coincide with the cubic stochastic object $\scubic$. The three cases are treated using either explicit stochastic objects (Section \ref{section:analytic}), linear random operators (Subsection \ref{section:lincub}), or bilinear random operators (Section \ref{section:bilinear}), respectively. 
    
    \begin{proposition}[Terms involving no linear stochastic object]\label{ansatz:prop-Y-no}
    For all $A\geq 1$, there exists an $A$-certain event $\Event_A \in \mathcal{E}$ such that the following estimates hold: For all $N\geq 1$, $\zeta^{(1)},\zeta^{(2)}, \zeta^{(3)} \in \Symb_0^p$, and $\zeta \in \Symb_0^b$, it holds that 
    \begin{equation}\label{ansatz:eq-Y-no-1}
    \begin{aligned}
        &\bigg\| \Pi^\ast_{\leq N} \Big( P_{\leq N} \zeta^{(1)}, P_{\leq N} \zeta^{(2)}, P_{\leq N} \zeta^{(3)}  \Big) \bigg\|_{X^{-1/2+\delta_2,b_+-1}(\Jc)} \\
        \leq&\,  A T^\alpha \Big( 1 + \| v_{\leq N} \|_{X^{-1,b}(\Jc)}^6 + \| Y_{\leq N} \|_{X^{1/2+\delta_2,b}(\Jc)}^3 \Big) 
    \end{aligned}
    \end{equation}
    and 
    \begin{align}
    \Big\| \mathfrak{C}_{\leq N} \zeta_{\leq N} \Big\|_{X^{-1/2+\delta_2,b_+-1}(\Jc)} &\leq A T^\alpha\Big( 1 + \| v_{\leq N} \|_{X^{1/2-\delta_1,b}(\Jc)} \Big), \label{ansatz:eq-Y-no-2}\\ 
    \Big\| \big( \gamma_{\leq N} - \Gamma_{\leq N} \big) \zeta_{\leq N} \Big\|_{X^{-1/2+\delta_2,b_+-1}(\Jc)} &\leq A T^\alpha \big( 1+ \| v_{\leq N} \|_{X^{1/2-\delta_1,b}(\Jc)} \big). 
    \label{ansatz:eq-Y-no-3}
    \end{align}
    \end{proposition}
    
    The most difficult case in Proposition \ref{ansatz:prop-Y-no} is the first estimate \eqref{ansatz:eq-Y-no-1} for $\zeta^{(1)}=\zeta^{(2)}=\zeta^{(3)}=\scubic$. As was already discussed in the introduction, the corresponding nonic stochastic object is rather involved and estimated using the molecules from \cite{DH21}. 
    
        \subsection{Proof of local well-posedness}\label{section:ansatz-proof-lwp}
    
    In this subsection, we prove the qualitative local well-posedness (Proposition \ref{ansatz:prop-lwp-qualitative}) and the nonlinear smoothing estimate (Proposition \ref{ansatz:prop-nonlinear-smoothing}). To this end, 
    we first state the following quantitative local well-posedness result.

    \begin{proposition}[Quantitative local well-posedness]\label{ansatz:prop-lwp-quantitative}
    Let $A\geq 1$ and $M\geq 1$. Then, there exists an event $E_{M,A} \in \mathcal{E}\otimes \mathcal{Z}$ such that the following properties hold: 
    \begin{enumerate}[label=(\roman*)]
    \item (High probability) \label{ansatz:item-quantitative-1} The event $E_{M,A}$ is $A$-certain with respect to $\mathbb{P}\otimes \mathbb{Q}$. 
    \item (Para-controlled solutions) \label{ansatz:item-quantitative-2} For all $0<\tau\leq A^{-\Theta}$ and all $N\leq M$, the solution $u_{\leq N}$ of the caloric initial value problem \eqref{ansatz:eq-frequency-truncated-NLW} and the solutions $(v_{\leq N}, \Xone_{\leq N}, \Xtwo_{\leq N}, Y_{\leq N})$ of the para-controlled nonlinear wave equations (as in Remark \ref{ansatz:rem-para-formulation}) exist on $[-\tau,\tau]\times \T^3$. Furthermore, they satisfy the estimates 
    \begin{alignat*}{3}
    \big\|  u_{\leq N} \big\| _{X^{-1/2-\epsilon,b}([-\tau,\tau])} \leq&\, CA, 
    &\quad 
    \big\|  v_{\leq N} \big\| _{X^{1/2-\delta_1,b}([-\tau,\tau])} \leq& \, CA, \\
     \max_{j=1,2} \big\|  \Xj_{\leq N} \big\|_{(L_t^\infty \C_x^{1/2-\delta_1} \cap X^{1/2-\delta_1,b})([-\tau,\tau])} \leq& \, CA, 
    &\qquad 
       \big\|  Y_{\leq N} \big\| _{X^{1/2+\delta_2,b}([-\tau,\tau])} \leq& \, CA. 
    \end{alignat*}
    \item (Difference estimates) \label{ansatz:item-quantitative-3} For all $0<\tau \leq A^{-\Theta}$ and all $N_1,N_2 \leq M$, the differences satisfy 
    \begin{align*}
    \big\| u_{\leq N_1} - u_{\leq N_2} \big\| _{X^{-1/2-\epsilon,b}([-\tau,\tau])} &\leq \min\big(N_1,N_2\big)^{-\theta}, \\
    \big\| v_{\leq N_1} - v_{\leq N_2} \big\| _{X^{1/2-\delta_1,b}([-\tau,\tau])} &\leq \min\big(N_1,N_2\big)^{-\theta}, \\
    \max_{j=1,2} \big\| \Xj_{\leq N_1} - \Xj_{\leq N_2} \big\| _{(L_t^\infty \C_x^{1/2-\delta_1} \cap X^{1/2-\delta_1,b})([-\tau,\tau])} &\leq \min\big(N_1,N_2\big)^{-\theta}, \\
    \big\| Y_{\leq N_1} - Y_{\leq N_2} \big\| _{X^{1/2+\delta_2,b}([-\tau,\tau])} &\leq \min\big(N_1,N_2\big)^{-\theta}. 
    \end{align*}
    \end{enumerate}
    \end{proposition}
    
    To avoid confusion, we note that the constant $C$ in Proposition \ref{ansatz:prop-lwp-quantitative} cannot be removed by re-adjusting $A$. The reason is that the remainder $\initial{red}[][\leqM]=\initial{red}[][\leqM](A,\phi^{\cos}_{\leq M})$ of the caloric representation is $A$-dependent, which would then have to be re-adjusted as well.

    Before we prove the quantitative version (Proposition \ref{ansatz:prop-lwp-quantitative}), we prove that it implies the qualitative version (Proposition \ref{ansatz:prop-lwp-qualitative}). 
    
    \begin{proof}[Proof of Proposition \ref{ansatz:prop-lwp-qualitative}] As discussed in Remark \ref{ansatz:rem-ucal}, we now introduce the solutions $u_{\leq M}^{\phi}$ and $u_{\leq M}^{\textup{cal}}$ of the general and caloric initial value problem, respectively. To this end,
    let $M\geq 1$ be arbitrary, let $(\mathfrak{Z},\mathcal{Z},\mathbb{Q})$ be a probability space, and let $(\phi^{\cos}_{\leq M},\phi^{\sin}_{\leq M})$ be random functions satisfying 
    \begin{equation*}
    \Law_{\mathbb{Q}} \Big( \big(\phi^{\cos}_{\leq M}, \phi^{\sin}_{\leq M}\big) \Big) = \mu_{\leq M}. 
    \end{equation*}
    Furthermore, let $A\geq 1$ remain to be chosen and let $\initial{blue}$, $\initial{green}[][\leqM]$, and $\initial{red}[][\leqM]=\initial{red}[][\leqM]\big( A,\phi^{\cos}_{\leq M}\big)$ be the caloric random initial data from Proposition \ref{ansatz:prop-caloric}. In particular, it holds that 
    \begin{equation}\label{ansatz:eq-lwp-qualitative-p1}
    \Law_{\mathbb{P}\otimes \mathbb{Q}}\Big( 
    \initial{blue}-\initial{green}[][\leqM]+\initial{red}[][\leqM] \Big)
    = \Law_{\mathbb{Q}} \Big( \big(\phi^{\cos}_{\leq M}, \phi^{\sin}_{\leq M}\big) \Big) = \mu_{\leq M}. 
    \end{equation}
   Let  $N\leq M$ be a frequency-truncation parameter. We denote by $u_{\leq N}^{\phi}$ the unique global solution of  \eqref{ansatz:eq-frequency-truncated-NLW-raw}, where the initial data $(\phi^{\cos}, \phi^{\sin}\big)$ denotes any element of the Sobolev space $\mathscr{H}_x^{-1/2-\epsilon}$. Furthermore, we denote by $u_{\leq N}^{\textup{cal}}$ the unique global solution of \eqref{ansatz:eq-frequency-truncated-NLW}, where the initial data is of the caloric form $\initial{blue}-\initial{green}[][\leqM]+\initial{red}[][\leqM]$. 
    
   Since solutions of the frequency-truncated nonlinear wave equation \eqref{ansatz:eq-frequency-truncated-NLW-raw} are uniquely determined by (and depend continuously on) their initial data, \eqref{ansatz:eq-lwp-qualitative-p1} implies that 
    \begin{equation}\label{ansatz:eq-lwp-qualitative-p2}
    \begin{aligned}
    &\mu_{\leq M} \bigg( \bigcap_{N_1,N_2=1}^M 
    \Big\{ \big\| u_{\leq N_1}^{\phi} - u_{\leq N_2}^{\phi} \big\|_{L_t^\infty \mathscr{H}_x^{-1/2-\epsilon}([-\tau,\tau]\times \T^3)} \leq \min(N_1,N_2)^{-\theta}\Big\} \bigg) \\
    =& \big( \mathbb{P} \otimes \mathbb{Q} \big)  \bigg( \bigcap_{N_1,N_2=1}^M 
    \Big\{ \big\| u_{\leq N_1}^{\textup{cal}} - u_{\leq N_2}^{\textup{cal}} \big\|_{L_t^\infty \mathscr{H}_x^{-1/2-\epsilon}([-\tau,\tau]\times \T^3)} \leq \min(N_1,N_2)^{-\theta}\Big\} \bigg). 
    \end{aligned}
    \end{equation}
    We now choose the parameter $A$, which was previously unspecified, as $A:= \tau^{-\theta}$. Due to Proposition \ref{ansatz:prop-lwp-quantitative}, there exists an $A$-certain event $E_A \in \mathcal{E} \otimes \mathcal{Z}$ such that, on this event,  the estimate  
    \begin{equation*}
    \big\| u_{\leq N_1}^{\textup{cal}} - u_{\leq N_2}^{\textup{cal}} \big\|_{X^{-1/2-\epsilon,b}([-\tau,\tau]\times \T^3)} \leq  \min(N_1,N_2)^{-\theta}
    \end{equation*}
    holds for all $N_1,N_2 \leq M$. Due to the embedding $X^{-1/2-\epsilon,b}\hookrightarrow L_t^\infty \mathscr{H}_x^{-1/2-\epsilon}$,  our choice of $A$, and \eqref{ansatz:eq-lwp-qualitative-p2}, it follows that 
    \begin{equation}\label{ansatz:eq-lwp-qualitative-p3}
    \begin{aligned}
    &\mu_{\leq M} \bigg( \bigcap_{N_1,N_2=1}^M 
    \Big\{ \big\| u_{\leq N_1}^{\phi} - u_{\leq N_2}^{\phi} \big\|_{L_t^\infty \mathscr{H}_x^{-1/2-\epsilon}([-\tau,\tau]\times \T^3)} \leq \min(N_1,N_2)^{-\theta}\Big\} \bigg) \\
    \geq& \, 1 - c_1^{-1} \exp\big(-c_1\tau^{c_1}\big).
    \end{aligned}
    \end{equation}
    In order to prove the desired conclusion, it only remains to pass to the (limiting) Gibbs measure $\mu$. To this end, we recall from \eqref{intro:eq-Gibbs-limit} that the frequency-truncated Gibbs measures $\mu_{\leq M}$ converge weakly to the limiting Gibbs measure $\mu$. For any fixed $N\geq 1$, it follows that 
    \begin{align*}
    &\mu \bigg( \bigcap_{N_1,N_2=1}^N 
    \Big\{ \big\| u_{\leq N_1}^{\phi} - u_{\leq N_2}^{\phi} \big\|_{L_t^\infty \mathscr{H}_x^{-1/2-\epsilon}([-\tau,\tau]\times \T^3)} \leq \min(N_1,N_2)^{-\theta}\Big\} \bigg) \\
    \geq&\, \lim_{M\rightarrow \infty} \mu_{\leq M} \bigg( \bigcap_{N_1,N_2=1}^N 
    \Big\{ \big\| u_{\leq N_1}^{\phi} - u_{\leq N_2}^{\phi} \big\|_{L_t^\infty \mathscr{H}_x^{-1/2-\epsilon}([-\tau,\tau]\times \T^3)} \leq \min(N_1,N_2)^{-\theta}\Big\} \bigg) \\
    \geq& \, \lim_{M\rightarrow \infty} \mu_{\leq M} \bigg( \bigcap_{N_1,N_2=1}^N 
    \Big\{ \big\| u_{\leq N_1}^{\phi} - u_{\leq N_2}^{\phi} \big\|_{L_t^\infty \mathscr{H}_x^{-1/2-\epsilon}([-\tau,\tau]\times \T^3)} \leq \min(N_1,N_2)^{-\theta}\Big\} \bigg) \\ 
    \geq& \, 1 - c_1^{-1} \exp\big(-c_1\tau^{c_1}\big).
    \end{align*}
    Next, we want to let $N$ tend to infinity. To this end, we note that the intersection over all frequency-scales $N_1$ and $N_2$ satisfying $1\leq N_1,N_2 \leq N$ is decreasing in $N$. Since all probability measures are continuous from above, it follows that 
    \begin{align*}
         &\mu \bigg( \bigcap_{N_1,N_2=1}^\infty
    \Big\{ \big\| u_{\leq N_1}^{\phi} - u_{\leq N_2}^{\phi} \big\|_{L_t^\infty \mathscr{H}_x^{-1/2-\epsilon}([-\tau,\tau]\times \T^3)} \leq \min(N_1,N_2)^{-\theta}\Big\} \bigg) \\
     =&\, \lim_{N\rightarrow \infty} \mu \bigg( \bigcap_{N_1,N_2=1}^N 
    \Big\{ \big\| u_{\leq N_1}^{\phi} - u_{\leq N_2}^{\phi} \big\|_{L_t^\infty \mathscr{H}_x^{-1/2-\epsilon}([-\tau,\tau]\times \T^3)} \leq \min(N_1,N_2)^{-\theta}\Big\} \bigg) \\
    \geq& \, 1 - c_1^{-1} \exp\big(-c_1\tau^{c_1}\big). 
    \end{align*}
    By choosing 
    \begin{equation*}
        \operatorname{Local}_\tau := \bigcap_{N_1,N_2=1}^\infty
    \Big\{ \big\| u_{\leq N_1}^{\phi} - u_{\leq N_2}^{\phi} \big\|_{L_t^\infty \mathscr{H}_x^{-1/2-\epsilon}([-\tau,\tau]\times \T^3)} \leq \min(N_1,N_2)^{-\theta}\Big\}, 
    \end{equation*}
    we obtain the desired conclusion. 
    \end{proof}
    
    We now present the proof of the quantitative local well-posedness result (Proposition \ref{ansatz:prop-lwp-quantitative}). 
    
    \begin{proof}[Proof of Proposition \ref{ansatz:prop-lwp-quantitative}] 
    The argument is based on our main estimates (Subsection \ref{section:main-estimates}) and the contraction mapping principle. We first define the ball 
    \begin{align*}
    \mathbb{B}_A := \Big\{& \big(v_{\leq N}, \Xone_{\leq N}, \Xtwo_{\leq N}, Y_{\leq N} \big) \colon 
    \big\| v_{\leq N} \big\|_{X^{1/2-\delta_1,b}([-\tau,\tau])} \leq CA,  \, \\
    & \big\| \Xone_{\leq N} \big\|_{(L_t^\infty \C_x^{1/2-\delta_1} \cap X^{1/2-\delta_1,b})([-\tau,\tau])} \leq CA, \, 
     \big\| \Xtwo_{\leq N} \big\|_{(L_t^\infty \C_x^{1/2-\delta_1} \cap X^{1/2-\delta_1,b})([-\tau,\tau])} \leq CA, \, \\
    &\big\| Y_{\leq N} \big\|_{X^{1/2+\delta_2,b}([-\tau,\tau])} \leq CA \Big\}.
    \end{align*}
    Furthermore, we define the map\footnote{To avoid confusion,  we note that $\Upsilon_{\leq N}$ actually only depends on $v_{\leq N}$ and $Y_{\leq N}$, but not does not depend on $\Xone_{\leq N}$ or $\Xtwo_{\leq N}$. This is because in all terms of the evolution equation for $Y_{\leq N}$ previously involving the para-controlled components $\Xone_{\leq N}$ and $\Xtwo_{\leq N}$, we utilized the double Duhamel-trick to replace them with the para-controlled operators $\XXone$ and $\XXtwo$.}
    \begin{equation*}
    \Upsilon_{\leq N}\big[v_{\leq N}, \Xone_{\leq N}, \Xtwo_{\leq N}, Y_{\leq N}\big]= \Big( \Upsilon_{\leq N}^v, \, \Upsilon_{\leq N}^{X^{(1)}}, \, \Upsilon_{\leq N}^{X^{(2)}}, \Upsilon_{\leq N}^Y \Big)\big[v_{\leq N}, \Xone_{\leq N}, \Xtwo_{\leq N}, Y_{\leq N}\big]
    \end{equation*}
    by
    \begin{alignat*}{3}
    \Upsilon_{\leq N}^v &:= \Upsilon_{\leq N}^{X^{(1)}} + \Upsilon_{\leq N}^{X^{(2)}} + \Upsilon_{\leq N}^Y, 
    \quad&\quad 
    \Upsilon_{\leq N}^{X^{(1)}} &:= \XXone[v_{\leq N},Y_{\leq N}], \\
    \Upsilon_{\leq N}^{X^{(2)}} &:= \XXtwo[v_{\leq N},], 
    \quad& \text{and} \quad \quad
    \Upsilon_{\leq N}^{Y} &:= \slinear[red][\leqM] + \Duh \Big[ \eqref{ansatz:eq-YN-1} + \eqref{ansatz:eq-YN-2} + \hdots + \eqref{ansatz:eq-YN-8} \Big]. 
    \end{alignat*}
   Similar as for $\initial{blue}[][]$ and $\initial{green}[][\leqM]$, $\slinear[red][\leqM]$ denotes the linear evolution with initial data $\initial{red}[][\leqM]$. 
   
   We now prove that $\Upsilon_{\leq N}$ maps $\mathbb{B}_A$ back into itself, i.e., that $\Upsilon_{\leq N}$ is a self-map on $\mathbb{B}_A$. To this end, we first restrict to an $A$-certain event $E_{M,A} \in \mathcal{E}\otimes \mathcal{Z}$ on which the statements in Proposition \ref{ansatz:prop-caloric}, Proposition \ref{ansatz:prop-Xj}, Proposition \ref{ansatz:prop-Y-two}, Proposition \ref{ansatz:prop-Y-one}, and Proposition \ref{ansatz:prop-Y-no} are satisfied. Then, we pick an arbitrary element 
   \begin{equation*}
   \big(v_{\leq N}, \Xone_{\leq N}, \Xtwo_{\leq N}, Y_{\leq N}\big) \in \mathbb{B}_A.
   \end{equation*}
   Due to the triangle inequality and the definition of $\Upsilon_{\leq N}^v$, it suffices to estimate $\Upsilon_{\leq N}^{X^{(1)}}$, $\Upsilon_{\leq N}^{X^{(2)}}$, and $\Upsilon_{\leq N}^Y$. Using Proposition \ref{ansatz:prop-Xj} and the upper bound on $\tau$, $\Upsilon_{\leq N}^{X^{(1)}}$ and $\Upsilon_{\leq N}^{X^{(2)}}$ can be estimated by
   \begin{equation}\label{ansatz:eq-lwp-quantitative-lwp-p1}
   \begin{aligned}
      & \max_{j=1,2} \big\| \Upsilon_{\leq N}^{X^{(j)}} \big\|_{(L_t^\infty \C_x^{1/2-\delta_1} \cap X^{1/2-\delta_1,b})([-\tau,\tau])} \\
       \leq& \,  A \tau^{b_+-b} \Big( 1 + \| v_{\leq N} \|_{X^{-1,b}([-\tau,\tau])}^2 + \| Y_{\leq N} \|_{X^{1/2+\delta_2,b}([-\tau,\tau])} \Big) \\
        \leq& \, 3 C^2 A^3 \tau^{b_+-b} \leq CA/4.  
   \end{aligned}
   \end{equation}
    Thus, it remains to estimate $\Upsilon_{\leq N}^Y$. Since $C\geq 1$ is a sufficiently large, the contribution of $\initial{red}[][\leqM]$ can be estimated using Lemma \ref{prep:lem-xsb}, which yields
    \begin{equation}\label{ansatz:eq-lwp-quantitative-lwp-p2}
    \Big\| \, \slinear[red][\leqM] \Big\|_{X^{1/2+\delta_2,b}([-\tau,\tau])} \leq \, C/4 \big\| \,  \initial{red}[][\leqM]\big\|_{\mathscr{H}_x^{1/2+\delta_2}} \leq CA/4. 
    \end{equation}
    From our main estimates (Proposition \ref{ansatz:prop-Y-two}, Proposition \ref{ansatz:prop-Y-one}, and Proposition \ref{ansatz:prop-Y-no}), it also follows that 
    \begin{equation}\label{ansatz:eq-lwp-quantitative-lwp-p3}
    \begin{aligned}
    &\Big\| \Duh \Big[ \eqref{ansatz:eq-YN-1} + \eqref{ansatz:eq-YN-2} + \hdots + \eqref{ansatz:eq-YN-8} \Big] \Big\|_{X^{1/2+\delta_2,b}([-\tau,\tau])} \\
    \leq&\, CA  \tau^{b_+-b} \Big( 1 + \| v_{\leq N} \|_{X^{-1,b}([-\tau,\tau])}^6 + \| Y_{\leq N} \|_{X^{1/2+\delta_2,b}([-\tau,\tau])}^3 \Big) \\
    \leq&\, 3C^7 A^7 \tau^{b_+-b}  \leq CA/4.  
    \end{aligned}
    \end{equation}
    By combining \eqref{ansatz:eq-lwp-quantitative-lwp-p2} and \eqref{ansatz:eq-lwp-quantitative-lwp-p3}, it follows that 
    \begin{equation}\label{ansatz:eq-lwp-quantitative-lwp-p4}
    \Big\| \Upsilon_{\leq N}^Y \Big\|_{X^{1/2+\delta_2,b}([-\tau,\tau])}\leq CA/2. 
    \end{equation}
    This completes our proof of the self-mapping property of $\Upsilon_{\leq N}$ on $\mathbb{B}_A$.  \\
    
    The contraction property of $\Upsilon_{\leq N}$ on $\mathbb{B}_A$ follows from similar arguments as the self-mapping property. It only requires to replace the main estimates (Proposition \ref{ansatz:prop-Xj}, Proposition \ref{ansatz:prop-Y-two}, Proposition \ref{ansatz:prop-Y-one}, and Proposition \ref{ansatz:prop-Y-no}) by minor generalizations. For example, instead of an estimate of
    \begin{equation*}
        \Big( \Pi^\ast_{\leq N} - \HLL\Big)\Big( \slinear[blue][\leqN], Y_{\leq N}, Y_{\leq N} \Big)
    \end{equation*}
    as in Proposition \ref{ansatz:prop-Y-one}, we require an estimate of 
     \begin{equation*}
        \Big( \Pi^\ast_{\leq N} - \HLL\Big)\Big( \slinear[blue][\leqN], Y_{\leq N}^{(2)}, Y_{\leq N}^{(3)} \Big).
    \end{equation*}
    Since the validity of each minor generalization is an easy consequence of our arguments (in Section \ref{section:diagrams}-Section \ref{section:analytic2}), we omit the standard (but notationally tedious) details. This completes the contraction mapping argument and, as a result, the proof of \ref{ansatz:item-quantitative-2}. \\
    
    The difference estimates from \ref{ansatz:item-quantitative-3} also follows from a minor modification of our previous argument. In all of the frequency-localized estimates in the body of our article (Section \ref{section:analytic}-Section \ref{section:analytic2}), we exhibit a gain in the maximal frequency-scale. After possibly restricting to another $A$-certain event, it follows\footnote{The right-hand side contains a gain in $\tau$ since the linear parts of the evolution do not depend on the frequency-truncation parameters $N_1$ and $N_2$.} that 
    \begin{align*}
    &\big\| u_{\leq N_1} - u_{\leq N_2} \big\| _{X^{-1/2-\epsilon,b}([-\tau,\tau])} + 
    \big\| v_{\leq N_1} - v_{\leq N_2} \big\| _{X^{1/2-\delta_1,b}([-\tau,\tau])} \\
    +& \sum_{j=1,2} \big\| \Xj_{\leq N_1} - \Xj_{\leq N_2} \big\|_{(L_t^\infty \C_x^{1/2-\delta_1} \cap X^{1/2-\delta_1,b})([-\tau,\tau])}  + 
    \big\| Y_{\leq N_1} - Y_{\leq N_2} \big\| _{X^{1/2+\delta_2,b}([-\tau,\tau])}  \\
    \leq& \, C^7 A^7 \tau^{b_+-b} \Big( \min\big(N_1,N_2\big)^{-\theta} + \| v_{\leq N_1} - v_{\leq N_2} \|_{X^{-1,b}} + \| Y_{\leq  N_1} - Y_{\leq N_2} \|_{X^{1/2+\delta_2,b}} \Big). 
    \end{align*}
    Since $0<\tau \leq A^{-\Theta}$, a continuity argument yields the desired conclusion. 
    \end{proof}
    
    Finally, we present the proof of our nonlinear smoothing estimate (Proposition \ref{ansatz:prop-nonlinear-smoothing}). The argument mostly relies on ideas already present in the proofs of Proposition \ref{ansatz:prop-lwp-qualitative} and Proposition \ref{ansatz:prop-lwp-quantitative}. 
    
    \begin{proof}[Proof of Proposition \ref{ansatz:prop-nonlinear-smoothing}]
    
    As in the proof of Proposition \ref{ansatz:prop-lwp-qualitative} and as discussed in Remark \ref{ansatz:rem-ucal}, we first convert the statement regarding the Gibbs measure $\mu_{\leq M}$ into a statement regarding the caloric initial data. To this end,  let $A \geq 1$ remain to be chosen, let $(\mathfrak{Z},\mathcal{Z},\mathbb{Q})$ be a probability space, and let $(\phi^{\cos}_{\leq M},\phi^{\sin}_{\leq M})\colon \mathfrak{Z} \rightarrow \mathscr{H}_x^{-1/2-\epsilon}$ be random functions satisfying 
    \begin{equation*}
        \Law_{\mathbb{Q}}\Big(\big(\phi^{\cos}_{\leq M},\phi^{\sin}_{\leq M}\big) \Big) = \mu_{\leq M}. 
    \end{equation*}
    Furthermore, let 
    $\initial{blue}$, $\initial{green}[][\leqM]$, and $\initial{red}[][\leqM]=\initial{red}[][\leqM](A,\phi^{\cos}_{\leq M})$ be the caloric initial data from Proposition \ref{ansatz:prop-caloric}. Finally, let $u_{\leq N}^{\textup{cal}}$ be the unique global solution of the caloric initial value problem \eqref{ansatz:eq-frequency-truncated-NLW}. Then, the nonlinear smoothing estimate is equivalent to the fact that the event
    \begin{align*}
    \NSE_{M,\tau,A}^{\textup{cal}} := \bigcap_{T\geq 1} \bigcap_{N=1}^M \bigg\{ 
     \Big\| P_{\leq N} \Duh \Big[ \mathbf{1}\big\{ 0\leq t\leq \tau \big\} \Big( \lcol (P_{\leq N}u_{\leq N}^{\textup{cal}})^3 \rcol + \gamma_{\leq N} u_{\leq N}^{\textup{cal}} \Big) \Big] \Big\|_{\NSN([-T,T])} \leq T^\alpha \bigg\}
    \end{align*}
    satisfies
    \begin{equation}\label{ansatz:eq-nonlinear-smoothing-p1}
        \Big( \mathbb{P} \otimes \mathbb{Q} \Big)\big( \NSE_{M,\tau,A}^{\textup{cal}} \big) \geq 1- c_1^{-1} \exp\big(-c_1\tau^{c_1}\big). 
    \end{equation}
    We now choose  
    \begin{equation}\label{ansatz:eq-nonlinear-smoothing-A}
    A:= \tau^{-\theta}.     
    \end{equation}
    Then, Proposition \ref{ansatz:prop-lwp-quantitative} implies that there exists an $A$-certain event $E_A \in \mathcal{E}\otimes \mathcal{Z}$ such that, on this event and on the time-interval $[-\tau,\tau]$, the identities 
    \begin{align}\label{ansatz:eq-nonlinear-smoothing-p2}
    u_{\leq N}^{\textup{cal}} = \slinear[blue] - \scubic[\leqN] - \slinear[green][\leqM] + 3 \squintic[\leqN] + v_{\leq N} \quad 
    \text{and} \quad v_{\leq N}= \Xone_{\leq N} + \Xtwo_{\leq N} + Y_{\leq N}
    \end{align}
    are satisfied, $\big(v_{\leq N},\Xone_{\leq N}, \Xtwo_{\leq N}, Y_{\leq N}\big)$ solves the para-controlled nonlinear wave equations, and the estimates in Proposition \ref{ansatz:prop-lwp-quantitative}.\ref{ansatz:item-quantitative-2} are satisfied. Inserting the decomposition \eqref{ansatz:eq-nonlinear-smoothing-p2} into the nonlinearity and arguing as in the derivation of our Ansatz (Subsection \ref{section:ansatz-explicit}-Subsection \ref{section:ansatz-para}), we have that
    \begin{align}
     &P_{\leq N}  \Big[  \Big( \lcol (P_{\leq N}u_{\leq N}^{\textup{cal}})^3 \rcol + \gamma_{\leq N} u_{\leq N}^{\textup{cal}} \Big) \Big] \notag \\
     =& - \scubicnl[\leqN] + P_{\leq N} \bigg[ 3 \squadratic[\leqN] \scubic[\leqN]  - \Gamma_{\leq N} \slinear[blue][\leqN] \bigg] \label{ansatz:eq-nonlinear-smoothing-p3}\\
     +& \Big\{ \eqref{ansatz:eq-X11-a}+ \eqref{ansatz:eq-X11-b} + \eqref{ansatz:eq-X11-c} + \eqref{ansatz:eq-X11-d} + \eqref{ansatz:eq-X21} \Big\} \label{ansatz:eq-nonlinear-smoothing-p4} \\
     +& \Big\{ \eqref{ansatz:eq-YN-1} + \hdots  + \eqref{ansatz:eq-YN-8} \Big\}. \label{ansatz:eq-nonlinear-smoothing-p5} 
    \end{align}
    Using Lemma \ref{sttime}, Lemma \ref{analytic:lem-cubic}, Proposition \ref{analytic:prop-quintic} , and \eqref{ansatz:eq-nonlinear-smoothing-A},  it follows that 
    \begin{equation}\label{ansatz:eq-nonlinear-smoothing-p6}
    \begin{aligned}
    &\bigg\|  \Duh \bigg[ \mathbf{1}\big\{ 0\leq  t \leq \tau \big\}  \scubicnl[\leqN] \bigg] \bigg\|_{(L_t^\infty \C_x^{-\epsilon} \cap X^{-\epsilon,b})([-T,T])} \\
    +& \bigg\|  P_{\leq N} \Duh \bigg[ \mathbf{1}\big\{ 0\leq   t \leq \tau \big\} \Big( 3 \squadratic[\leqN]  \scubic[\leqN]  - \Gamma_{\leq N} \slinear[blue][\leqN] \Big) \bigg] \bigg\|_{(L_t^\infty \C_x^{-\epsilon} \cap X^{-\epsilon,b})([-T,T])} \\
    \leq&\, A \tau^{b_+-b} T^{\alpha} \leq T^{\alpha}/3. 
    \end{aligned}
    \end{equation}
    We note that the quintic object can be estimated at regularity $1/2-\epsilon$ instead of $-\epsilon$, but this is not used in our nonlinear smoothing estimate\footnote{Of course, the better regularity of the quintic object is heavily used in the proof of our quantitative local well-posedness result (Proposition \ref{ansatz:prop-lwp-quantitative})}. From a (minor variant of) Proposition \ref{ansatz:prop-Xj}, the estimates in Proposition \ref{ansatz:prop-lwp-quantitative}.\ref{ansatz:item-quantitative-2}, and \eqref{ansatz:eq-nonlinear-smoothing-A}, it follows that 
    \begin{equation}\label{ansatz:eq-nonlinear-smoothing-p7}
    \begin{aligned}
    &\Big\| \Duh\Big[ \mathbf{1}\big\{ 0\leq   t\leq \tau \big\} \, \eqref{ansatz:eq-nonlinear-smoothing-p4} \Big] \Big\|_{(L_t^\infty \C_x^{-\epsilon} \cap X^{-\epsilon,b})([-T,T])} \\
    \leq& \,  A \tau^{b_+ - b} T^{\alpha} \Big( 1+ \| v_{\leq N}\|_{X^{-1,b}([0,\tau])}^2 + \| Y_{\leq N} \|_{X^{1/2+\delta_2,b}([0,\tau])}\Big) \\
    \leq& \,  T^{\alpha}/3. 
    \end{aligned}
    \end{equation}
    Finally, from our main estimates (Proposition \ref{ansatz:prop-Y-two}, Proposition \ref{ansatz:prop-Y-one}, and Proposition \ref{ansatz:prop-Y-no}), the estimates in Proposition \ref{ansatz:prop-lwp-quantitative}.\ref{ansatz:item-quantitative-2}, and \eqref{ansatz:eq-nonlinear-smoothing-A}, it follows that 
      \begin{equation}\label{ansatz:eq-nonlinear-smoothing-p8}
    \begin{aligned}
    &\Big\| \Duh\Big[ \mathbf{1}\big\{ 0\leq   t\leq \tau \big\} \, \eqref{ansatz:eq-nonlinear-smoothing-p5} \Big] \Big\|_{X^{1/2+\delta_2,b}([-T,T])} \\
    \leq& \,  A \tau^{b_+ - b} T^{\alpha} \Big( 1+ \| v_{\leq N}\|_{X^{-1,b}([0,\tau])}^6 + \| Y_{\leq N} \|_{X^{1/2+\delta_2,b}([0,\tau])}^3\Big) \\
    \leq& \,  T^{\alpha}/3.
   \end{aligned}
    \end{equation}
    By combining \eqref{ansatz:eq-nonlinear-smoothing-p6}, \eqref{ansatz:eq-nonlinear-smoothing-p7}, and \eqref{ansatz:eq-nonlinear-smoothing-p8}, we obtain that the event $ \NSE_{M,\tau,A}^{\textup{cal}}$ is $A$-certain. Due to our choice of $A$ from  \eqref{ansatz:eq-nonlinear-smoothing-A}, this implies \eqref{ansatz:eq-nonlinear-smoothing-p1}. 
    \end{proof}

\section{Global well-posedness and invariance}\label{section:global}

In this section, we prove our main theorem (Theorem \ref{intro:thm-rigorous}). 
In Subsection
\ref{section:global-bounds}, we obtain global bounds for the para-controlled components. The main ingredients are our main estimates (as in Subsection \ref{section:main-estimates}) and a global version of our nonlinear smoothing estimate (Proposition \ref{ansatz:prop-nonlinear-smoothing}). The latter is proved using a variant of Bourgain's globalization argument \cite{B94}. 
In Subsection \ref{section:global-stability}, we prove a stability estimate. In comparison to \cite{B20II,OOT20}, the main novelty is that the stability estimate is stated in terms of the individual components of the caloric initial data. In Subsection \ref{section:global-proof}, we combine the global bounds and stability estimates to prove our main theorem (Theorem \ref{intro:thm-rigorous}). 

\subsection{Global bounds}\label{section:global-bounds}

The stability theory presented in the next subsection is based on our quantitative local well-posedness result (Proposition \ref{ansatz:prop-lwp-quantitative}). In particular, it involves the caloric initial data (from Proposition \ref{ansatz:prop-caloric}) and the para-controlled components $\big(\Xone_{\leq M},\Xtwo_{\leq M},Y_{\leq M}\big)$. In order to insert the global bounds into our stability theory, they need to control the para-controlled components. This is the subject of the next proposition. 

\begin{proposition}[Para-controlled global bounds]\label{global:prop-para-controlled}
Let $M\geq 1$, $A\geq 1$, and let $\initial{blue}$, $\initial{green}[][\leqM]$, and $\initial{red}[][\leqM](A)$ be the caloric initial data from Proposition \ref{ansatz:prop-caloric}. Then, there exists an event $\operatorname{Glb}_{M,A}\in \mathcal{E}\otimes \mathcal{Z}$ such that the following properties are satisfied:
\begin{enumerate}[label=(\roman*)]
    \item\label{global:item-bounds-1} (High probability) It holds that 
    \begin{equation*}
     \big( \mathbb{P} \otimes \mathbb{Q} \big)\Big( \operatorname{Glb}_{M,A} \Big) \geq 1 - c_3^{-1} \exp\big(-c_3 A^{c_3} \big). 
    \end{equation*}
    \item\label{global:item-bounds-2} (Para-controlled global bounds) On the event $\operatorname{Glb}_{M,A}$, the solution $u_{\leq M}$ of the caloric initial value problem \eqref{ansatz:eq-frequency-truncated-NLW-orig} and the solutions $(v_{\leq M}, \Xone_{\leq M},\Xtwo_{\leq M},Y_{\leq M})$ of the para-controlled nonlinear wave equations \eqref{ansatz:eq-X11-a}-\eqref{ansatz:eq-X21-initial} and \eqref{ansatz:eq-YN-1}-\eqref{ansatz:eq-YN-8} exist globally. Furthermore, they satisfy
    \begin{align*}
    &\max \bigg(\big\|  u_{\leq M} \big\| _{X^{-1/2-\epsilon,b}([0,T])}, \, \big\|  v_{\leq M} \big\| _{X^{1/2-\delta_1,b}([0,T])}, 
    \big\|  \Xone_{\leq M} \big\|_{(L_t^\infty \C_x^{1/2-\delta_1} \cap X^{1/2-\delta_1,b})([0,T])}, \\ 
    &\hspace{7ex}
    \big\|  \Xtwo_{\leq M} \big\|_{(L_t^\infty \C_x^{1/2-\delta_1} \cap X^{1/2-\delta_1,b})([0,T])}, \, 
    \big\|  Y_{\leq M} \big\| _{X^{1/2+\delta_2,b}([0,T])} \bigg) \\
    \leq& \, C_3 \exp\big( C_3 (AT)^{C_3} \big)
    \end{align*}
    for all $T\geq 1$. 
\end{enumerate}
\end{proposition}

\begin{remark}\label{global:rem-caloric} Proposition \ref{global:prop-para-controlled}, whose proof utilizes the invariance of the frequency-truncated Gibbs measure, establishes the \emph{global} propagation of randomness. More precisely, it proves that the decomposition 
\begin{equation}\label{global:eq-global-propagation-randomness}
u_{\leq M}=\slinear[blue] - \scubic[\leqM] - \slinear[green][\leqM] + 3 \squintic[\leqM] + \Xone_{\leq M} + \Xtwo_{\leq M} + Y_{\leq M}
\end{equation}
holds for all time and that, again for all time, all terms in the decomposition have the same regularity as in the local theory. We emphasize that the global random structure in \eqref{global:eq-global-propagation-randomness} depends on the caloric representation of the Gibbs measure, since the three components   $\initial{blue}$, $\initial{green}[][\leqM]$, and $\initial{red}[][\leqM]$ appear individually\footnote{The red caloric initial data $\initial{red}[][\leqM]$ is hidden in the initial condition $Y_{\leq M}[0]=\initial{red}[][\leqM]$.}.

In the earlier works \cite{B20II} and \cite{OOT21},  the random initial data was represented using the drift measure \cite{BG20} or the variational approach \cite{BG18} of Barashkov and Gubinelli, respectively. Similar as in this article, the random structure in the local theory of \cite{B20II,OOT21} relies on a representation of the Gibbs measure which consists of a Gaussian free field and a random shift. In contrast to this article, however, the random structure in the global theory is stated without explicit reference to the representation of the Gibbs measure (see e.g. \cite[Proposition 3.1]{B20II} or \cite[Proposition 6.5]{OOT21}). This requires a significant computational effort (see e.g. \cite[Section 9]{B20II} or \cite[Proposition 6.5 and Remark 6.6]{OOT21}) which, due to the larger number of stochastic objects, makes this approach practically infeasible in our setting.
\end{remark}

The proof of Proposition \ref{global:prop-para-controlled} occupies the remainder of this subsection. As a stepping stone, we first prove a global version of our previous nonlinear smoothing estimate (Proposition \ref{ansatz:prop-nonlinear-smoothing}). 

\begin{lemma}[Global nonlinear smoothing estimate]\label{global:lem-nonlinear-smoothing}
Let $M\geq 1$ and let $A\geq 1$. Then, there exists a Borel measurable event $\NSEg_{M,A} \subseteq \mathscr{H}_x^{-1/2-\epsilon}$ such that the following two properties hold:
\begin{enumerate}[label=(\roman*)]
    \item\label{global:item-nonlinear-smoothing-1} (High probability) It holds that 
    \begin{equation*}
        \mu_{\leq M}\big( \NSEg_{M,A}\big) \geq 1- c_2^{-1} \exp\big( -c_2 A^{c_2} \big).  
    \end{equation*}
    \item\label{global:item-nonlinear-smoothing-2} (Nonlinear smoothing) For all $(\phi^{\cos},\phi^{\sin}) \in \NSEg_{M,A}$ and all $T\geq 1$, the solution $u_{\leq M}$ of \eqref{ansatz:eq-frequency-truncated-NLW-raw} satisfies 
    \begin{equation*}
        \Big\| P_{\leq M} \Duh \Big[  \Big( \lcol (P_{\leq M}u_{\leq M})^3 \rcol + \gamma_{\leq M} u_{\leq M} \Big) \Big] \Big\|_{\NSN([0,T])} \leq A T^\alpha. 
    \end{equation*}
\end{enumerate}
\end{lemma}

\begin{remark}
The basic idea behind the proof of Lemma \ref{global:lem-nonlinear-smoothing}, which is to prove global bounds using local bounds and invariance, goes back to Bourgain \cite{B94}. 
\end{remark}

\begin{remark}[Comparison with Proposition \ref{ansatz:prop-nonlinear-smoothing} and Proposition \ref{global:prop-para-controlled}]
In contrast our earlier nonlinear smoothing estimate (Proposition \ref{ansatz:prop-nonlinear-smoothing}), the estimate in Lemma \ref{global:lem-nonlinear-smoothing} does not require the restriction to times $0\leq t \leq \tau$ in the integrand. However, while the (local in time) nonlinear smoothing estimate from Proposition \ref{ansatz:prop-nonlinear-smoothing} holds for all frequency-truncation parameters $N\leq M$, Lemma \ref{global:lem-nonlinear-smoothing} only holds for $N=M$. \\

In contrast to Proposition \ref{global:prop-para-controlled}, Lemma \ref{global:lem-nonlinear-smoothing} refers to neither the caloric initial data nor the para-controlled components. As a result, its proof is more tractable through invariance-based methods. The drawback is that, since it does not refer to the para-controlled components, Lemma \ref{global:lem-nonlinear-smoothing} yields less precise information on the solution. In particular, it only controls the nonlinear part of the evolution at regularity $0-$.  
\end{remark}

\begin{proof}[Proof of Lemma \ref{global:lem-nonlinear-smoothing}:]
Since the unspecified parameter $\alpha$ occurs in both Proposition \ref{ansatz:prop-nonlinear-smoothing} and Lemma \ref{global:lem-nonlinear-smoothing}, we denote the parameter from Proposition \ref{ansatz:prop-nonlinear-smoothing} by $\alpha^\prime$. Throughout this proof, we denote the time-variable of the Duhamel integral by $t$ and the time-variable of the integrand by $t^\prime$. \\

We first define the event 
\begin{equation*}
\NSEg_{M,A} := \bigcap_{T \geq 1} \NSEg_{M,A,T}\,,
\end{equation*}
where 
\begin{equation}\label{global:eq-nonlinear-smoothing-p0}
\NSEg_{M,A,T} := \bigg\{ \big(\phi^{\cos}, \phi^{\sin}\big) \colon 
 \Big\| P_{\leq M} \Duh \Big[  \Big( \lcol (P_{\leq M}u_{\leq M})^3 \rcol + \gamma_{\leq M} u_{\leq M} \Big) \Big] \Big\|_{\NSN([0,T])} \leq A T^\alpha \bigg\}. 
\end{equation}
From our definition in \eqref{global:eq-nonlinear-smoothing-p0}, we directly obtain the nonlinear smoothing estimate in \ref{global:item-nonlinear-smoothing-2}. As a result, it only remains to prove the that $\NSEg_{M,A,T}$ has high probability, i.e., the estimate in \ref{global:item-nonlinear-smoothing-1}. 
To this end, we fix any $T\geq 1$, let $J\in \mathbb{N}$ remain to be chosen, and set $\tau:=T/J$. Then, we decompose
\begin{equation}\label{global:eq-nonlinear-smoothing-p1}
\begin{aligned}
 &P_{\leq M} \Duh \Big[  \Big( \lcol (P_{\leq M}u_{\leq M})^3 \rcol + \gamma_{\leq M} u_{\leq M} \Big)(t^\prime) \Big](t) \\
 =& \sum_{j=0}^{J-1} P_{\leq M} \Duh \Big[ \mathbf{1}\big\{ j\tau \leq t^\prime \leq (j+1) \tau \big\} \Big( \lcol (P_{\leq M}u_{\leq M})^3 \rcol + \gamma_{\leq M} u_{\leq M} \Big)(t^\prime) \Big](t). 
\end{aligned}
\end{equation}
From a change of variables, it follows that 
\begin{equation}\label{global:eq-nonlinear-smoothing-p2}
\begin{aligned}
     &P_{\leq M} \Duh \Big[ \mathbf{1}\big\{ j\tau \leq t^\prime \leq (j+1) \tau \big\} \Big( \lcol (P_{\leq M}u_{\leq M})^3 \rcol + \gamma_{\leq M} u_{\leq M} \Big)(t^\prime) \Big](t) \\
     =& P_{\leq M} \Duh \Big[ \mathbf{1}\big\{ 0 \leq t^\prime \leq \tau \big\} \Big( \lcol (P_{\leq M}u_{\leq M})^3 \rcol + \gamma_{\leq M} u_{\leq M} \Big)(t^\prime+j\tau) \Big](t-j\tau). 
\end{aligned}
\end{equation}
Using \eqref{global:eq-nonlinear-smoothing-p1}, \eqref{global:eq-nonlinear-smoothing-p2}, and the time-translation invariance of $L_t^q L_x^r$ and $X^{s,b}$-norms, we obtain that 
\begin{equation}\label{global:eq-nonlinear-smoothing-p3}
\begin{aligned}
&\bigg\| P_{\leq M} \Duh \Big[  \Big( \lcol (P_{\leq M}u_{\leq M})^3 \rcol + \gamma_{\leq M} u_{\leq M} \Big)(t^\prime) \Big](t) \bigg\|_{\NSN([0,T])} \\
\leq & \sum_{j=0}^{J-1} \, \bigg\| P_{\leq M} \Duh \Big[ \mathbf{1}\big\{ 0 \leq t^\prime \leq \tau \big\} \Big( \lcol (P_{\leq M}u_{\leq M})^3 \rcol + \gamma_{\leq M} u_{\leq M} \Big)(t^\prime+j\tau) \Big](t-j\tau) \bigg\|_{\NSN([0,T])} \\
\leq & \sum_{j=0}^{J-1} \, \bigg\| P_{\leq M} \Duh \Big[ \mathbf{1}\big\{ 0 \leq t^\prime \leq \tau \big\} \Big( \lcol (P_{\leq M}u_{\leq M})^3 \rcol + \gamma_{\leq M} u_{\leq M} \Big)(t^\prime+j \tau) \Big](t) \bigg\|_{\NSN([-T,T])}. 
\end{aligned}
\end{equation}
Using  \eqref{global:eq-nonlinear-smoothing-p3}, a union bound, and the invariance of $\mu_{\leq M}$ under $u_{\leq M}$, we obtain that 
\begin{align}
&\mu_{\leq M}\big( \mathscr{H}_x^{-1/2-\epsilon}\backslash \NSEg_{M,A,T} \big) \notag \\
\leq& \sum_{j=0}^{J-1} \mu_{\leq M} \bigg(  \, 
\bigg\| P_{\leq M} \Duh \Big[ \mathbf{1}\big\{ 0 \leq t^\prime \leq \tau \big\} \Big( \lcol (P_{\leq M}u_{\leq M})^3 \rcol + \gamma_{\leq M} u_{\leq M} \Big)(t^\prime +j \tau) \Big](t) \bigg\|_{\NSN([-T,T])} > A J^{-1}T^\alpha \bigg) \notag \\
\leq& \, J \cdot  \mu_{\leq M} \bigg(  \, 
\bigg\| P_{\leq M} \Duh \Big[ \mathbf{1}\big\{ 0 \leq t^\prime \leq \tau \big\} \Big( \lcol (P_{\leq M}u_{\leq M})^3 \rcol + \gamma_{\leq M} u_{\leq M} \Big)(t^\prime) \Big](t) \bigg\|_{\NSN([-T,T])} > A J^{-1}T^\alpha \bigg). \label{global:eq-nonlinear-smoothing-p4}
\end{align}
As mentioned at the beginning of the proof, we denote the parameter from Proposition \ref{ansatz:prop-nonlinear-smoothing} by $\alpha^\prime$. We now choose $\alpha=\alpha^\prime+2$ and $J= AT^2$.   Using Proposition \ref{ansatz:prop-nonlinear-smoothing}, it follows that 
\begin{align*}
&\mu_{\leq M} \bigg(  \, 
\bigg\| P_{\leq M} \Duh \Big[ \mathbf{1}\big\{ 0 \leq t^\prime \leq \tau \big\} \Big( \lcol (P_{\leq M}u_{\leq M})^3 \rcol + \gamma_{\leq M} u_{\leq M} \Big)(t^\prime) \Big](t) \bigg\|_{\NSN([-T,T])} > A J^{-1}T^\alpha \bigg) \\
=&  \mu_{\leq M} \bigg(  \, 
\bigg\| P_{\leq M} \Duh \Big[ \mathbf{1}\big\{ 0 \leq t^\prime \leq \tau \big\} \Big( \lcol (P_{\leq M}u_{\leq M})^3 \rcol + \gamma_{\leq M} u_{\leq M} \Big)(t^\prime) \Big](t) \bigg\|_{\NSN([-T,T])} > T^{\alpha^\prime} \bigg) \\
\leq& c_1^{-1} \exp\big( c_1^{-1} \tau^{-c_1} \big) = c_1^{-1} \exp\big(-c_1 T^{c_1} A^{c_1} \big).
\end{align*}
Inserting this back into \eqref{global:eq-nonlinear-smoothing-p4}, we obtain that 
\begin{equation}\label{global:eq-nonlinear-smoothing-p5}
   \mu_{\leq M}\big( \mathscr{H}_x^{-1/2-\epsilon}\backslash \NSEg_{M,A,T} \big) \leq c_1^{-1} A T^2 \exp\big( -c_1 T^{c_1} A^{c_1}\big).  
\end{equation}

By choosing $c_2=c_2(c_1)$ sufficiently small, \eqref{global:eq-nonlinear-smoothing-p5} implies that 
\begin{equation*}
   \mu_{\leq M}\big( \mathscr{H}_x^{-1/2-\epsilon}\backslash \NSEg_{M,A} \big) \leq 
    \sum_{T\geq 1} \mu_{\leq M}\big( \mathscr{H}_x^{-1/2-\epsilon}\backslash \NSEg_{M,A,T} \big) \leq c_2^{-1} \exp\big( -c_2 A^{c_2}\big). 
\end{equation*}
This completes the proof. 
\end{proof}

The next ingredient in our proof of the para-controlled global bounds is the following nonlinear estimate, which is a modification of Proposition \ref{ansatz:prop-Y-one} and Proposition \ref{ansatz:prop-Y-no}. 

\begin{lemma}[Terms with one or zero linear stochastic objects revisited]\label{global:lem-one-zero-revisited}
For all $A\geq 1$, there exists an $A$-certain event $E_A \in \mathcal{E}$ such that the following estimates hold: Let $N\geq 1$, let $T\geq 1$, let $0 \in \Jc \subseteq [0,T]$ be a closed interval, and let $v_{\leq N},Y_{\leq N}\colon [0,T] \times \T^3 \rightarrow \R$ satisfy
\begin{equation}\label{global:eq-one-zero-decomp}
v_{\leq N}= \XXone[v_{\leq N},Y_{\leq N}] + \XXtwo[v_{\leq N}] + Y_{\leq N}. 
\end{equation}
Let $0 \in \mathcal{J} \subseteq [0,T]$ be a closed interval. If $\zeta^{(1)}=\slinear[blue]$ and $\zeta^{(2)},\zeta^{(3)}\in \Symb_0^b$, it holds that
\begin{equation}\label{global:eq-one-zero-a}
\begin{aligned}
&\Big\| \Big( \Pi^\ast_{\leq N} - \HLL \Big) \Big( \slinear[blue][\leqN], \zeta^{(2)}_{\leq N}, \zeta^{(3)}_{\leq N} \Big) \Big\|_{X^{-1/2+\delta_2,b_+-1}(\mathcal{J})} \\
\leq& \, A T^\alpha \Big( 1+ \| v_{\leq N} \|_{\NSN([0,T])}^3 \Big) \Big( 1+ \| Y_{\leq N} \|_{X^{1/2+\delta_2,b}(\mathcal{J})} \Big). 
\end{aligned}
\end{equation}
If $\zeta^{(1)},\zeta^{(2)},\zeta^{(3)}\in \Symb_0^b$, it holds that 
\begin{equation}\label{global:eq-one-zero-b}
\begin{aligned}
&\Big\| \ \Pi^\ast_{\leq N} \Big( \zeta^{(1)}_{\leq N}, \zeta^{(2)}_{\leq N}, \zeta^{(3)}_{\leq N} \Big) \Big\|_{X^{-1/2+\delta_2,b_+-1}(\mathcal{J})} \\
\leq& \, A T^\alpha \Big( 1+ \| v_{\leq N} \|_{\NSN([0,T])}^4 \Big) \Big( 1+ \| Y_{\leq N} \|_{X^{1/2+\delta_2,b}(\mathcal{J})} \Big). 
\end{aligned}
\end{equation}
\end{lemma}

\begin{remark}
The most important aspect of \eqref{global:eq-one-zero-a} and \eqref{global:eq-one-zero-b} is that the right-hand side only grows linearly in $\| Y_{\leq N}\|_{X^{1/2+\delta_2,b}}$.
\end{remark}

Since the argument is extremely similar to the proof of Proposition \ref{ansatz:prop-Y-one} and Proposition \ref{ansatz:prop-Y-no}, it is postponed until Section \ref{section:proof-main-estimates}. 
Equipped with Lemma \ref{global:lem-nonlinear-smoothing} and Lemma \ref{global:lem-one-zero-revisited}, we are now ready to prove Proposition \ref{global:prop-para-controlled}. 

\begin{proof}[Proof of Proposition \ref{global:prop-para-controlled}]
In this proof, we will only restrict to events from Lemma \ref{global:lem-nonlinear-smoothing} and $A$-certain events, which both satisfy the probability estimate in \ref{global:item-bounds-1}. As a result, we focus our attention on the verification of the estimates in \ref{global:item-bounds-2}.\\

Using Lemma \ref{global:lem-nonlinear-smoothing}, we obtain (with sufficiently high probability) that 
\begin{equation*}
\Big\| u_{\leq M} - \Big( \slinear[blue][] - \slinear[green][\leqM] + \slinear[red][\leqM] \Big) \Big\|_{\NSN([0,T])} \leq A T^{\alpha}
\end{equation*}
for all $T\geq 1$. After possibly restricting to an $A$-certain event, it follows from our estimates for explicit stochastic objects  (Proposition \ref{ansatz:prop-caloric}, Lemma \ref{analytic:lem-linear}, Lemma \ref{analytic:lem-cubic}, and Proposition \ref{analytic:prop-quintic}) that
\begin{equation}\label{global:eq-bounds-p1}
\begin{aligned}
&\Big\| v_{\leq M} \Big \|_{\NSN([0,T])}  \\
\leq& \,
\Big\| u_{\leq M} - \Big( \slinear[blue][] - \slinear[green][\leqM] + \slinear[red][\leqM] \Big) \Big\|_{\NSN([0,T])}  
+ \Big\| \, \slinear[red][\leqM] \Big\|_{\NSN([0,T])}
+  \Big\|  \, \scubic[\leqM] \Big\|_{\NSN([0,T])} 
+  3 \,  \Big\|  \, \squintic[\leqM] \Big\|_{\NSN([0,T])} \\
\leq& \, 2 A T^\alpha
\end{aligned}
\end{equation}
for all $T\geq 1$. 

We now fix $T\geq 1$, let $J\in \mathbb{N}$, and set $\tau := T/J$. For any $0\leq j \leq J$, we define 
\begin{equation*}
A_j := A + \big\| Y_{\leq M} \big\|_{X^{1/2+\delta_2,b}([0,\,j\tau])} 
+ \big\| \big( \partial_t^2 +1 - \Delta \big) Y_{\leq M} \big\|_{X^{-1/2+\delta_2,b-1}([0,\,j\tau])}.
\end{equation*}
For all $0\leq j\leq J-1$, we obtain from the gluing lemma (Lemma \ref{prep:lem-gluing}) and time-localization (Lemma \ref{sttime}) that 
\begin{equation}\label{global:eq-bounds-p2}
\begin{aligned}
A_{j+1}
\lesssim&\, A+ \big\| \, \initial{red}[][\leqM] \big\|_{\mathscr{H}_x^{1/2+\delta_2}}
+ \big\| \big( \partial_t^2 +1 - \Delta \big) Y_{\leq M} \big\|_{X^{-1/2+\delta_2,b-1}([0,(j+1)\tau])} \\
\lesssim&\, A 
+ \big\| \big( \partial_t^2 +1 - \Delta \big) Y_{\leq M} \big\|_{X^{-1/2+\delta_2,b-1}([0,j\tau])} + \big\| \big( \partial_t^2 +1 - \Delta \big) Y_{\leq M} \big\|_{X^{-1/2+\delta_2,b-1}([j\tau,(j+1)\tau])} \\
\lesssim&\, A_j + \tau^{b_+-b} 
\big\| \big( \partial_t^2 +1 - \Delta \big) Y_{\leq M} \big\|_{X^{-1/2+\delta_2,b_+-1}([0,(j+1)\tau])}. 
\end{aligned}
\end{equation}
Arguing as in the derivation of \eqref{ansatz:eq-YN-1}-\eqref{ansatz:eq-YN-8} but without (immediately) inserting the para-controlled Ansatz in terms with only one or zero linear stochastic objects, the evolution equation for $Y_{\leq M}$ can be written as 
\begin{align}
&\big( \partial_t^2 + 1 - \Delta \big) Y_{\leq M} \\
=&   - P_{\leq M} \bigg[ 9 \squadratic[\leqM] \squintic[\leqM] - \Gamma_{\leq M}  \scubic[\leqM] - 18 \mathfrak{C}^{(1,5)}_{\leq M} P_{\leq M} \slinear[blue][\leqM] \label{global:eq-YM-1} \\
    &\hspace{8ex}- 9 \Big( 2 \HLLM + \HHLM \Big) \Big( \, \slinear[blue][\leqM], \slinear[blue][\leqM], \squintic[\leqM] \Big) \bigg] \notag \allowdisplaybreaks[3] \\
    +&3 P_{\leq M} \bigg[ \squadratic[\leqM] \slinear[green][\leqM] - \Big( 2 \HLLM + \HHLM \Big) \Big( \slinear[blue][\leqM] , \slinear[blue][\leqM] , \slinear[green][\leqM] \Big) \bigg] \\
    -&3 P_{\leq M} \bigg[ \squadratic[\leqM] \XXoneM - 
    \Big( 2 \HLLM + \HHLM \Big) \Big( \slinear[blue][\leqM] , \slinear[blue][\leqM], \XXoneM \Big) \bigg] \allowdisplaybreaks[3] \label{global:eq-YM-2} \\
    -& P_{\leq M} \bigg[ 3 \squadratic[\leqM] \XXtwoM - 
    \Big( 6 \HLLM + 3 \HHLM \Big) \Big( \slinear[blue][\leqM] , \slinear[blue][\leqM], \XXtwoM \Big)
    + \Gamma_{\leq M} \Big( 3 \squintic[\leqM] - \slinear[green][\leqM]+ v_{\leq M} \Big) \bigg] \allowdisplaybreaks[3] \label{global:eq-YM-3} \\
    -& 3 P_{\leq M} \bigg[ \squadratic[\leqM] Y_{\leq M} - 
    \Big( 2 \HLLM + \HHLM + \RESM \Big) \Big( \slinear[blue][\leqM] , \slinear[blue][\leqM], Y_{\leq M} \Big) \bigg] \allowdisplaybreaks[3] \label{global:eq-YM-4} \\
    +& \sum_{\zeta^{(2)},\zeta^{(3)} \in \Symb_0^b } \Athree\big(\, \slinear, \zeta^{(2)}, \zeta^{(3)} \big)  P_{\leq M} \Big( \Pi^\ast_{\leq M}- \HLLM \Big) \Big( \slinear[blue][\leqM], \zeta^{(2)}_{\leq M}, \zeta^{(3)}_{\leq M} \Big) \allowdisplaybreaks[3] \label{global:eq-YM-5} \\
    +& \sum_{\zeta^{(1)},\zeta^{(2)},\zeta^{(3)} \in \Symb_0^b} \Athree \big( \zeta^{(1)}, \zeta^{(2)}, \zeta^{(3)} \big) P_{\leq M} \Pi^\ast_{\leq M} \Big( P_{\leq M} \zeta^{(1)}, P_{\leq M} \zeta^{(2)}, P_{\leq M} \zeta^{(3)}  \Big) \allowdisplaybreaks[3] \label{global:eq-YM-6} \\
    +& \mathfrak{C}_{\leq M} \sum_{\zeta \in \Symb_0^b} \Aone(\zeta) P_{\leq M} \zeta_{\leq M}
    + (\gamma_{\leq M}-\Gamma_{\leq M}) \sum_{\zeta \in \Symb_0^b} \Atone(\zeta) P_{\leq M} \zeta_{\leq M}. \label{global:eq-YM-7}
\end{align}
After possible restricting to a further $A$-certain event, we now assume that our main estimates (Proposition \ref{ansatz:prop-Y-two} and Proposition \ref{ansatz:prop-Y-no}) and the modified version (from Lemma \ref{global:lem-one-zero-revisited}) are satisfied. We can then control \eqref{global:eq-YM-1}-\eqref{global:eq-YM-4} using the estimates in Proposition \ref{ansatz:prop-Y-two}, \eqref{global:eq-YM-5}-\eqref{global:eq-YM-6} using the estimates in Lemma \ref{global:lem-one-zero-revisited}, and \eqref{global:eq-YM-7} using the estimates in Proposition \ref{ansatz:prop-Y-no}. As a result, we obtain that 
\begin{equation}\label{global:eq-bounds-p3}
\begin{aligned}
&\big\| \big( \partial_t^2 +1 - \Delta \big) Y_{\leq M} \big\|_{X^{-1/2+\delta_2,b-1}([0,(j+1)\tau])} \\
\leq& \,  
A T^\alpha \Big( 1 + \big\| v_{\leq M} \big\|_{\NSN([0,T])}^4 \Big) \Big( 1+ \big\| Y_{\leq M}\big\|_{X^{1/2+\delta_2,b}([0,(j+1)\tau])} \Big). 
\end{aligned}
\end{equation}
By combining \eqref{global:eq-bounds-p1}, \eqref{global:eq-bounds-p2}, and \eqref{global:eq-bounds-p3}, we obtain 
\begin{equation}\label{global:eq-bounds-p4}
A_{j+1} \leq C A_j + C \tau^{b_+-b} A^5 T^{5\alpha} \big(1+A_j\big), 
\end{equation}
where $C=C(\delta_1,\delta_2,\eta,\nu,b,b_+)\gg 1$. 
After choosing 
\begin{equation*}
\tau := \Big( 4 C A^5 T^{5 \alpha} \Big)^{-1/(b_+-b)},  
\end{equation*}
iterating \eqref{global:eq-bounds-p4} and recalling the definition of $A_j$ yields
\begin{equation}\label{global:eq-bounds-p5}
 \big\| Y_{\leq M} \big\|_{X^{1/2+\delta_2,b}([0,\,T])} \leq C_3 \exp\bigg( \frac{1}{10} C_3 (AT)^{C_3} \bigg). 
\end{equation}
This proves (a slightly better) version of the estimate for $Y_{\leq M}$ in \ref{global:item-bounds-2}. The desired estimate for $\Xone_{\leq M}$ and $\Xtwo_{\leq M}$ now follows directly from Proposition \ref{ansatz:prop-Xj}, \eqref{global:eq-bounds-p1}, and \eqref{global:eq-bounds-p5}. Finally, the desired estimates for $u_{\leq M}$ and $v_{\leq M}$ follow directly from the triangle inequality and the estimates for $\Xone_{\leq M}$, $\Xtwo_{\leq M}$, and $Y_{\leq M}$.
\end{proof}

\begin{remark} In the proof of Proposition \ref{global:prop-para-controlled}, we first obtained bounds on the nonlinear remainder $v_{\leq M}$ using invariance (as in the proof of Lemma \ref{global:lem-nonlinear-smoothing}) and then upgraded the bounds using our main estimates (as in Subsection \ref{section:main-estimates}). A similar approach was previously used in \cite[Proposition 6.9]{OOT21}, see e.g. the invariance-based bound \cite[(6.95)]{OOT21} and the following analysis of \cite[(6.96)]{OOT21}. The main difference between our proof of Proposition \ref{global:prop-para-controlled} and the proof of \cite[Proposition 6.9]{OOT21} is that our nonlinear smoothing estimate utilizes multilinear dispersive effects, whereas \cite[(6.94)]{OOT21} only controls the Wick-ordered square pointwise in time. 
\end{remark}

\subsection{Stability theory}\label{section:global-stability}

In this subsection, we prove a stability estimate (Proposition \ref{global:prop-stability}). Together with the global bounds from Subsection \ref{section:global-bounds}, it will be used to prove the main theorem.

\begin{proposition}[Global stability estimate]\label{global:prop-stability}
Let $M\geq 1$, let $A\geq 1$, and let $\initial{blue}$, $\initial{green}[][\leqM]$, and $\initial{red}[][\leqM](A)$ be the caloric initial data from Proposition \ref{ansatz:prop-caloric}. Then, there exists an event $\Stbl_{M,A} \in \mathcal{E}\otimes \mathcal{Z}$ which satisfies
\begin{equation*}
\big( \mathbb{P} \otimes \mathbb{Q} \big)\Big( \Stbl_{M,A} \Big) 
\geq 1- c_4^{-1} \exp\big( -c_4 A^{c_4} \big)
\end{equation*}
and such that, on this event, the following implication hold for all $T\geq 1$:

Let  $\big(v_{\leq M},\Xone_{\leq M},\Xtwo_{\leq M},Y_{\leq M}\big)$ be the para-controlled components with frequency-truncation parameter $M$, let $B\geq 1$, and assume that
\begin{align*}
    &\max \bigg(\big\|  u_{\leq M} \big\| _{X^{-1/2-\epsilon,b}([0,T])}, \, \big\|  v_{\leq M} \big\| _{X^{1/2-\delta_1,b}([0,T])}, 
    \big\|  \Xone_{\leq M} \big\|_{(L_t^\infty \C_x^{1/2-\delta_1} \cap X^{1/2-\delta_1,b})([0,T])}, \\ 
    &\hspace{7ex}
    \big\|  \Xtwo_{\leq M} \big\|_{(L_t^\infty \C_x^{1/2-\delta_1} \cap X^{1/2-\delta_1,b})([0,T])}, \, 
    \big\|  Y_{\leq M} \big\| _{X^{1/2+\delta_2,b}([0,T])} \bigg) 
    \leq \, B. 
    \end{align*}
    Furthermore, let $1\leq N\leq M$ be a frequency-scale satisfying 
    \begin{equation*}
    C_4 \exp\Big( C_4 \big( AB T \big)^{C_4} \Big) N^{-\theta} \leq 1. 
    \end{equation*}
    Then, the solution $u_{\leq N}$ and para-controlled components $\big(v_{\leq N}, \Xone_{\leq N}, \Xtwo_{\leq N}, Y_{\leq N}\big)$ with frequency-truncation parameter $N$ satisfy
    \begin{align*}
    &\max \bigg(\big\|  u_{\leq M} -  u_{\leq N} \big\| _{X^{-1/2-\epsilon,b}([0,T])}, \, \big\|  v_{\leq M} - v_{\leq N} \big\| _{X^{1/2-\delta_1,b}([0,T])}, \\
    &\hspace{7ex}
    \big\|  \Xone_{\leq M} -  \Xone_{\leq N} \big\|_{(L_t^\infty \C_x^{1/2-\delta_1} \cap X^{1/2-\delta_1,b})([0,T])}, \,
    \big\|  \Xtwo_{\leq M} -  \Xtwo_{\leq N} \big\|_{(L_t^\infty \C_x^{1/2-\delta_1} \cap X^{1/2-\delta_1,b})([0,T])}, \\  
    &\hspace{7ex}\big\|  Y_{\leq M} - Y_{\leq N} \big\| _{X^{1/2+\delta_2,b}([0,T])} \bigg) \leq C_4 \exp\Big( C_4 (ABT)^{C_4}\Big) N^{-\theta}. 
    \end{align*}
\end{proposition}

\begin{proof}
Fix any $T\geq 1$, let $J\in \mathbb{N}$ remain to be chosen, and set $\tau:= T/J$. 
We now define a sequence by  $D_0 := N^{-\theta}$ and, for all $0\leq j \leq J$, 
\begin{equation}\label{global:eq-stability-p1}
\begin{aligned}
D_j :=& \, N^{-\theta}
+ \big\| v_{\leq M} - v_{\leq N} \big\|_{X^{1/2-\delta_1,b}([0,j\tau])} 
+ \big\| (\partial_t^2 +1 - \Delta) (v_{\leq M} - v_{\leq N}) \big\|_{X^{-1/2-\delta_1,b-1}([0,j\tau])} \\
&+ \big\| Y_{\leq M} - Y_{\leq N} \big\|_{X^{1/2+\delta_2,b}([0,j\tau])}
+ \big\| (\partial_t^2 +1 - \Delta) (Y_{\leq M} - Y_{\leq N}) \big\|_{X^{1/2+\delta_2,b-1}([0,j\tau])}.
\end{aligned}
\end{equation}
In the following, we assume that 
\begin{equation}\label{global:eq-stability-p2}
D_j \leq B
\end{equation}
for all $0\leq j \leq J$. Since our final estimate \eqref{global:eq-stability-p7} yields a stronger estimate on $D_j$, the assumption \eqref{global:eq-stability-p2} can be verified by using the same estimates as below and a continuity argument\footnote{Alternatively, we could have defined a stopping time $T_\ast \leq T$ such that the upper bound by $B$ holds on $[0,T_\ast]$, replaced all intervals $[0,j\tau]$ by $[0,j\tau] \cap [0,T_\ast]$, and verified through our estimates below that $T_\ast=T$. However, this would be a significant burden on our notation.}.

For any $0\leq j \leq J-1$, it follows from  Lemma \ref{prep:lem-xsb}, Lemma \ref{sttime}, and Lemma \ref{prep:lem-gluing} that 
\begin{equation}\label{global:eq-stability-p3}
\begin{aligned}
D_{j+1} \lesssim& \,  
N^{-\theta} 
+ \big\| (\partial_t^2 +1 - \Delta) (v_{\leq M} - v_{\leq N}) \big\|_{X^{-1/2-\delta_1,b-1}([0,(j+1)\tau])} \\
&+\big\| (\partial_t^2 +1 - \Delta) (Y_{\leq M} - Y_{\leq N}) \big\|_{X^{-1/2+\delta_2,b-1}([0,(j+1)\tau])} \\
\lesssim&\, D_j 
+ \big\| (\partial_t^2 +1 - \Delta) (v_{\leq M} - v_{\leq N}) \big\|_{X^{-1/2-\delta_1,b-1}([j\tau,(j+1)\tau])} \\
&+\big\| (\partial_t^2 +1 - \Delta) (Y_{\leq M} - Y_{\leq N}) \big\|_{X^{-1/2+\delta_2,b-1}([j\tau,(j+1)\tau])} \\ 
\lesssim&\, D_j 
+ \tau^{b_+-b}\big\| (\partial_t^2 +1 - \Delta) (v_{\leq M} - v_{\leq N}) \big\|_{X^{-1/2-\delta_1,b_+-1}([j\tau,(j+1)\tau])} \\
&+ \tau^{b_+-b} \big\| (\partial_t^2 +1 - \Delta) (Y_{\leq M} - Y_{\leq N}) \big\|_{X^{-1/2+\delta_2,b_+-1}([j\tau,(j+1)\tau])} \\ 
\lesssim&\, D_j 
+ \tau^{b_+-b}\big\| (\partial_t^2 +1 - \Delta) (v_{\leq M} - v_{\leq N}) \big\|_{X^{-1/2-\delta_1,b_+-1}([0,(j+1)\tau])} \\
&+ \tau^{b_+-b} \big\| (\partial_t^2 +1 - \Delta) (Y_{\leq M} - Y_{\leq N}) \big\|_{X^{-1/2+\delta_2,b_+-1}([0,(j+1)\tau])}. \\ 
\end{aligned}
\end{equation}
We recall that the evolution equations of $Y_{\leq M}$ and $Y_{\leq N}$ are given by \eqref{ansatz:eq-YN-1}-\eqref{ansatz:eq-YN-8} with frequency-truncation parameters $M$ and $N$, respectively. Due to the identities  
\begin{equation*}
v_{\leq M}= \Xone_{\leq M} + \Xtwo_{\leq M} + Y_{\leq M} \quad \text{and} \quad 
v_{\leq N}= \Xone_{\leq N} + \Xtwo_{\leq N} + Y_{\leq N},
\end{equation*}
the evolution equations of $v_{\leq M}$ and $v_{\leq N}$ are given by the sums of \eqref{ansatz:eq-X11-a}-\eqref{ansatz:eq-X11-d}, \eqref{ansatz:eq-X21}, and \eqref{ansatz:eq-YN-1}-\eqref{ansatz:eq-YN-8} with frequency-truncation parameters $M$ and $N$, respectively.  As already described in the proof of Proposition \ref{ansatz:prop-lwp-quantitative}, a minor modification of our main estimates yields an $A$-certain set $E_A \in \mathcal{E}$ on which 
\begin{equation}\label{global:eq-stability-p4}
\begin{aligned}
&\big\| (\partial_t^2 +1 - \Delta) (v_{\leq M} - v_{\leq N}) \big\|_{X^{-1/2-\delta_1,b_+-1}([0,(j+1)\tau])} \\
&+  \big\| (\partial_t^2 +1 - \Delta) (Y_{\leq M} - Y_{\leq N}) \big\|_{X^{-1/2+\delta_2,b_+-1}([0,(j+1)\tau])} \\
\lesssim& \, A T^{\alpha} \Big( 1+ 
\| v_{\leq M} \|_{X^{1/2-\delta_1,b}([0,(j+1)\tau])}^5 
+\| v_{\leq N} \|_{X^{1/2-\delta_1,b}([0,(j+1)\tau])}^5 
 +  \| Y_{\leq M} \|_{X^{1/2+\delta_2,b}([0,(j+1)\tau])}^5 \\
 &\hspace{6ex} +  \| Y_{\leq N} \|_{X^{1/2+\delta_2,b}([0,(j+1)\tau])}^5 \Big) \\
 &\times \Big( N^{-\theta} 
 + \big\| v_{\leq M} - v_{\leq N} \big\|_{X^{1/2-\delta_1,b}([0,(j+1)\tau])} 
 +\big\| Y_{\leq M} - Y_{\leq N} \big\|_{X^{1/2+\delta_2,b}([0,(j+1)\tau])} \Big). 
\end{aligned}
\end{equation}
By combining the global bounds from our assumption, \eqref{global:eq-stability-p2}, \eqref{global:eq-stability-p3}, and \eqref{global:eq-stability-p4}, we obtain that 
\begin{equation}\label{global:eq-stability-p5}
D_{j+1} \leq C D_j + C A T^{\alpha} \tau^{b_+-b} B^5 D_{j+1},
\end{equation}
where $C=C(\delta_1,\delta_2,\eta,\nu,b_+,b)\gg 1$ is a constant. We now choose 
\begin{equation}\label{global:eq-stability-p6}
\tau := \Big( 2 C A T^{\alpha} B^5 \Big)^{-1/(b_+-b)}. 
\end{equation}
After inserting our choice of $\tau$ into \eqref{global:eq-stability-p5} and iterating, we obtain that
\begin{equation}\label{global:eq-stability-p7}
\begin{aligned}
D_J &\leq (2C)^J D_0 \\
&\leq \exp\Big( \ln(2C) T \big( 2C A T^\alpha B^5\big)^{1/(b_+-b)} \Big) N^{-\theta} \\
&\leq C_4 \exp\bigg( \frac{C_4}{10} (ABT)^{C_4}\bigg) N^{-\theta}. 
\end{aligned}
\end{equation}
Due to the definition of $D_J$, this yields the desired difference estimate for $v_{\leq M}-v_{\leq N}$ and $Y_{\leq M}-Y_{\leq N}$. The difference estimate for $\Xone_{\leq M}-\Xone_{\leq N}$ and $\Xtwo_{\leq M}-\Xtwo_{\leq N}$ then follows from (a minor modification of) Proposition \ref{ansatz:prop-Xj}. Finally, the difference estimate for $u_{\leq M}-u_{\leq N}$ follows from estimates for explicit stochastic objects (Lemma \ref{analytic:lem-cubic} and Lemma \ref{analytic:prop-quintic}) and the estimates for the $v$, $\Xone$, $\Xtwo$, and $Y$-components. \end{proof}

\subsection{Proof of global well-posedness and invariance}\label{section:global-proof}

Equipped with the para-controlled global bounds (Proposition \ref{global:prop-para-controlled}) and stability theory (Proposition \ref{global:prop-stability}), we now prove our main theorem (Theorem \ref{intro:thm-rigorous}). 

\begin{proof}[Proof of Theorem \ref{intro:thm-rigorous}:]
We first prove the global well-posedness, i.e., the global existence of the limiting dynamics. Using time-reversal symmetry, it suffices to prove for all $T\geq 1$ that 
\begin{equation}\label{global:eq-global-limit}
\begin{aligned}
\lim_{K \rightarrow \infty} \mu \bigg( \Big \{& 
\big( \phi^{\cos}, \phi^{\sin} \big) \in \mathscr{H}_x^{-1/2-\epsilon}\colon 
\big\| u_{\leq N_1} - u_{\leq N_2} \big\|_{C_t^0 \mathscr{H}_x^{-1/2-\epsilon}([0,T] \times \T^3)} \leq \min\big(N_1,N_2\big)^{-\theta/2}\\ &\text{for all } N_1,N_2\geq K\Big\} \bigg)=1. 
\end{aligned}
\end{equation}
Here, $u_{\leq N}$ denotes the unique global solution of the frequency-truncated cubic wave equation with initial data $(\phi^{\cos},\phi^{\sin})$. For any frequency-scales $K\leq L$, we define 
\begin{equation}
\begin{aligned}
\Diff_{K,L,T} := \Big\{& 
\big( \phi^{\cos}, \phi^{\sin} \big) \in \mathscr{H}_x^{-1/2-\epsilon}\colon 
\big\| u_{\leq N_1} - u_{\leq N_2} \big\|_{C_t^0 \mathscr{H}_x^{-1/2-\epsilon}([0,T] \times \T^3)} \leq \min\big(N_1,N_2\big)^{-\theta/2} \\
&\text{for all } K \leq N_1,N_2 \leq L \Big\} 
\end{aligned}
\end{equation}
Using the continuity from above of general probability measures, using the weak convergence from \eqref{intro:eq-Gibbs-limit}, and using that $\Diff_{K,L,T}$ is closed in $\mathscr{H}_x^{-1/2-\epsilon}$, it follows that 
\begin{equation}
\text{LHS of } \eqref{global:eq-global-limit} = \lim_{K\rightarrow \infty} \lim_{L\rightarrow \infty} \mu\Big( \Diff_{K,L,T} \Big) \geq \lim_{K\rightarrow \infty} \lim_{L\rightarrow \infty} \limsup_{M\rightarrow \infty} \mu_{\leq M} \Big( \Diff_{K,L,T} \Big). 
\end{equation}
As a result, it suffices to prove that for all $A\geq 1$ and all $T\geq 1$, there exists a frequency-scale $K_0=K_0(A,T)$ such that 
\begin{equation}\label{global:eq-global-limit-reduced}
\mu_{\leq M} \Big( \Diff_{K,L,T} \Big) \geq 1- 2 c_4 \exp\Big(-c_4 A^{c_4} \Big) 
\end{equation}
for all frequency-scales $M\geq L \geq K \geq K_0$. \\

Similar as in our proof of local well-posedness (Proposition \ref{ansatz:prop-lwp-qualitative}), our proof of \eqref{global:eq-global-limit-reduced} relies on the caloric initial data. To this end, let $\initial{blue}$, $\initial{green}[][\leqM]$, and $\initial{red}[][\leqM]$ be as in Proposition \ref{ansatz:prop-caloric}. Then, it holds that 
\begin{equation*}
\mu_{\leq M} \Big( \Diff_{K,L,T} \Big) = \Big( \mathbb{P} \otimes \mathbb{Q} \Big) \Big(  \Big\{ \initial{blue} - \initial{green}[][\leqM] + \initial{red}[][\leqM] \in \Diff_{K,L,T} \Big\} \Big). 
\end{equation*}
We now let the events $\operatorname{Glb}_{M,A}$ and $\Stbl_{M,A}$ be as in Proposition \ref{global:prop-para-controlled} and Proposition \ref{global:prop-stability}, respectively. Using the corresponding probability estimates, it follows that 
\begin{equation*}
\big( \mathbb{P} \otimes \mathbb{Q} \big) \Big( \operatorname{Glb}_{M,A} \medcap \Stbl_{M,A} \Big) \geq 1 - 2 c_4 \exp\Big( -c_4 A^{c_4} \Big). 
\end{equation*}
As a result, it suffices to prove the inclusion
\begin{equation}
\operatorname{Glb}_{M,A} \medcap \Stbl_{M,A} \subseteq \Big\{ \initial{blue} - \initial{green}[][\leqM] + \initial{red}[][\leqM] \in \Diff_{K,L,T} \Big\}
\end{equation}
for all frequency-scales $M\geq L \geq K \geq K_0$. As long is $K_0=K_0(A,T)$ is sufficiently large, this follows directly from the global bounds (from Proposition \ref{global:prop-para-controlled}.\ref{global:item-bounds-2}), the stability estimates (from Proposition \ref{global:prop-stability}), and the triangle inequality. This completes the proof of global well-posedness.\\

It now only remains to prove the invariance of the Gibbs measure. To this end, let $f\colon \mathscr{H}_x^{-1/2-\epsilon} \rightarrow \R$ be bounded and globally Lipschitz, let $T\geq 1$, and let $t\in [0,T]$. In order to prove invariance, it suffices to prove that 
\begin{equation}\label{global:eq-invariance}
\E_{\mu} \Big[ f\big( u[t] \big) \Big] = \E_{\mu} \Big[ f\big( u[0] \big) \Big]. 
\end{equation}
Using the definition of the limiting dynamics $u$, it follows that
\begin{equation*}
    \E_{\mu} \Big[ f\big( u[t] \big) \Big] = \lim_{N\rightarrow \infty}  \E_{\mu} \Big[ f\big( u_{\leq N}[t] \big) \Big]. 
\end{equation*}
For any fixed $N\geq 1$, the frequency-truncated solution $u_{\leq N}[t]$ depends continuously on the initial data. Together with the weak convergence of the Gibbs measures from \eqref{intro:eq-Gibbs-limit}, it follows that
\begin{equation*}
    \lim_{N\rightarrow \infty}  \E_{\mu} \Big[ f\big( u_{\leq N}[t] \big) \Big]
    = \lim_{N\rightarrow \infty} \lim_{M\rightarrow \infty} \E_{\mu_{\leq M}} \Big[ f\big( u_{\leq N}[t] \big) \Big].
\end{equation*}
Using our estimate from the proof of global well-posedness, i.e., \eqref{global:eq-global-limit-reduced}, it follows that 
\begin{equation*}
    \lim_{N\rightarrow \infty} \lim_{M\rightarrow \infty} \E_{\mu_{\leq M}} \Big[ f\big( u_{\leq N}[t] \big) \Big] = \lim_{M\rightarrow \infty} \E_{\mu_{\leq M}} \Big[ f\big( u_{\leq M}[t] \big) \Big]
\end{equation*}
Since the frequency-truncated solution $u_{\leq M}[t]$ preserves the Gibbs measure, it follows that 
\begin{equation*}
 \lim_{M\rightarrow \infty} \E_{\mu_{\leq M}} \Big[ f\big( u_{\leq M}[t] \big) \Big]
 =  \lim_{M\rightarrow \infty} \E_{\mu_{\leq M}} \Big[ f\big( u_{\leq M}[0] \big) \Big]
 = \lim_{M\rightarrow \infty} \E_{\mu_{\leq M}} \Big[ f \big( (\phi^{\cos},\phi^{\sin})\big)\Big]. 
 \end{equation*}
 After again using the weak convergence of the Gibbs measures from \eqref{intro:eq-Gibbs-limit}, this completes the proof of \eqref{global:eq-invariance} and hence the proof of invariance.
\end{proof}

\setlength\parindent{15pt}

\section{Integer lattice counting and basic tensors estimates}
In this section we make some preparations for the 
proof, namely we introduce some counting estimates for integer lattice points and prove the basic random tensors estimates for the base tensor (Subsection \ref{base}), the cubic stochastic tensor (Subsection \ref{cubic-tensor}), the quintic stochastic tensor (Subsection \ref{quintic-tensor}) and the sine-cancellation tensor (Subsection \ref{sine-tensor}). Once established, these estimates will be used as a black box throughout the proofs of the estimates in Sections \ref{section:analytic}--\ref{section:analytic2} below. 

\label{section:counting}
\subsection{Lattice point counting estimates}
We start with some basic integer lattice point counting bounds which have appeared already in the literature (see for example \cite{B20II}) but which we record for the sake of completeness in the following three lemmas.

\begin{lemma}[A basic counting lemma]\label{counting:lem1}
Given dyadic numbers $A, N$, and $a\in \mathbb{Z}^3$ satisfying $|a|_{\infty}\sim A$, we have the following lattice point counting bounds:
\begin{equation}\label{counting:eq1}
\sup_{m\in \mathbb{Z}} \# \{n\in \mathbb{Z}^3 :  |n|_{\infty}\sim N,\, |\langle a+n\rangle\pm\langle n\rangle -m|\leq 1\} \lesssim \min (A, N)^{-1}  N^3.
\end{equation}
\begin{equation}\label{counting:eq2}
\sup_{m\in \mathbb{Z}} \# \{n\in \mathbb{Z}^3 :  |n|_{\infty}\sim N, \,|\langle a+n\rangle+\langle n\rangle -m|\leq 1\} \lesssim N^2.
\end{equation}
\end{lemma}
\begin{proof}
Since the $\ell^2$ and $\ell^\infty$-norms on $\mathbb{Z}^3$ are comparable, the conditions $|a|_\infty \sim A$ and $|n|_\infty \sim N$ can be replaced by $|a| \sim A$ and $|n|\sim N$ respectively{\footnote{We will use this fact throughout the paper when convenient without any further explanation.}}.
The proof of the estimate (\ref{counting:eq1}) is same as Lemma 4.15 in \cite{B20II}. The estimates (\ref{counting:eq2})  can also be proved using the same argument, except the case when $A\ll N$, so we only consider (\ref{counting:eq2})    with $A\ll N$.
Since $\big|\langle\xi\rangle-|\xi|\big|\leq 1$ for all $\xi\in \mathbb{R}^3$, the $\langle\cdot \rangle$ in (\ref{counting:eq2})   can be replaced by $|\cdot|$ after increasing the implicit constant.

Since the integer vectors are 1-separated, it follows that 
\[
 \# \{n\in \mathbb{Z}^3 : |n|\sim N,\,  || a+n|+| n| -m|\lesssim 1\}\lesssim \mathrm{Leb}\left( \{\xi\in \mathbb{R}^3\,:\,|\xi|\sim N, \, || a+\xi|+|\xi| -m|\lesssim 1\}\right).
\]
We now decompose the above Lebesgue measure in the following way:
\begin{equation*}
\begin{aligned}
&\mathrm{Leb}\left( \{\xi\in \mathbb{R}^3 :  |\xi|\sim N, \, || a+\xi|+|\xi| -m|\leq 1\}\right)\\
\lesssim &\sum_{
\substack {m_1,m_2\in \mathbb{Z}\\ |m_1+m_2-m|\lesssim 1
\\|m_1-m_2|\lesssim |a|\sim A}} \mathrm{Leb}\left( \{\xi\in \mathbb{R}^3 :  |\xi|\sim N,\, |a+\xi|=m_1+\mathcal{O}(1), \, |\xi|=m_2+\mathcal{O}(1)\}\right)\\
\lesssim&\, A \cdot\sup_{m_1, m_2\in \mathbb{Z}} \mathrm{Leb}\left( \{\xi\in \mathbb{R}^3 : |\xi|\sim N, \, |a+\xi|=m_1+\mathcal{O}(1),\, |\xi|=m_2+\mathcal{O}(1)\}\right).
\end{aligned}
\end{equation*}
Using the above decomposition, it will suffice to prove that
\begin{equation}\label{counting:eq2''}
\mathrm{Leb}\left( \{\xi\in \mathbb{R}^3 : |\xi|\sim N,\, |a+\xi|=m_1+\mathcal{O}(1), \,|\xi|=m_2+\mathcal{O}(1)\}\right)\lesssim \min(A,N)^{-1} N^2
\end{equation}
which directly follows the same argument in the proof of Lemma 4.15 in \cite{B20II}; we omit the details.\end{proof}

\begin{lemma}[A box-counting lemma]\label{counting:lem-annuli-box} Given dyadic numbers $A, N$, and $a\in \mathbb{Z}^3$ such that $|a|_\infty\sim A$,  let $\zeta \in \Z^3$ be arbitrary. Then, we have the following lattice point counting estimate:
\begin{equation}\label{counting:eq2'}
\sup_{m\in \mathbb{Z}} \# \{n\in \mathbb{Z}^3: |n|_\infty \sim N, \, |n-\zeta|_\infty\lesssim A,  \, |\langle a+n\rangle\pm\langle n\rangle -m|\leq 1\} \lesssim N^2.
\end{equation}
\end{lemma}
\begin{proof}
By our assumptions we have that $|a| \sim A$, $|n|\sim N$, and $|n-\zeta|\lesssim A$. We then argue essentially as in the proof of Lemma \ref{counting:lem1}, but use the following estimate for the Lebesgue measure: 
\begin{align*}
&\mathrm{Leb}\left( \{\xi\in \mathbb{R}^3 :  |\xi|\sim N, \,  |\xi - \zeta| \lesssim A,  \, || a+\xi|\pm|\xi| -m|\leq 1\}\right)\\
\lesssim &\sum_{
\substack {m_1,m_2\in \mathbb{Z}\\ |m_1 \pm m_2-m|\lesssim 1
\\|m_2 - |\zeta||\lesssim A}} \mathrm{Leb}\left( \{\xi\in \mathbb{R}^3 :  |\xi|\sim N,\, |a+\xi|=m_1+\mathcal{O}(1), \, |\xi|=m_2+\mathcal{O}(1)\}\right)\\
\lesssim&\, A \cdot\sup_{m_1, m_2\in \mathbb{Z}} \mathrm{Leb}\left( \{\xi\in \mathbb{R}^3 : |\xi|\sim N, \, |a+\xi|=m_1+\mathcal{O}(1),\, |\xi|=m_2+\mathcal{O}(1)\}\right).\qedhere
\end{align*}
\end{proof}

The next lemma, improves upon \eqref{counting:eq1} provided the two vectors $n$ and $n+a$ satisfy what in \cite{DNY19, DNY21} is called the $\Gamma$-condition. The $\Gamma$-condition needed in our case is with respect of the $| \cdot |_{\infty}$-norm but we still call it $\Gamma$-condition. More precisely, given $\Gamma\in\Rb$, we say that two vectors $n_1,\,n_2\in\Zb^3$ satisfy the $\Gamma$-condition if 
\begin{equation}\label{counting:def-gamma}
\text{either}\quad |n_1|_{\infty} \leq \Gamma \leq  |n_2|_{\infty} \quad \text{or} \quad |n_2|_{\infty}\leq \Gamma \leq  |n_1|_{\infty}.  
\end{equation}

\begin{lemma}[The $\Gamma$-condition counting lemma]\label{counting:lem1''}
Given $\Gamma\in \mathbb{R}$,  dyadic numbers $A, N$, and $a\in \mathbb{Z}^3$ such that $|a|_\infty\sim A$, we have the following lattice point counting bounds:
\begin{equation}\label{counting:gamma}
\sup_{m\in \mathbb{Z}} \# \{n\in \mathbb{Z}^3: |n|_{\infty}\sim N, \, |\langle a+n\rangle\pm\langle n\rangle -m|\leq 1,\, |n|_\infty\geq \Gamma \geq |n+a|_\infty\} \lesssim  N^2 \log N.
\end{equation} The same bound holds if one assumes $|n|_\infty\leq \Gamma \leq |n+a|_\infty$.
\end{lemma}

\begin{proof} 
If $A\geq N/100$ then (\ref{counting:gamma}) follows directly from (\ref{counting:eq1}); so below we will assume $A\leq N/100$. In particular $\langle n\rangle\sim N\sim \langle a+n\rangle$, so by symmetry we only need to prove (\ref{counting:gamma}).

Denote by $n=(x,y,z)$ and $a=(x_a,y_a,z_a)$ two vectors in $\mathbb{Z}^3$. Without loss of generality, suppose $|n|_\infty=|x|$, then $|n|_\infty\geq \Gamma \geq |n+a|_\infty$ implies that $|x|\geq \Gamma\geq |x+x_a|\geq |x|-|x_a|$. Now define $\delta \in 2^\Zb\cap[N^{-10},1]$ (with the understanding that the $\delta< N^{-10}$ case is included in $\delta=N^{-10}$) such that
\[\left|\frac{n}{|n|} - \frac{n+a}{|n+a|}\right|=2 \sin \big(\frac{ \angle(n,n+a)}{2}\big)\in[\delta,2\delta],
\]where $\angle(n,n+a)$ is the angle between the vectors $n$ and $n+a$. Since $|a|\sim A \leq  N/100 \sim |n|/100$, we have that $\delta\leq  1/10$.
Define the unit vectors $n/|n|=n'=(x',y',z')$ and $(a+n)/|a+n| = n''= (x'',y'',z'')$,  we claim that 
\begin{equation}\label{counting:gammaclaim}\max(| y'- y''|,|z'-z''|)\geq \delta/100.\end{equation}
In fact, suppose the contrary, then $|y'-y''|<\delta/100$ and $|z'-z''|<\delta/100$. Since $|n'-n''|\geq \delta$, we have $|x'-x''|\geq \delta/2$. But $x'$ and $x''$ must have the same sign because $|x|=|n|_{\infty}>10A\geq 5|x_a|$, and we also have $|x'|\geq 1/2$ because  $|x|=|n|_{\infty}\geq |n|/2$. Since $(x',y',z')$ and $(x'',y'',z'')$ are both unit vectors, we then have
\[\delta/2\leq |x'-x''|=\frac{|(x')^2-(x'')^2|}{|x'|+|x''|}=\frac{|(y')^2+(z')^2-(y'')^2-(z'')^2|}{|x'|+|x''|}\leq 4(|y'-y''|+|z'-z''|)\leq \delta/25,\] which is a contradiction.

Now, using \eqref{counting:gammaclaim}, without loss of generality, we may assume $|y'-y''|\geq \delta/100$. Let $f(n):= |n|-|n+a|$ so $\nabla f(n) = \frac{n}{|n|} - \frac{n+a}{|n+a|} = n'-n''=(x'-x'',y'-y'',z'-z'')$ with $|y'-y''|\geq \delta/100$.
Since by, assumption, $|f(n)-m|\lesssim 1$ and $|\frac{\partial f(n)}{\partial y}|=|y'-y''|\geq \delta/100$, we know that for any fixed values of $x$ and $z$, the number of choices for $y$ is $\lesssim 1/\delta$. 

Now it remains to count the choices of $x$ and $z$. By the law of sines we have
\begin{equation}\label{eq:xzbound}
|xz_a - zx_a|\leq |n\times a|= |n|\cdot|a|\cdot\sin (\angle(n,a))\sim 
N \cdot |a|\cdot \frac{|n+a|}{|a|} \sin (\angle(n,n+a))\lesssim N^2\cdot\delta.
\end{equation}
By assumption $|x|\geq \Gamma\geq |x|-|x_a|$, so the number of the choices of $x$ is at most $|x_a|$; when $x$ is fixed, by \eqref{eq:xzbound} we know that the number of choices of $z$ is at most $N^2 \delta/|x_a|$. Thus the number of choices for $(x,y,z)$ is 
\[\lesssim |x_a|\cdot\frac{N^2\delta}{|x_a|}\cdot\frac{1}{\delta}=N^2.\]
In summary, we have that the number of the choices of $(x,y,z)$ $\lesssim N^2$ with any fixed $\delta$. Summing over $\delta$ yields the desired bound in (\ref{counting:gamma}). 
\end{proof}

In the following two lemmas we prove the basic atomic estimates that will be used in later proofs, especially in Section \ref{section:analytic2}. Lemma \ref{counting:lem2} contains the estimates for a single nonlinear interaction, which corresponds to a single ``atom" in the sense of Section \ref{section:analytic2} (and as in \cite{DH21}), and may involve 2, 3, or 4 unknown vectors. Similarly, Lemma \ref{counting:lem3} contains estimates corresponding to two connected atoms, which involve at least 2 and at most 7 unknown vectors. They are listed separately, because in some cases they provide better bounds than if one simply applies Lemma \ref{counting:lem2} twice.

\begin{lemma}[Lattice point counting I]\label{counting:lem2}
Given $2\leq q\leq 4$, $\pm_j\in\{\pm\}$ and dyadic numbers $N_j\geq 1$ (for $1\leq j\leq q$), and $(n_{\mathrm{ex}},m)\in\mathbb{Z}^3\times\mathbb{R}$, consider the set
\begin{equation}\label{counting:eq3}
\mathcal{M}_q=\bigg\{(n_1,\cdots,n_q)\in(\mathbb{Z}^3)^q:\langle n_j\rangle\sim N_j,\,\,\,\sum_{j=1}^q(\pm_j)n_j=n_{\mathrm{ex}},\,\,\,\bigg|\sum_{j=1}^q(\pm_j)\langle n_j\rangle-m\bigg|\leq 1\bigg\}.
\end{equation} Assume that $\langle n_{\mathrm{ex}}\rangle\sim M$, and that $n_{\mathrm{ex}}=0\,(M=1)$ when $q=4$, also denote $N^{(1)}\geq\cdots\geq N^{(q)}$ be the decreasing rearrangement of $N_j$, and define $\pm^{(j)}$ correspondingly. Then we have the following bounds.
\begin{enumerate}
\item If $q=2$ we have
\begin{equation}\label{counting:eq4}\#\mathcal{M}_2\lesssim \left\{
\begin{aligned}&\max((N^{(2)})^2,(N^{(2)})^3M^{-1}),&\textrm{if }\pm^{(1)}=\mp^{(2)};\\
&(N^{(2)})^2,&\textrm{if }\pm^{(1)}=\pm^{(2)}.
\end{aligned}
\right.
\end{equation} If we moreover assume $\pm_1=\mp_2$ and that the two vectors
$n_1$ and $n_2$ satisfy $\Gamma$-condition (\ref{counting:def-gamma}), then we have \begin{equation}\label{counting:eq4+}\#\Mc_2\lesssim (N^{(2)})^2M.\end{equation}

\item If $q=3$ we have 
\begin{equation}\label{counting:eq5}\#\mathcal{M}_3\lesssim (N^{(2)})^3(N^{(3)})^3(\mathrm{med}(N^{(2)},N^{(3)},M))^{-1}\lesssim (N^{(2)})^3(N^{(3)})^2.
\end{equation}

\item If $q=4$ we have
\begin{equation}\label{counting:eq6}\#\mathcal{M}_4\lesssim (N^{(2)})^3(N^{(3)})^2(N^{(4)})^3.
\end{equation} If we moreover assume $|(\pm_1)n_1+(\pm_2)n_2|\lesssim L$, then we have 
\begin{equation}\label{counting:eq6+}\#\mathcal{M}_4\lesssim L(N_1\cdots N_4)^2(\max(N_1,N_2))^{-1}(\max(N_3,N_4))^{-1}.
\end{equation}
\item Summarizing (\ref{counting:eq4}), (\ref{counting:eq5}) and (\ref{counting:eq6}), we have
\begin{align}\label{counting:eqsum1}
\#\Mc_q&\lesssim (N_1\cdots N_q)^2\cdot (N^{(1)})^{-1},&\mathrm{if\ }&q\leq 3;\\
\label{counting:eqsum2}\#\Mc_q&\lesssim (N_1\cdots N_q)^2\cdot N^{(4)}(N^{(1)})^{-1},& \mathrm{if\ }&q=4.
\end{align}
\end{enumerate}
\end{lemma}
\begin{proof} Except in the proof of (\ref{counting:eq6+}), by symmetry we may assume $N_1\geq N_2\geq\cdots$, so $N^{(j)}=N_j$. First (\ref{counting:eq4}) follows directly from (\ref{counting:eq1}) and (\ref{counting:eq2}) with $N$ replaced by $N_2$ and $A$ replaced by $M$. As for (\ref{counting:eq4+}), by the definition of $|\cdot|_\infty$, we may assume $|(n_1)^1|\geq \Gamma \geq |(n_2)^1|$ or $|(n_2)^1|\geq \Gamma \geq |(n_1)^1|$ where $(\cdot)^1$ represents the first coordinate. In either case we have $|(n_2)^1|\in [\Gamma-O(M),\Gamma+O(M)]$, so $(n_2)^1$ has $\lesssim M$ choices. The other coordinates of $n_2$ each has $\lesssim N_2$ choices, so (\ref{counting:eq4+}) is true.

Next we prove (\ref{counting:eq5}). Let $|n_{\mathrm{ex}}-(\pm_3)n_3|\sim R$, then for fixed $n_3$, the number of choices for $(n_1,n_2)$ is $\lesssim N_2^3(\min(N_2,R))^{-1}$ by (\ref{counting:eq4}), while the number of choices for $n_3$ is $\lesssim\min(N_3,R)^3$. Now, if $M\lesssim N_3$, then $R\lesssim N_3$ and $\mathrm{med}(N^{(2)},N^{(3)},M)\sim N_3$, so 
\[\#\mathcal{M}_3\lesssim\sum_{R\lesssim N_3} R^3\cdot N_2^3R^{-1}\lesssim N_2^3N_3^2;\] if $M\gg N_3$, then $R\sim M$ and $\mathrm{med}(N^{(2)},N^{(3)},M)\sim \min(N_2,R)$, so
\[\#\mathcal{M}_3\lesssim N_3^3\cdot N_2^3(\min(N_2,R))^{-1}.\] Either way we have proved (\ref{counting:eq5}). Now (\ref{counting:eq6}) is a consequence of (\ref{counting:eq5}) because we may first fix $n_4$ which has $N_4^3$ choices, then apply (\ref{counting:eq5}) to get
\[\#\mathcal{M}_4\lesssim N_4^3\cdot N_2^3N_3^2.\] Clearly (\ref{counting:eqsum1}) also follows from (\ref{counting:eq4}), (\ref{counting:eq5}) and (\ref{counting:eq6}).

Finally we prove (\ref{counting:eq6+}). By symmetry we may assume $N_1\geq N_2$ and $N_3\geq N_4$. First we may fix the value of $(\pm_1)n_1+(\pm_2)n_2$, which has at most $L^3$ choices. Once these are fixed, we then apply (\ref{counting:eq4}) for $(n_1,n_2)$ and $(n_3,n_4)$ separately to get
\[\#\Mc_4\lesssim L^3\cdot(N_2^2+N_2^3L^{-1})(N_4^2+N_4^3L^{-1})\lesssim L^3\cdot(N_1N_2^2L^{-1})\cdot(N_3N_4^2L^{-1})=L(N_1N_2N_3N_4)^2(N_1N_3)^{-1}\] since $N_1\gtrsim\max(N_2,L)$ and $N_3\gtrsim\max(N_4,L)$. This proves (\ref{counting:eq6+}).
\end{proof}
\begin{lemma}[Lattice point counting II]\label{counting:lem3}
Given $q,r\geq 1$ such that $q+r\leq 4$, $\pm_j,\pm_j'\in\{\pm\}$ (for $1\leq j\leq q$) and $\pm_j''\in\{\pm\}$ (for $1\leq j\leq r$), and dyadic numbers $N_j,M_j\geq 1$ (for $1\leq j\leq q$) and $L_j\geq 1$ (for $1\leq j\leq r$), and $(n_{\mathrm{ex}},n_{\mathrm{ex}}',m,m')\in(\Zb^3)^2\times\Rb^2$, consider the set
\begin{multline}
\label{counting:eq7}\Mc_{q,r}=\bigg\{(n_1,\cdots,n_q,\ell_1,\cdots,\ell_r,n_1',\cdots,n_q')\in(\Zb^3)^{q+r+q}:\langle n_j\rangle\sim N_j,\,\,\langle n_j'\rangle\sim M_j,\,\,\langle \ell_j\rangle\sim L_{j},\,\,
\\\sum_{j=1}^q(\pm_j)n_j+\sum_{j=1}^r(\pm_{j}'')\ell_j=n_{\mathrm{ex}},\,\,\sum_{j=1}^q(\pm_j')n_j'+\sum_{j=1}^r(\pm_{j}'')\ell_j=n_{\mathrm{ex}}',\\
\bigg|\sum_{j=1}^q(\pm_j)\langle n_j\rangle+\sum_{j=1}^r(\pm_{j}'')\langle \ell_j\rangle-m\bigg|\leq 1,\,\,\bigg|\sum_{j=1}^q(\pm_j')\langle n_j'\rangle+\sum_{j=1}^r(\pm_{j}'')\langle \ell_j\rangle-m'\bigg|\leq 1\bigg\}.
\end{multline} Denote by $N^{(1)}\geq\cdots\geq N^{(q)}$ the decreasing rearrangements of $N_j$, and similarly define $M^{(j)}$ and $L^{(j)}$. Then we have that:
\begin{align}
\label{counting:eqsum3}
\#\Mc_{q,r}&\lesssim(N_1\cdots N_qM_1\cdots M_qL_1\cdots L_r)^2\cdot (N^{(1)})^{-1}(M^{(1)})^{-1},&\mathrm{if\ }&q+r=4,\,\,q\geq 1,\,r\geq 2;\\
\label{counting:eqsum4}
\#\Mc_{q,r}&\lesssim (N_1\cdots N_qM_1\cdots M_qL_1\cdots L_r)^2\cdot (N^{(1)})^{-1}(M^{(1)})^{-1}L^{(1)},&\mathrm{if\ }&q=3,\,\,r=1;\\
\label{counting:eqsum5}
\#\Mc_{q,r}&\lesssim (N_1\cdots N_qM_1\cdots M_qL_1\cdots L_r)^2\cdot (N^{(1)})^{-1}(M^{(1)})^{-1}(L^{(1)})^{-1},
&\mathrm{if\ }&q+r\leq 3,\,\,q,r\geq 1.
\end{align}
\end{lemma}
\begin{proof} First, we shall reduce to the case where $N_j=M_j$, $\pm_j=\pm_j'$ and $n_{\mathrm{ex}}=n_{\mathrm{ex}}'$, $m=m'$. In fact, for fixed choice of $\ell:=(\ell_1,\cdots,\ell_r)$, the number of choices for $(n_1,\cdots,n_q)$ depends on the value of $\ell$ and $(N_j,\pm_j,n_{\mathrm{ex}},m)$, so we may denote it by $\Gc_{\ell,N_j,\pm_j,n_{\mathrm{ex}},m}$; then we have
\[\#\Mc_{q,r}=\sum_{\ell}\Gc_{\ell,N_j,\pm_j,n_{\mathrm{ex}},m}\cdot \Gc_{\ell,M_j,\pm_j',n_{\mathrm{ex}}',m'}\leq\bigg(\sum_{\ell}\Gc_{\ell,N_j,\pm_j,n_{\mathrm{ex}},m}^2\bigg)^{1/2}\bigg(\sum_{\ell}\Gc_{\ell,M_j,\pm_j',n_{\mathrm{ex}}',m'}^2\bigg)^{1/2}.\] Therefore, if we can prove (\ref{counting:eqsum3})--(\ref{counting:eqsum5}) with $M_j$ replaced by $N_j$ etc., then by symmetry we also have the same estimates with $N_j$ replaced by $M_j$ etc., and so by taking geometric average we obtain (\ref{counting:eqsum3})--(\ref{counting:eqsum5}) with both $N_j$ etc. and $M_j$ etc.

Now we may assume $M_j=N_j$ etc. By symmetry we may assume $N_1\geq\cdots\geq N_q$ and $L_1\geq\cdots\geq L_r$, so $N^{(j)}=N_j$ and $L^{(j)}=L_j$. Most of the desired bounds directly follow from applying (\ref{counting:eqsum1}) or (\ref{counting:eqsum2}) once or twice. Indeed, if $q=r=1$, then (\ref{counting:eqsum1}) implies $\#\Mc_{1,1}\lesssim N_1^2L_1$ because $n_1'$ is uniquely fixed once $n_1$ and $\ell_1$ are fixed. Similarly, if $(q,r)=(1,2)$ then (\ref{counting:eqsum1}) implies $\#\Mc_{1,2}\lesssim (N_1L_2)^2L_1$, and if $(q,r)=(1,3)$ then (\ref{counting:eqsum2}) implies $\#\Mc_{1,3}\lesssim (N_1L_1L_2L_3)^2$. Moreover, if $(q,r)=(3,1)$ then we first apply (\ref{counting:eqsum2}) for $(n_1,n_2,n_3,\ell_1)$ and then apply (\ref{counting:eqsum1}) for $(n_1',n_2',n_3')$ to get $\#\Mc_{3,1}\lesssim N_1(N_2N_3)^2L_1^3(N_2N_3)^2N_1=N_1^2(N_2N_3)^4L_1^3$.

We are left with the cases where $q=2$ and $r\in\{1,2\}$. If $r=1$, let $\langle n_{\mathrm{ex}}-(\pm_3)\ell_1\rangle\sim R$, then either $R\lesssim N_2$ or $R\sim N_1$. Moreover, for fixed $\ell_1$, the number of choices for $(n_1,n_2,n_1',n_2')$ is $\lesssim N_2^6(\min(N_2,R))^{-2}$ due to (\ref{counting:eq4}), while the number of choices for $\ell_1$ is $\lesssim(\min(L_1,R))^3$. Therefore we have
\[\#\Mc_{2,1}\lesssim\sum_{R\lesssim N_2\textrm{ or }R\sim N_1}N_2^6(\min(N_2,R))^{-2}\cdot(\min(L_1,R))^3\lesssim N_1^2N_2^4L_1,\] which can be directly verified by enumerating the cases (namely, when $R\sim N_1$, or $L_1\lesssim R\lesssim N_2$, or $R\lesssim\min(L_1,N_2)$). Finally, if $q=2$ and $r=2$, then we may first fix $\ell_2$, which has $\lesssim L_2^3$ choices, and reduce to the case $(q,r)=(2,1)$. This gives the bound
\[\#\Mc_{2,2}\lesssim N_1^2N_2^4L_1L_2^3\lesssim N_2^4(N_1L_1L_2)^2.\] In any case this provides the needed bounds for (\ref{counting:eqsum3})--(\ref{counting:eqsum5}).
\end{proof}

\subsection{Tensors and tensor norms}\label{def_tensor} Before proving our basic tensors estimates in the next subsections,  we recall the definition of tensors and $L^2$-based tensor norms that will be used in this paper. In what follows the we will use capital letters  such as $A$ to denote a set of indices and we will denote by $n_A=(n_j:j\in A)$ tensor variables.  The presentation follows \cite{DNY20} where these were first introduced.

\begin{definition}\label{tensornorms}  A \emph{tensor} $h=h_{n_A}$ is a function $(\Zb^d)^A\to \Cb$, with $n_A$ being the input variables. The \emph{support} of $h$ is the set of $n_A$ such that $h_{n_A}\neq 0$. These tensors may depend on other parameters, such as $t$ or $(\lambda,\lambda_j)$, in which case we will write respectively $h_{n_A}=h_{n_A}(t)$ or $h_{n_A}=h_{n_A}(\lambda,\lambda_j)$. Note that unlike \cite{DNY20}, the tensors we will use here are not random; however the norms in which they are estimated are the same as in \cite{DNY20}.

A \emph{partition} of $A$ is a pair of sets $(B,C)$ such that $B\cup C=A$ and $B\cap C=\varnothing$. For such $(B,C)$ define the norm $\|\cdot\|_{n_B\to n_C}$ by
 \[\|h\|_{n_B\to n_C}^2=\sup\bigg\{\sum_{n_C}\bigg|\sum_{n_B}h_{n_A}\cdot z_{n_B}\bigg|^2:\sum_{n_B}|z_{n_B}|^2=1\bigg\}.\]By duality we have that 
 \begin{equation}\label{duality}\|h\|_{n_B\to n_C}=\sup\bigg\{\bigg|\sum_{n_B,n_C}h_{n_A}\cdot z_{n_B}\cdot y_{n_C}\bigg|:\sum_{n_B}|z_{n_B}|^2=\sum_{n_C}|y_{n_C}|^2=1\bigg\},
 \end{equation} hence $\|h\|_{n_B\to n_C}=\|h\|_{n_C\to n_B}$. If $B=\varnothing$ or $C=\varnothing$ we get the Hilbert-Schmidt norm $\|\cdot\|_{n_A}$ defined by
 \[\|h\|_{n_A}^2=\sum_{n_A}|h_{n_A}|^2.\] Note that trivially $\|h\|_{n_B\to n_C}\leq \|h\|_{n_A}$. Finally, for a tensor $h_{n_A}(\lambda)$ depending on $\lambda$, we define the norm $\|\cdot\|_{L^2_\lambda[n_B\to n_C]}$ by
\begin{equation}\label{def:L2lambda}
\|h_{n_A}(\lambda)\|_{L^2_\lambda[n_B\to n_C]} =\big\|\|h_{n_A}(\lambda)\|_{n_B\to n_C}\big\|_{L^2_\lambda}.
\end{equation}
\end{definition}

\subsection{Base tensors estimates}\label{base}
In this subsection, we introduce the base tensor $h^b_{nn_1n_2n_3} $, which describes the basic relations among the spatial frequencies $(n_j)_{j=1,2,3}$ of each component in the nonlinearity and the spatial frequency $n$ of the solution $u$, and estimate the corresponding Hilbert-Schmidt and operator norms. 
The base tensor $h^b_{nn_1n_2n_3} $ is given by 
\begin{equation}\label{base-tensor}
h^b_{nn_1n_2n_3}:= 1_{N}(n)\cdot \prod_{j=1}^3 1_{N_j}(n_j) \cdot \mathbf{1}{\big\{n=n_{123}\big\}} \cdot \mathbf{1}{\big\{|\Omega(n_1,n_2,n_3)-m|\leq 1\big\}},
\end{equation} where  $N$ and $N_j$'s ($j\in \{1,2,3\}$)  are dyadic numbers, $m\in \mathbb{R}$, and \[\Omega(n_1,n_2,n_3):= \langle n_{123}\rangle + \sum_{j=1}^3 (\pm_j) \langle n_j \rangle.\] For convenience we often relabel $n=:- n_0, \,  N=:N_0$   and rewrite the base tensor \eqref{base-tensor} in a symmetrized form:
\begin{equation}
\label{symmetrized-tensor}
h^b_{n_0n_1n_2n_3}= \prod_{j=0}^3 1_{N_j}(n_j) \cdot \mathbf{1}{\big\{n_0+n_{123}=0\big\}} \cdot \mathbf{1}{\big\{|\Omega-m|\leq 1\big\}},
\end{equation}
where here $\Omega$ has been rewritten as $\sum_{j=0}^3 (\pm_j) \langle n_j \rangle.$

\begin{lemma}[Base tensors estimates]\label{counting:basetensor_est}
We consider the base tensor $h^b$ defined as in (\ref{symmetrized-tensor}) with the dyadic numbers  $N_j$'s  ($j\in \{0,1,2,3\}$) and $m\in \mathbb{R}$. 
Then we have the following estimates:
\begin{enumerate}
\item Let $J\subseteq \{0,1,2,3\}$ satisfy $\#J=3$. Then we have
\begin{equation}\label{counting:base1}
\|h^b\|^2_{n_0n_1n_2n_3}\lesssim (\underset{j\in J}{\mathrm{med}}\, N_j)^{-1} \cdot \prod_{j\in J} N_j^3,
\end{equation} and
\begin{equation}\label{counting:base1'}
\|h^b\|^2_{n_0n_1n_2n_3}\lesssim (N_0N_1N_2N_3)^2\cdot \frac{\min{(N_0, N_1, N_2, N_3)}}{ \max{(N_0, N_1, N_2, N_3)}}.
\end{equation}
\smallskip
 
\item Let $J\subseteq \{j_1,j_2,j_3\}$ satisfy $\#J=2$, where $\{j_1,j_2,j_3,j_4\}= \{0,1,2,3\}$. Then we have
\begin{equation}\label{counting:base2}
\|h^b\|^2_{n_{j_1}n_{j_2}n_{j_3}\to n_{j_4}}\lesssim (\underset{j\in J\cup\{j_4\}}{\mathrm{med}}\, N_j)^{-1} \cdot \prod_{j\in J} N_j^3 \lesssim (\min_{j\in J} N_j)^{-1} \cdot \prod_{j\in J} N_j^3, 
\end{equation}
and
\begin{equation}\label{counting:base2'}
\|h^b\|^2_{n_{j_1}n_{j_2}n_{j_3}\to n_{j_4}}\lesssim (N_{j_1}N_{j_2}N_{j_3})^2\cdot (\max_{j\in \{j_1,j_2,j_3\}} N_j)^{-1}.
\end{equation}
\smallskip

\item Let us further localize the tensor to $|n_{j_1j_2}|=|n_{j_3j_4}|\sim N_{j_1j_2}=N_{j_3j_4}$. Then we have
\begin{equation}\label{counting:base3}
\|h^b\|^2_{n_{j_1}n_{j_2}\to n_{j_3}n_{j_4}}\lesssim \min (N_{j_1}, N_{j_1j_2})^{-1}\cdot \min(N_{j_3}, N_{j_3j_4})^{-1} \cdot(N_{j_1}N_{j_3})^3.
\end{equation}
\end{enumerate}
\end{lemma}
\smallskip

\begin{proof} 
First, (\ref{counting:base1}) and (\ref{counting:base1'})
follow from  (\ref{counting:eq6}) and (\ref{counting:eqsum2}) respectively. Then (\ref{counting:base2}) and (\ref{counting:base2'}) follow from (\ref{counting:eq5}) and (\ref{counting:eqsum1}) after applying Schur's test to the tensor norm.
Now consider (\ref{counting:base3}). Without loss of generality, we suppose that $(j_1,j_2,j_3,j_4)=(1,2,3,0)$ and $N_2=\max (N_1,N_2)\leq \max (N_3, N_0)=N_0$. Since $x\mapsto \min (x, y)^{-1}\cdot x^3$ is monotonically increasing in $x\geq 1$, it suffices to prove that 
\begin{equation*}
\|\mathbf{1}{\big\{|n_{12}|\sim N_{12}\big\}}h^b\|^2_{n_1n_2\to n_3n_0}\lesssim \min (N_1, N_{12})^{-1} \cdot \min (N_3, N_{03})^{-1}\cdot (N_1 N_3)^3.
\end{equation*}
To this end, using Schur's test we have that
\begin{multline*}
\|\mathbf{1}{\big\{|n_{12}|\sim N_{12}\big\}}h^b\|^2_{n_1n_2\to n_3n_0}\lesssim
\big(\sup_{n_0,n_3} \sum_{\substack{n_1,n_2\\n_{0123}=0}}\mathbf{1}{\big\{|n_{03}|\sim N_{03}\big\}}\cdot h^b_{n_0n_1n_2n_3} \big)\\\times \big(\sup_{n_1,n_2} \sum_{\substack{n_1,n_2\\n_{0123}=0}}\mathbf{1}{\big\{|n_{12}|\sim N_{12}\big\}}\cdot h^b_{n_0n_1n_2n_3} \big).
\end{multline*}
For the first factor, by  the lattice counting estimate (\ref{counting:eq1}) in  Lemma \ref{counting:lem1}, we have
\[\begin{aligned}
\sup_{n_0,n_3} \sum_{\substack{n_1,n_2\\n_{0123}=0}}\mathbf{1}{\big\{|n_{03}|\sim N_{03}\big\}}\cdot h^b_{n_0n_1n_2n_3} & \lesssim 
\sup_{m'\in \mathbb Z}\,\sup_{\substack{n_0,n_3\\|n_{03}|\sim N_{03}}} \sum_{|n_1|\sim N_1}\, \mathbf{1}{\big\{|\pm_1\langle n_1\rangle\pm_2\langle n_{013}\rangle-m'|\leq 1\big\}} \\
&\lesssim \min (N_{12}, N_1)^{-1} \cdot N_1^3.
\end{aligned}\]
Similarly, the second factor can be bounded by $\min (N_{03}, N_3)^{-1} \cdot N_3^3$, which together with the bound of the first factor yields (\ref{counting:base3}).
\end{proof}
\begin{remark} 
In Lemma \ref{counting:basetensor_est} we can minimize the right hand sides of (\ref{counting:base2}) and (\ref{counting:base1}) by making $J$  contain the two indices $j$'s in (\ref{counting:base2}) and respectively the three indices $j$'s in (\ref{counting:base1}) such that $N_j$'s are the smallest, in which case (\ref{counting:base1}) and (\ref{counting:base2}) are equivalent to the counting estimates (\ref{counting:eq6}) and (\ref{counting:eq5}) in Lemma \ref{counting:lem2}.
\end{remark}

\subsection{The cubic tensor}\label{cubic-tensor}
In this subsection, we prove estimates for cubic tensors, which correspond to stochastic objects such as $\scubic$ (see Subsection \ref{section:diagram-wave}). In Lemma \ref{counting:cubic_tensor} below, the cubic tensors $h$ and $H$ correspond to the cubic nonlinearity and its Duhamel integral, respectively.

\begin{lemma}[The cubic tensor estimates]\label{counting:cubic_tensor}
Suppose that $N_1,N_2,N_3,N_{123}$ are dyadic numbers and $\lambda_1, \lambda_2, \lambda_3$ are real numbers. Let $N_{\max} = \max (N_1,N_2,N_3)$, $N_{\min}=\min (N_{123}, N_1, N_2, N_3)$ and $\chi(t)$ be a Schwartz function.  Define the tensors $h_{nn_1n_2n_3}(t, \lambda_1, \lambda_2,\lambda_3)$ and $H_{nn_1n_2n_3}(t, \lambda_1, \lambda_2,\lambda_3)$ by 
\begin{align}
\label{def:cubic_tensor1}
h_{nn_1n_2n_3}(t,\lambda_1,\lambda_2,\lambda_3) 
&:= \mathbf{1}{\big\{ n=n_{123} \big\}}\cdot {1_{N_{123}}(n)} \Big( \prod_{j=1}^3 \frac{{1_{N_j}(n_j)}}{\langle n_j\rangle} \Big)
 \cdot\chi(t) e^{i(\pm\langle n_1\rangle\pm\langle n_2\rangle\pm\langle n_3\rangle+\lambda_1+\lambda_2+\lambda_3)t},\\
\label{def:cubic_tensor2}
H_{nn_1n_2n_3}(t,\lambda_1,\lambda_2,\lambda_3) 
&:= \mathbf{1}{\big\{ n=n_{123} \big\}}\cdot \frac{1_{N_{123}}(n)}{\langle n\rangle} \Big( \prod_{j=1}^3 \frac{1_{N_j}(n_j)}{\langle n_j\rangle} \Big)\\&\times
\int_0^t \dt'\,\chi(t)\chi(t')\cdot \sin((t-t^\prime)\langle n \rangle) e^{i(\pm\langle n_1\rangle\pm\langle n_2\rangle\pm\langle n_3\rangle+\lambda_1+\lambda_2+\lambda_3)t'}.
\end{align}
Then there exist two functions $A_j=A_j(\lambda,\lambda_1,\lambda_2,\lambda_3)$, $j= 1, 2$, such that 
\begin{equation}\label{counting:cubic_coef}
\| \langle \lambda \rangle^{b_{+}-1}  A_1(\lambda,\lambda_1,\lambda_2,\lambda_3)\|_{L^2_{\lambda}}\lesssim N_{\max}^{\epsilon},\qquad
\| \langle \lambda \rangle^{b_+}  A_2(\lambda,\lambda_1,\lambda_2,\lambda_3)\|_{L^2_{\lambda}}\lesssim N_{\max}^{\epsilon},
\end{equation} and that we have the following bounds for the shifted space time Fourier transform of \eqref{def:cubic_tensor1} and \eqref{def:cubic_tensor2}, 
\begin{align}
\label{counting:cubic_bound1}
\|\widetilde{h}_{nn_1n_2n_3}(\lambda, \lambda_1,\lambda_2,\lambda_3)\|_{nn_1n_2n_3}&\lesssim A_1(\lambda,\lambda_1,\lambda_2,\lambda_3)N_{123}\cdot \big(\frac{N_{\min}}{N_{\max}}\big)^{\frac{1}{2}},
\\\label{counting:cubic_bound2}
\|\widetilde{H}_{nn_1n_2n_3}(\lambda, \lambda_1,\lambda_2,\lambda_3)\|_{nn_1n_2n_3}&\lesssim
A_2(\lambda,\lambda_1,\lambda_2,\lambda_3)\cdot \big(\frac{N_{\min}}{N_{\max}}\big)^{\frac{1}{2}},
\\\label{counting:cubic_bound2mid}
\|\widetilde{h}_{nn_1n_2n_3}(\lambda, \lambda_1,\lambda_2,\lambda_3)\|_{nn_B\to n_C}&\lesssim
A_1(\lambda,\lambda_1,\lambda_2,\lambda_3)N_{123}\cdot N_{\max}^{-\frac{1}{2}},
\\\label{counting:cubic_bound2'}
\|\widetilde{H}_{nn_1n_2n_3}(\lambda, \lambda_1,\lambda_2,\lambda_3)\|_{nn_B\to n_C}&\lesssim
A_2(\lambda,\lambda_1,\lambda_2,\lambda_3)\cdot N_{\max}^{-\frac{1}{2}},
\end{align} for all partitions $(B,C)$ of $\{1,2,3\}$ where $C\neq \varnothing$. Furthermore, the $\partial_\lambda$ and $\partial_{\lambda_j}$ derivatives of $\widetilde{h}$ and $\widetilde{H}$ satisfy the same estimates as $\widetilde{h}$ and $\widetilde{H}$ themselves.
\end{lemma}
\begin{proof}
First, by the definition (\ref{def:cubic_tensor1}) of $h$, and definition \eqref{tilted} we have that
\begin{equation}\label{eq:htilde}
\begin{aligned}
\widetilde{h}_{nn_1n_2n_3}(\lambda,\lambda_1,\lambda_2,\lambda_3) 
&= \mathbf{1}{\big\{ n=n_{123} \big\}}\cdot {1_{N_{123}}(n)} \cdot \widehat{\chi}(\lambda-\lambda_1-\lambda_2-\lambda_3-\Omega) \cdot \Big( \prod_{j=1}^3 {1_{N_j}(n_j)} \Big),
\end{aligned}
\end{equation}
where $\Omega =\pm \langle n\rangle\pm \langle n_1\rangle\pm\langle n_2\rangle\pm \langle n_3\rangle$.
Then by the definition (\ref{base-tensor}) of the base tensor and the Minkowski's inequality, we have
\begin{equation}
\| \widetilde{h}_{nn_1n_2n_3}(\lambda,\lambda_1,\lambda_2,\lambda_3)\|_{nn_1n_2n_3} \lesssim A_1(\lambda,\lambda_1,\lambda_2,\lambda_3) (N_1N_2N_3)^{-1}\cdot \sup_{\substack{m\in \mathbb{Z}\\|m|\lesssim N_{\max}}} \|h^b\|_{nn_1n_2n_3},
\end{equation}
where 
\begin{equation}\label{cubic:A1}
A_1(\lambda,\lambda_1,\lambda_2,\lambda_3):= \sum_{\substack{m\in \mathbb{Z}\\|m|\lesssim N_{\max}}} |\widehat{\chi}(\lambda-\lambda_1-\lambda_2-\lambda_3-m)|;
\end{equation} note that here we have replaced $\langle n_j\rangle^{-1}1_{N_j}(n_j)$ by $N_j^{-1}$, a simplification which will be also used multiple times below. Define $\Lambda:=\lambda_1+\lambda_2+\lambda_3$, then (\ref{cubic:A1}) implies that
\begin{equation}\label{cubic:A1-2}|A_1|\lesssim \min(1,N_{\mathrm{max}}\cdot\langle\lambda-\Lambda\rangle^{-1}),
\end{equation}
whence
\begin{equation}\label{cubic:A1-3}
\| \langle \lambda \rangle^{b_{+}-1}  A_1(\lambda,\lambda_1,\lambda_2,\lambda_3)\|_{L^2_{\lambda}}\lesssim N_{\mathrm{max}}^{\frac{3}{2}(b^+-\frac{1}{2})}\cdot\| \langle \lambda \rangle^{b_{+}-1}\langle \lambda-\Lambda\rangle^{-\frac{3}{2}(b^+-\frac{1}{2})}\|_{L_\lambda^2}\lesssim N_{\mathrm{max}}^\epsilon,
\end{equation}which is the first inequality in (\ref{counting:cubic_coef}). Now, without loss of the generality, we may assume $N_1\geq N_2\geq N_3$. By applying the base tensor estimate (\ref{counting:base1}) in Lemma \ref{counting:basetensor_est}, we have \begin{equation*}
 \|h^b\|^2_{nn_1n_2n_3} \lesssim {\mathrm{med}}(N_2, N_3, N_{123})^{-1} \cdot ( N_2N_3N_{123})^3\lesssim N_2^3 N_3^2 N_{123}^3
\leq N_{123}\cdot (N_{123}N_1N_2N_3)^2 \cdot N_{\textup{max}}^{-1}.
\end{equation*}
By the symmetry of $(n,n_1,n_2,n_3)$ in the base tensor $h^b$, the above inequality implies
\[
 \|h^b\|^2_{nn_1n_2n_3} \lesssim N_{\min}\cdot (N_{123}N_1N_2N_3)^2 \cdot N_{\textup{max}}^{-1}.
\]
which proves (\ref{counting:cubic_bound1}).

Second, by the definition (\ref{def:cubic_tensor2}) of $H$, definition \eqref{tilted} and  Lemma \ref{lin} we can write
\begin{equation*}
\begin{aligned}
\widetilde{H}_{nn_1n_2n_3}(\lambda,\lambda_1,\lambda_2,\lambda_3) 
&= K(\lambda, \lambda_1+\lambda_2+\lambda_3+\Omega)\cdot \mathbf{1}{\big\{ n=n_{123} \big\}}\cdot \frac{1_{N_{123}}(n)}{\langle n\rangle}   \cdot \Big( \prod_{j=1}^3 \frac{1_{N_j}(n_j)}{\langle n_j\rangle}\Big),
\end{aligned}
\end{equation*}
where $\Omega =\pm \langle n\rangle\pm \langle n_1\rangle\pm \langle n_2\rangle\pm \langle n_3\rangle$ and $K(\lambda,\lambda_1+\lambda_2+\lambda_3+\Omega)$ is as defined in Lemma \ref{lin} and Remark \ref{rmk:notation}. 
Based on the above expression of $\widetilde{H}_{nn_1n_2n_3}(\lambda,\lambda_1,\lambda_2,\lambda_3)$ and (\ref{esti}), we can do{\footnote{This type of level-set decomposition of $\Omega$ will be frequently used later.}} a \emph{level-set decompositon} of $\Omega$ to obtain that, 
\begin{multline*}
\| \widetilde{H}\|_{nn_1n_2n_3} \lesssim  \sum_{\substack{m\in \mathbb{Z}\\|m|\lesssim N_{\max}}} \frac{1}{\langle \lambda\rangle\langle\lambda\pm (\lambda_1+\lambda_2+\lambda_3+m)\rangle} \\\times\Big\|  \mathbf{1}{\big\{n=n_{123}\big\}}\cdot \mathbf{1}{\big\{|\Omega -m|\leq 1\big\}} \cdot\frac{1_{N_{123}}(n)}{\langle n\rangle}   \cdot \Big( \prod_{j=1}^3 \frac{1_{N_j}(n_j)}{\langle n_j\rangle} \Big)\Big\|_{nn_1n_2n_3}.
\end{multline*}
Defining
\[A_2(\lambda,\lambda_1,\lambda_2,\lambda_3):= \sum_{\substack{m\in \mathbb{Z}\\|m|\lesssim N_{\max}}} \frac{1 }{\langle \lambda\rangle\langle\lambda\pm (\lambda_1+\lambda_2+\lambda_3+m)\rangle},\]
we then have that
\begin{equation}\label{counting:cubic_bound2''}
\| \widetilde{H}\|_{nn_1n_2n_3}  \lesssim \,\, |A_2(\lambda,\lambda_1,\lambda_2,\lambda_3)|\cdot \sup_{\substack{m\in \mathbb{Z}\\|m|\lesssim N_{\max}}} \|h^b\|_{nn_1n_2n_3} \cdot (N_{123}N_1N_2N_3)^{-1}
\end{equation}
which proves (\ref{counting:cubic_bound2}) by the base tensor estimate (\ref{counting:base1}) in Lemma \ref{counting:basetensor_est}. Furthermore, it is elementary to show that
\begin{equation}\label{eq:Abound}
A_2(\lambda,\lambda_1,\lambda_2,\lambda_3)\lesssim 
\log (N_{\max}) \cdot \langle\lambda\rangle^{-1}\cdot \min \big(1, N_{\max} \langle\lambda\pm\Lambda\rangle^{-1}\big)
\end{equation} (where again $\Lambda:=\lambda_1+\lambda_2+\lambda_3$), which is the same as (\ref{cubic:A1-2}) for $A_1$ but with an extra $\langle \lambda\rangle^{-1}$ factor on the right hand side and at most a $\log (N_{\mathrm{max}})$ loss. Then we can prove the second inequality in (\ref{counting:cubic_coef}) in the same way as the first one, by repeating (\ref{cubic:A1-3}).

Finally, the proofs of (\ref{counting:cubic_bound2mid}) and (\ref{counting:cubic_bound2'}) follow in the same manner, where one uses  (\ref{counting:base2})--(\ref{counting:base3}) in Lemma \ref{counting:basetensor_est}, instead of  (\ref{counting:base1}). The $\partial_\lambda$ and $\partial_{\lambda_j}$ estimates also follow in the same way because all derivatives of $\widehat{\chi}$ and $K$ satisfy the same bounds as these functions themselves.
\end{proof}

\begin{corollary}\label{counting:cor-cubic}
Suppose $\varphi_j \in \{\sin, \cos\}$ for $j=1,\cdots, 3$ and $N_1,\cdots,N_3,N_{123}$ are dyadic numbers, and let $N_{\textup{max}} =\max (N_1, N_2, N_3)$. Consider the cubic tensor $H_{nn_1n_2n_3}(t)$, which arises from the cubic stochastic object $\scubic$ and is defined as
\begin{equation}\label{def:cubic_tensor3}
\begin{aligned}
H_{nn_1n_2n_3}(t) =& H_{nn_1n_2n_3}[N_1,N_2,N_3,N_{123},\varphi_1,\varphi_2,\varphi_3](t) \\
=& \mathbf{1}\{ n=n_{123} \}\cdot \frac{1_{N_{123}}(n)}{\langle n\rangle}\cdot \Big( \prod_{j=1}^3 \frac{1_{N_j}(n_j)}{\langle n_j\rangle} \Big) \\
&\times\Big(\int\displaylimits_{0}^{t} \dt^\prime\, \chi(t)\chi(t')\cdot \sin((t-t^\prime)\langle n \rangle) \prod_{j=1}^3 \varphi_j(t^\prime \langle n_j \rangle) \Big).
\end{aligned}
\end{equation}
Then we have that, 
\begin{align}\label{counting:cubic_bound3}
 \sum_{ \substack{\varphi_1,\varphi_2,\varphi_3 \\ \in \{ \cos,\sin\} }} \| \langle\lambda \rangle^{b_+}\widetilde{H}_{nn_1n_2n_3}(\lambda)\|_{L^2_\lambda[nn_1n_2n_3]}&\lesssim N_{123}^{1/2}\cdot N_{\textup{max}}^{-1/2+\epsilon},\\ 
 \label{counting:cubic_bound4}
 \sum_{ \substack{\varphi_1,\varphi_2,\varphi_3 \\ \in \{ \cos,\sin\} }} \| \langle\lambda \rangle^{b_+}\widetilde{H}_{nn_1n_2n_3}(\lambda)\|_{L^2_\lambda[nn_B\to n_C]}&\lesssim N_{\mathrm{max}}^{-1/2+\epsilon}.
\end{align}
for any partition $(B,C)$ of $\{1,2,3\}$ with $C\neq\varnothing$. Here the $L^2_\lambda[nn_B\to n_C]$ norm is defined in \eqref{def:L2lambda}.
\end{corollary}
\begin{proof}
By the definition as in (\ref{def:cubic_tensor3}), we see that $H_{nn_1n_2n_3}(t)$ can be written as a linear combination of $H_{nn_1n_2n_3}(t,\lambda_1,\lambda_2,\lambda_3)$ as in (\ref{def:cubic_tensor2}), by setting $\lambda_1=\lambda_2=\lambda_3=0$ and choosing the $\pm$ signs. The desired bounds (\ref{counting:cubic_bound3}) and (\ref{counting:cubic_bound4}) then directly follow from (\ref{counting:cubic_coef}), (\ref{counting:cubic_bound2}) and (\ref{counting:cubic_bound2'}).
\end{proof}
\subsection{The quintic tensor}\label{quintic-tensor}
In this subsection, we will study the tensor estimates for the quintic stochastic objects such as $\squintic$ (see details in Section \ref{section:diagrams}). Similar as in Subsection \ref{cubic-tensor}, the tensors $h$ and $H$ in Lemma \ref{counting:quintic_tensor} correspond to the quintic nonlinearity and its Duhamel integral, respectively.
\begin{lemma}[The quintic tensor estimates]\label{counting:quintic_tensor}
Let $N_0,\cdots,N_5,N_{234}$ be dyadic numbers and $\lambda, \lambda_1,\cdots,\lambda_5$ be real numbers. Let $N_{\max} = \max (N_1,\cdots,N_5)$. Define the tensors
 \begin{equation}
\label{def:quintic_tensor1}
\begin{aligned}
h_{n_0 n_1 \hdots n_5}(t,\lambda_1,\cdots, \lambda_5) &= \mathbf{1}{\big\{ n_0 = n_{12345} \big\}} 1_{N_0}(n_0) \frac{1_{N_{234} }(n_{234})}{\langle n_{234} \rangle} \Big( \prod_{j=1}^5 \frac{1_{N_j}(n_j)}{\langle n_j \rangle} \Big) 
 \cdot e^{it(\pm\langle n_1\rangle\pm\langle n_5\rangle+\lambda_1+\lambda_5)}\\
&\times  \Big( \int\displaylimits_0^{t} \dt^{\prime}\, \chi(t) \chi(t') \sin\big((t-t^{\prime}) \langle n_{234} \rangle\big) e^{it'(\pm\langle n_2\rangle\pm\langle n_3\rangle\pm\langle n_4\rangle+\lambda_2+\lambda_3+\lambda_4)} \Big),
\end{aligned}
\end{equation}
\begin{equation}\label{def:quintic_tensor2}\begin{aligned}
H_{n_0 n_1 \hdots n_5}(t,\lambda_1,\cdots, \lambda_5) &=\, \mathbf{1}{\big\{ n_0 = n_{12345} \big\}} \cdot \frac{1_{N_{234}}(n_{234})}{\langle n_{234} \rangle} \cdot \Big( \prod_{j=0}^5 \frac{1_{N_j}(n_j)}{\langle n_j \rangle} \Big) 
\\
& \times \bigg( \int\displaylimits_0^{t} \dt^{\prime}\, \chi(t) \chi(t') \sin\big((t-t') \langle n_{0} \rangle\big)  e^{it'(\pm\langle n_1\rangle\pm\langle n_5\rangle+\lambda_1+\lambda_5)}  \\
&\times   \Big( \int\displaylimits_0^{t'} \dt^{\prime\prime}\, \chi(t') \chi(t'') \sin\big((t'-t'') \langle n_{234} \rangle\big)  e^{it''(\pm\langle n_2\rangle\pm\langle n_3\rangle\pm\langle n_4\rangle+\lambda_2+\lambda_3+\lambda_4)} \Big)\bigg).
\end{aligned}
\end{equation} Then there exist two functions $B_j=B_j(\lambda,\lambda_1,\cdots,\lambda_5)$, \,$ j= 1,2 $, such that
\begin{equation}
\label{counting:quintic_coef}
\| \langle \lambda \rangle^{b_+-1}  B_1(\lambda,\lambda_1,\cdots,\lambda_5)\|_{L^2_{\lambda}}\lesssim N_{\max}^{\epsilon},\qquad \| \langle \lambda \rangle^{b_+}  B_2(\lambda,\lambda_1,\cdots,\lambda_5)\|_{L^2_{\lambda}}\lesssim N_{\max}^{\epsilon},
\end{equation} and that we have
\begin{align}
\label{counting:quintic_bound1}
\left\|\widetilde{h}(\lambda, \lambda_1,\cdots,\lambda_5)\right\|_{n_0n_1\cdots n_5}&\lesssim  \frac{N_0\cdot\min{(N_2,N_3,N_4,N_{234})}^{\frac{1}{2}}}{\max(N_0,N_1,N_5)^{\frac{1}{2}}\max(N_2,N_3,N_4)^{\frac{1}{2}}} B_1(\lambda,\lambda_1,\cdots,\lambda_5)\\\nonumber
&\lesssim N_0 \cdot  N_{\max}^{-\frac{1}{2}}\cdot  B_1(\lambda,\lambda_1,\cdots,\lambda_5),\\
\label{counting:extra_HS1}
\left\|\widetilde{h}(\lambda, \lambda_1,\cdots,\lambda_5)\right\|_{n_0n_1\cdots n_5}&\lesssim  \frac{N_0 \min (N_0,N_1,N_2,N_5)^{\frac{1}{2}} (N_3N_4)^{\frac{1}{2}}}{N_2\max (N_0,N_1,N_2,N_5)^{\frac{1}{2}}} \cdot  B_1(\lambda,\lambda_1,\cdots,\lambda_5)\quad \text{when } N_2\sim N_{234},\\
\label{counting:quintic_bound2}
\left\|\widetilde{H}(\lambda, \lambda_1,\cdots,\lambda_5)\right\|_{n_0n_1\cdots n_5}&\lesssim   N_{\max}^{-\frac{1}{2}}\cdot B_2(\lambda,\lambda_1,\cdots,\lambda_5),\\
\nonumber
\big\|\widetilde{h}(\lambda, \lambda_1,\cdots,\lambda_5)\big\|_{n_0n_A\to n_Bn_5}
&\lesssim N_0^{\frac{1}{2}}{N_5}^{-\frac{1}{2}}\cdot \big\{\max (N_0, N_2, N_3, N_4)^{-\frac{1}{2}} +\max (N_2, N_3, N_4, N_5)^{-\frac{1}{2}}\big\}\\
\label{counting:quintic_bound1'}
&\times B_1(\lambda,\lambda_1,\cdots,\lambda_5),
\\
\nonumber
\big\|\widetilde{H}(\lambda, \lambda_1,\cdots,\lambda_5)\big\|_{n_0n_A\to n_Bn_5}
&\lesssim (N_0N_5)^{-\frac{1}{2}}\cdot \big\{\max (N_0, N_2, N_3, N_4)^{-\frac{1}{2}} +\max (N_2, N_3, N_4, N_5)^{-\frac{1}{2}}\big\}\\
\label{counting:quintic_bound2'}
&\times B_2(\lambda,\lambda_1,\cdots,\lambda_5),
\end{align} for any partition $(A,B)$ of $\{1,2,3,4\}$. The same bounds hold for all  $\partial_\lambda$ and $\partial_{\lambda_j}$ derivatives of $\widetilde{h}$ and $\widetilde{H}$.
\end{lemma}
\begin{proof}
First consider (\ref{counting:quintic_bound1}). We may use $\sin z=\frac{e^{iz}-e^{-iz}}{2i}$ to rewrite $\sin((t-t^{\prime}) \langle n_{234} \rangle)$. Then by definition \eqref{tilted}, Lemma \ref{lin} and \eqref{def:quintic_tensor1}, we have up to a linear combination that 
\begin{multline}
\widetilde{h}_{n_0 n_1 \hdots n_5}(\lambda,\lambda_1,\cdots, \lambda_5) =\, \mathbf{1}{\big\{ n_0 = n_{12345} \big\}} \cdot 1_{N_0}(n_0)\cdot \frac{1_{N_{234} }(n_{234})}{\langle n_{234} \rangle} \cdot\Big( \prod_{j=1}^5 \frac{1_{N_j}(n_j)}{\langle n_j \rangle} \Big) \\
\times K(\lambda-\lambda_1-\lambda_5 -\Omega', \lambda_2+\lambda_3+\lambda_4+\Omega'').
\end{multline}
Similarly, by expanding also $\sin((t-t^{\prime}) \langle n_{0} \rangle)$, we have up to a linear combination that
\begin{multline}
\widetilde{H}_{n_0 n_1 \hdots n_5}(\lambda,\lambda_1,\cdots, \lambda_5) =\, \mathbf{1}{\big\{ n_0 = n_{12345} \big\}} \cdot \frac{1_{N_{234} }(n_{234})}{\langle n_{234} \rangle} \cdot \Big( \prod_{j=0}^5 \frac{1_{N_j}(n_j)}{\langle n_j \rangle} \Big) 
 \\\times \int_{\mathbb{R}}\ K(\lambda, \sigma)
 \cdot K(\sigma-\lambda_1-\lambda_5 -\Omega', \lambda_2+\lambda_3+\lambda_4+\Omega'')\,\mathrm{d}\sigma,
\end{multline}
where in the above formulas \[
\Omega' = \pm\langle n_0\rangle\pm\langle n_1\rangle\pm\langle n_5\rangle\pm_{234}\langle n_{234}\rangle,\quad\Omega^{\prime\prime} = \mp_{234}\langle n_{234}\rangle\pm\langle n_2\rangle \pm\langle n_3\rangle\pm\langle n_4\rangle.
\]
By a level set decomposition and using Lemma \ref{lin} (and replacing $\langle n_j\rangle^{-1}\mathbf{1}_{N_j}(n_j)$ by $N_j^{-1}$ as above), we have  that
\begin{equation}\label{quintic_eq1}
\begin{aligned}
\|\widetilde{h}(\lambda,\lambda_1,\cdots, \lambda_5)\|_{n_0n_1\cdots n_5} \lesssim&\, 
\sum_{\substack{m_1,m_2\in \mathbb{Z}\\
|m_1|,|m_2|\lesssim N_{\max}}} \frac{1 }{\langle \lambda-\lambda_1-\lambda_5- m_1 \rangle \langle \lambda-\sum_{j=1}^5 \lambda_j -m_2\rangle}\\
&\times (\prod_{j=1}^5 N_j)^{-1}\cdot N_{234}^{-1}\cdot
 \bigg\|\sum_{n_{234}} h^b_{n_0n_1n_{234}n_5} \cdot h^b_{n_{234}n_2n_3n_4}\bigg\|_{n_0\cdots n_5},
\end{aligned}
\end{equation}
where $h^b_{n_0n_1n_{234}n_5}$ and $h^b_{n_{234}n_2n_3n_4}$ are two base tensors as defined in (\ref{base-tensor}) with the conditions $|m_1-\Omega'|\leq 1$ and $|m_2-m_1-\Omega''|\leq 1$ respectively. 
We define 
\[B_1(\lambda,\lambda_1,\cdots,\lambda_5):= \sum_{\substack{m_1,m_2\in \mathbb{Z}\\
|m_1|,|m_2|\lesssim N_{\max}}} \frac{1 }{\langle \lambda-\lambda_1-\lambda_5- m_1 \rangle \langle \lambda-\sum_{j=1}^5 \lambda_j -m_2\rangle},\]
and then
\begin{equation}\label{quintic_extra_eq1}
\| \widetilde{h}\|_{n_0n_1\cdots n_5}  \lesssim \,\, |B_1(\lambda,\lambda_1,\cdots,\lambda_5)|\cdot (\prod_{j=1}^5 N_j)^{-1}\cdot N_{234}^{-1}\cdot\sup_{m_1,m_2}
 \bigg\|\sum_{n_{234}} h^b_{n_0n_1n_{234}n_5} \cdot h^b_{n_{234}n_2n_3n_4}\bigg\|_{n_0\cdots n_5}.
\end{equation} By Lemma \ref{counting:basetensor_est} and the merging estimate (Lemma \ref{counting:lem-merging}) we have that
\begin{equation}\label{quintic_eq3}
\begin{aligned}
&\,\,\,\,\,\,\,\,(\prod_{j=1}^5 N_j)^{-1}\cdot N_{234}^{-1}\cdot\bigg\|\sum_{n_{234}} h^b_{n_0n_1n_{234}n_5} \cdot h^b_{n_{234}n_2n_3n_4}\bigg\|_{n_0\cdots n_5}\\
&\lesssim (\prod_{j=1}^5 N_j)^{-1}\cdot N_{234}^{-1}\cdot \|h^b_{n_0n_1n_{234}n_5}\|_{n_0n_1n_5\to n_{234}}\cdot \|h^b_{n_{234}n_2n_3n_4}\|_{n_2n_3n_4 n_{234}}\\
&\lesssim \frac{ N_0 \cdot \min (N_2,N_3,N_4,N_{234})^{1/2}}{\max (N_0,N_1,N_5)^{1/2}\cdot\max(N_2,N_3,N_4)^{1/2}}\lesssim N_0\cdot\max(N_0,\cdots,N_5)^{-1/2},
\end{aligned}
\end{equation} which proves (\ref{counting:quintic_bound1}). Furthermore, it is elementary to show that
\[|B_1|\lesssim (\log N_{\mathrm{max}})^2\cdot\min(1,N_{\mathrm{max}}\langle \lambda-\Lambda_1\rangle^{-1}),\] where $\Lambda_1:=\lambda_1+\lambda_5$, which is again similar to (\ref{cubic:A1-2}), so by repeating (\ref{cubic:A1-3}) we obtain the first inequality in (\ref{counting:quintic_coef}).

The proof of \eqref{counting:extra_HS1} is similar to the above proof of \eqref{counting:quintic_bound1}.   By $N_2\sim N_{234}$ and \eqref{counting:base1'} in Lemma \ref{counting:basetensor_est}, we get that $\|h^b_{n_0n_1n_{234}n_5}\|_{n_0n_1n_{234}n_5}\lesssim(N_0N_1N_2N_5)\cdot \min (N_0,N_1,N_2,N_5)^{1/2}\cdot\max (N_0,N_1,N_2,N_5)^{-1/2}$. By Schur's test and using the fact that $(n_2,n_3,n_4)$ has at most $(N_3N_4)^3$ choices in the support of $h_{n_{234}n_2n_3n_4}^b$ once $n_{234}$ is fixed, we obtain that $\|h^b_{n_{234}n_2n_3n_4}\|_{n_2n_3n_4\to n_{234}}\lesssim(N_3N_4)^{3/2}$.
Therefore, the estimate \eqref{quintic_extra_eq1} and the first inequality in \eqref{quintic_eq3} prove \eqref{counting:extra_HS1}.

Next, similarly as with \eqref{quintic_eq1}, for the tensor $\widetilde{H}(\lambda,\lambda_1,\cdots, \lambda_5)$, we have that 
\begin{equation}\label{quintic_eq2}
\begin{aligned}
\|\widetilde{H}(\lambda,\lambda_1,\cdots, \lambda_5)\|_{n_0n_1\cdots n_5} &\lesssim
|B_2(\lambda,\lambda_1,\cdots,\lambda_5)| \cdot (\prod_{j=0}^5 N_j)^{-1} N_{234}^{-1}\\&\times
 \bigg\|\sum_{n_{234}} h^b_{n_0n_1n_{234}n_5} \cdot h^b_{n_{234}n_2n_3n_4}\bigg\|_{n_0\cdots n_5}
\end{aligned}
\end{equation}
where $h^b_{n_0n_1n_{234}n_5}$ and $h^b_{n_{234}n_2n_3n_4}$ are the same as in \eqref{quintic_eq1}, and 
\[B_2(\lambda,\lambda_1,\cdots,\lambda_5):= \sum_{\substack{m_1,m_2\in \mathbb{Z}\\
|m_1|,|m_2|\lesssim N_{\max}}} \int_\mathbb{R}\frac{\mathrm{d}\sigma }{\langle\lambda\rangle\langle\lambda-\sigma\rangle\langle \sigma-\lambda_1-\lambda_5- m_1 \rangle \langle \sigma-\sum_{j=1}^5 \lambda_j -m_2\rangle}.\] The same arguments above then prove \eqref{counting:quintic_bound2}, and it is elementary to show that
\[|B_2|\lesssim (\log N_{\mathrm{max}})^3 \langle\lambda\rangle^{-1}\cdot\min(1,N_{\mathrm{max}}\langle \lambda-\Lambda_1\rangle^{-1}),\] where again $\Lambda_1:=\lambda_1+\lambda_5$, which then implies the second inequality in (\ref{counting:quintic_coef}).

Now we turn to (\ref{counting:quintic_bound1'}) and (\ref{counting:quintic_bound2'}). By the same argument above, it suffices to prove that
\begin{equation}\label{counting:quintic_eq3}
\begin{aligned}\bigg\|\sum_{n_{234}} h^b_{n_0n_1n_{234}n_5} \cdot h^b_{n_{234}n_2n_3n_4}\bigg\|_{n_0n_A\to n_Bn_5}&\lesssim \prod_{j=0}^5N_j\cdot N_{234}(N_0N_5)^{-\frac{1}{2}}\\&\times \big\{\max (N_0, N_2, N_3, N_4)^{-\frac{1}{2}} +\max (N_2, N_3, N_4, N_5)^{-\frac{1}{2}}\big\}.
\end{aligned}\end{equation} Here we shall apply the merging estimate (Lemma \ref{counting:lem-merging}). By choosing the order of $(h^{(1)}, h^{(2)})$ in Lemma \ref{counting:lem-merging}, we are able to get 
\begin{multline}\label{counting:quintic_eq4}\textrm{LHS\ of\ }(\ref{counting:quintic_eq3})\lesssim\min\big(\|h_{n_0n_1n_{234}n_5}^b\|_{n_X\to n_Y}\cdot\|h_{n_{234}n_2n_3n_4}\|_{n_Z\to n_W},\\\|h_{n_0n_1n_{234}n_5}^b\|_{n_{X'}\to n_{Y'}}\cdot\|h_{n_{234}n_2n_3n_4}\|_{n_{Z'}\to n_{W'}}\big)\end{multline} where $0\in X\cap X'$, $5\in Y\cap Y'$, and $Z\in\{\varnothing,\{2\},\{3\},\{4\}\}$ and $Z'=Z\cup\{234\}$. By Lemma \ref{counting:basetensor_est} we have
\begin{equation}\label{counting:quintic_eq5}\|h_{n_0n_1n_{234}n_5}^b\|_{n_X\to n_Y}\lesssim N_0N_1N_{234}N_5\cdot (N_0N_5)^{-\frac{1}{2}},\end{equation} and moreover we have
\begin{equation}\label{counting:quintic_eq6}\|h_{n_0n_1n_{234}n_5}^b\|_{n_X\to n_Y}\lesssim N_0N_1N_{234}N_5\cdot (N_0N_5)^{-\frac{1}{2}}\big(N_0^{-\frac{1}{2}}+N_5^{-\frac{1}{2}}\big)\end{equation} if $\min(N_0,N_5)\gg \max(N_2,N_3,N_4)$; the same holds for $(X',Y')$. In fact, if $|X|=1$ (or $|Y|=1$) then the stronger bound (\ref{counting:quintic_eq6}) directly follows from Lemma \ref{counting:basetensor_est}; if $|X|=|Y|=2$, then Lemma \ref{counting:basetensor_est} already implies (\ref{counting:quintic_eq5}), and when $\min(N_0,N_5)\gg N_{234}$ we have the improved estimate (\ref{counting:quintic_eq6}), because we are in the case $A\gtrsim N$ when applying Lemma \ref{counting:lem1} after using Schur's bound, which then gains an extra factor of either $N_0^{-1}$ or $N_5^{-1}$.

Next consider the relevant norms of $h_{n_{234}n_2n_3n_4}$. If $Z=\varnothing$ then $Z'=\{234\}$, and hence
\begin{equation}\label{counting:quintic_eq7}\|h_{n_{234}n_2n_3n_4}\|_{n_{Z'}\to n_{W'}}\lesssim N_2N_3N_4\cdot \max(N_2,N_3,N_4)^{-\frac{1}{2}}.\end{equation} If not, say $Z=\{2\}$, then we have
\begin{equation}\label{counting:quintic_eq8}\|h_{n_{234}n_2n_3n_4}\|_{n_Z\to n_W}\lesssim N_{234}N_3N_4\cdot \max(N_{234},N_3,N_4)^{-\frac{1}{2}}\lesssim N_2N_3N_4\cdot \max(N_2,N_3,N_4)^{-\frac{1}{2}}\end{equation} if $N_{234}\lesssim N_2$, and
\begin{equation}\label{counting:quintic_eq9}\|h_{n_{234}n_2n_3n_4}\|_{n_{Z'}\to n_{W'}}\lesssim N_2N_3N_4\cdot \max(N_3,N_4)^{-\frac{1}{2}}\sim N_2N_3N_4\cdot \max(N_2,N_3,N_4)^{-\frac{1}{2}}\end{equation}  if $N_{234}\gg N_2$ (this last inequality (\ref{counting:quintic_eq9}) is due to the same reason as above when applying Lemma \ref{counting:lem1}). Now putting together the above bounds by choosing between $(Z,W)$ and $(Z',W')$ appropriately, we can get that
\[\textrm{LHS\ of\ }(\ref{counting:quintic_eq3})\lesssim (N_0\cdots N_5)\cdot N_{234}\cdot (N_0N_5)^{-\frac{1}{2}}\cdot\big(\max\big(\max(N_2,N_3,N_4),\min(N_0,N_5)\big)\big)^{-\frac{1}{2}},\] which then proves (\ref{counting:quintic_eq3}). The $\partial_{\lambda}$ and $\partial_{\lambda_j}$ derivative estimates follow just like in Lemma \ref{counting:cubic_tensor}.
\end{proof}
\begin{corollary}\label{counting:cor-quintic}
Suppose $\varphi_j \in \{\sin, \cos\}$ for $j=1,\cdots, 5$ and $N_0,\cdots,N_5,N_{234}$ are dyadic numbers, let $N_{\textup{max}} =\max (N_1,\cdots,N_5)$. Consider the quintic tensor $H_{n_0\cdots n_5}(t)$, which arises from the quintic stochastic object $\squintic$ and is defined as
\begin{equation}\label{def:quintic_tensor3}
\begin{aligned}
H_{n_0 n_1 \hdots n_5}(t) &=\, \mathbf{1}{\big\{ n_0 = n_{12345} \big\}}  \cdot  \prod_{j=0}^5 \frac{1_{N_j}(n_j)}{\langle n_j \rangle}\cdot\int\displaylimits_0^{t} \dt^{\prime}\, \chi(t) \chi(t') \sin\big((t-t') \langle n_{0} \rangle\big)  \varphi_1(t'\langle n_1\rangle)\\
&\times \frac{1_{N_{234}}(n_{234})}{\langle n_{234} \rangle} \varphi_5(t')\Big( \int\displaylimits_0^{t'} \dt^{\prime\prime}\, \chi(t') \chi(t'') \sin\big((t'-t'') \langle n_{234} \rangle\big)  \prod_{j=2}^4\varphi_j(t''\langle n_j\rangle) \Big).
\end{aligned}
\end{equation}
Then we have
\begin{align}\label{counting:quintic_bound3}
 \| \langle\lambda \rangle^{b_+}\widetilde{H}_{n_0\cdots n_5}(\lambda)\|_{L^2_\lambda[n_0\cdots n_5]}&\lesssim  N_{\textup{max}}^{-1/2+\epsilon},\\
\| \langle\lambda \rangle^{b_+}\widetilde{H}_{n_0\cdots n_5}(\lambda)\|_{L^2_\lambda[n_0n_A\to n_Bn_5]}&\lesssim N_{\mathrm{max}}^\epsilon(N_0N_5)^{-\frac{1}{2}}\cdot \big\{\max (N_0, N_2, N_3, N_4)^{-\frac{1}{2}}
\nonumber\\&\qquad\qquad\qquad+\max (N_2, N_3, N_4, N_5)^{-\frac{1}{2}}\big\}, \label{counting:quintic_bound4}
\end{align}
for any choice of $\varphi_j$ and any partition $(A,B)$ of $\{1,2,3,4\}$.
\end{corollary}
\begin{proof} This follows from Lemma \ref{counting:quintic_tensor} in the same way that Corollary \ref{counting:cor-cubic} follows from Lemma \ref{counting:cubic_tensor}.
\end{proof}
\subsection{The sine-cancellation kernel and tensor}\label{sine-tensor} 
In this subsection, we study the sine-cancellation tensor, which occurs in the quintic stochastic object with one pairing (see e.g. Subsection \ref{section:quintic-diagram}). The estimates of this tensor rely on a cancellation originating from the sine-function in the Duhamel integral which was first observed in \cite{GKO18}. Before introducing the sine-cancellation tensor, we first introduce the simpler $\Sine$-kernel.
\begin{definition}[The sine-cancellation kernel]\label{def:sine_kernel}
For any frequency-scales $K$ and $L$ and any $r\in \Z^3$, we define the sine-cancellation kernel $\Sine\colon \R \rightarrow \R$ by 
\begin{equation}
    \Sine(t,r) = \Sine[K,L](t,r) := \sum_{\substack{n_1,n_2 \in \Z^3\colon \\ n_{12} = r}} 1_{K}(n_1) 1_L(n_2) \frac{\sin\big(t \langle n_1 \rangle\big)}{\langle n_1 \rangle} \frac{\cos\big( t \langle n_2 \rangle\big)}{\langle n_2 \rangle^2}.
\end{equation}
\end{definition}

\begin{remark}
The choice of using $K$ and $L$ as notation for our frequency-scales instead of $N_1$ and $N_2$ is deliberate, since using $N_1$ and $N_2$ would be confusing in the symmetrization below.
\end{remark}

\begin{lemma}[Symmetrization of the sine-cancellation kernel]\label{counting:lem-Sine-symmetrization}
For any frequency-scales $K$ and $L$ and any $r\in \Z^3$, it holds that 
\begin{align}
\Sine[K,L](t,r) 
&= \frac{1}{4} \sum_{\substack{n_1,n_2 \in \Z^3 \colon \\ n_{12}= r }} 
\bigg[ 1_{K}(n_1) 1_{L}(n_2) \frac{\langle n_1 \rangle - \langle n_2 \rangle}{\langle n_1 \rangle^2 \langle n_2 \rangle^2} \sin\big( t (\langle n_1 \rangle - \langle n_2 \rangle) \big) \bigg] \label{counting:eq-Sine-i1}\\
&+ \frac{1}{4}  \sum_{\substack{n_1,n_2 \in \Z^3 \colon \\ n_{12}= r }} 
\bigg[ 1_{K}(n_1) 1_{L}(n_2) \frac{\langle n_1 \rangle + \langle n_2 \rangle}{\langle n_1 \rangle^2 \langle n_2 \rangle^2} \sin\big( t (\langle n_1 \rangle + \langle n_2 \rangle) \big) \bigg]  \label{counting:eq-Sine-i2} \\
&- \frac{1}{2}\sum_{\substack{n_1,n_2 \in \Z^3 \colon \\ n_{12}= r }} 
\bigg[ \Big(  1_{K}(n_1) 1_{L}(n_2) - 1_L (n_1) 1_K(n_2) \Big) \frac{\cos\big( t \langle n_1 \rangle\big)}{\langle n_1 \rangle^2} \frac{\sin\big( t \langle n_2 \rangle\big)}{\langle n_2 \rangle} \bigg]. \label{counting:eq-Sine-i3}
\end{align}
Furthermore, on the support of the indicator function $1_{K}(n_1) 1_{L}(n_2) - 1_L (n_1) 1_K(n_2)$, the vectors $n_1$ and $n_2$ satisfy the $\Gamma$-condition \eqref{counting:def-gamma}.
\end{lemma}

\begin{remark}
As part of our proof of Lemma \ref{counting:lem-Sine-symmetrization}, we utilize a hidden total derivative in the $\Sine$-kernel, which was not observed in the earlier literature \cite{B20II,GKO18}.
\end{remark}

\begin{proof} The identity for the $\Sine$-kernel follows from direct calculations. As this calculation is important in our proof we shall detail it below. From the identity $\cos^\prime=-\sin$, we obtain that 
\begin{align}
\Sine(t,r) &=  \sum_{\substack{n_1,n_2 \in \Z^3\colon \\ n_{12} = r}} 1_{K}(n_1) 1_L(n_2) \frac{\sin\big(t \langle n_1 \rangle\big)}{\langle n_1 \rangle} \frac{\cos\big( t \langle n_2 \rangle\big)}{\langle n_2 \rangle^2} \notag \\
&= -  \sum_{\substack{n_1,n_2 \in \Z^3\colon \\ n_{12} = r}} \, 1_{K}(n_1) 1_L(n_2) \, \partial_t \bigg( \frac{\cos\big(t \langle n_1 \rangle\big)}{\langle n_1 \rangle^2} \bigg) \frac{\cos\big( t \langle n_2 \rangle\big)}{\langle n_2 \rangle^2}. 
\label{counting:eq-Sine-p1}
\end{align}
By symmetrizing in $n_1$ and $n_2$, it follows that 
\begin{align}
&\eqref{counting:eq-Sine-p1} \notag \\
=& - \frac{1}{2} 
\sum_{\substack{n_1,n_2 \in \Z^3\colon \\ n_{12} = r}}
\bigg[ 
\, 1_{K}(n_1) 1_L(n_2) \, \partial_t \bigg( \frac{\cos\big(t \langle n_1 \rangle\big)}{\langle n_1 \rangle^2} \bigg) \frac{\cos\big( t \langle n_2 \rangle\big)}{\langle n_2 \rangle^2}  \notag \\
&\hspace{19ex}
+ \, 1_{L}(n_1) 1_K(n_2) \, \frac{\cos\big(t \langle n_1 \rangle\big)}{\langle n_1 \rangle^2} \, \partial_t \bigg(  \frac{\cos\big( t \langle n_2 \rangle\big)}{\langle n_2 \rangle^2} \bigg) \bigg] \notag \\
=& - \frac{1}{2} \sum_{\substack{n_1,n_2 \in \Z^3\colon \\ n_{12} = r}}
\bigg[ 1_{K}(n_1) 1_L(n_2) \, \bigg( \partial_t \bigg( \frac{\cos\big(t \langle n_1 \rangle\big)}{\langle n_1 \rangle^2} \bigg) \frac{\cos\big( t \langle n_2 \rangle\big)}{\langle n_2 \rangle^2}  +  \frac{\cos\big(t \langle n_1 \rangle\big)}{\langle n_1 \rangle^2} \, \partial_t \bigg(  \frac{\cos\big( t \langle n_2 \rangle\big)}{\langle n_2 \rangle^2} \bigg) \bigg] \label{counting:eq-Sine-p2} \\
+& \frac{1}{2} \sum_{\substack{n_1,n_2 \in \Z^3\colon \\ n_{12} = r}}
\bigg[ \Big( 1_{K}(n_1) 1_L(n_2) - 1_L(n_1) 1_K(n_2) \Big) \frac{\cos\big(t \langle n_1 \rangle\big)}{\langle n_1 \rangle^2} \, \partial_t \bigg(  \frac{\cos\big( t \langle n_2 \rangle\big)}{\langle n_2 \rangle^2} \bigg) \bigg]. \label{counting:eq-Sine-p3}
\end{align}
After performing the $t$-derivative, the term \eqref{counting:eq-Sine-p3} contributes \eqref{counting:eq-Sine-i3} in the desired identity. Thus, it remains to consider \eqref{counting:eq-Sine-p2}. Due to the product rule, we can convert the summand into a total derivative, which yields 
\begin{equation}
\eqref{counting:eq-Sine-p2} =  - \frac{1}{2} \sum_{\substack{n_1,n_2 \in \Z^3\colon \\ n_{12} = r}}
\bigg[ 1_{K}(n_1) 1_L(n_2) \, \partial_t \bigg( 
 \frac{\cos\big(t \langle n_1 \rangle\big)}{\langle n_1 \rangle^2}  \frac{\cos\big( t \langle n_2 \rangle\big)}{\langle n_2 \rangle^2} \bigg) \bigg].
\end{equation}
Using a product-to-sum rule for trigonometric functions, we further obtain that 
\begin{align*}
     &- \frac{1}{2} \sum_{\substack{n_1,n_2 \in \Z^3\colon \\ n_{12} = r}}
\bigg[ 1_{K}(n_1) 1_L(n_2) \, \partial_t \bigg( 
 \frac{\cos\big(t \langle n_1 \rangle\big)}{\langle n_1 \rangle^2}  \frac{\cos\big( t \langle n_2 \rangle\big)}{\langle n_2 \rangle^2} \bigg) \bigg] \\
 =& - \frac{1}{4}  \sum_{\substack{n_1,n_2 \in \Z^3\colon \\ n_{12} = r}}
\bigg[ 1_{K}(n_1) 1_L(n_2) \, 
 \ \frac{1}{\langle n_1 \rangle^2} \frac{1}{\langle n_2 \rangle^2} \, \partial_t \Big( \cos\big( t (\langle n_1 \rangle - \langle n_2 \rangle) \big) + \cos \big( t ( \langle n_1 \rangle + \langle n_2 \rangle ) \big) \Big) \bigg]. 
\end{align*}
After performing the $t$-derivative, this leads to the terms \eqref{counting:eq-Sine-i1} and \eqref{counting:eq-Sine-i2} in our desired identity.

Now we consider the support of the indicator function $1_{K}(n_1) 1_{L}(n_2) - 1_L (n_1) 1_K(n_2)$. If $K=L$ this is clearly zero by symmetry; suppose $K\neq L$, say $K>L$ by symmetry, then $1_{K}(n_1) 1_{L}(n_2)$ is supported in $|n_1|_{\infty}\geq K\geq |n_2|_\infty$, and $ 1_L (n_1) 1_K(n_2)$ is supported in $|n_2|_{\infty}\geq K\geq |n_1|_\infty$, which is what we need.
\end{proof}

We now prove a direct estimate of the $\Sine$-kernel which, while not as powerful as Lemma \ref{counting:sine_tensor}, will sometimes be easier to use.

\begin{lemma}[Direct estimate of the $\Sine$-kernel]\label{counting:lem-Sine-estimate}
For all frequency-scales $K$ and $L$, all $r\in \Z^3$, all times $t\in \R$, and all $\lambda\in \R$, it holds that 
\begin{equation}\label{counting:eq-Sine-estimate}
\Big| \int_0^t \dt^\prime \, \Sine[K,L](t-t^\prime,r) \, e^{i\lambda t^\prime} \,\Big| 
\lesssim \langle t \rangle \, \frac{\log^2\big( 2 + \max(K,L,|\lambda|)\big)}{\max(K,L,|\lambda|)}. 
\end{equation}
\end{lemma}

\begin{proof}
We distinguish the cases $|\lambda|\gg \max(K,L)$ and $|\lambda|\lesssim \max(K,L)$. If $|\lambda|\gg \max(K,L)$, we obtain from the definition of the $\Sine$-kernel that 
\begin{align*}
 \Big| \int_0^t \dt^\prime \, \Sine[K,L](t-t^\prime,r) \, e^{i\lambda t^\prime} \,\Big|
 &\lesssim \langle t \rangle \, \langle \lambda \rangle^{-1} \sum_{\substack{n_1,n_2 \in \Z^3\colon \\ n_{12}=r}} 1_{K}(n_1) 1_L(n_2) \frac{1}{\langle n_1 \rangle} \frac{1}{\langle n_2 \rangle^2}  \\
 &\lesssim \langle t \rangle \, \langle \lambda \rangle^{-1} K^{-1} L^{-2} \min(K,L)^3 \\
 &\lesssim \langle t \rangle \, \langle \lambda \rangle^{-1}. 
\end{align*}
Thus, it remains to treat the case $|\lambda|\lesssim \max(K,L)$. After using Lemma \ref{counting:lem-Sine-symmetrization} and performing the $t^\prime$-integral, it follows that 
\begin{align}
  & \Big| \int_0^t \dt^\prime \, \Sine[K,L](t-t^\prime,r) \, e^{i\lambda t^\prime} \,\Big| \notag \\
  \lesssim& \, \langle t \rangle  \langle r \rangle K^{-2} L^{-2}  \sum_{\pm} \sum_{\substack{n_1,n_2 \in \Z^3\colon \\ n_{12}=r}} 1_{K}(n_1) 1_L(n_2) \frac{1}{1+|\langle n_1 \rangle - \langle n_2 \rangle \pm \lambda|} \label{counting:eq-Sine-estimate-p1} \\ 
  +& \, \langle t \rangle K^{-2} L^{-2} \max(K,L) \sum_{\substack{n_1,n_2 \in \Z^3\colon \\ n_{12}=r}} 1_{K}(n_1) 1_L(n_2) \frac{1}{1+|\langle n_1 \rangle + \langle n_2 \rangle \pm \lambda|} \label{counting:eq-Sine-estimate-p2} \\ 
  +& \, \langle t \rangle \min(K,L)^{-1} K^{-1} L^{-1} 
   \sum_{\pm_1, \pm_2} \sum_{\substack{n_1,n_2 \in \Z^3\colon \\ n_{12}=r}} \bigg[ \big| 1_{K}(n_1) 1_{L}(n_2) - 1_L (n_1) 1_K(n_2) \big|  \label{counting:eq-Sine-estimate-p3}\\
   &\hspace{42ex}\times \frac{1}{1+ |\langle n_1 \rangle \pm_1 \langle n_2 \rangle \pm_2 \lambda|}  \bigg]. \notag 
\end{align}
Since all three summands \eqref{counting:eq-Sine-estimate-p1}, \eqref{counting:eq-Sine-estimate-p2}, and \eqref{counting:eq-Sine-estimate-p3} are symmetric in $K$ and $L$, we now assume that $K\geq L$. Using Lemma \ref{counting:lem1} and a level-set decomposition of $\langle n_1\rangle - \langle n_2\rangle$, the first summand is estimated by 
\begin{align*}
\eqref{counting:eq-Sine-estimate-p1} 
&\lesssim 
\langle t \rangle \log(2+K) \langle r \rangle K^{-2} L^{-2} \min\big( \langle r \rangle, L\big)^{-1} L^3 \\
&\lesssim \langle  t \rangle \log(2+K) \max\big( \langle r \rangle, L \big)  K^{-2} \\
&\lesssim \langle t \rangle \log(2+K) K^{-1}. 
\end{align*}
Using Lemma \ref{counting:lem1}, the good sign in $\langle n_1\rangle + \langle n_2\rangle$, and a level-set decomposition of $\langle n_1\rangle + \langle n_2\rangle$, the second summand is estimated by 
\begin{align*}
    \eqref{counting:eq-Sine-estimate-p2} \lesssim 
    \langle  t \rangle \log(2+K) K^{-1} L^{-2} L^2  \lesssim \langle  t \rangle \log(2+K) K^{-1}. 
\end{align*}
It remains to treat the third summand \eqref{counting:eq-Sine-estimate-p3}.  We further distinguish the two cases $\langle  r \rangle \gtrsim L$ and $\langle r \rangle \ll L$. In the case $\langle r \rangle \gtrsim L$, we utilize Lemma  \ref{counting:lem1} and a level-set decomposition of $\langle n_1 \rangle \pm_1 \langle n_2 \rangle$, which yield 
\begin{align*}
\eqref{counting:eq-Sine-estimate-p3} &\lesssim
 \, \langle t \rangle  K^{-1} L^{-2} 
   \sum_{ \pm_1, \pm_2} \sum_{\substack{n_1,n_2 \in \Z^3\colon \\ n_{12}=r}} \bigg[ \Big( 1_{K}(n_1) 1_{L}(n_2) + 1_L (n_1) 1_K(n_2) \Big)   \frac{1}{1+ | \langle n_1 \rangle \pm_1 \langle n_2 \rangle \pm_2 \lambda|}  \bigg] \\
   &\lesssim \langle t \rangle \log(2+K) K^{-1}. 
\end{align*}
In the case $\langle r \rangle \ll L$, we first note that $L\leq K \lesssim \max(L,\langle r \rangle) \lesssim L$. As a result, it follows that $K\sim L$. 
Together with the support property in Lemma \ref{counting:lem-Sine-symmetrization}, it follows that 
\begin{equation}\label{counting:eq-Sine-estimate-p4}
\begin{aligned}
\eqref{counting:eq-Sine-estimate-p3} &\lesssim
\langle t \rangle K^{-3} \sum_{\pm_1, \pm_2} \sum_{\substack{n_1,n_2 \in \Z^3\colon \\ n_{12}=r}}  
 \mathbf{1}\Big\{|n_1|_\infty\geq \Gamma\geq |n_1-r|_\infty\mathrm{\ or\ }|n_1|_\infty\leq\Gamma\leq |n_1-r|_\infty\Big\}\\&\,\times\frac{1}{1+ | \langle n_1 \rangle \pm_1 \langle n_2 \rangle \pm_2 \lambda|}.
 \end{aligned}
\end{equation}
 By Lemma \ref{counting:lem1''} and a similar level-set decomposition, it follows that
\begin{align*}
\eqref{counting:eq-Sine-estimate-p4} \lesssim \langle t \rangle \log^2(2+K) K^{-1},
\end{align*}
which completes the estimate of the third summand \eqref{counting:eq-Sine-estimate-p3}. 
\end{proof}

\begin{lemma}[The sine-cancellation tensor estimates]\label{counting:sine_tensor}
 Suppose $N_0, N_2,\cdots,N_5,N_{234}$ are dyadic numbers and $\lambda, \lambda_3,\lambda_4,\lambda_5$ are real numbers, define $N_{\max} = \max (N_0, N_2,\cdots,N_5)$. Recall the Sine kernel defined in Definition \ref{def:sine_kernel}. Define the tensors
 \begin{equation}
\label{def:sine_tensor1}
\begin{aligned}
h^{\mathrm{sine}}_{n_0 n_3n_4 n_5}(t,\lambda_3,\lambda_4, \lambda_5) &=\, \mathbf{1}{\big\{ n_0 = n_{345}\big\}} \cdot 1_{N_0}(n_0)\cdot  \Big( \prod_{j=3}^5 \frac{1_{N_j}(n_j)}{\langle n_j \rangle} \Big) 
  e^{it(\pm\langle n_5\rangle+\lambda_5)}\\
&\times  \Big( \int\displaylimits_0^{t} \dt^{\prime}\, \chi(t) \chi(t') \cdot \Sine[N_{234},N_2](t-t^\prime,n_{34}) e^{it'(\pm\langle n_3\rangle\pm\langle n_4\rangle+\lambda_3+\lambda_4)} \Big),
\end{aligned}
\end{equation}
\begin{equation}\label{def:sine_tensor2}\begin{aligned}
H^{\mathrm{sine}}_{n_0 n_3n_4 n_5}(t,\lambda_3,\lambda_4, \lambda_5) &= \mathbf{1}{\big\{ n_0 = n_{345}\big\}} \frac{1_{N_0}(n_0)}{\langle n_0 \rangle}\prod_{j=3}^5 \frac{1_{N_j}(n_j)}{\langle n_j \rangle} \int\displaylimits_0^{t} \dt^{\prime}\, \chi(t)\chi(t')\cdot \sin\big((t-t')\langle n_0\rangle\big) \\
&\times e^{it'(\pm\langle n_5\rangle+\lambda_5)} \Big( \int\displaylimits_0^{t'} \dt^{\prime\prime}\, \chi(t') \chi(t'') \cdot \Sine[N_{234},N_2](t'-t^{\prime\prime},n_{34}) \, e^{it''(\pm\langle n_3\rangle\pm\langle n_4\rangle+\lambda_3+\lambda_4)} \Big).
\end{aligned}
\end{equation} Then there exist two functions $C_j=C_j(\lambda,\lambda_3,\lambda_4,\lambda_5)$ for $1\leq j\leq 2$, such that
\begin{equation}
\label{counting:sine_coef}
\| \langle \lambda \rangle^{b_+-1}  C_1(\lambda,\lambda_3,\lambda_4,\lambda_5)\|_{L^2_{\lambda}}\lesssim N_{\max}^{\epsilon},\qquad \| \langle \lambda \rangle^{b_+}  C_2(\lambda,\lambda_3,\lambda_4,\lambda_5)\|_{L^2_{\lambda}}\lesssim N_{\max}^{\epsilon},
\end{equation} and that we have
\begin{align}
\label{counting:sine_bound1}
\|\widetilde{h^{\mathrm{sine}}}(\lambda, \lambda_3,\lambda_4,\lambda_5)\|_{n_0n_3n_4 n_5}&\lesssim  N_0 N_{\max}^{-\frac{1}{2}}\cdot  C_1(\lambda,\lambda_3,\lambda_4,\lambda_5),\\
\label{counting:sine_extra_HS}
\|\widetilde{h^{\mathrm{sine}}}(\lambda, \lambda_3,\lambda_4,\lambda_5)\|_{n_0n_3n_4 n_5}&\lesssim \frac{\min(N_0,N_5)^{\frac{3}{2}}(N_3N_4)^{\frac{1}{2}}}{N_5\cdot\max(N_2,N_{234})}   C_1(\lambda,\lambda_3,\lambda_4,\lambda_5),\\
\label{counting:sine_bound2}
\|\widetilde{H^{\mathrm{sine}}}(\lambda, \lambda_3,\lambda_4,\lambda_5)\|_{n_0n_3n_4 n_5}&\lesssim  N_{\max}^{-\frac{1}{2}}\cdot  C_2(\lambda,\lambda_3,\lambda_4,\lambda_5),\\
\label{counting:sine_bound1'}
\big\|\widetilde{h^{\mathrm{sine}}}(\lambda, \lambda_3,\lambda_4,\lambda_5)\big\|_{n_0n_A\to n_Bn_5}
&\lesssim \,\,N_0^{\frac{1}{2}} N_5^{-\frac{1}{2}}\cdot N_{\max}^{-\frac{1}{2}} \cdot  C_1(\lambda,\lambda_3,\lambda_4,\lambda_5),\\
\label{counting:sine_bound2'}
\big\|\widetilde{H^{\mathrm{sine}}}(\lambda, \lambda_3,\lambda_4,\lambda_5)\big\|_{n_0n_A\to n_Bn_5}
&\lesssim \,\,(N_0 N_5)^{-\frac{1}{2}}\cdot N_{\max}^{-\frac{1}{2}} \cdot  C_2(\lambda,\lambda_3,\lambda_4,\lambda_5),\\
\label{counting:sine_bound1''}
\big\|\widetilde{h^{\mathrm{sine}}}(\lambda, \lambda_3,\lambda_4,\lambda_5)\big\|_{n_0n_5\to n_3n_4}
&\lesssim \,\,N_0\big(\max (N_0, N_5)\cdot\max (N_3,N_4)\big)^{-\frac{1}{2}}\cdot N_2^{-1}\cdot  C_1(\lambda,\lambda_3,\lambda_4,\lambda_5),\\
\label{counting:sine_bound2''}
\big\|\widetilde{H^{\mathrm{sine}}}(\lambda, \lambda_3,\lambda_4,\lambda_5)\big\|_{n_0n_5\to n_3n_4}
&\lesssim \,\, \big(\max (N_0, N_5)\cdot\max (N_3,N_4)\big)^{-\frac{1}{2}}\cdot N_2^{-1}\cdot  C_2(\lambda,\lambda_3,\lambda_4,\lambda_5),
\end{align} for any partition $(A,B)$ of $\{3,4\}$. The same bounds hold for all $\partial_{\lambda}$ and $\partial_{\lambda_j}$ derivatives of these tensors.
\end{lemma}
\begin{proof}
In the proof below we will also fix a dyadic number $N_{34}$ and insert the indicator function $\mathbf{1}_{N_{34}}(n_3+n_4)$; note that summing over ${N_{34}}$ introduces a $\log(N_{\mathrm{max}})$ factor, but this can easily be absorbed by slightly adjusting the definition of $C_1$. Similarly, we will ignore all possible $\log (N_{\mathrm{max}})$ factors that may occur below. 

By \eqref{tilted} and Lemma \ref{counting:lem-Sine-symmetrization}, we can decompose the tensor $(\widetilde{h^{\mathrm{sine}}})_{n_0n_3n_4n_5}(\lambda,\lambda_3,\lambda_4,
\lambda_5)$ into $h^{(-)}$, $h^{(+)}$ and $h^{(0)}$, where $h_{n_0n_3n_4n_5}^{(*)}=h_{n_0n_3n_4n_5}^{(*)}(\lambda,\lambda_3,\lambda_4,\lambda_5)$ for $(*)\in\{(-),(+),(0)\}$, and
\begin{multline}\label{eq:h-}
h^{(-)}_{n_0n_3n_4n_5}=\sum_{n_2\in \mathbb{Z}^3}  K\big(\lambda-\lambda_5\pm\langle n_5\rangle\pm\langle n_0\rangle\mp_{234} (\langle n_{234}\rangle-\langle n_2\rangle), \lambda_3+\lambda_4\pm_{234} (\langle n_{234}\rangle-\langle n_2\rangle)\pm \langle n_3\rangle\pm\langle n_4\rangle)\big)\\
\times \mathbf{1}{\big\{ n_0 = n_{345}\big\}}  1_{N_0}(n_0)1_{N_{34}}(n_3+n_4) 1_{N_{234}}(n_{234})\cdot 1_{N_2}(n_2)\cdot \frac{\langle n_{234}\rangle-\langle n_2\rangle}{\langle n_{234}\rangle^2\langle n_2\rangle^2}  \Big( \prod_{j=3}^5 \frac{1_{N_j}(n_j)}{\langle n_j \rangle} \Big),
\end{multline}
\begin{multline*}
h^{(+)}_{n_0n_3n_4n_5}=\sum_{n_2\in \mathbb{Z}^3}  K\big(\lambda-\lambda_5\pm\langle n_5\rangle\pm\langle n_0\rangle\mp_{234} (\langle n_{234}\rangle+\langle n_2\rangle), \lambda_3+\lambda_4\pm_{234} (\langle n_{234}\rangle+\langle n_2\rangle)\pm \langle n_3\rangle\pm\langle n_4\rangle)\big)\\
\times\mathbf{1}{\big\{ n_0 = n_{345}\big\}}  1_{N_0}(n_0)1_{N_{34}}(n_3+n_4) \cdot 1_{N_{234}}(n_{234})\cdot 1_{N_2}(n_2)\cdot \frac{\langle n_{234}\rangle+\langle n_2\rangle}{\langle n_{234}\rangle^2\langle n_2\rangle^2}  \Big( \prod_{j=3}^5 \frac{1_{N_j}(n_j)}{\langle n_j \rangle} \Big),  
\end{multline*}
\begin{multline*}
h^{(0)}_{n_0n_3n_4n_5}=\sum_{n_2\in \mathbb{Z}^3}  K\big(\lambda-\lambda_5\pm\langle n_5\rangle\pm\langle n_0\rangle\mp_{234} (\langle n_{234}\rangle\pm\langle n_2\rangle), \lambda_3+\lambda_4\pm_{234} (\langle n_{234}\rangle\pm\langle n_2\rangle)\pm \langle n_3\rangle\pm\langle n_4\rangle)\big)\\
\times \mathbf{1}{\big\{ n_0 = n_{345}\big\}}  1_{N_0}(n_0) 1_{N_{34}}(n_3+n_4)\cdot \frac{ 1_{N_{234}}(n_{234})\cdot 1_{N_2}(n_2)- 1_{N_{2}}(n_{234})\cdot 1_{N_{234}}(n_2)}{\langle n_{234}\rangle^2\langle n_2\rangle}  \Big( \prod_{j=3}^5 \frac{1_{N_j}(n_j)}{\langle n_j \rangle} \Big).
\end{multline*} Similarly we can decompose $\widetilde{H^{\mathrm{sine}}}$ into $H^{(-)}$, $H^{(+)}$ and $H^{(0)}$ components. We will focus on the discussion of $h^{(-)}$ and $H^{(-)}$ in parts (1)--(3) below; for the other two components we will only point out the points of difference in the proof in part (4). 

(1) Start with the formula (\ref{eq:h-}) for $h^{(-)}$. Similar to the proof of Lemma \ref{counting:quintic_tensor}, define 
\begin{equation}
C_1(\lambda,\lambda_3,\lambda_4,\lambda_5):= \sum_{\substack{m_1,m_2\in \mathbb{Z}\\
|m_1|,|m_2|\lesssim N_{\max}}} \frac{1 }{\langle \lambda-\lambda_5- m_1 \rangle \langle \lambda-\sum_{j=3}^5 \lambda_j -m_2\rangle},
\end{equation} then in the same way as (\ref{counting:quintic_coef}) we can prove the first inequality in (\ref{counting:sine_coef}). By level set decomposition we then have
\begin{multline}\label{counting:levelset1}
\|h^{(-)}\|_{n_0n_3n_4n_5}
\lesssim \, (N_3N_4N_5)^{-1}
C_1(\lambda,\lambda_3,\lambda_4,\lambda_5)\cdot
\sup_{m_1,m_2\in \mathbb{Z}} \bigg\|\mathbf{1}\big\{n_0=n_{345}\big\}1_{N_0}(n_0) 1_{N_{34}}(n_3+n_4)\prod_{j=3}^5\mathbf{1}_{N_j}(n_j)\\\times \mathbf{1}\big\{|\Omega-m_2|\leq 1\big\}\cdot\bigg(\sum_{n_2}1_{N_2}(n_2)1_{N_{234}}(n_{234}) \frac{|\langle n_{234}\rangle-\langle n_2\rangle|}{\langle n_{234}\rangle^2\,\langle n_2\rangle^2} \cdot \mathbf{1}\big\{|\Omega'-m_1|\leq 1\big\}\bigg)\bigg\|_{n_0n_3n_4n_5}.
\end{multline}
 where \begin{equation}\label{counting:resfactors}
\Omega =\pm\langle n_0\rangle\pm\langle n_3\rangle\pm\langle n_4\rangle\pm\langle n_5\rangle,\quad \Omega' = \pm\langle n_0\rangle\pm\langle n_5\rangle\mp_{234}(\langle n_{234}\rangle-\langle n_2\rangle).
\end{equation} For simplicity, we will define the tensor inside the $\|\cdot\|$ norm sign to be $\Hc_{n_0n_3n_4n_5}^{(-)}$. The idea is to examine the support of $\Hc^{(-)}$ together with the size $|\Hc^{(-)}|$.

Clearly $\Hc^{(-)}$ is supported in the set where $\langle n_j\rangle\sim N_j\,(j\in\{0,3,4,5\})$, $|n_3+n_4|\sim N_{34}$, and $n_0=n_{345}$ and $|\Omega-m_2|\lesssim 1$. Moreover, by Lemma \ref{counting:lem1} and the restriction $|\Omega'-m_1|\leq 1$, we know that
\begin{equation}\label{counting:sine-sum}|\Hc_{n_0n_3n_4n_5}^{(-)}|\lesssim N_2^{-2}N_{234}^{-2}N_{34}\cdot\big(\min(N_2,N_{234})^2+\min(N_2,N_{234})^3N_{34}^{-1}\big)\lesssim \max(N_2,N_{234})^{-1}.\end{equation} As such, we have that
\begin{equation}\label{counting:sine_extra_eq1}
\|h^{(-)}\|_{n_0n_3n_4n_5}^2\lesssim \bigl((N_3N_4N_5)^{-2} C_1^2 \max(N_2,N_{234})^{-2}\bigr) \cdot\sum_{(n_0,n_3,n_4,n_5)}1,\end{equation} where the summation is taken over the support of $\Hc^{(-)}$ described above. 
This summation can be trivially bounded by $(N_3N_4\min(N_0,N_5))^3$ which together with the first factor in \eqref{counting:sine_extra_eq1} proves \eqref{counting:sine_extra_HS} for $h^{(-)}$.

Let $i\in\{0,3,4,5\}$ be such that $N_i=\min(N_0,N_3,N_4,N_5)$. If $\max(N_2,N_{234})\gtrsim N_i$ then Lemma \ref{counting:lem2} already implies
\[\|h^{(-)}\|_{n_0n_3n_4n_5}^2\lesssim(N_3N_4N_5)^{-2} C_1^2\cdot\max(N_2,N_{234})^{-2}\cdot (N_0N_3N_4N_5)^2\frac{N_i}{\max(N_0,N_3,N_4,N_5)}\lesssim N_0^2N_{\mathrm{max}}^{-1}\cdot C_1^2,\] which is what we need. Now suppose $\max(N_2,N_{234})\ll N_i$, then we shall fix $(n_i,n_j)$, where $\{i,j\}\in\{\{0,5\},\{3,4\}\}$. Using the bound $|n_i\pm n_j|\sim N_{34}$, we apply Lemma \ref{counting:lem1} to get that
\[\|h^{(-)}\|_{n_0n_3n_4n_5}^2\lesssim(N_3N_4N_5)^{-2} C_1^2\cdot\max(N_2,N_{234})^{-2}\cdot N_i^3N_{34}^3\cdot N_{\mathrm{max}}^3N_{34}^{-1}\lesssim N_0^2N_i^{-1}N_{\mathrm{max}}^{-1}\cdot C_1^2,\] which is more than what we need. This proves (\ref{counting:sine_bound1}) for $h^{(-)}$. The proof of (\ref{counting:sine_extra_HS}) is similar, where one counts the number of choices of $(n_0,n_3,n_4,n_5)$ by first fixing either $n_0$ or $n_5$ depending on whether $N_0\lesssim N_5$ or $N_5\lesssim N_0$, and then fixing $(n_3,n_4)$.

The proof of (\ref{counting:sine_bound2}) for $h^{(-)}$ (plus the second inequality in (\ref{counting:sine_coef})) also follows similarly. In fact, we shall define
\[C_2(\lambda,\lambda_3,\lambda_4,\lambda_5)=\sum_{\substack{m_1,m_2\in \mathbb{Z}\\
|m_1|,|m_2|\lesssim N_{\max}}} \int_\mathbb{R}\frac{\mathrm{d}\sigma }{\langle\lambda\rangle\langle\lambda-\sigma\rangle\langle \sigma-\lambda_5- m_1 \rangle \langle \sigma-\sum_{j=3}^5 \lambda_j -m_2\rangle},\] which then satisfies the second inequality in (\ref{counting:sine_coef}) by repeating the arguments in the proof of Lemma \ref{counting:quintic_tensor}, and also gives the bound (\ref{counting:sine_bound2}) for $h^{(-)}$ by making the level set decomposition similar to (\ref{counting:levelset1}) and applying the bound for $\Hc^{(-)}$ obtained above.

(2) Now we turn to (\ref{counting:sine_bound1'}); the proof of (\ref{counting:sine_bound2'}) follows the same lines as above, which we will omit. By the same reduction above, it suffices to prove that
\begin{equation}\label{counting:sine-ex1}
\|\Hc^{(-)}\|_{n_0n_A\to n_Bn_5}\lesssim N_0N_3N_4N_5\cdot (N_0 N_5)^{-\frac{1}{2}}\cdot N_{\max}^{-\frac{1}{2}}.
\end{equation} Here again we will use the upper bound (\ref{counting:sine-sum}) of $\Hc^{(-)}$ and Schur's lemma. If $A=\{3\},\,B=\{4\}$ (the other case $A=\{4\},B=\{3\}$ being similar), then Schur's bound and the assumption $|n_3+n_4|\sim N_{34}$ implies that
\[
\begin{aligned}\|\Hc^{(-)}\|_{n_0n_3\to n_4n_5}&\lesssim \max(N_2,N_{234})^{-1}\min(N_0,N_3,N_{34})^{3/2}\min(N_5,N_4,N_{34})^{3/2}\\
&\lesssim \max(N_2,N_{234})^{-1}\min(N_0,N_3)^{1/2}N_3^{1/2}N_{34}^{1/2}\cdot\min(N_5,N_4)^{1/2}N_4^{1/2}N_{34}^{1/2}\\
&=N_0N_3N_4N_5\cdot (N_0N_5)^{-1/2}\cdot\max(N_0,N_3)^{-1/2}\max(N_5,N_4)^{-1/2}\frac{N_{34}}{\max(N_2,N_{234})}\\
&\lesssim N_0N_3N_4N_5\cdot (N_0N_5)^{-1/2}\cdot N_{\mathrm{max}}^{-1/2},
\end{aligned}\] which is what we need. If instead $A=\varnothing,\,B=\{3,4\}$  (the other case $A=\{3,4\},\,B=\varnothing$ is similar), then Lemma \ref{counting:lem2} already gives
\[\|\Hc^{(-)}\|_{n_0\to n_3n_4n_5}\lesssim \max(N_2,N_{234})^{-1}N_0N_3N_4N_5\cdot N_0^{-1}N_5^{-1/2},\] which is enough if $N_{\mathrm{max}}\sim N_0\mathrm{\ or\ }N_2\mathrm{\ or\ }N_{234}$; otherwise we have $N_0,N_2,N_{234}\ll N_{\mathrm{max}}$, then it must be that $N_3\sim N_4\sim N_{\mathrm{max}}$, hence in applying Schur's bound and counting $(n_3,n_4,n_5)$ we may first fix $n_5$ and then apply Lemma \ref{counting:lem1} to get
\[\|\Hc^{(-)}\|_{n_0\to n_3n_4n_5}\lesssim N_{34}^{-1}\cdot N_{34}^{3/2}\cdot N_3^{3/2}N_{34}^{-1/2}=N_3^{3/2}\lesssim N_0N_3N_4N_5\cdot (N_0N_5)^{-1/2}\cdot N_{\mathrm{max}}^{-1/2},\] which is also adequate. This finishes the proof of (\ref{counting:sine_bound1'}) and (\ref{counting:sine_bound2'}) for $h^{(-)}$.

(3) Next consider (\ref{counting:sine_bound1''}) (together with (\ref{counting:sine_bound2''})). Again it suffices to prove that
\begin{equation}\label{counting:sine-ex2}
\|\Hc^{(-)}\|_{n_0n_5\to n_3n_4}\lesssim N_0N_3N_4N_5\cdot  \big(\max (N_0, N_5)\cdot\max (N_3,N_4)\big)^{-\frac{1}{2}}\cdot N_2^{-1}.
\end{equation} This follows directly from the pointwise bound (\ref{counting:sine-sum}) and the fact that the base tensor $h_{n_0n_3n_4n_5}^b$ satisfies that
\[\|h^b\|_{n_0n_5\to n_3n_4}\lesssim N_0N_3N_4N_5\cdot \big(\max (N_0, N_5)\cdot\max (N_3,N_4)\big)^{-\frac{1}{2}},\] which is proved using Schur's bound and Lemma \ref{counting:lem1}. This proves (\ref{counting:sine_bound1''}) and (\ref{counting:sine_bound2''}) for $h^{(-)}$.

(4) Finally, we demonstrate how the same proof applies to the terms $h^{(+)}$ and $h^{(0)}$. In fact, the only difference is that the tensor $\Hc^{(-)}$ will be replaced by $\Hc^{(+)}$ and $\Hc^{(0)}$, which are defined as below. The tensor $\Hc^{(+)}$ is defined as 
\begin{equation}\label{counting:defh+}
\begin{aligned}
\Hc_{n_0n_3n_4n_5}^{(+)}&:=\mathbf{1}\big\{n_0=n_{345}\big\}1_{N_0}(n_0)1_{N_{34}}(n_3+n_4)\prod_{j=3}^5 1_{N_j}(n_j)\\&\times \mathbf{1}\big\{|\Omega-m_2|\leq 1\big\}\cdot\bigg(\sum_{n_2}1_{N_2}(n_2)1_{N_{234}}(n_{234}) \frac{\langle n_{234}\rangle+\langle n_2\rangle}{\langle n_{234}\rangle^2\,\langle n_2\rangle^2} \cdot \mathbf{1}\big\{|\Omega^{(+)}-m_1|\leq 1\big\}\bigg)
\end{aligned}
\end{equation} and $\Hc^{(0)}$ is defined as
\begin{align}
\Hc_{n_0n_3n_4n_5}^{(0)}&:=\mathbf{1}\big\{n_0=n_{345}\big\}1_{N_0}(n_0)1_{N_{34}}(n_3+n_4)\prod_{j=3}^5 1_{N_j}(n_j)\cdot\mathbf{1}\big\{|\Omega-m_2|\leq 1\big\}\cdot\bigg(\sum_{n_2}1_{N_2}(n_2)1_{N_{234}}(n_{234})\nonumber\\&\times \frac{ 1_{N_{234}}(n_{234})\cdot 1_{N_2}(n_2)- 1_{N_{2}}(n_{234})\cdot 1_{N_{234}}(n_2)}{\langle n_{234}\rangle^2\langle n_2\rangle} \cdot \mathbf{1}\big\{|\Omega^{(0)}-m_1|\leq 1\big\}\bigg)\label{counting:defh0}
\end{align}where $\Omega$ is as in (\ref{counting:resfactors}), but $\Omega^{(+)}$ and $\Omega^{(0)}$ are defined as
\[\Omega^{(+)} = \pm\langle n_0\rangle\pm\langle n_5\rangle\mp_{234}(\langle n_{234}\rangle+\langle n_2\rangle),\quad \Omega^{(0)}=\pm\langle n_0\rangle\pm\langle n_5\rangle\pm_{234} (\langle n_{234}\rangle\pm\langle n_2\rangle) \] with some arbitrary choices of $\pm$ signs. Clearly $\Hc^{(+)}$ and $\Hc^{(0)}$ have the same support as $\Hc^{(-)}$, so it suffices to prove that they also satisfy the same pointwise bounds (\ref{counting:sine-sum}) as $\Hc^{(-)}$, which then guarantees that the same proof will go through.

The estimate for $\Hc^{(+)}$ follows from Lemma \ref{counting:lem1}, because of the favorable sign in $\langle n_{234}\rangle+\langle n_2\rangle$; in fact by (\ref{counting:eq2}) we have
\[|\Hc_{n_0n_3n_4n_5}^{(+)}|\lesssim N_2^{-2}N_{234}^{-2}(N_2+N_{234})\cdot\min(N_2,N_{234})^2\lesssim\max(N_2,N_{234})^{-1},\] which is just (\ref{counting:sine-sum}). The estimate for $\Hc^{(0)}$, on the other hand, follows from Lemma \ref{counting:lem1''}. In fact, exploiting the $\Gamma$ condition in the support of $1_{N_{234}}(n_{234})\cdot 1_{N_2}(n_2)- 1_{N_{2}}(n_{234})\cdot 1_{N_{234}}(n_2)$, and using (\ref{counting:gamma}), we have
\[|\Hc_{n_0n_3n_4n_5}^{(0)}|\lesssim N_2^{-1}N_{234}^{-2}\cdot \min(N_2,N_{234})^2\log (N_{\mathrm{max}}) \lesssim \max(N_2,N_{234})^{-1}\log (N_{\mathrm{max}}),\] which is just (\ref{counting:sine-sum}), where the $\log N_{\mathrm{max}}$ loss can be absorbed by slightly adjusting the definition of $C_1$ (or $C_2$). This completes the proof (the $\partial_{\lambda}$ and $\partial_{\lambda_j}$ derivatives are dealt with trivially as before).
\end{proof}
\begin{corollary}\label{counting:cor-sine-cancel}
Suppose $\varphi_j \in \{\sin, \cos\}$ for $3\leq j\leq 5$ and $N_0,N_2,\cdots,N_5,N_{234}$ are dyadic numbers, let $N_{\textup{max}} =\max (N_0, N_2, \cdots,N_5)$. Consider the Sine tensor $H^{\mathrm{sine}}_{n_0n_3n_4n_5}(t)$, which is defined as
\begin{equation}\label{def:sine_tensor3}
\begin{aligned}
H^{\mathrm{sine}}_{n_0 n_3n_4 n_5}(t,\lambda_3,\lambda_4, \lambda_5) &=\, \mathbf{1}\big\{ n_0 = n_{345}\big\} \cdot  \frac{1_{N_0}(n_0)}{\langle n_0 \rangle}\prod_{j=3}^5 \frac{1_{N_j}(n_j)}{\langle n_j \rangle} \int\displaylimits_0^{t} \dt^{\prime}\, \chi(t)\chi(t')\cdot \sin\big((t-t')\langle n_0\rangle\big)\, \varphi_5(t'\langle n_5\rangle)\\
&\times  \Big( \int\displaylimits_0^{t'} \dt^{\prime\prime}\, \chi(t') \chi(t'') \cdot \Sine[N_{234},N_2](t'-t^{\prime\prime},n_{34}) \, \varphi_3(t''\langle n_3\rangle)\varphi_4(t''\langle n_4\rangle) \Big).
\end{aligned}
\end{equation}
Then we have that
\begin{align}\label{counting:sine_bound3}
\| \langle\lambda \rangle^{b_+}\widetilde{H}_{n_0n_3n_4n_5}(\lambda)\|_{L^2_\lambda[n_0n_3n_4n_5]}&\lesssim N_{\textup{max}}^{-1/2+\epsilon},\\
 \label{counting:sine_bound4}
\| \langle\lambda \rangle^{b_+}\widetilde{H}_{n_0n_3n_4n_5}(\lambda)\|_{L^2_\lambda[n_0n_A\to n_Bn_5]}&\lesssim (N_0 N_5)^{-\frac{1}{2}}\cdot N_{\mathrm{max}}^{-1/2+\epsilon},\\
 \label{counting:sine_bound5}
\| \langle\lambda \rangle^{b_+}\widetilde{H}_{n_0n_3n_4n_5}(\lambda)\|_{L^2_\lambda[n_0n_5\to n_3n_4]}&\lesssim N_{\mathrm{max}}^{\epsilon}\cdot \big(\max (N_0, N_5)\cdot\max (N_3,N_4)\big)^{-\frac{1}{2}}\cdot N_2^{-1},
\end{align}
for any choice of $\varphi_j$ and any partition $(A,B)$ of $\{3,4\}$.
\end{corollary}
\begin{proof} This follows from Lemma \ref{counting:sine_tensor} in the same way that Corollary \ref{counting:cor-cubic} follows from Lemma \ref{counting:cubic_tensor}.
\end{proof}

\subsection{Tensor and $p$-moment estimates reductions}\label{prep:remark-reduction}  We start by explaining how a multilinear estimate in $X^{s,b}$-spaces can be reduced to a tensor estimate. This (standard) reduction will be used in the proof of several estimates (see e.g. Lemma \ref{bilinear:lem-bilinear}, Lemma \ref{linear:lem-quad-dyadic} and Proposition \ref{linear:prop-lincub}); but only explained in detail here.
To illustrate this reduction, we focus on a bilinear estimate 
\begin{equation}\label{reduction:eq1}
{X^{1/2,b}\times X^{1/2,b}\rightarrow X^{-1/2,b_+-1}},\,
(w_2,w_3) \mapsto P_{N_0}  \Big[ \slinear[blue][N_1] P_{N_2}w_2 P_{N_3}w_3 \Big].
\end{equation}
In a first step, we decompose $w_2$ and $w_3$ into a superposition of time-modulated linear waves and  insert the decomposition into the Duhamel integral. More precisely, we write
\begin{equation}\label{reduction:eq2}
w_j(t,x)=\sum_{\pm_j}\sum_{n_j\in \mathbb{Z}^3}\int_{\mathbb{R}}\mathrm{d}\lambda_j e^{i(\pm_j \langle n_j\rangle+\lambda_j)t} e^{i\langle n_j, x\rangle} \widetilde{w_j}^{\pm_j}(\lambda_j, n_j),
\end{equation}
where 
$
\max_{\pm_j} \|\langle\lambda_j\rangle^b \langle n_j\rangle^{1/2} \widetilde{w_j}^{\pm_j}(\lambda_j, n_j)\|_{L^2_{\lambda_j}\ell^2_{n_j}} \sim 
\|w_j\|_{X^{1/2, b}}.
$
In order to closely match (\ref{reduction:eq2}), we also write
\begin{equation}\label{reduction:eq3}
\slinear[blue][] =\sum_{\pm_1}\sum_{n_1\in \mathbb{Z}^3} 
\frac{g_{n_1}^{\pm_1}}{\langle n_1\rangle} e^{\pm_1 it\langle n_1\rangle} e^{i\langle n_1,x\rangle}.
\end{equation}
Equipped with \eqref{reduction:eq2} and \eqref{reduction:eq3}, we write
\begin{align}
\langle\nabla\rangle^{-\frac{1}{2}} P_{N_0}  \Big[ \slinear[blue][N_1] P_{N_2}w_2 P_{N_3}w_3 \Big]=&\sum_{\pm_1,\pm_2,\pm_3}\int_\mathbb{R}\mathrm{d}\lambda_2\int_\mathbb{R}\mathrm{d}\lambda_3 \sum_{n_0,n_1,n_2,n_3\in\mathbb{Z}^3}\Big[ \mathbf{1}\big\{n_0=n_{123}\big\}\cdot \bigg( \prod_{j=0}^3 1_{N_j} (n_j)\bigg) \nonumber\\
&\times \langle n_0\rangle^{-1/2} \langle n_1\rangle^{-1}\cdot e^{i(\pm_1\langle n_1\rangle\pm_2\langle n_2\rangle\pm_3\langle n_3\rangle +\lambda_2+\lambda_3)t} \cdot e^{i\langle n_0,x\rangle}\nonumber\\
&\times g_{n_1}^{\pm_1}\cdot\widetilde{w_2}^{\pm_2}(\lambda_2, n_2)\cdot\widetilde{w_3}^{\pm_3}(\lambda_3, n_3)\Big]\label{reduction:eq4}.
\end{align}
In a second step, we bring the expression in (\ref{reduction:eq4}) into a tensor form. To this end, we use the cubic tensor $h_{n_0n_1n_2n_3}(t,\lambda_2,\lambda_3)$ defined in above \eqref{def:cubic_tensor1} with $\lambda_1=0$ and then 
\eqref{reduction:eq4} can be written as
\begin{equation*}
\begin{aligned}
\langle\nabla\rangle^{-\frac{1}{2}} P_{N_0}  \Big[ \slinear[blue][N_1] P_{N_2}w_2 P_{N_3}w_3 \Big]=&\sum_{\pm_1,\pm_2,\pm_3}\int_\mathbb{R}\mathrm{d}\lambda_2\int_\mathbb{R}\mathrm{d}\lambda_3 \sum_{n_0,n_1,n_2,n_3\in\mathbb{Z}^3}\Big[ \langle n_0\rangle^{-1/2} \langle n_1\rangle^{-1} \\
&\times  h_{n_0n_1n_2n_3}(t,\lambda_2,\lambda_3) \cdot e^{i\langle n_0,x\rangle}\cdot g_{n_1}^{\pm_1}\cdot\widetilde{w_2}^{\pm_2}(\lambda_2, n_2)\cdot\widetilde{w_3}^{\pm_3}(\lambda_3, n_3)\Big].
\end{aligned}
\end{equation*}
Using $\widetilde{h}_{n_0n_1n_2n_3}(\lambda,\lambda_2,\lambda_3)$ as in \eqref{eq:htilde} with $\lambda_1=0$, $n_0=-n$ and $N_0=N_{123}$, except now that we the fix signs $\pm$ in $\Omega$ by $\pm_j\, (j=1,2,3)$, we obtain that
\begin{multline}\label{reduction:eq5}
\mathcal{F}\Big(\langle\nabla\rangle^{-\frac{1}{2}} P_{N_0}  \Big[ \slinear[blue][N_1] P_{N_2}w_2 P_{N_3}w_3 \Big]\Big)(\pm_0\langle n_0\rangle+\lambda, n_0)
=\sum_{\pm_0,\pm_1,\pm_2,\pm_3}\int_\mathbb{R}\mathrm{d}\lambda_2\int_\mathbb{R}\mathrm{d}\lambda_3 \sum_{n_0,n_1,n_2,n_3\in\mathbb{Z}^3}\Big[ \langle n_0\rangle^{-1/2} \langle n_1\rangle^{-1} \\
\times  \widetilde{h}_{n_0n_1n_2n_3}(\lambda,\lambda_2,\lambda_3) \cdot g_{n_1}^{\pm_1}\cdot\widetilde{w_2}^{\pm_2}(\lambda_2, n_2)\cdot\widetilde{w_3}^{\pm_3}(\lambda_3, n_3)\Big].
\end{multline}

In a third and final step, we now estimate the expression in \eqref{reduction:eq5}. We may fix the choices of $\pm_j\,(0\leq j\leq 3)$ as there are only $2^4$ choices, then we have
\begin{equation*}
\begin{aligned}
&\Bigg\| \langle \lambda\rangle^{b_+-1} \mathcal{F}\Big(\langle\nabla\rangle^{-\frac{1}{2}} P_{N_0}  \Big[ \slinear[blue][N_1] P_{N_2}w_2 P_{N_3}w_3 \Big]\Big)(\pm_0\langle n_0\rangle+\lambda, n_0)\Bigg\|_{L_\lambda^2\ell_{n_0}^2}\\
\leq &\int_{\mathbb{R}^2}\mathrm{d}\lambda_2\mathrm{d}\lambda_3 \Bigg\|\sum_{n_0,\cdots,n_3\in\mathbb{Z}^3}\Big[\langle \lambda\rangle^{b_+-1} \langle n_0\rangle^{-1/2} \langle n_1\rangle^{-1}\widetilde{h}_{n_0n_1n_2n_3}(\lambda,\lambda_2,\lambda_3) \cdot g_{n_1}^{\pm_1}\widetilde{w_2}^{\pm_2}(\lambda_2, n_2)\widetilde{w_3}^{\pm_3}(\lambda_3, n_3)\Big]\Bigg\|_{L_\lambda^2\ell_{n_0}^2}\\
\lesssim&\, N_1^{-1}(N_0N_2N_3)^{-\frac{1}{2}}\cdot\sup_{\lambda_2,\lambda_3\in\mathbb{R}} \Big[(\langle\lambda_2\rangle\langle \lambda_3\rangle)^{-(b_--\frac{1}{2})}\Big\|\sum_{n_1} \langle\lambda\rangle^{b_+-1}\, \widetilde{h}_{n_0n_1n_2n_3}(\lambda,\lambda_2,\lambda_3) g_{n_1}\Big\|_{\ell^2_{n_2}\times\ell^2_{n_3}\to L^2_\lambda\ell^2_{n_0}}
\Big]\\
&\times \int_{\mathbb{R}^2}\mathrm{d}\lambda_2\mathrm{d}\lambda_3 \prod_{j=2}^3\langle\lambda_j\rangle^{b_--\frac{1}{2}}\big\|\langle n_j\rangle^{\frac{1}{2}} |\widetilde{w}^{\pm_j}_j(\lambda_j,n_j)|\big\|_{\ell^2_{n_j }}\\
\lesssim&\,\sup_{\lambda_2,\lambda_3\in\mathbb{R}} \Big[(\langle\lambda_2\rangle\langle \lambda_3\rangle)^{-(b_--\frac{1}{2})}
\Big\|\sum_{n_1} \langle\lambda\rangle^{b_+-1}\, \widetilde{h}_{n_0n_1n_2n_3}(\lambda,\lambda_2,\lambda_3) g_{n_1}\Big\|_{\ell^2_{n_2}\times\ell^2_{n_3}\to L^2_\lambda\ell^2_{n_0}}
\Big]\\
&\times N_1^{-1}(N_0N_2N_3)^{-\frac{1}{2}}\cdot \prod_{j=2}^3 \|w_j\|_{X^{\frac{1}{2}, b}}.
\end{aligned}
\end{equation*}

In order to estimate the operator norm in \eqref{reduction:eq1}, it therefore suffices to control
\begin{equation}\label{reduction:goal}
\Big\|
(\langle\lambda_2\rangle\langle\lambda_3\rangle)^{-(b_--\frac{1}{2})}
\Big\|
\sum_{n_1} \langle\lambda\rangle^{b_+-1}\, \widetilde{h}_{n_0n_1n_2n_3}(\lambda,\lambda_2,\lambda_3) g_{n_1}
\Big\|_{L^2_\lambda(\ell^2_{n_2}\times \ell^2_{n_3}\to\ell^2_{n_0})}
\Big\|_{L_{\lambda_2,\lambda_3}^\infty}.
\end{equation} Note that we have now reduced the estimate for \eqref{reduction:eq1} to the estimate for (\ref{reduction:goal}), which involves the $\ell^2_{n_2}\times \ell^2_{n_3}\to\ell^2_{n_0}$ bilinear operator norm of certain tensors with parameters $(\lambda,\lambda_2,\lambda_3)$. 

\medskip

Next, we explain how to estimate the $L_\omega^p$ moments of suprema or integrals of random variables over a continuous set of parameters such as (\ref{reduction:goal}). This argument is equivalent to the \emph{meshing argument} in \cite{DNY20} but is more suitable for our purpose here as most of our estimates are stated in terms of $L_\omega^p$ moments instead of $A$-certain events like in \cite{DNY20}.

Take (\ref{reduction:goal}) as an example. The goal is to estimate (see for example Lemma \ref{bilinear:lem-bilinear})
\begin{equation}\label{meshing:goal}
\mathbb{E}\big(|\eqref{reduction:goal}|^p\big)^{1/p}\end{equation}
for all large enough $p$; here the main obstacle is the $L_{\lambda_2,\lambda_3}^\infty$ norm in (\ref{reduction:goal}) because Minkowski's inequality does not allow us to control the $L_\omega^pL_{\lambda_2,\lambda_3}^\infty$ norm by the $L_{\lambda_2,\lambda_3}^\infty L_\omega^p$ norm.

However, let $q:=(b_--\frac{1}{2})^{-5}$, then by Sobolev embedding, we have
\begin{equation}\label{counting:sob}\|F(\lambda,\lambda_2,\lambda_3)\|_{L_{\lambda_2,\lambda_3}^\infty}\lesssim \|F(\lambda, \lambda_2,\lambda_3)\|_{L_{\lambda_2,\lambda_3}^q}+\|\nabla_{\lambda_2,\lambda_3} F(\lambda, \lambda_2,\lambda_3)\|_{L_{\lambda_2,\lambda_3}^q}\end{equation}for any complex valued or Banach space valued function $F$. This basically allows us to control the $L_{\lambda_2,\lambda_3}^\infty$ norm by the $L_{\lambda_2,\lambda_3}^q$ norm with one extra derivative, but since any $(\lambda,\lambda_2,\lambda_3)$ derivative of $\widetilde{h}$ satisfies the same estimates as $\widetilde{h}$ itself, in fact we only need to treat the contribution corresponding to the first term on the right hand side of (\ref{counting:sob}). 

Now, assuming the above considerations, we take  $p>q$ and apply Minkowski inequality to get, 
\begin{align}\mathbb{E}\big(|\eqref{reduction:goal}|^p\big)^{1/p}&\lesssim 
\Big\| (\langle\lambda_2\rangle\langle\lambda_3\rangle)^{-(b_--\frac{1}{2})}\Big\|
\sum_{n_1}\langle \lambda\rangle^{b_+-1}\widetilde{h}_{n_0n_1n_2n_3}(\lambda,\lambda_2,\lambda_3) g_{n_1}
\Big\|_{\ell^2_{n_2}\times \ell^2_{n_3}\to \ell^2_{n_0}}
\Big\|_{L_\omega^pL_{\lambda_2,\lambda_3}^qL_\lambda^2}\nonumber\\ 
&\lesssim \Big\| (\langle\lambda_2\rangle\langle\lambda_3\rangle)^{-(b_--\frac{1}{2})}\Big\|
\sum_{n_1}\langle \lambda\rangle^{b_+-1}\widetilde{h}_{n_0n_1n_2n_3}(\lambda,\lambda_2,\lambda_3) g_{n_1}
\Big\|_{\ell^2_{n_2}\times \ell^2_{n_3}\to \ell^2_{n_0}}
\Big\|_{L_{\lambda_2,\lambda_3}^qL_\lambda^2L_\omega^p}\nonumber\\
&\lesssim  \sup_{\lambda_2,\lambda_3}\Big\|\langle \lambda\rangle^{b_+-1}\Big(\Eb\Big\|
\sum_{n_1}\widetilde{h}_{n_0n_1n_2n_3}(\lambda,\lambda_2,\lambda_3) g_{n_1}
\Big\|_{\ell^2_{n_2}\times \ell^2_{n_3}\to \ell^2_{n_0}}^p\Big)^{1/p}\Big\|_{L_\lambda^2}.\label{counting:mesh2}
\end{align}
Note that on the right hand side of (\ref{counting:mesh2}) we are taking the $L_\omega^p$ moments for \emph{fixed} $(\lambda,\lambda_2,\lambda_3)$, instead of taking $L_\omega^p$ moments of some \emph{supremum} in $(\lambda_2,\lambda_3)$. This is the reduction we need, since then by the level-set decomposition of $\widetilde{h}$  (as in the first part of the proof of Lemma \ref{counting:cubic_tensor}), the $L_\omega^p$ moment for \emph{fixed} $(\lambda,\lambda_2,\lambda_3)$ in the last line of \eqref{counting:mesh2} can be controlled by $A_1(\lambda,0,\lambda_2,\lambda_3)$ (defined in \eqref{cubic:A1}) times the same $L_\omega^p$ moment but with the base tensor $h^b$ replacing $\widetilde{h}$.  Since $\|\langle \lambda\rangle^{b_+-1}A_1(\lambda,0,\lambda_2,\lambda_3)\|_{L^2_{\lambda}} \lesssim N_{\max}^{\epsilon}$ by \eqref{counting:cubic_coef}, one is finally reduced to controlling the $L_\omega^p$ moment \[\mathbb{E}\Big[\Big\|
\sum_{n_1}  (h^b)_{n_0n_1n_2n_3}\, g_{n_1}
\Big\|^p_{\ell^2_{n_2}\times \ell^2_{n_3}\to \ell^2_{n_0}}\Big]^{1/p}.\label{counting:eq-final-reduction} \]
This completes the reductions of this subsection since proof of the tensor estimates such as \eqref{counting:eq-final-reduction} is the subject of other sections of this article  (see  Lemma \ref{bilinear:lem-bilinear} for an estimate of \eqref{counting:eq-final-reduction} and Lemma \ref{linear:lem-quad-dyadic} and Proposition \ref{linear:prop-lincub} for similar estimates).

\section{Algebraic and graphical aspects of stochastic diagrams}\label{section:diagrams}

In this section, we describe explicit stochastic objects, which are defined directly as Picard iterates of the nonlinear wave equation. In contrast to the para-controlled components and nonlinear remainder, the explicit stochastic objects do not require the solution of a fixed point problem.

In Subsection \ref{section:diagram-parabolic}, we discuss the Gibbs measure and the associated caloric initial data (from Proposition \ref{ansatz:prop-caloric}). To this end, we discuss the stochastic quantization of the Gibbs measure through the (stochastic) nonlinear heat equation, i.e., the parabolic $\Phi^4_3$-model (cf. \cite{CC18,GIP15,GH21,H14,HM18,MW17,MW20}). In this context, we also introduce our notation for the stochastic diagrams of the  nonlinear heat equation. In the remaining subsections, we describe the stochastic diagrams of the cubic nonlinear wave equations. As already mentioned in the introduction, we will utilize two different sets of graphical notations, which are called shorthand and longhand diagrams, respectively. The shorthand diagrams for the wave equation, which have already been described in \mbox{Section \ref{section:ansatz},}  will include \mbox{$\slinear[blue]$, $\scubic$, and $\squintic$}.  From a type-setting perspective, these diagrams are relatively small and will be used in our Ansatz and contraction arguments (see Section \ref{section:ansatz}). In order to keep the notation in our Ansatz manageable, however, the shorthand diagrams will actually contain sub-terms which are not reflected in the shorthand notation. For instance, as already seen in \eqref{ansatz:eq-quintic}, the quintic diagram $\squintic$ will also contain terms with one or two pairings. In contrast, the longhand diagrams contain all information on both pairings and renormalizations. For example, the longhand diagrams 
\begin{equation*}
\quintic[][][0], \quad \quintic[][][1], \quad \text{and} \quad \quintic[][][2].
\end{equation*}
describe the terms in $\squintic$ with zero, one, or two pairings. From a type-setting perspective, it is difficult to include the longhand diagrams in mathematical formulas. As a result, they have not been used in our Ansatz (Section \ref{section:ansatz}). However, the longhand diagrams are incredibly useful in precise estimates of higher-order iterates. Without this notation, it would be much harder to exhibit the \oftt-cancellation (Proposition \ref{diagram:prop-1533-cancellation}), which is the main result of this section.

\subsection{Parabolic stochastic quantization}\label{section:diagram-parabolic} 

In this subsection, we discuss the properties of the caloric initial data (as in Proposition \ref{ansatz:prop-caloric}) and the Gibbs measure. 
As described in the introduction, the Gibbs measure is given by the product of the $\Phi^4_3$-measure and the Gaussian free field $\cg$. Both measures can be viewed as invariant measures for nonlinear and linear stochastic heat equations, respectively. Since the stochastic heat equations will be used to represent the initial position and velocity of the wave equation, which are propagated through the Fourier-multipliers  $\cos(t\langle \nabla \rangle)$ and $\sin(t\langle \nabla \rangle)$, we label the two stochastic heat equations using the super-scripts ``$\cos$" and ``$\sin$". 

In the following, we denote the time-variable in the heat equations by $s\in \R$, since the variable $t$ is reserved for the time-variable of nonlinear wave equations. 

\subsubsection{The linear and nonlinear stochastic heat equations}
We let $(\Omega,\mathcal{E},\mathbb{P})$ be the probability space from Subsection \ref{section:prelim-multiple-stochastic} and let 
let $(\Wp[\sin][s][n])_{n\in \Z^3}$ be an independent copy of the Gaussian process $(\Wp[][s][n])_{n\in \Z^3}$ from Subsection \ref{section:prelim-multiple-stochastic}. 
We  also define the Gaussian process
\begin{equation}
\Wp[\sin][s][x] := \sum_{n\in \Z^3} \Wp[\sin][s][n] e^{inx}. 
\end{equation}

The It\^{o}-derivative $\dW[sin][s][x]$ yields a real-valued space-time white noise. The (ancient) linear stochastic heat equation with space-time white noise forcing is then given by 
\begin{equation}\label{diagram:eq-SLH}
\big(\partial_s+ 1 - \Delta \big) \phi^{\sin} = \sqrt{2}~ \dW[\sin] \quad (s,x) \in \R \times \T^3.\\ 
\end{equation}
The reason for the $\sqrt{2}$-factor lies in Lemma \ref{diagram:lemma-linear-heat-invariance} below. In order to solve \eqref{diagram:eq-SLH}, we introduce the heat propagator 
$(e^{-s (1-\Delta)})_{s\geq 0}$. On the torus $\T^3$, the heat propagator is given by a Fourier multiplier with symbol $\exp(-s\langle n \rangle^2)$, i.e., it is defined by
\begin{equation}\label{diagram:eq-Hlin}
\mathcal{F}\Big(e^{-s (1-\Delta)}f\Big)(n) := e^{-s\langle n \rangle^2} \widehat{f}(n).
\end{equation}
Then, a stationary solution{\footnote{That is, the solution is a stationary process, its distribution constant in time.}} of \eqref{diagram:eq-SLH} is given explicitly by
\begin{equation}\label{diagram:eq-phi-sin}
    \Phi^{\sin}(s,x) = \sqrt{2} \sum_{n\in \Z^3} \bigg[
    \Big( \int_{-\infty}^s e^{-(s-\tau) \langle n \rangle^2} \dW[sin][\tau][n]\Big) e^{inx} \bigg]. 
\end{equation}
It is easy to see that \eqref{diagram:eq-phi-sin} is well-defined as a space-time distribution. For us, the most important property of $\Phi^{\sin}(s,x)$ is its law for any fixed time $s \in \R$.

\begin{lemma}[Law of $\Phi^{\sin}$] \label{diagram:lemma-linear-heat-invariance}
For all $s\in \R$, it holds that $\Law_{\mathbb{P}}(\Phi^{\sin}(s,\cdot))= \cg$. 
\end{lemma}

\begin{proof} Let $s \in \R$ be arbitrary but fixed. Then, we can write $\Phi^{\sin}(s,x)= \sum_{n\in \Z^3} \frac{g_n}{\langle n \rangle} e^{i\langle n, x \rangle}$, where 
\begin{equation*}
g_n := \sqrt{2} \langle n \rangle^2 \int_{-\infty}^{s} e^{-(s-\tau)\langle n \rangle^2} \dW[sin][\tau][n].
\end{equation*}
Using the properties of $\Wp[sin]$, it is easy to see that 
\begin{enumerate}[label={(\roman*)},leftmargin=7mm]
\item $g_0$ is a real-valued Gaussian and $g_n$ is a complex-valued Gaussian with independent and identically distributed real and imaginary parts for all $n\neq 0$,
\item $g_n$ and $g_m$ are independent as long as $n\neq \pm m$,
\item and $\overline{g_n}=g_{-n}$.
\end{enumerate}
Thus, it only remains to prove that $\E[|g_n|^2]=1$. Indeed, it follows from It\^{o}'s isometry that 
\begin{equation*}
\E[|g_n|^2] = 2 \langle n \rangle^2 \int_{-\infty}^{s} e^{-2 (s-\tau)\langle n \rangle^2} = 1. \qedhere
\end{equation*}
\end{proof}

This completes our analysis of the linear stochastic heat equation \eqref{diagram:eq-SLH} and we now turn to the nonlinear stochastic heat equation. To this end, we let $\Wp[\cos][s](x)$ be an independent copy of the Gaussian process $\Wp[\sin][s][x]$, which was introduced above. 
Furthermore, we let $s_0 \in (-\infty,0)$ be an initial time. The upper bound $s_0<0$ stems from Definition \ref{diagram:definition-caloric} below, where we evaluate the nonlinear heat equation at $s=0$. On a formal level, the nonlinear stochastic heat equation is then given by 
\begin{equation}\label{diagram:eq-SNLH}
\begin{cases}
\big(\partial_s+ 1 - \Delta \big) \Phi^{\cos} = - (\Phi^{\cos})^3 - \infty\cdot \Phi^{\cos} + \sqrt{2} ~\dW[\cos]  \qquad (s,x) \in (s_0,\infty) \times \T^3,\\ 
\Phi^{\cos}|_{s=s_0} = \phi^{\cos}. 
\end{cases}
\end{equation}
Here, the term ``$\infty \cdot \Phi^{\cos}\,$'' represents a renormalization. In order to give a rigorous meaning to the renormalized equation \eqref{diagram:eq-SNLH}, we introduce the frequency-truncated nonlinear stochastic heat equation
\begin{equation}\label{diagram:eq-truncated-SNLH}
\begin{cases}
\big(\partial_s + 1 - \Delta \big) \Phi_{\leq N}^{\cos}  = -  P_{\leq N}  \Big( \lcol \big(P_{\leq N} \Phi_{\leq N}^{\cos} \big)^3 \rcol  + \gamma_{\leq N}  \Phi^{\cos}_{\leq N} \Big)+ \sqrt{2} \dW[\cos] \quad  (s,x) \in (s_0,\infty) \times \T^3,\\
\Phi^{\cos}_{\leq N} \big|_{s=s_0}= \phi^{\cos}. 
\end{cases}
\end{equation}

Here, the term $\lcol \big(P_{\leq N} \Phi^{\cos} \big)^3\rcol$ is the Wick-ordered cubic nonlinearity.  The additional (logarithmic) renormalization constant $\gamma_{\leq N}$ and multiplier $\Gamma_{\leq N}$, which were already used in Section \ref{section:ansatz}, and are defined in the following definition. 
\begin{definition}[The renormalization constant and multipier]\label{sec6:def:Gamma}
For all $N\geq 1$, we define the renormalization constant $\gamma_{\leq N}$ and renormalization multiplier $\Gamma_{\leq N}$ by 
\begin{equation}\label{diagram:eq-Gamma-heat}
\begin{aligned}
\gamma_{\leq N} &= \Gamma_{\leq N}(0), \\ 
\Gamma_{\leq N}(n) &=   6 \cdot 1_{\leq N}(n) \sum_{\substack{m_1,m_2,m_3 \in \Z^3\colon \\ m_{123}=n}}   \Big( \prod_{j=1}^{3} 1_{\leq N}(m_j) \langle m_j \rangle^{-2} \Big). 
\end{aligned}
\end{equation}
\end{definition}
We emphasize that the nonlinearity in \eqref{diagram:eq-truncated-SNLH} coincides with the nonlinearity in the frequency-truncated wave equation \eqref{diagram:eq-wave-truncated} below (and encountered in \eqref{ansatz:eq-frequency-truncated-NLW} above). This is important for our analysis because otherwise the (frequency-truncated) Gibbs measures would not agree. 
 We now define the energy function 
\begin{equation}\label{diagram:eq-heat-energy}
    E_{\leq N}(\phi) := \frac{1}{2} \int_{\T^3} \dx \Big( \phi(x)^2 + |\nabla \phi(x)|^2 \Big) + \frac{1}{4} \int_{\T^3} \dx \lcol (P_{\leq N} \phi)^4 (x) \rcol +  \frac{1}{2}  \gamma_{\leq N} \int_{\T^3} \dx \,  (P_{\leq N} \phi)^2(x). 
\end{equation}
For future use, we also isolate the potential energy $V_{\leq N}$ in \eqref{diagram:eq-heat-energy}, which is defined as 
\begin{equation}
     V_{\leq N}(\phi) :=
     \frac{1}{4} \int_{\T^3} \dx \lcol (P_{\leq N} \phi)^4 (x) \rcol +\frac{1}{2}  \gamma_{\leq N} \int_{\T^3} \dx \,  (P_{\leq N} \phi)^2(x).
\end{equation}
Using the definition of $E_{\leq N}$, we can rewrite the nonlinear stochastic heat equation \eqref{diagram:eq-truncated-SNLH} as the Langevin equation corresponding to $E_{\leq N}$. More precisely,  \eqref{diagram:eq-truncated-SNLH} can be written as 
\begin{equation}\label{diagram:eq-truncated-Langevin}
\begin{cases}
\partial_s \Phi_{\leq N}^{\cos}  = -  \nabla_{L^2} E_{\leq N}(\Phi_{\leq N}^{\cos}) + \sqrt{2} \dW[\cos] \quad  (t,x) \in (t_0,\infty) \times \T^3,\\
\Phi^{\cos}_{\leq N} \big|_{t=t_0}= \phi^{\cos}. 
\end{cases}
\end{equation}
Here, $\nabla_{L^2}$ refers to the gradient with respect to the $L^2$-metric.
In addition to the evolution equation \eqref{diagram:eq-truncated-Langevin}, the energy $E_{\leq N}$ yields a frequency-truncated version of the $\Phi^4_3$-measure. On a formal level, it is given by the expression
\begin{equation*} 
    \mathrm{d}\Phi^4_{3;\leq N}(\phi) = \text{``} \mathcal{Z}^{-1}_{\leq N} \exp\big( - E_{\leq N}(\phi) \big) \mathrm{d}\phi \text{''}. 
\end{equation*}
Despite issues related to the infinite-dimensional Lebesgue measure $\mathrm{d}\phi$, which can not be given a rigorous meaning, we can give a rigorous definition of $\Phi^4_{3;\leq N}$ based on the Gaussian free field and the potential energy $V_{\leq N}$.

\begin{definition}[The frequency-truncated $\Phi^4_3$-measure]
For any dyadic scale $N$, we define
\begin{equation*}
    \mathrm{d}\Phi^4_{3;\leq N}(\phi) := \mathcal{Z}^{-1}_{\leq N} \exp\big( - V_{\leq N}(\phi) \big) \mathrm{d}\cg(\phi).
\end{equation*}
\end{definition}

The following lemma follows directly from standard methods for stochastic differential equations (see e.g. \cite{KS91,O03}).

\begin{lemma}[Invariance of the frequency-truncated $\Phi^4_3$-measure]\label{diagram:lemma-nonlinear-heat-invariance}
For any $N\geq 1$, the frequency-truncated nonlinear stochastic heat equation \eqref{diagram:eq-truncated-SNLH} is globally well-posed. In addition, the $\Phi^4_{3;\leq N}$-measure is invariant under the flow. 
\end{lemma}

While Lemma \ref{diagram:lemma-nonlinear-heat-invariance} has important consequences for a fixed frequency-scale $N$, it does not yield any information regarding the limiting behavior as $N\rightarrow \infty$. Before proceeding with the analysis of this limit, we need to introduce our graphical notation for the stochastic heat equation. 

\subsubsection{Graphical notation for stochastic heat equations}\label{section:graphical-heat}
We now introduce our graphical notation and start with the linear stochastic heat equation \eqref{diagram:eq-SLH}. We use a black dot to represent the stochastic forcing term in \eqref{diagram:eq-SLH}. More precisely, we write
\begin{equation}
\initial{black}[sin] \overset{\text{def}}{=} \sqrt{2} ~  \dW[\sin]. 
\end{equation}
Utilizing the dot notation, we represent the (stationary) solution $\Phi^{\sin}(s,x)$ as 
\begin{equation}\label{diagram:eq-heat-sin-x}
\heat[][sin][][(s,x)] := \Phi^{\sin}(s,x) = \sqrt{2} \bigg( \int\displaylimits_{-\infty}^s e^{-(s-\tau)(1-\Delta)} ~ \dW[\sin][\tau] \bigg)(x). 
\end{equation}
The densely dotted arrow represents the Duhamel integral for the linear heat equation. The orange dot together with its subscript represents an evaluation in space-time at $(s,x)$. In the following analysis, it is often more convenient to work with the spatial Fourier transform of $\phi^{\sin}$, which is represented as 
\begin{equation}\label{diagram:eq-heat-sin-n}
\heat[][sin][][(s,n)]:=  \mathcal{F}\big( \Phi^{\sin}(s,\cdot)\big)(n) = \sqrt{2} \int\displaylimits_{-\infty}^s e^{-(s-\tau) \langle n \rangle^2} \dW[\sin][\tau][n].
\end{equation}
We emphasize that our change from physical to frequency space is only marked by the subscript of the orange dot, which requires no additional components in our graphical notation. We now turn to the nonlinear stochastic heat equation \eqref{diagram:eq-truncated-SNLH}. While the solution of the linear equation \eqref{diagram:eq-SLH} is  exactly given by the stochastic object in \eqref{diagram:eq-heat-sin-n}, the solution of the nonlinear equation \eqref{diagram:eq-truncated-SNLH} will only be approximated by stochastic objects. Similar as above, the stochastic forcing term is represented using a black dot and we write  
\begin{equation}
\initial{black}[cos] \overset{\text{def}}{=} \sqrt{2} ~  \dW[\cos]. 
\end{equation}
Using our previous graphical representations for the Duhamel integral and evaluations, we can represent the zero-th Picard iterate by 
\begin{equation}
\heat[][cos][][(s,n)]:= \sqrt{2} \int\displaylimits_{-\infty}^s e^{-(s-\tau) \langle n \rangle^2} \dW[\cos][\tau][n].
\end{equation}
In preparation for Subsection \ref{section:quintic-diagram} below, we prove the following covariance identity.

\begin{lemma}[Covariance identity for the linear heat evolution]\label{diagram:lemma-covariance-heat}
For all $\varphi,\varphi^\prime \in \{ \cos, \sin\}$, $s,s^\prime \in \R$, and $n,n^\prime \in \Z^3$, it holds that 
\begin{equation}
\E\left[ \heat[][\varphi][][(s,n)] \heat[][\varphi^\prime][][(s^\prime,n^\prime)] \right] = \delta_{n+n^\prime=0}\, \delta_{\varphi=\varphi^\prime} \,\frac{1}{\langle n \rangle^2} e^{-|s-s^\prime| \langle n\rangle^2}.
\end{equation}
\end{lemma}

The argument is a minor variant of the proof of Lemma \ref{diagram:lemma-linear-heat-invariance}. 

\begin{proof}
Using It\^{o}'s isometry, it holds that 
\begin{align*}
\E\left[ \heat[][\varphi][][(s,n)] \heat[][\varphi^\prime][][(s^\prime,n^\prime)] \right] 
&= 2  \E \bigg[ \bigg( \int_{-\infty}^s \dW[\varphi][\tau][n] e^{-(s-\tau) \langle n \rangle^2} \bigg) 
\cdot \bigg( \int_{-\infty}^{s^\prime} \dW[\varphi^\prime][\tau][n^\prime] e^{-(s^\prime-\tau) \langle n^\prime \rangle^2} \bigg)\bigg] \\
&= 2  \delta_{n+n^\prime=0}\, \delta_{\varphi=\varphi^\prime}  \int_{-\infty}^{\min(s,s^\prime)} e^{-(s-\tau) \langle n \rangle^2}e^{-(s^\prime-\tau) \langle n \rangle^2} \\
&= \delta_{n+n^\prime=0}\, \delta_{\varphi=\varphi^\prime} \,\frac{1}{\langle n \rangle^2} e^{-|s-s^\prime| \langle n\rangle^2}.
\end{align*}
This completes the proof. 
\end{proof}

Since the nonlinear stochastic heat equation \eqref{diagram:eq-truncated-SNLH} contains the Littlewood-Paley operators $P_{\leq N}$, it is convenient to include frequency-truncations directly inside the stochastic diagrams. To this end, we define
\begin{equation}\label{diagram:heat-linear-cos}
\heat[][cos][\leqN][] := P_{\leq N} \heat[][cos][][] \quad. 
\end{equation}
In the cubic Picard iterate of \eqref{diagram:eq-truncated-SNLH}, we only take into account the Wick-ordered cubic nonlinearity, since the $\gamma_{\leq N}$-term is designed to cancel double-resonances, which only occur at higher orders. Then, we write 
\begin{equation}\label{diagram:eq-heat-cubic}
\begin{aligned}
\cubicheat[][\leqN][(s,x)] &:= \int\displaylimits_{-\infty}^s e^{-(s-\tau) (1-\Delta)} P_{\leq N} \left( \biglcol\Bigg(  ~ \heat[][cos][\leqN][(\tau)]\Bigg)^3 \bigrcol \right)(x) ~\mathrm{d}\tau \quad. 
\end{aligned}
\end{equation}

Here, the multiplication symbol at the joint in \eqref{diagram:eq-heat-cubic} represents the Wick-ordered cubic nonlinearity. We emphasize that, just like the linear stochastic object, the cubic stochastic object is stationary in time. 
In frequency-space, the cubic Picard iterate of the nonlinear heat equation is given by 
\begin{equation}\label{diagram:eq-heat-cubic-frequency}
\begin{aligned}
\cubicheat[][\leqN][(s,n)] 
&= 2^{\frac{3}{2}} 1_{\leq N}(n) \sum_{\substack{n_1,n_2,n_3\in \Z^3 \colon \\ n_{123}=n}} \bigg[ \Big( \prod_{j=1}^3 1_{\leq N}(n_j) \Big) \\
&\times  
\int\displaylimits_{[-\infty,s]^3} \bigotimes_{j=1}^3 \dW[\cos][s_j][n_j]   \Big( \hspace{-1ex} \int\displaylimits_{\max(s_1,s_2,s_3)}^{s} \hspace{-2ex} \mathrm{d}\tau ~ e^{-(s-\tau) \langle n \rangle^2} \prod_{j=1}^3 e^{-(\tau-s_j) \langle n_j \rangle^2} \Big)  \bigg] . 
\end{aligned}
\end{equation}

In order to concisely state the evolution equation for the nonlinear remainder, we now introduce shorthand diagrams for the linear and cubic longhand diagrams in \eqref{diagram:heat-linear-cos} and \eqref{diagram:eq-heat-cubic}. We let
\begin{equation}\label{diagram:eq-heat-shorthand}
\shlinear[\leqN][cos] (s,x) = \heat[][cos][\leqN][(s,x)] 
\qquad \text{and} \qquad
\shcubic[\leqN][cos] (s,x) = \cubicheat[][\leqN][(s,x)]. 
\end{equation}

In order to solve \eqref{diagram:eq-truncated-SNLH}, we now make the Ansatz
\begin{equation}\label{diagram:eq-heat-Ansatz}
\Phi_{\leq N}^{\cos}(s,x) = \shlinear[][cos](s,x) - \shcubic[\leqN][cos](s,x) + \Psi^{\cos}_{\leq N}(s,x). 
\end{equation}
Then, the nonlinear remainder $ \Psi^{\cos}_{\leq N}(s,x)$ is a solution to the forced nonlinear stochastic heat equation
\begin{equation}\label{diagram:eq-heat-N-forced}
\begin{aligned}
\big(\partial_s + 1 - \Delta \big) \Psi^{\cos}_{\leq N} &= - P_{\leq N} \left[ 
\lcol   \left( \shlinear[\leqN][cos] - \shcubic[\leqN][cos] + P_{\leq N} \Psi^{\cos}_{\leq N}   \right)^3 \rcol 
-\lcol \left( \shlinear[\leqN][cos] \right)^3 \rcol  \right.\\
&+\left. \gamma_{\leq N}   \left( \shlinear[\leqN][cos] - \shcubic[\leqN][cos] + P_{\leq N} \Psi^{\cos}_{\leq N}  \right) \right]
 \end{aligned}
 \end{equation}
and has initial data 
\begin{equation}
\Psi^{\cos}_{\leq N} \big|_{s=s_0} = \phi^{\cos}  - \shlinear[][cos](s_0) + \shcubic[\leqN][cos](s_0). 
\end{equation}
Here, $\phi^{\cos}$ is the initial data from \eqref{diagram:eq-truncated-SNLH}

As is well-known, \eqref{diagram:eq-heat-N-forced} cannot be solved (uniformly in $N$) through a direct contraction argument. The properties of the caloric initial data (as in Proposition \ref{ansatz:prop-caloric}), which will be defined momentarily, therefore can only be derived from recent advances in singular SPDEs \cite{CC18,GIP15,H14}.

\subsubsection{Caloric initial data} 
We now use the solutions of the stochastic heat equations to define the caloric initial data. 
For this, we need to evaluate the solutions at a single time, which we have chosen as $s=0$. Of course, this choice is purely notational, and any other fixed time (such as $s=1$) would also have worked. Before we define the caloric initial data, we introduce additional notation. We define the combination of sine and cosine-based linear stochastic objects by 
\begin{equation}\label{gau:proc}
\initial{black} := \Big( \initial{black}[cos], \initial{black}[sin]\Big) \qquad \text{and} \qquad \heat = \left( \heat[][cos], \heat[][sin] \right). 
\end{equation}
Second, we define the evolution at $s=0$ by 
\begin{equation}
\evaluation[box][(x)] = \evaluation[][(0,x)] \qquad \text{and} \qquad \evaluation[box][(n)] = \evaluation[][(0,n)].
\end{equation} 
In our analysis of the wave Picard iterates, the ``box'' will primarily appear as an internal node (rather than at the bottom) of our stochastic objects (see e.g. \eqref{diagram:eq-quintic-mixed}), which is our reason for choosing black over orange. 
We now define the blue and green caloric initial data.

\begin{definition}[Blue and green caloric initial data]\label{diagram:definition-caloric}
For all $N\geq 1$, we define the linear and cubic Gaussian chaoses by 
\begin{equation}\label{bluegreen}
\initial{blue}(x) := \heat[box][][][(x)] \qquad \text{and} \qquad \initial{green}[][\leqN] (x) := \left( \cubicheat[box][\leqN][(x)], ~ 0 \right). 
\end{equation}
\end{definition}

In other words, the blue and green caloric initial data are obtained by evaluating the linear and cubic stochastic objects of the heat equation at time $s=0$. Since the linear and cubic stochastic objects are defined explicitly in terms of $\Wp[\cos][][]$ and $\Wp[\sin][][]$, but do not depend on the initial data $\phi^{\cos}$ in \eqref{diagram:eq-truncated-Langevin}, the same holds for the blue and green caloric initial data.

 The elementary regularity properties of $\initial{blue}$ and $\initial{green}[][\leqN]$ are recorded in the next lemma. However, the regularity information alone is far from sufficient for our analysis, and more intricate estimates will be proved throughout the rest of the article. 

\begin{lemma}[Regularity of blue and green caloric initial data]\label{diagram:lem-initial}
For all $p\geq 2$, it holds that 
\begin{equation}\label{diagram:eq-initial-1}
\E \Big[ \big\| \, \initial{blue}[][] \big\|_{\C_x^{-1/2-\epsilon} \times \C_x^{-1/2-\epsilon}}^p \Big]^{1/p} \lesssim p^{1/2}
\end{equation}
and 
\begin{equation}\label{diagram:eq-initial-2}
\E \Big[  \sup_{N} \big\| \, \initial{green}[][\leqN] \big\|_{\C_x^{1/2-\epsilon} \times \C_x^{1/2-\epsilon} }^p \Big]^{1/p} \lesssim p^{3/2}.
\end{equation}
\end{lemma}

\begin{remark}
As previously described in \eqref{ansatz:eq-initial-convention}, our initial conditions for the nonlinear wave equation \eqref{ansatz:eq-frequency-truncated-NLW-raw} are phrased in terms of $(u,\langle \nabla \rangle^{-1} \partial_t u)|_{t=0}$, which explains why the cosine and sine-components in Lemma \ref{diagram:lem-initial} are controlled at the same regularity.
\end{remark}

\begin{proof}
From Lemma \ref{diagram:lemma-linear-heat-invariance}, which also applies to the linear evolution in $\Phi^{\cos}$, we obtain that 
\begin{equation*}
\Law_{\mathbb{P}} \big( \, \initial{blue}[cos] \big) = \Law_{\mathbb{P}} \big( \, \initial{blue}[sin] \big) = \cg. 
\end{equation*}
For any $\varphi \in \{\cos,\sin\}$, it then follows from Gaussian hypercontractivity that
\begin{equation*}
\E \Big[ \big\|  \, \initial{blue}[\varphi] \big\|_{\C_x^{-1/2-\epsilon}}^p\Big]^{2/p}
\lesssim p \sum_{n \in \Z^3} \langle n \rangle^{-3+\epsilon} \lesssim p. 
\end{equation*}
This completes the proof of \eqref{diagram:eq-initial-1} and it remains to prove \eqref{diagram:eq-initial-2}. From Definition \ref{diagram:definition-caloric} and \eqref{diagram:eq-heat-cubic-frequency}, we obtain that $\initial{green}[sin][\leqN]=0$ and 
\begin{align*}
\initial{green}[cos][\leqN](x) 
=&  2^{\frac{3}{2}} \sum_{N_0,N_1,N_2,N_3 \leq N} \sum_{\substack{n_0,n_1,n_2,n_3\in \Z^3 \colon \\ n_{0}=n_{123}}} \bigg[ \Big( \prod_{j=0}^3 1_{N_j}(n_j) \Big) \\
\times&  
\int\displaylimits_{[-\infty,0]^3} \bigotimes_{j=1}^3 \dW[\cos][s_j][n_j]   \Big( \hspace{-1ex} \int\displaylimits_{\max(s_1,s_2,s_3)}^{0} \hspace{-2ex} \mathrm{d}\tau ~ e^{\tau \langle n_0 \rangle^2} \prod_{j=1}^3 e^{-(\tau-s_j) \langle n_j \rangle^2} \Big)  e^{i\langle n_0, x \rangle} \bigg].
\end{align*}
To simplify the notation, we set $N_{\textup{max}}=\max(N_0,N_1,N_2,N_3)$. 
Using Gaussian hypercontractivity, we obtain that 
\begin{align*}
&\E \bigg[ \,  \bigg\| \sum_{\substack{n_0,n_1,n_2,n_3\in \Z^3 \colon \\ n_{0}=n_{123}}} \bigg[ \Big( \prod_{j=0}^3 1_{N_j}(n_j) \Big) \\
&\hspace{10ex}\times   
\int\displaylimits_{[-\infty,0]^3} \bigotimes_{j=1}^3 \dW[\cos][s_j][n_j]   \Big( \hspace{-1ex} \int\displaylimits_{\max(s_1,s_2,s_3)}^{0} \hspace{-2ex} \mathrm{d}\tau ~ e^{\tau \langle n \rangle^2} \prod_{j=1}^3 e^{-(\tau-s_j) \langle n_j \rangle^2} \Big)  e^{i\langle n_0, x \rangle} \bigg] \bigg\|_{\C_x^{1/2-\epsilon}}^p \bigg]^{2/p} \\
\lesssim& \,p^3  N_0^{1-\epsilon} 
\sum_{\substack{n_0,n_1,n_2,n_3\in \Z^3 \colon \\ n_{0}=n_{123}}} \bigg[ \Big( \prod_{j=0}^3 1_{N_j}(n_j) \Big) 
\int_{[-\infty,0]^3} \ds_1 \ds_2 \ds_3 \bigg( \int\displaylimits_{\max(s_1,s_2,s_3)}^0 \dtau
e^{\tau \langle n_0 \rangle^2} \prod_{j=1}^3 e^{-(\tau-s_j) \langle n_j \rangle^2} \bigg)^2  \bigg] \\
\lesssim& \, p^3 N_0^{1-\epsilon} N_{\textup{max}}^{-2} \Big( \prod_{j=0}^3 N_j^{-2} \Big) \sum_{\substack{n_0,n_1,n_2,n_3\in \Z^3 \colon \\ n_{0}=n_{123}}}  \Big( \prod_{j=0}^3 1_{N_j}(n_j) \Big) \\
\lesssim& \, p^3 N_{\textup{max}}^{-\epsilon}. 
\end{align*}
This frequency-localized estimate directly implies \eqref{diagram:eq-initial-2}.
\end{proof}

It only remains to define the red caloric initial data. To this end, we let $A\geq 1$
be a parameter, which dictates the bound on $\initial{red}[][\leqN](A,\phi^{\cos})$. We now choose the initial time as 
\begin{equation}\label{diagram:eq-s0}
s_0 = s_0(A) := - A^{-C} <0,
\end{equation}
where $C\gg 1$ is any large\footnote{
The value of $C$ is irrelevant for our argument and essentially arbitrary. As explained in the proof of Proposition \ref{ansatz:prop-caloric}, the choice of $s_0=s_0(A)$ is linked to the local theory of the nonlinear stochastic heat equation \cite{CC18,GIP15,H14}.
One can require all (parabolic) random objects to be of size $A^k$ in their respective norms, where $k$ denotes the degree of the random object. Then, the local theory dictates that $C$ is sufficiently large, say, $C=10$, but yields local well-posedness on events with probability $\sim \exp(-A^2)$. Alternatively, one can require that all random objects in the local theory are of size $A^\delta$, where $0<\delta\ll 1$. Then, it is possible to take $C=1$, but the event now only has probability $\sim \exp(-A^c)$, where $0<c=c(\delta)\ll 1$.} constant.

\begin{definition}[The red caloric initial data]\label{diagram:def-red}
For all $A\geq 1$, $N\geq 1$, and $\phi^{\cos} \in \C_x^{-1/2-\epsilon}$, we define 
\begin{equation*}
\initial{red}[][\leqN] (x;A,\phi^{\cos}):= \Big( \Psi^{\cos}_{\leq N}(0,x;s_0,\phi^{\cos}), 0 \Big), 
\end{equation*}
where $\Psi^{\cos}_{\leq N}$ is the solution of \eqref{diagram:eq-heat-N-forced} with initial time $s_0=s_0(A)$ and initial data $\phi^{\cos}$. 
\end{definition}

We note that, since $s$ has been set to zero, $\initial{red}[][\leqN] (x;A,\phi^{\cos})$ is a function of only the spatial variable $x\in \T^3$. The red caloric initial data also depends on the parameter $A\geq 1$ (through the initial time $s_0$) and on the initial data $\phi^{\cos}$ from \eqref{diagram:eq-truncated-Langevin}.

 It remains to  prove Proposition \ref{ansatz:prop-caloric}, which contains the most relevant properties of the caloric initial data.

\begin{proof}[Proof of Proposition \ref{ansatz:prop-caloric}:]
By definition, it holds that 
\begin{equation*}
\initial{blue} - \initial{green}[][\leqM]  + \initial{red}[][\leqM](A,\phi^{\cos}_{\leq M}) = \Big( \Phi^{\cos}_{\leq M}(0,x), \Phi^{\sin}(0,x) \Big),
\end{equation*}
where $\Phi^{\cos}_{\leq M}$ is a solution of \eqref{diagram:eq-truncated-SNLH} with frequency-truncation parameter $M$, initial time $s_0=s_0(A)$, and initial data $\phi^{\cos}_{\leq M}$ and $\Phi^{\sin}$ is as in \eqref{diagram:eq-phi-sin}. By using Lemma \ref{diagram:lemma-nonlinear-heat-invariance} and Lemma \ref{diagram:lemma-linear-heat-invariance}, it follows that 
\begin{equation*}
    \Law_{\mathbb{P}\otimes \mathbb{Q}} \Big( \initial{blue} - \initial{green}[][\leqM]  + \initial{red}[][\leqM](A,\phi^{\cos}_{\leq M}) \Big) 
    = \Law_{\mathbb{P} \otimes \mathbb{Q}} \Big( \big( \Phi^{\cos}_{\leq M}(0,x), \Phi^{\sin}(0,x) \big) \Big) = \mu_{\leq M}. 
\end{equation*}
The properties \ref{ansatz:item-caloric-1} and \ref{ansatz:item-caloric-2} simply refer to the definitions of $\initial{blue}$ and $\initial{green}[][\leqM]$. As a result, the only non-trivial part is the estimate of $\initial{red}[][\leqM](A,\phi^{\cos}_{\leq M})$. This estimate essentially follows from the established well-posedness theory  for the nonlinear stochastic heat equation \cite{CC18,GIP15}, which represents $\Psi^{\cos}_{\leq M}$ as a para-controlled component at spatial regularity $1-$ and a smooth nonlinear remainder at spatial regularity $3/2-$. The only technical difference lies in the sharp Fourier-cutoffs in \eqref{diagram:eq-truncated-SNLH} and the necessary modifications are sketched in the appendix (Appendix \ref{section:heat-appendix}). 
\end{proof}

\subsection{Linear and cubic diagrams}\label{section:diagram-wave}

In this subsection, we study the linear and cubic stochastic diagrams for the nonlinear wave equation. As explained at the beginning of this section, we will introduce both the shorthand and longhand diagrams.

\subsubsection{The linear evolution} We first recall that the linear wave equation is given by
\begin{equation}\label{diagram:eq-wave-linear-deterministic}
\begin{cases}
(\partial_t^2 +1 - \Delta) u = 0 \qquad (t,x) \in \R \times \T^3, \\
u \big|_{t=0} = \phi^{\cos} , \quad \partial_t u\big|_{t=0} = \langle \nabla \rangle \phi^{\sin}. 
\end{cases}
\end{equation}
As discussed in Section \ref{section:ansatz}, the initial condition is phrased in terms of $(u, \langle \nabla \rangle^{-1} \partial_t u)$, which allows a more symmetric treatment of $\phi^{\cos}$ and $\phi^{\sin}$.  The solution of the linear wave equation \eqref{diagram:eq-wave-linear-deterministic} is then given by
\begin{equation}\label{diagram:eq-wave-linear-propagator}
u = \cos(t\langle \nabla \rangle) \phi^{\cos} + \sin(t\langle \nabla \rangle) \phi^{\sin}.
\end{equation}
Due to our representation of the Gibbs measure through the caloric initial data (Proposition \ref{ansatz:prop-caloric}), we are interested in solving the linear wave equation with initial data given by $\initial{blue}$, $\initial{green}[][\leqN]$, and $\initial{red}[][\leqN]$. Similar as in Subsection \ref{section:diagram-parabolic}, we represent the linear evolutions by the stochastic diagrams 
\begin{equation}
\linear{blue}, \quad \linear{green}[\leqN][], \quad \text{and} \quad \linear{red}[\leqN][].  
\end{equation}
To be more precise, we define 
\begin{align}
\linear{blue}[][(t,n)] &= \cos(t\langle n\rangle) \initial{blue}[cos]\nclose(n) + \sin(t\langle n \rangle) \initial{blue}[sin]\nclose(n), \label{diagram:eq-linear-fourier} \\
\linear{green}[\leqN][(t,n)] &= \cos(t\langle n\rangle) \initial{green}[cos][\leqN] \nclose(n),  \label{diagram:eq-linear-green-fourier} \\
\linear{red}[\leqN][(t,n)] &= \cos(t\langle n\rangle) \initial{red}[cos][\leqN] \nclose(n).  \label{diagram:eq-linear-red-fourier}
\end{align}
We remark that both the green and red linear evolution only contain a cosine-term, since the sine-portion of the initial data vanishes (see Definition \ref{diagram:definition-caloric} and Definition \ref{diagram:def-red}). The three diagrams in \eqref{diagram:eq-linear-fourier}, \eqref{diagram:eq-linear-green-fourier}, and \eqref{diagram:eq-linear-red-fourier} are the longhand diagrams for the linear evolution. In addition to the longhand diagrams, we also work with the shorthand diagrams 
\begin{equation}
\slinear[blue](t,x) := \linear{blue}[][(t,x)], \quad \slinear[green][\leqN](t,x) := \linear{green}[\leqN][(t,x)], \quad \text{and} \quad 
\slinear[red][\leqN](t,x):=\linear{red}[\leqN][(t,x)],
\end{equation}
which were previously introduced in Subsection \ref{section:ansatz-explicit}. 

In the analysis of the nonlinear wave equation, we will encounter probabilistic resonances in products containing multiple factors of $\slinear[blue]$. To address the probabilistic resonances, we will make use of the following covariance identity.

\begin{lemma}[Covariance of linear evolution]\label{diagram:lemma-covariance}
For all $t,t^\prime \in \R$ and $n,n^\prime \in \Z^3$, it holds that
\begin{equation}\label{diagram:eq-covariance}
    \E\bigg[ \linear{blue}[][(t,n)]  \linear{blue}[][(t^\prime,n^\prime)] \bigg]
    = \delta_{n+n^\prime=0} \frac{\cos((t-t^\prime)\langle n \rangle)}{\langle n \rangle^2}. 
\end{equation}
\end{lemma}
\begin{proof}
The argument is a minor modification of the proof of Lemma \ref{diagram:lemma-linear-heat-invariance}. From the definition of the stochastic diagram, we obtain that 
\begin{align*}
    \linear{blue}[][(t,n)]
    &= \sqrt{2} \sum_{\varphi=\cos,\sin} \varphi(t \langle n \rangle) 
    \int_{-\infty}^0 \dW[\varphi][s][n] e^{s\langle n \rangle^2}, \\
       \linear{blue}[][(t^\prime,n^\prime)] 
       &= \sqrt{2} \sum_{\varphi^\prime=\cos,\sin} \varphi^\prime(t^\prime \langle n^\prime \rangle) \int_{-\infty}^0 \dW[\varphi^\prime][s][n^\prime] e^{s\langle n^\prime \rangle^2}.
\end{align*}
It follows that 
\begin{align*}
&\E\bigg[ \linear{blue}[][(t,n)]  \linear{blue}[][(t^\prime,n^\prime)] \bigg]\\
=&2 \sum_{\varphi,\varphi^\prime=\cos,\sin} \bigg( \varphi(t \langle n \rangle)
\varphi^\prime(t^\prime \langle n^\prime \rangle) 
\E \Big[ \big( \int_{-\infty}^0 \dW[\varphi][s][n] e^{s \langle n \rangle^2}\big) \cdot  \big( \int_{-\infty}^0 \dW[\varphi^\prime][s][n^\prime] e^{s\langle n^\prime \rangle^2} \big) \Big] \bigg)\\
=& \frac{\delta_{n+n^\prime=0}}{\langle n \rangle^2} \sum_{\varphi,\varphi^\prime=\cos,\sin} \delta_{\varphi=\varphi^\prime} \varphi(t \langle n \rangle)
\varphi^\prime(t^\prime \langle n \rangle) \\
=& \delta_{n+n^\prime=0} \frac{1}{\langle n \rangle^2}  \Big( 
 \cos(t \langle n \rangle) \cos(t^\prime \langle n \rangle) 
 + \sin(t \langle n \rangle) \sin(t^\prime \langle n \rangle) \Big) \\
 =&  \delta_{n+n^\prime=0} \frac{\cos((t-t^\prime)\langle n \rangle)}{\langle n \rangle^2}. 
\end{align*}
This yields the desired identity \eqref{diagram:eq-covariance}. 
\end{proof}

 \subsubsection{The cubic iterate} 
We recall from \eqref{ansatz:eq-frequency-truncated-NLW} that the frequency-truncated nonlinear wave equation with caloric initial data is given by  
\begin{equation}\label{diagram:eq-wave-truncated}
\begin{cases}
(\partial_t^2 + 1 - \Delta) u_{ \leq N} = - P_{\leq N} \Big( \lcol  (P_{\leq N} u_{\leq N})^3 \rcol + \Gamma_{\leq N} u_{\leq N} \Big) - (\gamma_{\leq N} - \Gamma_{\leq N}) P_{\leq N} u_{\leq N} \\
u_{\leq N}[0]= \initial{blue}[][] - \initial{green}[][\leqM] + \initial{red}[][\leqM].
\end{cases}
\end{equation}
In Subsection \ref{section:ansatz-explicit}, we previously introduced the short-hand diagram $\scubic$. The corresponding long-hand diagram is given by   
\begin{equation}\label{diagram:eq-cubic-physical}
\cubic[\leqN][(t,x)] = \int\displaylimits_{0}^{t} \dt^\prime \,  \frac{\sin((t-t^\prime)\langle \nabla \rangle)}{\langle \nabla \rangle} P_{\leq N} \biglcol \left(  \linear{blue}[\leqN][(t^\prime,x)]\hspace{-1ex} \right)^3 \bigrcol ~. 
\end{equation}
By inserting the representation of $\initial{blue}$ from Definition \ref{diagram:definition-caloric}, we obtain that 
\begin{equation}\label{diagram:eq-cubic-expanded}
\begin{aligned}
 \cubic[\leqN][(t,n)] 
 =&1_{\leq N}(n) \langle n \rangle^{-1} 
\sum_{ \substack{\varphi_1,\varphi_2,\varphi_3 \\ \in \{ \cos,\sin\} }} \sum_{\substack{ n_1,n_2,n_3 \in \Z^3 \colon \\ n=n_{123}}} 
\bigg[ \Big( \prod_{j=1}^3 \frac{1_{\leq N}(n_j)}{\langle n_j \rangle} \Big)  \\
\times &
\Big(\int\displaylimits_{0}^{t} \dt^\prime \sin((t-t^\prime)\langle n \rangle) \prod_{j=1}^3 \varphi_j(t^\prime \langle n_j \rangle) \Big)
\SI[n_j,\varphi_j\colon 1 \leq j \leq 3] \bigg]. 
\end{aligned}
\end{equation}
Here, the multiple stochastic integral $\SI[n_j,\varphi_j\colon 1 \leq j \leq 3]$ is as in \eqref{diagram:eq-multiple-stochastic}. We again emphasize that the longhand diagram is only a different notation for the shorthand diagram, i.e., 
\begin{equation}
\scubic[\leqN] (t,x) := \cubic[\leqN][(t,x)]. 
\end{equation}

\subsection{Quintic diagram}\label{section:quintic-diagram}
We first recall from \eqref{ansatz:eq-quintic} that the shorthand diagram $\squintic[\leqN]$ is defined as the solution of 
\begin{equation}\label{diagram:eq-quintic-shorthand}
\begin{cases}
3 \, (\partial_t^2 +1 - \Delta) \squintic[\leqN] = P_{\leq N}  \Big( 3 \squadratic[\leqN]  \scubic[\leqN]  -  \Gamma_{\leq N} \slinear \Big) \qquad (t,x) \in \R \times \T^3, \\
\squintic \hspace{0ex}[0]=0. 
\end{cases}
\end{equation}
To avoid confusion, we emphasize that a combinatorial factor of three has been included on the left-hand side of \eqref{diagram:eq-quintic-shorthand}, which is consistent with the $3 \squadratic[\leqN] \scubic[\leqN]$-term.   In order to utilize the renormalization in \eqref{diagram:eq-quintic-shorthand}, it is simpler to first compute just the nonlinearity (without the Duhamel integral). Using longhand diagrams, we write
\begin{align}\label{diagram:eq-quintic-pure}
3  P_{\leq N} \Big( \squadratic[\leqN]  \scubic[\leqN] \Big)
&= 6  ~ \quinticnl[\leqN][][0]+ 18 ~ \quinticnl[\leqN][][1] +  18 ~ \quinticnl[\leqN][][2] . 
\end{align}
In \eqref{diagram:eq-quintic-pure}, the pairings of two leafs are marked by combining the two vertices and adjusting the corresponding arrows. In addition, the resulting ``double'' vertex is colored in full. The first, second, and third summand in \eqref{diagram:eq-quintic-pure} contain zero, one, and two pairings, respectively.  The equation also contains two combinatorial factors of six, which are due to the possible number of equivalent pairings. Taking into account the covariance identity (Lemma \ref{diagram:lemma-covariance}), the longhand diagrams are given by the following formulas. For the quintic diagram with zero pairings, we have that 
\begin{equation}\label{diagram:eq-quintic-pure-zero}
\begin{aligned}
\quinticnl[\leqN][(t,n)][0] 
=& 1_{\leq N}(n)  \sum_{ \substack{\varphi_1,\hdots,\varphi_5 \\ \in \{ \cos, \sin\}} }  
\sum_{ \substack{n_1,\hdots,n_5 \in \Z^3\colon \\ n=n_{12345} } }  
\bigg[ \frac{1_{\leq N}(n_{234})}{\langle n_{234} \rangle}  \Big( \prod_{j=1}^5 \frac{1_{\leq N}(n_j)}{\langle n_j \rangle} \Big) 
 \Big( \prod_{j=1,5} \varphi_j(t \langle n_j \rangle) \Big)  \\
&\times \bigg( \int\displaylimits_{0}^{t} \dt^\prime \sin\big((t-t^\prime) \langle n_{234} \rangle \big) \Big( \prod_{j=2,3,4} \varphi_j(t^\prime \langle n_j \rangle) \Big) \bigg)    \\
&\times  \SI[n_j,\varphi_j \colon 1 \leq j \leq 5]  \bigg].
\end{aligned}
\end{equation}
In \eqref{diagram:eq-quintic-pure-zero}, the frequency $n_1$ is associated with the lower left vertex, the three frequencies $n_2,n_3$, and $n_4$ are associated with the three higher vertices (from left to right), and the frequency $n_5$ is associated with the lower right vertex. A similar numbering will also be used in all following diagrams. For the purely hyperbolic quintic diagram with one pairing, we have that 
\begin{equation}\label{diagram:eq-quintic-pure-one}
\begin{aligned}
\quinticnl[\leqN][(t,n)][1] 
=& 1_{\leq N}(n) \sum_{ \substack{\varphi_3,\varphi_4,\varphi_5 \\ \in \{ \cos, \sin\}} }  
\sum_{ \substack{n_3,n_4,n_5 \in \Z^3\colon \\ n=n_{345} } } \sum_{n_2 \in \Z^3}   
\bigg[ \frac{1_{\leq N}(n_{234})}{\langle n_{234} \rangle}  \frac{1_{\leq N}(n_2)}{\langle n_2 \rangle^2} \Big( \prod_{j=3}^5 \frac{1_{\leq N}(n_j)}{\langle n_j \rangle} \Big) 
  \\
&\times \varphi_5(t\langle n_5 \rangle)  \bigg( \int\displaylimits_{0}^{t} \dt^\prime \sin\big((t-t^\prime) \langle n_{234} \rangle \big) \cos((t-t^\prime) \langle n_2 \rangle) \Big( \prod_{j=3,4} \varphi_j(t^\prime \langle n_j \rangle) \Big) \bigg) \\
&\times   \SI[n_j,\varphi_j \colon 3 \leq j \leq 5]  \bigg].
\end{aligned}
\end{equation}
Finally, for the  quintic diagram with two pairings, we have that 
\begin{equation}\label{diagram:eq-quintic-pure-two}
\begin{aligned}
&\quinticnl[\leqN][(t,n)][2] \\
=& 1_{\leq N}(n) \sum_{ \substack{\varphi_3  \in \{ \cos, \sin\}} }  
\sum_{ \substack{n_3\in \Z^3\colon \\ n=n_{3} } } \sum_{n_2,n_4 \in \Z^3}   
\bigg[ \frac{1_{\leq N}(n_{234})}{\langle n_{234} \rangle}  \frac{1_{\leq N}(n_2)}{\langle n_2 \rangle^2}  \frac{1_{\leq N}(n_3)}{\langle n_3 \rangle} 
 \frac{1_{\leq N}(n_4)}{\langle n_4 \rangle^2}
  \\
&\times \bigg( \int\displaylimits_{0}^{t} \dt^\prime \sin\big((t-t^\prime) \langle n_{234} \rangle \big) \cos((t-t^\prime) \langle n_2 \rangle) \varphi_3(t^\prime \langle n_3 \rangle) \cos((t-t^\prime)\langle n_4 \rangle)  \bigg)  \\
&\times  \SI[n_3,\varphi_3]  \bigg]. 
\end{aligned}
\end{equation}

We now take a deeper look at the difference of the resonant part \eqref{diagram:eq-quintic-pure-two} and the renormalization in \eqref{diagram:eq-quintic-shorthand}. To this end, we first make the following definition. 

\begin{definition}[$t$-dependent version of $\Gamma_{\leq N}$]\label{diagram:definition-Gamma-t}
For any $t\in \R$ and $n\in \Z^3$, we define
\begin{equation}
\Gamma_{\leq N}(n,t) := 6 \cdot 1_{\leq N}(n) \sum_{\substack{n_1,n_2,n_3 \in \Z^3 \colon \\ n_{123}=n }} \bigg[  \prod_{j=1}^3 \frac{1_{\leq N}(n_j)}{\langle n_j \rangle^2} \cos\big( t \langle n_j \rangle \big)  \bigg]. 
\end{equation}
We note that $\Gamma_{\leq N}(n,0)$ agrees with $\Gamma_{\leq N}$ from \eqref{diagram:eq-Gamma-heat}. 
\end{definition}

We now prove the following explict formula for the renormalized resonant part of the quintic nonlinearity.

\begin{lemma}[The renormalized resonant part of the quintic object]\label{diagram:lem-quintic-resonant}
For all $N\geq 1$, it holds that 
\begin{equation}\label{diagram:eq-quintic-resonant}
\begin{aligned}
18 ~ \quinticnl[\leqN][(t,n)][2] - \Gamma_{\leq N}(n) \linear{blue}[][(t,n)]  &= - \Gamma_{\leq N}(n,t) \langle n \rangle^{-1} \SI[n,\cos] \\
&- \sum_{\varphi \in \{ \cos, \sin \}} \bigg[ \bigg( \int_0^t \dt^\prime \Gamma_{\leq N}(n,t-t^\prime) \big( \partial_t \varphi \big)\big( t^\prime \langle n \rangle \big) \bigg) \SI[n,\varphi] \bigg]. 
\end{aligned}
\end{equation}
\end{lemma}

\begin{proof}
The desired identity \eqref{diagram:eq-quintic-resonant} is equivalent to 
\begin{equation}\label{diagram:eq-resonant-part-p1}
\begin{aligned}
18 ~~ \quinticnl[\leqN][(t,n)][2] 
    =& \Gamma_{\leq N}(n) \linear{blue}[][(t,n)] - \Gamma_{\leq N}(n,t) \langle n \rangle^{-1} \SI[n,\cos]  \\
-& \sum_{\varphi \in \{\cos,\sin\}} \bigg[ \bigg( \int_{0}^t \dt^\prime\,   \Gamma_{\leq N}(n,t-t^\prime)  (\partial_t \varphi)(t^\prime \langle n \rangle) \bigg) \SI[n,\varphi] \bigg]. 
\end{aligned}
\end{equation}
We first recall the expression for the quintic object with two pairings, which is given by  \eqref{diagram:eq-quintic-pure-two}. 
Due to the constraint on $n_3\in \Z^3$, we simply replace $n_3$ by $n$ and denote $\varphi_3$ by $\varphi$. We also introduce the auxiliary frequency $n_6:= -n_{234} = - n_{24} -n$ and then switch variables $(n_2,n_4,n_6)\rightarrow (-n_2,-n_4,-n_6)$. This yields that 
\begin{align}
&18 ~ \quinticnl[\leqN][(t,n)][2]  \notag \\
=&
18 \cdot 1_{\leq N}(n) \langle n \rangle^{-1} \sum_{\varphi \in \{\cos,\sin\}} \sum_{\substack{n_2,n_4,n_6 \in \Z^3 \colon \\ n_{246}=n}} \bigg[ \bigg( \prod_{j=2,4} \frac{1_{\leq N}(n_j)}{\langle n_j \rangle^2} \bigg) \frac{1_{\leq N}(n_6)}{\langle n_6 \rangle}\notag \\
&\times \Big( \int_{0}^t \dt^\prime\,  \cos((t-t^\prime) \langle n_2 \rangle) 
\cos((t-t^\prime) \langle n_4 \rangle) 
\sin((t-t^\prime) \langle n_6 \rangle) 
\varphi(t^\prime \langle n \rangle) \Big) \SI[n,\varphi] \bigg]\notag \\
=&18 \cdot \frac{1_{\leq N}(n)}{\langle n \rangle} \sum_{\varphi \in \{\cos,\sin\}} \sum_{\substack{n_2,n_4,n_6 \in \Z^3 \colon \\ n_{246}=n}} \bigg[ \Big( \prod_{j=2,4,6} \frac{1_{\leq N}(n_j)}{\langle n_j \rangle^2}  \Big) \notag\\
&\times \bigg( \int_{0}^t \dt^\prime\,  \cos((t-t^\prime) \langle n_2 \rangle) 
\cos((t-t^\prime) \langle n_4 \rangle) 
\partial_{t^\prime} \Big(\cos((t-t^\prime) \langle n_6 \rangle) \Big)
\varphi(t^\prime \langle n \rangle) \bigg) \SI[n,\varphi] \bigg].\label{diagram:eq-resonant-part-p3}
\end{align}

By symmetrizing the sum in the triple $(n_2,n_4,n_6)$, we can perform the replacement
\begin{equation*}
\cos((t-t^\prime) \langle n_2 \rangle) 
\cos((t-t^\prime) \langle n_4 \rangle) 
\partial_{t^\prime} \Big(\cos((t-t^\prime) \langle n_6 \rangle) \Big)
\quad \rightarrow \quad \frac{1}{3} \,  \partial_{t^\prime}\Big( \prod_{j=2,4,6} \cos((t-t^\prime) \langle n_j \rangle) \Big). 
\end{equation*}
By recalling Definition \ref{diagram:definition-Gamma-t} and integrating by parts, we obtain that 
\begin{align*}
\eqref{diagram:eq-resonant-part-p3}
=& \langle n \rangle^{-1} \sum_{\varphi \in \{\cos,\sin\}} \bigg[ \int_{0}^t \dt^\prime\,  \partial_{t^\prime}\Big( \Gamma_{\leq N}(n,t-t^\prime) \Big) \varphi(t^\prime \langle n \rangle) \SI[n,\varphi] \bigg] \\
=&  \langle n \rangle^{-1} \sum_{\varphi\in \{\cos,\sin\}} \Gamma_{\leq N}(n) \varphi(t \langle n \rangle) \SI[n,\varphi]
 -  \langle n \rangle^{-1} \sum_{\varphi\in \{\cos,\sin\}} \Gamma_{\leq N}(n,t) \varphi(0) \SI[n,\varphi] \\
&-  \sum_{\varphi \in \{\cos,\sin\}} \bigg[ \int_{0}^t \dt^\prime\,   \Gamma_{\leq N}(n,t-t^\prime)  (\partial_t \varphi)(t^\prime \langle n \rangle) \SI[n,\varphi] \bigg] \\
=& \Gamma_{\leq N}(n) \linear{blue}[][(t,n)] - \Gamma_{\leq N}(n,t)  \langle n \rangle^{-1} \SI[n,\cos] 
-  \sum_{\varphi \in \{\cos,\sin\}} \bigg[ \int_{0}^t \dt^\prime\,   \Gamma_{\leq N}(n,t-t^\prime)  (\partial_t \varphi)(t^\prime \langle n \rangle) \SI[n,\varphi] \bigg],
\end{align*} which completes the proof.
\end{proof}

In our previous Lemma \ref{diagram:lem-quintic-resonant}, we obtained a detailed description of the resonant part of the quintic nonlinearity. In the following analysis, however, we will usually work with the Duhamel integral of the resonant part \eqref{diagram:eq-quintic-resonant} and not the resonant part itself. For notational purposes, it is convenient to depict this Duhamel integral using a new stochastic diagram. 

\begin{definition}[The resistor]\label{diagram:definition-resistor}
We define the resistor as 
\begin{equation}\label{diagram:eq-resistor}
 18~   \resistor[][\leqN][] := \Duh \left[ 18 ~ \quinticnl[\leqN][][2]  -  \Gamma_{\leq N} \linear{blue}[][] \right]. 
\end{equation}
In addition to the longhand diagram for the resistor in \eqref{diagram:eq-resistor}, we also represent the right-hand side of \eqref{diagram:eq-resistor} as $18\,\sresistor$.
\end{definition}
\begin{remark}
We call the stochastic diagram in \eqref{diagram:eq-resistor} the resistor, and use the same arrow as for resistors in electrical circuits, because this term is already renormalized. As a result, the Gaussian process $\initial{black}$ defined in \eqref{gau:proc} wants to "resist" being propagated, but is only partially successful. 
\end{remark}

The definition \eqref{diagram:eq-resistor} is pleasant from a conceptual perspective. From a practical perspective, however, it is unwieldy. Instead, we will be working with the following representation, which is a direct consequence of Lemma \ref{diagram:lem-quintic-resonant}.

\begin{corollary}[Exact representation of the resistor]\label{diagram:corollary-resistor}
For all $N\geq 1$, it holds that 
\begin{align*}
18~   \resistor[][\leqN][(t,n)] 
&= - \langle n \rangle^{-2}  \Big( \int_0^t \dt^\prime \sin \big( (t-t^\prime) \langle n \rangle \big) \Gamma_{\leq N}(n,t^\prime) \Big) \SI[n,\cos] \\
&- \langle n \rangle^{-1} \sum_{\varphi \in \{ \cos, \sin \}} \bigg[  \bigg( \int_0^t \dt^\prime \int_0^{t^\prime} \dt^{\prime \prime} \sin\big( (t-t^\prime) \langle n \rangle \big) \Gamma_{\leq N}(n,t^\prime -t^{\prime \prime}) \big(\partial_t \varphi\big)\big( t^{\prime \prime} \langle n \rangle \big) \bigg) \SI[n,\varphi] \bigg].
\end{align*}
\end{corollary}

\subsection{Heat-wave quintic diagram}\label{section:diagram-mixed} 

In this subsection, we examine heat-wave stochastic diagrams. In the previous subsections, our stochastic objects contained either heat propagators (Subsection \ref{section:diagram-parabolic}) or wave propagators (Subsection \ref{section:diagram-wave} and Subsection \ref{section:quintic-diagram}), but not both\footnote{Strictly speaking, the stochastic integrals $\SI$ in Subsection \ref{section:diagram-wave} and \ref{section:quintic-diagram} contain heat propagators, but this does not enter into our analysis.}. In contrast, the heat-wave stochastic objects simultaneously contain both heat and wave propagators. The heat-wave stochastic objects occur naturally due to our caloric representation of the Gibbs measure which, in turn, is needed due to the singularity of the Gibbs measure. For the analytical aspects of the heat-wave stochastic objects, we refer the reader to Subsection \ref{section:analytic-quintic-mixed} below.

 We now let $M \geq N$ and examine the heat-wave quintic term  $\squadratic[\leqN] P_{\leq N} \slinear[green][\leqM]$. 
Using the definition of the caloric initial data $\initial{blue}$ and $\initial{green}$ (from Definition \ref{diagram:definition-caloric}), we obtain that
\begin{equation}\label{diagram:eq-quintic-mixed}
P_{\leq N} \Big( \squadratic[\leqN] P_{\leq N} \slinear[green][\leqM] \Big)
= \scalebox{0.9}{\quinticmixednl[\leqN,\leqM][][0] \hspace{-2ex}+ 6 \hspace{-4ex} \quinticmixednl[\leqN,\leqM][][1] + 6  \hspace{-5ex} \quinticmixednl[\leqN][][2].} 
\end{equation}
Just as in \eqref{diagram:eq-quintic-pure}, the pairings of two leafs are marked by combining two vertices and adjusting the corresponding arrows. In the last diagram of \eqref{diagram:eq-quintic-mixed}, we only indicate the dependence on $N$, since the corresponding expression will be the same for all $M\geq N$ (see \eqref{diagram:eq-quintic-mixed-two} below). The heat-wave longhand diagrams are given by the following expressions. For the heat-wave quintic object with zero pairings, we have that 
\begin{equation}\label{diagram:eq-quintic-mixed-nopairing}
\begin{aligned}
& \scalebox{0.8}{\quinticmixednl[\leqN,\leqM][(t,n)][0]}  \allowdisplaybreaks[3] \\
=& 2^{\frac{5}{2}} 1_{\leq N}(n) 
\sum_{\substack{n_1,\hdots,n_5\in \Z^3\colon \\ n=n_{12345}  }} 
\sum_{\substack{ \varphi_1,\hdots,\varphi_5 \\ \in \{ \cos, \sin \} } } 
\bigg[ \Big( \prod_{j=1,5 }1_{\leq N}(n_j) \Big) 
1_{\leq N}(n_{234}) \Big( \prod_{j=2,3,4} 1_{\leq M}(n_j) \Big) 
 \\
&\times  \Big( \prod_{j=2}^4 \mathbf{1}\{\varphi_j=\cos\}  \Big) \varphi_1(t\langle n_1 \rangle) \cos(t \langle n_{234} \rangle) \varphi_5(t\langle n_5 \rangle)  \\
&\times \int\displaylimits_{[-\infty,0]^5} \bigotimes_{j=1}^5 \dW[\varphi_j][s_j](n_j)  \hspace{-2ex} \int\displaylimits_{\max(s_2,s_3,s_4)}^0  \hspace{-4ex} \dtau 
\Big(\prod_{j=1,5} e^{s_j \langle n_j \rangle^2} \Big)  e^{\tau \langle n_{234} \rangle^2} \Big( \prod_{j=2,3,4} e^{-(\tau-s_j) \langle n_j \rangle^2} \Big) \bigg] . 
\end{aligned}
\end{equation}
For the heat-wave quintic object with one pairing, we have that 
\begin{equation}\label{diagram:eq-quintic-mixed-onepairing}
\begin{aligned}
& \scalebox{0.8}{\quinticmixednl[\leqN,\leqM][(t,n)][1]}  \allowdisplaybreaks[3] \\
=&  2^{\frac{3}{2}} 1_{\leq N}(n) 
\sum_{\substack{n_3,n_4,n_5\in \Z^3\colon \\ n=n_{345}  }} \sum_{n_2 \in \Z^3}
\sum_{\substack{ \varphi_3,\varphi_4,\varphi_5 \\ \in \{ \cos, \sin \} } } 
\bigg[ 1_{\leq N}(n_2) 1_{\leq N}(n_{234}) \Big( \prod_{j=3}^4 1_{\leq M}(n_j) \Big) 1_{\leq N}(n_5)   \langle n_2 \rangle^{-2} \\
&\times  \Big( \prod_{j=3}^4 \mathbf{1}\{\varphi_j=\cos\}  \Big) \cos(t\langle n_2 \rangle) \cos(t \langle n_{234} \rangle) \varphi_5(t\langle n_5 \rangle)  \\
&\times \int\displaylimits_{[-\infty,0]^3} \bigotimes_{j=3}^5 \dW[\varphi_j][s_j](n_j)  \hspace{-2ex} 
 \int\displaylimits_{\max(s_3,s_4)}^0  \hspace{-2ex} \dtau ~ 
  e^{\tau \big( \langle n_2 \rangle^2 + \langle n_{234} \rangle^2\big)} \Big( \prod_{j=3,4} e^{-(\tau-s_j) \langle n_j \rangle^2} \Big) e^{s_5 \langle n_5 \rangle^2}  \bigg] . 
\end{aligned}
\end{equation}
Finally, the heat-wave quintic object with two pairings is given by 
\begin{equation}\label{diagram:eq-quintic-mixed-two}
\begin{aligned}
 \scalebox{0.8}{\quinticmixednl[\leqN][(t,n)][2] }
=&  2^{\frac{1}{2}} 1_{\leq N}(n) 
\sum_{\substack{n_3\in \Z^3\colon \\ n=n_{3}  }} \sum_{n_2,n_4 \in \Z^3}
\bigg[ 1_{\leq N}(n_{234}) \Big( \prod_{j=2,3,4} 1_{\leq N}(n_j) \Big)  \langle n_2 \rangle^{-2}  \langle n_4 \rangle^{-2} \\
&\times  \cos(t\langle n_2 \rangle) \cos(t \langle n_{234} \rangle) \cos(t\langle n_4 \rangle)  \\
&\times \int_{-\infty}^{0} \dW[\cos][s_3](n_3)  
 \int\displaylimits_{s_3}^0  \dtau ~ 
  e^{\tau \big( \langle n_2 \rangle^2 + \langle n_{234} \rangle^2 + \langle n_4 \rangle^2 \big)}  e^{-(\tau-s_3) \langle n_3 \rangle^2}   \bigg] . 
\end{aligned}
\end{equation}
In \eqref{diagram:eq-quintic-mixed-two}, we were able to replace $1_{\leq M}(n_3)$ by $1_{\leq N}(n_3)$ due to the restriction $n_3=n$. In particular, the resonant part \eqref{diagram:eq-quintic-mixed-two} no longer depends on $M$ (as long as $M\geq N$). 

Similar as in Lemma \ref{diagram:lem-quintic-resonant} above, we rewrite the resonant part in \eqref{diagram:eq-quintic-mixed-two} using the $\Gamma_{\leq N}$-multiplier. For this, however, we need the following final generalization of $\Gamma_{\leq N}$. 

\begin{definition}[$(t,s)$-dependent version of $\Gamma_{\leq N}(n)$]\label{diagram:definition-Gamma}
For any $t\in \R$, $s\leq 0$, and $n\in \Z^3$, we define
\begin{equation}
\Gamma_{\leq N}(n,t,s):= 6\cdot 1_{\leq N}(n) \sum_{\substack{n_1,n_2,n_3\in \Z^3\colon \\ n_{123}=n}} \bigg[ \prod_{j=1}^3 
\Big( \frac{1_{\leq N}(n_j)}{\langle n_j \rangle^2} 
\cos(t \langle n_j \rangle) e^{s \langle n_j \rangle^2} \Big) \bigg].
\end{equation}
We note that $\Gamma_{\leq N}(n,t,0)$ agrees with $\Gamma_{\leq N}(n,t)$ from Definition \ref{diagram:definition-Gamma-t}. 
\end{definition}

\begin{lemma}[Resonant part of the heat-wave quintic diagram]\label{diagram:lem-resonant-mixed-quintic}
For all $N \geq 1$, it holds that 
\begin{align*}
18 ~ \quinticmixednl[\leqN][(t,n)][2] &= \Gamma_{\leq N}(n,t)  \langle n \rangle^{-1} \SI[n,\cos]
- \sqrt{2} 
\int_{-\infty}^{0} \dW[\cos][s][n]\Gamma_{\leq N}(n,t,s)  \\
&+\sqrt{2}  \langle n \rangle^2 \int_{-\infty}^{0} \dW[\cos][s][n] \Big( \int_s^0 \dtau \,  \Gamma_{\leq N}(n,t,\tau)  e^{-(\tau-s) \langle n\rangle^2} \Big).
\end{align*}
\end{lemma}
\begin{remark}
Surprisingly, the first summand in Lemma \ref{diagram:lem-resonant-mixed-quintic} exactly cancels the first summand in Lemma \ref{diagram:lem-quintic-resonant}. However, this cancellation turns out to not be essential for our argument and will not be used. 
\end{remark}
\begin{proof} We recall the representation of the heat-wave quintic object from \eqref{diagram:eq-quintic-mixed-two}. We then note that $n_3 \in \Z^3$ is fixed as $n_3 = n$. As before, we introduce the auxiliary frequency-variable $n_6:= -n_{234}=-n_{24}-n$ and then switch variables $(n_2,n_4,n_6)\rightarrow (-n_2,-n_4,-n_6)$. For notational convenience, we also relabel the $s_3$-variable as $s$. This yields  
\begin{align*}
&18 \quinticmixednl[\leqN][(t,n)][2] \allowdisplaybreaks[3]\\
=& \sqrt{2}\cdot 18 \cdot 1_{\leq N}(n) 
\sum_{\substack{n_2,n_4,n_6\in \Z^3\colon \\ n_{246}=n}} \bigg[
\Big( \prod_{j=2,4,6} 1_{\leq N}(n_j) \cos(t\langle n_j \rangle) \Big) 
\langle n_2 \rangle^{-2} \langle n_4 \rangle^{-2} \\
&\times \int_{-\infty}^{0} \dW[\cos][s][n] 
\Big( \int_s^0 \dtau \, e^{\tau (\langle n_2 \rangle^2+\langle n_4 \rangle^2)}
e^{\tau \langle n_6 \rangle^2} e^{-(\tau-s) \langle n \rangle^2} \Big)
\bigg] \\
=& \sqrt{2}\cdot 18 \cdot 1_{\leq N}(n) 
\sum_{\substack{n_2,n_4,n_6\in \Z^3\colon \\ n_{246}=n}} \bigg[
\Big( \prod_{j=2,4,6} \frac{1_{\leq N}(n_j)}{\langle n_j \rangle^2} \cos(t\langle n_j \rangle) \Big) 
 \\
&\times \int_{-\infty}^{0} \dW[\cos][s][n] 
\Big( \int_s^0 \dtau \, e^{\tau (\langle n_2 \rangle^2+\langle n_4 \rangle^2)}
\partial_\tau \Big( e^{\tau \langle n_6 \rangle^2} \Big) e^{-(\tau-s) \langle n \rangle^2} \Big)
\bigg].
\end{align*}
By symmetrizing in the triple $(n_2,n_4,n_6)$, we can perform the replacement 
\begin{equation*}
e^{\tau (\langle n_2 \rangle^2+\langle n_4 \rangle^2)}
\partial_\tau \Big( e^{\tau \langle n_6 \rangle^2} \Big)
\quad \rightarrow \quad 
\frac{1}{3} \partial_\tau \Big( \prod_{j=2,4,6} e^{\tau\langle n_j \rangle^2} \Big). 
\end{equation*}
Using Definition \ref{diagram:definition-Gamma} and integration by parts, this leads to 
\begin{align*}
& \sqrt{2}\cdot 18 \cdot 1_{\leq N}(n) 
\sum_{\substack{n_2,n_4,n_6\in \Z^3\colon \\ n_{246}=n}} \bigg[
\Big( \prod_{j=2,4,6} \frac{1_{\leq N}(n_j)}{\langle n_j \rangle^2} \cos(t\langle n_j \rangle) \Big) 
 \\
&\times \int_{-\infty}^{0} \dW[\cos][s][n] 
\Big( \int_s^0 \dtau \, e^{\tau (\langle n_2 \rangle^2+\langle n_4 \rangle^2)}
\partial_\tau \Big( e^{\tau \langle n_6 \rangle^2} \Big) e^{-(\tau-s) \langle n \rangle^2} \Big)
\bigg] \\
=& \sqrt{2} \cdot \int_{-\infty}^{0} \dW[\cos][s][n] \Big( \int_s^0 \dtau \, \partial_\tau \Big( \Gamma_{\leq N}(n,t,\tau) \Big) e^{-(\tau-s) \langle n\rangle^2} \Big) \\
=& \sqrt{2} \cdot  \Gamma_{\leq N}(n,t,0) 
\int_{-\infty}^{0} \dW[\cos][s][n] e^{s \langle n \rangle^2} \\
-& \sqrt{2} \cdot  
\int_{-\infty}^{0} \dW[\cos][s][n]\Gamma_{\leq N}(n,t,s)  \\
+&\sqrt{2} \cdot  \langle n \rangle^2 \int_{-\infty}^{0} \dW[\cos][s][n] \Big( \int_s^0 \dtau \,  \Gamma_{\leq N}(n,t,\tau)  e^{-(\tau-s) \langle n\rangle^2} \Big).
\end{align*}
After recalling the definition of $\SI[n,\cos]$ from \eqref{diagram:eq-multiple-stochastic}, this completes the proof.
\end{proof}
\subsection{\protect{Sextic diagrams and the \oftt-cancellation}}\label{section:sextic-diagrams}

In this subsection, we examine the resonant parts of sextic diagrams. Our main objective is to exhibit the \oftt-cancellation, whose significance was already discussed in Subsection \ref{section:ansatz-remainder}.  Recalling the notation and definitions in  \eqref{ansatz:eq-mathfrakC} in  Subsection \ref{section:linearquintic}, we first derive an explicit formula for $\mathfrak{C}^{(1,5)}_{\leq N}$. In Subsection \ref{section:cubiccubic}, we then derive an explicit formula for $\mathfrak{C}^{(3,3)}_{\leq N}$. In Subsection \ref{section:1533}, we combine both explicit formulas and exhibit the \oftt-cancellation. Finally, in Subsection \ref{section:mathfrakC}, we control $\mathfrak{C}_{\leq N}$.

\subsubsection{The linear$\times$quintic-object}\label{section:linearquintic}
We now derive an explicit formula for $\mathfrak{C}^{(1,5)}_{\leq N}$. As a result, we are interested in the expectation (or resonant part) of $\slinear[blue][\leqN] \cdot \squintic[\leqN]$.
This expectation is not affected by the quintic and cubic chaos in \eqref{diagram:eq-quintic-pure}, which, even after multiplication with a Gaussian random variable, always have zero expectation. 
As a result, we are only interested in the product of the linear evolution $\slinear$ and the resistor $\sresistor$. Keeping the same pre-factor as in Definition \ref{diagram:definition-resistor}, we write this term as 
\begin{equation}\label{diagram:eq-product-linear-resistor}
18 \, \linear{blue}[\leqN] \resistor[][\leqN][] = 18 \,  \linearandresistor[\leqN][][0] + 18\, \linearandresistor[\leqN][][1].
\end{equation}

The first summand in \eqref{diagram:eq-product-linear-resistor}, which is defined as the non-resonant part of the product, can be expressed using multiple stochastic integrals. However, since its expectation is zero, this explicit representation is not needed here. Indeed, we obtain that 
\begin{equation}\label{diagram:eq-C15-linear-resistor}
3 \mathfrak{C}^{(1,5)}_{\leq N} = 3 \mathbb{E} \Big[ \slinear[blue][\leqN] \squintic[\leqN] \Big] 
= 18 \linearandresistor[\leqN][][1].
\end{equation}
In particular, we are only interested in the resonant part of the linear$\times$resistor-object from \eqref{diagram:eq-product-linear-resistor}. 

\begin{lemma}[Explicit formula for $\mathfrak{C}^{(1,5)}_{\leq N}$]\label{diagram:lemma-linear-quintic}
For all $N\geq 1$, it holds that 
\begin{equation}\label{diagram:eq-linear-quintic}
\begin{aligned}
    3 \, \mathfrak{C}^{(1,5)}_{\leq N} &= 
    - \frac{1}{2} \sum_{n\in \Z^3} \langle n \rangle^{-2} \int_0^t \dt^{\prime} \int_0^{t} \dt^{\prime\prime} 
    \sin\big( (t-t^\prime) \langle n \rangle\big)  \sin\big( (t-t^{\prime\prime}) \langle n \rangle\big) 
    \Gamma_{\leq N}(n,t^\prime-t^{\prime\prime}) \\
    &- \sum_{n \in \Z^3} \bigg[ \langle n \rangle^{-3} \cos(t\langle n \rangle) \int_0^t \dt^\prime \sin\big( (t-t^\prime) \langle n \rangle \big) \Gamma_{\leq N}(n,t^\prime) \bigg].  
\end{aligned}
\end{equation}
\end{lemma}

\begin{proof}
Due to \eqref{diagram:eq-C15-linear-resistor}, it remains to analyze the resonant part of the linear$\times$resistor-object. Using Corollary \ref{diagram:corollary-resistor}, it can be explicitly written as 
\begin{align}
    &18\, \linearandresistor[\leqN][][1]  \notag \\
    =&  - \sum_{n \in \Z^3} \langle n \rangle^{-3} \cos(t\langle n \rangle) \int_0^t \dt^\prime \sin\big( (t-t^\prime) \langle n \rangle \big) \Gamma_{\leq N}(n,t^\prime) \label{diagram:eq-linear-quintic-p0} \\   
    -&  \sum_{\varphi\in \{ \cos, \sin \}} \sum_{n\in \Z^3} \langle n \rangle^{-2} \int_0^t \dt^{\prime} \int_0^{t^\prime} \dt^{\prime\prime} 
    \sin\big( (t-t^\prime) \langle n \rangle\big) \varphi\big(t\langle n\rangle\big) (\partial_t \varphi)\big( t^{\prime\prime} \langle n \rangle \big)
    \Gamma_{\leq N}(n,t^\prime-t^{\prime\prime}) \label{diagram:eq-linear-quintic-p1}
\end{align}
The first term \eqref{diagram:eq-linear-quintic-p0} coincides with the second term in \eqref{diagram:eq-linear-quintic}. Regarding the second term \eqref{diagram:eq-linear-quintic-p1}, we first note that 
\begin{align*}
\sum_{\varphi\in \{ \cos, \sin\}} \varphi\big( t \langle n \rangle\big) 
(\partial_t \varphi) \big( t^{\prime\prime} \langle n \rangle\big) 
&= - \cos\big( t \langle n \rangle\big) \sin\big( t^{\prime\prime} \langle n \rangle\big)
+ \sin\big( t \langle n \rangle\big) \cos\big( t^{\prime\prime} \langle n \rangle\big) \\
&= \sin\big( (t-t^{\prime\prime}) \langle n \rangle\big). 
\end{align*}
By symmetrizing in $t^\prime$ and $t^{\prime\prime}$, we can rewrite \eqref{diagram:eq-linear-quintic-p1} as
\begin{align*}
    &-  \sum_{\varphi\in \{ \cos, \sin \}} \sum_{n\in \Z^3} \langle n\rangle^{-2} \int_0^t \dt^{\prime} \int_0^{t^\prime} \dt^{\prime\prime} 
    \sin\big( (t-t^\prime) \langle n \rangle\big) \varphi\big(t\langle n\rangle\big) (\partial_t \varphi)\big( t^{\prime\prime} \langle n \rangle \big)
    \Gamma_{\leq N}(n,t^\prime-t^{\prime\prime}) \\
    &=- \frac{1}{2} \sum_{n\in \Z^3} \langle n\rangle^{-2}  \int_0^t \dt^{\prime} \int_0^{t} \dt^{\prime\prime} 
    \sin\big( (t-t^\prime) \langle n \rangle\big)  \sin\big( (t-t^{\prime\prime}) \langle n \rangle\big) 
    \Gamma_{\leq N}(n,t^\prime-t^{\prime\prime}). 
\end{align*}
This yields the first term in \eqref{diagram:eq-linear-quintic} and therefore has the desired form. 
\end{proof}

\subsubsection{\protect{The cubic$\times$cubic-object}}\label{section:cubiccubic} 
We now turn to the explicit formula for $\mathfrak{C}^{(3,3)}_{\leq N}$, which was previously defined as
\begin{equation*}
\mathfrak{C}^{(3,3)}_{\leq N}= \mathbb{E} \bigg[ \Big( \scubic[\leqN] \Big)^2 \bigg].
\end{equation*}
In order to analyze $\mathfrak{C}^{(3,3)}_{\leq N}$, we utilize the longhand diagrams. In this notation, we can decompose the square of  the cubic object $\scubic$ as 
\begin{equation}\label{diagram:eq-square-cubic}
\begin{aligned}
&\quad \left( \cubic[\leqN][(t,x)] \right)^2   \\
&=\cubiccubic[\leqN][(t,x)][0] 
+ 9 \cubiccubic[\leqN][(t,x)][1] 
+ 18  \, \cubiccubic[\leqN][(t,x)][2] 
+ 6 \, \cubiccubic[\leqN][(t,x)][3]. 
\end{aligned}
\end{equation}
The first three summands in \eqref{diagram:eq-square-cubic}, which contain zero, one, or two pairings, can be represented using multiple stochastic integrals (as for the linear, cubic, and quintic diagrams). However, since we are currently interested in only the expectation of \eqref{diagram:eq-square-cubic}, we do not need to analyze them here. After taking expectations, we obtain that 
\begin{equation*}
  \mathfrak{C}^{(3,3)}_{\leq N} = \mathbb{E} \left[  \left( \cubic[\leqN][(t,x)] \right)^2 \right] = 
   6 \, \cubiccubic[\leqN][(t,x)][3],
\end{equation*}
which implies that only the last stochastic object in \eqref{diagram:eq-square-cubic} is relevant here. 
By introducing the auxiliary variable $n:= n_{123}$, where $n_1,n_2$, and $n_3$ are the three frequencies of the vertices, we obtain that 
\begin{equation}\label{diagram:eq-cubiccubic-1}
    \begin{aligned}
     6 \, \cubiccubic[\leqN][(t,x)][3] 
    =& 6 \sum_{\substack{n,n_1,n_2,n_3\in \Z^3 \colon \\ n_{123}=n}} 
    \bigg[  \frac{1_{\leq N}(n)}{\langle n \rangle^2}
    \prod_{1\leq j \leq 3} \frac{1_{\leq N}(n_j)}{\langle n_j \rangle^2} 
      \bigg( \int_0^t \dt^\prime \int_0^t \dt^{\prime\prime} \sin\Big((t-t^\prime) \langle n \rangle\Big)\\
    &\times\sin\Big((t-t^{\prime\prime}) \langle n \rangle \Big) \prod_{j=1,2,3} \cos\Big( (t^\prime - t^{\prime\prime}) \langle n_j \rangle \Big) \bigg) \bigg].  
    \end{aligned}
\end{equation}

Using the definition of $\Gamma_{\leq N}$, \eqref{diagram:eq-cubiccubic-1} leads to the following lemma.

\begin{lemma}[Explicit formula for $\mathfrak{C}^{(3,3)}_{\leq N}$]\label{diagram:lemma-cubic-cubic}
For all $N\geq 1$, it holds that 
\begin{equation*}
   \mathfrak{C}^{(3,3)}_{\leq N} =
   \sum_{n\in \Z^3} \bigg[ \frac{1}{\langle n \rangle^2}
   \int_0^t \dt^\prime \int_0^t \dt^{\prime\prime} \sin\Big((t-t^\prime) \langle n \rangle\Big)\sin\Big((t-t^{\prime\prime}) \langle n \rangle\Big) 
   \Gamma_{\leq N}(n,t^\prime-t^{\prime\prime}) \bigg]. 
\end{equation*}
\end{lemma}

\subsubsection{\protect{The \oftt-cancellation}}\label{section:1533} 
Equipped with the explicit formulas for $\mathfrak{C}^{(1,5)}_{\leq N}$ and $\mathfrak{C}^{(3,3)}_{\leq N}$, we now exhibit the \oftt-cancellation. 

\begin{proposition}[The \oftt-cancellation]\label{diagram:prop-1533-cancellation}
For all $N\geq 1$, it holds that 
\begin{equation}\label{diagram:eq-main-cancellation-identity}
\begin{aligned}
    \mathfrak{C}_{\leq N}(t) 
    &= - 2 \sum_{n \in \Z^3} \bigg[ \langle n \rangle^{-3} \cos(t\langle n \rangle) \int_0^t \dt^\prime \sin\big( (t-t^\prime) \langle n \rangle \big) \Gamma_{\leq N}(n,t^\prime) \bigg]. 
\end{aligned}
\end{equation}
\end{proposition}

\begin{remark}
We emphasize that the \oftt-cancellation leading to \eqref{diagram:eq-main-cancellation-identity} is quite different from a renormalization, which requires a modification of the equation. We refer to Proposition \ref{diagram:prop-1533-cancellation} as a cancellation because $\mathfrak{C}_{\leq N}$ is controlled uniformly in $N$ (see Lemma \ref{diagrams:lem-mathfrakC}), but the individual terms $\mathfrak{C}^{(1,5)}_{\leq N}$ and $\mathfrak{C}^{(3,3)}_{\leq N}$ from \eqref{ansatz:eq-mathfrakC} are unbounded as $N \to \infty$ (Lemma \ref{diagram:lem-C33}). 
\end{remark}

\begin{proof} 
By definition, it holds that 
\begin{equation*}
\mathfrak{C}_{\leq N}= 6 \mathfrak{C}^{(1,5)}_{\leq N} + \mathfrak{C}^{(3,3)}_{\leq N}.
\end{equation*}
By inserting the explicit formulas for $\mathfrak{C}^{(1,5)}_{\leq N}$ and $\mathfrak{C}^{(3,3)}_{\leq N}$ (Lemma \ref{diagram:lemma-linear-quintic} and Lemma \ref{diagram:lemma-cubic-cubic}), it follows that 
\begin{equation}\label{diagram:eq-main-cancellation-p1}
\begin{aligned}
 \mathfrak{C}_{\leq N}(t) =& 
    -  \sum_{n\in \Z^3} \langle n\rangle^{-2} \int_0^t \dt^{\prime} \int_0^{t} \dt^{\prime\prime} 
    \sin\big( (t-t^\prime) \langle n \rangle\big)  \sin\big( (t-t^{\prime\prime}) \langle n \rangle\big) 
    \Gamma_{\leq N}(n,t^\prime-t^{\prime\prime}) \\
    &-2 \sum_{n \in \Z^3}  \langle n \rangle^{-3} \cos(t\langle n \rangle) \int_0^t \dt^\prime \sin\big( (t-t^\prime) \langle n \rangle \big) \Gamma_{\leq N}(n,t^\prime) \\
    &+\sum_{n\in \Z^3}  \langle n\rangle^{-2}  \int_0^t \dt^{\prime} \int_0^{t} \dt^{\prime\prime} 
    \sin\big( (t-t^\prime) \langle n \rangle\big)  \sin\big( (t-t^{\prime\prime}) \langle n \rangle\big) 
    \Gamma_{\leq N}(n,t^\prime-t^{\prime\prime}) .
\end{aligned}
\end{equation}
We now observe that the first and third summand in \eqref{diagram:eq-main-cancellation-p1} exactly cancel each other. As a result, we obtain the desired identity \eqref{diagram:eq-main-cancellation-identity}. 
\end{proof}

\subsubsection{\protect{The cancellation-constant $\mathfrak{C}_{\leq N}$}}\label{section:mathfrakC}

At the end of this subsection, we examine the behavior of  $\mathfrak{C}_{\leq N}$,  $\mathfrak{C}^{(1,5)}_{\leq N}$, and $\mathfrak{C}^{(3,3)}_{\leq N}$. While the next lemma is rather analytic, and therefore conflicts with the algebraic spirit of this section, we like to keep it close to our discussion of the \oftt-cancellation. 

\begin{lemma}[Control of $\mathfrak{C}_{\leq N}$]\label{diagrams:lem-mathfrakC}
For all $\chi \in C^\infty_c(\R)$, it holds that 
\begin{equation}\label{diagrams:eq-mathfrakC}
\| \chi(t) \mathfrak{C}_{\leq N}(t) \|_{H_t^{1-\epsilon}} \lesssim_{\chi} 1. 
\end{equation}
\end{lemma}

\begin{proof}
By inserting the definition of $\Gamma_{\leq N}(n,t)$ and relabeling $n$ as $n_0$, we obtain that 
\begin{align*}
\mathfrak{C}_{\leq N}(t) 
&= -2 \sum_{\substack{N_0,N_1,N_2,N_3 \colon \\ N_{\textup{max}} \leq N}} 
\sum_{\substack{n_0,n_1,n_2,n_3 \in \Z^3 \colon \\ n_0 = n_{123} }} \bigg[
\Big( \prod_{j=0}^3 1_{N_j}(n_j) \Big) \langle n_0 \rangle^{-3} \langle n_1 \rangle^{-2} \langle n_2 \rangle^{-2} \langle n_3 \rangle^{-2}  \\
&\hspace{2ex} \times  \cos(t\langle n \rangle)  \int_0^t \dt^\prime \sin \big( (t-t^\prime) \langle n_0 \rangle\big) 
\prod_{j=1}^3 \cos\big( t^\prime \langle n_j \rangle \big) \bigg] \\
&=: \sum_{\substack{N_0,N_1,N_2,N_3 \colon \\ N_{\textup{max}} \leq N}} 
\mathfrak{C}[N_0,N_1,N_2,N_3](t). 
\end{align*}
To obtain the desired estimate, it suffices to prove the stronger estimate
\begin{equation}\label{diagrams:eq-mathfrakC-p1}
\big| \mathfrak{C}[N_0,N_1,N_2,N_3](t)\big| \lesssim \langle t \rangle \,  N_{\textup{max}}^{-1+\epsilon/2}. 
\end{equation}
By performing the $t^\prime$-integral and then using the level-set decomposition
\begin{equation*}
\bigg( 1 + \Big| \sum_{j=0}^3 (\pm_j) \langle n_j \rangle\Big|\bigg)^{-1} 
\leq \sum_{m\in \Z} \frac{1}{\langle m \rangle} \,  \mathbf{1}\bigg\{ \Big|
\sum_{j=0}^3 (\pm_j) \langle n_j \rangle -m \Big| \leq 1 \bigg\},
\end{equation*}
we obtain that 
\begin{align}
&\big| \mathfrak{C}[N_0,N_1,N_2,N_3](t)\big| \notag \\
\lesssim&\, N_0^{-3} N_1^{-2} N_2^{-2} N_3^{-3} 
\sum_{\substack{\pm_0,\pm_1, \\ \pm_2,\pm_3}}
\sum_{\substack{n_0,n_1,n_2,n_3 \in \Z^3 \colon \\ n_0 = n_{123} }}
\bigg[ \Big( \prod_{j=0}^3 1_{N_j}(n_j) \Big) \bigg( 1 + \Big| \sum_{j=0}^3 (\pm_j) \langle n_j \rangle\Big|\bigg)^{-1} \bigg] \notag  \\
\lesssim&\, N_{\textup{max}}^{\epsilon/2} N_0^{-3} N_1^{-2} N_2^{-2} N_3^{-3}  \sup_{m\in \Z}
\sum_{\substack{\pm_0,\pm_1, \\ \pm_2,\pm_3}}
\sum_{\substack{n_0,n_1,n_2,n_3 \in \Z^3 \colon \\ n_0 = n_{123} }}
\bigg[ \Big( \prod_{j=0}^3 1_{N_j}(n_j) \Big) \mathbf{1}\Big\{ \Big| \sum_{j=0}^3 (\pm_j) \langle n_j \rangle - m\Big| \leq 1 \Big\} \bigg]. \label{diagrams:eq-mathfrakC-p2}
\end{align}
Using our counting estimate \eqref{counting:base1'}, we obtain that
\begin{equation*}
\eqref{diagrams:eq-mathfrakC-p2} \lesssim N_{\textup{max}}^{\epsilon/2} N_0^{-1} N_{\textup{min}} N_{\textup{max}}^{-1} \lesssim N_{\textup{max}}^{-1+\epsilon/2}.  
\end{equation*}
This completes the proof of \eqref{diagrams:eq-mathfrakC-p1} and hence the proof of the Lemma. 
\end{proof}

In the final lemma of this section, we prove the divergence of $\mathfrak{C}^{(3,3)}_{\leq N}$. While this lemma will not be used in this paper, it shows that the \oftt-cancellation is essential to our argument. 

\begin{lemma}[Divergence of $\mathfrak{C}^{(3,3)}_{\leq N}$]\label{diagram:lem-C33}
Let $\rho \colon \R \rightarrow [0,1]$ be any smooth function satisfying $\rho|_{[1/4,3/4]}=1$ and $\rho|_{\R\backslash [1/8,7/8]}=0$. Then, we obtain for all $N\geq 1$ that 
\begin{equation}\label{diagram:eq-C33}
\int_{\R} \mathrm{d}t \rho(t) \mathfrak{C}_{\leq N}^{(3,3)}(t) \geq c \log(N) - c^{-1}, 
\end{equation}
where $0<c\ll 1$ is a sufficiently small absolute constant. In particular, $\mathfrak{C}_{\leq N}^{(3,3)}$ diverges in the sense of distributions. 
\end{lemma}

\begin{remark}[Connection of $\mathfrak{C}^{(3,3)}_{\leq N}$ to the $L_t^4 L_x^4$-Strichartz estimate]\label{diagram:rem-strichartz}
The $L_t^4L_x^4$-Strichartz estimate for the periodic wave equation is given by 
\begin{equation}\label{diagram:eq-strichartz}
\Big\| \sum_{n \in \Z^3} \cos\big( t \langle n \rangle \big) \widehat{\phi}^{\cos}(n) e^{i\langle n ,x \rangle} \Big\|_{L_t^4 L_x^4([-1,1]\times \T^3)} \lesssim \| \phi^{\cos} \|_{H_x^{1/2}(\T^3)}. 
\end{equation}
As will be shown below, the $H_x^{1/2}$-norm in \eqref{diagram:eq-strichartz} cannot be replaced by the weaker $B^{1/2}_{2,\infty}$-Besov norm. This will turn out to be reason for the divergence of $\mathfrak{C}_{\leq N}^{(3,3)}$. 
\end{remark}

\begin{proof}[Proof of Lemma \ref{diagram:lem-C33}:]
From Lemma \ref{diagram:lemma-cubic-cubic} and trigonometric identities, it follows that 
\begin{align*}
\mathfrak{C}^{(3,3)}_{\leq N}(t)
&=  \sum_{n\in \Z^3} \bigg[ \frac{1}{\langle n \rangle^2}
   \int_0^t \dt^\prime \int_0^t \dt^{\prime\prime} \sin\Big((t-t^\prime) \langle n \rangle\Big)\sin\Big((t-t^{\prime\prime}) \langle n \rangle\Big) 
   \Gamma_{\leq N}(n,t^\prime-t^{\prime\prime}) \bigg] \\
 &=  \frac{1}{2} \sum_{n\in \Z^3} \bigg[ \frac{1}{\langle n \rangle^2}
   \int_0^t \dt^\prime \int_0^t \dt^{\prime\prime} \cos\Big((t^\prime-t^{\prime\prime}) \langle n \rangle\Big) 
   \Gamma_{\leq N}(n,t^\prime-t^{\prime\prime}) \bigg] \\
   &- \frac{1}{2} \sum_{n\in \Z^3} \bigg[ \frac{1}{\langle n \rangle^2}
   \int_0^t \dt^\prime \int_0^t \dt^{\prime\prime} 
   \cos\Big((2t-t^\prime-t^{\prime\prime}) \langle n \rangle\Big) 
   \Gamma_{\leq N}(n,t^\prime-t^{\prime\prime}) \bigg] \\
  &=: \mathfrak{C}_{\leq N}^{(3,3),a}(t) + \mathfrak{C}_{\leq N}^{(3,3),b}(t). 
\end{align*}
We first treat the contribution of $\mathfrak{C}_{\leq N}^{(3,3),a}$, which is the divergent part. By breaking the symmetry in $t^\prime$ and $t^{\prime \prime}$ and direct calculations, we obtain that 
\begin{align}
\mathfrak{C}_{\leq N}^{(3,3),a}(t)
&=   \sum_{n\in \Z^3} \bigg[ \frac{1}{\langle n \rangle^2}
   \int_0^t \dt^\prime \int_0^{t^\prime} \dt^{\prime\prime} \cos\Big((t^\prime-t^{\prime\prime}) \langle n \rangle\Big) 
   \Gamma_{\leq N}(n,t^\prime-t^{\prime\prime}) \bigg] \notag \\
 &=   \sum_{n\in \Z^3} \bigg[ \frac{1}{\langle n \rangle^2}
   \int_0^t \dt^\prime \int_0^{t^\prime} \dt^{\prime\prime} \cos\Big(t^{\prime\prime} \langle n \rangle\Big) 
   \Gamma_{\leq N}(n,t^{\prime\prime}) \bigg] \notag \\
 &=\sum_{n\in \Z^3} \bigg[ \frac{1}{\langle n \rangle^2}
    \int_0^{t} \dt^{\prime\prime} (t-t^{\prime \prime}) \cos\Big(t^{\prime\prime} \langle n \rangle\Big) 
   \Gamma_{\leq N}(n,t^{\prime\prime}) \bigg] \notag \\    
   &=\sum_{n\in \Z^3} \bigg[ \frac{1}{\langle n \rangle^2}
    \int_0^{t} \dt^{\prime} (t-t^{\prime }) \cos\Big(t^{\prime} \langle n \rangle\Big) 
   \Gamma_{\leq N}(n,t^{\prime}) \bigg]. \label{diagram:eq-C33-q1} 
\end{align}
By inserting the definition of $\Gamma_{\leq N}$ into \eqref{diagram:eq-C33-q1}, it follows that
\begin{align}
\mathfrak{C}_{\leq N}^{(3,3),a}(t) 
&= \int_0^t \dt^\prime (t-t^\prime) \bigg[ 
\sum_{n \in \Z^3} \sum_{\substack{n_1,n_2,n_3 \in \Z^3\colon \\ n=n_{123}}} 
1_{\leq N}(n) \frac{\cos( t^\prime \langle n \rangle)}{\langle n \rangle^2} 
\bigg( \prod_{j=1,2,3} 1_{\leq N}(n_j) \frac{\cos( t^\prime \langle n_j \rangle )}{\langle n_j \rangle^2} \bigg) \bigg]  \notag \\
&= \int_0^t \dt^\prime (t-t^\prime) 
\bigg\| \sum_{n\in \Z^3} 1_{\leq N}(n) \frac{\cos( t^\prime \langle n \rangle)}{\langle n \rangle^2}  e^{i \langle n,x \rangle} \bigg\|_{L_x^4(\T^3)}^4. \label{diagram:eq-C33-q2}
\end{align}
By inserting \eqref{diagram:eq-C33-q2} into the integral against the test-function $\rho$ from the statement of this lemma, it follows that
\begin{align*}
\int_{\R} \dt  \, \rho(t) \mathfrak{C}_{\leq N}^{(3,3),a}(t)  
&= \int_0^\infty \dt \int_0^{t} \dt^\prime  \rho(t)  (t-t^\prime) 
\bigg\| \sum_{n\in \Z^3} 1_{\leq N}(n) \frac{\cos( t^\prime \langle n \rangle)}{\langle n \rangle^2}  e^{i \langle n,x \rangle} \bigg\|_{L_x^4(\T^3)}^4 \\
&= \int_0^\infty \dt^\prime  \Big( \int_{t^\prime}^\infty \dt \, \rho(t) (t-t^\prime) \Big) 
\bigg\| \sum_{n\in \Z^3} 1_{\leq N}(n) \frac{\cos( t^\prime \langle n \rangle)}{\langle n \rangle^2}  e^{i \langle n,x \rangle} \bigg\|_{L_x^4(\T^3)}^4 \\
&\gtrsim \bigg\| \sum_{n\in \Z^3} 1_{\leq N}(n) \frac{\cos( t \langle n \rangle)}{\langle n \rangle^2}  e^{i \langle n,x \rangle} \bigg\|_{L_t^4 L_x^4([0,1/4]\times \T^3)}^4. 
\end{align*}
As a result, it only remains to prove that 
\begin{equation}\label{diagram:eq-C33-q3} 
\bigg\| \sum_{n\in \Z^3} 1_{\leq N}(n) \frac{\cos( t \langle n \rangle)}{\langle n \rangle^2}  e^{i \langle n,x \rangle} \bigg\|_{L_t^4 L_x^4([0,1/4]\times \T^3)}^4 \gtrsim c \log(N). 
\end{equation}
As already mentioned in Remark \ref{diagram:rem-strichartz}, \eqref{diagram:eq-C33-q3} is connected with (Besov-space versions of) the $L_t^4 L_x^4$-Strichartz estimate. Indeed, the sequence of function $\phi^{\cos}_{\leq N}$ given by 
\begin{equation*}
\widehat{\phi}^{\cos}_{\leq N}(n)= 1_{\leq N}(n) \langle n \rangle^{-2} 
\end{equation*}
is unbounded in $H_x^{1/2}$ but uniformly bounded in  the Besov space $B^{1/2}_{2,\infty}$. In order to prove \eqref{diagram:eq-C33-q3}, we let $(P_K^{\textup{sm}})_{K}$ be a sequence of smooth Littlewood-Paley operators (as opposed to the sharp frequency-cutoffs used in the rest of this paper). By first using the Littlewood-Paley square function estimate and then the embedding $\ell_K^2 \hookrightarrow \ell_K^4$, it follows that 
\begin{align*}
&\bigg\| \sum_{n\in \Z^3} 1_{\leq N}(n) \frac{\cos( t \langle n \rangle)}{\langle n \rangle^2}  e^{i \langle n,x \rangle} \bigg\|_{L_t^4 L_x^4([0,1/4]\times \T^3)}^4 \\ 
\gtrsim& \, \bigg\| P^{\textup{sm}}_K \bigg( \sum_{n\in \Z^3} 1_{\leq N}(n) \frac{\cos( t \langle n \rangle)}{\langle n \rangle^2}  e^{i \langle n,x \rangle}  \bigg) \bigg\|_{L_t^4 L_x^4 \ell_K^2 ([0,1/4]\times \T^3 \times 2^{\mathbb{N}_0})}^4 \\ 
\gtrsim& \, \bigg\| P^{\textup{sm}}_K \bigg( \sum_{n\in \Z^3} 1_{\leq N}(n) \frac{\cos( t \langle n \rangle)}{\langle n \rangle^2}  e^{i \langle n,x \rangle}  \bigg) \bigg\|_{L_t^4 L_x^4 \ell_K^4 ([0,1/4]\times \T^3 \times 2^{\mathbb{N}_0})}^4 \\
=& \, \sum_{K} \int_0^{1/4} \dt \int_{\T^3}  \dx \, \bigg| 
P^{\textup{sm}}_K \bigg( \sum_{n\in \Z^3} 1_{\leq N}(n) \frac{\cos( t \langle n \rangle)}{\langle n \rangle^2}  e^{i \langle n,x \rangle}  \bigg) \bigg|^4. 
\end{align*}
For all frequency-scales $K\ll N$, it holds that 
\begin{equation*}
 \bigg| 
P^{\textup{sm}}_K \bigg( \sum_{n\in \Z^3} 1_{\leq N}(n) \frac{\cos( t \langle n \rangle)}{\langle n \rangle^2}  e^{i \langle n,x \rangle}  \bigg) \bigg| \gtrsim \mathbf{1}\big\{ |t|,|x| \ll K^{-1} \big\} K. 
\end{equation*}
As a result, 
\begin{align*}
    &\sum_{K} \int_0^{1/4} \dt \int_{\T^3}  \dx \, \bigg| 
P^{\textup{sm}}_K \bigg( \sum_{n\in \Z^3} 1_{\leq N}(n) \frac{\cos( t \langle n \rangle)}{\langle n \rangle^2}  e^{i \langle n,x \rangle}  \bigg) \bigg|^4 \\
\gtrsim&\,  \sum_{K \ll N} K^4 \int_0^{1/4} \dt \int_{\T^3}  \dx \, \mathbf{1}\big\{ |t|,|x| \ll K^{-1} \big\} \\
\gtrsim&\, \sum_{K \ll N} 1 
\gtrsim \, \log(N). 
\end{align*}
This completes the proof of \eqref{diagram:eq-C33-q3}. It only remains to bound the contribution of $\mathfrak{C}_{\leq N}^{(3,3),b}$. By first breaking the symmetry in $t^\prime$ and $t^{\prime\prime}$ and then using the change of variables $t^{\prime\prime} \rightarrow t^\prime -t^{\prime\prime}$, we obtain that 
\begin{align}
 \mathfrak{C}_{\leq N}^{(3,3),b}(t)&=   - \sum_{n\in \Z^3} \bigg[ \frac{1}{\langle n \rangle^2}
   \int_0^t \dt^\prime \int_0^{t^\prime} \dt^{\prime\prime} 
   \cos\big((2t-t^\prime-t^{\prime\prime}) \langle n \rangle\big) 
   \Gamma_{\leq N}(n,t^\prime-t^{\prime\prime}) \bigg] \notag \\ 
   &= - \sum_{n\in \Z^3} \bigg[ \frac{1}{\langle n \rangle^2}
   \int_0^t \dt^\prime \int_0^{t^\prime} \dt^{\prime\prime} 
   \cos\big((2t-2t^\prime+t^{\prime\prime}) \langle n \rangle\big) 
   \Gamma_{\leq N}(n,t^{\prime\prime}) \bigg] \label{diagram:eq-C33-q4}. 
\end{align}
By performing the $t^\prime$-integral in \eqref{diagram:eq-C33-q4}, we obtain that
\begin{align*}
\eqref{diagram:eq-C33-q4} =  \frac{1}{2} \sum_{n\in \Z^3} \bigg[ \frac{1}{\langle n \rangle^3}
   \int_0^t \dt^{\prime\prime} 
   \Big( \sin\big(t^{\prime\prime}  \langle n \rangle\big)  - \sin\big( (2t-t^{\prime \prime}) \langle n \rangle \big) \Big)
   \Gamma_{\leq N}(n,t^{\prime\prime}) \bigg]. 
\end{align*}
By arguing exactly as in the proof of Lemma \ref{diagrams:lem-mathfrakC}, it follows that
\begin{equation*}
\sup_{N}\big| \mathfrak{C}^{(3,3),b}_{\leq N}(t)\big| \lesssim \langle t \rangle, 
\end{equation*}
which controls the  contribution of this term to \eqref{diagram:eq-C33}. 
\end{proof}

\section{Analytic aspects of basic stochastic diagrams}\label{section:analytic}

In Section \ref{section:diagrams}, 
we addressed algebraic aspects of the stochastic diagrams. More precisely, we derived explicit formulas and exhibited the \oftt-cancellation. In this section, we start our treatment of analytic aspects, such as regularity estimates. We first focus on basic stochastic objects, i.e., stochastic objects with degree less than or equal to six. The treatment of all higher-order objects, which requires some of the tools in Section \ref{section:linear} and Section \ref{section:para}, will be postponed until Section \ref{section:analytic2} below.

\subsection{The renormalization constant and multiplier}

Before we start with estimates of the stochastic objects themselves, we first examine the renormalization constant $\gamma_{\leq N}$ and the renormalization multiplier $\Gamma_{\leq N}$, which are as in Definition \ref{diagram:definition-Gamma-t}. Our estimates of $\gamma_{\leq N}$ and $\Gamma_{\leq N}$ will be used in the treatment of the quintic stochastic object (see Proposition \ref{analytic:prop-quintic}), the para-controlled operators (see Section \ref{section:para}), and higher-order objects (see Section \ref{section:analytic2}). 

In the next lemma, we prove that the difference of $\gamma_{\leq N}-\Gamma_{\leq N}$ remains bounded uniformly in $N$, even though $\gamma_{\leq N}$ and $\Gamma_{\leq N}$ individually diverge logarithmically. 

\begin{lemma}[Estimate of $\gamma_{\leq N}-\Gamma_{\leq N}$]\label{analytic:lem-difference-gamma}
For all frequency-scales $N$ and all frequencies $n\in \Z^3$, it holds that 
\begin{equation}\label{analytic:eq-difference-gamma}
\Big| \gamma_{\leq N} - \Gamma_{\leq N}(n) \Big| \lesssim \langle n \rangle^\epsilon. 
\end{equation}
\end{lemma}

\begin{proof}
We recall from Definition \ref{diagram:definition-Gamma-t} that 
\begin{align*}
\Gamma_{\leq N}(n) &= 6 \cdot 1_{\leq N}(n) 
\sum_{\substack{n_0,n_1,n_2,n_3\in \Z^3 \colon \\ n_{123}=n }}
\bigg[ \prod_{j=1}^3 1_{\leq N}(n_j) \langle n_j \rangle^{-2} \bigg], \\
\gamma_{\leq N} &= \Gamma_{\leq N}(0). 
\end{align*}
It trivially holds that $|\Gamma_{\leq N}(n)| \lesssim \log(N)$ for all $n \in \Z^3$. In particular, \eqref{analytic:eq-difference-gamma} easily holds for all $|n|\gtrsim N$, and it remains to treat the case $|n| \ll N$. To this end, we first decompose 
\begin{align}
&\Gamma_{\leq N}(n) - \gamma_{\leq N} \notag \\
=& 6 \sum_{N_1,N_2,N_3 \leq N}  
\sum_{\substack{n_1,n_2,n_3\in \Z^3}}
\bigg[ \big( \mathbf{1}\big\{ n_{123}=n \big\} - \mathbf{1} \big\{n_{123}=0 \big\}\big) 
\Big( \prod_{j=1}^3 1_{N_j}(n_j) \langle n_j \rangle^{-2} \Big) \bigg].  \label{analytic:eq-gamma-p1}
\end{align}
We treat the dyadic components in \eqref{analytic:eq-gamma-p1} separately. By symmetry, it suffices to treat the case $N_1 \geq N_2 \geq N_3$. Using the change of variables $n_3 \rightarrow n_3+n$ for the $\mathbf{1}\big\{ n_{123}=n\big\}$-term, we obtain that
\begin{align}
&\sum_{\substack{n_1,n_2,n_3\in \Z^3}}
\bigg[ \big( \mathbf{1}\big\{ n_{123}=n \big\} - \mathbf{1} \big\{n_{123}=0 \big\}\big) 
\Big( \prod_{j=1}^3 1_{N_j}(n_j) \langle n_j \rangle^{-2} \Big) \bigg] \notag \\
=& \sum_{\substack{n_1,n_2,n_3\in \Z^3}}
\bigg[  \mathbf{1} \big\{n_{123}=0 \big\} 
\Big( \prod_{j=1}^2 1_{N_j}(n_j) \langle n_j \rangle^{-2} \Big) \Big( 1_{N_3}(n_3-n)\langle n_3 - n \rangle^{-2} - 1_{N_3}(n_3) \langle n_3 \rangle^{-2}\Big)  \bigg]. \label{analytic:eq-gamma-p2}
\end{align}
Using Young's convolution inequality, it follows that 
\begin{align*}
\Big| \eqref{analytic:eq-gamma-p2} \Big| 
&\lesssim \big\| 1_{N_1}(n_1) \langle n_1 \rangle^{-2} \big\|_{\ell_{n_1}^\infty} 
\big\| 1_{N_2}(n_2) \langle n_2 \rangle^{-2} \big\|_{\ell_{n_2}^1}
\big\| 1_{N_3}(n_3-n)\langle n_3 - n \rangle^{-2} - 1_{N_3}(n_3) \langle n_3 \rangle^{-2} \big\|_{\ell_{n_3}^1}  \\
&\lesssim N_1^{-2} N_2  \langle n \rangle \lesssim N_1^{-1} \langle n \rangle. 
\end{align*}
Since the contribution of the dyadic component is only non-trivial for $|n| \lesssim N_1$, this yields the desired estimate. 
\end{proof}

In the following, we will often estimate time-integrals of the time-dependent renormalization multiplier from Definition \ref{diagram:definition-Gamma-t}. In these estimates, it is convenient to work with the following decomposition. 

\begin{definition}[Dyadic components of time-dependent renormalization multiplier]
\label{analytic:def-Gamma-dyadic}
For all frequency-scales $N_0,N_1,N_2,N_3$, we define 
\begin{equation}
\Gamma [N_\ast](n_0,t) := 6 \cdot 1_{N_0}(n_0) 
\sum_{\substack{n_0,n_1,n_2,n_3 \in \Z^3\colon \\ n_0 =n_{123}}}
\bigg[  \prod_{j=1}^3 \frac{1_{N_j}(n_j)}{\langle n_j \rangle^2} \cos\big( t \langle n_j \rangle \big) \bigg]. 
\end{equation}
\end{definition}
From Definition \ref{diagram:definition-Gamma-t} and Definition \ref{analytic:def-Gamma-dyadic}, it directly follows that 
\begin{equation*}
\Gamma_{\leq N}(n_0,t) = \sum_{N_0,N_1,N_2,N_3 \leq N} \Gamma [N_\ast](n_0,t). 
\end{equation*}

Equipped with Definition \ref{analytic:def-Gamma-dyadic}, we now state and prove the following integral estimates. 

\begin{lemma}[\protect{Estimate of $\Gamma [N_\ast]$}]\label{analytic:lem-Gamma-dyadic}
For all frequency-scales $N_0,N_1,N_2,N_3$, $n_0 \in \Z^3$, $t\in \R$, and $\lambda \in \R$, it holds that
\begin{equation}\label{analytic:eq-Gamma-dyadic-1}
\Big| \int_0^t \dt^\prime \, \Gamma [N_\ast](n_0,t-t^\prime) 
e^{i \lambda t^\prime} \Big| 
\lesssim \langle t \rangle \log(N_{\textup{max}}) \max\big( N_{\textup{max}}, \langle \lambda \rangle \big)^{-1}. 
\end{equation}
Furthermore, we have for all $\chi \in C^\infty_c(\R)$ that 
\begin{equation}\label{analytic:eq-Gamma-dyadic-2}
\Big| \int_\R \dt \, \chi(t) \Gamma [N_\ast](n_0,t) 
e^{i \lambda t} \Big| 
\lesssim_\chi \log(N_{\textup{max}}) \max\big( N_{\textup{max}}, \langle \lambda \rangle \big)^{-1}. 
\end{equation}
\end{lemma}

\begin{proof}
We only prove \eqref{analytic:eq-Gamma-dyadic-1}, since the proof of  \eqref{analytic:eq-Gamma-dyadic-2} is extremely similar. By inserting the definition of $\Gamma [N_\ast]$, computing the $t^\prime$-integral, and using a level-set decomposition,
it holds that 
\begin{align*}
&\Big| \int_0^t \dt^\prime \, \Gamma [N_\ast](n_0,t-t^\prime) 
e^{i \lambda t^\prime} \Big| \\
\lesssim& \, \langle t \rangle N_1^{-2} N_2^{-2} N_3^{-2}
\sum_{m\in \Z} \sum_{\substack{\pm_0,\pm_1,\\ \pm_2,\pm_3}}
\sum_{\substack{n_1,n_2,n_3 \in \Z^3\colon \\ n_0 = n_{123} }}
\bigg[ \Big( \prod_{j=0}^3 1_{N_j}(n_j) \Big) 
\mathbf{1}\Big\{\big | \sum_{j=0}^3 (\pm_j) \langle n_j \rangle - m \big|\leq 1 \Big\} 
(1+|\lambda-m|)^{-1} \bigg] \\
\lesssim& \, \langle t \rangle N_1^{-2} N_2^{-2} N_3^{-2} 
\Big( \sum_{\substack{m \in \Z \colon \\ |m| \lesssim N_{\textup{max}}}} (1+|\lambda-m|)^{-1} \Big) \\
&\times \sup_{m\in \Z} \sum_{\substack{\pm_0,\pm_1,\\ \pm_2,\pm_3}}
\sum_{\substack{n_1,n_2,n_3 \in \Z^3\colon \\ n_0 = n_{123} }}
\bigg[ \Big( \prod_{j=0}^3 1_{N_j}(n_j) \Big) 
\mathbf{1}\Big\{\big | \sum_{j=0}^3 (\pm_j) \langle n_j \rangle - m \big|\leq 1 \Big\} \bigg].
\end{align*}
The sum over $m\in \Z$ can be estimated by
\begin{equation*}
\sum_{\substack{m \in \Z \colon \\ |m| \lesssim N_{\textup{max}}}} (1+|\lambda-m|)^{-1}
\lesssim \log(N_{\textup{max}}) \min\Big( 1, N_{\textup{max}} \langle \lambda \rangle^{-1}\Big). 
\end{equation*}
Using Lemma \ref{counting:lem2}, the sum over $n_1,n_2,n_3 \in \Z^3$ can be estimated by 
\begin{equation*}
N_1^{-2} N_2^{-2} N_3^{-2} \sum_{\substack{\pm_0,\pm_1,\\ \pm_2,\pm_3}}
\sum_{\substack{n_1,n_2,n_3 \in \Z^3\colon \\ n_0 = n_{123} }}
\bigg[ \Big( \prod_{j=0}^3 1_{N_j}(n_j) \Big) 
\mathbf{1}\Big\{\big | \sum_{j=0}^3 (\pm_j) \langle n_j \rangle - m \big|\leq 1 \Big\} \bigg] 
\lesssim N_{\textup{max}}^{-1}. 
\end{equation*}
The combined estimate yields the desired conclusion. 
\end{proof}

\subsection{The linear and cubic stochastic objects}
In this subsection, we obtain regularity and Strichartz estimates for $\slinear[blue]$, $\slinear[green][\leqM]$, and $\scubic$, which verify the regularity claims from Figure \ref{figure:structure-NLW}.(A) in the introduction. We start with estimates for the linear evolutions $\slinear[blue]$ and $\slinear[green][\leqM]$. 

\begin{lemma}[Regularity of linear evolutions]\label{analytic:lem-linear}
For any $T\geq 1$ and $p\geq 2$, it holds that  
\begin{align}
\E \bigg[ \sup_N \Big\| \, \slinear[blue][\leqN] \Big\|_{(L_t^\infty \C_x^{-1/2-\epsilon} \cap X^{-1/2-\epsilon,b})([-T,T])}^p \bigg]^{1/p} \lesssim p^{1/2} T^\alpha. 
\end{align}
Furthermore, it holds that 
\begin{align}
\E \bigg[ \sup_{\substack{M,N\colon \\ M\geq N}}  \Big\| P_{\leq N} \slinear[green][\leqM] \Big\|_{(L_t^\infty \C_x^{1/2-\epsilon} \cap X^{1/2-\epsilon,b})([-T,T])}^p \bigg]^{1/p} \lesssim p^{3/2} T^\alpha. 
\end{align}
\end{lemma}

\begin{proof} 
From Lemma \ref{prep:lem-xsb}, it follows that 
\begin{equation*}
 \Big\| \, \slinear[blue][\leqN] \Big\|_{X^{-1/2-\epsilon,b}([-T,T])} 
 \lesssim T^\alpha \big\| \, \initial{blue}[][] \big\|_{\mathscr{H}_x^{-1/2-\epsilon}} \quad \text{and} \quad 
 \Big\| \, P_{\leq N} \slinear[green][\leqM] \Big\|_{X^{1/2-\epsilon,b}([-T,T])} 
 \lesssim T^\alpha \big\| \, \initial{green}[][\leqM] \big\|_{\mathscr{H}_x^{-1/2-\epsilon}}.
 \end{equation*}
Using Lemma \ref{diagram:lem-initial}, it follows that
\begin{align*}
\E \bigg[ \sup_N \Big\| \, \slinear[blue][\leqN] \Big\|_{ X^{-1/2-\epsilon,b}([-T,T])}^p \bigg]^{1/p} \lesssim p^{1/2} T^\alpha \quad \text{and} \quad \E \bigg[ \sup_{\substack{M,N\colon \\ M\geq N}} \Big\| \, P_{\leq N}\slinear[green][\leqM] \Big\|_{ X^{1/2-\epsilon,b}([-T,T])}^p \bigg]^{1/p} \lesssim p^{3/2} T^\alpha,
\end{align*}
which yields the desired estimates in $X^{s,b}$-spaces. 
Similar as in the proof of Lemma \ref{diagram:lem-initial}, the desired bound in $L_t^\infty \C_x^s$-spaces follows from a standard argument involving Gaussian hypercontractivity and translation-invariance. For this reason, we only treat the blue linear evolution $\slinear[blue]$. We let $q=q(\epsilon)\gg 1$ to be chosen and obtain, using Sobolev-embedding in space-time, that 
\begin{equation*}
\Big\| \, \slinear[blue][N] \Big\|_{L_t^\infty \C_x^{-1/2-\epsilon}([-T,T]\times \T^3)} \lesssim  T N^{4/q} \Big\| \langle \nabla \rangle^{-1/2-\epsilon}\,  \slinear[blue][N] \Big\|_{L_t^q L_x^q([-T,T]\times \T^3)}.
\end{equation*}
Using Minkowski's integral inequality and Gaussian hypercontractivity, we obtain for all $p\geq q$ that 
\begin{align*}
T N^{4/q}  \E \bigg[ \Big\| \langle \nabla \rangle^{-1/2-\epsilon}\,  \slinear[blue][N] \Big\|_{L_t^q L_x^q([-T,T]\times \T^3)}^p \bigg]^{1/p}
&\lesssim T N^{4/q}  \Big\| \langle \nabla \rangle^{-1/2-\epsilon}\,  \slinear[blue][N] \Big\|_{L_t^q L_x^q L_\omega^p([-T,T]\times \T^3 \times \Omega)} \\
&\lesssim p^{1/2} T N^{4/q}  \Big\| \langle \nabla \rangle^{-1/2-\epsilon}\,  \slinear[blue][N] \Big\|_{L_t^q L_x^q L_\omega^2([-T,T]\times \T^3 \times \Omega)}.
\end{align*}
Using spatial translation-invariance, Minkowski's integral inequality, and the embedding $X^{0,b}\hookrightarrow L_t^\infty L_x^2$, it follows that 
\begin{align*}
   p^{1/2} T N^{4/q}  \Big\| \langle \nabla \rangle^{-1/2-\epsilon}\,  \slinear[blue][N] \Big\|_{L_t^q L_x^q L_\omega^2([-T,T]\times \T^3 \times \Omega)}
   &= p^{1/2} T N^{4/q}  \Big\| \langle \nabla \rangle^{-1/2-\epsilon}\,  \slinear[blue][N] \Big\|_{L_t^q L_x^2 L_\omega^2([-T,T]\times \T^3 \times \Omega)} \\
   &\lesssim  p^{1/2} T N^{4/q}   \E \bigg[ \Big\| \langle \nabla \rangle^{-1/2-\epsilon}\,  \slinear[blue][N] \Big\|_{L_t^q L_x^2 ([-T,T]\times \T^3)}^2 \bigg]^{1/2}\\
   &\lesssim p^{1/2} T^{1+\frac{1}{q}} N^{4/q}   \E \bigg[ \Big\|\,  \slinear[blue][N] \Big\|_{X^{-1/2-\epsilon,b}([-T,T])}^2 \bigg]^{1/2}.
\end{align*}
By choosing $q=8\epsilon^{-1}$, using our previous $X^{s,b}$-estimate (with $\epsilon$ replaced by $\epsilon/2$), and summing over all dyadic scales, this yields the desired conclusion.
\end{proof}

We now turn from the linear evolutions to the cubic stochastic object.

\begin{lemma}[Regularity of cubic stochastic object]\label{analytic:lem-cubic}
For all $T\geq 1 $ and $p\geq 2$, it holds that 
\begin{align}\label{analytic:eq-cubic}
\E \bigg[ \sup_N \Big\| \, \scubic[\leqN] \Big\|_{(L_t^\infty \C_x^{-\epsilon} \cap X^{-\epsilon,b_+})([-T,T])}^p \bigg]^{1/p} \lesssim p^{3/2} T^\alpha. 
\end{align}
\end{lemma}

\begin{proof} We only treat the case $T=1$, since the general case follows from minor modifications, and we now restrict all statements to the space-time interval $[-1,1]\times \T^3$. We only prove the $X^{1/2-\epsilon,b_+}$-estimate since, similar as in the proof of Lemma \ref{analytic:lem-linear}, the $L_t^\infty \C_x^{1/2-\epsilon}$-estimate then follows from Gaussian hypercontractivity and translation-invariance.

Similar as in \eqref{diagram:eq-cubic-expanded}, the cubic stochastic object is given by 
\begin{equation}\label{analytic:eq-cubic-p1}
\begin{aligned}
  \scubic[\leqN] 
  =&
  \sum_{\substack{\varphi_1,\varphi_2,\varphi_3 \in \\ \{ \cos, \sin \}}} 
  \sum_{\substack{n_0,n_1,n_2,n_3 \in \Z^3 \colon \\ n_0=n_{123}}} 
  \bigg[
  \Big( \prod_{j=0}^3 1_{\leq N}(n_j) \langle n_j \rangle^{-1}
    \Big) e^{i\langle n_0 ,x \rangle} \\
    &\times \bigg( \int_0^t \sin\big((t-t^\prime) \langle n_0 \rangle\big)
    \prod_{j=1}^3 \varphi_j\big(t^\prime \langle n_j \rangle\big) \bigg)
  \SI[n_j, \varphi_j \colon 1 \leq j \leq 3 ] \bigg]. 
\end{aligned}
\end{equation}
Using a dyadic decomposition and the tensor from Corollary \ref{counting:cor-cubic}, we can write (\ref{analytic:eq-cubic-p1}) as
\begin{equation*}
\scubic[\leqN]  = \sum_{N_0,N_1,N_2,N_3 \leq N} \scubic \hspace{0ex} [N_\ast].
\end{equation*}
In this expression,
\begin{equation}\label{analytic:eq-cubic-p2}
\begin{aligned}
  \scubic \hspace{0ex} [N_\ast] 
  =& \sum_{\substack{\varphi_1,\varphi_2,\varphi_3 \in \\ \{ \cos, \sin \}}} 
  \sum_{\substack{n_0,n_1,n_2,n_3 \in \Z^3 }} 
  \bigg[
  H_{n_0 n_1 n_2 n_3}[N_\ast,\varphi_\ast](t) e^{i\langle n_0 , x \rangle}
  \SI[n_j, \varphi_j \colon 1 \leq j \leq 3 ] \bigg],
  \end{aligned}
\end{equation}

$N_\ast$ denotes the dependence on the frequency-scales $N_0,N_1,N_2$, and $N_3$, and $\varphi_\ast$ denotes the dependence on the phase-functions $\varphi_1$, $\varphi_2$, and $\varphi_3$. Using \eqref{analytic:eq-cubic-p2} and Corollary \ref{counting:cor-cubic}, it follows that 
\begin{align*}
\E \bigg[ \Big\| \scubic \hspace{0ex} [N_\ast] \Big\|_{X^{-\epsilon,b_+}}^2 \bigg] &\lesssim N_0^{-2\epsilon} \max_{\pm} \sum_{\substack{\varphi_1,\varphi_2,\varphi_3 \in \\ \{ \cos, \sin \}}} 
  \sum_{\substack{n_0,n_1,n_2,n_3 \in \Z^3 }} 
  \int_{\R} \dlambda \langle \lambda \rangle^{2b_+} \big| \widetilde{H}^{\pm}_{n_0 n_1 n_2 n_3}[N_\ast,\varphi_\ast](\lambda)\big|^2 \\
  &\lesssim N_{\textup{max}}^{-2\epsilon}. 
\end{align*}
Using Gaussian hypercontractivity, this directly implies the desired result. 
\end{proof}

\subsection{The product of blue and green linear evolutions}

We now turn to the product of $\slinear[blue]$ and $\slinear[green][\leqM]$, which is a quartic stochastic object. While this product does not occur in our Ansatz for the solution, it naturally occurs in the nonlinearity. As we saw in Lemma \ref{analytic:lem-linear}, the regularities of $\slinear[blue]$ and $\slinear[green][\leqM]$ are $-1/2-\epsilon$ and $1/2-\epsilon$, respectively. While the product $\slinear[blue]\cdot \slinear[green][\leqM]$ cannot therefore  be controlled using only the individual regularities, it can be controlled using more detailed information on the linear evolutions.
\begin{lemma}[\protect{Product estimate for $\slinear[blue]\, \cdot \, \slinear[green][\leqM]$}]\label{analytic:lem-resonant-blue-green}
For all $T\geq 1$, $p\geq 2$, and frequency-scales $K_1$ and $K_{2}$, it holds that 
\begin{align*}
\E \bigg[ \sup_M \Big\| \slinear[blue][K_1] \, P_{K_2} \slinear[green][\leqM] \Big\|_{L_t^\infty \C_x^{-1/2-\epsilon}([-T,T])}^p \bigg]^{1/p}
\lesssim p^2 T^\alpha K_2^{-1/2+\epsilon}. 
\end{align*}
\end{lemma}

Our choice to use $K_1$ and $K_2$ to denote the frequency-scales is deliberate. In the proof, we will replace $K_1$ and $K_2$ by $N_1$ and $N_{234}$, but in applications, we will replace $K_1$ and $K_2$ by $N_1$ and $N_2$.

\begin{proof}
We only treat the case $T=1$, since the general case follows from minor modifications. We also set $N_1:=K_1$ and $N_{234}:=K_2$. Using the definition of $\initial{blue}$ and $\initial{green}[][\leqM]$ (Definition \ref{diagram:definition-caloric}), we obtain that 
\begin{align}
&\slinear[blue][K_1] P_{K_2} \slinear[green][\leqM] \notag\\
=& \,  \sum_{\varphi_1 \in \{ \cos, \sin\}} \sum_{n_1 \in \Z^3} 
\bigg[ 1_{N_1}(n_1) \varphi_1\big( t \langle n_1 \rangle\big) 
e^{i \langle n_1,x \rangle} \SI[n_1,\varphi_1] \bigg] \notag \\
\times&  
\sum_{n_2,n_3,n_4 \in \Z^3} 
\bigg[ 1_{N_{234}}(n_{234})  \Big( \prod_{j=2}^4 1_{\leq M}(n_j) \Big) 
\cos\big( t \langle n_{234} \rangle\big) 
e^{i\langle n_{234}, x \rangle} \notag \\
&\hspace{10ex} \times \int_{(-\infty,0]^3} \bigotimes_{j=2}^4 \dW[\cos][s_j][n_j]
\int_{\max(s_2,s_3,s_4)}^0 \dtau e^{\tau \langle n_{234} \rangle^2} 
\Big( \prod_{j=2}^4 e^{-(\tau-s_j) \langle n_j \rangle^2} \Big) \bigg] \notag \allowdisplaybreaks[3]\\
=&  \sum_{\varphi_1 \in \{ \cos, \sin\}} \sum_{n_1,n_2,n_3,n_4 \in \Z^3}\bigg[
 1_{N_{234}}(n_{234})  1_{N_1}(n_1) \Big( \prod_{j=2}^4 1_{\leq M}(n_j) \Big) 
 \varphi_1\big( t \langle n_1 \rangle\big) \cos\big( t \langle n_{234} \rangle\big) 
 e^{i\langle n_{1234}, x \rangle} \label{analytic:eq-blue-green-p1} \\
 &\hspace{1ex}\times  \int\displaylimits_{(-\infty,0]^4} \dW[\varphi_1][s_1][n_1]\dW[\cos][s_2][n_2]
\dW[\cos][s_3][n_3]\dW[\cos][s_4][n_4] \, e^{s_1 \langle n_1 \rangle^2} \hspace{-3ex}
\int\displaylimits_{\max(s_2,s_3,s_4)}^0 \hspace{-3ex} \dtau e^{\tau \langle n_{234} \rangle^2} 
\Big( \prod_{j=2}^4 e^{-(\tau-s_j) \langle n_j \rangle^2} \Big) \bigg] \allowdisplaybreaks[3] \notag \\ 
+& 3  \sum_{n_1,n_2,n_3,n_4 \in \Z^3}\bigg[ \mathbf{1}\big\{ n_{12} = 0 \big\}
 1_{N_{234}}(n_{234})  1_{N_1}(n_1) \Big( \prod_{j=2}^4 1_{\leq M}(n_j) \Big) \langle n_2 \rangle^{-2}
 \cos\big( t \langle n_1 \rangle\big) \cos\big( t \langle n_{234} \rangle\big)  
\label{analytic:eq-blue-green-p2} \\
 &\hspace{1ex}\times  e^{i\langle n_{34}, x \rangle}  \int\displaylimits_{(-\infty,0]^2} \dW[\cos][s_3][n_3]
 \dW[\cos][s_4][n_4] 
\int\displaylimits_{\max(s_3,s_4)}^0 \hspace{-3ex} \dtau e^{\tau (\langle n_2 \rangle^2 +\langle n_{234} \rangle^2)} 
\Big( \prod_{j=3}^4 e^{-(\tau-s_j) \langle n_j \rangle^2} \Big) \bigg] \allowdisplaybreaks[3] \notag.
\end{align}
We now separately estimate the non-resonant and resonant parts \eqref{analytic:eq-blue-green-p1} and \eqref{analytic:eq-blue-green-p2} separately. \\

\emph{The non-resonant part \eqref{analytic:eq-blue-green-p1}:} We first dyadically localize the frequencies $n_2$, $n_3$, and $n_4$ in \eqref{analytic:eq-blue-green-p1}, which leads to 
\begin{equation}\label{analytic:eq-blue-green-p3}
\begin{aligned}
    &\sum_{\varphi_1 \in \{ \cos, \sin\}} \sum_{n_1,n_2,n_3,n_4 \in \Z^3}\bigg[
 1_{N_{234}}(n_{234})  \Big( \prod_{j=1}^4 1_{N_j}(n_j) \Big) 
 \varphi_1\big( t \langle n_1 \rangle\big) \cos\big( t \langle n_{234} \rangle\big) 
 e^{i\langle n_{1234}, x \rangle}  \\
 &\hspace{1ex}\times  \int\displaylimits_{(-\infty,0]^4} \dW[\varphi_1][s_1][n_1]\dW[\cos][s_2][n_2]
 \dW[\cos][s_3][n_3] \dW[\cos][s_4][n_4] \, e^{s_1 \langle n_1 \rangle^2} \hspace{-3ex}
\int\displaylimits_{\max(s_2,s_3,s_4)}^0 \hspace{-3ex} \dtau e^{\tau \langle n_{234} \rangle^2} 
\Big( \prod_{j=2}^4 e^{-(\tau-s_j) \langle n_j \rangle^2} \Big) \bigg]. 
\end{aligned}
\end{equation}
Using Gaussian hypercontractivity, it follows that 
\begin{align}
&\E \Big[ \big\| \eqref{analytic:eq-blue-green-p3} \big\|_{L_t^\infty \C_x^{-1/2-\epsilon}([-1,1])}^p \Big]^{2/p} \notag  \\
\lesssim& \,  p^4 N_{\textup{max}}^\epsilon \sum_{n_1,n_2,n_3,n_4 \in \Z^3} \bigg[ \langle n_{1234} \rangle^{-1-2\epsilon} 1_{N_{234}}(n_{234}) \Big( \prod_{j=1}^4 1_{N_j}(n_j) \Big) \notag \\
&\hspace{1ex} \times \int_{(-\infty,0]^4} \ds_1 \hdots \ds_4 \, e^{2s_1 \langle n_1 \rangle^2} \bigg( \int_{\max(s_2,s_3,s_4)}^0 \dtau e^{\tau \langle n_{234} \rangle^2} \prod_{j=2}^4 e^{-(\tau-s_j) \langle n_j \rangle^2} \bigg)^2 \bigg]. \label{analytic:eq-blue-green-p4} 
\end{align}
We note that 
\begin{equation*}
 \int_{\max(s_2,s_3,s_4)}^0 \dtau e^{\tau \langle n_{234} \rangle^2} \prod_{j=2}^4 e^{-(\tau-s_j) \langle n_j \rangle^2} 
 \lesssim  \max_{j=2,3,4} \int_{s_j}^0 \dtau  e^{-(\tau-s_j) \langle n_j \rangle^2} \lesssim \max(N_2,N_3,N_4)^{-2},
\end{equation*}
which implies
\begin{equation}\label{analytic:eq-blue-green-integral-estimate}
\begin{aligned}
&\bigg( \int_{\max(s_2,s_3,s_4)}^0 \dtau e^{\tau \langle n_{234} \rangle^2} \prod_{j=2}^4 e^{-(\tau-s_j) \langle n_j \rangle^2} \bigg)^2 \\
\lesssim& \max(N_2,N_3,N_4)^{-2} \bigg( \int_{\max(s_2,s_3,s_4)}^0 \dtau e^{\tau \langle n_{234} \rangle^2} \prod_{j=2}^4 e^{-(\tau-s_j) \langle n_j \rangle^2} \bigg).
\end{aligned}
\end{equation}
By first using \eqref{analytic:eq-blue-green-integral-estimate}, then performing the $s_1,s_2,s_3,s_4$ and $\tau$-integrals, summing in $n_1$, and finally summing in $n_2,n_3$, and $n_4$, we obtain that 
\begin{align*}
\eqref{analytic:eq-blue-green-p4} 
&\lesssim p^4 N_{\textup{max}}^\epsilon \max(N_2,N_3,N_4)^{-2} N_{234}^{-2}\, \Big( \prod_{j=1}^4 N_j^{-2} \Big) \\
&\times \sum_{n_1,n_2,n_3,n_4 \in \Z^3} \bigg[ \langle n_{1234} \rangle^{-1-2\epsilon} 1_{N_{234}}(n_{234}) \Big( \prod_{j=1}^4 1_{N_j}(n_j) \Big)  \bigg] \\
&\lesssim  p^4 N_{\textup{max}}^\epsilon N_1^{-2\epsilon} \max(N_2,N_3,N_4)^{-2} N_{234}^{-2}\, \Big( \prod_{j=2}^4 N_j^{-2} \Big) 
 \sum_{n_2,n_3,n_4 \in \Z^3} \bigg[ 1_{N_{234}}(n_{234}) \Big( \prod_{j=2}^4 1_{N_j}(n_j) \Big)  \bigg] \\
 &\lesssim  p^4 N_{\textup{max}}^\epsilon N_1^{-2\epsilon} \max(N_2,N_3,N_4)^{-1} \\
 &\lesssim  p^4  N_1^{-\epsilon} \max(N_2,N_3,N_4)^{-1+\epsilon}.
\end{align*}
This yields the desired estimate for the non-resonant portion \eqref{analytic:eq-blue-green-p4}. \\

\emph{The resonant part \eqref{analytic:eq-blue-green-p2}:} We first replace $n_1$ by $-n_2$ and use a dyadic localization in $n_2$, $n_3$, $n_4$, and $n_{34}$, which leads to 
\begin{equation}\label{analytic:eq-blue-green-p5}
\begin{aligned}
& \sum_{n_2,n_3,n_4 \in \Z^3}\bigg[
 1_{N_{234}}(n_{234})  1_{N_{34}}(n_{34}) \Big( \prod_{j=2}^4 1_{N_j}(n_j) \Big) \langle n_2 \rangle^{-2}
 \cos\big( t \langle n_2 \rangle\big) \cos\big( t \langle n_{234} \rangle\big)  
\\
 &\hspace{1ex}\times  e^{i\langle n_{34}, x \rangle}  \int\displaylimits_{(-\infty,0]^2} \dW[\cos][s_3][n_3]
 \dW[\cos][s_4][n_4] 
\int\displaylimits_{\max(s_3,s_4)}^0 \hspace{-3ex} \dtau e^{\tau (\langle n_2 \rangle^2 +\langle n_{234} \rangle^2)} 
\Big( \prod_{j=3}^4 e^{-(\tau-s_j) \langle n_j \rangle^2} \Big) \bigg]. 
\end{aligned}
\end{equation}
Due to the symmetry in $n_3$ and $n_4$, we can assume that $N_3\geq N_4$.  Using Gaussian hypercontractivity, we obtain that 
\begin{align}
&\E \Big[ \big\| \eqref{analytic:eq-blue-green-p5} \big\|_{L_t^\infty \C_x^{-1/2-\epsilon}([-1,1])}^p \Big]^{2/p} \notag \\
\lesssim& \, p^4 N_{\textup{max}}^\epsilon N_2^{-4} N_{34}^{-1-2\epsilon} 
\sum_{n_3,n_4 \in \Z^3}
\bigg[ 1_{N_{34}}(n_{34}) \Big( \prod_{j=3}^4 1_{N_j}(n_j)\Big) \notag \\
&\times 
\int_{(-\infty,0]^2} \ds_3 \ds_4 \bigg( \sum_{n_2\in \Z^3}
 1_{N_{234}}(n_{234})  1_{N_2}(n_2) \int_{\max(s_3,s_4)}^0 \dtau e^{\tau (\langle n_2 \rangle^2 + \langle n_{234} \rangle^2)} \prod_{j=3}^4 e^{-(\tau-s_j) \langle n_j \rangle^2} \bigg)^2 \bigg] \label{analytic:eq-blue-green-p6}
\end{align}
In the exponential factors in \eqref{analytic:eq-blue-green-p6}, we first use $\langle n_2 \rangle \geq N_2/2$ and $\langle n_{234} \rangle\geq N_{234}/2$. Then, the square of the sum over $n_2$ can be estimated by 
\begin{align*}
&\bigg( \sum_{n_2\in \Z^3}
 1_{N_{234}}(n_{234})  1_{N_2}(n_2) \int_{\max(s_3,s_4)}^0 \dtau e^{\tau (\langle n_2 \rangle^2 + \langle n_{234} \rangle^2)} \prod_{j=3}^4 e^{-(\tau-s_j) \langle n_j \rangle^2} \bigg)^2 \\
 \lesssim&\, N_2^6  \,\bigg( \int_{\max(s_3,s_4)}^0 \dtau e^{\frac{1}{4} \tau (N_2^2 + N_{234}^2)}  \prod_{j=3}^4 e^{-(\tau-s_j) \langle n_j \rangle^2}\bigg)^2 \\
 \lesssim&   \,N_2^6  \, \bigg( \int_{-\infty}^0 \dtau e^{\frac{1}{4} \tau N_2^2} \bigg) \bigg( \int_{\max(s_3,s_4)}^0 \dtau e^{\frac{1}{4} \tau (N_2^2 + N_{234}^2)}  \prod_{j=3}^4 e^{-(\tau-s_j) \langle n_j \rangle^2}\bigg) \\
  \lesssim&   \,N_2^4 \,  \bigg( \int_{\max(s_3,s_4)}^0 \dtau e^{\frac{1}{4} \tau (N_2^2 + N_{234}^2)}  \prod_{j=3}^4 e^{-(\tau-s_j) \langle n_j \rangle^2}\bigg).
\end{align*}
By integrating in $s_3$, $s_4$, and $\tau$, and finally summing over $n_3$ and $n_4$, it follows that 
\begin{align*}
\eqref{analytic:eq-blue-green-p6}
&\lesssim p^4 N_{\textup{max}}^\epsilon \max(N_2,N_{234})^{-2} N_{34}^{-1-2\epsilon} N_3^{-2}N_4^{-2} 
\sum_{n_3,n_4} \bigg[ 1_{N_{34}}(n_{34}) \Big( \prod_{j=3}^4 1_{N_j}(n_j)\Big) \bigg] \\
&\lesssim p^4 N_{\textup{max}}^\epsilon \max(N_2,N_{234})^{-2} N_{34}^{2-2\epsilon} N_3^{-1}. 
\end{align*}
Since $N_{34}\lesssim \max(N_2,N_{234})$ and, due to our assumption $N_3\geq N_4$, $N_{34}\lesssim N_3$, we obtain that 
\begin{align*}
    N_{\textup{max}}^\epsilon \max(N_2,N_{234})^{-2} N_{34}^{2-2\epsilon} N_3^{-1} 
    \lesssim \, N_{\textup{max}}^\epsilon \max(N_2,N_{234},N_3)^{-1} 
    \lesssim \, N_{\textup{max}}^{-1+\epsilon}.
\end{align*}
This contribution is (better than) acceptable. 
\end{proof}

\subsection{The quintic stochastic objects} 

In the previous subsections, we already estimated linear, cubic, and quartic stochastic objects. In this subsection, we estimate the quintic stochastic object
 $\squintic[\leqN]$.

\begin{proposition}[Regularity of the quintic stochastic object]\label{analytic:prop-quintic}
For all $T\geq 1 $ and $p\geq 2$, it holds that 
\begin{align}
\E \bigg[ \sup_N \Big\| \, \squintic[\leqN] \Big\|_{(L_t^\infty \C_x^{1/2-\epsilon} \cap X^{1/2-\epsilon,b})([-T,T])}^p \bigg]^{1/p} \lesssim p^{5/2} T^\alpha. 
\end{align}
\end{proposition}

In our estimate of the quintic stochastic object $\squintic[\leqN]$, we need to distinguish between the zero, one, and two pairing cases. To this end, we use the longhand diagrams from Subsection \ref{section:quintic-diagram}, 
and write the quintic stochastic diagram as 
\begin{align}
3 \squintic[\leqN](t,x) &= 3 ~ \scalebox{\medscale}{\quintic[\leqN][(t,x)][0]}+ 18 ~ \scalebox{\medscale}{\quintic[\leqN][(t,x)][1]}  + 18~   \scalebox{\medscale}{\resistor[][\leqN][(t,x)]}.  \label{analytic:eq-quintic-decomp-pure}
\end{align}
We recall that the third summand in \eqref{analytic:eq-quintic-decomp-pure} is called the resistor, which is also represented using the
 shorthand diagram $\sresistor$
 and has been defined in Definition \ref{diagram:definition-resistor}. 
 Our argument now splits into three steps, which address the three terms in \eqref{analytic:eq-quintic-decomp-pure}.

\begin{lemma}[No pairings]\label{analytic:lem-quintic-pure-nopairing}
For all $T\geq 1$ and $p\geq 2$, it holds that 
\begin{equation}\label{analytic:eq-quintic-pure-nopairing}
\E \left[ \sup_N \,  \left\| \, \scalebox{\medscale}{\quintic[\leqN][(t,x)][0]} \, \right\|_{( L_t^\infty \C_x^{1/2-\epsilon} \cap X^{1/2-\epsilon,b} )([-T,T])}^p \right]^{1/p} \lesssim \, p^{5/2} T^\alpha. 
\end{equation}
\end{lemma}

\begin{proof} 
We only treat the case $T=1$, since the general case follows from minor modifications and we restrict all of the statements below to the space-time interval $[-1,1]\times \T^3$. We only prove the $X^{1/2-\epsilon,b_+}$-estimate since, similar as in the proof of Lemma \ref{analytic:lem-linear}, the $L_t^\infty \C_x^{1/2-\epsilon}$-estimate then follows from Gaussian hypercontractivity and translation-invariance.

Similar as in \eqref{diagram:eq-quintic-pure-zero}, we write 
\begin{equation}\label{analytic:eq-quintic-pure-nopairing-p1}
\begin{aligned}
\scalebox{\medscale}{\quintic[\leqN][(t,x)][0]}
 =& \sum_{\substack{\varphi_1,\hdots, \varphi_5 \in \\ \{ \cos, \sin \}}}
 \sum_{\substack{n_0,n_1,\hdots,n_5 \in \Z^3 \colon \\ n_0 =n_{12345}}}
 \bigg[ \frac{1_{\leq N}(n_{234})}{\langle n_{234} \rangle} 
 \Big( \prod_{j=0}^5 \frac{1_{\leq N}(n_j)}{\langle n_j \rangle} \Big) e^{i \langle n_0 , x \rangle}      \\
 &\times \bigg( \int_0^t \dt^\prime
 \int_0^{t^\prime} \dt^{\prime\prime} 
 \sin\big((t-t^\prime) \langle n_0 \rangle \big)
 \sin \big( (t^\prime-t^{\prime\prime}) \langle n_{234} \rangle\big) 
 \varphi_1\big( t^\prime \langle n_1 \rangle \big) \\
 &\times  \Big(  
 \prod_{j=2,3,4} \varphi_j\big( t^{\prime\prime} \langle n_j \rangle\big)  \Big)
 \varphi_5 \big( t^\prime \langle n_5 \rangle \big) \bigg) \, \SI[n_j,\varphi_j \colon 1 \leq j \leq 5] \bigg]. 
\end{aligned}
\end{equation}
Using a dyadic decomposition and the tensor from Corollary \ref{counting:cor-quintic}, we write 
\begin{equation}
\scalebox{\medscale}{\quintic[\leqN][(t,x)][0]} = \sum_{N_0,N_1,\hdots,N_5 , N_{234}\leq N} 
\scalebox{\medscale}{\quintic[][(t,x)][0]} \hspace{1ex} [N_\ast],
\end{equation}
where 
\begin{align}\label{analytic:dyadic-quintic-tree}
\scalebox{\medscale}{\quintic[][(t,x)][0]} \hspace{1ex} [N_\ast]
 =& \sum_{\substack{\varphi_1,\hdots, \varphi_5 \in \\ \{ \cos, \sin \}}}
 \sum_{\substack{n_0,n_1,\hdots,n_5 \in \Z^3 }}
 \bigg[ H_{n_0 n_1 \hdots n_5}[N_\ast,\varphi_\ast](t) e^{i\langle n_0,x\rangle} 
 \SI[n_j,\varphi_j \colon 1 \leq j \leq 5] \bigg], 
\end{align}
$N_\ast$ denotes the dependence on $N_0,N_1,\hdots,N_5$, and $N_{234}$, and $\varphi_\ast$ denotes the dependence on $\varphi_1,\hdots,\varphi_5$. Using Corollary \ref{counting:cor-quintic}, we obtain that 
\begin{align}\label{analytic:eq-quintic-pure-nopairing-p3}
    \E \left[  \left\| \, \scalebox{\medscale}{\quintic[][(t,x)][0]} \hspace{1ex} [N_\ast] \, \right\|_{X^{1/2-\epsilon,b_+}}^2 \right] \lesssim \, N_0^{1-2\epsilon} N_{\textup{max}}^{-1+\epsilon} \lesssim N_{\textup{max}}^{-\epsilon}. 
\end{align} Together with Gaussian hypercontractivity, this implies the desired estimate. 
\end{proof}

We now turn to the  one-pairing case in which me make crucial use of the sine-cancellation.

\begin{lemma}[One pairing]\label{analytic:lem-quintic-pure-one-pairing}
For all $T\geq 1$ and $p\geq 2$, it holds that 
\begin{equation}\label{analytic:eq-quintic-pure-one-pairing}
\E \left[ \sup_N \,  \left\| \, \scalebox{\medscale}{\quintic[\leqN][(t,x)][1]} \, \right\|_{( L_t^\infty \C_x^{1/2-\epsilon} \cap X^{1/2-\epsilon,b_+} )([-T,T])}^p \right]^{1/p} \lesssim \, p^{3/2} T^\alpha. 
\end{equation}
\end{lemma}

\begin{proof} 
We only treat the case $T=1$, since the general case follows from minor modifications and we restrict all statements below to the space-time interval $[-1,1]\times \T^3$. 
We only prove the $X^{1/2-\epsilon,b_+}$-estimate since, similar as in the proof of Lemma \ref{analytic:lem-linear}, the $L_t^\infty \C_x^{1/2-\epsilon}$-estimate then follows from Gaussian hypercontractivity and translation-invariance.

The argument is essentially the same as in the proof of Lemma \ref{analytic:lem-quintic-pure-nopairing}, but uses the sine-cancellation tensor instead of the quintic tensor. Similar as in \eqref{diagram:eq-quintic-pure-one}, we write 
\begin{equation}\label{analytic:eq-quintic-pure-one-pairing-p1}
\begin{aligned}
 \scalebox{\medscale}{\quintic[\leqN][(t,x)][1]}
 =& \sum_{\substack{\varphi_3,\varphi_4, \varphi_5 \in \\ \{ \cos, \sin \}}}
 \sum_{\substack{n_0,n_1,\hdots,n_5 \in \Z^3 \colon \\ n_0 =n_{345}}}
 \bigg[ \mathbf{1}\big\{n_{12}=0 \big\} \frac{1_{\leq N}(n_{234})}{\langle n_{234} \rangle}  
 \Big( \prod_{j=0}^5 \frac{1_{\leq N}(n_j)}{\langle n_j \rangle} \Big) e^{i \langle n_0 , x \rangle}  \\
 &\times 
 \bigg( \int_0^t \dt^\prime \int_0^{t^\prime} \dt^{\prime\prime}
 \sin \big( (t-t^\prime) \langle n_{0} \rangle\big)
 \sin \big( (t^\prime-t^{\prime\prime}) \langle n_{234} \rangle\big) 
 \cos\big( (t^\prime-t^{\prime\prime}) \langle n_2 \rangle \big) \\
 &\times \Big( \prod_{j=3,4} \varphi_j\big( t^{\prime\prime} \langle n_j \rangle\big)
 \Big)  \varphi_5 \big( t^\prime \langle n_5 \rangle \big)  \bigg)
 \SI[n_j,\varphi_j \colon 3 \leq j \leq 5] \bigg]. 
\end{aligned}
\end{equation}
Using a dyadic decomposition and the tensor from Corollary \ref{counting:cor-sine-cancel}, we write 
\begin{equation}
 \scalebox{\medscale}{\quintic[\leqN][(t,x)][1]} = \sum_{N_0,N_1,\hdots,N_5 , N_{234}\leq N} 
 \scalebox{\medscale}{\quintic[][(t,x)][1]} \hspace{1ex} [N_\ast],
\end{equation}
where 
\begin{equation}\label{analytic:eq-quintic-pure-one-pairing-p2}
\begin{aligned}
\scalebox{\medscale}{\quintic[][(t,x)][1]} \hspace{1ex} [N_\ast]
 =&  \sum_{\substack{\varphi_3,\varphi_4, \varphi_5 \in \\ \{ \cos, \sin \}}}
 \sum_{\substack{n_0,n_1,\hdots,n_5 \in \Z^3 \colon \\ n_0 =n_{345}}}
 \bigg[ H^{\textup{sine}}_{n_0 n_3 n_4 n_5}[N_\ast,\varphi_\ast](t) e^{i \langle n_0 , x \rangle}  
 \SI[n_j,\varphi_j \colon 3 \leq j \leq 5] \bigg], 
\end{aligned}
\end{equation}
$N_\ast$ denotes the dependence on $N_0,N_1,\hdots,N_5$, and $N_{234}$, and $\varphi_\ast$ denotes the dependence on $\varphi_3,\varphi_4$, and $\varphi_5$. Using Corollary \ref{counting:cor-sine-cancel}, it follows that 
\begin{align}\label{analytic:eq-quintic-pure-one-pairing-p3}
    \E \left[  \left\| \,  \scalebox{\medscale}{\quintic[\leqN][(t,x)][1]} \, \right\|_{X^{1/2-\epsilon,b_+}}^2 \right] \lesssim \, N_0^{1-2\epsilon} N_{\textup{max}}^{-1+\epsilon} \lesssim N_{\textup{max}}^{-\epsilon}. 
\end{align}
Together with Gaussian hypercontractivity, this implies the desired estimate. 
\end{proof}

It now only remains to treat the resistor $\sresistor$, which is the (renormalized) two-pairing term in the quintic object. 

\begin{lemma}[\protect{The resistor $\sresistor[\leqN]$}]\label{analytic:lem-resistor}
For all $T\geq 1$ and $p\geq 2$, it holds that 
\begin{equation}
\E \bigg[ \sup_N \Big\| \hspace{1ex} \sresistor[\leqN] \Big\|_{(L_t^\infty \C_x^{1/2-\epsilon} \cap X^{1/2-\epsilon,b})([-T,T])}^p \bigg]^{1/p} 
\lesssim p^{1/2} T^\alpha. 
\end{equation}
\end{lemma}

\begin{proof}
We only treat the case $T=1$, since the general case follows from minor modifications. We only prove the $X^{1/2-\epsilon,b_+}$-estimate since, similar as in the proof of Lemma \ref{analytic:lem-linear}, the $L_t^\infty \C_x^{1/2-\epsilon}$-estimate then follows from Gaussian hypercontractivity and translation-invariance.

From the definition of $\sresistor[\leqN]$ (Definition \ref{diagram:definition-resistor}), it follows that 
\begin{equation*}
18 \, \, \sresistor[\leqN] = \Duh \left[ 18\, \, \quinticnl[\leqN][][2] - \Gamma_{\leq N} \,  \linear{blue}[\leqN] \right]. 
\end{equation*}
From Lemmma \ref{diagram:lem-quintic-resonant} and Definition \ref{analytic:def-Gamma-dyadic}, it follows that 
\begin{align}
&18\, \, \quinticnl[\leqN][][2] - \Gamma_{\leq N} \,  \linear{blue}[\leqN] \notag \\
=& - \sum_{n_0 \in \Z^3} \bigg[ \langle n_0 \rangle^{-1} \Gamma_{\leq N}(n_0,t) e^{i\langle n_0 , x \rangle} \SI[n_0,\cos] \bigg] \notag \\
-& \sum_{\varphi \in \{ \cos, \sin \}} \sum_{n_0 \in \Z^3} 
\bigg[ \Big( \int_0^t \dt^\prime \Gamma_{\leq N}(n_0,t-t^\prime) (\partial_t \varphi) \big( t^\prime \langle n_0 \rangle\big) \Big) 
e^{i \langle n_0 ,x \rangle} \SI[n_0,\varphi] \bigg]  \notag \\
=& - \sum_{N_0,N_1,N_2,N_3\leq N} \sum_{n_0 \in \Z^3} \bigg[ \langle n_0 \rangle^{-1} \Gamma[N_\ast](n_0,t) e^{i\langle n_0 , x \rangle} \SI[n_0,\cos] \bigg] \label{analytic:eq-resistor-p1} \\ 
-& \sum_{\varphi \in \{ \cos, \sin \}}  \sum_{N_0,N_1,N_2,N_3\leq N}
\sum_{n_0 \in \Z^3} 
\bigg[ \Big( \int_0^t \dt^\prime \Gamma[N_\ast](n_0,t-t^\prime) (\partial_t \varphi) \big( t^\prime \langle n_0 \rangle\big) \Big) 
e^{i \langle n_0 ,x \rangle} \SI[n_0,\varphi] \bigg]. \label{analytic:eq-resistor-p2}
\end{align}
Due to Gaussian hypercontractivity and Lemma \ref{prep:lem-xsb}, it suffices to prove that 
\begin{align}
& \E \bigg[ \Big\| \sum_{n_0 \in \Z^3} \Big[ \langle n_0 \rangle^{-1} \Gamma[N_\ast](n_0,t) e^{i\langle n_0 , x \rangle} \SI[n_0,\cos] \Big] \Big\|_{X^{-1/2-\epsilon,b_+-1}}^2 \bigg] \label{analytic:eq-resistor-p3} \\
+&  \E \bigg[ \Big\| \sum_{n_0 \in \Z^3} 
\Big[ \Big( \int_0^t \dt^\prime \Gamma[N_\ast](n_0,t-t^\prime) (\partial_t \varphi) \big( t^\prime \langle n_0 \rangle\big) \Big) 
e^{i \langle n_0 ,x \rangle} \SI[n_0,\varphi] \Big] \Big\|_{X^{-1/2-\epsilon,b_+-1}}^2 \bigg] 
\label{analytic:eq-resistor-p4} \\
\lesssim& \, N_{\textup{max}}^{-\epsilon}. \notag 
\end{align}
Using Lemma \ref{analytic:lem-Gamma-dyadic}, the first term \eqref{analytic:eq-resistor-p3} is estimated by 
\begin{align*}
\eqref{analytic:eq-resistor-p3}
&\lesssim \sum_{\pm_0} \sum_{n_0 \in \Z^3}
\langle n_0 \rangle^{-3-2\epsilon} \int_\R \dlambda \, \langle \lambda \rangle^{2(b_+-1)}
\Big| \int_{\R} \chi(t) \Gamma[N_\ast](n_0,t) e^{i(\pm_0 \langle n_0 \rangle + \lambda) t}
\Big|^2 \\
&\lesssim N_{\textup{max}}^{-2+\epsilon} N_0^{-2\epsilon},
\end{align*}
which is acceptable. Similarly, using Lemma \ref{analytic:lem-Gamma-dyadic}, the second term \eqref{analytic:eq-resistor-p4} is estimated by
\begin{align*}
\eqref{analytic:eq-resistor-p4}
&\lesssim \E \bigg[ \Big\| \sum_{n_0 \in \Z^3} 
\Big[ \Big( \int_0^t \dt^\prime \Gamma[N_\ast](n_0,t-t^\prime) (\partial_t \varphi) \big( t^\prime \langle n_0 \rangle\big) \Big) 
e^{i \langle n_0 ,x \rangle} \SI[n_0,\varphi] \Big] \Big\|_{L_t^2 H_x^{-1/2-\epsilon}}^2 \bigg] \\ 
&\lesssim \sum_{n_0 \in \Z^3}
\langle n_0 \rangle^{-1-2\epsilon} \Big| \int_0^t \dt^\prime \Gamma[N_\ast](n_0,t-t^\prime) (\partial_t \varphi) \big( t^\prime \langle n_0 \rangle\big) \Big) \Big|^2 \\
&\lesssim N_0^{2-2\epsilon} N_{\textup{max}}^{-2+\epsilon}  \\ 
&\lesssim N_{\textup{max}}^{-\epsilon}, 
\end{align*}
which is acceptable. 
\end{proof}

Equipped with the estimates of the three summands in \eqref{analytic:eq-quintic-decomp-pure}, we now prove the main estimate of this subsection. 

\begin{proof}[Proof of Proposition \ref{analytic:prop-quintic}:]
This follows directly from the decomposition in \eqref{analytic:eq-quintic-decomp-pure} and the estimates in Lemma \ref{analytic:lem-quintic-pure-nopairing}, Lemma \ref{analytic:lem-quintic-pure-one-pairing}, and Lemma \ref{analytic:lem-resistor}. 
\end{proof}

\subsection{The heat-wave quintic object}\label{section:analytic-quintic-mixed}

In this subsection, we treat the heat-wave quintic object, which relates to (frequency-localized) versions of $\squadratic[\leqN] P_{\leq N} \slinear[green][\leqM]$. 
As discussed in Section \ref{section:introduction} and Subsection \ref{section:diagram-mixed}, the reason for referring to $\squadratic[\leqN] P_{\leq N} \slinear[green][\leqM]$ as a heat-wave stochastic object is as follows: The green caloric initial data $\initial{green}[][\leqM]$ is defined as the cubic Picard iterate of the stochastic heat equation. As a result, the product $\squadratic[\leqN] P_{\leq N} \slinear[green][\leqM]$  contains nonlinear interactions involving both heat and wave propagators. 

\begin{proposition}[Heat-wave quintic object]\label{analytic:prop-mixed}
For all $T\geq 1$ and $p\geq 2$, it holds that 
\begin{align*}
&\E \bigg[ \sup_{\substack{M,N\colon \\ M \geq N}} \bigg\| P_{\leq N}  \bigg[  \squadratic[\leqN] P_{\leq N} \slinear[green][\leqM] - \Big( 2 \HLL + \HHL \Big) \Big( \slinear[blue][\leqN], \slinear[blue][\leqN], P_{\leq N} \slinear[green][\leqM] \Big) \bigg] \bigg\|_{X^{-1/2+\delta_2,b_+-1}([-T,T])}^p \bigg]^{1/p} \\
\lesssim&\, p^{5/2} T^\alpha. 
\end{align*}
\end{proposition}

The following computations are more involved than in previous subsections. The reason is that, since the interplay of heat and wave propagators is only needed in this subsection, we did not previously set up all of the corresponding tensor estimates. \\

In our proof of Lemma \ref{analytic:prop-mixed}, we utilize  the long-hand diagrams from Section \ref{section:diagrams}. To this end, we first decompose
\begin{align}
&P_{\leq N}  \bigg[  \squadratic[\leqN] P_{\leq N} \slinear[green][\leqM] - \Big( 2 \HLL + \HHL \Big) \Big( \slinear[blue][\leqN], \slinear[blue][\leqN], P_{\leq N} \slinear[green][\leqM] \Big) \bigg] \notag \\
=& \sum_{\substack{N_0,N_1,N_{234},N_5 \leq N\colon \\ N_{234} > \max(N_1,N_5)}} 
P_{N_0} \bigg[ \lcol \slinear[blue][N_1] \slinear[blue][N_5] \rcol P_{N_{234}} \slinear[green][\leqM] \bigg]. 
\end{align}
Arguing similarly as in Subsection \ref{section:diagram-mixed}, we then write
\begin{equation}\label{analytic:eq-mixed-decomposition-012}
\begin{aligned}
&\sum_{\substack{N_0,N_1,N_{234},N_5 \leq N\colon \\ N_{234} > \max(N_1,N_5)}} 
P_{N_0} \bigg[ 
\lcol \slinear[blue][N_1] \slinear[blue][N_5] \rcol 
P_{N_{234}} \slinear[green][\leqM] \bigg] \\
=& 
\sum_{\substack{N_0,N_1,N_{234},N_5 \leq N\colon \\ N_{234} > \max(N_1,N_5)}}  \sum_{N_2,N_3,N_4 \leq M}  
\left[ 
\scalebox{0.8}{\quinticmixednl[][][0]} \hspace{0ex} [N_\ast] \right. \\
&+ \left. 6\scalebox{0.8}{\quinticmixednl[][][1]} \hspace{0ex} [N_\ast]  
+ 6 \scalebox{0.8}{\quinticmixednl[][][2]} \hspace{0ex} [N_\ast] 
\right],
\end{aligned}
\end{equation}
where the frequency-localized stochastic diagrams are defined as follows: 
The heat-wave quintic object with no pairing is defined as
\begin{equation}\label{analytic:eq-mixed-0-algebraic}
\begin{aligned}
&\scalebox{0.8}{\quinticmixednl[][(t,x)][0]} \hspace{0ex} [N_\ast] \allowdisplaybreaks[4] \\
=& 2^{5/2} 
\sum_{ \substack{ \varphi_1, \hdots, \varphi_5 \\ \in \{ \cos, \sin\}}}
\sum_{\substack{n_0,n_1,\hdots,n_5 \colon \\ n_0 = n_{12345} }} 
 \bigg[ 1_{N_{234}}(n_{234}) \Big( \prod_{j=0}^5 1_{N_j}(n_j) \Big) \Big( \prod_{j=2}^4 \mathbf{1}\{ \varphi_j = \cos \} \Big) e^{i \langle n_0 ,x \rangle} \  \\
&\times \varphi_1\big( t \langle n_1 \rangle\big) \cos\big( t \langle n_{234} \rangle \big) \varphi_5\big( t \langle n_5 \rangle\big) \\
&\times \int\displaylimits_{(-\infty,0]^5} \bigotimes_{j=1}^5 \dW[\varphi_j][s_j][n_j] \int\displaylimits_{\max(s_2,s_3,s_4)}^0 \mathrm{d}\tau \Big( \prod_{j=1,5} e^{s_j \langle n_j \rangle^2} \Big) e^{\tau  \langle n_{234} \rangle^2} \prod_{j=2,3,4} e^{-(\tau-s_j) \langle n_j \rangle^2} \bigg].
\end{aligned}
\end{equation}
The heat-wave quintic object with one pairing is defined as
\begin{equation}\label{analytic:eq-mixed-1-algebraic}
\begin{aligned}
&\scalebox{0.8}{\quinticmixednl[][(t,x)][1]} \hspace{0ex} [N_\ast] \allowdisplaybreaks[4] \\
=& 2^{3/2} \mathbf{1}\big\{ N_1 = N_2 \big\} 
\sum_{ \substack{ \varphi_3,\varphi_4, \varphi_5 \\ \in \{ \cos, \sin\}}}
\sum_{\substack{n_0,n_3,n_4,n_5 \colon \\ n_0 = n_{345} }} \sum_{n_1,n_2 \in \Z^3} 
 \bigg[ \mathbf{1}\big\{ n_{12}=0 \big\} 1_{N_{234}}(n_{234}) \Big( \prod_{j=0}^5 1_{N_j}(n_j) \Big) \\
&\times  e^{i \langle n_0 ,x \rangle} \Big( \prod_{j=3}^4 \mathbf{1}\{ \varphi_j = \cos \} \Big)  \cos\big( t \langle n_1 \rangle\big) \cos\big( t \langle n_{234} \rangle \big) \varphi_5\big( t \langle n_5 \rangle\big) \\
&\times \int\displaylimits_{(-\infty,0]^3} \bigotimes_{j=3}^5 \dW[\varphi_j][s_j][n_j] \int\displaylimits_{\max(s_3,s_4)}^0 \mathrm{d}\tau e^{\tau (\langle n_2 \rangle^2 + \langle n_{234} \rangle^2)} \prod_{j=3,4} e^{-(\tau-s_j) \langle n_j \rangle^2} e^{s_5\langle n_5 \rangle^2} \bigg].
\end{aligned}
\end{equation}
Finally, the heat-wave quintic object with two pairings is defined as 
\begin{equation}\label{analytic:eq-mixed-2-algebraic}
\begin{aligned}
&\scalebox{0.8}{\quinticmixednl[][(t,x)][2]} \hspace{0ex} [N_\ast] \allowdisplaybreaks[3] \\
=& 2^{1/2} \mathbf{1}\big\{ N_0 = N_3 \big\} \mathbf{1}\big\{ N_1 = N_2 \big\} \mathbf{1} \big\{ N_4 = N_5 \big\}\\
&\times 
\sum_{\substack{n_0,n_1,\hdots n_5 \in \Z^3 \colon \\ n_0 =n_3}} \bigg[ \mathbf{1}\big\{ n_{12}=n_{45}=0\big\} 1_{N_{234}}(n_{234}) \Big( \prod_{j=0}^5 1_{N_j}(n_j) \Big) \langle n_2 \rangle^{-2} \langle n_4 \rangle^{-2} e^{i\langle n_0 ,x \rangle} 
 \\
&\times \cos\big( t \langle n_1 \rangle \big) \cos\big( t \langle n_{234}\rangle \big) 
\cos\big( t \langle n_5 \rangle\big) \int_{-\infty}^0 \dW[\cos][s_3][n_3] 
\int_{s_3}^0 \dtau e^{\tau (\langle n_2 \rangle^2 + \langle n_{234} \rangle^2 + \langle n_4 \rangle^2)} e^{-(\tau-s_3) \langle n_3 \rangle^2} \bigg] . 
\end{aligned}
\end{equation}

We estimate the contributions of \eqref{analytic:eq-mixed-0-algebraic}, \eqref{analytic:eq-mixed-1-algebraic}, and \eqref{analytic:eq-mixed-2-algebraic} separately. The first lemma addresses the heat-wave quintic object with no pairing.  

\begin{lemma}[Heat-wave quintic object with no pairing]\label{analytic:lem-mixed-0}
For all $T\geq 1$, $p\geq 2$, and frequency-scales $N_0,N_1,\hdots,N_5,N_{234}$ satisfying $N_{234} \geq \max(N_1,N_5)^\eta$, it holds that 
\begin{equation}\label{analytic:eq-mixed-0-estimate}
\E \left[ \, \left\| \scalebox{0.8}{\quinticmixednl[][(t,x)][0]} \hspace{0ex} [N_\ast] 
 \right\|_{X^{-1/2+\delta_2,b_+-1}([-T,T])}^p \right]^{1/p} 
 \lesssim p^{5/2} T^\alpha N_{\textup{max}}^{-\epsilon}. 
\end{equation}
\end{lemma}

\begin{proof}
We only treat the case $T= 1$, since the general case follows from minor modifications. By Gaussian hypercontractivity, it suffices to prove \eqref{analytic:eq-mixed-0-estimate} for $p=2$. From \eqref{analytic:eq-mixed-0-algebraic}, we obtain for all $t\in [-T,T]$ that 
\begin{align*}
&\scalebox{0.8}{\quinticmixednl[][(t,x)][0]} \hspace{0ex} [N_\ast] \\
&=     2^{5/2} \sum_{\substack{n_0,n_1,\hdots,n_5 \colon \\ n_0 = n_{12345} }} 
\sum_{ \substack{ \varphi_1, \hdots, \varphi_5 \\ \in \{ \cos, \sin\}}} \bigg[  \Big( \prod_{j=2}^4 1_{N_j}(n_j) \Big) \Big( \prod_{j=2}^4 \mathbf{1}\{ \varphi_j = \cos \} \Big)  
e^{i \langle n_0 ,x \rangle} h_{n_0 n_1 n_{234}n_5}(t) \\
&\times \int\displaylimits_{(-\infty,0]^5} \bigotimes_{j=1}^5 \dW[\varphi_j][s_j][n_j] \int\displaylimits_{\max(s_2,s_3,s_4)}^0 \mathrm{d}\tau \Big( \prod_{j=1,5} e^{s_j \langle n_j \rangle^2} \Big) e^{\tau  \langle n_{234} \rangle^2} \prod_{j=2,3,4} e^{-(\tau-s_j) \langle n_j \rangle^2}  \bigg],
\end{align*}
where 
\begin{align*}
&h_{n_0 n_1 n_{234}n_5}(t) \\
=& h_{n_0 n_1 n_{234}n_5}[N_0,N_1,N_{234},N_5,\varphi_1,\varphi_5](t) \\
=& 1_{N_{234}}(n_{234}) \Big( \prod_{j=0,1,5} 1_{N_j}(n_j) \Big) \mathbf{1}\big\{ n_0 = n_1 + n_{234} +n_5\big\} \chi\big(t\big) \varphi_1\big( t \langle n_1 \rangle\big) \cos\big( t \langle n_{234} \rangle \big) \varphi_5\big( t \langle n_5 \rangle\big). 
\end{align*}
We also let $\widetilde{h}^{\pm_0}_{n_0 n_1 n_{234} n_5}$ be the twisted Fourier transform from   \eqref{prelim:eq-twisted-FT-h}.
Due to the orthogonality of the multiple stochastic integrals, we obtain that
\begin{align}
&\E \left[ \, \left\| \scalebox{0.8}{\quinticmixednl[][(t,x)][0]} \hspace{0ex} [N_\ast] 
 \right\|_{X^{-1/2+\delta_2,b_+-1}([-1,1])}^2 \right]  \notag \\
 \lesssim&\, N_0^{-1+2\delta_2} \sum_{\pm_0} \sum_{ \substack{\varphi_1,\varphi_5 \\ \in \{ \cos, \sin \}}} \sum_{n_0,n_1,\hdots,n_5 \in \Z^3} \bigg[ \prod_{j=2}^4 1_{N_j}(n_j) \int_{\R} \dlambda \bigg(  \langle \lambda \rangle^{2 (b_+-1)} \Big| \widetilde{h}^{\pm_0}_{n_0 n_1 n_{234} n_5}(\lambda)\Big|^2 \notag \\
 &\times \int_{(-\infty,0]^5} \ds_1 \hdots \ds_5 \prod_{j=1,5} e^{2 s_j \langle n_j \rangle^2} \Big( \int\displaylimits_{\max(s_2,s_3,s_4)}^0 \dtau \prod_{j=2,3,4} e^{-(\tau-s_j) \langle n_j \rangle^2} e^{\tau \langle n_{234} \rangle^2} \Big)^2 \bigg) \bigg] \label{analytic:eq-mixed-0-p1}
\end{align}
By first using \eqref{analytic:eq-blue-green-integral-estimate} and then integrating in $s_1,\hdots,s_5$ and $\tau$, we obtain that 
\begin{equation}\label{analytic:eq-mixed-0-p2}
\begin{aligned}
   & \int_{(-\infty,0]^5} \ds_1 \hdots \ds_5 \prod_{j=1,5} e^{2 s_j \langle n_j \rangle^2} \Big( \int\displaylimits_{\max(s_2,s_3,s_4)}^0 \dtau \prod_{j=2,3,4} e^{-(\tau-s_j) \langle n_j \rangle^2} e^{\tau \langle n_{234} \rangle^2} \Big)^2  \\
   \lesssim&\,  \max(N_2,N_3,N_4)^{-2} N_{234}^{-2} \prod_{j=1}^5 N_j^{-2}. 
\end{aligned}
\end{equation}
After fixing $n_{234}$, the remaining possible choices of $n_2,n_3$, and $n_4$ are bounded by 
\begin{equation*}
    \med(N_2,N_3,N_4)^3 \min(N_2,N_3,N_4)^3 \lesssim N_2^2 N_3^2 N_4^2.
\end{equation*} 
Together with \eqref{analytic:eq-mixed-0-p2} and Lemma \ref{counting:cubic_tensor}, this implies 
\begin{align*}
\eqref{analytic:eq-mixed-0-p1} 
&\lesssim \,  \max(N_2,N_3,N_4)^{-2} N_{234}^{-2} N_0^{-1+2\delta_2} \prod_{j=1}^5 N_j^{-2} \\
&\times \sum_{\pm_0} \sum_{ \substack{\varphi_1,\varphi_5 \\ \in \{ \cos, \sin \}}} \sum_{n_0,n_1,\hdots,n_5 \in \Z^3} \bigg[ \prod_{j=2}^4 1_{N_j}(n_j) \int_{\R} \dlambda  \langle \lambda \rangle^{2 (b_+-1)} \Big| \widetilde{h}^{\pm_0}_{n_0 n_1 n_{234} n_5}(\lambda)\Big|^2 \bigg] \\
&\lesssim  \max(N_2,N_3,N_4)^{-2} N_{234}^{-2} N_0^{-1+2\delta_2}  N_1^{-2} N_5^{-2}\\
&\times \sum_{\pm_0} \sum_{ \substack{\varphi_1,\varphi_5 \\ \in \{ \cos, \sin \}}} \sum_{n_0,n_1,n_{234},n_5 \in \Z^3} \bigg[ \prod_{j=2}^4 1_{N_j}(n_j) \int_{\R} \dlambda   \langle \lambda \rangle^{2 (b_+-1)} \Big| \widetilde{h}^{\pm_0}_{n_0 n_1 n_{234} n_5}(\lambda)\Big|^2  \bigg] \\
&\lesssim N_{\textup{max}}^{2(\delta_2+\epsilon)} \min(N_0,N_1,N_{234},N_5) \max(N_2,N_3,N_4)^{-2} \\
&\lesssim N_{\textup{max}}^{2(\delta_2+\epsilon)} \max(N_2,N_3,N_4)^{-1}.
\end{align*}
Due to our assumption, it holds that $\max(N_2,N_3,N_5) \gtrsim N_{234} \gtrsim \max(N_1,N_5)^\eta$. Since $\eta \gg \delta_2,\epsilon$, this yields the desired estimate. 
\end{proof}

We now address the heat-wave quintic object with one pairing, i.e., \eqref{analytic:eq-mixed-1-algebraic}. 

\begin{lemma}[Heat-wave quintic object with one pairing]
\label{analytic:lem-mixed-1}
For all $T\geq 1$, $p\geq 2$, and frequency-scales $N_0,N_1,\hdots,N_5,N_{234}$ satisfying $N_{234} \geq \max(N_1,N_5)^\eta$, it holds that 
\begin{equation}\label{analytic:eq-mixed-1-estimate}
\E \left[ \, \left\| \scalebox{0.8}{\quinticmixednl[][(t,x)][1]} \hspace{0ex} [N_\ast] 
 \right\|_{X^{-1/2+\delta_2,b_+-1}([-T,T])}^p \right]^{1/p} 
 \lesssim p^{3/2} T^\alpha N_{\textup{max}}^{-\epsilon}. 
\end{equation}
\end{lemma}

\begin{proof} 
We only treat the case $T= 1$ and $p=2$,  since the general case follows from either minor modifications or  Gaussian hypercontractivity,  respectively. Due to the pairing in $n_1$ and $n_2$, we implicitly restrict to $N_1=N_2$. From \eqref{analytic:eq-mixed-1-algebraic}, we obtain for all $t\in [-1,1]$ that 
\begin{align*}
    &\scalebox{0.8}{\quinticmixednl[][(t,x)][1]} \hspace{0ex} [N_\ast] \allowdisplaybreaks[4] \\
=& 2^{3/2}  
\sum_{ \substack{ \varphi_3,\varphi_4, \varphi_5 \\ \in \{ \cos, \sin\}}} 
\sum_{\substack{n_0,n_3,n_4,n_5 \in \Z^3 \colon \\ n_0 = n_{345} }} 
\sum_{n_1,n_2 \in \Z^3} 
\bigg[ \mathbf{1}\big\{ n_{12}=0 \big\}  \Big( \prod_{j=3}^4 1_{N_j}(n_j) \Big) \\
&\times  e^{i \langle n_0 ,x \rangle} \Big( \prod_{j=3}^4 \mathbf{1}\{ \varphi_j = \cos \} \Big)  h_{n_0 n_1 n_{234} n_5}(t) \\
&\times \int\displaylimits_{(-\infty,0]^3} \bigotimes_{j=3}^5 \dW[\varphi_j][s_j][n_j] \int\displaylimits_{\max(s_3,s_4)}^0 \mathrm{d}\tau e^{\tau (\langle n_2 \rangle^2 + \langle n_{234} \rangle^2)} \prod_{j=3,4} e^{-(\tau-s_j) \langle n_j \rangle^2} e^{s_5\langle n_5 \rangle^2} \bigg],
\end{align*}
where 
\begin{align*}
&h_{n_0n_1 n_{234}n_5}(t) \\
=& h_{n_0n_1 n_{234}n_5}[N_0,N_1,N_{234},N_5,\varphi_5](t) \\
=&\, 1_{N_{234}}(n_{234}) \Big(\prod_{j=0,1,5}1_{N_j}(n_j) \Big)
\mathbf{1}\big\{ n_0 = n_1 + n_{234} +n_5 \big\} \chi(t) \cos\big( t \langle n_1 \rangle\big) \cos\big( t \langle n_{234} \rangle \big) \varphi_5\big( t \langle n_5 \rangle\big). 
\end{align*}
We also let $\widetilde{h}^{\pm_0}_{n_0 n_1 n_{234} n_5}$ be the twisted Fourier transform from   \eqref{prelim:eq-twisted-FT-h}.
Using the orthogonality of the multiple stochastic integrals, it holds that 
\begin{align}
&\E \left[ \, \left\| \scalebox{0.8}{\quinticmixednl[][(t,x)][1]} \hspace{0ex} [N_\ast] 
 \right\|_{X^{-1/2+\delta_2,b_+-1}([-1,1])}^2 \right]  \notag \\
 \lesssim& \,  N_0^{-1+2\delta_2} 
 \sum_{\pm_0} \sum_{\varphi_5 \in \{ \cos, \sin \}} 
 \sum_{n_0,n_3,n_4,n_5 \in \Z^3} \bigg[ \Big( \prod_{j=3}^4 1_{N_j}(n_j) \Big) \notag \\
 \times& \int_{\R} \dlambda
  \int\displaylimits_{(-\infty,0]^3} \ds_3\ds_4 \ds_5 \, \langle \lambda \rangle^{2(b_+-1)} 
 \bigg( \sum_{n_1,n_2 \in \Z^3} \mathbf{1}\big\{ n_{12}=0 \big\} \langle n_2 \rangle^{-2} 
 \big|\widetilde{h}^{\pm_0}_{n_0 n_1 n_{234} n_5}(\lambda)\big| \label{analytic:eq-mixed-1-p1}\\
 & \hspace{4ex}\times
 \int\displaylimits_{\max(s_3,s_4)}^0 \dtau e^{\tau (\langle n_2 \rangle^2 + \langle n_{234} \rangle^2)} \prod_{j=3,4} e^{-(\tau-s_j) \langle n_j \rangle^2} e^{s_5\langle n_5 \rangle^2} \bigg)^2 \bigg]. \notag 
\end{align}
Under the frequency-restrictions in \eqref{analytic:eq-mixed-1-p1}, it holds that 
\begin{align*}
&\int\displaylimits_{\max(s_3,s_4)}^0 \dtau e^{\tau (\langle n_2 \rangle^2 + \langle n_{234} \rangle^2)} \prod_{j=3,4} e^{-(\tau-s_j) \langle n_j \rangle^2} e^{s_5\langle n_5 \rangle^2} \lesssim \, \int\displaylimits_{\max(s_3,s_4)}^0 \dtau e^{\frac{\tau}{4} (N_2^2 + N_{234}^2)} \prod_{j=3,4} e^{-\frac{(\tau-s_j)}{4} N_j^2} e^{ \frac{s_5}{4} N_5^2},
\end{align*}
where the right-hand side only depends on the frequency-scales (but not the frequency-variables). This allows us pull the $\tau$-integral out of the sum over $n_1,n_2 \in \Z^3$ and perform the $s_3,s_4,s_5$, and $\tau$-integrals. As a result, we obtain that 
\begin{align}
\eqref{analytic:eq-mixed-1-p1}
&\lesssim \,  \max(N_2,N_3,N_4)^{-2} \max(N_2,N_{234})^{-2} N_0^{-1+2\delta_2} N_3^{-2} N_4^{-2} N_5^{-2} \notag \\
&\times
 \sum_{\pm_0} \sum_{\varphi_5 \in \{ \cos, \sin \}} 
 \sum_{n_0,n_3,n_4,n_5 \in \Z^3} \bigg[  \Big( \prod_{j=3}^4 1_{N_j}(n_j) \Big) \notag \\
 &\times \int_{\R} \dlambda
   \langle \lambda \rangle^{2(b_+-1)} 
 \bigg( \sum_{n_1,n_2 \in \Z^3} \mathbf{1}\big\{ n_{12}=0 \big\} \langle n_2 \rangle^{-2} 
 \big|\widetilde{h}^{\pm_0}_{n_0 n_1 n_{234} n_5}(\lambda)\big|\bigg)^2 \bigg] \label{analytic:eq-mixed-1-p2}. 
\end{align}
We now use Cauchy-Schwarz\footnote{This estimate can be slightly improved by using dispersive effects, but would only yield a gain in $|n_{34}|$, which would not significantly affect the rest of the argument.} in $n_2\in \Z^3$, which costs us a factor of $N_2^3$, and yields that 
\begin{align}
\eqref{analytic:eq-mixed-1-p2}
&\lesssim \,  \max(N_2,N_3,N_4)^{-2} \max(N_2,N_{234})^{-2} N_0^{-1+2\delta_2} N_2^{-1} N_3^{-2} N_4^{-2} N_5^{-2}  \notag \\
&\times
 \sum_{\pm_0} \sum_{\varphi_5 \in \{ \cos, \sin \}} 
 \sum_{\substack{n_0,n_1,n_2,n_3,\\n_4,n_5 \in \Z^3}} \bigg[ \mathbf{1}\big\{ n_{12}=0 \big\}  \Big( \prod_{j=3}^4 1_{N_j}(n_j) \Big)  \notag \\
 &\times \int_{\R} \dlambda
   \langle \lambda \rangle^{2(b_+-1)} 
 \big|\widetilde{h}^{\pm_0}_{n_0 n_1 n_{234} n_5}(\lambda)\big|^2  \bigg]. 
 \label{analytic:eq-mixed-1-p3}
\end{align}
We now first perform the $\lambda$-integral and then finally insert the constraint $n_1 =-n_2$. Together with a level-set decomposition of the multilinear dispersive symbol as in Section \ref{section:counting}, this leads to 
\begin{align}
 \eqref{analytic:eq-mixed-1-p3}
 &\lesssim \,  N_{\textup{max}}^{2\epsilon} \max(N_2,N_3,N_4)^{-2} \max(N_2,N_{234})^{-2} N_0^{-1+2\delta_2} N_2^{-1} N_3^{-2} N_4^{-2} N_5^{-2}  \notag \\
&\times \sup_{m\in\Z}
\sum_{\substack{\pm_0, \pm_2, \\ \pm_{234}, \pm_5 }} 
 \sum_{\substack{n_0,n_2,n_3,\\n_4,n_5 \in \Z^3}}
 \bigg[ \mathbf{1}\big\{n_0=n_{345}\big\} 1_{N_{234}}(n_{234}) 
 \Big( \prod_{j=0,2,3,4,5} 1_{N_j}(n_j) \Big)  \notag \\
 &\times 
 \mathbf{1}\big\{ \big|\pm_0 \langle n_0 \rangle \pm_2 \langle n_2 \rangle \pm_{234} \langle n_{234} \rangle \pm_5 \langle n_5 \rangle - m \big|\leq 1\big\} \bigg].\label{analytic:eq-mixed-1-p4}
\end{align}
By symmetry in $n_3$ and $n_4$, we may now assume that $N_3\geq N_4$. Furthermore, as long as our final estimate has room in $N_{\textup{max}}$, we may further restrict \eqref{analytic:eq-mixed-1-p4} to the regime $|n_{34}|\sim N_{34}$. By using Lemma \ref{counting:lem1} to perform the summation over $n_2 \in \Z^3$, it follows that 
\begin{align}
\eqref{analytic:eq-mixed-1-p4}\Big|_{|n_{34}|\sim N_{34}} 
&\lesssim  \,  N_{\textup{max}}^{2\epsilon} \max(N_2,N_3,N_4)^{-2} \max(N_2,N_{234})^{-2} N_0^{-1+2\delta_2} N_2^{-1} N_3^{-2} N_4^{-2} N_5^{-2}  \notag \\
&\times \min(N_{2},N_{34})^{-1}N_2^3  
 \sum_{\substack{n_0,n_3,\\n_4,n_5 \in \Z^3}}
 \bigg[ \mathbf{1}\big\{n_0=n_{345}\big\}  1_{N_{34}}(n_{34})
 \Big( \prod_{j=0,3,4,5} 1_{N_j}(n_j) \Big) \bigg]. \label{analytic:eq-mixed-1-p5}
\end{align}
The  remaining sum over $n_0,n_3,n_4,$ and $n_5$ in \eqref{analytic:eq-mixed-1-p5} is estimated by
\begin{align*}
     \sum_{\substack{n_0,n_3,\\n_4,n_5 \in \Z^3}}
 \bigg[ \mathbf{1}\big\{n_0=n_{345}\big\}  1_{N_{34}}(n_{34})
 \Big( \prod_{j=0,3,4,5} 1_{N_j}(n_j) \Big) \bigg] \lesssim N_{34}^3 N_4^3 \min(N_0,N_5)^3. 
\end{align*}
As a result, we obtain that
\begin{align}
\eqref{analytic:eq-mixed-1-p5} 
\lesssim& \,  N_{\textup{max}}^{2(\delta_2+\epsilon)} \max(N_2,N_3,N_4)^{-2} \max(N_2,N_{234})^{-2}  \min(N_2,N_{34})^{-1} N_2^{2} N_{34}^3 N_3^{-1} \notag  \\
\lesssim& \,  N_{\textup{max}}^{2(\delta_2+\epsilon)} \max(N_2,N_3,N_4)^{-2} \max(N_2,N_{234})^{-2}  \max(N_2,N_{34}) N_2 N_{34}^2 N_3^{-1}. 
\label{analytic:eq-mixed-1-p6}
\end{align}
In our estimate of \eqref{analytic:eq-mixed-1-p6}, we now use that 
\begin{alignat*}{3}
\max(N_2,N_{34}) \max(N_2,N_{234})^{-1}&\lesssim 1, 
\quad& N_{34} \max(N_2,N_{234})^{-1}&\lesssim 1, \\
 N_{34} N_3^{-1} &\lesssim 1, \quad \text{and} \quad&  N_2 \max(N_2,N_3,N_4)^{-1} &\lesssim 1.
\end{alignat*}
As a result, it follows that 
\begin{equation*}
\eqref{analytic:eq-mixed-1-p6} \lesssim N_{\textup{max}}^{2(\delta_2+\epsilon)} \max(N_2,N_3,N_4)^{-1}.
\end{equation*}
Due to our assumption $N_{234} \gtrsim \max(N_1,N_5)^\eta$, this yields the desired estimate. 
\end{proof}

We now address the heat-wave quintic object with two pairings, i.e., \eqref{analytic:eq-mixed-2-algebraic}. 

\begin{lemma}[Heat-wave quintic object with two pairings]
\label{analytic:lem-mixed-2}
For all $T\geq 1$, $p\geq 2$, and frequency-scales $N_0,N_1,\hdots,N_5,N_{234}$, it holds that 
\begin{equation}\label{analytic:eq-mixed-2-estimate}
\E \left[ \, \left\| \scalebox{0.8}{\quinticmixednl[][(t,x)][2]} \hspace{0ex} [N_\ast] 
 \right\|_{X^{-1/2+\delta_2,b_+-1}([-T,T])}^p \right]^{1/p} 
 \lesssim p^{1/2} T^\alpha N_{\textup{max}}^{-\epsilon}. 
\end{equation}
\end{lemma}

In contrast to Lemma \ref{analytic:lem-mixed-0} and Lemma \ref{analytic:lem-mixed-1}, Lemma \ref{analytic:lem-mixed-2} does not require a lower bound on $N_{234}$. The reason is that both factors originating from the blue caloric initial data have been paired. 

\begin{proof}
It suffices to treat the case $T=1$ and $p=2$, since the general case follows from minor modifications or Gaussian hypercontractivity, respectively.

In the explicit formula \eqref{analytic:eq-mixed-2-algebraic}, we first remove the sum over $n_0,n_1$, and $n_5$ by inserting the restrictions 
 $n_0=n_3$, $n_1=-n_2$, and $n_5=-n_4$. To simplify the notation, we further fix $N_0=N_3$, $N_1=N_2$, and $N_5=N_4$. For all $t\in [-1,1]$, this implies 
 \begin{align*}
  \scalebox{0.8}{\quinticmixednl[][(t,x)][2]} \hspace{0ex} [N_\ast] 
  =& \, 2^{1/2} \sum_{n_2,n_3,n_4 \in \Z^3} \bigg[ \langle n_2 \rangle^{-2} \langle n_4 \rangle^{-2} e^{i\langle n_3 , x\rangle} h_{n_3 n_2 n_{234} n_4}(t) \\
  &\times 
   \int_{-\infty}^0 \dW[\cos][s_3][n_3] 
\int_{s_3}^0 \dtau e^{\tau (\langle n_2 \rangle^2 + \langle n_{234} \rangle^2 + \langle n_4 \rangle^2)} e^{-(\tau-s_3) \langle n_3 \rangle^2} \bigg],
 \end{align*}
 where 
 \begin{align*}
    &h_{n_3 n_2 n_{234} n_4}(t) \\
    =&\, 1_{N_{234}}(n_{234}) \Big( \prod_{j=2}^4 1_{N_j}(n_j) \Big) \mathbf{1}\big\{ n_{234} =n_2 + n_3 + n_4 \big\} \chi(t)  \cos\big( t \langle n_2 \rangle \big) \cos\big( t \langle n_{234}\rangle \big) 
\cos\big( t \langle n_4 \rangle\big). 
 \end{align*}
 We also let     $\widetilde{h}^{\pm_3}_{n_3 n_2 n_{234} n_4}(\lambda_3)$ be the twisted Fourier transform from   \eqref{prelim:eq-twisted-FT-h}.
 It follows that 
 \begin{align}
  &  \E \left[ \, \left\| \scalebox{0.8}{\quinticmixednl[][(t,x)][2]} \hspace{0ex} [N_\ast] 
 \right\|_{X^{-1/2+\delta_2,b_+-1}}^2 \right] \notag \\
 \lesssim&\, N_2^{-4} N_3^{-1+2\delta} N_4^{-4} \sum_{\pm_3} \sum_{n_3\in \Z^3} \int_{-\infty}^0 \ds_3 \int_{\R} \dlambda_3 \bigg[ \langle \lambda_3 \rangle^{2(b_+-1)} \notag \\
&\times  \bigg( \sum_{n_2,n_4 \in \Z^3} \big|  
 \widetilde{h}^{\pm_3}_{n_3 n_2 n_{234} n_4}(\lambda_3) \big| 
 \int_{s_3}^0 \dtau  e^{\tau (\langle n_2 \rangle^2 + \langle n_{234} \rangle^2 + \langle n_4 \rangle^2)} e^{-(\tau-s_3) \langle n_3 \rangle^2} \bigg)^2 \bigg]. 
 \label{analytic:eq-mixed-2-p1}
 \end{align}
 We now estimate the $\tau$-integral by 
 \begin{equation*}
      \int_{s_3}^0 \dtau  e^{\tau (\langle n_2 \rangle^2 + \langle n_{234} \rangle^2 + \langle n_4 \rangle^2)} e^{-(\tau-s_3) \langle n_3 \rangle^2}
      \lesssim  \int_{s_3}^0 \dtau  e^{\frac{\tau}{4} (N_2^2 + N_{234}^2 +N_4^2)} e^{-\frac{(\tau-s_3)}{4} N_3^2}, 
 \end{equation*}
 which depends only on the frequency-scales but not the frequency-variables. As a result, we can pull the $\tau$-integral out of the sum over $n_2,n_4 \in \Z^3$ and perform the $\tau$ and $s_3$-integrals. Since
 \begin{align*}
  &\int_{-\infty}^0 \ds_3 \bigg( \int_{s_3}^0 \dtau  e^{\frac{\tau}{4} (N_2^2 + N_{234}^2 +N_4^2)} e^{-\frac{(\tau-s_3)}{4} N_3^2} \bigg)^2 \\
  \lesssim& \, \max(N_2,N_{234},N_4)^{-2} \int_{-\infty}^0 \ds_3  \int_{s_3}^0 \dtau  e^{\frac{\tau}{4} (N_2^2 + N_{234}^2 +N_4^2)} e^{-\frac{(\tau-s_3)}{4} N_3^2} \\
  \lesssim&\, \max(N_2,N_{234},N_4)^{-4} N_3^{-2}
 \end{align*}
 and $\max(N_2,N_{234},N_4)\sim \max(N_2,N_3,N_4)$, it follows that 
 \begin{align}
    \eqref{analytic:eq-mixed-2-p1} 
    &\lesssim \max(N_2,N_3,N_4)^{-4} N_2^{-4} N_3^{-3+2\delta_2} N_4^{-4} \notag \\
     &\times \sum_{\pm_3} \sum_{n_3\in \Z^3} \int_{\R} \dlambda_3 \bigg[ \langle \lambda_3 \rangle^{2(b_+-1)} 
 \bigg( \sum_{n_2,n_4 \in \Z^3} \big|  
 \widetilde{h}^{\pm_3}_{n_3 n_2 n_{234} n_4}(\lambda_3) \big| \bigg)^2 \bigg].
 \label{analytic:eq-mixed-2-p2}
 \end{align}
 We now use Cauchy-Schwarz\footnote{While this application of Cauchy-Schwarz is rather crude and a better estimate can be obtained through dispersive effects, it is sufficient for the proof of this lemma.} in $n_2,n_4\in \Z^3$ and then apply Lemma \ref{counting:cubic_tensor} (see also Lemma \ref{counting:basetensor_est}), which yields 
 \begin{align}
    \eqref{analytic:eq-mixed-2-p2} &\lesssim  \max(N_2,N_3,N_4)^{-4} N_2^{-1} N_3^{-3+2\delta_2} N_4^{-1} \notag \\
    &\times  \sum_{\pm_3} \sum_{n_2,n_3,n_4\in \Z^3} \int_{\R} \dlambda_3 \langle \lambda_3 \rangle^{2(b_+-1)} \big| \widetilde{h}^{\pm_3}_{n_3 n_2 n_{234} n_4}(\lambda_3)  \big|^2 \notag \\
    &\lesssim  \max(N_2,N_3,N_4)^{-4} N_2^{-1} N_3^{-3+2\delta_2} N_4^{-1} 
  \times \Big( \med(N_2,N_3,N_4)^{-1} N_2^3 N_3^3 N_4^3 \Big) \notag \\
  &\lesssim \max(N_2,N_3,N_4)^{-4} \med(N_2,N_3,N_4)^{-1} N_2^2 N_3^{2\delta_2} N_4^2 \label{analytic:eq-mixed-2-p3}. 
 \end{align}
 Since $\med(N_2,N_3,N_4) \gtrsim \min(N_2,N_4)$, we obtain that
 \begin{equation*}
    \eqref{analytic:eq-mixed-2-p3} \lesssim  \max(N_2,N_3,N_4)^{-4} \max(N_2,N_4) N_2 N_3^{2\delta_2} N_4 \lesssim \max(N_2,N_3,N_4)^{-1+2\delta_2}. 
 \end{equation*}
 This yields the desired estimate. 
\end{proof}

Equipped with our previous lemmas, we can now complete our analysis of the heat-wave stochastic object. 

\begin{proof}[Proof of Proposition \ref{analytic:prop-mixed}]
The desired estimate directly follows from the decomposition \eqref{analytic:eq-mixed-decomposition-012} as well as Lemma \ref{analytic:lem-mixed-0}, Lemma \ref{analytic:lem-mixed-1}, and Lemma \ref{analytic:lem-mixed-2}.
\end{proof}

\subsection{A sextic stochastic object} In this subsection we prove the following propositions, which concerns the regularity of two sextic stochastic objects, namely the products of two cubic objects and one linear and one quintic object.
\begin{proposition}\label{analytic:prop-sextic1} For any dyadic scales $N_1\leq N_2$, we have
\begin{equation}\label{analytic:eq-sextic1}\Big\|P_{N_1}\scubic[\leqN]\cdot P_{N_2}\scubic[\leqN]-\mathfrak{C}_{\leq N}^{(3,3)}[N_1,N_2]\Big\|_{L_\omega^p(C_t^0C_x^{-1/10}([-T,T]))}\lesssim p^{3}T^\alpha N_2^{-100\nu}.
\end{equation}
\end{proposition}
\begin{proof} Define $\Mc_{3,3}:=P_{N_1}\scubic[\leqN]\cdot P_{N_2}\scubic[\leqN]-\mathfrak{C}_{\leq N}^{(3,3)}[N_1,N_2]$. Start by choosing a parameter between $1$ and $\nu$ (say $\sqrt{\nu}$) and considering the projections $P_{\leq N_2^{\sqrt{\nu}}}\Mc_{3,3}$ and $P_{\geq N_2^{\sqrt{\nu}}}\Mc_{3,3}$. By Lemma \ref{analytic:lem-cubic}, we have
\[\big\|P_{N_j}\scubic[\leqN]\big\|_{L_\omega^{2p}(C_t^0C_x^0([-T,T]))}\lesssim (2p)^{3/2}T^\alpha N_2^\epsilon\] for $1\leq j\leq 2$, and therefore 
\[\Big\|P_{N_1}\scubic[\leqN]\cdot P_{N_2}\scubic[\leqN]\Big\|_{L_\omega^p(C_t^0C_x^0([-T,T]))}\lesssim p^3T^\alpha N_2^\epsilon,\]which then implies (\ref{analytic:eq-sextic1}) for $P_{\geq N_2^{\sqrt{\nu}}}\Mc_{3,3}$, thus below we will consider $P_{\leq N_2^{\sqrt{\nu}}}\Mc_{3,3}$ only (hence $N_1\sim N_2$). By using the reduction arguments in Subsection \ref{prep:remark-reduction}, it suffices to bound $\Eb|\Mc_{3,3}(t,x)|^2\lesssim N_{2}^{-100\sqrt{\nu}}$ for fixed values of $(t,x)$.

As in Lemma \ref{analytic:lem-cubic}, we may decompose $P_{N_1}\scubic[\leqN]$ into $\scubic\,[M_1,M_2,M_3,M_{123}=N_1]$ and decompose $P_{N_2}\scubic[\leqN]$ into $\scubic\,[M_4,M_5,M_6,M_{456}=N_2]$. We may assume $M_+:=\max(M_1,\cdots M_6)\leq N_2^{100}$, otherwise we can gain an extra power of $M_+$ using (\ref{counting:cubic_bound3}) (with $N_{123}$ and $N_{\mathrm{max}}$ replaced by $N_2$ and $M_+$) and argue as in the above paragraph. Now recall from (a dyadic version of) (\ref{diagram:eq-cubic-expanded}) and Corollary \ref{counting:cor-cubic} that we have
\begin{equation}\label{analytic:eq-sextic2}\scubic\,[M_1,M_2,M_3,M_{123}=N_1](t,x)=\sum_{n_0,n_1,n_2,n_3}H_{n_0n_1n_2n_3}(t)e^{in_0\cdot x}\Sc\Ic(n_j,\varphi_j:1\leq j\leq 3),\end{equation} where $H_{n_0n_1n_2n_3}(t)$ satisfies
\begin{equation}\label{analytic:eq-sextic3}\sup_t\|H\|_{n_0n_1n_2n_3}\lesssim M_{123}^{\epsilon},\quad \sup_t\|H\|_{n_A\to n_B}\lesssim \max(M_1,M_2,M_3)^{-1/2}
\end{equation} for any partition $(A,B)$ of $\{0,1,2,3\}$ with $A,B\neq\varnothing$. Similarly we have
\begin{equation}\label{analytic:eq-sextic4}\scubic\,[M_4,M_5,M_6,M_{456}=N_2](t,x)=\sum_{n_0',n_4,n_5,n_6}(H')_{n_0'n_4n_5n_6}(t)e^{in_0\cdot x}\Sc\Ic(n_j,\varphi_j:4\leq j\leq 6),\end{equation} where $(H')_{n_0'n_4n_5n_6}(t)$ satisfies
\begin{equation}\label{analytic:eq-sextic5}\sup_t\|h'\|_{n_0'n_4n_5n_6}\lesssim M_{456}^{\epsilon},\quad \sup_t\|H'\|_{n_A\to n_B}\lesssim \max(M_4,M_5,M_6)^{-1/2}
\end{equation} for any partition $(A,B)$ of $\{0,4,5,6\}$ with $A,B\neq\varnothing$ (and $n_0$ replaced by $n_0'$). Therefore
\begin{equation}\label{analytic:eq-sextic6}
\begin{aligned}&P_{\leq N_2^{\sqrt{\nu}}}\Big(P_{N_1}\scubic[\leqN] P_{N_2}\scubic[\leqN]\Big)(t,x)\\&=\sum_\Pc\sum_{\substack{n_0,n_0',n_1,\cdots,n_6\\n_i+n_j=0\,(\forall \{i,j\}\in\Pc)}}e^{i(n_0+n_0')\cdot x}\mathbf{1}\big\{|n_0+n_0'|\leq N_2^{\sqrt{\nu}}\big\}H_{n_0n_1n_2n_3}(t)(H')_{n_0'n_4n_5n_6}(t)\Sc\Ic(n_j,\varphi_j:j\in O),
\end{aligned}
\end{equation} where $\Pc$ is a collection of pairings, i.e. disjoint two-element subsets $\{i,j\}$ of $\{1,\cdots,6\}$, that does not contain any subset in $\{1,2,3\}$ or in $\{4,5,6\}$, and $O$ is the set of indices in $\{1,\cdots,6\}$ not in a subset in $\Pc$. Moreover  the renormalization term $\mathfrak{C}_{\leq N}^{(3,3)}[N_1,N_2]$, upon decomposing in $(M_1,\cdots,M_6)$, exactly corresponds to the cases where $\Pc$ contains three pairings, we know that $\Mc_{3,3}$ can be written in the same form as (\ref{analytic:eq-sextic6}) but with $\Pc$ containing at most two pairings.

(1) Suppose $\Pc=\varnothing$, then we have
\begin{equation}\label{analytic:eq-sextic7}\begin{aligned}
P_{\leq N_2^{\sqrt{\nu}}}\Mc_{3,3}(t,x)=&\sum_{n_0,n_0',n_1,\cdots,n_6}e^{i(n_0+n_0')\cdot x}\mathbf{1}\big\{|n_0+n_0'|\leq N_2^{\sqrt{\nu}}\big\}H_{n_0n_1n_2n_3}(t)\\
&\times
(H')_{n_0'n_4n_5n_6}(t)\Sc\Ic(n_j,\varphi_j:1\leq j\leq 6),
\end{aligned}
\end{equation}
hence by Lemma \ref{counting:lem-merging} we have
\begin{equation}\label{analytic:eq-sextic8}
\begin{aligned}\Eb|P_{\leq N_2^{\sqrt{\nu}}}\Mc_{3,3}(t,x)|^2&\lesssim\Big\|\sum_{n_0,n_0'}e^{i(n_0+n_0')\cdot x}\mathbf{1}\big\{|n_0+n_0'|\leq N_2^{\sqrt{\nu}}\big\}H_{n_0n_1n_2n_3}(t)(H')_{n_0'n_4n_5n_6}(t)\Big\|_{n_1\cdots n_6}^2\\
&\lesssim\|H\|_{n_0n_1n_2n_3}^2\cdot\|e^{i(n_0+n_0')\cdot x}\mathbf{1}\big\{|n_0+n_0'|\leq N_2^{\sqrt{\nu}}\big\}\|_{n_0'\to n_0}^2\cdot\|H'\|_{n_4n_5n_6\to n_0'}^2\\&\lesssim N_2^\epsilon (N_2^{\sqrt{\nu}})^3(\max(M_4,M_5,M_6))^{-1},
\end{aligned}
\end{equation} which proves (\ref{analytic:eq-sextic1}).

(2) Suppose $|\Pc|=1$, say $\Pc=\{\{1,4\}\}$ by symmetry, then we have
\begin{equation}\label{analytic:eq-sextic9}\begin{aligned}P_{\leq N_2^{\sqrt{\nu}}}\Mc_{3,3}(t,x)=&\sum_{n_0,n_0',n_1}e^{i(n_0+n_0')\cdot x}\mathbf{1}\big\{|n_0+n_0'|\leq N_2^{\sqrt{\nu}}\big\}H_{n_0n_1n_2n_3}(t)\\&\times(H')_{n_0',-n_1,n_5n_6}(t)\Sc\Ic(n_j,\varphi_j:j\in\{2,3,5,6\}),\end{aligned}
\end{equation}
hence by Lemma \ref{counting:lem-merging} we have
\begin{equation}\label{analytic:eq-sextic10}
\begin{aligned}\Eb|P_{\leq N_2^{\sqrt{\nu}}}\Mc_{3,3}(t,x)|^2&\lesssim\Big\|\sum_{n_0,n_0',n_1}e^{im\cdot x}\mathbf{1}\big\{m=n_0+n_0'\big\}H_{n_0n_1n_2n_3}(t)(H')_{n_0',-n_1,n_5n_6}(t)\Big\|_{n_2n_3n_5n_6}^2\\
&\lesssim\|H\|_{n_0n_1n_2n_3}^2\cdot\|e^{i(n_0+n_0')\cdot x}\mathbf{1}\big\{|n_0+n_0'|\leq N_2^{\sqrt{\nu}}\big\}\|_{n_0'\to n_0}^2\cdot\|H'\|_{n_5n_6\to n_0'n_4}^2\\&\lesssim N_2^\epsilon (N_2^{\sqrt{\nu}})^3(\max(M_4,M_5,M_6))^{-1},
\end{aligned}
\end{equation} which proves (\ref{analytic:eq-sextic1}).

(3) Suppose $|\Pc|=2$, say $\Pc=\{\{1,4\},\{2,5\}\}$, then we have
\begin{equation}\label{analytic:eq-sextic11}\begin{aligned}P_{\leq N_2^{\sqrt{\nu}}}\Mc_{3,3}(t,x)=&\sum_{n_0,n_0',n_1,n_2}e^{i(n_0+n_0')\cdot x}\mathbf{1}\big\{|n_0+n_0'|\leq N_2^{\sqrt{\nu}}\big\}H_{n_0n_1n_2n_3}(t)\\&\times (H')_{n_0',-n_1,-n_2,n_6}(t)\Sc\Ic(n_3,n_6,\varphi_3,\varphi_6),\end{aligned}
\end{equation}
hence by Lemma \ref{counting:lem-merging} we have
\begin{equation}\label{analytic:eq-sextic12}
\begin{aligned}\Eb|P_{\leq N_2^{\sqrt{\nu}}}\Mc_{3,3}(t,x)|^2&\lesssim\Big\|\sum_{n_0,n_0',n_1,n_2}e^{im\cdot x}\mathbf{1}\big\{m=n_0+n_0'\big\}h_{n_0n_1n_2n_3}(t)(h')_{n_0',-n_1,-n_2,n_6}(t)\Big\|_{n_3n_6}^2\\
&\lesssim\|H\|_{n_0n_1n_2n_3}^2\cdot\|e^{i(n_0+n_0')\cdot x}\mathbf{1}\big\{|n_0+n_0'|\leq N_2^{\sqrt{\nu}}\big\}\|_{n_0'\to n_0}^2\cdot\|H'\|_{n_6\to n_0'n_4n_5}^2\\&\lesssim N_2^\epsilon (N_2^{\sqrt{\nu}})^3(\max(M_4,M_5,M_6))^{-1},
\end{aligned}
\end{equation} which proves (\ref{analytic:eq-sextic1}).
\end{proof}
\begin{proposition}\label{analytic:prop-sextic2} For any dyadic scale $N_1$ and $N_2$ such that $N_1\lesssim N_2$, we have
\begin{equation}\label{analytic:eq-sextic13}\Big\|P_{N_1}\slinear[blue][\leqN]\cdot P_{N_2}\squintic[\leqN]-\mathfrak{C}_{\leq N}^{(1,5)}[N_1,N_2]\Big\|_{L_\omega^p(C_t^0C_x^{-1/10}([-T,T]))}\lesssim p^{3}T^\alpha N_2^{-100\nu}.
\end{equation}
\end{proposition}
\begin{proof} Define the object $\Mc_{1,5}:=P_{N_1}\slinear[blue][\leqN]\cdot P_{N_2}\squintic[\leqN]-\mathfrak{C}_{\leq N}^{(1,5)}[N_1,N_2]$. Start by considering the projections $P_{\leq N_2^{\sqrt{\nu}}}\Mc_{1,5}$ and $P_{\geq N_2^{\sqrt{\nu}}}\Mc_{1,5}$; by Lemma \ref{analytic:lem-linear} and Proposition \ref{analytic:prop-quintic}, like in the proof of Proposition \ref{analytic:prop-sextic2} we can deduce that
\[\Big\|P_{N_1}\slinear[blue][\leqN]\cdot P_{N_2}\squintic[\leqN]\Big\|_{L_\omega^p(C_t^0C_x^0([-T,T]))}\lesssim p^3T^\alpha N_1^{1/2+\epsilon}N_2^{-1/2+\epsilon},\] which implies (\ref{analytic:eq-sextic13}) for $P_{\geq N_2^{\sqrt{\nu}}}\Mc_{1,5}$ as $N_1\lesssim N_2$. Now we only need to consider $P_{\leq N_2^{\sqrt{\nu}}}\Mc_{1,5}$, so we may assume $N_1\sim N_2$.

Recall the decomposition of $\squintic[\leqN]$ in (\ref{analytic:eq-quintic-decomp-pure}); we shall refer to the first two terms in (\ref{analytic:eq-quintic-decomp-pure}) as $(\squintic[\leqN])_0$ and $(\squintic[\leqN])_1$.

(1) First consider the term $(\squintic[\leqN])_0$. As in Lemma \ref{analytic:lem-quintic-pure-nopairing}, we may decompose $P_{N_2}\squintic[\leqN]$ into \[(\squintic[\leqN])_{0}[M_0=N_2,M_1,\cdots M_5,M_{234}].\] By the dyadic version expression \eqref{analytic:dyadic-quintic-tree} of the term $(\squintic[\leqN])_0$, we write
\[(\squintic[\leqN])_{0}[M_*]=\sum_{n_0',n_1,\cdots, n_5}H_{n_0'n_1\cdots n_5}(t)e^{in_0'\cdot x}\Sc\Ic(n_j,\varphi_j:1\leq j\leq 5),\] where $M_*=(M_0=N_2,M_1,\cdots,M_5,M_{234})$ and by Corollary \ref{counting:cor-quintic}, $H$ is a tensor satisfying
\[\sup_t\|H\|_{n_0'\cdots n_5}\lesssim (\max(M_1,\cdots,M_5))^{-1/2+\epsilon},\] \[\sup_t\|H\|_{n_0'n_A\to n_Bn_5}\lesssim (M_0M_5)^{-1/2}\big( M_0^{-\frac{1}{2}} +\max (M_2, M_3, M_4, M_5)^{-\frac{1}{2}}\big)\] for any partition $(A,B)$ of $\{1,2,3,4\}$.

Now we consider two cases depending on whether there is a pairing between $\slinear[blue]$ and $\squintic$. If there is no pairing, then we have
\begin{multline*}P_{\leq N_2^{\sqrt{\nu}}}\Mc_{1,5}(t,x)=\sum_{n_0,n_0',n_1,\cdots,n_5}\mathbf{1}\big\{|n_0+n_0'|\leq N_2^{\sqrt{\nu}}\big\}e^{i(n_0+n_0')\cdot x}\langle n_0\rangle^{-1}\varphi_0(\langle n_0\rangle t)\\\times H_{n_0'n_1\cdots n_5}(t)\Sc\Ic(n_j,\varphi_j:0\leq j\leq 5),
\end{multline*} thus we can estimate
\[
\begin{aligned}\Eb|P_{\leq N_2^{\sqrt{\nu}}}M_{1,5}(t,x)|^2&\lesssim\bigg\|\sum_{n_0'}\mathbf{1}\big\{|n_0+n_0'|\leq N_2^{\sqrt{\nu}}\big\}e^{i(n_0+n_0')\cdot x}\langle n_0\rangle^{-1}\varphi_0(\langle n_0\rangle t)\cdot H_{n_0'n_1\cdots n_5}(t)\bigg\|_{n_0n_1\cdots n_5}^2\\
&\lesssim\big\|\mathbf{1}\big\{|n_0+n_0'|\leq N_2^{\sqrt{\nu}}\big\}e^{i(n_0+n_0')\cdot x}\langle n_0\rangle^{-1}\varphi_0(\langle n_0\rangle t)\big\|_{n_0\to n_0'}\cdot\|H\|_{n_0'n_1\cdots n_5}\\
&\lesssim N_1^{-1}(N_2^{\sqrt{\nu}})^3(\max(M_3,M_4,M_5))^{-1/2},
\end{aligned}\] which proves (\ref{analytic:eq-sextic13}).

The case of one pairing is a bit trickier. Assume $n_0+n_1=0$ (other cases are either similar or easier), then we have
\begin{equation}\label{analytic:5tensorexp}
\begin{aligned}P_{\leq N_2^{\sqrt{\nu}}}\Mc_{1,5}(t,x)=&\sum_{n_0,n_0',n_2,\cdots,n_5}\mathbf{1}\big\{|n_0+n_0'|\leq N_2^{\sqrt{\nu}}\big\}e^{i(n_0+n_0')\cdot x}\langle n_0\rangle^{-1}\varphi_0(\langle n_0\rangle t)\\&\times H_{n_0',-n_0,n_2,\cdots n_5}(t)\Sc\Ic(n_j,\varphi_j:2\leq j\leq 5),\end{aligned}\end{equation} thus we can estimate
\begin{equation}\label{analytic:eq-sextic-extra}
\begin{aligned}\Eb|P_{\leq N_2^{\sqrt{\nu}}}M_{1,5}(t,x)|^2&\lesssim\bigg\|\sum_{n_0'}\mathbf{1}\big\{|n_0+n_0'|\leq N_2^{\sqrt{\nu}}\big\}e^{i(n_0+n_0')\cdot x}\langle n_0\rangle^{-1}\varphi_0(\langle n_0\rangle t)\cdot H_{n_0',-n_0,\cdots n_5}(t)\bigg\|_{n_2\cdots n_5}^2\\
&\lesssim\big\|\mathbf{1}\big\{|n_0+n_0'|\leq N_2^{\sqrt{\nu}}\big\}e^{i(n_0+n_0')\cdot x}\langle n_0\rangle^{-1}\varphi_0(\langle n_0\rangle t)\big\|_{n_0n_0'}\cdot\|H\|_{n_2\cdots n_5\to n_0'n_1}\\
&\lesssim N_1^{-1}(N_2^{\sqrt{\nu}})^{3/2}N_1^{3/2}\cdot(N_2M_5)^{-1/2}\big( N_2^{-1/2} +\max (M_2, M_3, M_4, M_5)^{-1/2}\big).
\end{aligned}\end{equation} This implies (\ref{analytic:eq-sextic13}) if $\max(M_2,\cdots,M_5)\geq N_2^{1/1000}$.

Now, suppose $\max(M_2,\cdots,M_5)\leq N_2^{1/1000}$, then instead of (\ref{analytic:5tensorexp}), we shall examine the full expression of $\Mc_{1,5}$ directly and exploit the cancellation structure as in Lemma \ref{counting:lem-Sine-symmetrization}. By losing at most an $N_2^{1/50}$ factor, we may also fix the values of $n_2,\cdots,n_5$. This allows us to write $\Mc_{1,5}$ in terms of the sine cancellation kernel, namely, up to a linear combination we have
\[\Fc_x\Lc_{1,5}(t,k)=\int_0^t\mathrm{Sine}[M_0,M_1](t-t',k)e^{i\lambda t'}\,\mathrm{d}t'\] with $\lambda$ being a parameter depending on the choice of $(n_2,\cdots,n_5)$ and $|\lambda|\lesssim N_2^{1/10}$. Lemma \ref{counting:lem-Sine-estimate} then immediately implies that
\[\sup_k|\Fc\Mc_{1,5}(t,k)|\lesssim N_2^{-9/10}\] for each $|k|\lesssim N_2^{\sqrt{\nu}}$, which then easily proves (\ref{analytic:eq-sextic13}).

(2) Now consider the term $(\squintic[\leqN])_{1}$. Similar to (1) we can make dyadic decompositions, and now reduce to \[(\squintic[\leqN])_{1}[M_*]=\sum_{n_0',n_3,n_4, n_5}H_{n_0'n_3n_4n_5}(t)e^{in_0'\cdot x}\Sc\Ic(n_j,\varphi_j:3\leq j\leq 5),\] where $M_*=(M_0=N_2,M_1,\cdots,M_5,M_{234})$ and by using \eqref{counting:sine_bound3} and \eqref{counting:sine_bound5} in Corollary \ref{counting:cor-sine-cancel}, $H$ is a tensor satisfying
\[\sup_t\|H\|_{n_0'\cdots n_5}\lesssim (\max(M_3,M_4,M_5))^{-1/2+\epsilon},\quad \sup_t\|H\|_{n_3n_4\to n_0'n_5}\lesssim M_0^{-1/2+\epsilon}(\max(M_3,M_4))^{-1/2+\epsilon}.\] Similar to (1), we can consider the no-pairing case where
\begin{multline*}P_{\leq N_2^{\sqrt{\nu}}}\Mc_{1,5}(t,x)=\sum_{n_0,n_0',n_3,n_4,n_5}\mathbf{1}\big\{|n_0+n_0'|\leq N_2^{\sqrt{\nu}}\big\}e^{i(n_0+n_0')\cdot x}\langle n_0\rangle^{-1}\varphi_0(\langle n_0\rangle t)\\\times H_{n_0'n_3n_4n_5}(t)\Sc\Ic(n_j,\varphi_j:j\in\{0,3,4,5\})\end{multline*} and we have
\[
\begin{aligned}\Eb|P_{\leq N_2^{\sqrt{\nu}}}M_{1,5}(t,x)|^2&\lesssim\bigg\|\sum_{n_0'}\mathbf{1}\big\{|n_0+n_0'|\leq N_2^{\sqrt{\nu}}\big\}e^{i(n_0+n_0')\cdot x}\langle n_0\rangle^{-1}\varphi_0(\langle n_0\rangle t)\cdot H_{n_0'n_3n_4n_5}(t)\bigg\|_{n_0n_3n_4n_5}^2\\
&\lesssim\big\|\mathbf{1}\big\{|n_0+n_0'|\leq N_2^{\sqrt{\nu}}\big\}e^{i(n_0+n_0')\cdot x}\langle n_0\rangle^{-1}\varphi_0(\langle n_0\rangle t)\big\|_{n_0\to n_0'}\cdot\|H\|_{n_0'n_3n_4n_5}\\
&\lesssim N_1^{-1}(N_2^{\sqrt{\nu}})^3(\max(M_3,M_4,M_5))^{-1/2},
\end{aligned}\] which proves (\ref{analytic:eq-sextic13}). 

As for the one-pairing case, we only consider the hardest case $n_0+n_5=0$. Then we have \begin{multline*}P_{\leq N_2^{\sqrt{\nu}}}\Mc_{1,5}(t,x)=\sum_{n_0,n_0',n_3,n_4}\mathbf{1}\big\{|n_0+n_0'|\leq N_2^{\sqrt{\nu}}\big\}e^{i(n_0+n_0')\cdot x}\langle n_0\rangle^{-1}\varphi_0(\langle n_0\rangle t)\\\times H_{n_0',n_3n_4,-n_0}(t)\Sc\Ic(n_3,n_4,\varphi_3,\varphi_4)\end{multline*} and we have
\[
\begin{aligned}\Eb|P_{\leq N_2^{\sqrt{\nu}}}M_{1,5}(t,x)|^2&\lesssim\bigg\|\sum_{n_0,n_0'}\mathbf{1}\big\{|n_0+n_0'|\leq N_2^{\sqrt{\nu}}\big\}e^{i(n_0+n_0')\cdot x}\langle n_0\rangle^{-1}\varphi_0(\langle n_0\rangle t)\cdot H_{n_0',n_3n_4,-n_0}(t)\bigg\|_{n_3n_4}^2\\
&\lesssim\big\|\mathbf{1}\big\{|n_0+n_0'|\leq N_2^{\sqrt{\nu}}\big\}e^{i(n_0+n_0')\cdot x}\langle n_0\rangle^{-1}\varphi_0(\langle n_0\rangle t)\big\|_{n_0 n_0'}\cdot\|H\|_{n_3n_4\to n_0'n_5}\\
&\lesssim N_1^{-1}(N_2^{\sqrt{\nu}})^{3/2}N_1^{3/2}\cdot M_0^{-1/2+\epsilon}(\max(M_3,M_4))^{-1/2+\epsilon}M_1^{-1+\epsilon},
\end{aligned}\] which proves (\ref{analytic:eq-sextic13}) if $\max(M_1,M_3,M_4)\geq N_2^{1/1000}$. If $\max(M_1,M_3,M_4)\leq N_2^{1/1000}$, then we can rewrite $\Mc_{1,5}$ using the sine cancellation kernel in essentially the same way as (1), and using Lemma \ref{counting:lem-Sine-estimate} we can again prove (\ref{analytic:eq-sextic13}).

(3) Finally consider the resistor term. Since the renormalization term $\mathfrak{C}_{\leq N}^{(1,5)}[N_1,N_2]$ exactly equals the contribution where the linear object $\slinear[blue]$ pairs with the resistor $\sresistor$, in studying $\Mc_{1,5}$ we may assume there is no pairing. As such, by Corollary \ref{diagram:corollary-resistor} and Lemma \ref{analytic:lem-resistor} we have
\[P_{\leq N_2^{\sqrt{\nu}}}\Mc_{1,5}(t,x)=\sum_{n_0,n_0'}\mathbf{1}\big\{|n_0+n_0'|\leq N_2^{\sqrt{\nu}}\big\}e^{i(n_0+n_0')\cdot x}\langle n_0\rangle^{-1}\varphi_0(\langle n_0\rangle t)\cdot F_{\leq N}^{\varphi_0'}(t,n_0')\Sc\Ic(n_0,n_0',\varphi_0,\varphi_0')\] where $F_{\leq N}^{\varphi_0'}$ is a function satisfying
\[\sup_t\|F_{\leq N}^{\varphi_0'}\|_{\ell_{n_0'}^2}\lesssim\|F_{\leq N}^{\varphi_0'}\|_{X^{0,b}}\lesssim N_2^{-1/2+\epsilon}.\] This easily implies that
\[
\begin{aligned}\Eb|P_{\leq N_2^{\sqrt{\nu}}}M_{1,5}(t,x)|^2&\lesssim\big\|\mathbf{1}\big\{|n_0+n_0'|\leq N_2^{\sqrt{\nu}}\big\}e^{i(n_0+n_0')\cdot x}\langle n_0\rangle^{-1}\varphi_0(\langle n_0\rangle t)\cdot R_{n_0'}(t)\big\|_{n_0n_0'}^2\\
&\lesssim N_2^{-1+2\epsilon}\cdot\sup_{n_0'}\sum_{n_0:|n_0+n_0'|\leq N_2^{\sqrt{\nu}}}\langle n_0\rangle^{-2}\lesssim N_2^{-1+2\epsilon}N_0^{-2}(N_2^{\sqrt{\nu}})^3
\end{aligned}\] which proves (\ref{analytic:eq-sextic13}).
\end{proof}

We finish this subsection with a crude estimate on the cancellation constants $\mathfrak{C}^{(1,5)}_{\leq N}$ and $\mathfrak{C}^{(3,3)}_{\leq N}$, which will be useful in the analysis of  the para-controlled terms $\XXone$ and $\XXtwo$ (see Section \ref{section:para}). 

\begin{lemma}[\protect{Crude estimates of $\mathfrak{C}^{(1,5)}_{\leq N}$ and $\mathfrak{C}^{(3,3)}_{\leq N}$}]\label{analytic:lem-C15-C33}
Let $T\geq 1$, let $N_1,N_2,N$ be frequency-scales and let $\mathfrak{C}^{(1,5)}_{\leq N}[N_1,N_2]$ and $\mathfrak{C}^{(3,3)}_{\leq N}[N_1,N_2]$ be as in Definition \ref{ansatz:def-dyadic-modified-product}. Then, it holds that 
\begin{equation}\label{analytic:eq-C15-C33-a}
\big|\mathfrak{C}^{(1,5)}_{\leq N}[N_1,N_2](t)\big| 
+\big|\mathfrak{C}^{(3,3)}_{\leq N}[N_1,N_2](t)\big| 
\lesssim \max(N_1,N_2)^{2\epsilon} T^\alpha. 
\end{equation}
for all $t\in [-T,T]$. Furthermore, it holds that 
\begin{equation}\label{analytic:eq-C15-C33-b}
\big\| \chi^2(t/T) \mathfrak{C}^{(1,5)}_{\leq N}[N_1,N_2](t)\big\|_{H_t^b}
+\big\| \chi^2(t/T) \mathfrak{C}^{(3,3)}_{\leq N}[N_1,N_2](t)\big\|_{H_t^b}
\lesssim \max(N_1,N_2)^{2} T^\alpha. 
\end{equation}
\end{lemma}

The first estimate \eqref{analytic:eq-C15-C33-a} will only be applied when $N_1$ and $N_2$ are not the maximal frequency-scales and the second estimate \eqref{analytic:eq-C15-C33-b} will only be applied when $N_1$ and $N_2$ are extremely small compared to the maximal frequency-scale. As a result, the losses in $\max(N_1,N_2)$ in \eqref{analytic:eq-C15-C33-a} and \eqref{analytic:eq-C15-C33-b} pose no problems. 

\begin{proof}
We only prove the estimates for $\mathfrak{C}^{(1,5)}_{\leq N}[N_1,N_2]$, since the estimates for $\mathfrak{C}^{(3,3)}_{\leq N}[N_1,N_2]$ are similar. From the definition of $\mathfrak{C}^{(1,5)}_{\leq N}[N_1,N_2]$ (see also Subsection \ref{section:quintic-diagram}), it follows that $\mathfrak{C}^{(1,5)}_{\leq N}[N_1,N_2]\equiv 0$ unless $N_1=N_2$, which will now be assumed for the rest of this proof. Using translation-invariance, Lemma \ref{analytic:lem-linear}, and Proposition \ref{analytic:prop-quintic}, it follows that
\begin{align*}
\big| \mathfrak{C}^{(1,5)}_{\leq N}[N_1,N_2](t)\big| 
=& \bigg|\E \bigg[ \int_{\T^3} \dx \, \slinear[blue][N_1](t,x) \, 
P_{N_2} \squintic[\leqN] (t,x) \bigg] \bigg| \\
\lesssim& \, \E \Big[ \Big\| \, \slinear[blue][N_1] (t,x) \Big\|_{L_x^2}^2 \Big]^{1/2}
\E \Big[ \Big\| \, P_{N_2} \squintic[\leqN] (t,x) \Big\|_{L_x^2}^2 \Big]^{1/2} \\
\lesssim& \, N_1^{1/2+\epsilon} N_2^{-1/2+\epsilon} T^\alpha \\
\lesssim&\, \max(N_1,N_2)^{2\epsilon} T^\alpha. 
\end{align*}
This proves the first estimate \eqref{analytic:eq-C15-C33-a}. Similarly, we obtain from the algebra property of $H_t^b$ that 
\begin{align}
&\big\| \chi^2(t/T) \mathfrak{C}^{(1,5)}_{\leq N}[N_1,N_2](t)\big\|_{H_t^b} \notag \\
\leq& \, \E \bigg[ \int_{\T^3} \dx \, \Big\| \chi^2(t/T) \slinear[blue][N_1](t,x) \, 
P_{N_2} \squintic[\leqN] (t,x) \Big\|_{H_t^b} \bigg]  \notag \\
\lesssim&\, 
\Big( \E \int_{\T^3} \dx  \, \Big\| \chi(t/T) \slinear[blue][N_1](t,x) \Big\|_{H_t^b}^2 \Big)^{1/2}
\Big( \E \int_{\T^3} \dx  \, \Big\| \chi(t/T) P_{N_2} \squintic[\leqN](t,x) \Big\|_{H_t^b}^2 \Big)^{1/2}. \label{analytic:eq-C15-C33-p1}
\end{align}
Since $\langle \lambda \rangle^b \lesssim \langle n \rangle^b \langle \lambda \mp |n| \rangle^b$ for all $\lambda \in \R$ and $n \in \Z^3$, we obtain from  Lemma \ref{analytic:lem-linear} and Proposition \ref{analytic:prop-quintic} that
\begin{align*}
\eqref{analytic:eq-C15-C33-p1} 
&\lesssim \, \Big( \E  \Big\| \chi(t/T) \slinear[blue][N_1](t,x) \Big\|_{X^{b,b}}^2  \Big)^{1/2}
\Big( \E  \Big\| \chi(t/T) P_{N_2} \squintic[\leqN](t,x) \Big\|_{X^{b,b}}^2  \Big)^{1/2} \\
&\lesssim \, N_1^{b+1/2+\epsilon} N_2^{b-1/2+\epsilon} T^\alpha \\
&\lesssim \, N_1 N_2 T^\alpha. \qedhere
\end{align*}

\end{proof}

\section{Bilinear random operator}\label{section:bilinear}

In this section, we analyze (frequency-localized versions) of the bilinear operator
\begin{equation}\label{bilinear:eq-full-bilinear}
(w_2,w_3) \in X^{1/2,b} \times X^{1/2,b} \mapsto 
P_{\leq N}\Big[ \,\slinear[blue][\leqN]  P_{\leq N} w_2 P_{\leq N}w_3 \Big].
\end{equation}
In the method of random tensors \cite{DNY20}, operators such as \eqref{bilinear:eq-full-bilinear} are viewed as linear operators in the (frequency-space) tensor product $\widehat{w}_2 \otimes \widehat{w}_3$. As will be explained in Remark \ref{bilinear:rem-abstract-linear-bound} and Remark \ref{bilinear:rem-linear-bound} below, however, this linear treatment does not yield our desired estimates. Instead, the analysis of \eqref{bilinear:eq-full-bilinear} requires a bilinear treatment, which extends the ideas in \cite{DNY20}. 
We start this section with an abstract bilinear estimate, which will then be applied to \eqref{bilinear:eq-full-bilinear}. 

\begin{lemma}[Abstract bilinear estimate]\label{bilinear:lem-abstract}
Let $D\geq 1$, let $h=h_{abcd}\colon (\Z^D)^4 \rightarrow \mathbb{C}$ be a tensor supported on a finite set, and let $(g_b)_{b\in \Z^D}$ be a sequence of independent standard complex Gaussians. Define the bilinear random operator $B\colon \ell^2 \times \ell^2 \rightarrow \ell^2$ by 
\begin{equation*}
B(v,w)_a := \sum_{b,c,d\in \Z^D} h_{abcd} g_b v_c w_d. 
\end{equation*}
Then, we have for all $\epsilon>0$ and $p\geq 2$ that 
\begin{equation}\label{bilinear:eq-abstract}
\begin{aligned}
&\mathbb{E} \Big[ \big\| B \big\|_{\ell^2 \times \ell^2 \rightarrow \ell^2}^p \Big]^{1/p} \\ \lesssim_\epsilon& \big( \# \supp h \big)^\epsilon\max 
\Big( \| h\|_{ad\rightarrow bc} \| h \|_{ac\rightarrow bd},
\| h \|_{ad\rightarrow bc} \| h \|_{abc\rightarrow d}, 
\| h \|_{ac\rightarrow bd} \| h \|_{abd \rightarrow c},
\| h \|_{ab\rightarrow cd}^2 \Big)^{1/2}  p^{1/2} .
\end{aligned}
\end{equation}
\end{lemma}

\begin{remark}\label{bilinear:rem-abstract-linear-bound}
By treating $B$ as a linear operator acting on the tensor product $v\otimes w$, we could have estimated
\begin{equation*}
\| B \|_{\ell^2 \times \ell^2\rightarrow \ell^2}
\leq \Big\| \sum_{ b\in \Z^D} h_{abcd} g_b \Big\|_{a \rightarrow cd}. 
\end{equation*}
From the moment method (Proposition \ref{counting:prop-moment}), it follows that 
\begin{equation}\label{bilinear:eq-comparison}
\mathbb{E} \Big[ \Big\| \sum_{ b\in \Z^D} h_{abcd} g_b \Big\|_{a \rightarrow cd}^p \Big]^{1/p}
\lesssim \big( \# \supp h \big)^\epsilon \max \big( \| h \|_{a\rightarrow bcd}, \| h \|_{ab\rightarrow cd} \big). 
\end{equation}
While the right-hand side of \eqref{bilinear:eq-comparison} contains a full power of the $(3,1)$-tensor norm $\| h\|_{a\rightarrow bcd}$, the right-hand side of  \eqref{bilinear:eq-abstract} only contains  square-roots of $(3,1)$-tensor norms. In our application (see Remark \ref{bilinear:rem-linear-bound} below), the bilinear estimate therefore leads to an improvement over the linear tensor estimate. 
\end{remark}

The main ingredient in the proof of Lemma \ref{bilinear:lem-abstract} is a $TT^\ast$-argument. While a $TT^\ast$-argument is also at the heart of the linear tensor estimates in \cite{DNY20}, the bilinear structure of $B$ affords us greater flexibility in the choice of tensor norms. 

\begin{proof} 
For any $v,w \in \ell^2$, it holds that 
\begin{align*}
\big\| B(v,w) \big\|_{\ell^2}^2 
&= \sum_a \Big| \sum_{b,c,d} h_{abcd} g_b v_c w_d\Big|^2 \\
&= \sum_{c,c^\prime,d,d^\prime} \Big( \sum_{a,b,b^\prime} h_{abcd} \overline{h_{ab^\prime c^\prime d^\prime}} g_b \overline{g_{b^\prime}} \Big) v_c \overline{v_{c^\prime}} w_d \overline{w_{d^\prime}} \\
&=\sum_{c,c^\prime,d,d^\prime} \mathcal{B}_{cc^\prime d d^\prime} v_c \overline{v_{c^\prime}} w_d \overline{w_{d^\prime}},
\end{align*}
where we defined the tensor
\begin{equation*}
 \mathcal{B}_{cc^\prime d d^\prime} := 
 \sum_{a,b,b^\prime} h_{abcd} \overline{h_{ab^\prime c^\prime d^\prime}} g_b \overline{g_{b^\prime}}. 
\end{equation*}
As a result, we obtain that 
\begin{equation}\label{bilinear:eq-abstract-p1}
\| B \|_{\ell^2 \times \ell^2 \rightarrow \ell^2}^2 
= \sup_{\| v \|_{\ell^2} \leq 1} \sup_{\| w\|_{\ell^2} \leq 1} \| B(v,w)\|_{\ell^2}^2 
= \sup_{\| v \|_{\ell^2} \leq 1} \sup_{\| w\|_{\ell^2} \leq 1} \sum_{c,c^\prime,d,d^\prime} \mathcal{B}_{cc^\prime d d^\prime} v_c \overline{v_{c^\prime}} w_d \overline{w_{d^\prime}}. 
\end{equation}
In order to estimate the right-hand side of \eqref{bilinear:eq-abstract-p1}, we decompose the random tensor  $\mathcal{B}=\mathcal{B}_{cc^\prime d d^\prime}$ into its non-resonant and resonant components. More precisely, we decompose 
\begin{align*}
\mathcal{B}_{cc^\prime d d^\prime}
&= \sum_{a,b,b^\prime}  h_{abcd} \overline{h_{ab^\prime c^\prime d^\prime}}\big( g_b \overline{g_{b^\prime}} - \delta_{b=b^\prime} \big) + 
\sum_{a,b,b^\prime}  h_{abcd} \overline{h_{ab^\prime c^\prime d^\prime}} \delta_{b=b^\prime} \\
&=: \mathcal{B}_{cc^\prime d d^\prime}^{(2)} + \mathcal{B}_{cc^\prime d d^\prime}^{(0)}.
\end{align*}
Thus, the components $\mathcal{B}^{(2)}$ and $\mathcal{B}^{(0)}$ correspond to the second and zeroth-order Gaussian chaos in $\mathcal{B}$, respectively. Due to \eqref{bilinear:eq-abstract-p1} and our decomposition, it suffices to prove the two estimates
\begin{equation}\label{bilinear:eq-abstract-p2}
\begin{aligned}
&\mathbb{E} \Big[  \sup_{\| v\|_{\ell^2}\leq 1} \sup_{\| w\|_{\ell^2}\leq 1} \Big| 
\sum_{c,c^\prime,d,d^\prime} \mathcal{B}_{cc^\prime d d^\prime}^{(2)} v_c \overline{v_{c^\prime}} w_d \overline{w_{d^\prime}} \Big|^p \Big]^{1/p} \\
\lesssim&\big( \# \supp h\big)^\epsilon \max 
\Big( \| h\|_{ad\rightarrow bc} \| h \|_{ac\rightarrow bd},
\| h \|_{ad\rightarrow bc} \| h \|_{abc\rightarrow d}, 
\| h \|_{ac\rightarrow bd} \| h \|_{abd \rightarrow c} \Big)  p. 
\end{aligned}
\end{equation}
and 
\begin{equation}\label{bilinear:eq-abstract-p3}
\sup_{\| v\|_{\ell^2}\leq 1} \sup_{\| w\|_{\ell^2}\leq 1} \Big| 
\sum_{c,c^\prime,d,d^\prime} \mathcal{B}_{cc^\prime d d^\prime}^{(0)} v_c \overline{v_{c^\prime}} w_d \overline{w_{d^\prime}} \Big| \leq \| h \|_{ab\rightarrow cd}^2. 
\end{equation}
\emph{Estimate of the non-resonant part $\mathcal{B}^{(2)}$:}
Let $v,w\in \ell^2$ satisfy $\| v\|_{\ell^2}, \|w\|_{\ell^2}\leq 1$. Using  linear tensor norms, it holds that 
\begin{equation}\label{bilinear:eq-abstract-q1}
\Big| \sum_{c,c^\prime,d,d^\prime} \mathcal{B}_{cc^\prime d d^\prime}^{(2)} v_c \overline{v_{c^\prime}} w_d \overline{w_{d^\prime}} \Big| 
\leq \| \mathcal{B}_{cc^\prime d d^\prime}^{(2)} \|_{cc^\prime \rightarrow dd^\prime} \| v_c \overline{v_{c^\prime}} \|_{cc^\prime} \| w_d \overline{w_{d^\prime}} \|_{dd^\prime} \leq  \| \mathcal{B}_{cc^\prime d d^\prime}^{(2)} \|_{cc^\prime \rightarrow dd^\prime}. 
\end{equation}
Using the definition of $\mathcal{B}^{(2)}$ and the moment method (Proposition \ref{counting:prop-moment}), it follows that 
\begin{equation}\label{bilinear:eq-abstract-p4}
\begin{aligned}
\mathbb{E} \big[ \| \mathcal{B}_{cc^\prime d d^\prime}^{(2)} \|_{cc^\prime \rightarrow dd^\prime}^p \big]^{1/p}
\lesssim& \, (\# \supp h )^{\epsilon} \max\Big(
\big\| \sum_a h_{abcd}\overline{h_{ab^\prime c^\prime d^\prime}} \big\|_{bb^\prime cc^\prime \rightarrow dd^\prime}, 
\big\| \sum_a h_{abcd}\overline{h_{ab^\prime c^\prime d^\prime}} \big\|_{b cc^\prime \rightarrow b^\prime dd^\prime}, \\
&\big\| \sum_a h_{abcd}\overline{h_{ab^\prime c^\prime d^\prime}} \big\|_{b^\prime cc^\prime \rightarrow b dd^\prime},
\big\| \sum_a h_{abcd}\overline{h_{ab^\prime c^\prime d^\prime}} \big\|_{ cc^\prime \rightarrow bb^\prime dd^\prime}\Big) \cdot  p . 
\end{aligned}
\end{equation}
The four arguments of the maximum in \eqref{bilinear:eq-abstract-p4} are now controlled separately using the merging estimate (Lemma \ref{counting:lem-merging}). The first argument in  \eqref{bilinear:eq-abstract-p4} is estimated by 
\begin{equation}\label{bilinear:eq-abstract-p5}
\big\| \sum_a h_{abcd}\overline{h_{ab^\prime c^\prime d^\prime}} \big\|_{bb^\prime cc^\prime \rightarrow dd^\prime} 
\leq \| h_{abcd}\|_{abc\rightarrow d}
\| \overline{h_{ab^\prime c^\prime d^\prime}} \|_{b^\prime c^\prime \rightarrow ad^\prime} = \| h \|_{ad\rightarrow bc} \| h\|_{abc\rightarrow d}.
\end{equation}
The second argument in \eqref{bilinear:eq-abstract-p4} is estimated by 
\begin{equation}\label{bilinear:eq-abstract-p6}
\big\| \sum_a h_{abcd}\overline{h_{ab^\prime c^\prime d^\prime}} \big\|_{b cc^\prime \rightarrow b^\prime dd^\prime} 
\leq \| h_{abcd}\|_{bc\rightarrow ad} \| \overline{h_{ab^\prime c^\prime d^\prime}} \|_{ac^\prime \rightarrow b^\prime d^\prime} = \| h \|_{ad\rightarrow bc} \| h\|_{ac\rightarrow bd}. 
\end{equation}
Due to the symmetry in $b$ and $b^\prime$, the third argument in \eqref{bilinear:eq-abstract-p4} can be estimated similarly. Finally, the fourth argument in \eqref{bilinear:eq-abstract-p4} is estimated by 
\begin{equation}\label{bilinear:eq-abstract-p7}
    \big\| \sum_a h_{abcd}\overline{h_{ab^\prime c^\prime d^\prime}} \big\|_{ cc^\prime \rightarrow bb^\prime dd^\prime} 
    \leq \| h_{abcd}\|_{ac\rightarrow bd} \| \overline{h_{ab^\prime c^\prime d^\prime}}\|_{c^\prime\rightarrow ab^\prime d^\prime} 
    \leq \| h \|_{ac\rightarrow bd} \| h \|_{abd\rightarrow c}. 
\end{equation}
Since all three right-hand sides of \eqref{bilinear:eq-abstract-p5}, \eqref{bilinear:eq-abstract-p6}, and \eqref{bilinear:eq-abstract-p7} are acceptable, this completes the proof of \eqref{bilinear:eq-abstract-p2}. \\

\emph{Estimate of the resonant part $\mathcal{B}^{(0)}:$} 
Let $v,w\in \ell^2$ satisfy $\| v\|_{\ell^2}, \|w\|_{\ell^2}\leq 1$. Using  linear tensor norms, it holds that 
\begin{equation}\label{bilinear:eq-abstract-q2}
\Big| \sum_{c,c^\prime,d,d^\prime} \mathcal{B}_{cc^\prime d d^\prime}^{(0)} v_c \overline{v_{c^\prime}} w_d \overline{w_{d^\prime}} \Big| 
\leq \| \mathcal{B}_{cc^\prime d d^\prime}^{(0)} \|_{cd\rightarrow c^\prime d^\prime} 
\| v_c w_d \|_{cd} \| \overline{v_{c^\prime}} \overline{w_{d^\prime}} \|_{c^\prime d^\prime} \leq \| \mathcal{B}_{cc^\prime d d^\prime}^{(0)} \|_{cd\rightarrow c^\prime d^\prime} . 
\end{equation}
We emphasize that the linear tensor norms of $\mathcal{B}^{(2)}$ and $\mathcal{B}^{(0)}$ used in 
\eqref{bilinear:eq-abstract-q1} and \eqref{bilinear:eq-abstract-q2} are different, which is essential for the proof and only possible in the bilinear analysis. Using the merging estimates (Lemma \ref{counting:lem-merging}), it follows that 
\begin{align*}
     \| \mathcal{B}_{cc^\prime d d^\prime}^{(0)} \|_{cd\rightarrow c^\prime d^\prime} 
     = \big\| \sum_{a,b} h_{abcd} \overline{h_{abc^\prime d^\prime}} \big\|_{cd\rightarrow c^\prime d^\prime} 
     \leq \| h_{abcd} \|_{cd\rightarrow ab} \| \overline{h_{abc^\prime d^\prime}} \|_{ab\rightarrow c^\prime d^\prime} 
     = \| h \|_{ab\rightarrow cd}^2. 
\end{align*}
This completes the proof of \eqref{bilinear:eq-abstract-p3}. 
\end{proof}

We now apply the abstract bilinear estimate (Lemma \ref{bilinear:lem-abstract}) to our setting. 

\begin{proposition}[Bilinear random operator]\label{bilinear:lem-bilinear}
Let $p\geq 2$, let $N_0,N_1,N_2,N_3$ be frequency-scales satisfying $N_2,N_3 \ll N_1$, and let $T\geq 1$. Then, we have the bilinear estimate
\begin{equation}\label{bilinear:eq-bilinear}
\begin{aligned}
&\mathbb{E} \Big[ \sup_{\Jc} \Big\| (w_2,w_3) \mapsto P_{N_0}  \Big[ \slinear[blue][N_1] P_{N_2}w_2 P_{N_3}w_3 \Big]\Big\|_{X^{1/2,b}(\Jc)\times X^{1/2,b}(\Jc)\rightarrow X^{-1/2,b_+-1}(\Jc)}^p \Big]^{1/p} \\
\lesssim&\, p^{1/2} T^{\alpha} N_{\textup{max}}^\epsilon \big( N_1^{-1/4}+ \max(N_2,N_3)^{-1/3}\big), 
\end{aligned}
\end{equation}
where the supremum is taken over all closed intervals $0 \in \Jc \subseteq [-T,T]$. 
\end{proposition}

\begin{proof}

We utilize the reduction arguments in Subsection \ref{prep:remark-reduction}. Then, it suffices to prove the random tensor estimate
\begin{equation}\label{bilinear:eq-bilinear-p1}
\begin{aligned}
&N_1^{-3/2} N_2^{-1/2} N_3^{-1/2} 
\mathbb{E}\Big[ \Big\| (\widehat{w}_2,\widehat{w}_3) \mapsto \sum_{n_1,n_2,n_3} h_{n_0 n_1 n_2 n_3} g_{n_1} \widehat{w}_2(n_2) \widehat{w}_3(n_3) \Big\|_{\ell^2 \times \ell^2 \rightarrow \ell^2}^p  \Big]^{1/p}\\
\lesssim&\, N_{\textup{max}}^\epsilon \big( N_1^{-1/4}+ \max(N_2,N_3)^{-1/3}\big) p^{1/2},
\end{aligned}
\end{equation}
where $h$ is the base tensor from Section \ref{section:counting}, i.e., 
\begin{align}\label{bilinear:eq-bilinear-p2}
h_{n_0 n_1 n_2 n_3} := \Big(\prod_{j=0}^3 1_{N_j}(n_j)\Big) \mathbf{1}\{ n_0=n_{123} \} 
\mathbf{1}\{ | \Omega-m|\leq 1 \}, \quad 
\Omega:= \sum_{j=0}^3 (\pm_j) \langle n_ j \rangle, \quad \text{and} \quad m \in \Z. 
\end{align}

In the case $N_2 \sim N_3$, we further insert a dyadic localization to $|n_{23}|\sim N_{23}$. Using a standard box-localization argument in $n_2$ and $n_3$ (see e.g. \cite{Tao01}), we can further the insert the indicator function 
\begin{equation*}
 \prod_{j=2}^{3} 1_{Q_j}(n_j)
\end{equation*}
into our definition of $h$ in \eqref{bilinear:eq-bilinear-p2}, where $Q_2$ and $Q_3$ are boxes of sidelength $\sim N_{23}$ and at a distance $\lesssim N_2 \sim N_3$ from the origin. In total, we therefore replace\footnote{In the case $N_2 \not \sim N_3$, one can still insert the additional dyadic and box localizations, but they provide no new information.} the tensor $h$ by 
\begin{equation}\label{bilinear:eq-bilinear-p3}
\begin{aligned}
    h_{n_0 n_1 n_2 n_3} :=&
    \Big( \mathbf{1} \{ N_2 \not \sim N_3 \} + \mathbf{1}\{ N_2 \sim N_3 \} 1_{N_{23}}(n_{23}) \big( \prod_{j=2}^{3} 1_{Q_j}(n_j) \big) \Big) \\
    &\times\Big(\prod_{j=0}^3 1_{N_j}(n_j)\Big) \mathbf{1}\{ n_0=n_{123} \} 
\mathbf{1}\{ | \Omega-m|\leq 1 \}. 
\end{aligned}
\end{equation}
Due to the abstract bilinear estimate (Lemma \ref{bilinear:lem-abstract}), the random bilinear estimate \eqref{bilinear:eq-bilinear-p1} can be reduced to the deterministic tensor estimate 
\begin{equation}\label{bilinear:eq-bilinear-p4}
\begin{aligned}
   &N_1^{-3/2} N_2^{-1/2} N_3^{-1/2} \max\Big(
    \| h\|_{n_0 n_3 \rightarrow n_1 n_2} \| h\|_{n_0 n_2 \rightarrow n_1 n_3}, 
    \| h \|_{n_0 n_3 \rightarrow n_1 n_2} \| h \|_{n_0 n_1 n_2 \rightarrow n_3}, \\
    &\| h\|_{n_0 n_2 \rightarrow n_1 n_3} \| h \|_{n_0 n_1 n_3 \rightarrow n_2}, 
    \| h \|_{n_0 n_1 \rightarrow n_2 n_3}^2 \Big)^{1/2} \lesssim N_1^{-1/4} + \max(N_2,N_3)^{-1/3}. 
 \end{aligned}
\end{equation}
We now treat the four arguments in \eqref{bilinear:eq-bilinear-p4} separately. The most interesting contribution comes from the fourth argument, which will be treated last and is the only argument requiring the box-localization. \\

\emph{The first argument in \eqref{bilinear:eq-bilinear-p4}:} Using the assumption $N_2,N_3 \ll N_1$ and the base tensor estimate (Lemma \ref{counting:basetensor_est}), it holds that 
\begin{align*}
    &N_1^{-3/2} N_2^{-1/2} N_3^{-1/2} 
    \| h\|_{n_0 n_3 \rightarrow n_1 n_2}^{1/2} \| h\|_{n_0 n_2 \rightarrow n_1 n_3}^{1/2} \\
    \lesssim& \, N_1^{-3/2} N_2^{-1/2} N_3^{-1/2}  \big( N_2 N_3 \big)^{1/2} \big( N_2 N_3 \big)^{1/2} \\
    \lesssim& \, N_1^{-3/2} N_2^{1/2} N_3^{1/2} \lesssim N_1^{-1/2}. 
\end{align*}
This contribution is (better than) acceptable. \\

\emph{The second argument in \eqref{bilinear:eq-bilinear-p4}:} 
Using the assumption $N_2,N_3 \ll N_1$ and the base tensor estimate (Lemma \ref{counting:basetensor_est}), it holds that 
\begin{align*}
&N_1^{-3/2} N_2^{-1/2} N_3^{-1/2}
\| h \|_{n_0 n_3 \rightarrow n_1 n_2}^{1/2} \| h \|_{n_0 n_1 n_2 \rightarrow n_3}^{1/2}\\
\lesssim& \, N_1^{-3/2} N_2^{-1/2} N_3^{-1/2} \big( N_2 N_3 \big)^{1/2} \big( N_1^{3/2} N_2 \big)^{1/2} \\
\lesssim& \, N_1^{-3/4} N_2^{1/2} \lesssim N_1^{-1/4}. 
\end{align*}
This contribution is responsible for the $N_1^{-1/4}$-term in our final estimate \eqref{bilinear:eq-bilinear-p4}. \\

\emph{The third argument in \eqref{bilinear:eq-bilinear-p4}:} Due to symmetry in the $n_2$ and $n_3$-variables, the third argument in \eqref{bilinear:eq-bilinear-p4} can be treated exactly as the second argument in \eqref{bilinear:eq-bilinear-p4}. \\

\emph{The fourth argument in \eqref{bilinear:eq-bilinear-p4}:} For expository purposes, we distinguish the cases $N_2 \not \sim N_3$ and $N_2 \sim N_3$. In the easier case $N_2 \not\sim N_3$, it follows from our assumption $N_2,N_3 \ll N_1$ and the base tensor estimate (Lemma \ref{counting:basetensor_est}) that 
\begin{align*}
&N_1^{-3/2} N_2^{-1/2} N_3^{-1/2} \| h\|_{n_0 n_1 \rightarrow n_2 n_3} \\
\lesssim& N_1^{-3/2} N_2^{-1/2} N_3^{-1/2} \max(N_2,N_3)^{-1/2} N_1^{3/2} \min(N_2,N_3)\\
\lesssim& \max(N_2,N_3)^{-1/2}. 
\end{align*}
This contribution is better than the $\max(N_2,N_3)^{-1/3}$-term in our final estimate. Finally, we treat the case $N_2 \sim N_3$. Using Schur's test, the box-localization in \eqref{bilinear:eq-bilinear-p3}, and the box-counting lemma (Lemma \ref{counting:lem-annuli-box}), it follows that 
\begin{equation*}
\| h \|_{n_0 n_1\rightarrow n_2 n_3} \lesssim N_{23}^{-1/2} N_1^{3/2} \min\big(N_2,N_{23}^{3/2}\big) \lesssim N_1^{3/2} N_2^{2/3}. 
\end{equation*}
As a result, it follows that 
\begin{equation*}
N_1^{-3/2} N_2^{-1/2} N_3^{-1/2} \| h\|_{n_0 n_1 \rightarrow n_2 n_3}
\lesssim N_2^{-1/2-1/2+2/3} = N_2^{-1/3}. 
\end{equation*}
This contribution is responsible for the second summand $\max(N_2,N_3)^{-1/3}$ in our final estimate \eqref{bilinear:eq-bilinear-p4}. 
\end{proof}

\begin{remark}\label{bilinear:rem-linear-bound}
As mentioned above, the proof of Lemma \ref{bilinear:lem-bilinear} crucially relies on the abstract bilinear estimate (Lemma \ref{bilinear:lem-abstract}). If we instead try to estimate the left-hand side of \eqref{bilinear:eq-bilinear-p1} directly using the moment method (Proposition \ref{counting:prop-moment}), which only utilizes linear tensor analysis, it can only be bounded by 
\begin{equation*}
N_1^{-3/2} N_2^{-1/2} N_3^{-1/2} \max\big( \| h \|_{n_0 \rightarrow n_1 n_2 n_3}, \| h \|_{n_0 n_1\rightarrow n_2 n_3} \big). 
\end{equation*}
Unfortunately, the base tensor estimate (Lemma \ref{counting:basetensor_est}) only yields
\begin{align*}
N_1^{-3/2} N_2^{-1/2} N_3^{-1/2}\| h \|_{n_0 \rightarrow n_1 n_2 n_3} 
&\lesssim N_1^{-3/2} N_2^{-1/2} N_3^{-1/2} \max(N_2,N_3)^{3/2} \min(N_2,N_3) \\
&\lesssim N_1^{-3/2} \max(N_2,N_3) \min(N_2,N_3)^{1/2}. 
\end{align*}
Since the total power of all frequency-scales adds up to zero, this estimate is insufficient from a dimensional standpoint. To avoid confusion, we emphasize that this issue cannot be addressed through box localization, which yields no improvements when $N_{23}\sim N_2 \sim N_3$. 
\end{remark}

Using similar arguments as in the proof of Lemma \ref{bilinear:lem-bilinear}, we also prove a variant of our bilinear estimate, which concerns a different (and extremely specific) frequency-regime. This variant is only used in the proof of Corollary \ref{para:cor-blueXj} and we encourage the reader to skip it on first reading. 

\begin{lemma}[Variant of the bilinear estimate]\label{bilinear:lem-variant}
Let $T\geq 1$, let $p\geq 2$, let $N_0,N_1,N_2,N_3,N_{12}$ be frequency scales satisfying 
\begin{equation}\label{bilinear:eq-variant-frequency-regime}
N_{\textup{max}}\sim N_2, \quad N_{12} \gtrsim N_2^{1-1/100}, \quad \text{and} \quad N_0 \lesssim N_2^{1/100}. 
\end{equation}
Then, it holds that 
\begin{equation}\label{bilinear:eq-variant}
\begin{aligned}
&\mathbb{E} \Big[ \sup_{\Jc} \Big\| (w_2,w_3) \mapsto P_{N_0}  \Big[ P_{N_{12}} \Big( \slinear[blue][N_1] P_{N_2}w_2 \Big) P_{N_3}w_3 \Big]\Big\|_{X^{1/2,b}(\Jc)\times X^{0,b}(\Jc)\rightarrow X^{-1/2,b_+-1}(\Jc)}^p \Big]^{1/p} \\
\lesssim&\, p^{1/2} T^{\alpha} N_2^{-1/8}, 
\end{aligned}
\end{equation}
where the supremum is taken over all closed intervals $0 \in \Jc \subseteq [-T,T]$. 
\end{lemma}
We note that, in contrast to Lemma \ref{bilinear:lem-bilinear}, $w_3$ is placed in $X^{0,b}$ instead of $X^{1/2,b}$. 
The precise exponents $1/100$ and $1/8$ in \eqref{bilinear:eq-variant-frequency-regime} and \eqref{bilinear:eq-variant} are not important and we only used absolute numbers to avoid introducing a further small parameter. 

\begin{proof}[Proof of Lemma \ref{bilinear:lem-variant}:] 
By arguing as in the proof of Lemma \ref{bilinear:lem-bilinear}, it suffices to prove that 
\begin{equation}\label{bilinear:eq-variant-p1}
\begin{aligned}
   &N_0^{-1} N_1^{-2} N_2^{-1} \max\Big(
    \| h\|_{n_0 n_3 \rightarrow n_1 n_2} \| h\|_{n_0 n_2 \rightarrow n_1 n_3}, 
    \| h \|_{n_0 n_3 \rightarrow n_1 n_2} \| h \|_{n_0 n_1 n_2 \rightarrow n_3}, \\
    &\| h\|_{n_0 n_2 \rightarrow n_1 n_3} \| h \|_{n_0 n_1 n_3 \rightarrow n_2}, 
    \| h \|_{n_0 n_1 \rightarrow n_2 n_3}^2 \Big) \lesssim N_2^{-1/4-\epsilon}, 
\end{aligned}
\end{equation}
where 
\begin{align*}
h_{n_0 n_1 n_2 n_3} &:=
    1_{N_{12}}(n_{12}) 1_{N_{23}}(n_{23})  1_{N_{13}}(n_{13})
    \Big(\prod_{j=0}^3 1_{N_j}(n_j)\Big)  \Big( \prod_{j=2}^{3} 1_{Q_j}(n_j) \Big) \mathbf{1}\{ n_0=n_{123} \} 
\mathbf{1}\{ | \Omega-m|\leq 1 \}, \\
\Omega &:= \sum_{j=0}^3 (\pm_j) \langle n_ j \rangle, \qquad  m \in \Z,
\end{align*} 
and $Q_2$ and $Q_3$ are boxes of sidelength $\sim N_{23}$. The factor $1_{N_{13}}(n_{13})$ is only included to match the notation of Lemma \ref{counting:basetensor_est}, which will be used in all of the estimates below, since our frequency assumptions already imply $|n_{13}|\sim N_{13}$.  To simplify the notation, we set $\beta := 1/100$. \\

\emph{First argument in \eqref{bilinear:eq-variant-p1}:} It holds that
\begin{equation*}
\| h\|_{n_0 n_3 \rightarrow n_1 n_2}^2 \lesssim  \min(N_0,N_{12})^{-1} N_0^{3} \min(N_1,N_{12})^{-1} N_1^{3} \lesssim N_0^2 N_1^{2+\beta}
\end{equation*}
and 
\begin{equation*}
    \| h\|_{n_0 n_2 \rightarrow n_1 n_3}^2 
    \lesssim  \min(N_0, N_{13})^{-1} N_0^{3} \min(N_1,N_2)^{-1} N_1^{3}
    \lesssim N_0^2 N_1^2. 
\end{equation*}
As a result, 
\begin{align*}
N_0^{-1} N_1^{-2} N_2^{-1} 
    \| h\|_{n_0 n_3 \rightarrow n_1 n_2} \| h\|_{n_0 n_2 \rightarrow n_1 n_3} 
\lesssim \, N_0 N_1^{\beta/2} N_2^{-1} \lesssim N_2^{-1+3\beta/2}.
\end{align*}~\\

\emph{Second argument in \eqref{bilinear:eq-variant-p1}:} We have that 
\begin{equation*}
\| h \|_{n_0 n_3 \rightarrow n_1 n_2}^2 \lesssim \min(N_0,N_{12})^{-1} N_0^3 \min(N_{12},N_1)^{-1} N_1^3 \lesssim N_0^2 N_1^{2+\beta}
\end{equation*}
and 
\begin{equation*}
\| h \|_{n_0 n_1 n_2 \rightarrow n_3}^2 \lesssim \med\big( N_0,N_1,N_3 \big)^{-1} N_0^3 N_1^3 \lesssim N_0^3 N_1^{2+\beta}.
\end{equation*}
In the last inequality, we used that $N_0 \ll N_{12}$ implies that $N_{12}\sim N_3$. As a result, 
\begin{equation*}
N_0^{-1} N_1^{-2} N_2^{-1}\| h \|_{n_0 n_3 \rightarrow n_1 n_2} \| h \|_{n_0 n_1 n_2 \rightarrow n_3} 
\lesssim N_0^{3/2} N_1^{\beta} N_2^{-1} \lesssim N_2^{-1+5\beta/2}.
\end{equation*}~\\ 

\emph{Third argument in \eqref{bilinear:eq-variant-p1}:} It holds that
\begin{equation*}
\| h \|_{n_0 n_2 \rightarrow n_1 n_3}^2 \lesssim \min(N_0,N_{13})^{-1} N_0^3 \min(N_1,N_2)^{-1} N_1^3 \lesssim N_0^2 N_1^2
\end{equation*}
and 
\begin{equation*}
\| h \|_{n_0 n_1 n_3 \rightarrow n_2}^2 \lesssim \med(N_0,N_1,N_2)^{-1} N_0^3 N_1^3 \lesssim N_0^3 N_1^2. 
\end{equation*}
As a result, we obtain that 
\begin{equation*}
N_0^{-1} N_1^{-2} N_2^{-1} \| h\|_{n_0 n_2 \rightarrow n_1 n_3} \| h \|_{n_0 n_1 n_3 \rightarrow n_2} \lesssim N_0^2 N_2^{-1} \lesssim N_2^{-1+2\beta}. 
\end{equation*}

\emph{Fourth argument in \eqref{bilinear:eq-variant-p1}:} Using the box localization and the box counting lemma (Lemma \ref{counting:lem-annuli-box}), we obtain that
\begin{equation*}
\| h \|_{n_0 n_1 \rightarrow n_2 n_3}^2 \lesssim \min(N_0,N_{23})^{-1} N_0^3 \min(N_{23}^3, N_2^2)
\lesssim N_0^3 \min(N_{23}^3, N_2^2).
\end{equation*}
Since $N_{23} \lesssim \max(N_0,N_1)$, we obtain that
\begin{align*}
  N_0^{-1} N_1^{-2} N_2^{-1}   \| h \|_{n_0 n_1 \rightarrow n_2 n_3}^2
  &\lesssim N_0^2 N_1^{-2} N_2^{-1} \min\Big( \max(N_0,N_1)^3 , N_2^2 \Big) \\
  &\lesssim N_0^2 N_1^{-2} N_2^{-1} \big(\max(N_0,N_1)^{3}\big)^{2/3} \big(N_2^{2}\big)^{1/3}\\
  &\lesssim N_0^4 N_2^{-1/3} \lesssim N_2^{-1/3+4\beta}. 
\end{align*}
\end{proof}
\section{Linear random operators}\label{section:linear}

In this section, we analyze (frequency-localized versions of) the two random linear operators, 
\begin{equation}\label{linear:eq-motivation-quad}
w \mapsto P_{\leq N}  \Big[ \squadratic[\leqN] P_{\leq N} w \Big]
\end{equation}
and 
\begin{equation}\label{linear:eq-motivation-lincub}
w \mapsto P_{\leq N}  \Big[ \slinear[blue][\leqN] \scubic[\leqN] P_{\leq N} w \Big].
\end{equation}
In contrast to the bilinear estimate in Section \ref{section:bilinear}, both linear random operators are estimated directly using linear tensor norms. From an abstract perspective, this section is therefore closely related to similar linear tensor estimates in \cite{DNY20}. From a problem-specific perspective, however, the situation is rather different than in \cite{DNY20}. As was already discussed in our introduction (Subsection \ref{section:intro-main-ideas}) and in \cite[Subsection 9.2]{DNY20}, each dispersion relation requires different tensor estimates, and the tensor estimates entering into the analysis of \eqref{linear:eq-motivation-quad} and \eqref{linear:eq-motivation-lincub} are one of the novelties of our article. 

\subsection{Linear random operator involving the quadratic object}\label{section:quad}

We first treat the frequency-localized variants of the random linear operator \eqref{linear:eq-motivation-quad}. In our Ansatz (see Section \ref{section:ansatz}), we isolated several terms involving the quadratic stochastic object $\squadratic[\leqN]$. In \eqref{ansatz:eq-YN-1} and \eqref{ansatz:eq-YN-2}, we encountered the terms 
\begin{align*}
P_{\leq N} \bigg[ 9 \squadratic[\leqN] \squintic[\leqN] - \Gamma_{\leq N} \Big( \scubic[\leqN] + \slinear[green][\leqN] \Big) - 18 \mathfrak{C}^{(1,5)}_{\leq N} P_{\leq N} \slinear[blue][\leqN]  - 9 \Big( 2 \HLL + \HHL \Big) \Big( \, \slinear[blue][\leqN], \slinear[blue][\leqN], \squintic[\leqN] \Big) \bigg]
\end{align*}
and
\begin{equation*}
 3 P_{\leq N} \bigg[  \squadratic[\leqN] P_{\leq N} \slinear[green][\leqM] - \Big( 2 \HLL + \HHL \Big) \Big( \, \slinear[blue][\leqN], \slinear[blue][\leqN], P_{\leq N} \slinear[green][\leqM] \Big)\bigg]. 
\end{equation*}
Fortunately, both of them are explicit stochastic objects. The first object is treated in Section \ref{section:analytic2} below and requires additional ingredients. The second object has already been estimated in Section \ref{section:analytic}. In \eqref{ansatz:eq-YN-3} and \eqref{ansatz:eq-YN-4}, we encountered the terms 
\begin{equation}\label{linear:eq-motivation-quad-Xone}
P_{\leq N} \bigg[ \squadratic[\leqN] \XXone - 
    \Big( 2 \HLL + \HHL \Big) \Big( \slinear[blue][\leqN] , \slinear[blue][\leqN], \XXone \Big) \bigg]
\end{equation}
and 
\begin{equation}\label{linear:eq-motivation-quad-Xtwo}
P_{\leq N} \bigg[ 3 \squadratic[\leqN] \XXtwo - 
    \Big( 6 \HLL + 3 \HHL \Big) \Big( \slinear[blue][\leqN] , \slinear[blue][\leqN], \XXtwo \Big)
    + \Gamma_{\leq N} \Big( 3 \squintic[\leqN] + v_{\leq N} \Big) \bigg]. 
\end{equation}
While the methods of this section would control \eqref{linear:eq-motivation-quad-Xone} and \eqref{linear:eq-motivation-quad-Xtwo} at least for some (but not all) frequency-interactions, it is easier to directly use the para-controlled structure. As a result, the analysis of \eqref{linear:eq-motivation-quad-Xone} and \eqref{linear:eq-motivation-quad-Xtwo} is postponed entirely until Section \ref{section:para}. The only remaining term containing $\squadratic[\leqN]$, which is the subject of this subsection, is contained in \eqref{ansatz:eq-YN-5} and given by 
\begin{equation}\label{linear:eq-motivation-quad-Y}
P_{\leq N} \bigg[ \squadratic[\leqN] Y_{\leq N} - 
    \Big( 2 \HLL + \HHL + \RES \Big) \Big( \slinear[blue][\leqN] , \slinear[blue][\leqN], Y_{\leq N} \Big) \bigg].
\end{equation}
The main estimate of \eqref{linear:eq-motivation-quad-Y} is the content of the next proposition.  

\begin{proposition}[\protect{Random linear operator involving $\squadratic[\leqN]$}]\label{linear:prop-quad}
For all $p \geq 2$ and $T\geq 1$, it holds that
\begin{equation}
\begin{aligned}
&\mathbb{E} \bigg[ \sup_{N} \sup_{\Jc} \bigg\| Y \mapsto P_{\leq N}  \bigg[ \squadratic[\leqN] P_{\leq N} Y \\
&\hspace{15ex}- 
\Big( 2 \HLL + \HHL + \RES \Big) \Big( \slinear[blue][\leqN] , \slinear[blue][\leqN],  P_{\leq N} Y \Big) \bigg] \bigg\|_{X^{1/2+\delta_2,b}(\Jc) \rightarrow X^{-1/2+\delta_2,b_+-1}(\Jc)}^p \bigg]^{1/p} \\
    &\lesssim  T^{\alpha }p, 
\end{aligned}
\end{equation}
where the supremum is taken over all closed intervals $0 \in \Jc \subseteq [-T,T]$. 
\end{proposition}

The proof of Proposition \ref{linear:prop-quad} will be postponed until the end of the subsection.
From a notational standpoint, the formula for the random operator in \eqref{linear:eq-motivation-quad-Y} is rather tedious. To simplify the notation, we therefore make the following definitions. 

\begin{definition}[The $\Quad$-operators]\label{linear:def-quad-operator} 
We introduce the following two random operators: 
\begin{enumerate}
    \item For all $N\geq 1$, we denote the random operator in \eqref{linear:eq-motivation-quad-Y} by $\Quad_{\leq N}$, i.e., we define 
    \begin{equation*}
        \Quad_{\leq N}(Y):= P_{\leq N} \bigg[ \squadratic[\leqN] Y_{\leq N} - 
    \Big( 2 \HLL + \HHL + \RES \Big) \Big( \slinear[blue][\leqN] , \slinear[blue][\leqN], Y_{\leq N} \Big) \bigg].
    \end{equation*}
    \item For all $N_0,N_1,N_2,N_3,N_{13},N_{23}\geq 1$, we define the frequency-localized variant of $\Quad_{\leq N}$ by 
    \begin{align*}
    \Quad[N_\ast](Y) := \sum_{ \substack{n_0,n_1,n_2 ,n_3 \in \Z^3 \colon \\ n_0 = n_{123} }}
    \bigg[ 1_{N_{13}}(n_{13}) 1_{N_{23}}(n_{23}) \Big( \prod_{j=0}^3 1_{N_j}(n_j) \Big) \lcol  \slinear[blue][N_1](n_1) \slinear[blue][N_2](n_2) \rcol \widehat{Y}(n_3) e^{i \langle n_0 ,x \rangle} \bigg]. 
    \end{align*}
\end{enumerate}
\end{definition}
In the proof of Proposition \ref{linear:prop-quad} below, we will see that $\Quad_{\leq N}$ can be written as a linear combination of the frequency-localized operators $\Quad[N_\ast]$ (and additional easier terms). The main ingredient in the proof of Proposition \ref{linear:prop-quad} is the following lemma, which estimates the frequency-localized random operator $\Quad[N_\ast]$.

\begin{lemma}[\protect{Estimate of $\Quad [N_\ast]$}]\label{linear:lem-quad-dyadic}
Let $p\geq 2$, let $T\geq 1$, and let $N_0,N_1,N_2,N_3,N_{12},N_{13}$ be frequency-scales satisfying the non-resonant conditions 
\begin{equation}\label{linear:eq-quad-dyadic-nonres}
 \quad N_{13} \gtrsim N_2^\eta \quad \text{and} \quad N_{23}\gtrsim N_1^\eta. 
\end{equation}
Then, it holds that 
\begin{equation*}
\mathbb{E} \Big[  \sup_{\Jc} \big\|  \Quad[N_\ast] \big\|_{X^{1/2,b}(\Jc)\rightarrow X^{-1/2,b_+-1}(\Jc)}^p \Big]^{1/p} \lesssim N_{\textup{max}}^{\epsilon} \big( N_0^{-\eta/2} + N_3^{-1/2} \big)  T^{\alpha} p. 
\end{equation*}
If the frequency-scales instead satisfy the double-resonance conditions
\begin{equation}\label{linear:eq-quad-dyadic-doubleres}
 \quad N_{13} \lesssim N_2^\eta \quad \text{and} \quad N_{23}\lesssim N_1^\eta,
\end{equation}
then it  holds that 
\begin{equation*}
\mathbb{E} \Big[ \sup_{\Jc} \big\| \Quad[N_\ast] \big\|_{X^{1/2,b}(\Jc)\rightarrow X^{-1/2,b_+-1}(\Jc)}^p \Big]^{1/p} \lesssim N_{\textup{max}}^{-1+\epsilon}  T^\alpha p. 
\end{equation*}
\end{lemma}

\begin{remark}
The main estimate in Lemma \ref{linear:lem-quad-dyadic} concerns the non-resonant case \eqref{linear:eq-quad-dyadic-nonres}. As will be clear from the proof, the estimate in the double-resonance case \eqref{linear:eq-quad-dyadic-doubleres} is almost trivial. We emphasize that Lemma  \ref{linear:lem-quad-dyadic} does not cover the case of a single resonance such as $N_{13}\sim N_2$ and $N_{23}\lesssim N_1^\eta$.
\end{remark}

\begin{proof}
We first treat the non-resonant case \eqref{linear:eq-quad-dyadic-nonres}. Using the reduction arguments in Subsection \ref{prep:remark-reduction}, it suffices to prove the random tensor estimate
\begin{equation}
\begin{aligned}\label{linear:eq-quad-dyadic-p1}
&N_0^{-1/2} N_1^{-1} N_2^{-1} N_3^{-1/2} 
 \mathbb{E} \bigg[ \Big\| \, \widehat{Y} \mapsto \sum_{n_1,n_2,n_3 \in \Z^3} h_{n_0 n_1 n_2 n_3 } \lcol g_{n_1} g_{n_2} \rcol \widehat{Y}(n_3) \Big\|_{\ell^2\rightarrow \ell^2}^p \bigg]^{1/p} \\
\lesssim& N_{\textup{max}}^\epsilon \big( N_0^{-\eta/2} +N_3^{-1/2} \big) p. 
\end{aligned}
\end{equation}
Here, the tensor $h$ is given by 
\begin{align*}
h_{n_0 n_1 n_2 n_3} &= 1_{N_{13}}(n_{13}) 1_{N_{23}}(n_{23}) \Big( \prod_{j=0}^3 1_{N_j}(n_j) \Big) \mathbf{1}\{ n_0 = n_{123} \} \mathbf{1}\{ |\Omega-m |\leq 1 \}, \\
\Omega &= \sum_{j=0}^3 (\pm_j) \langle n_j \rangle, \qquad \text{and} \qquad 
m \in \Z. 
\end{align*}
Up to the additional dyadic localization of $n_{13}$ and $n_{23}$, $h$ agrees with the base tensor from Section \ref{section:counting}.
After utilizing the moment method (see Proposition \ref{counting:prop-moment} and Remark \ref{Bee3}), it then remains to prove that 
\begin{equation}\label{linear:eq-quad-dyadic-p2}
\begin{aligned}
&N_0^{-1/2} N_1^{-1} N_2^{-1} N_3^{-1/2}  
\max\Big( \| h \|_{n_0 n_1 n_2 \rightarrow n_3}, 
\| h \|_{n_0 n_1 \rightarrow n_2 n_3}, 
\| h \|_{n_0 n_2 \rightarrow n_1 n_3}, 
\| h \|_{n_0 \rightarrow n_1 n_2 n_3} \Big) \\
\lesssim& N_0^{-\eta/2} + N_3^{-1/2}. 
\end{aligned}
\end{equation}
The four arguments in the maximum of \eqref{linear:eq-quad-dyadic-p2} can all be treated using the base tensor estimate (Lemma \ref{counting:basetensor_est}). For expository purposes, we separate the remaining argument into four steps. 

\emph{Estimate of the first argument in \eqref{linear:eq-quad-dyadic-p2}:} Using the base tensor estimate, it holds that 
\begin{align*}
N_0^{-1/2} N_1^{-1} N_2^{-1} N_3^{-1/2}\| h \|_{n_0 n_1 n_2 \rightarrow n_3} 
&\lesssim N_0^{-1/2} N_1^{-1} N_2^{-1} N_3^{-1/2} \med(N_0,N_1,N_2)^{3/2} \min(N_0,N_1,N_2) \\ &\lesssim N_3^{-1/2}. 
\end{align*}

\emph{Estimate of the second argument in \eqref{linear:eq-quad-dyadic-p2}:} Using the base tensor estimate, it holds that 
\begin{align}
&N_0^{-1/2} N_1^{-1} N_2^{-1} N_3^{-1/2}  \| h \|_{n_0 n_1 \rightarrow n_2 n_3} \notag \\
\lesssim&\, N_0^{-1/2} N_1^{-1} N_2^{-1} N_3^{-1/2} \min(N_{23},N_0,N_1)^{-1/2} \min(N_0,N_1)^{3/2} 
  \label{linear:eq-quad-dyadic-p3} \\
  &\times \min(N_{23},N_2,N_3)^{-1/2} \min(N_2,N_3)^{3/2} \notag. 
\end{align}
Due to our assumption $N_{23}\gtrsim N_1^\eta$, it holds that 
$\min(N_{23},N_0,N_1)\gtrsim \min(N_0,N_1)^{\eta}$. Together with the trivial estimate 
$\min(N_{23},N_2,N_3)\gtrsim 1$, it follows that
\begin{equation*}
\eqref{linear:eq-quad-dyadic-p3} \lesssim N_0^{-1/2} N_1^{-1} N_2^{-1} N_3^{-1/2} 
\min(N_0,N_1)^{3/2-\eta/2} \min(N_2,N_3)^{3/2} \lesssim N_0^{-\eta/2}. 
\end{equation*}

\emph{Estimate of the third argument in \eqref{linear:eq-quad-dyadic-p2}:}
Due to symmetry in the $n_1$ and $n_2$-variables, the third argument in  \eqref{linear:eq-quad-dyadic-p2} can be treated exactly as the second argument. 

\emph{Estimate of the fourth argument in \eqref{linear:eq-quad-dyadic-p2}:} Using the base tensor estimate, it holds that 
\begin{align*}
   &N_0^{-1/2} N_1^{-1} N_2^{-1} N_3^{-1/2}   \| h \|_{n_0 \rightarrow n_1 n_2 n_3} \\
   &\lesssim N_0^{-1/2} N_1^{-1} N_2^{-1} N_3^{-1/2}  \med(N_1,N_2,N_3)^{3/2} \min(N_1,N_2,N_3)\\
   &\lesssim N_0^{-1/2}. 
\end{align*}
This contribution is better than the $N_0^{-\eta/2}$-term from the second and third argument. This completes our argument in the non-resonant case \eqref{linear:eq-quad-dyadic-nonres}.\\

It remains to treat the double-resonant case \eqref{linear:eq-quad-dyadic-doubleres}. Similar as above, it suffices to prove the random tensor estimate 
\begin{equation}
\begin{aligned}\label{linear:eq-quad-dyadic-p4}
&N_0^{-1/2} N_1^{-1} N_2^{-1} N_3^{-1/2} 
\mathbb{E} \bigg[ \Big\| \, \widehat{Y} \mapsto \sum_{n_1,n_2,n_3 \in \Z^3} h_{n_0 n_1 n_2 n_3 } \lcol g_{n_1} g_{n_2} \rcol \widehat{Y}(n_3) \Big\|_{\ell^2\rightarrow \ell^2}^p \bigg]^{1/p} \lesssim N_{\textup{max}}^{-1} p. 
\end{aligned}
\end{equation}
While we could estimate the left-hand side of \eqref{linear:eq-quad-dyadic-p4} using the moment method, we can also proceed using the following much simpler argument\footnote{The same estimate can also be derived by first using the moment method and then estimating all tensor norms by the Hilbert-Schmidt norm. However, we prefer the argument presented here, since it emphasizes that the moment method is not needed.}. Using Cauchy-Schwarz in the $n_3$-variable, it holds that 
\begin{align*}
&N_0^{-1/2} N_1^{-1} N_2^{-1} N_3^{-1/2} 
\mathbb{E} \bigg[ \Big\| \, \widehat{Y} \mapsto \sum_{n_1,n_2,n_3 \in \Z^3} h_{n_0 n_1 n_2 n_3 } \lcol g_{n_1} g_{n_2} \rcol \widehat{Y}(n_3) \Big\|_{\ell^2\rightarrow \ell^2}^p \bigg]^{1/p} \\
\lesssim&  N_0^{-1/2} N_1^{-1} N_2^{-1} N_3^{-1/2} \mathbb{E} \bigg[ \Big\| \sum_{n_1,n_2 \in \Z^3} h_{n_0 n_1 n_2 n_3 } \lcol g_{n_1} g_{n_2} \rcol
\Big\|_{n_0 n_3}^p \bigg]^{1/p}\\
\lesssim& N_0^{-1/2} N_1^{-1} N_2^{-1} N_3^{-1/2}  \big\| h \big\|_{n_0 n_1 n_2 n_3} p. 
\end{align*}
Under the constraint $n_0 =n_{123}$, the four frequencies $n_0,n_1,n_2$, and $n_3$ are uniquely determined by $n_{13}$, $n_{23}$, and either one of the frequencies $n_0,n_1,n_2,n_3$. Without utilizing any dispersive effects and only using the double-resonance conditions \eqref{linear:eq-quad-dyadic-doubleres}, it follows that 
\begin{align*}
    N_0^{-1/2} N_1^{-1} N_2^{-1} N_3^{-1/2}  \big\| h \big\|_{n_0 n_1 n_2 n_3} 
    \lesssim& \, N_0^{-1/2} N_1^{-1} N_2^{-1} N_3^{-1/2} N_{13}^{3/2} N_{23}^{3/2} \min(N_0,N_1,N_2,N_3)^{3/2} \\
    \lesssim&\, N_0^{-1/2} N_1^{-1+\frac{3\eta}{2}} N_2^{-1+\frac{3\eta}{2}} N_3^{-1/2} \min(N_0,N_1,N_2,N_3)^{3/2} \\
    \lesssim& \, N_{\textup{max}}^{-1}. \qedhere
    \end{align*}
\end{proof}

Equipped with the frequency-localized estimate (Lemma \ref{linear:lem-quad-dyadic}), we now prove the main proposition (Proposition \ref{linear:prop-quad}).
\begin{proof}[Proof of Proposition \ref{linear:prop-quad}]
From Definition \ref{ansatz:def-para-product} and Definition \ref{linear:def-quad-operator}, it follows that 
\begin{align}
P_{\leq N} \Big[ \squadratic[\leqN] P_{\leq N} Y \Big] 
&= \sum_{N_0,N_1,N_2,N_3 \leq N} \sum_{N_{13},N_{23}} \Quad[N_\ast](Y), \label{linear:eq-quad-p1}\\
\Big( 2 \HLL+\HHL\Big)\Big( \slinear[blue][\leqN], \slinear[blue][\leqN], P_{\leq N} Y \Big) 
&= \sum_{\substack{ N_0 ,N_1, N_2, N_3 \leq N \colon \\ N_3 \leq  \max(N_1,N_2)^\eta }} 
\sum_{N_{13},N_{23}} \Quad[N_\ast](Y). \label{linear:eq-quad-p2}
\end{align}
In the definition of $\RES$ (Definition \ref{ansatz:def-para-product}), we did not include the Wick-ordering of the product in $\slinear[blue][N_1]$ and $\slinear[blue][N_2]$. We  split $\RES$ into the Wick-ordered random operators and deterministic operators by writing
\begin{align}
\RES\Big( \slinear[blue][\leqN], \slinear[blue][\leqN], P_{\leq N} Y \Big) 
&= \sum_{\substack{ N_0 ,N_1, N_2, N_3 \leq N \colon \\ N_3 > \max(N_1,N_2)^\eta }} \sum_{N_{13},N_{23}}  \big( \mathbf{1}\{ N_{13} \leq N_2^\eta\} + \mathbf{1}\{N_{23}\leq N_1^\eta \} \big) \Quad[N_\ast](Y) \label{linear:eq-quad-p3} \\
&+ \sum_{\substack{ N_0 ,N_1, N_2, N_3 \leq N \colon \\ N_3 > \max(N_1,N_2)^\eta }} \sum_{N_{13},N_{23}}  \big( \mathbf{1}\{ N_{13} \leq N_2^\eta\} + \mathbf{1}\{N_{23}\leq N_1^\eta \} \big) \Quad^{\diamond}[N_\ast](Y), \label{linear:eq-quad-p4}
\end{align}
where we define the deterministic operator
\begin{equation}
\Quad^{\diamond}[N_\ast](Y) := \sum_{\substack{n_0,n_1,n_2,n_3 \in \Z^3 \colon \\ n_0 = n_{3}, \,  n_{12}= 0 }} \bigg[ 1_{N_{13}}(n_{13}) 1_{N_{23}}(n_{23}) \Big( \prod_{j=0}^3 1_{N_j}(n_j) \Big) \frac{1}{\langle n_1 \rangle^2} \widehat{Y}(n_3) e^{i\langle n_0 , x \rangle} \bigg]. 
\end{equation}
After subtracting \eqref{linear:eq-quad-p2}, \eqref{linear:eq-quad-p3}, and \eqref{linear:eq-quad-p4} from \eqref{linear:eq-quad-p1}, we obtain that 
\begin{align}
\Quad_{\leq N}(Y) 
&= \sum_{\substack{ N_0 ,N_1, N_2, N_3 \leq N \colon \\ N_3 > \max(N_1,N_2)^\eta }}
\sum_{\substack{N_{13}\colon \\ N_{13} > N_2^\eta }} \sum_{\substack{N_{23} \colon \\ N_{23}> N_1^\eta }} 
\Quad[N_\ast](Y) \label{linear:eq-quad-p5}\\
&- \sum_{\substack{ N_0 ,N_1, N_2, N_3 \leq N \colon \\ N_3 > \max(N_1,N_2)^\eta }}
\sum_{\substack{N_{13}\colon \\ N_{13} \leq N_2^\eta }} \sum_{\substack{N_{23} \colon \\ N_{23}\leq  N_1^\eta }} 
\Quad[N_\ast](Y)  \label{linear:eq-quad-p6}\\
&- \sum_{\substack{ N_0 ,N_1, N_2, N_3 \leq N \colon \\ N_3 > \max(N_1,N_2)^\eta }}
\sum_{N_{13},N_{23}} \big( \mathbf{1}\{ N_{13} \leq N_2^\eta\} + \mathbf{1}\{N_{23}\leq N_1^\eta \} \big) \Quad^{\diamond}[N_\ast](Y). \label{linear:eq-quad-p7}
\end{align}
We now  estimate the contributions of \eqref{linear:eq-quad-p5}, \eqref{linear:eq-quad-p6}, and \eqref{linear:eq-quad-p7} to $\Quad_{\leq N}$ separately. Using Lemma \ref{linear:lem-quad-dyadic}, the contribution of \eqref{linear:eq-quad-p5} can be estimated by 
\begin{equation}\label{linear:eq-quad-p8}
\begin{aligned}
&\sum_{\substack{ N_0 ,N_1, N_2, N_3 \leq N \colon \\ N_3 > \max(N_1,N_2)^\eta }}
\sum_{\substack{N_{13}\colon \\ N_{13} > N_2^\eta }} \sum_{\substack{N_{23} \colon \\ N_{23}> N_1^\eta }} \mathbb{E} \bigg[ \sup_{\Jc} \Big\|  \Quad[N_\ast](Y) \Big\|_{X^{1/2+\delta_2,b}(\Jc) \rightarrow X^{-1/2+\delta_2,b_+-1}(\Jc)}^p \bigg]^{1/p} \\
\lesssim p& T^{\alpha} \sum_{\substack{ N_0 ,N_1, N_2, N_3 \leq N \colon \\ N_3 > \max(N_1,N_2)^\eta }}
\sum_{\substack{N_{13}\colon \\ N_{13} > N_2^\eta }} \sum_{\substack{N_{23} \colon \\ N_{23}> N_1^\eta }} N_{\textup{max}}^\epsilon N_0^{\delta_2} \big( N_0^{-\eta/2} + N_3^{-1/2} \big) N_3^{-\delta_2}. 
\end{aligned}
\end{equation}
Under the frequency-restrictions in \eqref{linear:eq-quad-p8}, which guarantee that $N_3 \gtrsim N_{\textup{max}}^{\eta}$, we obtain 
\begin{equation*}
    N_{\textup{max}}^\epsilon N_0^{\delta_2} \big( N_0^{-\eta/2} + N_3^{-1/2} \big) N_3^{-\delta_2}
    \lesssim  N_{\textup{max}}^{\epsilon-\delta_2 \eta} 
    + N_{\textup{max}}^{\epsilon+\delta_2 - \big(1/2+\delta_2\big) \eta} 
    \lesssim N_{\textup{max}}^{-\epsilon}. 
\end{equation*}
The contribution of the double-resonance term \eqref{linear:eq-quad-p6} can be estimated directly from Lemma \ref{linear:lem-quad-dyadic} and requires no further frequency-scale considerations. Finally, we note that $\Quad^{\diamond}(Y)$ is a Fourier-multiplier with symbol 
\begin{equation*}
n_3 \mapsto \sum_{\substack{n_1,n_2 \in \Z^3 \colon \\ n_{12} = 0 }} \bigg[ 
 1_{N_{13}}(n_{13}) 1_{N_{23}}(n_{23}) \Big( \prod_{j=0}^3 1_{N_j}(n_j) \Big) \frac{1}{\langle n_1 \rangle^2} \bigg]. 
\end{equation*}
Under the frequency-restrictions in \eqref{linear:eq-quad-p7}, this symbol can be estimated by 
\begin{align*}
&\bigg| \sum_{\substack{n_1,n_2 \in \Z^3 \colon \\ n_{12} = 0 }} \bigg[ 
 1_{N_{13}}(n_{13}) 1_{N_{23}}(n_{23}) \Big( \prod_{j=0}^3 1_{N_j}(n_j) \Big) \frac{1}{\langle n_1 \rangle^2} \bigg] \bigg| \\
 \lesssim& \, \mathbf{1}\{ N_1 = N_2 \} N_1^{-2} \min(N_{13},N_{23})^3 \\
 \lesssim& \mathbf{1}\{ N_1 = N_2 \} N_1^{-2} N_3^{-1} \max(N_{1},N_{13})\min(N_{13},N_{23})^3 \\
 \lesssim& \mathbf{1}\{ N_1 = N_2 \} N_1^{-1+3\eta} N_3^{-1} \\
 \lesssim& N_{\textup{max}}^{-1+3\eta}. 
\end{align*}
As a result, the contribution of \eqref{linear:eq-quad-p7} can be estimated directly in $L_t^2 H_x^{-1/2+\delta_2}$, which is stronger than $X^{-1/2+\delta_2,b_+-1}$.
\end{proof}

\subsection{Linear random operator involving the linear and cubic objects}\label{section:lincub}
In this section, we treat (frequency-localized) versions of the linear random operator
\begin{equation}\label{linear:eq-lincub-motivation-1}
w \mapsto P_{\leq N} \Pi^\ast_{\leq N}\Big( \slinear[blue][\leqN],  \scubic[\leqN], P_{\leq N} w \Big). 
\end{equation}
In the evolution equation for $Y_{\leq N}$, i.e., in \eqref{ansatz:eq-YN-6}, we encountered several terms of this form. If $w$ is replaced by the cubic stochastic object $\scubic[\leqN]$, we obtain an explicit septic stochastic object, which will be estimated in Section \ref{section:analytic2} below. Thus, it remains to treat the cases when $w$ is of the form 
\begin{equation}\label{linear:eq-lincub-motivation-2}
\slinear[green][\leqM], \, \squintic[\leqN], \, \XXone, \, \XXtwo, \, \text{or} \, \, Y_{\leq N}. 
\end{equation}

In contrast to Subsection \ref{section:quad}, the main estimate of this subsection (Proposition \ref{linear:prop-lincub}) concerns all five terms in \eqref{linear:eq-lincub-motivation-2}. In the first four cases of \eqref{linear:eq-lincub-motivation-2}, however, it only covers most (but not all) relevant frequency-interactions. The remaining frequency-interactions will be treated  in Section \ref{section:proof-main-estimates} using additional ingredients from Section \ref{section:analytic} and Section \ref{section:para}. In order to concisely state our main estimate, we first make the following definition.

\begin{definition}[High$\times$low$\times$high-interactions]\label{linear:def-hlh}
For any $N\geq 1$ and function $w\colon \R \times \T^3 \rightarrow \R$, we define
\begin{equation}\label{linear:eq-hlh}
\HLH \Big( \slinear[blue][\leqN], \scubic[\leqN], P_{\leq N} w \Big) 
 = \sum_{ \substack{N_1,N_{234},N_5 \leq N\colon \\ N_{234}\leq N_1^{\nu}, \\ N_1 \sim N_5 >N_1^\eta}} 
 \slinear[blue][N_1] P_{N_{234}} \scubic[\leqN] P_{N_5} w. 
\end{equation}
\end{definition}
We encourage the reader to ignore the condition $N_5 > N_1^\eta$ in Definition \ref{linear:def-hlh}. It is only included to avoid double-counting and, if $N_1$ is sufficiently large, directly follows from $ N_1 \sim N_5$. We also point out that, due to the choice of our parameters, the condition $N_{234} \leq N_1^\nu$ is weaker than $N_{234}\leq N_1^\eta$. This weaker condition is needed for technical reasons. Otherwise, the proof of Proposition \ref{linear:prop-lincub} would require $\delta_1 \ll \eta$, which would be disadvantageous in the regularity estimates of Subsection \ref{section:para-regularity}, and therefore we do not impose this condition in Section \ref{section:prelim}.

In the quintic case $w=\squintic[\leqN]$, we eventually have to replace the product 
\begin{equation*}
    \slinear[blue][N_1] P_{N_5}w \qquad \text{in \eqref{linear:eq-hlh} by} \qquad \slinear[blue][N_1] P_{N_5} w - \mathfrak{C}^{(1,5)}_{\leq N}[N_1,N_5].
\end{equation*}
From a notational perspective, it is more convenient to perform this replacement later, as we prefer to not distinguish between the first four terms in \eqref{linear:eq-lincub-motivation-2}
throughout this subsection. 

We emphasize that, despite the similar notation, the operators $\HHL$ (from Definition \ref{ansatz:def-para-product}) and $\HLH$  in  \eqref{linear:eq-hlh} are rather different. While high$\times$high$\times$low-interactions always contains two linear stochastic objects, high$\times$low$\times$high-interactions always contains a linear and cubic stochastic object. Furthermore, while the frequency-restrictions in our high$\times$high$\times$low-interactions only impose $N_1^\eta<N_2$, the frequency-restrictions in our high$\times$low$\times$high-interactions impose $N_1 \sim N_5$, which is much stronger.

Equipped with Definition \ref{linear:def-hlh}, we now state our main estimate. 

\begin{proposition}[\protect{Linear random operator involving $\slinear[blue][\leqN] \scubic[\leqN]$}]\label{linear:prop-lincub}
Let $T\geq 1$ and $p\geq 2$. Then, it holds that 
\begin{equation}
\begin{aligned}
 &\E \bigg[ \sup_N \sup_{\Jc} \Big\| w \mapsto P_{\leq N}  \Big[ \slinear[blue][\leqN] \scubic[\leqN] P_{\leq N} w \\
&\hspace{12ex}- \Big( \HLL + \HLH \Big) \Big( \slinear[blue][\leqN], \scubic[\leqN], P_{\leq N} w \Big) \Big]
\Big\|_{X^{1/2-\delta_1,b}(\Jc)\rightarrow X^{-1/2+\delta_2,b_+-1}(\Jc)}^p \bigg]^{1/p} \\
\lesssim& \, T^{\alpha} p^2,
\label{linear:eq-lincub-w}
\end{aligned}
\end{equation}
where the supremum is taken over all closed intervals. Furthermore, it holds that  
\begin{equation}
\begin{aligned}
 &\E \bigg[ \sup_N \sup_{\Jc} \Big\| Y \mapsto P_{\leq N}  \Big[ \slinear[blue][\leqN] \scubic[\leqN] P_{\leq N} Y \\
&\hspace{12ex}-  \HLL  \Big( \slinear[blue][\leqN], \scubic[\leqN], P_{\leq N} Y \Big) \Big]
\Big\|_{X^{1/2+\delta_2,b}(\Jc)\rightarrow X^{-1/2+\delta_2,b_+-1}(\Jc)}^p \bigg]^{1/p} \\
\lesssim& \, T^{\alpha} p^2. 
\label{linear:eq-lincub-Y}
\end{aligned}
\end{equation}
\end{proposition}

We now briefly describe the strategy behind the proof of 
Proposition \ref{linear:prop-lincub}. We first separate the terms in $\slinear[blue][\leqN] \scubic[\leqN]$ containing one or zero pairings. 
Then, we utilize a dyadic decomposition and exhibit the \sine-cancellation (see e.g. Lemma \ref{counting:lem-Sine-estimate}). Finally, we use the quintic and \sine-cancellation tensor estimates.

Just as in Subsection \ref{section:quad}, we simplify the notation for the frequency-localized operators using additional definitions.

\begin{definition}[The $\LinCub$-operators]\label{linear:def-lincub}
For all frequency-scales $N_{234},N_0,\hdots,N_5$, we define 
\begin{equation}
\begin{aligned}
\LinCub^{(5)}[N_\ast](w) :=& 
\sum_{n_0, n_1, \hdots, n_5 \in \Z^3} 
\sum_{\substack{\varphi_1,\varphi_2,\varphi_3,\varphi_4 \\ \in \{\cos,\sin\}} } 
\bigg[ \mathbf{1}\{ n_0 = n_{12345} \} 1_{N_{234}}(n_{234}) \Big(\prod_{j=0}^5 1_{N_j}(n_j)\Big)  \\
&\times \langle n_{234} \rangle^{-1} \Big( \prod_{j=1}^4 \langle n_j \rangle^{-1} \Big) 
\varphi_1(t \langle n_1 \rangle)
\Big( \int_0^t \dt^\prime \sin\big( (t-t^\prime) \langle n_{234} \rangle\big) \prod_{j=2}^4 \varphi_j\big(t^\prime \langle n_j \rangle \big) \Big) \\
&\times \SI[n_j,\varphi_j:1\leq j \leq 4] \widehat{w}(t,n_5) e^{i \langle n_0 ,x \rangle} \bigg]. 
\end{aligned}
\end{equation}
Furthermore, we define
\begin{equation}
\begin{aligned}
&\LinCub^{\sin} [N_\ast](w)  \allowdisplaybreaks[3] \\
:=& \mathbf{1}\big\{ N_1 =N_2 \big\}  \sum_{n_0, n_3, n_4, n_5 \in \Z^3} 
\sum_{\substack{\varphi_3, \varphi_4 \in \\ \{ \cos, \sin \}}} 
\Bigg[ \mathbf{1}\{n_0=n_{345}\}   \Big( \prod_{j=0,3,4,5} 1_{N_j}(n_j) \Big) \langle n_3 \rangle^{-1} \langle n_4 \rangle^{-1}\\
&\times \bigg( \int_0^t \dt^\prime  \Sine[N_{234},N_2](t-t^\prime,n_{34}) \prod_{j=3}^4 \varphi_j\big( t^\prime \langle n_j \rangle \big) \bigg) \, \widehat{w}(t,n_5)  \\
&\times \SI[n_j, \varphi_j \colon 3\leq j \leq 4] \, e^{i\langle n_0 ,x \rangle} 
\Bigg],
\end{aligned}
\end{equation}
where the $\Sine$-kernel is as in Definition \ref{def:sine_kernel}. 
\end{definition}

\begin{remark}\label{linear:rem-lincub} The superscripts of $\LinCub^{(5)}$ and $\LinCub^{\sin}$ in Definition \ref{linear:def-lincub} are motivated by our notation for the tensors in Section \ref{section:counting}, which will be used in the estimates of $\LinCub^{(5)}$ and $\LinCub^{\sin}$. 

\end{remark}

Equipped with Definition \ref{linear:def-lincub}, we now state the following (algebraic) lemma.

\begin{lemma}[Decomposition using $\LinCub$-operators]\label{linear:lem-lincub-decomp}
For all frequency-scales $N\geq 1$ and functions $w\colon \R \times \T^3 \rightarrow \R$, it holds that 
\begin{align}
&P_{\leq N} \bigg[ \slinear[blue][\leqN] \scubic[\leqN] P_{\leq N} w 
- \Big( \HLL + \HLH \Big) \Big( \slinear[blue][\leqN], \scubic[\leqN], P_{\leq N} w \Big) \bigg] \notag \allowdisplaybreaks[3]\\
=& \sum_{\substack{N_0,N_1,N_2,N_3, \\ N_4,N_5,N_{234} \leq N}} \bigg[
\Big( \mathbf{1}\big\{ N_{234}> N_1^\eta \geq N_{5} \big\} 
+ \mathbf{1} \big\{ N_1\not \sim N_5 > N_1^\eta \big\} + \mathbf{1}\big\{ N_{234}> N_1^\nu, \, N_1 \sim N_5 > N_1^\eta \big\} \Big) \label{linear:eq-lincub-decomp-w1}  \\
& \hspace{14ex} \times \LinCub^{(5)}[N_\ast](w) \bigg]  \notag \\
+& 3 \sum_{\substack{N_0,N_1,N_2,N_3, \\ N_4,N_5,N_{234} \leq N}}  \bigg[ 
\Big( 1 - \mathbf{1} \big\{ N_{234},N_5 \leq N_1^\eta \big\} - \mathbf{1} \big\{ N_{234}\leq N_1^\nu, N_1 \sim N_5 > N_1^\eta\big\} \Big) 
\LinCub^{\sin}[N_\ast](w) \bigg]. \label{linear:eq-lincub-decomp-w2}
\end{align}
Furthermore, for all functions $Y \colon \R \times \T^3 \rightarrow \R$, it holds that 
\begin{align}
&P_{\leq N}  \bigg[ \slinear[blue][\leqN] \scubic[\leqN] P_{\leq N} Y 
-\HLL \Big( \slinear[blue][\leqN], \scubic[\leqN], P_{\leq N} Y \Big) \bigg] \notag \\
=&\sum_{\substack{N_0,N_1,N_2,N_3, \\ N_4,N_5,N_{234} \leq N}} 
\mathbf{1}\big\{ \max(N_{234},N_5)> N_1^\eta \big\} 
\LinCub^{(5)}[N_\ast](Y) \label{linear:eq-lincub-decomp-Y1} \\
+& 3 \sum_{\substack{N_0,N_1,N_2,N_3, \\ N_4,N_5,N_{234} \leq N}}
\mathbf{1}\big\{ \max(N_{234},N_5)> N_1^\eta \big\}  \LinCub^{\sin}[N_\ast](Y) \label{linear:eq-lincub-decomp-Y2}. 
\end{align}
\end{lemma}

Our reason for stating the (equivalent) frequency-restrictions in \eqref{linear:eq-lincub-decomp-w1} and \eqref{linear:eq-lincub-decomp-w2} differently is that the one-pairing case does not require the subtraction of the high$\times$low$\times$low and high$\times$low$\times$high-interactions  and the three frequency-regions in  \eqref{linear:eq-lincub-decomp-w2} can be treated separately.  
Due to this, it might have been more natural to subtract the high$\times$low$\times$low and high$\times$low$\times$high-interactions only in the zero-pairing case, but this would increase the notational complexity of our Ansatz. 

\begin{proof}
Let $N_1,N_{234},N_5 \leq N$. Using the definitions of the linear and cubic stochastic objects (see Subsection \ref{section:diagram-wave}), it holds that 
\begin{equation}\label{linear:eq-lincub-decomp-p1}
\begin{aligned}
&P_{\leq N} \bigg[ \slinear[blue][N_1] P_{N_{234}} \scubic[\leqN] P_{N_5} w \bigg] \allowdisplaybreaks[3]\\
=& \sum_{n_0,n_1,\hdots,n_5 \in \Z^3} 
\sum_{ \substack{\varphi_1,\varphi_2,\varphi_3,\varphi_4 \in \\ \in \{ \cos, \sin \}} } \bigg[ \mathbf{1}\{n_0=n_{12345}\}  1_{N_1}(n_1) 1_{N_{234}}(n_{234}) 1_{N_5}(n_5)  \Big( \prod_{j=0,3,4} 1_{\leq N}(n_j) \Big) \\
&\times \langle n_{234} \rangle^{-1} \Big( \prod_{j=1}^4 \langle n_j \rangle^{-1} \Big) 
\varphi_1\big( t \langle n_1 \rangle\big) \Big( \int_0^t \dt^\prime \sin\big( (t-t^\prime) \langle n_{234} \rangle\big) \prod_{j=2}^4 \varphi_j\big(t^\prime \langle n_j \rangle \big) \Big) \\
&\times\SI[n_1,\varphi_1] \SI[n_j,\varphi_j \colon 2\leq j\leq 4] \, \widehat{w}(t,n_5) e^{i \langle n_0, x\rangle} \bigg]. 
\end{aligned}
\end{equation}
From the product-formula for multiple stochastic integrals, we obtain the decomposition 
\begin{align}
\eqref{linear:eq-lincub-decomp-p1} 
=& \sum_{n_0,n_1,\hdots,n_5 \in \Z^3} 
\sum_{ \substack{\varphi_1,\varphi_2,\varphi_3,\varphi_4 \in \\ \in \{ \cos, \sin \}} } \bigg[ \mathbf{1}\{n_0=n_{12345}\}  1_{N_1}(n_1) 1_{N_{234}}(n_{234}) 1_{N_5}(n_5) \Big( \prod_{j=0,3,4} 1_{\leq N}(n_j) \Big) \notag \\
&\times \langle n_{234} \rangle^{-1} \Big( \prod_{j=1}^4 \langle n_j \rangle^{-1} \Big) 
\varphi_1\big( t \langle n_1 \rangle\big) \Big( \int_0^t \dt^\prime \sin\big( (t-t^\prime) \langle n_{234} \rangle\big) \prod_{j=2}^4 \varphi_j\big(t^\prime \langle n_j \rangle \big) \Big) \label{linear:eq-lincub-decomp-p2}\\
&\times\SI[n_j,\varphi_j \colon 1\leq j\leq 4] \,  \widehat{w}(t,n_5) e^{i \langle n_0, x\rangle} \bigg] \notag  \allowdisplaybreaks[3]\\
+&3 \sum_{n_0,n_1,\hdots,n_5 \in \Z^3} 
\sum_{ \substack{\varphi_1,\varphi_2,\varphi_3,\varphi_4 \in \\ \in \{ \cos, \sin \}} } \bigg[ \mathbf{1}\{n_0=n_{345}\} \mathbf{1}\{ n_{12}=0 \}  1_{N_1}(n_1) 1_{N_{234}}(n_{234}) 1_{N_5}(n_5)
\Big( \prod_{j=0,3,4} 1_{\leq N}(n_j) \Big) \notag \\
&\times \langle n_{234} \rangle^{-1} \Big( \prod_{j=1}^4 \langle n_j \rangle^{-1} \Big) 
\mathbf{1}\{ \varphi_1 = \varphi_2 \} \, 
\varphi_1\big( t \langle n_1 \rangle\big) \Big( \int_0^t \dt^\prime \sin\big( (t-t^\prime) \langle n_{234} \rangle\big) \prod_{j=2}^4 \varphi_j\big(t^\prime \langle n_j \rangle \big) \Big)\label{linear:eq-lincub-decomp-p3} \\
&\times\SI[n_j,\varphi_j \colon 3\leq j\leq 4] \,  \widehat{w}(t,n_5) e^{i \langle n_0, x\rangle} \bigg]. \notag 
\end{align}
We first discuss the non-resonant component \eqref{linear:eq-lincub-decomp-p2}, which is treated more easily. After inserting a dyadic decomposition in $n_0,n_2,n_3$, and $n_4$, it directly from Definition \ref{linear:def-lincub} that  
\begin{equation*}
\eqref{linear:eq-lincub-decomp-p2} = \sum_{N_0,N_2,N_3,N_4\leq N} \LinCub^{(5)}[N_\ast](w). 
\end{equation*}
After summing over all relevant frequency-scales $N_1,N_{234}$, and $N_5$, this directly leads to \eqref{linear:eq-lincub-decomp-w1} and \eqref{linear:eq-lincub-decomp-Y1}. We now discuss the resonant component \eqref{linear:eq-lincub-decomp-p3}, which, due to the \sine-cancellation, requires a more careful analysis. We first focus on the sum in $n_1, n_2, \varphi_1$, and $\varphi_2$. Using trigonometric identities and $n_{12}=0$, it holds that 
\begin{equation*}
\sum_{\varphi_1,\varphi_2 \{ \cos, \sin \}} \mathbf{1}\{ \varphi_1 = \varphi_2 \} \varphi_1\big( t \langle n_1 \rangle \big) \varphi_2 \big( t^\prime \langle n_2 \rangle \big)
= \cos\big( (t-t^\prime) \langle n_1 \rangle\big). 
\end{equation*}
From the definition of the $\Sine$-kernel (Definition \ref{def:sine_kernel}), it follows that 
\begin{align*}
&\sum_{n_1,n_2 \in \Z^3} \bigg[ \mathbf{1}\big\{ n_{12}=0 \big\} 1_{N_1}(n_1) 1_{N_{234}}(n_{234}) 1_{\leq N}(n_2) \frac{\sin\big( (t-t^\prime) \langle n_{234} \rangle\big)}{\langle n_{234} \rangle} \frac{\cos\big( (t-t^\prime) \langle n_1\big)}{\langle n_1 \rangle^2} \bigg] \\
=& \sum_{\substack{n_1,n_{234} \in \Z^3 \colon \\ n_1 + n_{234} = n_{34}}} 
\bigg[ 1_{N_{234}}(n_{234}) 1_{N_1}(n_1) \frac{\sin\big( (t-t^\prime) \langle n_{234} \rangle\big)}{\langle n_{234} \rangle} \frac{\cos\big( (t-t^\prime) \langle n_1\big)}{\langle n_1 \rangle^2} \bigg] \\
=& \Sine[N_{234},N_1](t-t^\prime,n_{34}). 
\end{align*} 
From a notational perspective, it is convenient to introduce the new frequency-scale $N_2$ and always impose the restriction $N_2=N_1$. Then, we obtain that 
\begin{align*}
\eqref{linear:eq-lincub-decomp-p3} &= 3 \cdot \mathbf{1}\big\{ N_1 = N_2 \big\} 
\sum_{n_0, n_3, n_4, n_5 \in \Z^3} 
\sum_{\substack{\varphi_3, \varphi_4 \in \\ \{ \cos, \sin \}}} 
\Bigg[ \mathbf{1}\{n_0=n_{345}\}   \Big( \prod_{j=0,3,4} 1_{\leq N}(n_j) \Big) 1_{N_5}(n_5) \langle n_3 \rangle^{-1} \langle n_4 \rangle^{-1}\\
&\times \bigg( \int_0^t \dt^\prime  \Sine[N_{234},N_2](t-t^\prime,n_{34}) \prod_{j=3}^4 \varphi_j\big( t^\prime \langle n_j \rangle \big) \bigg) \, \widehat{w}(t,n_5)   \SI[n_j, \varphi_j \colon 3\leq j \leq 4] \, e^{i\langle n_0 ,x \rangle} 
\Bigg] \\
&= 3 \sum_{N_0,N_3,N_4\leq N} \LinCub^{\sin}[N_\ast,w]. 
\end{align*}
After summing over the relevant frequency-scales $N_1$, $N_{234}$, and $N_5$, this leads to the desired contributions in \eqref{linear:eq-lincub-decomp-w2} and \eqref{linear:eq-lincub-decomp-Y2}. 
\end{proof}

Equipped with Lemma \ref{linear:lem-lincub-decomp} and our tensor estimates from Section \ref{section:counting}, we now prove Proposition \ref{linear:prop-lincub}. 

\begin{proof}[Proof of Proposition \ref{linear:prop-lincub}:]
Using Lemma \ref{linear:lem-lincub-decomp}, we can decompose the linear random operators from \eqref{linear:eq-lincub-w} and \eqref{linear:eq-lincub-Y} into $\LinCub^{(5)}$ and $\LinCub^{\sin}$-terms. We first treat the $\LinCub^{(5)}$-terms, which, even equipped with our tensor estimates from Section \ref{section:counting}, still require detailed case distinctions.  In the end, we treat the $\LinCub^{\sin}$-terms, whose estimates require no further case distinctions. \\

\emph{The $\LinCub^{(5)}$-terms:} We first treat the random operator in the $w$-variable from \eqref{linear:eq-lincub-decomp-w1}. Using the reduction arguments in Subsection \ref{prep:remark-reduction}, the moment method (Proposition \ref{counting:prop-moment}), and the quintic tensor estimate (Lemma \ref{counting:quintic_tensor}), it follows that 
\begin{align}
 &\mathbb{E} \bigg[ \sup_{\Jc} \Big\| w \mapsto  \eqref{linear:eq-lincub-decomp-w1} \Big\|_{X^{1/2-\delta_1,b}(\Jc)\rightarrow X^{-1/2+\delta_2,b_+-1}(\Jc)}^p \bigg]^{1/p}  \notag \\[1ex]
&\lesssim p^2 T^\alpha \sum_{\substack{N_0,N_1,N_2,N_3, \\ N_4,N_5,N_{234} \leq N}} \bigg[ \Big( 
\mathbf{1}\big\{ N_{234}> N_1^\eta \geq N_{5} \big\}+
\mathbf{1}\big\{ N_1\not \sim N_5 > N_1^\eta \big\}  \notag \\&\hspace{19ex}+ 
\mathbf{1}\big\{ N_{234}> N_1^\nu, \, N_1 \sim N_5 > N_1^\eta \big\} \Big) \label{linear:eq-lincub-p1} \\
&\hspace{2ex}\times N_{\textup{max}}^{\epsilon} N_0^{\delta_2} \Big( \max(N_0,N_2,N_3,N_4)^{-1/2} + \max(N_2,N_3,N_4,N_5)^{-1/2} \Big) N_5^{\delta_1} \bigg]. \notag 
\end{align}
Thus, it only remains to control the dyadic sum in \eqref{linear:eq-lincub-p1}. In this argument, we distinguish the contributions of the three indicator functions in \eqref{linear:eq-lincub-p1}. 
The contribution of the first indicator function  in \eqref{linear:eq-lincub-p1} is estimated by 
\begin{align*}
&\sum_{\substack{N_0,N_1,N_2,N_3, \\ N_4,N_5,N_{234} \leq N}} \bigg[ 
\Big( \mathbf{1}\big\{ N_{234}> N_1^\eta \geq N_{5} \big\}  \\
&\hspace{13ex}\times N_{\textup{max}}^{\epsilon} N_0^{\delta_2} \Big( \max(N_0,N_2,N_3,N_4)^{-1/2} + \max(N_2,N_3,N_4,N_5)^{-1/2} \Big) N_5^{\delta_1} \bigg] \\
\lesssim& \sum_{\substack{N_0,N_1,N_2,N_3, \\ N_4,N_5,N_{234} \leq N}} \bigg[ 
\Big( \mathbf{1}\big\{ N_{234}> N_1^\eta \geq N_{5} \big\} N_{\max}^{\epsilon + \eta \delta_1 + \delta_2 } \max(N_2,N_3,N_4)^{-1/2} \bigg] \\
\lesssim& \sum_{\substack{N_0,N_1,N_2,N_3, \\ N_4,N_5,N_{234} \leq N}} N_{\textup{max}}^{\epsilon+\eta \delta_1 + \delta_2 - \eta/2} \lesssim 1. 
\end{align*}
We now turn to the contribution of the second indicator function in \eqref{linear:eq-lincub-p1}.  Since $N_1 \not \sim N_5$, it holds that $\max(N_0,N_2,N_3,N_4) \gtrsim N_{\textup{max}}$. Furthermore, since $N_5 > N_1^\eta$, it also holds that 
\begin{equation*}
\max(N_2,N_3,N_4,N_5)^{-1/2} N_5^{\delta_1} \lesssim \max(N_2,N_3,N_4,N_5)^{-1/2+\delta_1} \lesssim N_{\max}^{- ( 1/2-\delta_1) \eta}. 
\end{equation*}
As a result, the contribution of the second indicator function in \eqref{linear:eq-lincub-p1} is bounded by 
\begin{align*}
   &\sum_{\substack{N_0,N_1,N_2,N_3, \\ N_4,N_5,N_{234} \leq N}} \bigg[ 
\Big( \mathbf{1}\big\{ N_1\not \sim N_5 > N_1^\eta \big\} \\
&\hspace{13ex}\times N_{\textup{max}}^{\epsilon} N_0^{\delta_2} \Big( \max(N_0,N_2,N_3,N_4)^{-1/2} + \max(N_2,N_3,N_4,N_5)^{-1/2} \Big) N_5^{\delta_1} \bigg] \\
\lesssim&\sum_{\substack{N_0,N_1,N_2,N_3, \\ N_4,N_5,N_{234} \leq N}} 
\Big( N_{\textup{max}}^{\epsilon+\delta_1+\delta_2-1/2} + N_{\textup{max}}^{\epsilon+\delta_2 - (1/2-\delta_1)\eta} \Big) 
\lesssim \,  1. 
\end{align*}
Finally, we treat the contribution of the third indicator function in \eqref{linear:eq-lincub-p1}. 
In this frequency-regime, it holds that 
\begin{align*}
   &\sum_{\substack{N_0,N_1,N_2,N_3, \\ N_4,N_5,N_{234} \leq N}} \bigg[ 
\Big( \mathbf{1}\big\{ N_{234}> N_1^\nu, \, N_1 \sim N_5 > N_1^\eta \big\}  \\
&\hspace{13ex}\times N_{\textup{max}}^{\epsilon} N_0^{\delta_2} \Big( \max(N_0,N_2,N_3,N_4)^{-1/2} + \max(N_2,N_3,N_4,N_5)^{-1/2} \Big) N_5^{\delta_1} \bigg] \\
\lesssim& \sum_{\substack{N_0,N_1,N_2,N_3, \\ N_4,N_5,N_{234} \leq N}} \bigg[ 
\mathbf{1}\big\{ N_{234}> N_1^\nu, \, N_1 \sim N_5 > N_1^\eta \big\} 
N_{\textup{max}}^{\epsilon+\delta_1+\delta_2} \max(N_2,N_3,N_4)^{-1/2}  \bigg]\\
\lesssim& \sum_{\substack{N_0,N_1,N_2,N_3, \\ N_4,N_5,N_{234} \leq N}} 
N_{\textup{max}}^{\epsilon+\delta_1+\delta_2-\nu/2} \lesssim 1.  
\end{align*}
This completes our estimate of the $w$-operator in \eqref{linear:eq-lincub-decomp-w1} and we now turn to the $Y$-operator from \eqref{linear:eq-lincub-decomp-Y1}. Arguing as for \eqref{linear:eq-lincub-p1}, we obtain that 
\begin{align}
 &\mathbb{E} \bigg[ \sup_{\Jc} \Big\| Y \mapsto  \eqref{linear:eq-lincub-decomp-Y1} \Big\|_{X^{1/2+\delta_2,b}(\Jc) \rightarrow X^{-1/2+\delta_2,b_+-1}(\Jc)}^p \bigg]^{1/p}  \notag \\[1ex]
&\lesssim p^2 T^\alpha  \sum_{\substack{N_0,N_1,N_2,N_3, \\ N_4,N_5,N_{234} \leq N}} \bigg[ 
\mathbf{1}\big\{ \max(N_{234},N_5)> N_1^\eta  \big\} \label{linear:eq-lincub-p2} \\
&\hspace{2ex}\times N_{\textup{max}}^{\epsilon} N_0^{\delta_2} \Big( \max(N_0,N_2,N_3,N_4)^{-1/2} + \max(N_2,N_3,N_4,N_5)^{-1/2} \Big) N_5^{-\delta_2} \bigg]. \notag 
\end{align}
While \eqref{linear:eq-lincub-p2} contains fewer frequency-restrictions than \eqref{linear:eq-lincub-p1}, the $N_5^{\delta_1}$-factor has been replaced by $N_5^{-\delta_2}$, which makes a significant difference. Using $\max(N_{234},N_5)>N_1^\eta$, it holds that 
\begin{equation*}
\max(N_2,N_3,N_4,N_5) \gtrsim \max(N_1,N_2,N_3,N_4,N_5)^\eta \gtrsim N_{\textup{max}}^{\eta}. 
\end{equation*}
As a result, we obtain that 
\begin{align*}
\eqref{linear:eq-lincub-p2} 
&\lesssim p^2 T^\alpha 
\sum_{\substack{N_0,N_1,N_2,N_3, \\ N_4,N_5,N_{234} \leq N}} \bigg[ 
\mathbf{1}\big\{ \max(N_{234},N_5)> N_1^\eta  \big\}\\
&\hspace{8ex}\times N_{\textup{max}}^{\epsilon}   \Big( \max(N_0,N_2,N_3,N_4)^{-1/2+\delta_2} N_5^{-\delta_2} 
+ \max(N_2,N_3,N_4,N_5)^{-1/2} \Big) \bigg] \\
&\lesssim p^2 T^\alpha  \sum_{\substack{N_0,N_1,N_2,N_3, \\ N_4,N_5,N_{234} \leq N}}    N_{\textup{max}}^{\epsilon} \Big( N_{\textup{max}}^{-\delta_2} + N_{\textup{max}}^{-\eta/2} \Big)
\lesssim p^2 T^\alpha . 
\end{align*}

\emph{The $\LinCub^{\sin}$-terms:} 
As in our estimate of $\LinCub^{(5)}$, we first use the reduction arguments in Subsection \ref{prep:remark-reduction}, then the moment method (Proposition \ref{counting:prop-moment}), and finally our sine-cancellation tensor estimates (Lemma \ref{counting:sine_tensor}).  Due to the  $N_{\textup{max}}^{-1/2}$-factor in Lemma \ref{counting:sine_tensor}, the remaining sum over all frequency-scales can easily be performed.
\end{proof} 
\section{Para-controlled calculus}\label{section:para}

In this section, we analyze nonlinear terms containing the para-controlled components $\SXXone$ and $\SXXtwo$. The estimates in the following three subsections are used to estimate nonlinear terms which, in addition to the para-controlled components, contain zero, one, or two factors of $\slinear[blue]$, respectively. 

In Subsection \ref{section:para-regularity}, we prove probabilistic Strichartz and regularity estimates for $\SXXj$. Despite their simplicity, these estimates will be sufficient to control nonlinear terms involving no factors of $\slinear[blue]$. 
In Subsection \ref{section:para-one-linear}, we estimate  products of $\slinear[blue]$ with one of the para-controlled components $\SXXone$ or $\SXXtwo$. Together with previous estimates from Section \ref{section:analytic}, Section \ref{section:bilinear}, and Section \ref{section:linear}, the product estimates control all terms involving exactly one factor of $\slinear[blue]$ and one para-controlled component.
Finally, in Subsection \ref{section:para-two-linear}, we treat the two random operators from \eqref{ansatz:eq-YN-3} and \eqref{ansatz:eq-YN-4}, which contain the quadratic stochastic object, i.e., two factors of $\slinear[blue]$. This part of our argument most heavily relies on the para-controlled structure of $\SXXone$ and $\SXXtwo$ and also requires the tensor estimates from Section \ref{section:counting}. 

Before starting with our estimates, we introduce additional notation. 
We recall from Definition \ref{ansatz:def-operators-para} that both $\XXone$ and $\XXtwo$ themselves consist of multiple different terms (see \eqref{ansatz:eq-X11-a}-\eqref{ansatz:eq-X11-d} and \eqref{ansatz:eq-X21}). Since some of these terms require different arguments, we isolate them using the following definition.

\begin{definition}[Decomposition of $\XXone$ and $\XXtwo$]\label{para:def-decomposition}~
Let $N$ and $N_0,N_1,N_2,N_3,N_{23}\leq N$ be frequency-scales.  Then, we define the following frequency-localized terms and operators:
\begin{enumerate}[label=(\roman*)] 
\item The high$\times$low$\times$low-portion of $\SXXone$ is defined as
\begin{align*}
 \XXhll[N_\ast,w_2,w_3] 
                := 1\big\{ N_2,N_3 \leq N_1^\eta\big\} P_{N_0}\Duh \Big[ \, \slinear[blue][N_1] P_{N_2} w_2 P_{N_3}w_3 \Big].
\end{align*} 
\item The resonant-portion of $\SXXone$ is defined as
\begin{equation*}
    \XXres[ N_\ast,w_2,w_3 ] 
                    := 1\big\{ N_{3} > \max(N_1,N_2)^\eta\big\} 
                    1\big\{ N_{23} \leq N_1^\eta \big\} P_{N_0} \Duh \Big[ \slinear[blue][N_1] P_{N_{23}} \Big( P_{N_2} w_2 P_{N_3} w_3 \Big) \Big].
\end{equation*}
\item The explicit portion of $\SXXone$ is defined as 
\begin{align*}
    \XXexpl &:= \Duh \bigg[ \Atone \big( \slinear[blue] \big) \big( \gamma_{\leq N} - \Gamma_{\leq N} \big) \, \slinear[blue][\leqN] +  \Aone\big( \slinear[blue] \big)  \mathfrak{C}_{\leq N}  \, \slinear[blue][\leqN] \\
    &+18 \sum_{\substack{N_0,N_1,N_2,N_3 \colon\\  N_{\textup{max}}\leq N, \\ \max(N_2,N_3) \leq N_1^\eta}} P_{N_0} \Big( \slinear[blue][N_1] \mathfrak{C}_{\leq N}^{(1,5)}[N_2,N_3](t) \Big) \\
    &- \Athree\Big( \slinear[blue], \scubic, \scubic \Big) \sum_{\substack{N_0,N_1,N_2,N_3 \colon\\  N_{\textup{max}}\leq N, \\ \max(N_2,N_3) \leq N_1^\eta}} P_{N_0} \Big( \slinear[blue][N_1] \mathfrak{C}_{\leq N}^{(3,3)}[N_2,N_3](t) \Big) \bigg],
\end{align*}
where $\Atone$, $\Aone$, and $\Athree$ are as in Lemma \ref{ansatz:lem-grouping}.
\item The operator version of $\SXXtwo$ is defined as 
\begin{align*}
\XXop[N_\ast,w] := - 3 P_{N_0} \Duh \Big[ \lcol \slinear[blue][N_1] \slinear[blue][N_2] \rcol P_{N_3} w \Big]. 
\end{align*}
\end{enumerate}
\end{definition}

In the following remark, we emphasize certain aspects of Definition \ref{para:def-decomposition}.
\begin{remark}
 In \eqref{ansatz:eq-X11-d}, we encountered $\RES\big(\slinear[blue][\leqN],\slinear[blue][\leqN], Y_{\leq N}\big)$, which will be estimated using the operator $\XXres[N_\ast,w_2,w_3]$ with $(w_2,w_3)=(\,\slinear[blue],Y)$.  The reason for treating this term as a bilinear operator in $w_2$ and $w_3$, rather than as a linear operator in $Y$, is as follows: Since $Y$ is viewed as a general element of $X^{1/2+\delta_2,b}$, it can behave like $\langle \nabla \rangle^{-1-\delta_2-\epsilon} \,  \slinear[blue]$. In this case, however, the high$\times$high$\rightarrow$low-interactions of $\slinear[blue]$ and $Y$ are no better than high$\times$high$\rightarrow$low-interactions of general elements $(w_2,w_3) \in X^{-1/2-\epsilon,b} \times X^{1/2+\delta_2,b}$.

We also emphasize that $\XXop$ contains the precise pre-factor $-3$, which is essential in Subsection \ref{section:para-two-linear} below, where a resonance originating from $\XXop$ has to cancel with the renormalization multiplier $\Gamma$. In all other terms, the precise pre-factors are not important.  
\end{remark}

Using Definition \ref{para:def-decomposition}, we can decompose $\XXone$ and $\XXtwo$ as follows: 

\begin{lemma}[Decomposition of $\XXone$ and $\XXtwo$]\label{para:lem-decomposition}
For all $N\geq 1$, it holds that 
\begin{align}
\XXone[v_{\leq N},Y_{\leq N}] 
&= - 6 \sum_{\substack{N_0,N_1,N_2,N_3 \colon \\ N_{\textup{max}} \leq N}}
\XXhll\Big[ N_\ast, \slinear[blue][\leqN], - P_{\leq N} \slinear[green][\leqM] + 3 \squintic[\leqN] + v_{\leq N} \Big] \allowdisplaybreaks[3]\\
&+ \sum_{\zeta^{(2)},\zeta^{(3)} \in \Symb_0^b} 
\sum_{\substack{N_0,N_1,N_2,N_3 \colon \\ N_{\textup{max}} \leq N}}
\Athree\Big( \slinear[blue], \zeta^{(2)}, \zeta^{(3)} \Big) \XXhll\big[ N_\ast, \zeta^{(2)}_{\leq N}, \zeta^{(3)}_{\leq N}\big] \allowdisplaybreaks[3] \\
&- 6  \sum_{\substack{N_0,N_1,N_2,N_3,N_{23} \colon \\ N_{\textup{max}} \leq N}} \XXres\big[N_\ast, \slinear[blue][\leqN], Y_{\leq N}\big] \allowdisplaybreaks[3] \label{para:eq-X11-decomp-c} \\
&+ \XXexpl. 
\end{align}
Furthermore, it holds that 
\begin{align*}
\XXtwo\big[v_{\leq N}\big] = \sum_{\substack{N_0,N_1,N_2,N_3 \colon \\ N_{\textup{max}} \leq N}} \XXop\Big[ N_\ast, - P_{\leq N} \slinear[green][\leqM] + 3 \squintic[\leqN] + v_{\leq N} \Big]. 
\end{align*}
\end{lemma}

\begin{proof}
The two identities follow directly from our definitions, i.e., Definition \ref{ansatz:def-para-product}, Definition \ref{ansatz:def-operators-para}, and Definition \ref{para:def-decomposition}. To avoid any possible confusion, we only emphasize the following two aspects:

First, the pre-factors in \eqref{ansatz:eq-X11-d} and the corresponding contribution \eqref{para:eq-X11-decomp-c} differ by a factor of two. The reason is that $\RES$ is the symmetric version of the interaction in the Duhamel integral of $\XXres$. 

Second, the $\mathfrak{C}_{\leq N}^{(1,5)}[N_2,N_3]$ and $\mathfrak{C}_{\leq N}^{(3,3)}[N_2,N_3]$-terms are included in $\XXexpl$ because $\XXhll$ and $\XXop$ are defined using the usual and not the modified product (from Definition \ref{ansatz:def-product}). 
\end{proof} 

Due to Lemma \ref{para:lem-decomposition}, the rest of this section focuses on the frequency-localized operators from  Definition \ref{para:def-decomposition}.

\subsection{Probabilistic Strichartz and regularity estimates}\label{section:para-regularity}

In this section, we control the $L_t^\infty \C_x^{1/2-\delta_1}$ and $X^{1/2-\delta_1,b}$-norms of $\SXXone$ and $\SXXtwo$. This eventually leads to a proof of Proposition \ref{ansatz:prop-Xj}, which is presented in Section \ref{section:proof-main-estimates} below. We start with our estimates for $\SXXone$. 

\begin{lemma}[Probabilistic Strichartz and regularity estimates for $\SXXone$]\label{para:lem-strichartz-XXone}~ Let $T\geq 1$ and $p\geq 2$. In the following, all suprema over $\mathcal{J}$ are taken over closed intervals $0 \in \mathcal{J} \subseteq [-T,T]$. 
\begin{enumerate}[label=(\roman*)]
    \item\label{para:item-strichartz-1} For all frequency-scales  $N_0,N_1,N_2,N_3$ satisfying $N_2,N_3 \leq N_1^\eta$, it holds that 
    \begin{align*}
    &\mathbb{E} \bigg[ \, \sup_{\mathcal{J}} |\mathcal{J}|^{-(b_+-b)p} \Big\| (w_2,w_3) \mapsto \XXhll[N_\ast,w_2,w_3] \Big\|_{X^{-1,b}(\mathcal{J})\times X^{-1,b}(\mathcal{J})\rightarrow X^{1/2-\delta_1,b}(\mathcal{J})}^p \bigg]^{1/p} \\
    &\lesssim p^{1/2} T^{\alpha} N_{\textup{max}}^{-\epsilon}. 
    \end{align*}
    Furthermore, it holds that 
    \begin{align*}
    &\mathbb{E} \bigg[ \sup_{\mathcal{J}} |\mathcal{J}|^{-(b_+-b)p} \Big\| (w_2,w_3) \mapsto \XXhll[N_\ast,w_2,w_3] \Big\|_{X^{-1,b}(\mathcal{J})\times X^{-1,b}(\mathcal{J})\rightarrow L_t^\infty \C_x^{1/2-\delta_1}(\mathcal{J})}^p  \bigg]^{1/p} \\
    &\lesssim p^{1/2}T^{\alpha} N_{\textup{max}}^{-\epsilon}. 
    \end{align*}
    \item \label{para:item-strichartz-2} For all frequency-scales  $N_0,N_1,N_2,N_3,N_{23}$ satisfying $N_3 > \max(N_1,N_2)^\eta$ and $N_{23}\leq N_1^\eta$, it holds that 
    \begin{align*}
   &\mathbb{E} \bigg[ \sup_{\mathcal{J}} |\mathcal{J}|^{-(b_+-b)p} \Big\| (w_2,w_3) \mapsto \XXres[N_\ast,w_2,w_3] \Big\|_{X^{-1/2-\epsilon,b}(\mathcal{J})\times X^{1/2+\delta_2,b}(\mathcal{J})\rightarrow X^{1/2-\delta_1,b}(\mathcal{J})}^p \bigg]^{1/p} \\ &\lesssim p^{1/2} T^{\alpha}N_{\textup{max}}^{-\epsilon}. 
    \end{align*}
    Furthermore, it holds that 
    \begin{align*}
    &\mathbb{E} \bigg[ \sup_{\mathcal{J}} |\mathcal{J}|^{-(b_+-b)p} \Big\| (w_2,w_3) \mapsto \XXres[N_\ast,w_2,w_3] \Big\|_{X^{-1/2-\epsilon,b}(\mathcal{J})\times X^{1/2+\delta_2,b}(\mathcal{J})\rightarrow L_t^\infty \C_x^{1/2-\delta_1}(\mathcal{J})}^p \bigg]^{1/p} \\ &\lesssim p^{1/2} T^{\alpha} N_{\textup{max}}^{-\epsilon}. 
    \end{align*}
    \item \label{para:item-strichartz-3} It holds that 
    \begin{align*}
    \E \bigg[ \sup_N \Big\| \XXexpl \Big\|_{(X^{1/2-\delta_1,b} \cap L_t^\infty \C_x^{1/2-\delta_1})([-T,T])}^p \bigg]^{1/p} 
    \lesssim p^{1/2} T^\alpha. 
    \end{align*}
\end{enumerate}
\end{lemma}

\begin{proof}
We only treat the case $T=1$, since the general case follows from minor modifications. We prove the  estimates for $\XXhll$, $\XXres$, and $\XXexpl$ separately. \\

\emph{Proof of \ref{para:item-strichartz-1}:} Using energy estimates, it holds that 
\begin{equation}\label{para:eq-strichartz-q1}
\begin{aligned}
\| \XXhll [ N_\ast, w_2,w_3] \|_{X^{1/2-\delta_1,b}(\mathcal{J})} 
&\leq N_1^{-1/2-\delta_1} \Big\| \slinear[blue][N_1] P_{N_2} w_2 P_{N_3} w_3 \Big\|_{L_t^2 L_x^2(\mathcal{J})} \\
&\lesssim |\mathcal{J}|^{1/2} N_1^{-1/2-\delta_1} \Big\| \slinear[blue][N_1] \Big\|_{L_t^\infty L_x^\infty(\mathcal{J})} 
\prod_{j=2}^3 \| P_{N_j} w_j \|_{L_t^\infty L_x^\infty(\mathcal{J})} \\
&\lesssim |\mathcal{J}|^{1/2} N_1^{\epsilon-\delta_1} N_2^{5/2} N_3^{5/2} \| \slinear[blue] \|_{L_t^\infty \C_x^{-1/2-\epsilon}(\mathcal{J})} \prod_{j=2}^3 \| w_j \|_{X^{-1,b}(\mathcal{J})}. 
\end{aligned}
\end{equation}
Due to the high$\times$low$\times$low-assumption, it holds that 
\begin{equation*}
    N_1^{\epsilon-\delta_1} N_2^{5/2} N_3^{5/2}  \lesssim N_1^{\epsilon + 5 \eta - \delta_1} \lesssim N_{\textup{max}}^{-\delta_1/2}. 
\end{equation*}
This completes our proof of the $X^{1/2-\delta_1,b}$-estimate and it remains to prove the $L_t^\infty \C_x^{1/2-\delta_1}$-estimate. 
Using the reduction arguments in Subsection \ref{prep:remark-reduction} and $\Jc\subseteq [-1,1]$, it suffices to prove that 
\begin{align*}
\sup_{\substack{n_2,n_3 \in \Z^3 \colon \\ |n_2|,|n_3|\leq N_1^\eta}} \sup_{t\in [-1,1]} \sup_{\lambda_2, \lambda_3 \in \R} \E \Big[  
\Big\| P_{N_0} \Duh \Big[ \slinear[blue][N_1] \prod_{j=2}^3 e^{i\langle n_j , x \rangle \pm_j i \langle n_j \rangle t+ i t \lambda_j }\Big]
\Big\|_{ \C_x^{1/2-\delta_1}}^2 \Big]^{1/2} \lesssim N_1^{-\epsilon}. 
\end{align*}
This estimate can be obtained exactly as in \eqref{para:eq-strichartz-q1}. \\

\emph{Proof of \ref{para:item-strichartz-2}:} The argument is similar to the proof of  \ref{para:item-strichartz-1} and we omit the details. The only additional ingredients are that the sum of the regularities of $w_2$ and $w_3$ is positive, i.e., $(-1/2-\epsilon)+(1/2+\delta_2)>0$, and a standard box localization argument (see e.g. \cite[Corollary 3.13]{Tao01}).  \\

\emph{Proof of \ref{para:item-strichartz-3}:} Using our estimate for $\mathfrak{C}_{\leq N}$ (Lemma \ref{diagrams:lem-mathfrakC}), $\mathfrak{C}_{\leq N}^{(1,5)}[N_2,N_3]$ and $\mathfrak{C}_{\leq N}^{(3,3)}[N_2,N_3]$ (Lemma \ref{analytic:lem-C15-C33}), and $\gamma_{\leq N}-\Gamma_{\leq N}$ (Lemma \ref{analytic:lem-difference-gamma}), the argument is again essentially as in \ref{para:item-strichartz-1}. We note that, due to translation invariance and Gaussian hypercontractivity, the estimate of $(\gamma_{\leq N}-\Gamma_{\leq N})\slinear[blue][\leqN]$ only requires operator bounds on $L_x^2$ and not on $L_x^\infty$. 
\end{proof}

\begin{lemma}[Probabilistic Strichartz and regularity estimates for $\SXXtwo$]\label{para:lem-strichartz-XXtwo}
Let $T\geq 1$ and $p\geq 2$. For all frequency-scales $N_0,N_1,N_2,N_3$ satisfying $N_3 \leq \max(N_1,N_2)^\eta < \min(N_1,N_2)$, it holds that 
\begin{equation}\label{para:eq-XXtwo-regularity}
\mathbb{E} \bigg[ \sup_{\mathcal{J}} |\mathcal{J}|^{-(b_+-b)p} \Big\| w \mapsto \XXop[N_\ast,w] \Big\|_{X^{-1,b}(\mathcal{J})\rightarrow X^{1/2-\delta_1,b}(\mathcal{J})}^p \bigg]^{1/p} \lesssim p T^\alpha N_{\textup{max}}^{-\epsilon},
\end{equation}
where the supremum is taken over all closed intervals $0 \in \mathcal{J}\subseteq [-T,T]$. 
Furthermore, it holds that 
\begin{equation}\label{para:eq-XXtwo-strichartz}
\mathbb{E} \bigg[ \sup_{\mathcal{J}} |\mathcal{J}|^{-(b_+-b) p } \Big\| w \mapsto \XXop[N_\ast,w] \Big\|_{X^{-1,b}(\mathcal{J})\rightarrow L_t^\infty \C_x^{1/2-\delta_1}(\mathcal{J})}^p \bigg]^{1/p} \lesssim p T^\alpha N_{\textup{max}}^{-\epsilon}.
\end{equation}
\end{lemma}

In contrast to Lemma \ref{para:lem-strichartz-XXone}, Lemma \ref{para:lem-strichartz-XXtwo} relies on multilinear dispersive effects. As a result, we can no longer (exclusively) rely on Strichartz-type estimates, but instead use our tensor estimates. 

\begin{proof} As before, we only treat the case $T=1$, since the general case follows from minor modifications. 
We first prove the regularity estimate \eqref{para:eq-XXtwo-regularity}. The argument has similarities with the proof of Lemma \ref{linear:lem-quad-dyadic}, but concerns a different frequency-regime. Using the reduction arguments in Subsection \ref{prep:remark-reduction}, it suffices to prove the random tensor estimate 
\begin{equation}\label{para:eq-strichartz-XXone-p1}
N_0^{-1/2-\delta_1} N_1^{-1}N_2^{-1} N_3 \, 
\mathbb{E} \bigg[ \Big\| \widehat{w} \mapsto \sum_{n_1,n_2} h_{n_0n_1 n_2 n_3} \lcol g_{n_1} g_{n_2} \rcol \widehat{w}(n_3) \Big\|_{\ell_{n_3}^2 \rightarrow \ell_{n_0}^2}^p \bigg]^{1/p} \lesssim p N_{\textup{max}}^{-\delta_1+4\eta}.
\end{equation}
Here, $h$ is the base tensor from Section \ref{section:counting}, i.e., 
\begin{align*}
h_{n_0 n_1 n_2 n_3} &= \Big( \prod_{j=0}^3 1_{N_j}(n_j) \Big) 1\{n_0 =n_{123}\} 1\{ |\Omega-m| \leq 1 \}, \\
\Omega &= \sum_{j=0}^3 (\pm_j) \langle n_j \rangle, \quad \text{and} \quad m \in \Z. 
\end{align*}
In fact, \eqref{para:eq-strichartz-XXone-p1} yields slightly stronger decay in $N_{\textup{max}}$ than our desired estimate \eqref{para:eq-XXtwo-regularity}. While the random tensor estimate could be proven using the moment method (Proposition \ref{counting:prop-moment}), it also follows from the following more elementary argument. 

Using Cauchy-Schwarz in the $n_3$-variable, it holds that 
\begin{align*}
    \mathbb{E} \bigg[ \Big\| \widehat{w} \mapsto \sum_{n_1,n_2} h_{n_0n_1 n_2 n_3} \lcol g_{n_1} g_{n_2} \rcol \widehat{w}(n_3) \Big\|_{\ell_{n_3}^2 \rightarrow \ell_{n_0}^2}^p \bigg]^{1/p} 
    &\lesssim \mathbb{E} \bigg[ \Big\|  \sum_{n_1,n_2} h_{n_0n_1 n_2 n_3} \lcol g_{n_1} g_{n_2}  \rcol  \Big\|_{n_0 n_3}^{p}  \bigg]^{1/p}  \\
    &\lesssim p \| h \|_{n_0 n_1 n_2 n_3}. 
\end{align*}
Using the Hilbert-Schmidt estimate from Lemma \ref{counting:basetensor_est}, it follows that 
\begin{align*}
&N_0^{-1/2-\delta_1} N_1^{-1}N_2^{-1} N_3 \| h \|_{n_0 n_1 n_2 n_3} \\
\lesssim& \, N_0^{-1/2-\delta_1} N_1^{-1}N_2^{-1} N_3 N_{\textup{min}}^{1/2} N_{\textup{max}}^{-1/2} N_0 N_1 N_2 N_3 \\
\lesssim&\, N_{\textup{max}}^{-\delta_1} N_3^{5/2} 
\lesssim \, N_{\textup{max}}^{-\delta_1+3\eta}. 
\end{align*}

This completes the proof of \eqref{para:eq-strichartz-XXone-p1} and therefore the proof of the regularity estimate \eqref{para:eq-XXtwo-regularity}. Using the proof of \eqref{para:eq-strichartz-XXone-p1}, the  Strichartz estimate \eqref{para:eq-XXtwo-strichartz} can now be obtained  exactly as in Lemma \ref{para:lem-strichartz-XXone} and we omit the details. 
\end{proof}

\subsection{Interactions with one linear stochastic object}\label{section:para-one-linear}

In the previous subsection, we proved that $\SXXone$ and $\SXXtwo$ have spatial regularity $1/2-\delta_1$, which is slightly less than $1/2$. From the discussion in Subsection \ref{section:ansatz-para}, it is clear that this cannot be improved. In particular, the products
\begin{equation}\label{para:eq-products-motivation}
\slinear[blue][\leqN] \XXone \qquad \text{and} \qquad \slinear[blue][\leqN] \XXtwo
\end{equation}
cannot be defined using only the regularity estimates of $\XXone$ and $\XXtwo$. As we will see in this subsection, however, the products in  \eqref{para:eq-products-motivation} can be defined using the para-controlled structure of $\XXone$ and $\XXtwo$.  

\begin{proposition}[\protect{Product estimates for $\slinear[blue] \, \SXXone$ and $\slinear[blue] \, \SXXtwo$}]\label{para:prop-resonant}
For all $A\geq 1$, there exists an $A$-certain event $E_A \in \mathcal{E}$ on which the following estimates hold: For all frequency-scales $N$, $K_1$, and $K_2$, all $T\geq 1$, closed intervals $0 \in \mathcal{J} \subseteq [-T,T]$, and $v_{\leq N},Y_{\leq N} \colon \mathcal{J} \times \T^3 \rightarrow \mathbb{R}$, it holds that 
\begin{equation}
\begin{aligned}
&\Big\| \slinear[blue][K_1] P_{K_2} \XXone[v_{\leq N}, Y_{\leq N}] \Big\|_{L_t^\infty \C_x^{-1/2-\epsilon}(\mathcal{J})} 
+ \Big\| \slinear[blue][K_1] P_{K_2} \XXtwo[v_{\leq N}] \Big\|_{L_t^\infty \C_x^{-1/2-\epsilon}(\mathcal{J})} \\
\leq& \, A T^\alpha K_2^{-1/2+10\eta} \Big( 1+ \| v_{\leq N}\|_{X^{-1,b}(\mathcal{J})}^2 + \| Y_{\leq N} \|_{X^{1/2+\delta_2,b}(\mathcal{J})} \Big). 
\end{aligned}
\end{equation}
\end{proposition}

We chose to denote the frequency-scales in Proposition \ref{para:prop-resonant} by $K_1$ and $K_2$ in order to not conflict with the frequency-scales in Definition \ref{para:def-decomposition}. 
The main ingredient in the proof of Proposition \ref{para:prop-resonant} is the \sine-cancellation (see e.g. Lemma \ref{counting:lem-Sine-symmetrization}). We emphasize that, in contrast to the product $\squadratic[\leqN] \XXtwo$ from the next subsection, the products in \eqref{para:eq-products-motivation} do not require a renormalization. 
We start by proving the product estimates for the different $\SXXone$-components. 

\begin{lemma}[\protect{Product estimates for $\slinear[blue] \, \SXXone$}]\label{para:lem-resonant-X1}
Let $T\geq 1$ and $p\geq 2$. In the following, all suprema over $\mathcal{J}$ are taken over closed intervals $0\in \mathcal{J} \subseteq [-T,T]$. 
\begin{enumerate}[label=(\roman*)]
\item\label{para:item-resonant-X1-hhl}
For all frequency-scales $N_1,N_2,N_3,N_4,N_{234}$ satisfying $\max(N_3,N_4)\leq N_2^\eta$, it holds that 
\begin{align*}
&\mathbb{E} \bigg[ \sup_{\mathcal{J}} \Big\| (w_3,w_4) \mapsto 
\slinear[blue][N_1] \XXhll[N_\ast;w_3,w_4] \Big\|_{X^{-1,b}(\mathcal{J}) \times X^{-1,b}(\mathcal{J}) \rightarrow L_t^\infty C_x^{-1/2-\epsilon}(\mathcal{J})}^p \bigg]^{1/p} \\
\lesssim& \, p T^\alpha N_{234}^{-1/2+10 \eta}. 
\end{align*}
\item\label{para:item-resonant-X1-res} For all frequency-scales $N_0,N_1,N_2,N_3,N_4,N_{34},N_{234}$ satisfying $N_4 > \max(N_2,N_3)^\eta$ and $N_{34}\leq N_2^\eta$, it follows that
\begin{align*}
&\mathbb{E} \bigg[ \sup_{\mathcal{J}} \Big\| (w_3,w_4) \mapsto 
\slinear[blue][N_1] \XXres[N_\ast;w_3,w_4] \Big\|_{X^{-1/2-\epsilon,b}(\mathcal{J}) \times X^{1/2+\delta_2,b}(\mathcal{J}) \rightarrow L_t^\infty C_x^{-1/2-\epsilon}(\mathcal{J})}^p \bigg]^{1/p} \\
\lesssim& \, p T^\alpha N_{234}^{-1/2+10 \eta}. 
\end{align*}
\item\label{para:item-resonant-X1-expl} 
For all frequency-scales $K_1$ and $K_2$, it holds that 
\begin{equation*}
\E \bigg[ \sup_N \Big\|\, \slinear[blue][K_1] P_{K_2} \XXexpl \Big\|_{L_t^\infty \C_x^{-1/2-\epsilon}([-T,T])}^p \bigg]^{1/p} \lesssim p T^\alpha K_2^{-1/2+10\eta}. 
\end{equation*}
\end{enumerate}
\end{lemma}

The choice of using $K_1$ and $K_2$ as frequency-scales in \ref{para:item-resonant-X1-expl} is deliberate, since, depending on the precise term in $\XXexpl$, $K_2$ is most naturally replaced by either $N_{234}$ or $N_2$.  

\begin{proof}
We only treat the case $T=1$ and $\mathcal{J}=[-1,1]$, since the general case follows from minor modifications. We prove the estimates for $\XXhll$, $\XXres$, and $\XXexpl$ separately.\\

\emph{Proof of \ref{para:item-resonant-X1-hhl}:} Using similar reduction arguments as in Subsection \ref{prep:remark-reduction}, it suffices to prove that 
\begin{equation}\label{para:eq-resonant-X1-p1}
\begin{aligned}
    &\sup_{t\in [-1,1]} \sup_{\pm_3,\pm_4} \sup_{\substack{n_3,n_4 \in \Z^3 \colon \\ |n_3|,|n_4| \leq N_2^\eta}} 
    \sup_{\lambda_3,\lambda_4 \in \R} \E \Bigg[ \bigg\| \sum_{n_1,n_2\in \Z^3} \bigg[ 1_{N_{234}}(n_{234}) \Big( \prod_{j=1}^2 1_{N_j}(n_j) \Big) e^{i\langle n_{1234}, x \rangle} \slinear[blue][N_1](t,n_1) \\
    &\times \Big( \int_0^t \dt^\prime 
    \frac{\sin\big( (t-t^\prime) \langle n_{234} \rangle\big)}{\langle n_{234} \rangle} \slinear[blue][N_2](t^\prime,n_2) 
    \Big( \prod_{j=3}^4 e^{i (\pm_j \langle n_j \rangle + \lambda_j ) t^\prime } \Big) 
    \Big) \bigg] \bigg\|_{H_x^{-1/2-\epsilon/2}}^2 \bigg] \lesssim N_1^{-\epsilon} N_{234}^{-1}. 
\end{aligned}
\end{equation}
In order to prove \eqref{para:eq-resonant-X1-p1}, we separate the non-resonant and resonant parts. Using the definition of the $\Sine$-kernel (Definition \ref{def:sine_kernel}), it holds that
\begin{align}
&\sum_{n_1,n_2\in \Z^3} \bigg[ 1_{N_{234}}(n_{234}) \Big( \prod_{j=1}^2 1_{N_j}(n_j) \Big) e^{i\langle n_{1234}, x \rangle} \slinear[blue][N_1](t,n_1) \notag  \\
    &\hspace{10ex}\times \Big( \int_0^t \dt^\prime 
    \frac{\sin\big( (t-t^\prime) \langle n_{234} \rangle\big)}{\langle n_{234} \rangle} \slinear[blue][N_2](t^\prime,n_2) 
    \Big( \prod_{j=3}^4 e^{i (\pm_j \langle n_j \rangle + \lambda_j ) t^\prime } \Big) 
    \Big) \bigg] \notag \\ 
=& \sum_{\substack{\varphi_1, \varphi_2 \\ \in \{ \cos, \sin \}}} \sum_{n_1,n_2\in \Z^3} \bigg[ 1_{N_{234}}(n_{234}) \Big( \prod_{j=1}^2 1_{N_j}(n_j) \Big) e^{i\langle n_{1234}, x \rangle} \langle n_1 \rangle^{-1} \langle n_{234} \rangle^{-1} \langle n_2 \rangle^{-1} \varphi_1\big( t \langle n_1 \rangle \big)  \label{para:eq-resonant-X1-p2}  \\
    &\hspace{10ex}\times \Big( \int_0^t \dt^\prime 
    \sin\big( (t-t^\prime) \langle n_{234} \rangle\big)  \varphi_2(t^\prime \langle n_2 \rangle) 
    \Big( \prod_{j=3}^4 e^{i (\pm_j \langle n_j \rangle + \lambda_j ) t^\prime } \Big) 
    \Big) \SI[n_j ,\varphi_j \colon 1\leq j \leq 2] \bigg] \notag \\
    +& \, e^{i\langle n_{34}, x\rangle} 1\big\{ N_1 = N_2 \big\} \bigg( \int_0^t \dt^\prime \, \Sine[N_{234},N_2](t-t^\prime,n_{34}) \Big( \prod_{j=3}^4 e^{i (\pm_j \langle n_j \rangle + \lambda_j ) t^\prime } \Big) \bigg).   \label{para:eq-resonant-X1-p3}  
\end{align}
We now estimate the non-resonant and resonant parts \eqref{para:eq-resonant-X1-p2} and \eqref{para:eq-resonant-X1-p3} separately. \\

For the non-resonant part \eqref{para:eq-resonant-X1-p2}, we have that 
\begin{align*}
\E \Big[ \big\| \eqref{para:eq-resonant-X1-p2} \big\|_{H_x^{-1/2-\epsilon/2}}^2 \Big] 
\lesssim&\, \sum_{n_1,n_2 \in \Z^3} \bigg[ 
1_{N_{234}}(n_{234}) \Big( \prod_{j=1}^2 1_{N_j}(n_j) \Big) \langle n_{1234} \rangle^{-1-\epsilon} \langle n_1 \rangle^{-2} \langle n_{234} \rangle^{-2} \langle n_2 \rangle^{-2} \bigg] \\
\lesssim&\, N_1^{-\epsilon} \sum_{n_2 \in \Z^3} \bigg[ 
1_{N_{234}}(n_{234})  1_{N_2}(n_2)  \langle n_{234} \rangle^{-2} \langle n_2 \rangle^{-2} \bigg] \\
\lesssim& \, N_1^{-\epsilon} N_{234}^{-1}. 
\end{align*}

For the resonant part, Lemma \ref{counting:lem-Sine-estimate} implies that 
\begin{equation*}
\big\| \eqref{para:eq-resonant-X1-p3} \big\|_{L_x^2}^2 \lesssim 1\big\{ N_1 = N_2 \big\} N_{234}^{-2+\epsilon},
\end{equation*}
which is more than acceptable. \\

\emph{Proof of \ref{para:item-resonant-X1-res}:}  
The argument is similar to the proof of  \ref{para:item-resonant-X1-hhl} and we omit the details. The only additional ingredients are that the sum of the regularities of $w_2$ and $w_3$ is positive, i.e., $(-1/2-\epsilon)+(1/2+\delta_2)>0$, and a standard box localization argument (see e.g. \cite[Corollary 3.13]{Tao01}). \\

\emph{Proof of \ref{para:item-resonant-X1-expl}:}  Together with Lemma \ref{diagrams:lem-mathfrakC} and Lemma \ref{analytic:lem-C15-C33}, the proof for the $\mathfrak{C}_{\leq N}$, $\mathfrak{C}_{\leq N}^{(1,5)}$, and $\mathfrak{C}_{\leq N}^{(3,3)}$-terms in $\XXexpl$ is as in \ref{para:item-resonant-X1-hhl}. Therefore, we only treat the $(\gamma_{\leq N}-\Gamma_{\leq N})$-term.

We now set $N_1 := K_1$ and $N_2 := K_2$. We first decompose the product into its non-resonant and resonant components. More precisely, we decompose 
\begin{align}
&\slinear[blue][N_1] \big( \gamma_{\leq N} - \Gamma_{\leq N} \big) 
\Duh \Big[ \, \slinear[blue][N_2] \Big] \notag \\
=& \sum_{\substack{\varphi_1, \varphi_2 \in \\ \{ \cos, \sin \} }}
\sum_{n_1,n_2 \in \Z^3} \bigg[ \Big( \prod_{j=1}^2 1_{N_j}(n_j) \Big) \langle n_1 \rangle^{-1} \langle n_2 \rangle^{-2} \big( \gamma_{\leq N} - \Gamma_{\leq N}(n_2) \big) \, e^{i \langle n_{12},x\rangle} \label{para:eq-resonant-X1-expl-p1}\\
&\times \varphi_1\big( t\langle n_1\rangle\big) 
\Big( \int_0^t \dt^\prime \sin\big( (t-t^\prime) \langle n_2 \rangle \big) \varphi_2\big( t^\prime \langle n_2 \rangle \big) \Big) 
\SI[n_j,\varphi_j\colon 1 \leq j \leq 2] \bigg] \notag \\
+& 
\sum_{n_1,n_2 \in \Z^3} \bigg[ 1\big\{ n_{12}=0 \big\} \Big( \prod_{j=1}^2 1_{N_j}(n_j) \Big) \langle n_1 \rangle^{-1} \langle n_2 \rangle^{-2} \big( \gamma_{\leq N} - \Gamma_{\leq N}(n_2) \big)  \label{para:eq-resonant-X1-expl-p2}\\
&\times
\Big( \int_0^t \dt^\prime \sin\big( (t-t^\prime) \langle n_2 \rangle \big) \cos\big( (t-t^\prime) \langle n_2 \rangle \big) \Big) 
 \bigg] \notag.
\end{align}
Together with Lemma \ref{analytic:lem-difference-gamma}, the non-resonant part \eqref{para:eq-resonant-X1-expl-p1} can be treated exactly as in \ref{para:item-resonant-X1-hhl}. Thus, it only remains to treat the (purely deterministic) resonant part \eqref{para:eq-resonant-X1-expl-p2}. This argument is similar as for the sine-cancellation (Lemma \ref{counting:lem-Sine-symmetrization}), but requires a minor modification due to the symbol of $(\gamma_{\leq N}-\Gamma_{\leq N})$. More precisely, we have that 
\begin{align}
&\hspace{1ex}\eqref{para:eq-resonant-X1-expl-p2} \notag \\
=& 1\big\{ N_1 = N_2 \big\} \sum_{n\in \Z^3}\bigg[ 
 \Big( \prod_{j=1}^2 1_{N_j}(n) \Big) \langle n \rangle^{-3} \big( \gamma_{\leq N} - \Gamma_{\leq N}(n) \big) 
\Big( \int_0^t \dt^\prime \sin\big( (t-t^\prime) \langle n \rangle \big) \cos\big( (t-t^\prime) \langle n \rangle \big) \Big) \bigg] \notag \\
=& \frac{1}{4} 1\big\{ N_1 = N_2 \big\}  
\sum_{n\in \Z^3}\bigg[ 
 \Big( \prod_{j=1}^2 1_{N_j}(n) \Big) \langle n \rangle^{-4} \big( \gamma_{\leq N} - \Gamma_{\leq N}(n) \big) \big( \cos\big( 2 t \langle n \rangle\big) -1 \big) \bigg].
 \label{para:eq-resonant-X1-expl-p3}
\end{align}
Using Lemma \ref{analytic:lem-difference-gamma}, we obtain that
\begin{equation*}
\big| \eqref{para:eq-resonant-X1-expl-p3} \big| \lesssim 1\big\{ N_1 = N_2\big\} N_2^{-1+\epsilon}, 
\end{equation*}
which yields the desired estimate.
\end{proof}

We now turn to the product estimate for the $\SXXtwo$-term.

\begin{lemma}[\protect{Product estimates for $\slinear[blue] \, \SXXtwo$}]\label{para:lem-resonant-X2}
Let $T\geq 1$ and $p\geq 2$. Let $N_1,N_2,N_3,N_4,N_{234}\geq 1$ be frequency scales such that $N_4 \leq \max(N_2,N_3)^\eta < \min(N_2,N_3)$, then it holds that 
\begin{align*}
\E \bigg[ \sup_{\mathcal{J}} \Big\| w_4 \mapsto \slinear[blue][N_1] \XXop[N_\ast,w_4] \Big\|_{X^{-1,b}(\mathcal{J})\rightarrow L_t^\infty \C_x^{-1/2-\epsilon}(\mathcal{J})}^p \bigg]^{1/p} \lesssim p T^{\alpha} \max(N_2,N_3)^{-1/2+5\eta}, 
\end{align*}
where the supremum is taken over all closed intervals $0 \in \mathcal{J} \subseteq [-T,T]$ and $N_\ast$ denotes the dependence on $N_{234},N_2,N_3$, and $N_4$. 
\end{lemma}

\begin{proof}
We only treat the case $T=1$ and $\mathcal{J}=[-1,1]$, since the general case follows from a minor modification.  \\

Using similar reduction arguments as in Subsection \ref{prep:remark-reduction}, it suffices to prove that 
\begin{equation}\label{para:eq-resonant-X2-p1}
\begin{aligned}
&\sup_{t\in [-1,1]} \sup_{\pm_4} \sup_{|n_4|\leq \max(N_2,N_3)^\eta} \sup_{\lambda_4 \in \R} \E \Bigg[ \bigg\| \sum_{n_1,n_2\in \Z^3} \bigg[ 1_{N_{234}}(n_{234}) \Big( \prod_{j=1}^3 1_{N_j}(n_j) \Big) e^{i\langle n_{1234}, x \rangle}  \\
    &\hspace{2ex} \times \slinear[blue][N_1](t,n_1) \Big( \int_0^t \dt^\prime 
    \frac{\sin\big( (t-t^\prime) \langle n_{234} \rangle\big)}{\langle n_{234} \rangle} \lcol \slinear[blue][N_2](t^\prime,n_2) \, \slinear[blue][N_3](t^\prime,n_3) \rcol  
     e^{i (\pm_4 \langle n_4 \rangle + \lambda_4 ) t^\prime } 
    \Big) \bigg] \bigg\|_{H_x^{-1/2-\epsilon/2}}^2 \bigg] \\
    &\lesssim N_1^{-\epsilon} \max(N_2,N_3)^{-1}. 
   \end{aligned}
\end{equation}
In order to prove \eqref{para:eq-resonant-X2-p1}, we separate the non-resonant and resonant parts. 
By inserting the definition of the $\Sine$-kernel (from Definition \ref{def:sine_kernel}), we obtain that 
\begin{align}
 &\sum_{n_1,n_2\in \Z^3} \bigg[ 1_{N_{234}}(n_{234}) \Big( \prod_{j=1}^3 1_{N_j}(n_j) \Big) e^{i\langle n_{1234}, x \rangle} \slinear[blue][N_1](t,n_1)   \notag \\
    &\hspace{2ex} \times  \Big( \int_0^t \dt^\prime 
    \frac{\sin\big( (t-t^\prime) \langle n_{234} \rangle\big)}{\langle n_{234} \rangle} \lcol \slinear[blue][N_2](t^\prime,n_2) \, \slinear[blue][N_3](t^\prime,n_3) \rcol  
     e^{i (\pm_4 \langle n_4 \rangle + \lambda_4 ) t^\prime } 
    \Big) \bigg]  \notag  \allowdisplaybreaks[3] \\
=& \sum_{\substack{\varphi_1,\varphi_2,\varphi_3 \\ \in \{ \cos, \sin \} }} 
\sum_{n_1, n_2 ,n_3 \in \Z^3} \bigg[ 1_{N_{234}}(n_{234}) \Big( \prod_{j=1}^3 1_{N_j}(n_j) \Big) \langle n_1 \rangle^{-1} \langle n_{234} \rangle^{-1} \langle n_2 \rangle^{-1} \langle n_3 \rangle^{-1}  e^{i \langle n_{1234} , x \rangle}   \notag  \\
&\hspace{2ex} \times \varphi_1\big( t \langle n_1 \rangle\big) 
\Big( \int_0^t \dt^\prime \, \sin\big( (t-t^\prime) \langle n_{234} \rangle \big) 
\varphi_2\big( t^\prime \langle n_2 \rangle \big) 
 \varphi_3\big( t^\prime \langle n_3 \rangle \big) 
 e^{i (\pm_4 \langle n_4 \rangle + \lambda_4 ) t^\prime } \Big)  \label{para:eq-resonant-X2-p2} \\
& \hspace{2ex} \times \SI[n_j,\varphi_j \colon 1 \leq j \leq 3] \bigg]   \notag  \allowdisplaybreaks[3]\\ 
+& 1\big\{ N_1 = N_2 \big\} \sum_{\substack{\varphi_3 \in \\ \{ \cos, \sin \}}} \sum_{n_3 \in \Z^3}
\bigg[ 1_{N_3}(n_3)  \langle n_3 \rangle^{-1} e^{i \langle n_{34}, x \rangle} \label{para:eq-resonant-X2-p3}   \\
&\hspace{2ex} \times \Big( \int_0^t 
\Sine[N_{234},N_2](t-t^\prime,n_{34})
\varphi_3\big(t^\prime \langle n_3 \rangle \big)
 e^{i (\pm_4 \langle n_4 \rangle + \lambda_4 ) t^\prime } \Big) 
 \SI[n_3,\varphi_3] \bigg]  \notag  \allowdisplaybreaks[3] \\
 +& 1\big\{ N_1 = N_3 \big\} \sum_{\substack{\varphi_2 \in \\ \{ \cos, \sin \}}} \sum_{n_2 \in \Z^3}
\bigg[ 1_{N_2}(n_2)  \langle n_2 \rangle^{-1} e^{i \langle n_{24}, x \rangle} \label{para:eq-resonant-X2-p4}   \\
&\hspace{2ex} \times \Big( \int_0^t 
\Sine[N_{234},N_3](t-t^\prime,n_{24})
\varphi_2\big(t^\prime \langle n_2 \rangle \big)
 e^{i (\pm_4 \langle n_4 \rangle + \lambda_4 ) t^\prime } \Big) 
 \SI[n_2,\varphi_2] \bigg]. \notag  
\end{align}
By symmetry in the $n_2$ and $n_3$-variables, it suffices to treat \eqref{para:eq-resonant-X2-p2} and \eqref{para:eq-resonant-X2-p3}. \\

\emph{Estimate of the non-resonant term \eqref{para:eq-resonant-X2-p2}:} 
By calculating the expectation and performing the sum in $n_1$, we obtain that 
\begin{align}
&\E \Big[ \Big\| \eqref{para:eq-resonant-X2-p2} \Big\|_{H_x^{-1/2-\epsilon/2}}^2 \Big]  \allowdisplaybreaks[3] \notag \\
\lesssim& \,  \sum_{\substack{\varphi_1,\varphi_2,\varphi_3 \\ \in \{ \cos, \sin \} }} 
\sum_{n_1, n_2 ,n_3 \in \Z^3} \bigg[ 1_{N_{234}}(n_{234}) \Big( \prod_{j=1}^3 1_{N_j}(n_j) \Big) \langle n_{1234} \rangle^{-1-\epsilon} \langle n_1 \rangle^{-2} \langle n_{234} \rangle^{-2} \langle n_2 \rangle^{-2} \langle n_3 \rangle^{-2}    \notag  \\
&\hspace{2ex} \times 
\bigg| \int_0^t \dt^\prime \, \sin\big( (t-t^\prime) \langle n_{234} \rangle \big) 
\varphi_2\big( t^\prime \langle n_2 \rangle \big) 
 \varphi_3\big( t^\prime \langle n_3 \rangle \big) 
 e^{i (\pm_4 \langle n_4 \rangle + \lambda_4 ) t^\prime } \bigg|^2 \bigg] \notag  \allowdisplaybreaks[3] \\
 \lesssim&\, N_1^{-\epsilon}  \sum_{\substack{\varphi_2,\varphi_3 \\ \in \{ \cos, \sin \} }} 
\sum_{n_2 ,n_3 \in \Z^3} \bigg[ 1_{N_{234}}(n_{234}) \Big( \prod_{j=2}^3 1_{N_j}(n_j) \Big) \langle n_{234} \rangle^{-2} \langle n_2 \rangle^{-2} \langle n_3 \rangle^{-2}   \notag  \\
&\hspace{2ex} \times 
\bigg| \int_0^t \dt^\prime \, \sin\big( (t-t^\prime) \langle n_{234} \rangle \big) 
\varphi_2\big( t^\prime \langle n_2 \rangle \big) 
 \varphi_3\big( t^\prime \langle n_3 \rangle \big) 
 e^{i (\pm_4 \langle n_4 \rangle + \lambda_4 ) t^\prime } \bigg|^2 \bigg].  \label{para:eq-resonant-X2-p5}
\end{align}
By performing the $t^\prime$-integral and using our counting estimate (Lemma \ref{counting:lem2}), 
we obtain that 
\begin{equation*}
\eqref{para:eq-resonant-X2-p5} \lesssim N_1^{-\epsilon} \max\big( N_{234},N_2,N_3\big)^{-1}, 
\end{equation*}
which is acceptable. \\ 

\emph{Estimate of the resonant term \eqref{para:eq-resonant-X2-p3}:} 
Using Lemma \ref{counting:lem-Sine-estimate}, we obtain that 
\begin{align*}
&\E \Big[ \Big\| \eqref{para:eq-resonant-X2-p3} \Big\|_{H_x^{-1/2-\epsilon/2}}^2 \Big] \\
\lesssim& \sum_{\substack{\varphi_3 \in \\ \{ \cos, \sin \}}}
\sum_{n_3 \in \Z^3} \bigg[ 1_{N_3}(n_3) 
\langle n_{34} \rangle^{-1-\epsilon} \langle n_3 \rangle^{-2} 
\bigg| \int_0^t \dt^\prime \Sine[N_{234},N_2](t-t^\prime,n_{34}) \varphi_3\big( t^\prime \langle n_{3} \rangle \big) e^{i (\pm_4 \langle n_4 \rangle + \lambda_4 ) t^\prime } \bigg|^2 \bigg] \\
\lesssim& \max\big(N_{234},N_2\big)^{-2+\epsilon} 
\sum_{n_3 \in \Z^3} \bigg[ 1_{N_3}(n_3) 
\langle n_{34} \rangle^{-1-\epsilon} \langle n_3 \rangle^{-2} \bigg] \\
\lesssim& \max\big(N_{234},N_2\big)^{-2+\epsilon}. 
\end{align*}
Since $N_4 \leq \max(N_2,N_3)^\eta$, it holds that $\max(N_{234},N_2) \sim \max(N_2,N_3)$. As a result, this contribution is (more than) acceptable. 
\end{proof}

It remains to prove the main result of this subsection (Proposition \ref{para:prop-resonant}). 

\begin{proof}[Proof of Proposition \ref{para:prop-resonant}]
Due to Lemma \ref{para:lem-decomposition}, $\XXone$ and $\XXtwo$ can be rewritten in terms of $\XXhll[N_\ast]$, $\XXres[N_\ast]$, $\XXexpl$,  and $\XXop[N_\ast]$. The desired estimate now follows directly from Lemma \ref{para:lem-resonant-X1} and Lemma \ref{para:lem-resonant-X2}. 
\end{proof}

At the end of this subsection, we prove the following corollary of the product estimate (Proposition \ref{para:prop-resonant}) and Lemma \ref{bilinear:lem-variant}. It will only be needed in the proof of Lemma \ref{global:lem-one-zero-revisited}, which is presented in Section \ref{section:proof-main-estimates}, and we encourage the reader to skip this corollary on first reading.

\begin{corollary}[\protect{The $\slinear[blue][\leqN] \XXj$-operator}]\label{para:cor-blueXj}
There exists an $A$-certain event $E_A \in \mathcal{E}$ such that, on this event, the following estimates hold: Let $T\geq 1$, let $0 \in \mathcal{J} \subseteq [-T,T]$ be a closed interval, and let $N,N_1,N_2,N_3$ be frequency-scales satisfying 
\begin{equation*}
N_1,N_2,N_3\leq N, \quad N_2 \geq N_1^\eta, \quad \text{and} \quad N_2 \geq N_3. 
\end{equation*}
Furthermore, let $v_{\leq N},w_{\leq N}, Y_{\leq N}\colon \mathcal{J} \times \T^3 \rightarrow \R$. Then,
\begin{align*}
&\max_{j=1,2} \bigg\| \slinear[blue][N_1] P_{N_2} \XXj [v_{\leq N},Y_{\leq N}] P_{N_3} w_{\leq N} \bigg\|_{X^{-1/2+\delta_2,b_+-1}(\mathcal{J})} \\
\leq& \, A T^\alpha \max\big(N_1,N_2\big)^{-\epsilon} 
\Big( 1+ \big\| v_{\leq N} \big\|_{X^{-1,b}(\mathcal{J})}^2 + \big\| Y_{\leq N} \big\|_{X^{1/2+\delta_2,b}(\mathcal{J})} \Big) \big\| w_{\leq N} \big\|_{X^{-\epsilon,b}(\mathcal{J})}. 
\end{align*}
\end{corollary}

\begin{remark}
Our notation for the frequency-scales in Corollary \ref{para:cor-blueXj} is chosen so that it agrees with our notation in Section \ref{section:proof-main-estimates}.  
\end{remark}

\begin{proof}
In the following, all estimates are made after possible restrictions to $A$-certain events. 
We only treat the case $T=1$ and $\mathcal{J}=[-1,1]$, since the general case follows from minor modifications. \\

To prove the desired estimate, we first decompose 
\begin{equation}\label{para:eq-blueXj-p1}
\begin{aligned}
&\slinear[blue][N_1] P_{N_2} \XXj [v_{\leq N},Y_{\leq N}] P_{N_3} w_{\leq N} \\
=& \sum_{N_0,N_{12}} P_{N_0} \bigg[ P_{N_{12}} \Big( \slinear[blue][N_1] P_{N_2} \XXj [v_{\leq N},Y_{\leq N}] \Big) P_{N_3} w_{\leq N}\bigg].
\end{aligned}
\end{equation}
We now proceed using two different estimates. Using $X^{-1/2+\delta_2,0}\hookrightarrow X^{-1/2+\delta_2,b_+-1}$, our product estimate (Proposition \ref{para:prop-resonant}), and $N_3 \leq N_2$, it holds that 
\begin{align}
&\bigg\|  P_{N_0} \bigg[ P_{N_{12}} \Big( \slinear[blue][N_1] P_{N_2} \XXj [v_{\leq N},Y_{\leq N}] \Big) P_{N_3} w_{\leq N}\bigg] \bigg\|_{X^{-1/2+\delta_2,b_+-1}}  \notag \\
\lesssim& \, \bigg\|  P_{N_0} \bigg[ P_{N_{12}} \Big( \slinear[blue][N_1] P_{N_2} \XXj [v_{\leq N},Y_{\leq N}] \Big) P_{N_3} w_{\leq N}\bigg] \bigg\|_{L_t^2 H_x^{-1/2+\delta_2}} \notag \\
\lesssim& \, N_0^{-1/2+\delta_2} \Big\| P_{N_{12}} \Big( \slinear[blue][N_1] P_{N_2} \XXj [v_{\leq N},Y_{\leq N}] \Big) \Big\|_{L_t^\infty L_x^\infty} \big\| w_{\leq N} \big\|_{L_t^\infty L_x^2} \notag \\
\lesssim& \, N_0^{-1/2+\delta_2} N_{12}^{1/2+\epsilon} N_2^{-1/2+10\eta} N_3^\epsilon
\Big( 1 + \| v_{\leq N} \|_{X^{-1,b}}^2 + \| Y_{\leq N}\|_{X^{1/2+\delta_2,b}} \Big)
\big\| w_{\leq N} \big\|_{X^{-\epsilon,b}}. \label{para:eq-blueXj-p2}
\end{align}
If $N_1 \gg N_2$, then $N_0 \sim N_{12} \sim N_1$. Together with our assumptions $N_2 \geq N_1^\eta$ and $N_2 \geq N_3$, the pre-factor in \eqref{para:eq-blueXj-p2} can then be bounded by 
\begin{equation*}
N_0^{-1/2+\delta_2} N_{12}^{1/2+\epsilon} N_2^{-1/2+10\eta} N_3^\epsilon \lesssim N_1^{\epsilon+\delta_2} N_2^{-1/2+10\eta+\epsilon} \lesssim N_1^{-\epsilon},
\end{equation*}
which is acceptable. Alternatively, we can simply rewrite the pre-factor in \eqref{para:eq-blueXj-p2} as 
\begin{equation*}
N_0^{-1/2+\delta_2} \Big( N_{12}/N_2 \Big)^{1/2+\epsilon} N_2^{\epsilon+10\eta}.
\end{equation*}
Thus, if either $N_0 \geq N_2^{1/100}$ or $N_{12} \leq N_2^{1-1/100}$, this also yields an acceptable contribution. Therefore, it only remains to treat the frequency-regime
\begin{equation*}
N_2 \gtrsim N_1,N_3, \quad N_0 \leq N_2^{1/100}, \quad \text{and} \quad N_{12} \geq N_2^{1-1/100}. 
\end{equation*}
In this case, the desired estimate follows directly from Lemma \ref{bilinear:lem-variant}.
\end{proof}

\subsection{Interactions with the quadratic stochastic object}\label{section:para-two-linear}

We now turn to the interactions of the para-controlled operators with the quadratic stochastic object $\squadratic[\leqN]$, i.e., two factors of $\slinear[blue]$. More precisely, we now estimate
\begin{gather}
\squadratic[\leqN] \XXone - \Big( 2 \HLL + \HHL \Big) \Big( \slinear[blue][\leqN], \slinear[blue][\leqN], \XXone \Big) \label{para:eq-motivation-11X1} \\
\text{and} \quad 3 \squadratic[\leqN] \XXtwo - \Big( 6 \HLL + 3 \HHL \Big) \Big( \slinear[blue][\leqN], \slinear[blue][\leqN], \XXtwo \Big) + \Gamma_{\leq N} \Big( 3 \squintic[\leqN] - P_{\leq N} \slinear[green][\leqM]  + v_{\leq N} \Big). \label{para:eq-motivation-11X2}
\end{gather}
The main estimates on \eqref{para:eq-motivation-11X1} and \eqref{para:eq-motivation-11X2} have been previously stated in Proposition \ref{ansatz:prop-Y-two}. In the following, we prove estimates for the individual components of $\XXone$ and $\XXtwo$ (as in Definition \ref{para:def-decomposition}), which will lead to the desired conclusions. 

\begin{lemma}[\protect{Estimate of $\squadratic[\leqN] \XXone$-terms}]\label{para:lem-11X1}
For all $T\geq 1$ and $p\geq 2$, we have the following estimates:
\begin{enumerate}[label=(\roman*)]
    \item\label{para:item-11X1-hll} For all frequency-scales $N_0,\hdots,N_5,N_{234}$ satisfying 
    \begin{equation*}
    N_{234} > \max(N_1,N_5)^\eta \qquad \text{and} \qquad N_3,N_4 \leq N_2^\eta,
    \end{equation*}
    it holds that 
    \begin{equation}\label{para:eq-11X1-hll}
    \begin{aligned}
    &\E \Big[ \sup_{\mathcal{J}} \Big\| 
    P_{N_0} \Big[ \lcol \slinear[blue][N_1] \slinear[blue][N_5] \rcol \XXhll[N_\ast,w_3,w_4] 
    \Big\|_{X^{-1,b}(\mathcal{J}) \times X^{-1,b}(\mathcal{J})\rightarrow X^{-1/2+\delta_2,b_+-1}(\mathcal{J})}^p \Big]^{1/p} \\
    \lesssim&\, p^{3/2} T^\alpha N_{\textup{max}}^{-\epsilon},
    \end{aligned}
    \end{equation}
    where the supremum is taken over all closed intervals $0\in \mathcal{J} \subseteq [-T,T]$ and $N_\ast$ denotes the dependence on $N_{234},N_2,N_3$, and $N_4$. 
    \item\label{para:item-11X1-res} For all frequency-scales $N_0,\hdots,N_5,N_{34},N_{234}$ satisfying 
    \begin{equation*}
    N_{234} > \max(N_1,N_5)^\eta, \quad N_4>N_2^\eta, \quad \text{and} \quad N_{34} \leq N_2^\eta, 
    \end{equation*}
    it holds that 
    \begin{equation}
    \begin{aligned}
    &\E \Big[ \sup_{\mathcal{J}} \Big\| 
    P_{N_0} \Big[ \lcol \slinear[blue][N_1] \slinear[blue][N_5] \rcol \XXres[N_\ast,w_3,w_4] 
    \Big\|_{X^{-1/2-\epsilon,b}(\mathcal{J}) \times X^{1/2+\delta_2,b}(\mathcal{J})\rightarrow X^{-1/2+\delta_2,b_+-1}(\mathcal{J})}^p \Big]^{1/p} \\
    \lesssim&\, p^{3/2} T^\alpha N_{\textup{max}}^{-\epsilon},
    \end{aligned}
    \end{equation}
    where the supremum is taken over all closed intervals $0\in \mathcal{J} \subseteq [-T,T]$ and $N_\ast$ denotes the dependence on $N_{234},N_{34},N_2,N_3$, and $N_4$. \item\label{para:item-11X1-expl} For all frequency-scales $K_0,K_1,K_2,K_3$ satisfying $K_3>\max(K_1,K_2)^\eta$, it holds that 
    \begin{align*}
    \E \bigg[ \sup_N \Big\| P_{K_0} \Big[ \lcol 
    \slinear[blue][K_1] \slinear[blue][K_2] \rcol P_{K_3} \XXexpl \Big] \Big\|_{X^{-1/2+\delta_2,b_+-1}([-T,T])}^p \bigg]^{1/p} \lesssim p^{3/2} T^\alpha K_{\textup{max}}^{-\epsilon}. 
    \end{align*} 
\end{enumerate}
\end{lemma}

\begin{proof}
We only treat the case $T=1$ and $\mathcal{J}=[-1,1]$, since the general case follows from minor modifications. \\

\emph{Proof of \ref{para:item-11X1-hll}:} We first decompose the argument in \eqref{para:eq-11X1-hll} into non-resonant and resonant terms. From the definition of $\XXhll$ and the definition of the $\Sine$-kernel (Definition \ref{def:sine_kernel}), we obtain that 
\begin{align}
      &P_{N_0} \Big[ \lcol \slinear[blue][N_1] \slinear[blue][N_5] \rcol \XXhll[N_\ast,w_3,w_4] \Big] \notag \allowdisplaybreaks[3] \\
      =&
      \sum_{\substack{n_0,\hdots,n_5 \in \Z^3 \colon \\ n_0 =n_{12345}}}
      \bigg[ 1_{N_{234}}(n_{234}) \Big( \prod_{j=0}^5 1_{N_j}(n_j) \Big) 
      e^{i \langle n_0 , x \rangle} 
      \lcol \slinear[blue][N_1](t,n_1) \slinear[blue][N_5](t,n_5)\rcol \notag \\
      &\times \Big( \int_0^t \dt^\prime \, \frac{\sin\big((t-t^\prime) \langle n_{234}\rangle\big)}{\langle n_{234} \rangle} \slinear[blue][N_2](t^\prime,n_2) \widehat{w}_3(t^\prime,n_3) \widehat{w}_4(t^\prime,n_4) \bigg] \notag \allowdisplaybreaks[3] \\
      =& \sum_{\substack{\varphi_1,\varphi_2,\varphi_5 \in \\ \{\cos, \sin\} }}
      \sum_{\substack{n_0,n_1,\hdots,n_5 \in \Z^3 \colon \\ n_0 = n_{12345} }}
      \bigg[ 1_{N_{234}}(n_{234}) \Big( \prod_{j=0}^5 1_{N_j}(n_j) \Big) 
      \langle n_{234} \rangle^{-1} \langle n_1 \rangle^{-1} \langle n_2 \rangle^{-1} 
      \langle n_5 \rangle^{-1} e^{i\langle n_0 , x\rangle} \label{para:eq-11X1-p1} \\
      &\times \varphi_1\big( t \langle n_1 \rangle\big) 
      \varphi_5\big( t \langle n_5 \rangle\big) 
      \Big( \int_0^t \dt^\prime \, \sin\big( (t-t^\prime) \langle n_{234} \rangle\big) \varphi_2\big( t^\prime \langle n_2 \rangle \big) 
      \widehat{w}_3(t^\prime,n_3) \widehat{w}_4(t^\prime,n_4) \Big)  \notag \\
      &\times \SI[n_j,\varphi_j \colon j=1,2,5] \bigg] \notag \allowdisplaybreaks[3] \\
      +& 1\big\{ N_1= N_2 \big\} \sum_{\substack{\varphi_5 \in \\ \{\cos, \sin\} }}
      \sum_{\substack{n_0,n_3,n_4,n_5 \in \Z^3 \colon \\ n_0 = n_{345} }}
      \bigg[  \Big( \prod_{j=0,3,4,5} 1_{N_j}(n_j) \Big) 
      \langle n_5 \rangle^{-1} e^{i\langle n_0 , x\rangle} \label{para:eq-11X1-p2} \\
      &\times  
      \varphi_5\big( t \langle n_5 \rangle\big) 
      \Big( \int_0^t \dt^\prime \, \Sine[N_{234},N_2](t-t^\prime,n_{34})  
      \widehat{w}_3(t^\prime,n_3) \widehat{w}_4(t^\prime,n_4) \Big)  
      \SI[n_5,\varphi_5] \bigg] \notag \allowdisplaybreaks[3] \\ 
       +& 1\big\{ N_2= N_5 \big\} \sum_{\substack{\varphi_1 \in \\ \{\cos, \sin\} }}
      \sum_{\substack{n_0,n_1,n_3,n_4 \in \Z^3 \colon \\ n_0 = n_{134} }}
      \bigg[  \Big( \prod_{j=0,1,3,4} 1_{N_j}(n_j) \Big) 
      \langle n_1 \rangle^{-1} e^{i\langle n_0 , x\rangle} \label{para:eq-11X1-p3} \\
      &\times  
      \varphi_1\big( t \langle n_1 \rangle\big) 
      \Big( \int_0^t \dt^\prime \, \Sine[N_{234},N_2](t-t^\prime,n_{34})  
      \widehat{w}_3(t^\prime,n_3) \widehat{w}_4(t^\prime,n_4) \Big)  
      \SI[n_1,\varphi_1] \bigg] \notag.
\end{align}
By symmetry in $n_1$ and $n_5$, it suffices to treat the non-resonant term \eqref{para:eq-11X1-p1} and the resonant term \eqref{para:eq-11X1-p2}. \\

\emph{The non-resonant term \eqref{para:eq-11X1-p1}:} We use the quintic tensor from Lemma \ref{counting:quintic_tensor},
\begin{equation}\label{para:eq-11X1-p4}
h_{n_0 n_1 \hdots n_5}(t,\lambda_3,\lambda_4) 
= h_{n_0 n_1 \hdots n_5}[N_0,\hdots,N_5,N_{234},\pm_1,\hdots,\pm_5](t,\lambda_3,\lambda_4),
\end{equation}
where we set $\lambda_1=\lambda_2=\lambda_5=0$. Using \eqref{para:eq-11X1-p4}, we can rewrite the non-resonant term as 
\begin{align}
\eqref{para:eq-11X1-p1} \notag 
=& \sum_{\substack{\pm_1,\pm_2,\pm_3,\\ \pm_4,\pm_5}} \int_{\R^2} \dlambda_3 \dlambda_4
\sum_{\substack{n_0,n_1,n_2 \\ n_3,n_4,n_5}} \bigg[  e^{i \langle n_0 ,x \rangle} 
h_{n_0 n_1 \hdots n_5}(t,\lambda_3,\lambda_4)  \notag  \\
&\times 
\widetilde{\langle \nabla \rangle w_3^{\pm_3}}(n_3,\lambda_3) \,
\widetilde{\langle \nabla \rangle w_4^{\pm_4}}(n_4,\lambda_4) \,
\SI[n_j,\pm_j\colon j=1,2,5] \bigg].\label{para:eq-11X1-p5}
\end{align}
The $\langle \nabla \rangle$-multiplier acting on $w_3$ and $w_4$ is due to the $\langle n_3 \rangle^{-1}$ and $\langle n_4 \rangle^{-1}$-factors in the definition of the quintic tensor. Due to the frequency assumptions $N_3,N_4 \leq N_2^\eta$, however, the $\langle \nabla \rangle$-multipliers are essentially irrelevant. Using the tensor estimate reduction argument in Subsection \ref{prep:remark-reduction},  it follows that 
\begin{align}
&\Big\| \eqref{para:eq-11X1-p5} \Big\|_{X^{-1,b}\times X^{-1,b} \rightarrow X^{-1/2+\delta_2,b_+-1}} \notag \\
\lesssim& \, N_0^{-1/2+\delta_2} N_3^2 N_4^2 
\max_{\pm_0,\pm_1,\hdots,\pm_5} \sup_{\lambda_3,\lambda_4 \in \R} \bigg[ 
\langle \lambda_3 \rangle^{-(b_--1/2)}\langle \lambda_4 \rangle^{-(b_--1/2)} \notag \\
&\times \Big\| \langle \lambda \rangle^{b_+-1} 
\Big\| \sum_{\substack{n_1,n_2,n_5}} 
\widetilde{h}^{\pm_0}_{n_0 n_1 \hdots n_5}(\lambda,\lambda_3,\lambda_4) 
\SI[n_j, \pm_j \colon j=1,2,5] \Big\|_{n_3 n_4 \rightarrow n_0} \Big\|_{L_\lambda^2} \bigg]. \label{para:eq-11X1-p6}
\end{align}
Since $N_3$ and $N_4$ are small compared to $N_2$, we simply bound\footnote{Of course, this bound is always possible, but is typically rather crude for large values of $N_3$ and $N_4$.} the $\| \cdot\|_{n_3n_4 \rightarrow n_0}$-norm by the Hilbert-Schmidt norm $\| \cdot \|_{n_0n_3n_4}$. Using the $p$-moment estimate reduction argument in Subsection \ref{prep:remark-reduction} and the Hilbert-Schmidt estimate \eqref{counting:extra_HS1}
in Lemma \ref{counting:quintic_tensor}, it follows that 
\begin{align*}
\E\Big[ \eqref{para:eq-11X1-p6}^p \Big]^{1/p}
\lesssim& \, p^{3/2} N_{\textup{max}}^\epsilon N_0^{-1/2+\delta_2} (N_3 N_4)^2 \\
&\times
\min\big( N_0,N_1,N_2,N_5)^{1/2} \max(N_0,N_1,N_2,N_5)^{-1/2} N_0 N_2^{-1} N_3^{1/2} N_4^{1/2} \\
\lesssim& \, p^{3/2} N_{\textup{max}}^{\epsilon+\delta_2} N_2^{-1/2} N_3^{5/2} N_4^{5/2}. 
\end{align*}
Due to our frequency assumptions, this is acceptable.

\emph{The resonant term \eqref{para:eq-11X1-p2}:} The argument is similar to the estimate of \eqref{para:eq-11X1-p1}, but with the quintic tensor replaced by the sine-cancellation tensor. To be precise, we use the sine-cancellation tensor
\begin{equation}\label{para:eq-11X1-p7}
h^{\textup{\sine}}_{n_0n_3n_4n_5}(t,\lambda_3,\lambda_4)
= h^{\textup{\sine}}_{n_0n_3n_4n_5}[N_0,N_1,\hdots,N_5,N_{234},\pm_3,\pm_4,\pm_5](t,\lambda_3,\lambda_4)
\end{equation}
from Lemma \ref{counting:sine_tensor}, where we set $\lambda_5=0$. Using \eqref{para:eq-11X1-p7}, we can rewrite the resonant term as 
\begin{align}
\eqref{para:eq-11X1-p2} &= 
\sum_{\pm_3,\pm_4,\pm_5} \int_{\R^2} \dlambda_3 \dlambda_4 \sum_{n_0,n_3,n_4,n_5} \bigg[ e^{i\langle n_0 ,x \rangle} h^{\textup{\sine}}_{n_0n_3n_4n_5}(t,\lambda_3,\lambda_4) \notag \\
&\times 
\widetilde{\langle \nabla \rangle w_3^{\pm_3}}(n_3,\lambda_3) \,
\widetilde{\langle \nabla \rangle w_4^{\pm_4}}(n_4,\lambda_4) \,
\SI[n_5,\pm_5] \bigg]. \label{para:eq-11X1-p8}
\end{align}
As before, the $\langle \nabla \rangle$-multipliers result from the definition of the sine-cancellation tensor, but are essentially irrelevant.  Use the tensor estimate reduction argument in Subsection \ref{prep:remark-reduction}, we obtain that\\
\begin{align}
&\Big\| \eqref{para:eq-11X1-p8} \Big\|_{X^{-1,b}\times X^{-1,b} \rightarrow X^{-1/2+\delta_2,b_+-1}} \notag \\
\lesssim& \, N_0^{-1/2+\delta_2} N_3^2 N_4^2 
\max_{\pm_0,\pm_1,\hdots,\pm_5} \sup_{\lambda_3,\lambda_4 \in \R} \bigg[ 
\langle \lambda_3 \rangle^{-(b_--1/2)}\langle \lambda_4 \rangle^{-(b_--1/2)} \notag \\
&\times \Big\| \langle \lambda \rangle^{b_+-1} 
\Big\| \sum_{\substack{n_5}} 
\widetilde{h}^{\textup{\sine},\pm_0}_{n_0 n_3 n_4 n_5}(\lambda,\lambda_3,\lambda_4) 
\SI[n_5,\pm_5] \Big\|_{n_3 n_4 \rightarrow n_0} \Big\|_{L_\lambda^2} \bigg]. \label{para:eq-11X1-p9}
\end{align}
As before, we estimate the operator norm $\| \cdot \|_{n_3 n_4 \rightarrow n_0}$ by the Hilbert-Schmidt norm $\| \cdot \|_{n_0 n_3 n_4}$. 
Using the $p$-moment estimate reduction argument in Subsection \ref{prep:remark-reduction} and the Hilbert-Schmidt estimate \eqref{counting:sine_extra_HS} in Lemma \ref{counting:sine_tensor}, we obtain that 
\begin{align*}
\E\Big[ \eqref{para:eq-11X1-p9}^p \Big]^{1/p}
\lesssim& \, p^{1/2} N_{\textup{max}}^\epsilon N_0^{-1/2+\delta_2} (N_3 N_4)^2 \\
&\times
\min(N_0,N_5)^{1/2} N_2^{-1} N_3^{1/2} N_4^{1/2} \\
\lesssim& \, p^{1/2} N_{\textup{max}}^{\epsilon+\delta_2} N_2^{-1} N_3^{5/2} N_4^{5/2}.
\end{align*}
Due to our assumptions on the frequency scales, this is acceptable.

\emph{Proof of \ref{para:item-11X1-res}:} 
The argument is similar to the proof of \ref{para:item-11X1-hll}. The only additional ingredient is a standard box-localization argument (see e.g. \cite[Corollary 3.13]{Tao01}) and we omit the standard details. \\

\emph{Proof of \ref{para:item-11X1-expl}:} Using Lemma \ref{diagrams:lem-mathfrakC} and Lemma \ref{analytic:lem-C15-C33}, the $\mathfrak{C}_{\leq N}$, $\mathfrak{C}_{\leq N}^{(1,5)}$, and $\mathfrak{C}_{\leq N}^{(3,3)}$-terms in $\XXexpl$ can be treated as in the proof of  \ref{para:item-11X1-hll}. Thus, it only remains to treat the $(\gamma_{\leq N}-\Gamma_{\leq N})$-term. To this end, we first set $N_j := K_j$, where $j=0,1,2,3$. Then, we decompose the product into its resonant and non-resonant components. More precisely, we decompose 
\begin{align}
&P_{N_0} \bigg[ \lcol \slinear[blue][N_1] \slinear[blue][N_2] \rcol 
\big( \gamma_{\leq N} - \Gamma_{\leq N} \big) \Duh \Big[ \, \slinear[blue][N_3] \Big] \bigg] \notag \allowdisplaybreaks[3] \\
=& \sum_{\substack{\varphi_1,\varphi_2,\varphi_3 \in \\ \{ \cos, \sin \}}}
 \sum_{\substack{n_0,n_1,n_2,n_3\in \Z^3 \colon \\ n_0=n_{123}}}
 \bigg[ \Big( \prod_{j=0}^3 1_{N_j}(n_j) \Big) \langle n_1 \rangle^{-1} \langle n_2 \rangle^{-1} \langle n_3 \rangle^{-2}
 \big( \gamma_{\leq N} - \Gamma_{\leq N}(n_3) \big) 
 e^{i\langle n_0 ,x \rangle} \label{para:eq-11X1-expl-p1} \\
 &\times 
 \varphi_1\big( t \langle n_1 \rangle\big) \varphi_2\big(t \langle n_2 \rangle \big)
 \Big( \int_0^t \dt^\prime \sin\big( (t-t^\prime) \langle n_{3} \rangle\big) \varphi_3\big( t\langle n_3 \rangle\big) \Big) \SI[n_j,\varphi_j \colon 1 \leq j \leq 3] \bigg] \notag \allowdisplaybreaks[3] \\
 +& 1\big\{ N_2 = N_3 \big\} \sum_{\substack{\varphi_1 \in \\ \{ \cos, \sin \}}} 
 \sum_{\substack{n_0,n_1,n_2,n_3\in \Z^3\colon \\ n_0=n_{1}}}
 \bigg[ 1\big\{ n_{23}=0 \big\}  \Big( \prod_{j=0}^3 1_{N_j}(n_j) \Big) \langle n_1 \rangle^{-1} \langle n_2 \rangle^{-1} \langle n_3 \rangle^{-2} 
\label{para:eq-11X1-expl-p2}  \\
&\times   
\big( \gamma_{\leq N} - \Gamma_{\leq N}(n_3) \big) e^{i\langle n_0 ,x \rangle}  \varphi_1\big( t \langle n_1 \rangle\big)
 \Big( \int_0^t \dt^\prime \sin\big( (t-t^\prime) \langle n_{3} \rangle\big) \cos\big( (t-t^\prime) \langle n_3 \rangle\big) \Big) \SI[n_1,\varphi_1] \bigg] \notag \allowdisplaybreaks[3] \\
  +& 1\big\{ N_1 = N_3 \big\} \sum_{\substack{\varphi_2 \in \\ \{ \cos, \sin \}}} 
 \sum_{\substack{n_0,n_1,n_2,n_3\in \Z^3\colon \\ n_0=n_{2}}}
 \bigg[ 1\big\{ n_{13}=0 \big\}  \Big( \prod_{j=0}^3 1_{N_j}(n_j) \Big) \langle n_1 \rangle^{-1} \langle n_2 \rangle^{-1} \langle n_3 \rangle^{-2} 
\label{para:eq-11X1-expl-p3}  \\
&\times   
\big( \gamma_{\leq N} - \Gamma_{\leq N}(n_3) \big) e^{i\langle n_0 ,x \rangle}  \varphi_2\big( t \langle n_2 \rangle\big)
 \Big( \int_0^t \dt^\prime \sin\big( (t-t^\prime) \langle n_{3} \rangle\big) \cos\big( (t-t^\prime) \langle n_3 \rangle\big) \Big) \SI[n_2,\varphi_2] \bigg] \notag.
\end{align}
By symmetry in $n_1$ and $n_2$, it suffices to treat the non-resonant component \eqref{para:eq-11X1-expl-p1} and the resonant component \eqref{para:eq-11X1-expl-p2}. Together with Lemma \ref{analytic:lem-difference-gamma}, the non-resonant component \eqref{para:eq-11X1-expl-p1} can be treated exactly as in \ref{para:item-11X1-hll}. For the resonant component \eqref{para:eq-11X1-expl-p2}, we argue similarly as in the proof of Lemma \ref{para:lem-resonant-X1}. More precisely, we obtain that 
\begin{align}
\eqref{para:eq-11X1-expl-p2} &= 
1\big\{ N_0 = N_1 \big\} 1\big\{ N_2=N_3 \big\}
 \sum_{\substack{\varphi_1 \in \\ \{ \cos, \sin \}}} 
 \sum_{\substack{n_1,n_3\in \Z^3\colon}}
 \bigg[ \Big( \prod_{j=1,3} 1_{N_j}(n_j) \Big) \langle n_1 \rangle^{-1}  \langle n_3 \rangle^{-3} \notag   \\
&\times   
\big( \gamma_{\leq N} - \Gamma_{\leq N}(n_3) \big) e^{i\langle n_1 ,x \rangle}  \varphi_1\big( t \langle n_1 \rangle\big)
 \Big( \int_0^t \dt^\prime \sin\big( (t-t^\prime) \langle n_{3} \rangle\big) \cos\big( (t-t^\prime) \langle n_3 \rangle\big) \Big) \SI[n_1,\varphi_1] \bigg] \notag \\
 &= -\frac{1}{4}\cdot 
1\big\{ N_0 = N_1 \big\} 1\big\{ N_2=N_3 \big\}
 \sum_{\substack{\varphi_1 \in \\ \{ \cos, \sin \}}} 
 \sum_{\substack{n_1,n_3\in \Z^3\colon}}
 \bigg[ \Big( \prod_{j=1,3} 1_{N_j}(n_j) \Big) \langle n_1 \rangle^{-1}  \langle n_3 \rangle^{-4} \notag   \\
&\times   
\big( \gamma_{\leq N} - \Gamma_{\leq N}(n_3) \big) e^{i\langle n_1 ,x \rangle}  \varphi_1\big( t \langle n_1 \rangle\big)
\Big( \cos\big( 2t \langle n_3 \rangle\big) -1 \Big) \SI[n_1,\varphi_1] \bigg]. \label{para:eq-11X1-expl-p4}
\end{align}
Using Gaussian hypercontractivity and Lemma \ref{analytic:lem-difference-gamma}, it follows that 
\begin{align*}
&\E \Big[ \big\| \eqref{para:eq-11X1-expl-p4} \big\|_{X^{-1/2+\delta_2,b_+-1}}^p \Big]^{2/p} \\
\lesssim&\, 1\big\{ N_0 = N_1 \big\} 1\big\{ N_2=N_3 \big\}
\sum_{n_1 \in \Z^3} 1_{N_1}(n_1) \langle n_1 \rangle^{-3+2\delta_2}
\Big( \sum_{n_3 \in \Z^3} 1_{N_3}(n_3) \langle n_3 \rangle^{-4} \big| \gamma_{\leq N}- \Gamma_{\leq N}(n_3)\big| \Big)^2 \\
\lesssim& 1\big\{ N_0 = N_1 \big\} 1\big\{ N_2=N_3 \big\} N_1^{2\delta_2} N_3^{-2}. 
\end{align*}
Due to our frequency-scale assumptions, this yields an acceptable contribution. 
\end{proof}

It remains to treat \eqref{para:eq-motivation-11X2}, i.e., the term involving $\squadratic$ and $\SXXtwo$. Aside from frequency-localizations, we further split \eqref{para:eq-motivation-11X2} into three further sub-terms, which correspond to the zero, one, and two-pairing case. The contribution of \eqref{para:eq-motivation-11X2} in the two-pairing case is partially but not completely cancelled by the renormalization multiplier $\Gamma$. In order to represent the remaining error, it is convenient to make the following definition.

\begin{definition}[Frequency-localized operator version of $\Gamma_{\leq N}$]\label{para:def-gamma-operator}
For all frequency-scales $K_0,K_1,K_2,K_3$ and $w\colon \R \times \T^3\rightarrow \R$,
we define
\begin{equation}\label{para:def-gamma-operator-eqn}
\begin{aligned}
\Gamma^{\textup{op}}[K_\ast](w) 
:=& 18 \sum_{\substack{k_0,k_1,k_2,k_3 \in \Z^3 \colon \\ k_0 = k_{123} }} 
\bigg[ \Big( \prod_{j=0}^3 1_{K_j}(k_j) \Big) e^{i \langle k_0 ,x \rangle} \\
&\times 
\int_0^t \dt^\prime \frac{\sin\big((t-t^\prime) \langle k_3 \rangle\big)}{\langle k_3 \rangle} \Big( \prod_{j=1}^2 \frac{\cos\big((t-t^\prime) \langle k_j \rangle\big)}{\langle k_j \rangle^2} \Big) \widehat{w}(t^\prime,k_0) \bigg]. 
\end{aligned}
\end{equation}
\end{definition}

\begin{remark}
The choice of using $K_j$ instead of $N_j$ to denote frequency-scales is deliberate, since otherwise the notation in Lemma \ref{para:lem-11X2} would conflict with the notation in the proof of Proposition \ref{analytic:prop-higher-order2} below. 
\end{remark}

Equipped with Definition \ref{para:def-gamma-operator}, we now obtain the following decomposition of the $\squadratic[\leqN] \XXtwo$-term.

\begin{lemma}[\protect{Decomposition of the $\squadratic[\leqN] \XXtwo$-term}]
\label{para:lem-11X2-decomposition}
For all $N\geq 1$, it holds that 
\begin{align}
&3 \, \squadratic[\leqN] \XXtwo - \Big( 6 \HLL + 3 \HHL \Big) \Big( \slinear[blue][\leqN], \slinear[blue][\leqN], \XXtwo \Big) + \Gamma_{\leq N} \Big( 3 \squintic[\leqN] - P_{\leq N} \slinear[green][\leqM]  + v_{\leq N} \Big)  \notag \\
=& \sum_{w_4} \sum_{\substack{N_0,\hdots,N_5,N_{234}\leq N\colon \\ 
N_{234}> \max(N_1,N_5)^\eta, \\ 
N_4 \leq \max(N_2,N_3)^\eta < \min(N_2,N_3) 
}}
P_{N_0} \bigg[ 3 \lcol \slinear[blue][N_1] \slinear[blue][N_5] \rcol \XXop [N_\ast,w_4]  
\label{para:eq-11X2-decomposition-1} \\
&\hspace{30ex} 
+ 1\big\{ N_0 = N_4, ~ N_1= N_2, ~ N_3=N_5 \big\} \, \Gamma^{\textup{op}}[N_4,N_2,N_3,N_{234}] w_4 \bigg]  \notag \\
+& \sum_{w_4} \bigg[ \Gamma_{\leq N} w_4 - 
\sum_{ \substack{ N_2,N_3,N_4,N_{234}\leq N: \\ 
N_{234} > \max(N_2,N_3)^\eta, \\
\min(N_2,N_3) > \max(N_2,N_3)^\eta, \\
N_4 \leq \max(N_2,N_3)^\eta.
}} \Gamma^{\textup{op}}[N_4,N_2,N_3,N_{234}] \, w_4\, \bigg], \label{para:eq-11X2-decomposition-2}
\end{align}
where the sum in $w_4$ is taken over $3 \squintic[\leqN]$, $-P_{\leq N} \slinear[green][\leqM]$, and $v_{\leq N}$. 
\end{lemma}

The purpose of the $\Gamma^{\textup{op}}$-terms in \eqref{para:eq-11X2-decomposition-1} is to exactly cancel the two-pairing terms. 

\begin{proof}
This follows directly from Lemma \ref{para:lem-decomposition} and the definitions of $\HLL$ and $\HHL$.
\end{proof}

We now separately treat the contributions of \eqref{para:eq-11X2-decomposition-1} and \eqref{para:eq-11X2-decomposition-2}. We start with \eqref{para:eq-11X2-decomposition-1}, which contains zero and one-pairing terms. 

\begin{lemma}[\protect{Estimate of the zero and one-pairing parts of $\squadratic[\leqN] \XXtwo$}]\label{para:lem-11X2}
For all $T\geq 1$ and $p\geq 2$, it holds that 
\begin{align*}
&\E \bigg[ \sup_N \sup_{\mathcal{J}} \Big\|
\sum_{\substack{N_0,\hdots,N_5,N_{234}\leq N\colon \\ 
N_{234}> \max(N_1,N_5)^\eta, \\ 
N_4 \leq \max(N_2,N_3)^\eta < \min(N_2,N_3) 
}}
P_{N_0} \bigg[ 3 \lcol \slinear[blue][N_1] \slinear[blue][N_5] \rcol \XXop [N_\ast,w_4]  
 \\
&\hspace{12ex} 
+ 1\big\{ N_0 = N_4, ~ N_1= N_2, ~ N_3=N_5 \big\} \, \Gamma^{\textup{op}}[N_4,N_2,N_3,N_{234}] w_4 \bigg] 
\Big\|_{X^{-1,b}(\mathcal{J})\rightarrow X^{-1/2+\delta_2,b_+-1}(\mathcal{J})}^p \bigg]^{1/p}\\
&\lesssim\, p^2 T^\alpha.
\end{align*}
\end{lemma}

The only reason for summing over all admissible frequency-scales $N_0,\hdots,N_5,N_{234}$ in Lemma \ref{para:lem-11X2} (instead of fixing the frequency-scales) is that it preserves the symmetry in $n_1$ and $n_5$ as well as in $n_2$ and $n_3$.

\begin{proof}
We only treat the case $T=1$ and $\mathcal{J}=[-1,1]$, since the general case follows from minor modifications. 
We first decompose the product into resonant and non-resonant terms. It holds that
\begin{align}
     &\sum_{\substack{N_0,\hdots,N_5,N_{234}\leq N\colon \\ 
N_{234}> \max(N_1,N_5)^\eta, \\ 
N_4 \leq \max(N_2,N_3)^\eta < \min(N_2,N_3) 
}} P_{N_0} \Big[ 3 \lcol \slinear[blue][N_1] \slinear[blue][N_5] \rcol \XXop[N_\ast,w_4] \Big] \notag \\
     =& - 9 \sum_{\substack{N_0,\hdots,N_5,N_{234}\leq N\colon \\ 
N_{234}> \max(N_1,N_5)^\eta, \\ 
N_4 \leq \max(N_2,N_3)^\eta < \min(N_2,N_3) 
}}  P_{N_0}\bigg[ \lcol \slinear[blue][N_1] \slinear[blue][N_5] \rcol
     P_{N_{234}} \Big[ \lcol \slinear[blue][N_2] \slinear[blue][N_3] \rcol P_{N_4}w_4 \Big] \bigg] \notag \allowdisplaybreaks[3] \\
     =& -9 
\sum_{\substack{N_0,\hdots,N_5,N_{234}\leq N\colon \\ 
N_{234}> \max(N_1,N_5)^\eta, \\ 
N_4 \leq \max(N_2,N_3)^\eta < \min(N_2,N_3) 
}} 
\sum_{\substack{\varphi_1,\varphi_2,\varphi_3,\varphi_4 \\ \in \{ \cos, \sin \}}} \sum_{\substack{n_0,n_1, \hdots, n_5 \in \Z^3 \colon \\ n_0 =n_{12345}}}
     \bigg[ 1_{N_{234}}(n_{234}) \Big( \prod_{j=0}^5 1_{N_j}(n_j) \Big)   \label{para:eq-11X2-p1} \\
     &\times  \langle n_{234} \rangle^{-1} \langle n_1 \rangle^{-1}\langle n_2 \rangle^{-1}\langle n_3 \rangle^{-1}\langle n_5 \rangle^{-1}     e^{i \langle n_0, x \rangle}   \varphi_1\big(t \langle n_1 \rangle\big) \varphi_5\big( t\langle n_5 \rangle\big) 
    \notag \\
     &\times  \Big( \int_0^t \dt^\prime \sin\big((t-t^\prime) \langle n_{234} \rangle\big) 
     \varphi_2\big( t^\prime \langle n_2 \rangle\big) 
     \varphi_3\big(t^\prime \langle n_3 \rangle\big)
     \widehat{w}_4(t^\prime,n_4) \Big)   \SI[n_j,\varphi_j \colon j=1,2,3,5] \bigg]\notag \allowdisplaybreaks[3] \\
     -& \, 36 
      \sum_{\substack{N_0,\hdots,N_5,N_{234}\leq N\colon \\ 
N_{234}> \max(N_1,N_5)^\eta, \\ 
N_4 \leq \max(N_2,N_3)^\eta < \min(N_2,N_3) 
}} 
     \sum_{\substack{\varphi_3,\varphi_5 \in \\ \{ \cos,\sin \}}}
     \sum_{\substack{n_0,n_3,n_4,n_5 \in \Z^3 \colon \\ n_0=n_{345}}} \bigg[  1\big\{ N_1 = N_2 \big\}
      1_{N_0}(n_0)  \frac{1_{N_3}(n_3)}{\langle n_3 \rangle} 
      \label{para:eq-11X2-p2} \\
     &\times    1_{N_4}(n_4)  
      \frac{1_{N_5}(n_5)}{\langle n_5 \rangle}   e^{i \langle n_0, x \rangle}   \varphi_5\big(t \langle n_5 \rangle\big) 
     \Big( \int_0^t \dt^\prime \Sine[N_{234},N_2](t-t^\prime,n_{34})
     \varphi_3\big(t^\prime \langle n_3 \rangle\big)
     \widehat{w}_4(t^\prime,n_4) \Big)  \notag \\
     &\times \SI[n_3,\varphi_3,n_5,\varphi_5] \bigg] \allowdisplaybreaks[3] \notag \\
     -&\, 18  
      \sum_{\substack{N_0,\hdots,N_5,N_{234}\leq N\colon \\ 
N_{234}> \max(N_1,N_5)^\eta, \\ 
N_4 \leq \max(N_2,N_3)^\eta < \min(N_2,N_3) 
}} 
     \sum_{\substack{n_0,n_1,\hdots n_5 \in \Z^3 \colon \\ n_0=n_{12345}}} \bigg[ 1\big\{ n_{12} = n_{35}=0 \big\} 
     1_{N_{234}}(n_{234}) \Big( \prod_{j=0}^5 1_{N_j}(n_j) \Big) 
     e^{i \langle n_0 ,x \rangle} \label{para:eq-11X2-p3} \\
     &\times \Big( \int_0^t \dt^\prime 
     \frac{\sin\big( (t-t^\prime) \langle n_{234}\rangle\big)}{\langle n_{234}\rangle}
     \frac{\cos\big( (t-t^\prime) \langle n_2 \rangle\big)}{\langle n_2 \rangle^2}
     \frac{\cos\big( (t-t^\prime) \langle n_3 \rangle\big)}{\langle n_3 \rangle^2}
     \widehat{w}_4(t^\prime,n_4) \Big) \bigg]. \notag 
\end{align}
In \eqref{para:eq-11X2-p2} and \eqref{para:eq-11X2-p3}, we already used symmetry to reduce to pairings between $n_1$ and $n_2$ or between $n_1$ and $n_2$ as well as $n_3$ and $n_5$, respectively.  Since the two-pairing terms \eqref{para:eq-11X2-p3} are exactly cancelled by  the $P_{N_0} \Gamma^{\textup{op}}[N_4,N_2,N_3,N_{234}] w_4$-terms, it remains to estimate the zero-pairing terms in \eqref{para:eq-11X2-p1} and the one-pairing terms in \eqref{para:eq-11X2-p2}. Since we already used symmetry to reduce to fixed pairings, we can now estimate each dyadic term separately (instead of the whole sum over all admissible frequency scales). \\

\emph{The zero-pairing terms in \eqref{para:eq-11X2-p1}:} The argument is similar as for the non-resonant term in Lemma \ref{para:lem-11X1}. 
We use the quintic tensor from Lemma \ref{counting:quintic_tensor},
\begin{equation}\label{para:eq-11X2-p4}
h_{n_0 n_1 \hdots n_5}(t,\lambda_4) 
= h_{n_0 n_1 \hdots n_5}[N_0,\hdots,N_5,N_{234},\pm_1,\hdots,\pm_5](t,\lambda_4),
\end{equation}
where we set $\lambda_1=\lambda_2=\lambda_3=\lambda_5=0$. Using \eqref{para:eq-11X2-p4}, we can rewrite a dyadic component of the zero-pairing term as 
\begin{align}
\eqref{para:eq-11X2-p1}\big|_{N_\ast}
=& \sum_{\substack{\pm_1,\pm_2,\pm_3,\\ \pm_4,\pm_5}} \int_{\R} \dlambda_4
\sum_{\substack{n_0,n_1,n_2 \\ n_3,n_4,n_5}} \bigg[  e^{i \langle n_0 ,x \rangle} 
h_{n_0 n_1 \hdots n_5}(t,\lambda_4)  \notag  \\
&\times 
\widetilde{\langle \nabla \rangle w_4^{\pm_4}}(n_4,\lambda_4) \,
\SI[n_j,\pm_j\colon j=1,2,3,5] \bigg].\label{para:eq-11X2-p5}
\end{align}
The $\langle \nabla \rangle$-multiplier acting on $w_4$ is due to the $\langle n_4 \rangle^{-1}$-factor in the definition of the quintic tensor. Due to the frequency assumptions $N_4 \leq \max(N_2,N_3)^\eta$, however, the $\langle \nabla \rangle$-multiplier is essentially irrelevant. Using the tensor estimate reduction argument in Subsection \ref{prep:remark-reduction},  it follows that 
\begin{align}
&\Big\| \eqref{para:eq-11X2-p5} \Big\|_{X^{-1,b}\times X^{-1,b} \rightarrow X^{-1/2+\delta_2,b_+-1}} \notag \\
\lesssim& \, N_0^{-1/2+\delta_2}  N_4^2 
\max_{\pm_0,\pm_1,\hdots,\pm_5} \sup_{\lambda_4 \in \R} \bigg[ 
\langle \lambda_4 \rangle^{-(b_--1/2)} \notag \\
&\times \Big\| \langle \lambda \rangle^{b_+-1} 
\Big\| \sum_{\substack{n_1,n_2,n_3,n_5}} 
\widetilde{h}^{\pm_0}_{n_0 n_1 \hdots n_5}(\lambda,\lambda_4) 
\SI[n_j, \pm_j \colon j=1,2,3,5] \Big\|_{n_4 \rightarrow n_0} \Big\|_{L_\lambda^2} \bigg]. \label{para:eq-11X2-p6}
\end{align}
Since $N_4$ is small compared to $\max(N_2,N_3)$, we simply bound the $\| \cdot\|_{n_4 \rightarrow n_0}$-norm by the Hilbert-Schmidt norm $\| \cdot \|_{n_0n_4}$. Using  the $p$-moment estimate reduction argument in Subsection \ref{prep:remark-reduction}  and the Hilbert-Schmidt estimate \eqref{counting:quintic_bound1} in Lemma \ref{counting:quintic_tensor}, it follows that 
\begin{align*}
\E\Big[ \eqref{para:eq-11X2-p6}^p \Big]^{1/p}
\lesssim& \, p^{2} N_{\textup{max}}^\epsilon N_0^{-1/2+\delta_2} N_4^2 \\
&\times
N_0 \min(N_2,N_3,N_4)^{1/2} \max(N_2,N_3,N_4)^{-1/2} \max(N_0,N_1,N_5)^{-1/2} \\
\lesssim& \, p^2 N_{\textup{max}}^{\delta_2+\epsilon} N_4^{5/2} \max(N_2,N_3,N_4)^{-1/2}. 
\end{align*}
Due to our frequency-scale assumptions $\max(N_2,N_3,N_4)\gtrsim \max(N_1,N_5)^\eta$ and $N_4 \leq \max(N_2,N_3)^\eta$, this yields an acceptable contribution. 
~\\ 

\emph{The one-pairing term \eqref{para:eq-11X2-p2}:} 
 The argument is similar as for the resonant term in Lemma \ref{para:lem-11X1} and utilizes the sine-cancellation tensor. 
To be precise, we use the sine-cancellation tensor
\begin{equation}\label{para:eq-11X2-p7}
h^{\textup{\sine}}_{n_0n_3n_4n_5}(t,\lambda_4)
= h^{\textup{\sine}}_{n_0n_3n_4n_5}[N_0,N_1,\hdots,N_5,N_{234},\pm_3,\pm_4,\pm_5](t,\lambda_4)
\end{equation}
from Lemma \ref{counting:sine_tensor}, where we set $\lambda_3=\lambda_5=0$. Using \eqref{para:eq-11X2-p7}, we can rewrite a dyadic component of the one-pairing term as 
\begin{align}
\eqref{para:eq-11X2-p2}\big|_{N_\ast} &= 
\sum_{\pm_3,\pm_4,\pm_5} \int_{\R} \dlambda_4 \sum_{n_0,n_3,n_4,n_5} \bigg[ e^{i\langle n_0 ,x \rangle} h^{\textup{\sine}}_{n_0n_3n_4n_5}(t,\lambda_4) \notag \\
&\times 
\widetilde{\langle \nabla \rangle w_4^{\pm_4}}(n_4,\lambda_4) \,
\SI[n_3,\pm_3,n_5,\pm_5] \bigg]. \label{para:eq-11X2-p8}
\end{align}
As before, the $\langle \nabla \rangle$-multiplier results from the definition of the sine-cancellation tensor, but is essentially irrelevant. Using  the tensor estimate reduction argument in Subsection \ref{prep:remark-reduction}, \\
\begin{align}
&\Big\| \eqref{para:eq-11X2-p8} \Big\|_{X^{-1,b}\times X^{-1,b} \rightarrow X^{-1/2+\delta_2,b_+-1}} \notag \\
\lesssim& \, N_0^{-1/2+\delta_2}  N_4^2 
\max_{\pm_0,\pm_1,\hdots,\pm_5} \sup_{\lambda_4 \in \R} \bigg[ 
\langle \lambda_4 \rangle^{-(b_--1/2)} \notag \\
&\times \Big\| \langle \lambda \rangle^{b_+-1} 
\Big\| \sum_{\substack{n_3,n_5}} 
\widetilde{h}^{\textup{\sine},\pm_0}_{n_0 n_3 n_4 n_5}(\lambda,\lambda_3,\lambda_4) 
\SI[n_3,\pm_3,n_5,\pm_5] \Big\|_{n_3 n_4 \rightarrow n_0} \Big\|_{L_\lambda^2} \bigg]. \label{para:eq-11X2-p9}
\end{align}
As before, we estimate the operator norm $\| \cdot \|_{n_4 \rightarrow n_0}$ by the Hilbert-Schmidt norm $\| \cdot \|_{n_0 n_4}$. 
Using the $p$-moment estimate reduction argument in Subsection \ref{prep:remark-reduction}
and the Hilbert-Schmidt estimate \eqref{counting:sine_extra_HS} in Lemma \ref{counting:sine_tensor}, we obtain that 
\begin{align}
\E\Big[ \eqref{para:eq-11X2-p9}^p \Big]^{1/p} 
\lesssim& \, p N_{\textup{max}}^\epsilon N_0^{-1/2+\delta_2} (N_3 N_4)^2 \notag \\
&\times
\min(N_0,N_5)^{3/2} \max(N_2,N_{234})^{-1} N_3^{1/2} N_4^{1/2} N_5^{-1}  \notag \\
\lesssim& \, p N_{\textup{max}}^{\epsilon+\delta_2} \max(N_2,N_{234})^{-1} N_3^{1/2} N_4^{5/2}. \label{para:eq-11X2-p10}
\end{align}

Since $N_4 \leq \max(N_2,N_3)^\eta$, it holds that 
$\max(N_2,N_{234}) \sim \max(N_2,N_3,N_4)$. As a result, we obtain
\begin{equation*}
\eqref{para:eq-11X2-p10} \lesssim p N_{\textup{max}}^{\epsilon+\delta_2} \max(N_2,N_3,N_4)^{-1/2} N_4^{5/2}. 
\end{equation*}
Due to our frequency-scale assumptions, this is acceptable. 
\end{proof}

It remains to estimate \eqref{para:eq-11X2-decomposition-2}, which contains the renormalized two-pairing terms. In contrast to \eqref{para:eq-11X2-decomposition-1}, \eqref{para:eq-11X2-decomposition-2} viewed as an operator in $w_4$ is entirely deterministic. 

\begin{lemma}[\protect{Estimate of the renormalized two-pairing terms in $\squadratic[\leqN] \XXtwo$}]\label{para:lem-11X2-renormalized}
For all $N\geq 1$, $T\geq 1$, and closed intervals $0 \in \Jc \subseteq [-T,T]$, it holds that 
\begin{align}
\bigg\| \,     
\sum_{ \substack{ N_2,N_3,N_4,N_{234}\leq N: \\ 
N_{234} > \max(N_2,N_3)^\eta, \\
\min(N_2,N_3) > \max(N_2,N_3)^\eta, \\
N_4 \leq \max(N_2,N_3)^\eta.
}} \Gamma^{\textup{op}}[N_4,N_2,N_3,N_{234}] \, w_4\, 
- \Gamma_{\leq N} w_4 \bigg\|_{X^{-1/2+\delta_2,b_+-1}(\Jc)} 
\lesssim T^\alpha \| w_4 \|_{X^{-\epsilon,b}(\Jc)}. 
\end{align}
\end{lemma}

\begin{proof}
We only treat the case $T=1$ and $\Jc=[-1,1]$, since the general case follows from minor modifications. Before starting with our estimates, we first rewrite the $\Gamma^{\textup{op}}$-terms. From Definition \ref{para:def-gamma-operator}, it follows that 
\begin{align}
&\sum_{ \substack{ N_2,N_3,N_4,N_{234}\leq N: \\ 
N_{234} > \max(N_2,N_3)^\eta, \\
\min(N_2,N_3) > \max(N_2,N_3)^\eta, \\
N_4 \leq \max(N_2,N_3)^\eta.
}} \Gamma^{\textup{op}}[N_4,N_2,N_3,N_{234}] \, w_4\, - 
\Gamma_{\leq N} w_4 \notag \\
=& \sum_{ \substack{ N_2,N_3,N_4,N_{234}\leq N }} 
\Gamma^{\textup{op}}[N_4,N_2,N_3,N_{234}] \, w_4\, - 
\Gamma_{\leq N} w_4 \label{para:eq-11X2-renormalized-p1} \\
+& \sum_{ \substack{ N_2,N_3,N_4,N_{234}\leq N}} \bigg[
\Big( 1 - 1\big\{ N_{234}>\max(N_2,N_3)^\eta, ~N_4 \leq \max(N_2,N_3)^\eta < \min(N_2,N_3) \big\} \Big)  \label{para:eq-11X2-renormalized-p2}  \\
&\hspace{17ex} \times  \Gamma^{\textup{op}}[N_4,N_2,N_3,N_{234}] \, w_4\,\bigg]. \notag 
\end{align}
We now treat \eqref{para:eq-11X2-renormalized-p1} and \eqref{para:eq-11X2-renormalized-p2} separately. \\

\emph{The contribution of \eqref{para:eq-11X2-renormalized-p1}:}
By first using  Definition \ref{para:def-gamma-operator} and then using the symmetry in $-n_{234},n_2$, and $n_3$, it holds that 
\begin{align}
\eqref{para:eq-11X2-renormalized-p1} 
&= 18 \sum_{\substack{n_2,n_3,n_4,n_{234}\in \Z^3 \colon \\ n_{234} = n_2 + n_3 + n_4}} 
\bigg[ 1_{\leq N}(n_{234}) \Big( \prod_{j=2}^4 1_{\leq N}(n_j) \Big) e^{i\langle n_4 , x \rangle} 
\notag \\
&\times 
\int_0^t \dt^\prime \partial_{t^\prime} \bigg( \frac{\cos\big((t-t^\prime) \langle n_{234}  \rangle\big)}{\langle n_{234} \rangle^2} \bigg) \Big( \prod_{j=2}^3 \frac{\cos\big((t-t^\prime) \langle n_j \rangle\big)}{\langle n_j \rangle^2} \Big) \widehat{w}(t^\prime,n_4) \bigg] - \Gamma_{\leq N} w_4 \notag  \\
&= 6  \sum_{\substack{n_2,n_3,n_4,n_{234}\in \Z^3 \colon \\ n_{234} = n_2 + n_3 + n_4}} 
\bigg[ 1_{\leq N}(n_{234}) \Big( \prod_{j=2}^4 1_{\leq N}(n_j) \Big) e^{i\langle n_4 , x \rangle} 
\notag \\
&\times 
\int_0^t \dt^\prime \partial_{t^\prime} \bigg( \frac{\cos\big((t-t^\prime) \langle n_{234}  \rangle\big)}{\langle n_{234} \rangle^2}  \Big( \prod_{j=2}^3 \frac{\cos\big((t-t^\prime) \langle n_j \rangle\big)}{\langle n_j \rangle^2} \Big) \bigg) \widehat{w}(t^\prime,n_4) \bigg] - \Gamma_{\leq N} w_4 \notag \\
&= \sum_{n_4 \in \Z^3} \bigg[ e^{i\langle n_4 ,x \rangle} \bigg( 
\int_0^t \dt^\prime \partial_{t^\prime}\Big( \Gamma_{\leq N}(n_4,t-t^\prime) \Big)  \widehat{w}(t^\prime,n_4) - \Gamma_{\leq N}(n_4) \widehat{w}(t,n_4) \bigg) \bigg]. \label{para:eq-11X2-renormalized-p3}
\end{align}
Using integration by parts, it follows that 
\begin{align}
\eqref{para:eq-11X2-renormalized-p3} &= 
-\sum_{n_4\in\Z^3} e^{i\langle n_4,x\rangle} \Gamma_{\leq N}(n_4,t) \widehat{w}_4(0,n_4) \label{para:eq-11X2-renormalized-p4} \\
&- \sum_{n_4\in \Z^3} e^{i\langle n_4 , x \rangle} \int_0^t \dt^\prime \Gamma_{\leq N}(n_4,t-t^\prime) \partial_{t^\prime} \widehat{w}_4(t^\prime,n_4). \label{para:eq-11X2-renormalized-p5} 
\end{align}
The first term \eqref{para:eq-11X2-renormalized-p4} can easily be estimated using a dyadic decomposition and \eqref{analytic:eq-Gamma-dyadic-1} from Lemma \ref{analytic:lem-Gamma-dyadic}. In order to estimate  
\eqref{para:eq-11X2-renormalized-p5}, we first decompose $w_4$ into a superposition of time-modulated linear waves. By inserting the decomposition into \eqref{para:eq-11X2-renormalized-p5}, it follows that 
\begin{equation*}
\eqref{para:eq-11X2-renormalized-p5} = 
i \sum_{\pm_4} \sum_{n_4\in \Z^3} \bigg[ e^{i\langle n_4 ,x \rangle} \bigg( \int_{\R} \dlambda_4 \big( \pm_4 \langle n_4 \rangle + \lambda_4 \big) \Big( \int_0^t \dt^\prime \Gamma_{\leq N}(n_4,t-t^\prime) e^{i (\pm_4 \langle n_4 \rangle + \lambda_4)t^\prime} \Big) \widetilde{w}_4^{\pm_4}(\lambda_4,n_4) \bigg)\bigg].
\end{equation*}
Using \eqref{analytic:eq-Gamma-dyadic-1} from Lemma \ref{analytic:lem-Gamma-dyadic} and a dyadic decomposition, it follows that
\begin{equation*}
\bigg|  \big( \pm_4 \langle n_4 \rangle + \lambda_4 \big) \Big( \int_0^t \dt^\prime \Gamma_{\leq N}(n_4,t-t^\prime) e^{i (\pm_4 \langle n_4 \rangle + \lambda_4)t^\prime} \Big) \bigg| \lesssim \big( N_4 + \langle \lambda_4 \rangle\big)^{b_--1/2}. 
\end{equation*}
After using Cauchy-Schwarz in $\lambda_4$, it follows that
\begin{equation*}
  \big\| \eqref{para:eq-11X2-renormalized-p5}  \big\|_{L_t^\infty H_x^{-2\epsilon}} \lesssim \| w_4 \|_{X^{-\epsilon,b}}, 
\end{equation*}
which yields the desired estimate. 
\\

\emph{The contribution of \eqref{para:eq-11X2-renormalized-p2}:}
Using Definition \ref{para:def-gamma-operator}, we see that the dyadic components in \eqref{para:eq-11X2-renormalized-p2} are given by 
\begin{equation}\label{para:eq-11X2-renormalized-p6}
\begin{aligned}
\sum_{n_2,n_3,n_4 \in \Z^3} \bigg[& \,   1_{\leq N}(n_{234}) \Big( \prod_{j=2}^4 1_{\leq N}(n_j) \Big) e^{i\langle n_4 , x \rangle} 
\notag \\
&\times 
\int_0^t \dt^\prime\,   \frac{\sin\big((t-t^\prime) \langle n_{234}  \rangle\big)}{\langle n_{234} \rangle} \Big( \prod_{j=2}^3 \frac{\cos\big((t-t^\prime) \langle n_j \rangle\big)}{\langle n_j \rangle^2} \Big) \widehat{w}(t^\prime,n_4) \bigg]. 
\end{aligned}
\end{equation}
The sum over frequency-scales in \eqref{para:eq-11X2-renormalized-p2} is supported on a set of frequency-scales where at least one of the three conditions 
\begin{equation}\label{para:eq-11X2-renormalized-p7}
\min\big(N_2,N_3\big) \leq \max(N_2,N_3)^\eta, \quad N_4 > \max(N_2,N_3)^\eta, \quad \text{or} \quad N_{234} < \max(N_2,N_3)^\eta
\end{equation}
is satisfied. 
By symmetry, we may also assume that $N_2\geq N_3$. We now further separate the argument into two sub-cases. \\

\emph{Case 1: $\min(N_2,N_3)\leq \max(N_2,N_3)^\eta$ or $N_4 > \max(N_2,N_3)^\eta$.}
In this case, we make use of the sine-cancellation. Using Definition \ref{def:sine_kernel}, we write
\begin{equation}
\begin{aligned}
\eqref{para:eq-11X2-renormalized-p6} =
\sum_{n_3,n_4 \in \Z^3} \bigg[&\, \langle n_3 \rangle^{-2} 1_{N_3}(n_3) 1_{N_4}(n_4) e^{i\langle n_4 ,x \rangle}  \\
&\times \int_0^t \dt^\prime \Sine[N_{234},N_2](t-t^\prime,n_{34}) \cos\big( (t-t^\prime) \langle n_3 \rangle\big) \widehat{w}_4(t^\prime,n_4) \bigg]. 
\label{para:eq-11X2-renormalized-p8}
\end{aligned}
\end{equation}
Together with the embedding $X^{-1/2+\delta_2,0}\hookrightarrow X^{-1/2+\delta_2,b_+-1}$, we obtain that 
\begin{align}
\Big\| \eqref{para:eq-11X2-renormalized-p8} \Big\|_{X^{-1/2+\delta_2,b_+-1}}^2 
\lesssim&\, \Big\| \eqref{para:eq-11X2-renormalized-p8} \Big\|_{L_t^2 H_x^{-1/2+\delta_2}}^2
\notag \\
\lesssim&\, \sup_{t\in [-1,1]} \sum_{n_4 \in \Z^3} 1_{N_4}(n_4) \langle n_4 \rangle^{-1+2\delta_2}
\bigg| \sum_{n_3\in \Z^3} \bigg[ 1_{N_3}(n_3) \langle n_3 \rangle^{-2} \label{para:eq-11X2-renormalized-p9} \\
&\hspace{2ex} \times  \int_0^t \dt^\prime \Sine[N_{234},N_2](t-t^\prime,n_{34}) \cos\big( (t-t^\prime) \langle n_3 \rangle\big) \widehat{w}_4(t^\prime,n_4)\bigg] \bigg|^2. \notag 
\end{align}
By first using the decomposition of $w_4$ into time-modulated linear waves (similar as the reduction arguments in Subsection \ref{prep:remark-reduction}) and then applying Lemma \ref{counting:lem-Sine-estimate}, it follows that 
\begin{align}
\eqref{para:eq-11X2-renormalized-p9} &\lesssim \max(N_{234},N_2)^{-2+2\epsilon} N_3^2 \big\| P_{N_4}w_4 \big\|_{X^{-1/2+\delta_2,b}}^2 \notag \\
&\lesssim \max(N_{234},N_2)^{-2+2\epsilon} N_3^2 N_4^{-1+2\delta_2+2\epsilon}
\big\| P_{N_4}w_4 \big\|_{X^{-\epsilon,b}}^2. \label{para:eq-11X2-renormalized-p10}
\end{align}
If $\min(N_2,N_3)\leq \max(N_2,N_3)^\eta$, it follows that 
\begin{align*}
\eqref{para:eq-11X2-renormalized-p10} \lesssim N_2^{-2+2\eta +2 \epsilon}  N_4^{-1+2\delta_2+2\epsilon}\big\| P_{N_4}w_4 \big\|_{X^{-\epsilon,b}}^2,
\end{align*}
which is acceptable. Alternatively, if $N_4 > \max(N_2,N_3)^\eta$, it follows that 
\begin{equation*}
  \eqref{para:eq-11X2-renormalized-p10} \lesssim \max(N_2,N_3,N_4)^{2\epsilon}  N_4^{-1+2\delta_2+2\epsilon}\big\| P_{N_4}w_4 \big\|_{X^{-\epsilon,b}}^2 
  \lesssim \max(N_2,N_3,N_4)^{-\eta+2\delta_2 + 4\epsilon}\big\| P_{N_4}w_4 \big\|_{X^{-\epsilon,b}}^2 , 
\end{equation*}
which is also acceptable. We note that, while the case $N_4 \geq \max(N_2,N_3)^\eta$ has been conveniently treated using the $\Sine$-kernel, it actually does not require the sine-cancellation, because one can gain $N_4^{-1/3}$ from the $X^{-1/2+\delta_2,b_+-1}$ norm of the output, while the sum over $(n_2,n_3)$ can be controlled using Lemma \ref{counting:lem2}. \\

\emph{Case 2: $N_{234}\leq \max(N_2,N_3)^\eta$.}
In this case, we write the dyadic component of \eqref{para:eq-11X2-renormalized-p2} as 
\begin{equation}\label{para:eq-11X2-renormalized-p11}
\begin{aligned}
\sum_{n_4 \in \Z^3} \int_0^t \dt^\prime &\bigg( 1_{N_4}(n_4)  \widehat{w}(t^\prime,n_4)  e^{i\langle n_4, x \rangle} 
\sum_{n_2,n_3 \in \Z^3} \bigg[ 1_{N_{234}}(n_{234}) \Big( \prod_{j=2}^3 1_{N_j}(n_j) \Big) \\
&\times
\frac{\sin\big( (t-t^\prime) \langle n_{234} \rangle\big)}{\langle n_{234} \rangle} 
\Big( \prod_{j=2}^3 \frac{\cos\big((t-t^\prime) \langle n_j \rangle\big)}{\langle n_j \rangle^2}
\Big) \bigg]\bigg). 
\end{aligned}
\end{equation}
Using the assumption $N_{234}\leq \max(N_2,N_3)^\eta$ and $N_2 \geq N_3$, it follows that 
\begin{align*}
  &\bigg| \sum_{n_2,n_3 \in \Z^3} \bigg[ 1_{N_{234}}(n_{234}) \Big( \prod_{j=2}^3 1_{N_j}(n_j) \Big) 
\frac{\sin\big( (t-t^\prime) \langle n_{234} \rangle\big)}{\langle n_{234} \rangle} 
\Big( \prod_{j=2}^3 \frac{\cos\big((t-t^\prime) \langle n_j \rangle\big)}{\langle n_j \rangle^2}
\Big) \bigg]\bigg| \\
\lesssim& \,  N_{234}^{-1} N_2^{-2} N_3^{-2} \times N_{234}^3 N_3^3 \lesssim N_2^{-1+2\eta}. 
\end{align*}
As a result, we obtain that 
\begin{align*}
\Big\| \eqref{para:eq-11X2-renormalized-p11} \Big\|_{X^{-1/2+\delta_2,b_+-1}}^2 
\lesssim&\, \Big\| \eqref{para:eq-11X2-renormalized-p11} \Big\|_{L_t^2 H_x^{-1/2+\delta_2}}^2
 \\
 \lesssim&\, N_2^{-2+4\eta} N_4^{-1+2\delta_2 +2\epsilon} \big\| P_{N_4}w_4 \big\|_{L_t^2 H_x^{-\epsilon}}^2 \\
\lesssim&\, N_2^{-2+4\eta} N_4^{-1+2\delta_2 +2\epsilon} \big\| P_{N_4}w_4 \big\|_{X^{-\epsilon,b}}^2,
\end{align*}
which is acceptable.

\end{proof}
\section{Analytic aspects of higher order stochastic diagrams} \label{section:analytic2}

In this section, we estimate the terms involving explicit stochastic objects which contribute to the remainder. In fact, some of them have already been treated in Sections \ref{section:bilinear}--\ref{section:para}, so here we only need to study the three (more difficult) remaining objects, namely the linear-cubic-cubic, the cubic-cubic-cubic, and the linear-linear-quintic interactions. We state the main estimates in the following two propositions. 

\begin{proposition}[Regularity of linear-cubic-cubic and cubic-cubic-cubic terms]\label{analytic:prop-higher-order}
Fix dyadic frequency-scales $K_0,K_1,K_2,K_3 \geq 1$, denote $K_{\mathrm{max}}:=\max(K_0,\cdots,K_3)$, let $p\geq 2$ and $T\geq 1$. Define the linear-cubic-cubic interaction
\[\Mc_{1,3,3}:=P_{K_0} \Pi_{\leq N}^\ast \big( P_{K_1}\slinear[blue][\leqN], P_{K_2}\scubic[\leqN], P_{K_3}\scubic[\leqN]\big)\] and the cubic-cubic-cubic interaction
\[\Mc_{3,3,3}:=P_{K_0} \Pi_{\leq N}^\ast \big( P_{K_1}\scubic[\leqN], P_{K_2}\scubic[\leqN], P_{K_3}\scubic[\leqN]\big),\] then we have the estimates
\begin{equation}\label{analytic:eq-higher-order}
\begin{aligned}
\| \Mc_{1,3,3}\|_{L^p_\omega X^{-1/2+\delta_2,b_+-1}([-T,T])} &\lesssim p^{15/2} T^\alpha (K_{\mathrm{max}})^{-\eta^2},\quad \mathrm{if }\,\max(K_2,K_3) \gtrsim (K_\mathrm{max})^\eta,
\\
\| \Mc_{3,3,3} \|_{L^p_\omega X^{-1/2+\delta_2,b_+-1}([-T,T])} &\lesssim p^{15/2} T^\alpha (K_{\mathrm{max}})^{-\eta^2}.
\end{aligned}
\end{equation}
\end{proposition}
\begin{proposition} [Regularity of the linear-linear-quintic term]\label{analytic:prop-higher-order2}Fix a dyadic frequency scale $N$, let $p\geq 2$ and $T\geq 1$. Define the linear-quintic-quintic interaction
\[\Mc_{1,1,5}:=P_{\leq N} \bigg[ 9 \squadratic[\leqN] \squintic[\leqN] - \Gamma_{\leq N} \cdot\scubic[\leqN] - 18 \mathfrak{C}^{(1,5)}_{\leq N} P_{\leq N} \slinear[blue][\leqN]  
   - 9 \Big( 2 \HLL + \HHL \Big) \Big( \, \slinear[blue][\leqN], \slinear[blue][\leqN], \squintic[\leqN] \Big) \bigg],\] then we have the estimate
\begin{equation}\label{analytic:eq-higher-order2}
            \| \Mc_{1,1,5}\|_{L_\omega^pX^{-1/2+\delta_2,b_+-1}([-T,T])}\lesssim p^{15/2}T^{\alpha}.
\end{equation}
\end{proposition}

In Subsection \ref{1+3+3} we will prove Proposition \ref{analytic:prop-higher-order} for $\Mc_{1,3,3}$, and in Subsection \ref{3+3+3} we will prove Proposition \ref{analytic:prop-higher-order} for $\Mc_{3,3,3}$. In Subsection \ref{1+1+5} we will study $\Mc_{1,1,5}$ and prove Proposition \ref{analytic:prop-higher-order2}. These objects are much more complicated than the quintic and lower order ones discussed in Section \ref{section:analytic}. While it is possible to treat them using only the graphical notation from Section \ref{section:diagrams} and tensor estimates, such an argument could be extremely lengthy. To be more efficient, we instead use the molecule formalism from \cite{DH21}. This method will be described in detail in Subsection \ref{1+3+3} below, and is also applied in Subsections \ref{3+3+3}--\ref{1+1+5}. Note that the estimate for $\Mc_{1,1,5}$ is stated and proved separately from the other two terms, because it has more complicated cancellation and renormalization structures, which will also be treated in Subsection \ref{1+1+5}.

\subsection{The linear-cubic-cubic stochastic object}\label{1+3+3} In this subsection we prove Proposition \ref{analytic:prop-higher-order} for $\Mc_{1,3,3}$.

\underline{\emph{Part I: reduction to counting estimates}}. We first reduce the estimate for $\Mc_{1,3,3}$ to a counting problem associated with the structure of this term.

\begin{figure}[h!]
  \includegraphics[scale=.45]{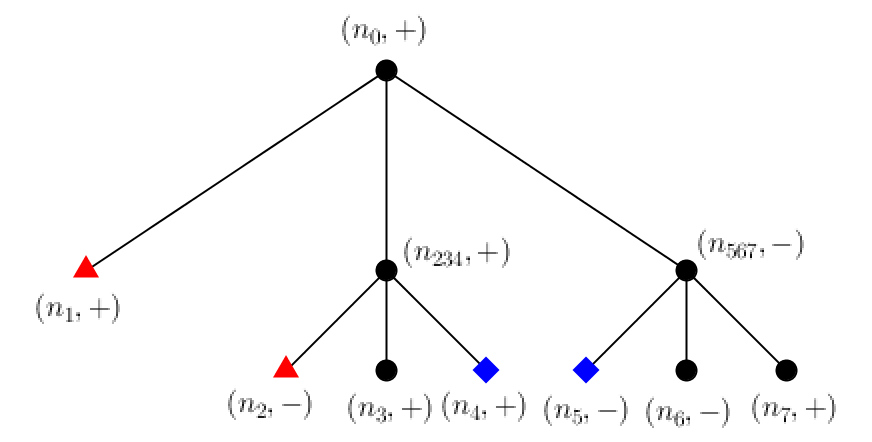}
  \caption{The diagram associated with $\Mc_{1,3,3}$ with the corresponding Fourier modes $n_j$ indicated. For each $j$, the sign $\pm_j$ is defined below, see (\ref{analytic:sept-coef4})); here we have made a specific choice of each, as written beside each $n_j$ (for example $\pm_{234}=+$ and $\pm_5=-$). Nodes of the same color and shape (other than black circled dot) are assumed to be paired in $\Pc$, for example here we have $\Pc=\{\{1,2\},\{4,5\}\}$.}
  \label{analytic:fig:sept-tree}
\end{figure}

Note that the definition \eqref{ansatz:eq-cubic} of $\scubic[\leqN]$ involves a time integral, which is restricted to $[-T,T]$ due to the time localization in (\ref{analytic:eq-higher-order}); by subdividing $[-T,T]$ into intervals of length $\sim 1$, and possibly losing factors of $T^\alpha$ in this process, we may reduce to the case $T=1$ (the estimates used below are invariant under time translations). This also means we can freely insert smooth cutoff functions $\chi=\chi(t)$ in the arguments below. Using Gaussian hypercontractivity (Lemma \ref{prelim:lem-hypercontractivity}) and Minkowski inequality as in the proof of Lemma \ref{analytic:lem-linear}, we may also reduce to the case $p=2$. Denote the wave number of the output $\Mc_{1,3,3}$ by $n_0$, and the wave numbers of the linear and two cubic inputs in $\Mc_{1,3,3}$ by $n_1$, $(n_2,n_3,n_4)$ and $(n_5,n_6,n_7)$ respectively; this is illustrated in Figure \ref{analytic:fig:sept-tree}, in which each node and corresponding Fourier mode, as well as some additional features (see below) are indicated.

Note that Figure \ref{analytic:fig:sept-tree} is drawn in the style of \cite{DH21} (except that we identify pairings by both color and shape); in the notation of the current paper it would correspond to the term \[\septicex[\leqN].\]

Without loss of generality, we may restrict to the region $\rho\geq 0$ in (\ref{defofxsb}), below we will replace the letter $\rho$ by $\xi$, so the $X^{-1/2+\delta_2,b_+-1}$ norm of $\Mc_{1,3,3}$ will contain the weight $\langle \xi-\langle n_0\rangle\rangle^{b_+-1}$, see (\ref{defofxsb}). Then, like in Section \ref{section:diagrams}, we can decompose
\begin{equation}\label{analytic:eq-septic2}
\Fc\Big(P_{K_1}\slinear[blue][\leqN]\cdot P_{K_2}\scubic[\leqN]\cdot P_{K_3}\scubic[\leqN]\Big)(t,n_0)=\sum_\Pc\sum_{\varphi_1,\cdots,\varphi_7\in\{\sin,\cos\}}\sum_{n_1,\cdots,n_7}(\Bc_{1,3,3})_{n_0\cdots n_7}(t)\cdot\Sc\Ic[n_j,\varphi_j:j\in O],
\end{equation} with the coefficient defined as
\begin{equation}\label{analytic:sept-coef}
\begin{aligned}(\Bc_{1,3,3})_{n_0\cdots n_7}(t)&=\prod_{j=1}^7\frac{1}{\langle n_j\rangle}\cdot\frac{1}{\langle n_{234}\rangle}\frac{1}{\langle n_{567}\rangle}\int_0^t\int_0^t\chi(t)\chi(t_1)\chi(t_2)\sin(\langle n_{234}\rangle (t-t_1))\\&\times\sin(\langle n_{567}\rangle (t-t_2))\varphi_1(\langle n_1\rangle t)\prod_{j=2}^4\varphi_j(\langle n_j\rangle t_1)\prod_{j=5}^7\varphi_j(\langle n_j\rangle t_2)\,\mathrm{d}t_1\mathrm{d}t_2.
\end{aligned}
\end{equation} Here in (\ref{analytic:eq-septic2}), $\chi$ is a smooth cutoff function, the sum in $\Pc$ is taken over all possible collections of \emph{pairings}, i.e. disjoint two-element subsets of $\{1,\cdots,7\}$, such that $\Pc$ does not contain any pairing within $\{2,3,4\}$ or $\{5,6,7\}$ (for example $\{2,3\}\not\in\Pc$; note that those $\Pc$ containing $\{2,3\}$ etc. are excluded due to the Wick ordering in the definition of $\scubic[]$). In the summation over $(n_1,\cdots,n_7)$ and $(\varphi_1,\cdots,\varphi_7)$ in (\ref{analytic:eq-septic2}) we assume that
\begin{equation}\label{analytic:sept-assump1}n_0=n_{1\cdots 7};\quad n_i+n_j=0\,\,\mathrm{and}\,\,\varphi_i=\varphi_j\,\,\mathrm{whenever}\,\,\{i,j\}\in\Pc,\end{equation}as well as\begin{equation}\label{analytic:sept-assump2}K_1/2<|n_1|_\infty\leq K_1,\quad K_2/2<|n_{234}|_\infty\leq K_2,\quad K_3/2< |n_{567}|_\infty\leq K_3,\end{equation} and $|n_\nu|_\infty\leq N\,(\forall \nu\in\{1,\cdots,7,234,567\})$. As a convention, in this section we will use $j$ to indicate nonzero single-digit indices (in $\{1,\cdots,7\}$ in the current case), and $\nu$ to indicate arbitrary indices (for example including $0$, $234$ and $567$). The stochastic integral $\Sc\Ic$ is defined as in (\ref{diagram:eq-multiple-stochastic}), where $O$ is the set of indices in $\{1,\cdots,7\}$ that are not in any subset in $\Pc$.

Note that the renormalization term $\mathfrak{C}_{\leq N}^{(3,3)}[K_2,K_3](t)\cdot P_{K_1}\slinear[blue][\leqN]$ (see Definition \ref{ansatz:def-dyadic-modified-product}) exactly corresponds to the cases where $\Pc$ consists of exactly three pairings between the sets $\{2,3,4\}$ and $\{5,6,7\}$ (for example $\Pc=\{\{2,5\},\{3,6\},\{4,7\}\}$). Thus, $\Mc_{1,3,3}$ can be decomposed in the same way as in (\ref{analytic:eq-septic2}), but restricted to $\Pc$ that does not contain any of the pairing within $\{2,3,4\}$ or $\{5,6,7\}$, and does not consist of exactly three pairings between $\{2,3,4\}$ and $\{5,6,7\}$. We call such $\Pc$ \emph{good}. From now on, we shall assume a good $\Pc$ is fixed.

Suppose $\Pc$ contains a subset $\{1,j\}$ where $2\leq j\leq 7$; by symmetry we may assume $j=2$. Then we must have $n_1+n_2=0$ and $\varphi_1=\varphi_2$ in (\ref{analytic:sept-coef}), which means that here we have a similar structure as the sine-cancellation kernel in Lemma \ref{counting:lem-Sine-symmetrization}. Therefore, we can perform the same reduction in Lemma \ref{counting:lem-Sine-symmetrization} and replace $\Bc_{1,3,3}$ in (\ref{analytic:sept-coef}) by $\Bc_{1,3,3}^{(\mathrm{s})}$, which is the superposition of three terms defined in (\ref{analytic:sept-coef2}) below.

With the above arguments, we may now fix a choice of good $\Pc$ and a choice of $\varphi_j\in\{\sin,\cos\}$ for $j\in O$, and consider the expression on the right hand side of (\ref{analytic:eq-septic2}) with summation over $n_j$ and $\varphi_i=\varphi_j$ for $\{i,j\}\in\Pc$ only, where the coefficient $\Bc_{1,3,3}$ is defined as in (\ref{analytic:sept-coef}). Moreover, if $\{1,2\}\in\Pc$ then we only sum over $\varphi_i=\varphi_j$ for $\{1,2\}\neq\{i,j\}\in\Pc$ and (\ref{analytic:sept-coef}) is replaced by
\begin{equation}\label{analytic:sept-coef2}
\begin{aligned}(\Bc_{1,3,3}^{(\mathrm{s})})_{n_0\cdots n_7}(t)&=\prod_{j=3}^7\frac{1}{\langle n_j\rangle}\cdot\frac{1}{\langle n_{567}\rangle}\int_0^t\int_0^t\chi(t)\chi(t_1)\chi(t_2) \\&\times\bigg[1_{K_1}(n_1)1_{K_2}(n_{234})\frac{\langle n_1\rangle-\langle n_{234}\rangle}{\langle n_1\rangle^2\langle n_{234}\rangle^2}\sin((t-t_1)(\langle n_1\rangle-\langle n_{234}\rangle))\\&+1_{K_1}(n_1)1_{K_2}(n_{234})\frac{\langle n_1\rangle+\langle n_{234}\rangle}{\langle n_1\rangle^2\langle n_{234}\rangle^2}\sin((t-t_1)(\langle n_1\rangle+\langle n_{234}\rangle))\\&-\frac{1_{K_1}(n_1)1_{K_2}(n_{234})-1_{K_2}(n_1)1_{K_1}(n_{234})}{2}\cdot\frac{\cos((t-t')\langle n_1\rangle)\sin((t-t')\langle n_{234}\rangle)}{\langle n_1\rangle^2\langle n_{234}\rangle}\bigg]\\&\times\sin(\langle n_{567}\rangle (t-t_2))\prod_{j=3}^4\varphi_j(\langle n_j\rangle t_1)\prod_{j=5}^7\varphi_j(\langle n_j\rangle t_2)\,\mathrm{d}t_1\mathrm{d}t_2.
\end{aligned}
\end{equation}

Next, we may apply the formulas $\sin z=(e^{iz}-e^{-iz})/(2i)$ and $\cos z=(e^{iz}+e^{-iz})/2$, to reduce the time integral in (\ref{analytic:sept-coef}) to a linear combination of integrals of form
\begin{equation}\label{analytic:sept-coef3}e^{i\langle n_0\rangle t}\int_0^t\int_0^t \chi(t)\chi(t_1)\chi(t_2)e^{i(\Omega_0t+\Omega_1t_1+\Omega_2t_2)}\,\mathrm{d}t_1\mathrm{d}t_2:=e^{i\langle n_0\rangle t}\cdot\Hc(t,\Omega_0,\Omega_1,\Omega_2),
\end{equation} where
\begin{equation}\label{analytic:sept-coef4}
\begin{aligned}
\Omega_0&:=-\langle n_0\rangle\pm_1\langle n_1\rangle\pm_{234}\langle n_{234}\rangle \pm_{567}\langle n_{567}\rangle,
\\\Omega_1&:=\mp_{234}\langle n_{234}\rangle\pm_2\langle n_2\rangle\pm_3\langle n_3\rangle\pm_4\langle n_4\rangle,\\\Omega_2&:=\mp_{567}\langle n_{567}\rangle\pm_5\langle n_5\rangle\pm_6\langle n_6\rangle\pm_7\langle n_7\rangle.
\end{aligned}
\end{equation} Here $\pm_\nu$ for $\nu\in\{1,\cdots,7,234,567\}$ are arbitrary signs and $\mp_\nu$ are the opposite signs; for convenience define also $\pm_0=+$. If $\{i,j\}\in\Pc$ then it is easy to see that the two occurrences of $\langle n_i\rangle=\langle n_j\rangle$ in the phase must be with opposite signs, thus $\pm_i=\mp_j$. As for (\ref{analytic:sept-coef2}) one has a similar expression; in addition, either $\pm_1=\pm_{234}$ and one has a factor $(\langle n_1\rangle+\langle n_{234}\rangle)/\langle n_1\rangle^2\langle n_{234}\rangle^2$ attached to this integral, or $\pm_1=\mp_{234}$ and one has a factor $(\langle n_1\rangle-\langle n_{234}\rangle)/\langle n_1\rangle^2\langle n_{234}\rangle^2$ attached to this integral, or one has the $\Gamma$-condition where $|n_1|_\infty\geq\Gamma\geq|n_{234}|_\infty$ or $|n_1|_\infty\leq\Gamma\leq|n_{234}|_\infty$ for some a fixed value $\Gamma$ (and one has a factor $\langle n_1\rangle^{-2}\langle n_{234}\rangle^{-1}$ attached to this integral).

For the function $\Hc(t,\Omega_0,\Omega_1,\Omega_2)$ occurring in (\ref{analytic:sept-coef3}), by Proposition \ref{counting:timeint} we have
\begin{equation}\label{analytic:sept-coefest}\int_{\Rb^3}(\langle\Omega_0\rangle\langle \Omega_1\rangle \langle\Omega_2\rangle)^{8(1/2-b_+)}|(\Fc_t\Hc)(\xi,\Omega_0,\Omega_1,\Omega_2)|\,\mathrm{d}\Omega_0\mathrm{d}\Omega_1\mathrm{d}\Omega_2\lesssim \langle \xi\rangle^{4(1/2-b_+)};
\end{equation} The same estimate also holds for all $\Omega_j$ derivatives of $\Hc$. Moreover, assume $\max_{1\leq j\leq 7}\langle n_j\rangle\sim K_\mathrm{max}^+$, then all $\langle \Omega_j\rangle\lesssim K_\mathrm{max}^+$ by  (\ref{analytic:sept-coef4}), so by taking supremum in each unit cube using Sobolev, we obtain that
\begin{equation}\label{analytic:sept-coefest2}|(\Fc_t\Hc)(\xi,\Omega_0,\Omega_1,\Omega_2)|\lesssim \widetilde{\Hc}(\xi,\lfloor \Omega_0\rfloor,\lfloor\Omega_1\rfloor,\lfloor\Omega_2\rfloor)
\end{equation} for some function $\widetilde{\Hc}$, such that
\begin{equation}\label{analytic:sept-coefest2+}\sum_{m_0,m_1,m_2}|\widetilde{\Hc}(\xi,m_0,m_1,m_2)|\lesssim (K_\mathrm{max}^+)^{O(b_+-1/2)}\langle \xi\rangle^{4(1/2-b_+)}.
\end{equation} Note also that $\overline{(\Fc_t\Hc)(\xi,\Omega_0,\Omega_1,\Omega_2)}=\Fc_t\Hc(-\xi,\-\Omega_0,-\Omega_1,-\Omega_2)$. If $K_\mathrm{max}\leq (K_\mathrm{max}^+)^{1/100}$, say $\langle n_2\rangle\sim K_\mathrm{max}^+$, then using the proof of Lemma \ref{analytic:lem-cubic} (in particular using (\ref{counting:cubic_bound3}) with $N_{123}$ and $N_{\mathrm{max}}$ replaced by $K_1$ and $K_{\mathrm{max}}^+$ respectively) we can easily bound the $L_t^\infty H^{3}$ norm of $P_{K_1}\scubic[\leqN]$ by $(K_\mathrm{max}^+)^{-1/10}$, and the $L_t^\infty H^{3}$ norm of the other two factors by $(K_\mathrm{max}^+)^{1/30}$, hence (\ref{analytic:eq-higher-order}) follows trivially. Thus, we may assume $K_\mathrm{max}\geq (K_\mathrm{max}^+)^{1/100}$, so $(K_\mathrm{max}^+)^{O(b_+-1/2)}=(K_\mathrm{max})^{O(b_+-1/2)}$ on the right hand side of (\ref{analytic:sept-coefest2+}), which is negligible in view of the $(K_\mathrm{max})^{-\eta^2}$ gain in (\ref{analytic:eq-higher-order}). 

At this point, we can reduce the estimate of $\Mc_{1,3,3}$ to a counting estimate. Recall that we may fix $\Pc$ and $\varphi_j\,(j\in O)$; further we may replace the time integrals in (\ref{analytic:sept-coef}) and (\ref{analytic:sept-coef2}) by integrals of form (\ref{analytic:sept-coef3}). Recall that we have fixed $p=2$. Since
\[\Eb\big[\Sc\Ic(n_j,\varphi_j:j\in O)\cdot\overline{\Sc\Ic(n_j',\varphi_j:j\in O)}\big]=0\] unless $(n_j':j\in O)$ is a permutation of $(n_j:j\in O)$, and $\Eb|\Sc\Ic(n_j,\varphi_j:j\in O)|^2\lesssim1$, we get that
\begin{multline}\label{analytic:sept-exp}\Eb\|\Mc_{1,3,3}\|_{X^{-1/2+\delta_2,b_+-1}}^2\lesssim\sum_{n_0}\langle n_0\rangle^{-1+2\delta_2}\int_\Rb\langle \xi\rangle^{2b_+-2}\,\mathrm{d}\xi\\\times\sum_{(n_j:j\in O)}\bigg|\sum_{(n_i=-n_j:\{i,j\}\in\Pc)}\alpha(n_0,\cdots,n_7)\Fc_t\Hc(\xi,\Omega_0,\Omega_1,\Omega_2)\bigg|^2,
\end{multline} where in the summation we still assume (\ref{analytic:sept-assump1})--(\ref{analytic:sept-assump2}), and $\alpha=\alpha(n_1,\cdots,n_7)$ is a function bounded by
\begin{equation}\label{analytic:sept-coefest3}|\alpha|\lesssim\prod_{j=1}^7\frac{1}{\langle n_j\rangle}\cdot\frac{1}{\langle n_{234}\rangle}\frac{1}{\langle n_{567}\rangle}
\end{equation} if $\{1,2\}\not\in\Pc$. If $\{1,2\}\in\Pc$, then either we have 
\begin{equation}\label{analytic:sept-coefest4}|\alpha|\lesssim\prod_{j=3}^7\frac{1}{\langle n_j\rangle}\cdot\frac{1}{\langle n_{567}\rangle}\cdot\left\{
\begin{aligned}
&\frac{\langle n_1\rangle+\langle n_{234}\rangle}{\langle n_1\rangle^2\langle n_{234}\rangle^2},&\textrm{if\ }\pm_1=\pm_{234}\textrm{\ in\ }(\ref{analytic:sept-coef4}),\\
&\frac{|\langle n_1\rangle-\langle n_{234}\rangle|}{\langle n_1\rangle^2\langle n_{234}\rangle^2},&\textrm{if\ }\pm_1=\mp_{234}\textrm{\ in\ }(\ref{analytic:sept-coef4}),
\end{aligned}
\right.
\end{equation} or we have (\ref{analytic:sept-coefest3}) but may further require the $\Gamma$-condition on $n_1$ and $n_{234}$ in the summation.

Since $\langle n_0\rangle\lesssim K_\mathrm{max}$ and $K_\mathrm{max}^+\leq (K_\mathrm{max})^{100}$, by losing at most $(K_\mathrm{max})^{O(\delta_2)}$ which is negligible in view of the $(K_\mathrm{max})^{-\eta^2}$ gain, we may replace $\langle n_0\rangle^{-1+2\delta_2}$ in (\ref{analytic:sept-exp}) by $\langle n_0\rangle^{-1}$, and also perform a dyadic decomposition and restrict that $\langle n_\nu\rangle\sim N_\nu\,(\nu\in\{0,1,\cdots,7,234,567\})$, and that $|n_1+n_{234}|\sim L$ if $\{1,2\}\in\Pc$. Note in particular that $N_0\sim K_0$, $N_1\sim K_1,\, N_{234}\sim K_2,\,N_{567}\sim K_3$. In either case, using (\ref{analytic:sept-coefest2}), we can estimate
\begin{equation*}
\begin{aligned}\Eb\|\Mc_{1,3,3}\|_{X^{1/2+\delta_2,b_+-1}}^2&\lesssim\sum_{n_0}\langle n_0\rangle^{-1}\int_\Rb\langle \xi\rangle^{2b_+-2}\,\mathrm{d}\xi\sum_{(n_j=n_j':j\in O)}\sum_{(n_i=-n_j:\{i,j\}\in\Pc)}\sum_{(n_i'=-n_j':\{i,j\}\in\Pc)}|\alpha(n_1,\cdots,n_7)|\\&\times|\alpha(n_1',\dots,n_7')|\cdot|\widetilde{\Hc}(\xi,\lfloor\Omega_0\rfloor,\lfloor\Omega_1\rfloor,\lfloor\Omega_2\rfloor)|\cdot|\widetilde{\Hc}(-\xi,\lfloor\Omega_0'\rfloor,\lfloor\Omega_1'\rfloor,\lfloor\Omega_2'\rfloor)|.
\end{aligned}
\end{equation*} Here the summation is taken over $(n_1,\cdots,n_7)$ and another vector $(n_1',\cdots,n_7')$, which satisfies the same equations, support requirements etc., as $(n_1,\cdots,n_7)$. Moreover, we require the equalities $n_j=n_j'$ for $j\in O$ and $n_i=-n_j$ and $n_i'=-n_j'$ for $(i,j)\in\Pc$, as indicated in the summation. The quantities $\Omega_j$ are defined from $(n_j)$ as in (\ref{analytic:sept-coef3}); $\Omega_j'$ are similarly defined from $(n_j')$, \emph{but we change $-\langle n_0\rangle$ to $\langle n_0\rangle$ and switch $\mp_\nu$ and $\pm_\nu$ in the definition of $\Omega_j'$.}

Now, using (\ref{analytic:sept-coefest3}), and the $\xi$ integrability of $\langle \xi\rangle^{(2b_+-2)+4(1/2-b_+)}$, we can bound
\begin{equation}\label{analytic:sept-counting}\Eb\|\Mc_{1,3,3}\|_{X^{1/2+\delta_2,b_+-1}}^2\lesssim\prod_{j=1}^7N_j^{-2}\cdot N_0^{-1}N_{234}^{-2}N_{567}^{-2}\cdot\Ac^{-2}\cdot\sup_{m_0,m_1,m_2,m_0',m_1',m_2'}(\#\Sigma),
\end{equation} where $\Sigma$ is the set defined by
\begin{multline}\label{analytic:sept-defset}\Sigma=\big\{(n_0,n_1,\cdots,n_7,n_1',\cdots,n_7'):n_0=n_{1\cdots 7}=n_{1\cdots 7}',\,\,n_i+n_j=n_i'+n_j'=0\,(\forall \{i,j\}\in\Pc),\\n_j=n_j'\,(\forall j\in O),\,\,
\Omega_j=m_j+O(1)\,(0\leq j\leq 2),\,\,\Omega_j'=m_j'+O(1)\,(0\leq j\leq 2)\big\},
\end{multline} note that in the definition of $\Sigma$ we also impose the above dyadic restrictions such as $\langle n_\nu\rangle\sim\langle n_\nu'\rangle\sim N_\nu$ and $|n_1+n_{234}|\sim  |n_1'+n_{234}'|\sim L$ etc. The quantity $\Ac$ equals $1$ if $\{1,2\}\not\in\Pc$, and if $\{1,2\}\in\Pc$ then either \begin{equation}\label{analytic:quantityA}
\Ac=\left\{
\begin{aligned}
&1,&\textrm{if\ }\pm_1=\pm_{234}\textrm{\ in\ }(\ref{analytic:sept-coef4}),\\
&N_{234}L^{-1},&\textrm{if\ }\pm_1=\mp_{234}\textrm{\ in\ }(\ref{analytic:sept-coef4})
\end{aligned}
\right.
\end{equation} (in view of (\ref{analytic:sept-coefest4})), or $\Ac=1$ but we we further require the $\Gamma$-condition for both $(n_j)$ and $(n_j')$ in (\ref{analytic:sept-defset}).

We have now reduced the estimate of $\Mc_{1,3,3}$ to the counting problem of estimating $\#\Sigma$ for fixed values of $(m_j,m_j')$, which we analyze in the next step using the notion of molecules.

\smallskip
\underline{\emph{Part II: reduction to molecules}}. The notion of \emph{molecules} is introduced in the recent work  of the second author with Z. Hani \cite{DH21}. With an algorithmic approach, it plays a fundamental role in the analysis of the diagrams occurring in \cite{DH21}, which are similar to Figure \ref{analytic:fig:sept-tree} but have arbitrarily high complexity.

Here, since we only need to deal with septic and nonic terms, our analysis will not be pure algorithmic, but rather an algorithmic-enumerative hybrid. In either case, we find it much easier (in terms of both calculations and presentations) than the traditional approach where one classifies all possibilities of pairings $\Pc$ and performs a case-by-case analysis.

\begin{figure}[h!]
  \includegraphics[scale=.4]{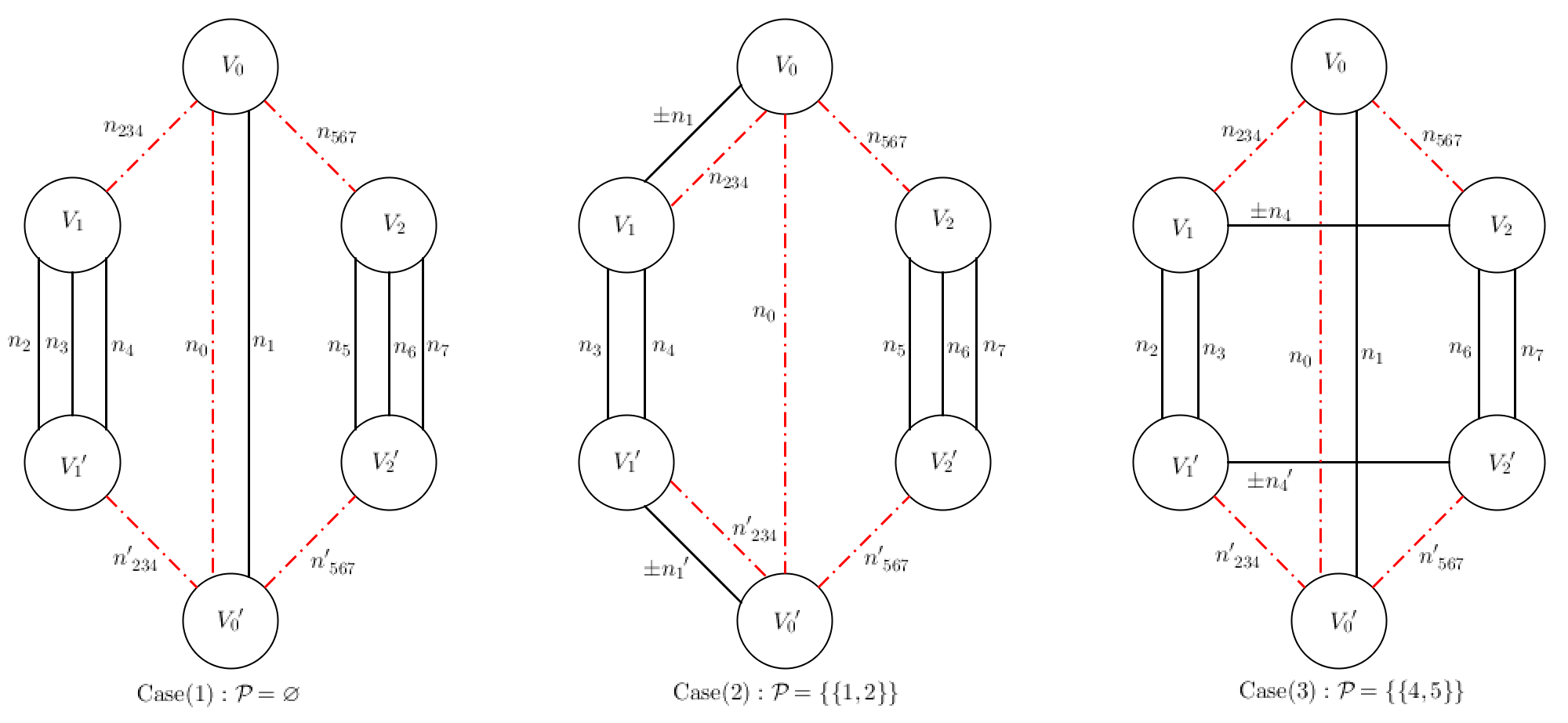}
  \caption{Cases (1)--(3) of the molecule $\Mb$ for the linear-cubic-cubic interaction. Here $V_0$ and $V_0'$ are the roots and $V_1,V_2,V_1',V_2'$ are their children. For each bond $e$, the value of $m_e$ is indicated as in Proposition \ref{analytic:prop-mol}; the directions of the bonds are not specified as in Remark \ref{analytic:rem1}. The $5$ bad bonds $e$ are those whose $m_e$ equals $\pm n_{234}, \pm n_{567}$, $\pm n_{234}',\pm n_{567}'$ or $n_0$, which has to do with the right hand side of (\ref{analytic:sept-counting}), see (\ref{analytic:molecount1}).}
  \label{analytic:fig:sept-mole1}
\end{figure}

The definition of molecules here is a special case of \cite{DH21}, although we will present it in a slightly different way in accordance with the notations chosen here.
\begin{definition}[Definition of molecules: linear-cubic-cubic case]\label{analytic:def-mole} Given a set $\Pc$ and the associated signs $\pm_\nu$ as above. Draw a ternary tree $\Tc$ such that $\Tc$ has root $0$ with three children $1$, $234$ and $567$, the nodes $234$ and $567$ have children $2,3,4$ and $5,6,7$ respectively, and nodes $1$--$7$ are all leaves (this is essentially the same tree in Figure \ref{analytic:fig:sept-tree}, which is the ternary tree corresponding to the linear-cubic-cubic interaction); suppose each node $\nu$ of the tree has sign $\pm_\nu$ (the root has sign $+$). Draw another identical tree $\Tc'$ whose nodes (denoted by $1'$, $234'$ etc.) have opposite signs with $\Tc$. We pair the leaves $(i,j)$ and $(i',j')$ for $\{i,j\}\in\Pc$, and also pair the leaves $(j,j')$ for $j\in O$; note that the signs of any two paired leaves must be opposite.

Now define a \emph{molecule} $\Mb$, which is a directed graph, as follows. Its vertices (called \emph{atoms}) are the branching nodes $V_j,V_j'\,(0\leq j\leq 2)$ of the two trees $\Tc$ and $\Tc'$. We connect two atoms $V$ and $\widetilde{V}$ by an directed edge $e$ (called \emph{bond}) if and only if (1) $\widetilde{V}$ is a child of $V$, in this case we let $e$ go from $V$ to $\widetilde{V}$ if $\widetilde{V}$ has $-$ sign, and from $\widetilde{V}$ to $V$ otherwise; (2) $V$ has a child leaf paired to a child leaf of $\widetilde{V}$, in this case we let $e$ go from the one whose child has $-$ sign to the one whose child has $+$ sign; (3) $V$ and $\widetilde{V}$ are the roots of $\Tc$ and $\Tc'$ respectively, in this case we let $e$ go from $V$ to $\widetilde{V}$. Note that multiple edges may be connected between two atoms. For example, in Figure \ref{analytic:fig:sept-tree} we have $\Pc=\{\{1,2\},\{4,5\}\}$ with the values of $\pm_i$ fixed as shown, then the corresponding molecule $\Mb$ is drawn in Case (4) of Figure \ref{analytic:fig:sept-mole2} below; here the bonds marked by $n_{234}$, $n_6$ and $n_0$ (for example) are constructed through cases (1)--(3) above respectively.
\end{definition}
\begin{proposition}\label{analytic:prop-mol} In the molecule $\Mb$ each atom has exactly degree $4$ (where we ignore the direction of bonds in the definition of degree); there is no bond between $V_i$ and $V_j'$ if $i\neq j$, and any (possible) bond from $V_i$ to $V_j$ corresponds to a unique bond from $V_j'$ to $V_i'$.

Now suppose $(n_0,n_1,\cdots,n_7,n_1',\cdots,n_7')\in \Sigma$. We assign a vector $m_q\in\Zb^3$ to each node $q\in\Tc\cup\Tc'$, such that if $q=\nu$ is a node of $\Tc$ then $m_q=\pm_\nu n_\nu$, and if $q=\nu'$ is a node of $\Tc'$ then $m_q=\pm_\nu n_\nu'$; note that any paired leaves must have the same $m_q$ value. Then, we assign a vector $m_e$ to each bond $e$ of the molecule $\Mb$ as follows. In case (1) in Definition \ref{analytic:def-mole} we define $m_e=m_q$ where $q$ is the node corresponding to the atom $V'$; in case (2) we define $m_e=m_q$ where $q$ is one of the leaves in the leaf pair; in case (3) we define $m_e=n_0$.

Then, for any atom $V\in\Mb$, we have
\begin{equation}\label{analytic:septic-moleeqn}\sum_{e}\pm_{V,e}m_e=0,\quad \sum_e\pm_{V,e}\langle m_e\rangle=m_{\mathrm{ex}}+O(1),
\end{equation} where $e$ runs over all bonds with one endpoint $V$, $\pm_{v,e}$ equals $+$ if $e$ is outgoing from $V$, and equal $-$ if $e$ is incoming at $v$, and $m_{\mathrm{ex}}$ is a constant vector which is one of $\pm m_j$ or $\pm m_j'$, see (\ref{analytic:sept-defset}). The mapping from $(n_0,n_1,\cdots,n_7,n_1',\cdots,n_7')\in \Sigma$ to $(m_e)$ satisfying (\ref{analytic:septic-moleeqn}) is a bijection.
\end{proposition}
\begin{proof} This can be checked using Definition \ref{analytic:def-mole}, the definitions of assignments $m_q$ and $m_e$, and the conditions in (\ref{analytic:sept-defset}) that define $\Sigma$. The first equation in (\ref{analytic:septic-moleeqn}) follows from the equalities $n_0=n_1+n_{234}+n_{567}$ and $n_{234}=n_2+n_3+n_4$ etc., and the second equation (\ref{analytic:septic-moleeqn}) follows from the condition $\Omega_j=m_j+O(1)$ and $\Omega_j'=m_j'+O(1)$ in (\ref{analytic:sept-defset}). The precise verifications are straightforward and are omitted here.
\end{proof}
\begin{remark}\label{analytic:rem1} We make a few observations. First, the molecule in Definition \ref{analytic:def-mole} is defined from the specific ternary tree $\Tc$ associated with the linear-cubic-cubic interaction. If one use any other ternary tree $\Tc$, then the rest of Definition \ref{analytic:def-mole} can be done in the same way to form a molecule $\Mb$. The assignments defined in Proposition \ref{analytic:prop-mol} are also the same and the conclusions of Proposition \ref{analytic:prop-mol} still hold. This will be used in Subsections \ref{3+3+3}--\ref{1+1+5} below, where $\Tc$ will be the ternary trees associated with the cubic-cubic-cubic and linear-linear-quintic interactions.

Second, the directions of the bonds in the molecule will not be too important in the counting estimates below, because they only affect the choice of signs $\pm_{V,e}$ in (\ref{analytic:septic-moleeqn}), but the counting estimates we will use, namely Lemma \ref{counting:lem2}--\ref{counting:lem3}, do not depend on the choices of these signs. As such, we will not specify the directions of bonds below unless necessary (when the direction affect the corresponding counting estimates).

Finally, since we only need to estimate each individual tree term (such as $\Eb\|\Mc_{1,3,3}\|_{X^{-1/2+\delta_2,b_+-1}}^2$) instead of the covariance of the different tree terms, here we may assume that the molecule $\Mb$ takes a special form where the only bonds between its upper and lower halves are vertical (i.e. no crossing bonds), see Figures \ref{analytic:fig:sept-mole1}--\ref{analytic:fig:sept-mole2}. Thus, the number of cases we need to consider is slightly smaller than that in \cite{DH21} for molecules of a fixed size.
\end{remark}
\underline{\emph{Part III: molecule counting estimates}}. By considering all the possibilities of $\Pc$ and exploiting symmetry, we can enumerate all the possibilities of the molecule $\Mb$ (without specifying bond directions, see Remark \ref{analytic:rem1}) corresponding to good $\Pc$ in Figures \ref{analytic:fig:sept-mole1}--\ref{analytic:fig:sept-mole2}. Note that $\{1,2\}\in\Pc$ if and only if a double bond is connected between $V_0$ and $V_1$, and $\pm_1=\pm_{234}$ if and only if these two bonds are in the same direction.
  \begin{figure}[h!]
  \includegraphics[scale=.4]{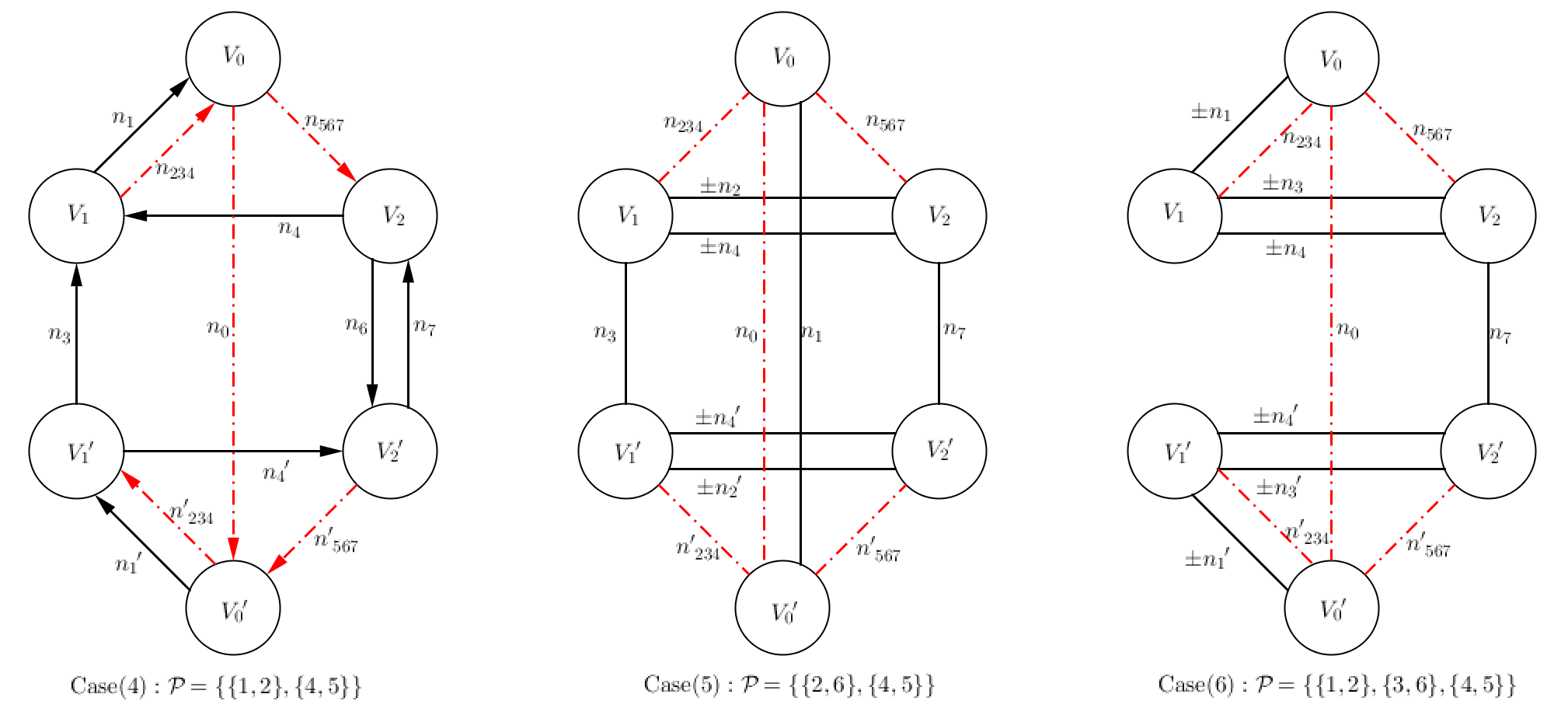}
  \caption{Cases (4)--(6) of the molecule $\Mb$ for the linear-cubic-cubic interaction. Here the directions of the bonds are specified only in Case (4), because this corresponds to the ternary tree in Figure \ref{analytic:fig:sept-tree}, in which the signs $\pm_j$ are fixed for the purpose of demonstration.}
  \label{analytic:fig:sept-mole2}
\end{figure}

By Proposition \ref{analytic:prop-mol}, the counting problem for $\#\Sigma$ can be reduced to the counting problem for the assignments $(m_e)$ where $e$ runs over all bonds in the molecule $\Mb$. For any $e$, since $m_e$ equals one of $\pm n_\nu$ or $\pm n_\nu'$, we can identify the corresponding dyadic number $N_e=N_\nu$ such that $N_e\sim \langle m_e\rangle=\langle n_\nu\rangle\sim \langle n_\nu'\rangle\sim N_\nu$. As such we can check that 
\begin{equation}\label{analytic:molecount1}\prod_{j=1}^7N_j^{-2}\cdot N_0^{-1}N_{234}^{-2}N_{567}^{-2}=\prod_{\textrm{good bonds }e}N_e^{-2}\cdot\prod_{\textrm{bad bonds }e}N_e^{-1}\end{equation} for the expression occurring on the right hand side of (\ref{analytic:sept-counting}), where the bad bonds are the ones indicated by red dotted lines in Figures \ref{analytic:fig:sept-mole1}--\ref{analytic:fig:sept-mole2}.

\smallskip
The main idea in estimating the number of assignments $(m_e)$ is to do it in steps. Each time we select one or two atoms in $\Mb$, and count the number of choices for $(m_e)$ where $e$ runs over all bonds connecting to one of the chosen atoms whose $m_e$ have not been fixed yet; after this step we fix the value of these $m_e$ and select the next atom(s) and so on. This can be viewed as a much simplified version of the counting algorithm in \cite{DH21}. In the case of one atom we can use Lemma \ref{counting:lem2}, and in the case of two atoms we can use Lemma \ref{counting:lem3} (or use Lemma \ref{counting:lem2} once or twice if $q=0$ or $r=0$ in Lemma \ref{counting:lem3}); note that in the two atoms case the role of $n_j$ and $n_j'$ in Lemma \ref{counting:lem3} is played by the $m_e$ where $e$ is connected to only one of the atoms, while the role of $\ell_j$ in Lemma \ref{counting:lem3} is played by the $m_e$ where $e$ is a bond connecting the two atoms. In each step we get an upper bound $\Xc$ for the number of choices, and multiply it by $\Yc$, which is the right hand side of (\ref{analytic:molecount1}) but only involving those bonds $e$ where $m_e$ is being counted. Let $\Xc\Yc=\Zc$, then in order to prove (\ref{analytic:eq-higher-order}) it suffices to show that
\begin{equation}\label{analytic:molestep}\prod_{\textrm{all steps}}\Zc\lesssim\Ac^2\cdot (K_\mathrm{max})^{-\eta^2}
\end{equation} in view of (\ref{analytic:sept-counting}). We will choose the atoms in some order that depends on the structure of the molecule; the basic guiding principle is to avoid loss in counting estimates in each individual step; for example we would prefer to choose two degree $4$ atoms when they are connected by a double bond, because the corresponding bound (\ref{counting:eqsum3}) saves a power. Precisely, we argue as follows.

(i) Consider Cases (1)--(4) in Figures \ref{analytic:fig:sept-mole1}--\ref{analytic:fig:sept-mole2}; note that there is at least a double bond between $V_2$ and $V_2'$ and at least a single bond between $V_1$ and $V_1'$. We may start by choosing the atoms $(V_2,V_2')$. By (\ref{counting:eqsum3}) and noting the two bad bonds connecting to $V_0$ and $V_0'$ respectively, we get $\Zc_1\lesssim 1$ in this first step. In fact, say we are in case (4), then $V_2$ and $V_2'$ are connected by two good bonds, say $e_1$ and $e_2$, and $V_2$ has one more bad bond $e_3$ and one more good bond $e_4$, while $V_2'$ has one more bad bond $e_5$ and one more good bond $e_6$, Then (\ref{counting:eqsum3}) implies that \[\Xc_1\lesssim (N_{e_1}\cdots N_{e_6})^{2}\max(N_{e_3},N_{e_4})^{-1}\max(N_{e_5},N_{e_6})^{-1},\quad \Yc_1=(N_{e_1}\cdots N_{e_6})^{-2}N_{e_3}N_{e_5},\] hence $\Zc_1\lesssim 1$.

Next we may choose the atoms $(V_0,V_0')$. Similarly, by (\ref{counting:eqsum5}) and noting the bad bond connecting $V_0$ and $V_0'$ and the two remaining bad bonds connecting to $V_1$ and $V_1'$, we get $\Zc_2\lesssim 1$. Finally we choose the atoms $(V_1,V_1')$ and use (\ref{counting:eqsum1}) to get \[\Zc_3\lesssim (N_{V_1,V_1'})^{-1};\quad N_{V_1,V_1'}:=\max_e N_e,\] where the maximum is taken over all bonds $e$ connecting $V_1$ and $V_1'$. Putting altogether, we see that (\ref{analytic:molestep}) is already proved if $N_{V_1,V_1'}\geq (K_{\mathrm{max}})^{10\eta^2}$.

Now suppose $N_{V_1,V_1'}\leq (K_\mathrm{max})^{10\eta^2}$, then we may first choose $(V_1,V_1')$, count the vectors $m_e$ trivially for bonds $e$ connecting $V_1$ and $V_1'$, and use (\ref{counting:eqsum1}) to get $\Zc_1\lesssim (K_\mathrm{max})^{40\eta^2}$ in this step. In the same way as above, we next choose $(V_0,V_0')$ and then $(V_2,V_2')$ to get $\Zc_2\lesssim 1$ and $\Zc_3\lesssim (N_{V_2,V_2'})^{-1}$ where $N_{V_2,V_2'}$ is defined similarly as $N_{V_1,V_1'}$ above. Thus, we see that (\ref{analytic:molestep}) is proved if $N_{V_2,V_2'}\geq (K_\mathrm{max})^{100\eta^2}$.

Next suppose $N_{V_j,V_j'}\leq (K_\mathrm{max})^{100\eta^2}$ for $1\leq j\leq 2$. Note that $\max(N_{234},N_{567})\gtrsim (K_\mathrm{max})^{\eta}$ by assumption; we shall assume $N_{234}\gtrsim (K_\mathrm{max})^{\eta}$, as the other case can be treated similarly (and is easier as there is no double bond between $V_0$ and $V_2$). Here we shall start with $(V_2,V_2')$ and get $\Zc_1\lesssim (K_\mathrm{max})^{O(\eta^2)}$ as above, then choose $(V_1,V_1')$ and count $m_e$ trivially for bonds $e$ connecting $V_1$ and $V_1'$. If $\{1,2\}\not\in\Pc$, then there is a single bond $e$ between $V_1$ and $V_0$, so the value of $m_e$ is uniquely fixed and we trivially get $\Zc_2\lesssim (K_\mathrm{max})^{O(\eta^2)}N_{234}^{-1}$; note that $N_{234}\gtrsim (K_\mathrm{max})^\eta$, and in the last step with $(V_0,V_0')$ we have $\Zc_3\lesssim 1$ due to (\ref{counting:eqsum1}), this also proves (\ref{analytic:molestep}).

Finally suppose $N_{V_j,V_j'}\leq (K_\mathrm{max})^{100\eta^2}$ for $1\leq j\leq 2$, and $\{1,2\}\in\Pc$, so there is a double bond between $V_0$ and $V_1$. Again we first choose $(V_2,V_2')$, then $(V_1,V_1')$ and finally $(V_0,V_0')$; as above we still have $\Zc_1\lesssim (K_\mathrm{max})^{O(\eta^2)}$ and $\Zc_3\lesssim 1$. In the second step, using (\ref{counting:eq4}) after counting $m_e$ trivially for bonds $e$ connecting $V_1$ and $V_1'$, and noticing that we have assumed $|n_1+n_{234}|\sim L$ (which plays the role of $M$ in (\ref{counting:eq4})), we get
\begin{equation}\label{analytic:eqn-z2new}\Zc_2\lesssim(K_\mathrm{max})^{O(\eta^2)}\cdot L^{-2};\quad \Zc_2\lesssim(K_\mathrm{max})^{O(\eta^2)}\cdot(\max(N_1,N_{234}))^{-2}\quad\textrm{if }\pm_1=\pm_{234}.
\end{equation} Moreover, if $\pm_1=\mp_{234}$ then either we have $\Ac=N_{234}L^{-1}$ due to (\ref{analytic:quantityA}), or we can use the $\Gamma$-condition, so using (\ref{counting:eq4+}) gives that $\Zc_2\lesssim (K_\mathrm{max})^{O(\eta^2)}\cdot L^2(\max(N_1,N_{234}))^{-2}$. In any case it is always true that \[\Zc_2\lesssim \Ac^2\cdot (K_\mathrm{max})^{O(\eta^2)}\cdot (\max(N_1,N_{234}))^{-1},\] which proves (\ref{analytic:molestep}) as $N_{234}\geq (K_\mathrm{max})^\eta$.

(ii) Consider Case (5) in Figure \ref{analytic:fig:sept-mole2}; note that there is a single bond $e_1$ between $V_1$ and $V_1'$, and a single bond $e_2$ between $V_2$ and $V_2'$. By symmetry we may assume $N_{e_1}\geq N_{e_2}$, so we shall start by choosing $(V_2,V_2')$. Using (\ref{counting:eqsum4}) we get $\Zc_1\lesssim N_{e_2}$, while as in (1) we next choose $(V_0,V_0')$ and get $\Zc_2\lesssim 1$ and finally choose $(V_1,V_1')$ to get $\Zc_3\lesssim N_{e_1}^{-2}$. Thus we are done if $N_{e_1}\gtrsim (K_\mathrm{max})^{100\eta^2}$. If $N_{e_2}\leq N_{e_1}\leq (K_\mathrm{max})^{100\eta^2}$, then we have $N_{V_j,V_j'}\leq (K_\mathrm{max})^{100\eta^2}$ for $1\leq j\leq 2$, so we can proceed as in the last part of (1) to conclude the proof.

(iii) Consider Case (6) in Figure \ref{analytic:fig:sept-mole2}. In this case we do something slightly different. First choose the atoms $(V_1',V_2')$; note that they are connected by a good double bond, and each has one good external bond and one bad external bond, so by (\ref{counting:eqsum3}) we get $\Zc_1\lesssim 1$. Then we choose $V_0'$, which now has only one bond remaining, and get $\Zc_2\lesssim N_0^{-1}$. Next we choose $V_2$ and use (\ref{counting:eqsum1}) to get $\Zc_3\lesssim 1$, and finally choose $V_1$ and repeat the last part of (1) (albeit for only one double bond instead of two) to get $\Zc_4\lesssim \Ac\cdot\max(N_1,N_{234})^{-1/2}$. This also proves (\ref{analytic:molestep}), because we always have $\max(N_0,N_1,N_{234})\gtrsim \max(N_{234},N_{567})\gtrsim (K_\mathrm{max})^\eta$ by assumption.

\smallskip
With the above discussions, we have now finished the proof of (\ref{analytic:molestep}), which then proves (\ref{analytic:eq-higher-order}) for $\Mc_{1,3,3}$.
\subsection{The cubic-cubic-cubic stochastic object}\label{3+3+3} In this subsection we prove Proposition \ref{analytic:prop-higher-order} for $\Mc_{3,3,3}$.

\underline{\emph{Part I: reduction to counting estimates}}. We shall reduce the estimate for $\Mc_{3,3,3}$ to a counting problem associated with the structure of this term; this part is largely identical to Part I of Subsection \ref{1+3+3}, so we will only list the main points.

Let the wave number of the output $\Mc_{3,3,3}$ be $n_0$, and the wave numbers of the three cubic inputs be $(n_1,n_2,n_3)$, $(n_4,n_5,n_6)$ and $(n_7,n_8,n_9)$. Let $\max_{j}\langle n_j\rangle\sim K_{\mathrm{max}}^+$, in the same way as in Subsection \ref{1+3+3} we may assume $K_{\mathrm{max}}^+\leq (K_\mathrm{max})^{100}$. Now, similar to (\ref{analytic:eq-septic2}), we have

\begin{equation}\label{analytic:eq-nonic2}
\Fc\Mc_{3,3,3}(t,n_0)=\sum_\Pc\sum_{\varphi_1,\cdots,\varphi_9\in\{\sin,\cos\}}\sum_{n_1,\cdots,n_9}(\Bc_{3,3,3})_{n_0\cdots n_9}(t)\cdot\Sc\Ic[n_j,\varphi_j:j\in O],
\end{equation} where the coefficients are
\begin{equation}\label{analytic:non-coef}
\begin{aligned}(\Bc_{3,3,3})_{n_0\cdots n_9}(t)&=\prod_{j=1}^9\frac{1}{\langle n_j\rangle}\cdot\frac{1}{\langle n_{123}\rangle}\frac{1}{\langle n_{456}\rangle}\frac{1}{\langle n_{789}\rangle}\int_0^t\int_0^t\int_0^t\chi(t)\chi(t_1)\chi(t_2)\chi(t_3)\\&\times\sin(\langle n_{123}\rangle (t-t_1))\sin(\langle n_{456}\rangle (t-t_2))\sin(\langle n_{789}\rangle (t-t_3))\\&\times\prod_{j=1}^3\varphi_j(\langle n_j\rangle t_1)\prod_{j=4}^6\varphi_j(\langle n_j\rangle t_2)\prod_{j=7}^9\varphi_j(\langle n_j\rangle t_3)\,\mathrm{d}t_1\mathrm{d}t_2\mathrm{d}t_3
\end{aligned}
\end{equation} instead of (\ref{analytic:sept-coef}). Here $\Pc$ is a collection of pairings of $\{1,\cdots,9\}$ and $O$ is the set of indices which is not in any subset in $\Pc$. Again $\Pc$ does not contain any pairing within $\{1,2,3\}$, $\{4,5,6\}$ or $\{7,8,9\}$ (so for example $\{1,2\}\not\in\Pc$), due to the renormalization in the definition of $\scubic$. Moreover, since the renormalization term \[\sum_{\textrm{cyclic}}\mathfrak{C}_{\leq N}^{(3,3)}[K_i,K_j](t)\cdot P_{K_k}\scubic[\leqN]\] (where the sum is taken over $(i,j,k)\in\{(1,2,3),(2,3,1),(3,1,2)\}$, see Definition \ref{ansatz:def-dyadic-modified-product}) exactly corresponds to the cases where $\Pc$ consists of exactly three pairings between two of the sets $\{1,2,3\}$, $\{4,5,6\}$ and $\{7,8,9\}$ (for example $\Pc=\{\{1,8\},\{2,7\},\{3,9\}\}$), we know that in (\ref{analytic:eq-nonic2}) we may also assume that $\Pc$ does not consist of exactly three such pairings. We will call such $\Pc$ \emph{good}.

Now we may fix a good $\Pc$ and $\varphi_j\,(j\in O)$, and reduce the time integral in (\ref{analytic:non-coef}) to a linear combination of integrals of form\begin{equation}\label{analytic:non-coef3}e^{i\langle n_0\rangle t}\int_0^t\int_0^t \int_0^t\chi(t)\chi(t_1)\chi(t_2)e^{i(\Omega_0t+\Omega_1t_1+\Omega_2t_2+\Omega_3t_3)}\,\mathrm{d}t_1\mathrm{d}t_2\mathrm{d}t_3:=e^{i\langle n_0\rangle t}\cdot\Hc(t,\Omega_0,\Omega_1,\Omega_2,\Omega_3),
\end{equation} where
\begin{equation}\label{analytic:non-coef4}
\begin{aligned}
\Omega_0&:=-\langle n_0\rangle\pm_{123}\langle n_{123}\rangle\pm_{456}\langle n_{456}\rangle \pm_{789}\langle n_{789}\rangle,
\\\Omega_1&:=\mp_{123}\langle n_{123}\rangle\pm_1\langle n_1\rangle\pm_2\langle n_2\rangle\pm_3\langle n_3\rangle,\\\Omega_2&:=\mp_{456}\langle n_{456}\rangle\pm_4\langle n_4\rangle\pm_5\langle n_5\rangle\pm_6\langle n_6\rangle,
\\\Omega_3&:=\mp_{789}\langle n_{789}\rangle\pm_7\langle n_7\rangle\pm_8\langle n_8\rangle\pm_9\langle n_9\rangle
\end{aligned}
\end{equation} with suitable signs $\pm_\nu$. Again, like in (\ref{analytic:sept-coefest2}) and (\ref{analytic:sept-coefest2+}), by Proposition \ref{counting:timeint} we have
\begin{equation}\label{analytic:non-coefest2}|(\Fc_t\Hc)(\xi,\Omega_0,\Omega_1,\Omega_2,\Omega_3)|\lesssim \widetilde{\Hc}(\xi,\lfloor \Omega_0\rfloor,\lfloor\Omega_1\rfloor,\lfloor\Omega_2\rfloor,\lfloor\Omega_3\rfloor)
\end{equation} for some function $\widetilde{\Hc}$, such that
\begin{equation}\label{analytic:non-coefest3}\sum_{m_0,m_1,m_2,m_3}|\widetilde{\Hc}(\xi,m_0,m_1,m_2,m_3)|\lesssim (K_\mathrm{max}^+)^{O(b_+-1/2)}\langle \xi\rangle^{4(1/2-b)},
\end{equation} where $(K_\mathrm{max}^+)^{O(b_+-1/2)}$ is again negligible. Let $\langle n_\nu\rangle\sim N_\nu$ for $\nu\in\{0,\cdots,9,123,456,789\}$ (note that $N_0\sim K_0$ and $N_{123}\sim K_1$ etc.), by repeating the arguments in Part I of Subsection \ref{1+3+3}, we can obtain, up to negligible factors, that
\begin{equation}\label{analytic:non-counting}\Eb\|\Mc_{3,3,3}\|_{X^{1/2+\delta_2,b_+-1}}^2\lesssim\prod_{j=1}^9N_j^{-2}\cdot N_0^{-1}N_{123}^{-2}N_{456}^{-2}N_{789}^{-2}\cdot\sup_{m_0,\cdots,m_3,m_0',\cdots,m_3'}(\#\Sigma),
\end{equation} where $\Sigma$ is the set defined by
\begin{multline}\label{analytic:non-defset}\Sigma=\big\{(n_0,n_1,\cdots,n_9,n_1',\cdots,n_9'):n_0=n_{1\cdots 9}=n_{1\cdots 9}',\,\,n_i+n_j=n_i'+n_j'=0\,\,(\forall \{i,j\}\in\Pc),\\n_j=n_j'\,\,(\forall j\in O),\,\,
\Omega_j=m_j+O(1)\,(0\leq j\leq 3),\,\,\Omega_j'=m_j'+O(1)\,\,(0\leq j\leq 3)\big\},
\end{multline} such that $(n_0,n_1',\cdots,n_9')$ satisfies the same equations, support requirements etc., as $(n_0,n_1,\cdots,n_9)$. The quantities $\Omega_j$ are defined from $(n_j)$ as in (\ref{analytic:non-coef3}), $\Omega_j'$ are similarly defined from $(n_j')$ but with opposite $\pm_\nu$ signs, and we assume $\langle n_\nu\rangle\sim\langle n_\nu'\rangle \sim N_\nu$.

\smallskip
\underline{\emph{Part II: reduction to molecules}}. Before moving to molecules, we first settle some simple cases. If $\{1,2,3\}\subseteq O$ (or if $\{4,5,6\}\subseteq O$ or $\{7,8,9\}\subseteq O$). By definition (\ref{analytic:non-defset}) of $\Sigma$, if $\{1,2,3\}\subseteq O$ then $n_j=n_j'$ for $1\leq j\leq 3$, hence $n_{123}=n_{123}'$. If we replace $n_{123}$ and $n_{123}'$ by two new variables $n_{\mathrm{ex}}$ and $n_{\mathrm{ex}}'$, then $(n_0,n_{\mathrm{ex}},n_4,\cdots,n_9,n_{\mathrm{ex}}',n_4',\cdots,n_9')$ belongs to a new set of form (\ref{analytic:sept-defset}). By the proof in Subsection \ref{1+3+3}, we see the number of choices for $(n_0,n_{\mathrm{ex}},n_4,\cdots,n_9,n_{\mathrm{ex}}',n_4',\cdots,n_9')$ is bounded by
\[N_0\cdot N_{123}^2N_{456}^2N_{789}^2(N_4\cdots N_9)^2\cdot
\left\{
\begin{aligned}&1,&&\textrm{always},\\
& (K_{\mathrm{max}})^{-\eta^2}, &&\textrm{if }\max(N_{456},N_{789})\geq (K_\mathrm{max})^\eta.
\end{aligned}
\right.\] Once these are fixed, we then count the number of choices of $(n_1=n_1',n_2=n_2',n_3=n_3')$, which is $\lesssim(N_1N_2N_3)^2N_{123}^{-1}$ using (\ref{counting:eqsum1}) and the fact that $\max(N_1,N_2,N_3)\gtrsim N_{123}$. Putting together, we see that these estimates imply (\ref{analytic:eq-higher-order}) no matter which one in $\{N_{123},N_{456},N_{789}\}$ is comparable to $K_\mathrm{max}$. From now on, we will assume that none of the sets $\{1,2,3\}$, $\{4,5,6\}$ or $\{7,8,9\}$ is contained in $O$.

Define the molecule $\Mb$ as in Definition \ref{analytic:def-mole}, but start with the ternary tree $\Tc$ associated with the cubic-cubic-cubic interaction; namely, $\Tc$ has root $0$ whose three children are $123$, $456$ and $789$, and each node $ijk$ has three children $i$, $j$ and $k$ which are all leaves. The rest of Definition \ref{analytic:def-mole} remains the same (for example the sign of each node $\nu$ is $\pm_\nu$, the atoms of $\Mb$ are branching nodes of $\Tc$ and $\Tc'$, etc.). Moreover, Proposition \ref{analytic:prop-mol} still holds, and its statement and proof remain exactly the same.

\begin{figure}[h!]
  \includegraphics[scale=.4]{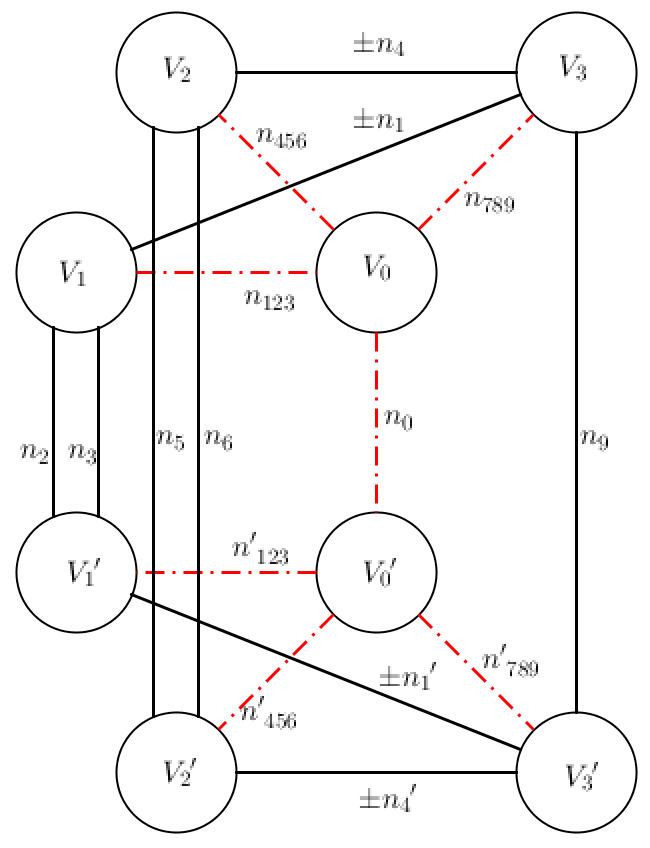}
  \caption{An example of molecule $\Mb$ for the cubic-cubic-cubic interaction corresponding to $\Pc=\{\{1,7\},\{4,8\}\}$. Here $V_0$ and $V_0'$ are the roots and $V_j,V_j'\,(1\leq j\leq 3)$ are their children. The other elements shown, including the bad bonds, are similar to Figure \ref{analytic:fig:sept-mole1}.}
  \label{analytic:fig:non-mole}
\end{figure}

We also point out two properties specific to the current molecule $\Mb$, which can be easily verified: it has no triple bond nor self connecting bond (because $\Pc$ is good, and none of the sets $\{1,2,3\}$, $\{4,5,6\}$ and $\{7,8,9\}$ is contained in $O$). In addition, if $V_0$ and $V_0'$ are the atoms corresponding to the roots of $\Tc$ and $\Tc'$, and $V_j,V_j'\,(1\leq j\leq 3)$ are their children, then there is a single bond $e_j$ between $V_0$ and each $V_j$, a single bond $e_j'$ between $V_0'$ and each $V_j'$, and a single bond $e_0$ between $V_0$ and $V_0'$. From the form of the product on the right hand side of (\ref{analytic:non-counting}), the bad bonds in $\Mb$ are exactly these $7$ single bonds (so good bonds are exactly those between $V_j,V_j'\,(1\leq j\leq 3)$. An example of such a molecule $\Mb$ (which has $8$ atoms) is depicted in Figure \ref{analytic:fig:non-mole}.

\underline{\emph{Part III: molecule counting estimates}}. By symmetry we shall assume $N_{123}\geq N_{456}\geq N_{789}$, and that they correspond to atoms $(V_1,V_2,V_3)$ and $(V_1',V_2',V_3')$ respectively. Same as in in Subsection \ref{1+3+3}, we only need to prove (which is much stronger than we need)
\begin{equation}\label{analytic:molestep2}(\textrm{number of choices for }(m_e))\cdot\Big(\prod_{\textrm{good bonds }e}N_e^{-2}\cdot\prod_{\textrm{bad bonds }e}N_e^{-1}\Big)\lesssim (K_\mathrm{max})^{-1};\end{equation} we do the counting in steps, with one or two atoms chosen at each step, and only need to estimate the value of $\Zc=\Xc\Yc$ in each step (as well as their product). In fact we will choose $V_0$ in the last step, before which we only work with $V_j,V_j'\,(1\leq j\leq 3)$. We also keep track of the number of (remaining) good bonds connected to each of these atoms, initially the values are $GB:=\{3,3,3,3,3,3\}$.

Consider $V_3$. Since there is no triple bonds, there are only two possibilities: when it has $3$ good single bonds, or when it has $1$ good single bond and $1$ good double bond.

(i) Suppose it has $3$ good single bonds, then we choose $V_3$ in the first step, by (\ref{counting:eqsum2}) we know $\Zc_1\lesssim N_{789}$, and after this step we have $GB=\{2,2,2,3,3\}$. Due to the lack of triple bonds, at least one atom with $2$ good bonds must be connected to an atom with $3$ good bonds. We next choose this atom with $2$ good bonds (which also has a bad bond connecting to $V_0$ or $V_0'$), by (\ref{counting:eqsum1}) we have $\Zc_2\lesssim 1$, and after this step we have $GB=\{1,2,2,3\}$. Now the $2$ atoms with $2$ good bonds cannot be connected by a double bond (since otherwise removing then would leave $GB=\{1,3\}$ which is impossible), so we next choose the one of these $2$ atoms that is not connected to the atom with $1$ good bond. By (\ref{counting:eqsum1}) we have $\Zc_3\lesssim 1$, and after this step we have $GB=\{1,1,2\}$. Next we choose the only atom with $2$ good bonds and get $\Zc_4\lesssim 1$ by (\ref{counting:eqsum1}).

Now there are two atoms (other than $V_0$ and $V_0'$) left, and each of them has only one bad bond. Next we choose these atoms and get $\Zc_{5}\lesssim N_{\nu_1}^{-1}$ and $\Zc_6\lesssim N_{\nu_2}^{-1}$, where $\nu_j\in\{123,456,789\}$ but only one of them is $789$ (as $V_3$ is chosen in the first step). Finally we choose $V_0$ to get $\Zc_7\lesssim N_0^{-1}$. Summing up, we get
\[\Zc_1\cdots \Zc_7\lesssim N_{789}\cdot N_{789}^{-1}\cdot N_{456}^{-1}\cdot N_0^{-1}\] which is bounded by $(N_\mathrm{max})^{-1}$, because the two maximal elements in $\{N_0,N_{123},N_{456},N_{789}\}$ must be comparable due to the equality $n_0=n_{123}+n_{456}+n_{789}$. As an example of this procedure, in the case of Figure \ref{analytic:fig:non-mole}, we have the sequence $V_3\to V_2\to V_1\to V_3'\to V_2'\to V_1'\to V_0$.

(ii) Suppose it has $1$ good single bond and $1$ good double bond, say the good double bond connects to another atom $\widetilde{V}$. In the first step we choose $(V_3,\widetilde{V})$; since they are connected by a good double bond, and each has one good external bond and one bad external bond, we can apply (\ref{counting:eqsum3}) to get $\Zc_1\lesssim 1$. Then we choose the next atoms arbitrarily, one at a time, such that each atom has at most $2$ good bonds when it is chosen. This is always possible because there is no triple bond, and by (\ref{counting:eqsum1}) we have $\Zc_j\lesssim 1$ at each step.

We continue until there are only two atoms (other than $V_0$ and $V_0'$) left, say $V_+$ and $V_-$. At least one of them is not $V_3'$; say $V_+\neq V_3'$. Then we choose $V_-$ and get $\Zc_5\lesssim 1$ using (\ref{counting:eqsum1}), and choose $V_+$, which now has only one bad bond, to get $\Zc_6\lesssim N_{\nu}^{-1}$ where $\nu\in\{123,456\}$ since $V_+\not\in\{V_3,V_3'\}$. Finally we choose $V_0$ to get $\Zc_7\lesssim N_0^{-1}$. Summing up, we get \[\Zc_1\cdots \Zc_7\lesssim N_{456}^{-1}\cdot N_0^{-1}\lesssim (K_\mathrm{max})^{-1},\] in the same way as (i).

\smallskip
With the above discussions, we have now finished the proof of (\ref{analytic:molestep2}), which then proves (\ref{analytic:eq-higher-order}) for $\Mc_{3,3,3}$. The proof of Proposition \ref{analytic:prop-higher-order} is now complete.
\subsection{The linear-linear-quintic stochastic object}\label{1+1+5} In this subsection we prove Proposition \ref{analytic:prop-higher-order2}. Compared to $\Mc_{1,3,3}$ and $\Mc_{3,3,3}$, the term $\Mc_{1,1,5}$ has more complicated form, so we first need some pre-processing. To be precise, in the argument below we will decompose $\Mc_{1,1,5}$ into two ``special" components $(\Mc_{1,1,5})_{sp1}$ and $(\Mc_{1,1,5})_{sp2}$, and two ``generic" components $(\Mc_{1,1,5})_{ge1}$ and $(\Mc_{1,1,5})_{ge2}$, and estimate them separately in Subsections \ref{generic1proof}--\ref{generic2proof}.

Recall from (\ref{analytic:eq-quintic-decomp-pure}) that $3\squintic[\leqN]$ can be decomposed into three terms, where the last term is the resistor $18\sresistor[\leqN]$. We shall denote the sum of the first two terms by $3(\squintic[\leqN])_{\mathrm{nr}}$, and will single out a first ``special" component of $\Mc_{1,1,5}$, namely
\begin{equation}\label{analytic:115sp1}(\Mc_{1,1,5})_{sp1}:=P_{\leq N} \bigg[ 54 \squadratic[\leqN] \sresistor[\leqN]- 18 \mathfrak{C}^{(1,5)}_{\leq N} P_{\leq N} \slinear[blue][\leqN]  
   - 54 ( 2 \HLL + \HHL \Big) \Big( \, \slinear[blue][\leqN], \slinear[blue][\leqN], \sresistor[\leqN] ) \bigg],\end{equation} so that
\begin{multline}\label{analytic:115rm1}\Mc_{1,1,5}-(\Mc_{1,1,5})_{sp1}=P_{\leq N} \bigg[ 9 \squadratic[\leqN] (\squintic[\leqN])_{\mathrm{nr}}- \Gamma_{\leq N} \cdot\scubic[\leqN]   
   \\
- 9 \Big( 2 \HLL + \HHL \Big) \Big( \, \slinear[blue][\leqN], \slinear[blue][\leqN], (\squintic[\leqN])_{\mathrm{nr}}\Big) \bigg].
\end{multline} Here $\HLL$ and $\HHL$ are defined as in Definition \ref{ansatz:def-para-product} but with $\squintic[\leqN]$ replaced by either $6\sresistor[\leqN]$ or $(\squintic[\leqN])_{\mathrm{nr}}$, except that for $\HLL$, the renormalization term $\mathfrak{C}_{\leq N}^{(1,5)}$ is included only in (\ref{analytic:115sp1}) but not in (\ref{analytic:115rm1}).

Now we dyadically decompose (\ref{analytic:115rm1}). By the definitions of $\HLL$ and $\HHL$ (see Definition \ref{ansatz:def-para-product}) we have
\begin{equation}\label{analytic:115rm1+} (\ref{analytic:115rm1})=\sum_{\substack{K_0,\cdots,K_3\leq N\\ K_3\geq (K_\mathrm{max})^\eta}}\sum_{L_1,L_2\lesssim N}\mathfrak{Quad}[K_*]\Big((\squintic[\leqN])_{\mathrm{nr}}\Big)-\Gamma_{\leq N} \cdot\scubic[\leqN],
\end{equation} where $K_*=(K_0,K_1,K_2,K_3,K_{13}=L_1,K_{23}=L_2)$ and $\mathfrak{Quad}$ is defined in Definition \ref{linear:def-quad-operator}. Denote also $K_{\mathrm{max}}:=\max(K_0,\cdots,K_3)$ similar to Proposition \ref{analytic:prop-higher-order}. 

We shall further decompose (\ref{analytic:115rm1}) and identify a first ``generic" component, namely
\begin{equation}\label{analytic:115ge1}(\Mc_{1,1,5})_{ge1}:=\sum_{K_*}^{(1)}\mathfrak{Quad}[K_*](\squintic[\leqN])_{\mathrm{nr}}.
\end{equation} Here $\sum^{(1)}$ includes those terms where $K_3\geq (K_\mathrm{max})^\eta$ and $\min(K_0,L_1,L_2)\geq (K_{\mathrm{max}})^{10\eta^2}$. Then
\begin{equation}\label{analytic:115rm2}
(\ref{analytic:115rm1})-(\Mc_{1,1,5})_{ge1}=\sum_{K_*}^{(2)}\mathfrak{Quad}[K_*](\squintic[\leqN])_{\mathrm{nr}}-\Gamma_{\leq N}\cdot\scubic[\leqN],
\end{equation} where $\sum^{(2)}$ includes those terms in (\ref{analytic:115rm1+}) that are not in $\sum^{(1)}$. 

Next, for fixed $K_*$, define the operator $\Gamma^{\mathrm{op}}[K_*]$ in the same way as in Definition \ref{para:def-gamma-operator}, but with the extra factor $1_{L_1}(n_{13})1_{L_2}(n_{23})$ inserted to (\ref{para:def-gamma-operator-eqn}). Then we can identify a second ``special" component from (\ref{analytic:115rm2}), namely
\begin{equation}\label{analytic:115sp2}
(\Mc_{1,1,5})_{sp2}:=\sum_{K_*}^{(2)}\Gamma^{\mathrm{op}}[K_*] \scubic[\leqN]-\Gamma_{\leq N}\scubic[\leqN].
\end{equation}Finally, the last component, which is the second ``generic" component, is given by
\begin{equation}\label{analytic:115ge2}(\Mc_{1,1,5})_{ge2}=\sum_{K_*}^{(2)}\bigg[\mathfrak{Quad}[K_*](\squintic[\leqN])_{\mathrm{nr}}-\Gamma^{\mathrm{op}}[K_*] \scubic[\leqN]\bigg],
\end{equation} so we have $\Mc_{1,1,5}=(\Mc_{1,1,5})_{sp1}+(\Mc_{1,1,5})_{sp2}+(\Mc_{1,1,5})_{ge1}+(\Mc_{1,1,5})_{ge2}$. It then suffices to prove (\ref{analytic:eq-higher-order2}) for each of these components. For convenience of presentation, these estimates will be proved  in the order $ge1\to sp2\to sp1\to ge2$, which is different from the natural order where they are introduced.
\subsubsection{The first generic component}\label{generic1proof} We start by proving (\ref{analytic:eq-higher-order2}) for the generic term $(\Mc_{1,1,5})_{ge1}$. This actually follows directly from Lemmas \ref{analytic:lem-quintic-pure-nopairing}--\ref{analytic:lem-quintic-pure-one-pairing} and Lemma \ref{linear:lem-quad-dyadic}. More precisely, if $K_3\geq (K_{\mathrm{max}})^\eta$ and $\min(K_0,L_1,L_2)\geq (K_{\mathrm{max}})^{10\eta^2}$, then Lemmas \ref{analytic:lem-quintic-pure-nopairing}--\ref{analytic:lem-quintic-pure-one-pairing} implies 
\[\Big\|P_{K_3}(\squintic[\leqN])_{\mathrm{nr}}\Big\|_{L_\omega^{2p}X^{1/2+\delta_2,b}([-T,T])}\lesssim (2p)^{5/2}K_3^{2\delta_2}T^\alpha,\] while (a slight modification of the proof of) Lemma \ref{linear:lem-quad-dyadic} implies that 
\[\big\|\mathfrak{Quad}[K_*]\big\|_{L_\omega^{2p}(X^{1/2+\delta_2,b}([-T,T])\to X^{-1/2+\delta_2,b_+-1}([-T,T]))}\lesssim 2p\cdot (K_{\mathrm{max}})^{O(\delta_2)}(K_0^{-\eta/2}+K_3^{-1/2})T^\alpha.\] Hence, by our choice of parameters, we have that
\[\Big\|\mathfrak{Quad}[K_*](\squintic[\leqN])_{\mathrm{nr}}\Big\|_{L_\omega^pX^{-1/2+\delta_2,b_+-1}([-T,T])}\lesssim p^{7/2}T^\alpha(K_{\mathrm{max}})^{-\eta^3},\] which implies (\ref{analytic:eq-higher-order2}) for $(\Mc_{1,1,5})_{ge1}$.

\subsubsection{The second special component} Next we prove (\ref{analytic:eq-higher-order2}) for the special term $(\Mc_{1,1,5})_{sp2}$, which is basically the same as Lemma \ref{para:lem-11X2-renormalized}. In fact, let $w_4:=\scubic[\leqN]$, then we can decompose (\ref{analytic:115sp2}) as
\begin{equation}\label{analytic:115sp2dec}(\ref{analytic:115sp2})=\bigg\{\sum_{K_*}\Gamma^{\mathrm{op}}[K_*] w_4-\Gamma_{\leq N}w_4\bigg\}-\sum_{K_*}^{(1)}\Gamma^{\mathrm{op}}[K_*] w_4:=\mathrm{I}-\sum_{K_*}\mathrm{II}[K_*],\end{equation} where in the first term no restriction is put on the summation in $K_*$. Now this term $\mathrm{I}$ is exactly equal to (\ref{para:eq-11X2-renormalized-p1}) since we are summing over all different choices of $L_1$ and $L_2$, so by the proof of Lemma \ref{para:lem-11X2-renormalized} we have
\[\mathbb{E}\|\mathrm{I}\|_{X^{-1/2+\delta_2,b_+-1}([-T,T])}^p\lesssim T^{p\alpha} \mathbb{E}\|w_4\|_{X^{-\epsilon,b}([-T,T]))}^p\lesssim T^{p\alpha} p^{3p/2}\] using also (\ref{analytic:eq-cubic}). As for $\mathrm{II}[K_*]$, we may fix a choice of $K_*$, and notice that $\sum^{(1)}$ requires that $K_0\geq K_{\mathrm{max}}^\eta$, so we are in \emph{Case (1)} of (\ref{para:eq-11X2-renormalized-p2}) as in the proof of Lemma \ref{para:lem-11X2-renormalized}. By repeating that proof (note also that the sine cancellation is not needed assuming $K_0\geq K_{\mathrm{max}}^\eta$), we see that
\[\mathbb{E}\|\mathrm{II}[K_*]\|_{X^{-1/2+\delta_2,b_+-1}([-T,T])}^p\lesssim (K_{\mathrm{max}})^{-\eta   p/4}T^{p\alpha} \mathbb{E}\|w_4\|_{X^{-\epsilon,b}([-T,T]))}^p\lesssim (K_{\mathrm{max}})^{-\eta p/4}T^{p\alpha} p^{3p/2}\] which is then summable in $K_*$ (noting that $\max(L_1,L_2)\lesssim K_{\mathrm{max}}$), and this proves (\ref{analytic:eq-higher-order2}) for $(\Mc_{1,1,5})_{sp2}$.
\subsubsection{The first special component}
Now we turn to $(\Mc_{1,1,5})_{sp1}$. Let $K_*=(K_0,\cdots, K_3,K_{13}=L_1,K_{23}=L_2)$, we make a dyadic decomposition by attaching the factor $\prod_{j=0}^3 1_{K_j}(n_j)\cdot 1_{L_1}(n_{13})1_{L_2}(n_{23})$ where as before $K_{\mathrm{max}}=\max(K_0\cdots,K_3)$, and note that, by the definitions of $\HLL$ and $\HHL$ in Definition \ref{ansatz:def-para-product}, we can always assume $K_3\geq (K_{\mathrm{max}})^\eta$ here. Note that now (this component of) $(\Mc_{1,1,5})_{sp1}$ can be written as $\mathfrak{Quad}[K_*](\sresistor[\leqN])$, so if $\min(K_0,L_1,L_2)\geq (K_{\mathrm{max}})^{10\eta^2}$, we can apply Lemma \ref{linear:lem-quad-dyadic} together with Lemma \ref{analytic:lem-resistor}, and argue as in Subsection \ref{generic1proof} to conclude that
\[\begin{aligned}&\quad\,\,\|\mathfrak{Quad}[K_*](\sresistor[\leqN])\|_{L_\omega^pX^{-1/2+\delta_2,b_+-1}([-T,T])}\\
&\lesssim\|P_{K_3}\sresistor[\leqN]\|_{L_\omega^{2p}X^{-1/2+\delta_2,b_+-1}([-T,T])}\cdot \big\|\mathfrak{Quad}[K_*]\big\|_{L_\omega^{2p}(X^{1/2+\delta_2,b}([-T,T])\to X^{-1/2+\delta_2,b_+-1}([-T,T]))}\\
&\lesssim p^{7/2}(K_{\mathrm{max}})^{O(\delta_2)}(K_0^{-\eta/2}+K_3^{-1/2})T^\alpha\lesssim p^{7/2}T^\alpha(K_{\mathrm{max}})^{-\eta^3},
\end{aligned}\] which proves (\ref{analytic:eq-higher-order2}) for $(\Mc_{1,1,5})_{sp1}$.

Now let us assume $\min(K_0,L_1,L_2)\leq (K_{\mathrm{max}})^{10\eta^2}$. By Corollary \ref{diagram:corollary-resistor}, we have
\begin{equation}\label{analytic:resistor-rep}\Fc_x\sresistor[\leqN](t,n)=F_{\leq N}^{\mathrm{sin}}(t,n)\Sc\Ic[n,\mathrm{sin}]+F_{\leq N}^{\mathrm{cos}}(t,n)\Sc\Ic[n,\mathrm{cos}]
\end{equation} for some explicitly defined functions $F_{\leq N}^\varphi$. As such, and noticing also that the terms with one pairing are exactly cancelled by the $\mathfrak{C}_{\leq N}^{(1,5)}$ factor, we see that $(\Mc_{1,1,5})_{sp1}=\sum_{K_*}\mathrm{III}[K_*]$ where the summation is taken over $K_*$ such that $K_3\geq (K_{\mathrm{max}})^\eta$ and $\min(K_0,L_1,L_2)\leq (K_{\mathrm{max}})^{10\eta^2}$, and
\begin{equation}\label{analytic:115sp1dec2}
\Fc_x\mathrm{III}[K_*](t,n_0)=\sum_{n_0=n_{123}}\sum_{(\varphi_j)}\prod_{j=0}^3 1_{K_j}(n_j)\cdot\frac{\varphi_1(t\langle n_1\rangle)\varphi_2(t\langle n_2\rangle)}{\langle n_1\rangle\langle n_2\rangle}F_{\leq N}^{\varphi_3}(t,n_3)\Sc\Ic[n_j,\varphi_j:1\leq j\leq 3].
\end{equation} By (the proof of) the $X^{1/2-\epsilon,b}$ bound in Lemma \ref{analytic:lem-resistor}, we can deduce the bound satisfied by the function $F_{\leq N}^\varphi$, namely
\begin{equation}\label{analytic:coef-bound}
\|P_{K_3}F_{\leq N}^\varphi\|_{L_t^\infty\ell_n^2}\lesssim\|P_{K_3}F_{\leq N}^\varphi\|_{X^{0,b}}\lesssim K_3^{-1/2+\epsilon}\lesssim K_3^{-1/4}.
\end{equation} Therefore, by (\ref{analytic:coef-bound}), hypercontractivity estimates and the reduction arguments in Subsection \ref{prep:remark-reduction} we have
\[
\begin{aligned}\|\mathrm{III}[K_*]\|_{L_\omega^pX^{-1/2+\delta_2,b_+-1}([-T,T])}&\lesssim T^{\alpha}p^{3/2}K_0^{-1/2+2\delta_2}\|\mathrm{III}[K_*]\|_{L_t^\infty\ell_n^2L_\omega^2}\\
&\lesssim T^\alpha p^{3/2} K_0^{-1/2+2\delta_2}(K_1K_2)^{-1}K_3^{-1/4}\cdot\bigg(\sup_{n_3}\sum_{(n_0,n_1,n_2)}1\bigg)^{1/2},
\end{aligned}\] where in the last sum over $(n_0,n_1,n_2)$ we require that $n_0=n_{123}$ and $\langle n_j\rangle \sim K_j\,(0\leq j\leq 3)$ as well as $\langle n_{13}\rangle\sim L_1,\,\langle n_{23}\rangle\sim L_2$. The value of this sum is clearly bounded by
\[(\min(K_0,L_1,L_2))^{3}(\mathrm{med}(K_0,K_1,K_2))^3\lesssim \big(K_0^{1/2-2\delta_2}K_1K_2\cdot (K_{\mathrm{max}})^{O(\eta^2)}\big)^2.\] This gives an acceptable contribution since $K_3\geq (K_{\mathrm{max}})^\eta$, so by summing over $K_*$, we have proved (\ref{analytic:eq-higher-order2}) for $(\Mc_{1,1,5})_{sp1}$.
\subsubsection{The second generic component}\label{generic2proof} Finally we prove (\ref{analytic:eq-higher-order2}) for the generic term $(\Mc_{1,1,5})_{ge2}$. The proof is similar to Sections \ref{1+3+3} and \ref{3+3+3}, and mostly relies on the molecular technology.

\underline{\emph{Part I: reduction to counting estimates}}. This part is largely identical to Part I in Subsection \ref{1+3+3} or \ref{3+3+3}, so we will only list the main points.

Let the wave number of the output $(\Mc_{1,1,5})_{ge2}$ be $n_0$, and the wave numbers of the two linear inputs be $n_1$ and $n_2$. Let the wave numbers of the inputs in the quintic term $(\squintic[\leqN])_{\mathrm{nr}}$ be $(n_3,\cdots,n_7)$, where $(n_5,n_6,n_7)$ are the wave numbers of the inputs in the cubic term $\scubic[\leqN]$. Let $\max_{j}\langle n_j\rangle\sim K_{\mathrm{max}}^+$, as before may assume $K_{\mathrm{max}}^+\leq (K_\mathrm{max})^{100}$. Note also that $|n_j|\sim K_j$ for $0\leq j\leq 2$, $|n_{34567}|\sim K_3$, and $|n_{j34567}|\sim L_j$ for $j\in\{1,2\}$. Now we have

\begin{equation}\label{analytic:eq-115term1}
\Fc(\Mc_{1,1,5})_{ge2}(t,n_0)=\sum_\Pc\sum_{\varphi_1,\cdots,\varphi_7\in\{\sin,\cos\}}\sum_{n_1,\cdots,n_7}(\Bc_{1,1,5})_{n_0\cdots n_7}(t)\cdot\Sc\Ic[n_j,\varphi_j:j\in O],
\end{equation} where the coefficients are
\begin{equation}\label{analytic:115-coef}
\begin{aligned}(\Bc_{1,1,5})_{n_0\cdots n_7}(t)&=\prod_{j=1}^7\frac{1}{\langle n_j\rangle}\cdot\frac{1}{\langle n_{567}\rangle}\frac{1}{\langle n_{34567}\rangle}\int_0^t\int_0^{t_1}\chi(t)\chi(t_1)\chi(t_2)\sin(\langle n_{34567}\rangle (t-t_1))\\&\times\sin(\langle n_{567}\rangle (t_1-t_2))\prod_{j=1}^2\varphi_j(\langle n_j\rangle t)\prod_{j=3}^4\varphi_j(\langle n_j\rangle t_1)\prod_{j=5}^7\varphi_j(\langle n_j\rangle t_2)\,\mathrm{d}t_2\mathrm{d}t_1.
\end{aligned}
\end{equation} Here $\Pc$ is a collection of pairings of $\{1,\cdots,7\}$ and $O$ the set of indices not in any subset in $\Pc$. Now $\Pc$ does not contain $\{1,2\}$, $\{3,4\}$ or any pairing within $\{5,6,7\}$, and contains at most one pair within the set $\{3,4,5,6,7\}$, due to our choice of $(\squintic[\leqN])_{\mathrm{nr}}$; moreover, the renormalization term
\[\Gamma^{\mathrm{op}}[K_*]\scubic[\leqN]\] exactly corresponds to the cases where $\Pc=\{\{1,3\},\{2,4\}\}$ or $\Pc=\{\{1,4\},\{2,3\}\}$, we know that in (\ref{analytic:eq-nonic2}) we may also assume that $\Pc$ is not $\{\{1,3\},\{2,4\}\}$ nor $\{\{1,4\},\{2,3\}\}$. We will call such $\Pc$ \emph{good}.

Now we may fix a good $\Pc$ and $\varphi_j\,(j\in O)$, and reduce the time integral in (\ref{analytic:115-coef}) to a linear combination of integrals of form\begin{equation}\label{analytic:non-coef3}e^{i\langle n_0\rangle t}\int_0^t\int_0^{t_1} \int_0^t\chi(t)\chi(t_1)\chi(t_2)e^{i(\Omega_0t+\Omega_1t_1+\Omega_2t_2)}\,\mathrm{d}t_2\mathrm{d}t_1:=e^{i\langle n_0\rangle t}\cdot\Hc(t,\Omega_0,\Omega_1,\Omega_2),
\end{equation} where
\begin{equation}\label{analytic:115-coef4}
\begin{aligned}
\Omega_0&:=-\langle n_0\rangle\pm_{1}\langle n_{1}\rangle\pm_{2}\langle n_{2}\rangle \pm_{34567}\langle n_{34567}\rangle,
\\\Omega_1&:=\mp_{34567}\langle n_{34567}\rangle\pm_3\langle n_3\rangle\pm_4\langle n_4\rangle\pm_{567}\langle n_{567}\rangle,\\\Omega_2&:=\mp_{567}\langle n_{567}\rangle\pm_5\langle n_5\rangle\pm_6\langle n_6\rangle\pm_7\langle n_7\rangle
\end{aligned}
\end{equation} with suitable signs $\pm_\nu$. Like before, by Proposition \ref{counting:timeint} we have
\begin{equation}\label{analytic:115-coefest2}|(\Fc_t\Hc)(\xi,\Omega_0,\Omega_1,\Omega_2)|\lesssim \widetilde{\Hc}(\xi,\lfloor \Omega_0\rfloor,\lfloor\Omega_1\rfloor,\lfloor\Omega_2\rfloor)
\end{equation} for some function $\widetilde{\Hc}$, such that
\begin{equation}\label{analytic:115-coefest3}\sum_{m_0,m_1,m_2}|\widetilde{\Hc}(\xi,m_0,m_1,m_2)|\lesssim (K_\mathrm{max}^+)^{O(b_+-1/2)}\langle \xi\rangle^{4(1/2-b)},
\end{equation} where $(K_\mathrm{max}^+)^{O(b_+-1/2)}$ is negligible. Let $\langle n_\nu\rangle\sim N_\nu$ for $\nu\in\{0,\cdots,7,34567,567\}$ (note that $N_0\sim K_0$ and $N_{34567}\sim K_3$ etc.), by repeating the arguments as before, we can obtain, up to negligible factors, that
\begin{equation}\label{analytic:115-counting}\Eb\|(\Mc_{1,1,5})_{ge2}\|_{X^{1/2+\delta_2,b_+-1}}^2\lesssim\prod_{j=1}^9N_j^{-2}\cdot N_0^{-1}N_{34567}^{-2}N_{567}^{-2}\Ac^{-2}\cdot\sup_{m_0,\cdots,m_3,m_0',\cdots,m_3'}(\#\Sigma),
\end{equation} where $\Sigma$ is the set defined by
\begin{multline}\label{analytic:115-defset}\Sigma=\big\{(n_0,n_1,\cdots,n_7,n_1',\cdots,n_7'):n_0=n_{1\cdots 7}=n_{1\cdots 7}',\,\,n_i+n_j=n_i'+n_j'=0\,(\forall \{i,j\}\in\Pc),\\n_j=n_j'\,(\forall j\in O),\,\,
\Omega_j=m_j+O(1)\,(0\leq j\leq 2),\,\,\Omega_j'=m_j'+O(1)\,(0\leq j\leq 2)\big\},
\end{multline} with the same notation as before; the dyadic support conditions we have are that $\langle n_\nu\rangle\sim N_\nu$ and $\langle n_{j34567}\rangle\sim L_j\,(j\in\{1,2\})$, and the same for $n_j'$. 

The quantity $\Ac$ in (\ref{analytic:115-defset}) is related to the sine-cancellation kernel, and is similar as in Subsection \ref{1+3+3} but slightly different, so we explain it here. First, we will assume $\{1,3\}\in\Pc$ (or $\{2,4\}\in\Pc$ etc. by symmetry, but we just switch the indices $1$ and $2$ if necessary) and $\pm_1=\mp_{34567}$; if this is not satisfied then $\Ac=1$. Now, in the definition (\ref{analytic:115ge2}) of $(\Mc_{1,1,5})_{ge2}$ we may fix $L_1$ and sum in $L_2$. By definition of $\sum^{(2)}$ (see \eqref{analytic:115rm2}), this summation in $L_2$ results in either a restriction $|n_{234567}|_\infty< (K_\mathrm{max})^{10\eta^2}$, or no restriction at all on $n_{234567}$ (in the latter case we must have $\min(K_0,L_1)<(K_\mathrm{max})^{10\eta^2}$). We shall symmetrize in $n_1$ and $n_{34567}$ as in Subsection \ref{1+3+3} and get the corresponding value of $\Ac$, except that we may have an extra term involving
\begin{equation}\label{analytic:symmetrization}\mathbf{1}\big\{|n_1+n_2|_\infty<(K_\mathrm{max})^{10\eta^2}\big\}-\mathbf{1}\big\{|n_{34567}+n_2|_\infty<(K_\mathrm{max})^{10\eta^2}\big\}\end{equation} due to the possible restriction $|n_{234567}|_\infty< (K_\mathrm{max})^{10\eta^2}$.

In summary, we have the following description of $\Ac$ and the set $\Sigma$ in (\ref{analytic:115-counting}) and (\ref{analytic:115-defset}). In the first scenario we have $\Ac=N_{34567}L_1^{-1}$ similar to (\ref{analytic:quantityA}). In the second scenario $\Ac=1$, but we have either $\{1,3\}\not\in\Pc$, or $\pm_1=\pm_{34567}$, or we can further require the $\Gamma$-condition in $\Sigma$, or we can require that \emph{exactly} one of $|n_{12}|_\infty$ and $|n_{234567}|_\infty$ is $<(K_\mathrm{max})^{10\eta^2}$ (in this scenario $L_1$ is still fixed but we are summing over $L_2$ so this parameter is absent in the definition of terms; also whatever conditions we have for $(n_j)$, the similar conditions will also hold for $(n_j')$).

\underline{\emph{Part II: reduction to molecules}}. Define the molecule $\Mb$ as in Definition \ref{analytic:def-mole}, but start with the ternary tree $\Tc$ associated with the linear-linear-quintic interaction; namely, $\Tc$ has root $0$ whose three children are $1$, $2$ and $34567$, and $34567$ has three children $3$, $4$ and $567$, and $567$ has three children $5$, $6$ and $7$, and nodes $1$--$7$ are all leaves. The rest of Definition \ref{analytic:def-mole} and Proposition \ref{analytic:prop-mol} remain the same.

\begin{figure}[h!]
  \includegraphics[scale=.4]{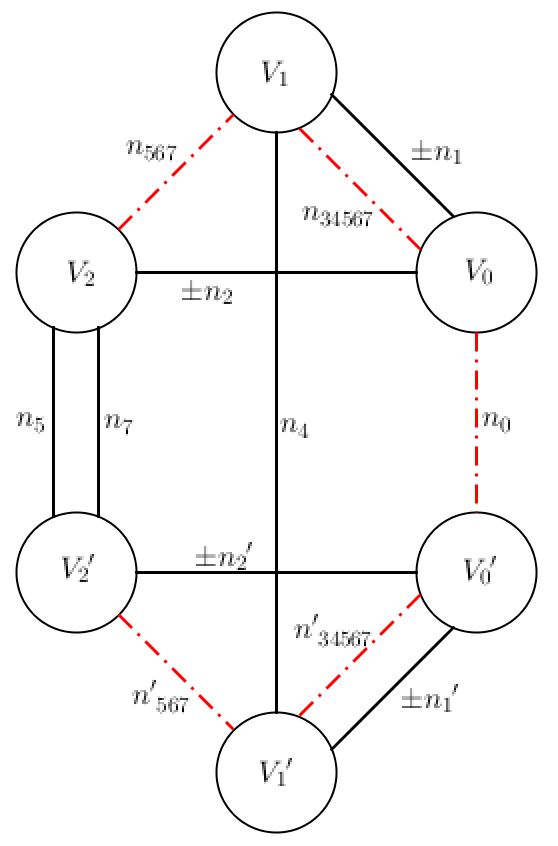}
  \caption{An example of molecule $\Mb$ for the linear-linear-quintic interaction corresponding to $\Pc=\{\{1,3\},\{2,6\}\}$. Here $V_0,V_0'$ are the roots, $V_1,V_1'$ are their children and $V_2,V_2'$ are the children of $V_1,V_1'$. The other elements shown, including the bad bonds, are similar to Figure \ref{analytic:fig:sept-mole1}.}
  \label{analytic:fig:115-mole}
\end{figure}

The following are easily verified for the molecule $\Mb$, thanks to $\Pc$ being good: it contains no self connecting bond and no triple bond between any $V_i$ and $V_j$ (or $V_i'$ and $V_j'$), here we assume $V_0$ is the atom corresponding to the root of $\Tc$, and $V_1$ corresponds to its child, and $V_2$ corresponds to its grandchild. Moreover, $\{1,3\}\in\Pc$ if and only if a double bond is connected between $V_0$ and $V_1$, and $\pm_1=\pm_{34567}$ if and only if these bonds have the same direction. The bad bonds in $\Mb$ are those bond $e$ whose $m_e$ equals $\pm n_{34567}$, $\pm n_{567}$ or $n_0$; there are $5$ of them, two connecting $V_1$ to $V_0$ and $V_2$, two connecting $V_1'$ to $V_0'$ and $V_2'$, and one connecting $V_0$ and $V_0'$. An example of such a molecule $\Mb$ is depicted in Figure \ref{analytic:fig:115-mole}.

\underline{\emph{Part III: molecule counting estimates}}. Recall from the definition of $\sum^{(2)}$ as in \eqref{analytic:115rm2} that we have $|n_{34567}|\sim N_{34567}\sim K_3\gtrsim (K_\mathrm{max})^{\eta}$, and that one of the values $|n_0|,|n_{134567}|,|n_{234567}|$ or (by symmetrization in (\ref{analytic:symmetrization})) $|n_{12}|$ is $\lesssim (K_\mathrm{max})^{10\eta^2}$. In any case, as in Subsections \ref{1+3+3} and \ref{3+3+3}, we only need to prove that
\begin{equation}\label{analytic:molestep3}(\textrm{number of choices for }(m_e))\cdot\Big(\prod_{\textrm{good bonds }e}N_e^{-2}\cdot\prod_{\textrm{bad bonds }e}N_e^{-1}\Big)\lesssim \Ac^2\cdot(K_\mathrm{max})^{-\eta^2},\end{equation} and again the left hand side can be written as a product of $\Zc$ over all steps.

Denote by $e_{\mathrm{bad}}$ the bad bond connecting $V_0$ and $V_1$; we have $N_{e_{\mathrm{bad}}}\sim N_{34567}\geq (K_\mathrm{max})^{\eta}$. We start by choosing $(V_0,V_0')$ and get a factor of $\Zc_1$; then we choose $(V_1,V_1')$ to get $\Zc_2$, and finally choose $(V_2,V_2')$ to get $\Zc_3$. By (\ref{counting:eqsum1}) we easily see that $\Zc_3\lesssim 1$. Also by (\ref{counting:eqsum5}), and since no bad bond is connected between $V_1$ and $V_1'$, we see that $\Zc_2\lesssim N_{V_1,V_1'}^{-1}$, where $N_{V_1,V_1'}$ is the maximum of $N_e$ for $e$ running over all bonds connecting $V_1$ and $V_1'$. Alternatively, if we count $m_e$ for those $e$ trivially in the second step, then by (\ref{counting:eq4}), we also get $\Zc_2\lesssim N_{V_1,V_1'}^9\cdot Q^{-1}$, where the quantity $Q$ is defined to be $Q=N_e$ if $V_1$ is connected to $V_0$ by a single bond $e$, and $Q=\max(N_{e_1},N_{e_2})$ if $V_1$ is connected to $V_0$ by a double bond $(e_1,e_2)$ of same directions, and $Q\sim |n_{e_1}-n_{e_2}|$ if $V_1$ is connected to $V_0$ by a double bond $(e_1,e_2)$ of opposite directions. As for $\Zc_1$, note that we may assume that one of $|n_0|,|n_{134567}|,|n_{234567}|,|n_{12}|$ is $\leq (K_\mathrm{max})^{10\eta^2}$. If $|n_0|\leq (K_\mathrm{max})^{10\eta^2}$ we count $n_0$ trivially and then apply (\ref{counting:eqsum1}) or (\ref{counting:eqsum5}) to get $\Zc_1\lesssim (K_\mathrm{max})^{O(\eta^2)}$; if $\min(|n_{134567}|,|n_{234567}|,|n_{12}|)\leq (K_\mathrm{max})^{10\eta^2}$, we may first choose $V_0'$ using (\ref{counting:eq6+}), and then choose $V_0$ using (\ref{counting:eqsum1}) to get the same result (strictly speaking, in the case $|n_{12}|\leq (K_\mathrm{max})^{10\eta^2}$ we get a potentially losing factor $N_1^{-1}N_0$ when choosing $V_0'$, but after $V_0'$ is chosen, the value of $n_{e_{\mathrm{bad}}}$ will be fixed up to $(K_\mathrm{max})^{10\eta^2}$, which leads to a factor $N_0^{-1}$ when choosing $V_0$, so in total we still have $\Zc_1\lesssim (K_\mathrm{max})^{O(\eta^2)}$).

In summary, we see that (\ref{analytic:molestep3}) is already proved, unless $N_{V_1,V_1'}\lesssim (K_\mathrm{max})^{O(\eta^2)}$, and $V_1$ is connected to $V_0$ by a double bond $(e_1,e_2)$ of opposite directions, and $|n_{e_1}-n_{e_2}|\lesssim (K_\mathrm{max})^{O(\eta^2)}$. Therefore, we may assume all these conditions hold; in particular we have $\{1,3\}\in\Pc$, $\pm_1=\mp_{34567}$ and $L_1\sim |n_{134567}|\lesssim (K_\mathrm{max})^{O(\eta^2)}$. Now if $\min(|n_{234567}|,|n_{12}|)\leq (K_\mathrm{max})^{O(\eta^2)}$ also, then each of $|n_0|,|n_1|,|n_2|,|n_{34567}|$ is within $(K_\mathrm{max})^{O(\eta^2)}$ of each other, which easily implies that $\Zc_1\lesssim(K_\mathrm{max})^{O(\eta^2)}\cdot N_{34567}^{-1}$ and hence (\ref{analytic:molestep3}); we may thus ignore the factor (\ref{analytic:symmetrization}), hence we have either $\Ac=N_{34567}L_1^{-1}$ which also implies (\ref{analytic:molestep3}) easily, or the $\Gamma$-condition $|n_1|_\infty\geq \Gamma\geq |n_{34567}|_\infty$ (or $|n_{34567}|_\infty\geq \Gamma\geq |n_1|_\infty$).

Now assume we are in the last case with the $\Gamma$-condition. In the first step, we may first choose $V_0'$ using (\ref{counting:eq6+}), then choose $V_0$ but using (\ref{counting:eq4+}) instead of (\ref{counting:eqsum1}) to get $\Zc_1\lesssim (K_\mathrm{max})^{O(\eta^2)}\cdot N_{34567}^{-1}$. Note that here we apply (\ref{counting:eq4+}) directly if $V_0$ has no bond connecting to $V_2$; otherwise $V_0$ has a single bond $e$ connecting to $V_2$ such that $n_e=\pm n_2$ (see Figure \ref{analytic:fig:115-mole}), then the value of $n_2$ must be fixed up to error $(K_\mathrm{max})^{O(\eta^2)}$, so we can apply (\ref{counting:eq4+}) when counting $(n_1,n_{34567})$. Therefore, in this case we also get (\ref{analytic:molestep3}).

\smallskip With the above discussions, we have now finished the proof of (\ref{analytic:molestep3}), which then proves (\ref{analytic:eq-higher-order2}) for $(\Mc_{1,1,5})_{ge2}$. The proof of Proposition \ref{analytic:prop-higher-order2} is now complete.

\section{Proof of main estimates}\label{section:proof-main-estimates}

In this section, we prove the main estimates from Subsection \ref{section:main-estimates} and the modified version from Lemma \ref{global:lem-one-zero-revisited}. Since all of the main ideas have already been presented in Sections \ref{section:diagrams}-\ref{section:analytic2}, each proof will primarily consist of references to earlier lemmas and propositions.

\begin{proof}[Proof of Proposition \ref{ansatz:prop-Xj}:] 
The desired estimate follows directly from the decomposition in Lemma \ref{para:lem-decomposition} and the estimates in Lemma \ref{para:lem-strichartz-XXone} and Lemma \ref{para:lem-strichartz-XXtwo}.
\end{proof}

We now turn to the nonlinear terms in the evolution equation for $Y_{\leq N}$. We first prove Proposition \ref{ansatz:prop-Y-two}, which controls terms containing two linear stochastic objects. 

\begin{proof}[Proof of Proposition \ref{ansatz:prop-Y-two}:] 
The three different estimates \ref{item:prop-Y-two-1}, \ref{item:prop-Y-two-2}, and \ref{item:prop-Y-two-3} in the proposition have already been proven exactly as stated in earlier sections. The two explicit stochastic objects in \ref{item:prop-Y-two-1} have been estimated in Proposition \ref{analytic:prop-mixed} and Proposition \ref{analytic:prop-higher-order2}. The two para-controlled random operators in \ref{item:prop-Y-two-2} have been estimated in Subsection \ref{section:para-two-linear}. The term involving $\XXone$ has been estimated in Lemma \ref{para:lem-11X1} and the term involving $\XXtwo$ has been estimated in Lemma \ref{para:lem-11X2-decomposition}, Lemma \ref{para:lem-11X2}, and Lemma \ref{para:lem-11X2-renormalized}.  Finally, the random operator in $Y$ from \ref{item:prop-Y-two-3} has been estimated in Proposition \ref{linear:prop-quad}. 
\end{proof}

We now prove Proposition \ref{ansatz:prop-Y-one}, which controls terms containing one linear stochastic object. Since Proposition \ref{ansatz:prop-Y-one} controls more terms than Proposition \ref{ansatz:prop-Y-two}, its proof will be longer. Nevertheless, it still primarily consists of a collection of earlier estimates.

\begin{proof}[Proof of Proposition \ref{ansatz:prop-Y-one}:] To simplify the notation, we focus only on $T=1$ and $\Jc=[-1,1]$, since the general case follows from minor modifications. Furthermore, we restrict all norms to $[-1,1]\times \T^3$. 
In the following, each estimate holds on $A$-certain events $E_A\in \mathcal{E}$, which will not be repeated below. For notational purposes, it is convenient to replace the $A$-factor in our desired estimate by $A^3$. By eventually replacing $E_A$ with $E_{A^{1/3}}$, this can be re-adjusted after our argument. In the statement of the proposition, we are concerned with the terms 
\begin{equation}\label{pmain:eq-one-recall}
\Big( \Pi^\ast_{\leq N}- \HLL \Big) \Big( \slinear[blue][\leqN], \zeta^{(2)}_{\leq N}, \zeta^{(3)}_{\leq N} \Big),
\end{equation}
where $\zeta^{(2)},\zeta^{(3)} \in \Symb_0^p$. In our estimate of \eqref{pmain:eq-one-recall}, we distinguish the cases when the arguments $\zeta^{(2)}$ and $\zeta^{(3)}$ contain two, one, or zero cubic stochastic objects. \\

\emph{Case 1: $\zeta^{(2)}=\zeta^{(3)}=\scubic$.} The resulting term is an explicit septic stochastic object, which has been estimated in Proposition \ref{analytic:prop-higher-order}. \\

\emph{Case 2: $\zeta^{(2)}=\scubic$.} In this case, we further distinguish the sub-cases
\begin{equation*}
\zeta^{(3)}= \squintic, \qquad \zeta^{(3)} \in \Big\{ \slinear[green][\leqM], \SXXone, \SXXtwo \Big\}, \qquad \text{and} \qquad \zeta^{(3)}= Y. 
\end{equation*}

\emph{Case 2.a: $\zeta^{(2)}=\scubic$ and $\zeta^{(3)}=\squintic$.} In this case, we 
prove that 
\begin{equation*}
\Big\| \Big( \Pi^\ast_{\leq N} - \HLL \Big)\Big( \slinear[blue][\leqN], \scubic[\leqN], \squintic[\leqN] \Big) \Big\|_{X^{-1/2+\delta_2,b_+-1}} \leq A. 
\end{equation*}
Using the $\HLH$-operator from Definition \ref{linear:def-hlh}, we first decompose\footnote{The $\HLL$-term contains no portion of $\mathfrak{C}^{(1,5)}_{\leq N}$ since $\mathfrak{C}^{(1,5)}_{\leq N}[N_1,N_3]$ is only non-zero for $N_1=N_3$.}
\begin{align}
&\Big( \Pi^\ast_{\leq N} - \HLL \Big)\Big( \slinear[blue][\leqN], \scubic[\leqN], \squintic[\leqN] \Big) \notag \\
=& \slinear[blue][\leqN] \scubic[\leqN] \squintic[\leqN] - \Big( \HLL + \HLH \Big) \Big( \slinear[blue][\leqN], \scubic[\leqN], \squintic[\leqN] \Big) \label{pmain:eq-135-p1}\\
+& \HLH \Big( \slinear[blue][\leqN], \scubic[\leqN], \squintic[\leqN] \Big) - \mathfrak{C}^{(1,5)}_{\leq N} \scubic[\leqN]. \label{pmain:eq-135-p2}
\end{align}
Using Proposition \ref{linear:prop-lincub} and Proposition \ref{analytic:prop-quintic},  the first summand \eqref{pmain:eq-135-p1} is controlled by 
\begin{align*}
    &\Big\| \slinear[blue][\leqN] \scubic[\leqN] \squintic[\leqN] - \Big( \HLL + \HLH \Big) \Big( \slinear[blue][\leqN], \scubic[\leqN], \squintic[\leqN] \Big) \Big\|_{X^{-1/2+\delta_2,b_+-1}} \\
    \leq& A \Big\| \squintic[\leqN] \Big\|_{X^{1/2-\delta_1,b}} 
    \leq A^2 .
\end{align*}
Thus, it remains to treat the second summand \eqref{pmain:eq-135-p2}. We emphasize that the estimate heaviliy relies on the $\mathfrak{C}^{(1,5)}_{\leq N}$-term, which originates from the \oftt-cancellation. From Definition \ref{linear:def-hlh}, we have that 
\begin{align}
&\HLH \Big( \slinear[blue][\leqN], \scubic[\leqN], \squintic[\leqN] \Big) - \mathfrak{C}^{(1,5)}_{\leq N} \scubic[\leqN] \notag \\
=& \sum_{\substack{N_1,N_2,N_3 \leq N \colon \\ N_2 \leq N_1^\nu, \\ N_1 \sim N_3 > N_1^\eta}} \bigg( \slinear[blue][N_1] P_{N_3} \squintic[\leqN] - \mathfrak{C}^{(1,5)}_{\leq N}[N_1,N_3] \bigg) P_{N_2} \scubic[\leqN] \label{pmain:eq-135-p3} \\
-& \bigg( \sum_{\substack{N_1,N_2,N_3 \leq N}} - \sum_{\substack{N_1,N_2,N_3 \leq N \colon \\ N_2 \leq N_1^\nu \\N_1 \sim N_3 > N_1^\eta}} \bigg)
\mathfrak{C}^{(1,5)}_{\leq N}[N_1,N_3] P_{N_2} \scubic[\leqN]. \label{pmain:eq-135-p4}
\end{align}
We first estimate the remaining main term, i.e., \eqref{pmain:eq-135-p3}. Using $X^{-1/2+\delta_2,0}\hookrightarrow X^{-1/2+\delta_2,b_+-1}$ and product estimates (see e.g. \cite[Appendix A]{Tao06}), we obtain that 
\begin{align}
    &\bigg\| \bigg( \slinear[blue][N_1] P_{N_3} \squintic[\leqN] - \mathfrak{C}^{(1,5)}_{\leq N}[N_1,N_3] \bigg) P_{N_2} \scubic[\leqN]  \bigg\|_{X^{-1/2+\delta_2,b_+-1}} \notag \\
\lesssim& \, \bigg\| \bigg( \slinear[blue][N_1] P_{N_3} \squintic[\leqN] - \mathfrak{C}^{(1,5)}_{\leq N}[N_1,N_3] \bigg) P_{N_2} \scubic[\leqN]  \bigg\|_{L_t^2 H_x^{-1/2+\delta_2}} \notag  \\
\lesssim& \,  \bigg\| \bigg( \slinear[blue][N_1] P_{N_3} \squintic[\leqN] - \mathfrak{C}^{(1,5)}_{\leq N}[N_1,N_3] \bigg)\bigg\|_{L_t^2 H_x^{-1/2+\delta_2}} 
\Big\| P_{N_2} \scubic[\leqN] \Big\|_{L_t^\infty \C_x^1} \label{pmain:eq-135-p5}
\end{align}
Using Proposition \ref{analytic:prop-sextic2}, Lemma \ref{analytic:lem-cubic}, and $N_1 \sim N_3$, it holds that 
\begin{equation*}
\eqref{pmain:eq-135-p5} \lesssim A^2 N_2^2  N_1^{-100\nu}. 
\end{equation*}
Due to the frequency-restrictions $N_2 \leq N_1^\nu$ in \eqref{pmain:eq-135-p3}, this yields an acceptable contribution. We now estimate the minor term \eqref{pmain:eq-135-p4}. Using $X^{-1/2+\delta_2,0}\hookrightarrow X^{-1/2+\delta_2,b_+-1}$, we have that 
\begin{align}
\Big\| \mathfrak{C}^{(1,5)}_{\leq N}[N_1,N_3] P_{N_2} \scubic[\leqN] \Big\|_{X^{-1/2+\delta_2,b_+-1}} 
&\lesssim \Big\| \mathfrak{C}^{(1,5)}_{\leq N}[N_1,N_3] P_{N_2} \scubic[\leqN] \Big\|_{L_t^2 H_x^{-1/2+\delta_2}} \notag \\
&\lesssim N_2^{-1/2+\delta_2+\epsilon} \big\| \mathfrak{C}^{(1,5)}_{\leq N}[N_1,N_3] \big\|_{L_t^\infty L_x^\infty} \Big\| P_{N_2} \scubic[\leqN] \Big\|_{L_t^\infty \C_x^{-\epsilon}}. 
\end{align}
Using Lemma \ref{analytic:lem-C15-C33}, Lemma \ref{analytic:lem-cubic}, and the frequency-restriction in \eqref{pmain:eq-135-p4}, plus the fact that $\mathfrak{C}^{(1,5)}_{\leq N}[N_1,N_3]=0$ unless $N_1\sim N_3$, the resulting contribution is acceptable.\\

\emph{Case 2.b: $\zeta^{(2)}=\scubic$ and $\zeta^{(3)} \in \big\{ \slinear[green][\leqM], \SXXone, \SXXtwo\big\}$:} In this case, we prove that
\begin{equation*}
\Big\|     \Big( \Pi^\ast_{\leq N} - \HLL \Big)\Big( \slinear[blue][\leqN], \scubic[\leqN], \zeta^{(3)} \Big) \Big\|_{X^{-1/2+\delta_2,b-1}} \leq A^3 \Big( 1 + \| v_{\leq N} \|_{X^{-1,b}}^2 + \| Y_{\leq N} \|_{X^{1/2+\delta_2},b} \Big). 
\end{equation*}
This case  slightly easier than Case 2.a, since it does not utilize the \oftt-cancellation. Using the $\HLH$-operator from Definition \ref{linear:def-hlh}, we first decompose 
\begin{align}
    &\Big( \Pi^\ast_{\leq N} - \HLL \Big)\Big( \slinear[blue][\leqN], \scubic[\leqN], \zeta^{(3)} \Big) \notag \\
=& \slinear[blue][\leqN] \scubic[\leqN] \zeta^{(3)}_{\leq N} - \Big( \HLL + \HLH \Big) \Big( \slinear[blue][\leqN], \scubic[\leqN], \zeta^{(3)}_{\leq N}  \Big) \label{pmain:eq-13z-p1}\\
+& \HLH \Big( \slinear[blue][\leqN], \scubic[\leqN], \zeta^{(3)}_{\leq N}  \Big). \label{pmain:eq-13z-p2}
\end{align}
The first summand \eqref{pmain:eq-13z-p1} can be estimated using Proposition \ref{linear:prop-lincub} and the regularity estimate for $\slinear[green][\leqM]$, $\SXXone$, or $\SXXtwo$ (Lemma \ref{analytic:lem-linear} or Proposition \ref{ansatz:prop-Xj}). For the second summand \eqref{pmain:eq-13z-p2}, we recall that 

\begin{equation*}
     \HLH \Big( \slinear[blue][\leqN], \scubic[\leqN], \zeta^{(3)}_{\leq N}  \Big) 
     = \sum_{\substack{N_1,N_{2},N_3 \colon \\ N_{2} \leq N_1^\nu \\ N_1 \sim N_3 >N_1^\eta}} \slinear[blue][N_1] P_{N_{2}} \scubic[\leqN] P_{N_3} \zeta^{(3)}_{\leq N}.
\end{equation*}
For the dyadic components, we have that 
\begin{align*}
\Big\|  \slinear[blue][N_1] P_{N_{2}} \scubic[\leqN] P_{N_3} \zeta^{(3)}_{\leq N} \Big\|_{X^{-1/2+\delta_2,b_+-1}} 
\lesssim& \, \Big\|  \slinear[blue][N_1] P_{N_{2}} \scubic[\leqN] P_{N_3} \zeta^{(3)}_{\leq N} \Big\|_{L_t^2 H_x^{-1/2+\delta_2}} \\
\lesssim& \Big\|  \slinear[blue][N_1]P_{N_3} \zeta^{(3)}_{\leq N} \Big\|_{L_t^2 H_x^{-1/2+\delta_2}} \Big\|  P_{N_{2}} \scubic[\leqN]  \Big\|_{L_t^\infty \C_x^{1}}. 
\end{align*}
Depending on $\zeta^{(3)} \in \big\{ \slinear[green][\leqM], \SXXone, \SXXtwo \big\}$, we use either Lemma \ref{analytic:lem-resonant-blue-green} or Proposition \ref{para:prop-resonant}, which yield 
\begin{equation*}
\Big\| \slinear[blue][N_1] P_{N_3} \zeta^{(3)}_{\leq N}\Big\|_{L_t^2 H_x^{-1/2+\delta_2}} \leq A^2 N_1^{-1/2+12\eta} \Big(  1+ \| v_{\leq N}\|_{X^{-1,b}}^2 + \| Y_{\leq N}\|_{X^{1/2+\delta_2,b}} \Big). 
\end{equation*}
Using Lemma \ref{analytic:lem-cubic}, we also obtain that 
\begin{equation*}
 \Big\|  P_{N_{2}} \scubic[\leqN]  \Big\|_{L_t^\infty \C_x^{1}} \leq A N_{2}^{2}. 
\end{equation*}
Under the frequency-restriction $N_{2}\leq N_1^\nu$, this yields the desired conclusion. \\

\emph{Case 2.c: $\zeta^{(2)}=\scubic$ and $\zeta^{(3)}=Y$.} From Proposition \ref{linear:prop-lincub}, we obtain that 
\begin{equation*}
\Big\|  \Big( \Pi^{\ast}_{\leq N} - \HLL \Big) \Big( \slinear[blue][N_1], \scubic[\leqN], Y_{\leq N} \Big) \Big\|_{X^{-1/2+\delta_2,b_+-1}} \leq A \| Y_{\leq N} \|_{X^{1/2+\delta_2,b}},
\end{equation*}
which is acceptable.\\

\emph{Case 3: $\zeta^{(2)},\zeta^{(3)} \in \Symb_{1/2}^p$.} We now treat the case when neither $\zeta^{(2)}$ nor $\zeta^{(3)}$ is given by the cubic stochastic object. We first decompose
\begin{equation}\label{pmain:eq-1zz-p1}
\Big( \Pi^\ast_{\leq N} - \HLL \Big) \Big( \slinear[blue][\leqN], \zeta^{(2)}_{\leq N}, \zeta^{(3)}_{\leq N}\Big) 
=\sum_{\substack{N_1,N_2,N_3 \colon \\ \max(N_2,N_3)>N_1^\eta}} \hspace{-2ex}
\Pi^\ast_{\leq N}  \Big( \slinear[blue][\leqN], P_{N_2} \zeta^{(2)}_{\leq N}, P_{N_3}\zeta^{(3)}_{\leq N}\Big).
\end{equation}
In the case $N_1^\eta < \max(N_2,N_3) \ll N_1$, it follows from the bilinear estimate (Lemma \ref{bilinear:lem-bilinear}) that 
\begin{align*}
  &\Big\| \Pi^\ast_{\leq N}  \Big( \slinear[blue][\leqN], P_{N_2} \zeta^{(2)}_{\leq N}, P_{N_3}\zeta^{(3)}_{\leq N}\Big) \Big\|_{X^{-1/2+\delta_2,b_+-1}} \\
  =& \Big\| \, \slinear[blue][\leqN]  P_{N_2} \zeta^{(2)}_{\leq N} P_{N_3}\zeta^{(3)}_{\leq N} \Big\|_{X^{-1/2+\delta_2,b_+-1}} \\
  \leq& A N_1^{\delta_2} \Big(N_1^{-1/4}+ \max(N_2,N_3)^{-1/3}\Big)
  N_2^{\delta_1} N_3^{\delta_1} \prod_{j=2}^3 \big\| \zeta^{(j)}_{\leq N} \big\|_{X^{1/2-\delta_1,b}} \\
  \leq& A N_1^{\delta_2 - \eta (1/3-2\delta_1)} \prod_{j=2}^3 \big\| \zeta^{(j)}_{\leq N} \big\|_{X^{1/2-\delta_1,b}}  \\
  \leq&  A N_1^{-\epsilon}\prod_{j=2}^3 \big\| \zeta^{(j)}_{\leq N} \big\|_{X^{1/2-\delta_1,b}}. 
\end{align*}
Using the regularity estimates for $\slinear[green][\leqM]$, $\SXXone$, $\SXXtwo$, or $Y$, we obtain that 
\begin{equation*}
    \prod_{j=2}^3 \big\| \zeta^{(j)}_{\leq N} \big\|_{X^{1/2-\delta_1,b}} \leq A^2  \Big(  1+ \| v_{\leq N}\|_{X^{-1,b}}^4 + \| Y_{\leq N}\|_{X^{1/2+\delta_2,b}}^2 \Big).
\end{equation*}
Since this yields an acceptable contribution, it remains to treat the case $\max(N_2,N_3)\gtrsim N_1$. By symmetry, we may assume that $N_2 = \max(N_2,N_3)$. In the following, we distinguish several different sub-cases. The case distinctions are driven by the following three aspects: 
\begin{itemize}
    \item[$\bullet$] In high$\times$high$\times$low-interactions between $\slinear[blue][\leqN]$, $\squintic[\leqN]$, and a third argument, we need to utilize the \oftt-cancellation.
    \item[$\bullet$] In high$\times$high$\times$low-interactions between $\slinear[blue][\leqN]$, an element of  $\big\{ \slinear[green][\leqM], \SXXone, \SXXtwo \big\}$, and a third argument, we need to use our product estimates (Lemma \ref{analytic:lem-resonant-blue-green} or Proposition \ref{para:prop-resonant}). 
    \item[$\bullet$] In comparison with $Y_{\leq N}$, the elements of $\big\{ \slinear[green][\leqM], \squintic, \SXXone, \SXXtwo \big\}$ have fewer derivatives but obey better Strichartz estimates. 
\end{itemize}~\\

\emph{Case 3.a: $\zeta^{(2)}=\squintic$ and $\zeta^{(3)} \in \Symb_{1/2}^p$.} For $\zeta^{(3)}\in \Symb_{1/2}^{p}$, it holds that 
\begin{align}
\Pi^\ast_{\leq N}\Big( \slinear[blue][N_1], P_{N_2} \squintic[\leqN], P_{N_3}\zeta^{(3)} \Big) 
&= \Big( \slinear[blue][N_1] P_{N_2} \squintic[\leqN] - \mathfrak{C}^{(1,5)}_{\leq N}[N_1,N_2]\Big)  P_{N_3}\zeta^{(3)} \label{pmain:eq-15z-p1} \\
&-  \mathbf{1}\Big\{ \zeta^{(3)} = \squintic \Big\} \, \mathfrak{C}^{(1,5)}_{\leq N}[N_1,N_3] P_{N_2}\squintic[\leqN]. \label{pmain:eq-15z-p2}
\end{align} 
While the \oftt-cancellation is essential in our estimate of \eqref{pmain:eq-15z-p1}, it is irrelevant in  the additional term in \eqref{pmain:eq-15z-p2}. Using $X^{-1/2+\delta_2,0}\hookrightarrow X^{-1/2+\delta_2,b_+-1}$ and product estimates, the first summand \eqref{pmain:eq-15z-p1} can be estimated by 
\begin{align}
&\Big\| \Big( \slinear[blue][N_1] P_{N_2} \squintic[\leqN] - \mathfrak{C}^{(1,5)}_{\leq N}[N_1,N_2]\Big)  P_{N_3}\zeta^{(3)} \Big\|_{X^{-1/2+\delta_2,b_+-1}} \notag \\
\leq& \, \Big\| \Big( \slinear[blue][N_1] P_{N_2} \squintic[\leqN] - \mathfrak{C}^{(1,5)}_{\leq N}[N_1,N_2]\Big)  P_{N_3}\zeta^{(3)} \Big\|_{L_t^2 H_x^{-1/2+\delta_2}} \notag \\
\leq&  \Big\|  \slinear[blue][N_1] P_{N_2} \squintic[\leqN] - \mathfrak{C}^{(1,5)}_{\leq N}[N_1,N_2] \Big\|_{L_t^\infty \C_x^{-1/2+\delta_1+\epsilon}} \big\| \zeta^{(3)} \big\|_{L_t^\infty H_x^{1/2-\delta_1}}. \label{pmain:eq-15z-p4} 
\end{align}
The first factor in \eqref{pmain:eq-15z-p4} is estimated using Proposition \ref{analytic:prop-sextic2}. The $\zeta^{(3)}$-factor in \eqref{pmain:eq-15z-p4} is estimated using the regularity estimate for either $\slinear[green][\leqM]$, $\SXXone$, $\SXXtwo$, or $Y$. As a result, we obtain that 
\begin{align*}
\eqref{pmain:eq-15z-p4} \lesssim A^3 N_2^{-1/200}
\Big(  1+ \| v_{\leq N}\|_{X^{-1,b}}^2 + \| Y_{\leq N}\|_{X^{1/2+\delta_2,b}} \Big),
\end{align*}
which is acceptable. 
We now turn to the additional term in the case $\zeta^{(3)}=\squintic$, i.e., the term in  \eqref{pmain:eq-15z-p2}. Using $X^{-1/2+\delta_2,0}\hookrightarrow X^{-1/2+\delta_2,b_+-1}$, it follows that 
\begin{align}
\Big\|  \mathfrak{C}^{(1,5)}_{\leq N}[N_1,N_3] P_{N_2}\squintic[\leqN] \Big\|_{X^{-1/2+\delta_2,b_+-1}} 
&\lesssim \Big\|  \mathfrak{C}^{(1,5)}_{\leq N}[N_1,N_3] P_{N_2}\squintic[\leqN] \Big\|_{L_t^2 H_x^{-1/2+\delta_2}} \notag \\
&\lesssim N_2^{-1+\delta_2+\epsilon} \| \mathfrak{C}^{(1,5)}_{\leq N}[N_1,N_3]\|_{L_t^\infty} \Big\| \squintic[\leqN] \Big\|_{L_t^\infty \C_x^{1/2-\epsilon}}. \label{pmain:eq-15z-p5}
\end{align}
Using Lemma \ref{analytic:lem-C15-C33}, Lemma \ref{analytic:prop-quintic}, and the frequency-condition $N_{2} \gtrsim N_1,N_3$, it follows that 
\begin{equation*}
    \eqref{pmain:eq-15z-p5} \leq A N_2^{-1+\delta_2+3\epsilon},
\end{equation*}
which is acceptable. \\ 

\emph{Case 3.b: $\zeta^{(2)} \in \big\{ \slinear[green][\leqM], \SXXone, \SXXtwo \big\}$ and $\zeta^{(3)} \in \Symb_{1/2}^p$.} Since this case does not require the \oftt-cancellation, it is slightly easier than Case 3.a. It holds that
\begin{equation}\label{pmain:eq-1Xz-p1}
\Pi^\ast_{\leq N} \Big( \slinear[blue][N_1], P_{N_2} \zeta^{(2)}, P_{N_3} \zeta^{(3)} \Big) 
= \slinear[blue][N_1] P_{N_2} \zeta^{(2)}_{\leq N} P_{N_3} \zeta^{(3)}_{\leq N} - 
\mathbf{1} \Big\{ \zeta^{(3)} = \squintic \Big\} \mathfrak{C}^{(1,5)}_{\leq N}[N_1,N_3] P_{N_2} \zeta^{(2)}_{\leq N}. 
\end{equation}
Since $\zeta^{(2)}$ has the highest frequency in \eqref{pmain:eq-1Xz-p1}, the \oftt-cancellation is not essential in \eqref{pmain:eq-1Xz-p1} and we estimate both summands separately. Using Strichartz and product estimates, the first summand in \eqref{pmain:eq-1Xz-p1} is estimated by 
\begin{align}
\Big\| \, \slinear[blue][N_1] P_{N_2} \zeta^{(2)}_{\leq N} P_{N_3} \zeta^{(3)}_{\leq N} \Big\|_{X^{-1/2+\delta_2,b_+-1}} 
&\lesssim \Big\| \, \slinear[blue][N_1] P_{N_2} \zeta^{(2)}_{\leq N} P_{N_3} \zeta^{(3)}_{\leq N} \Big\|_{L_t^2 H_x^{-1/2+\delta_2}}  \notag \\ 
&\lesssim \Big\| \, \slinear[blue][N_1] P_{N_2} \zeta^{(2)}_{\leq N}  \Big\|_{L_t^\infty \C_x^{-1/2+\delta_1+\epsilon}} \Big\| \zeta^{(3)} \Big\|_{L_t^\infty H_x^{1/2-\delta_1}}.  
\label{pmain:eq-1Xz-p2}
\end{align}
Depending on the choice of $\zeta^{(2)}$, the first factor in \eqref{pmain:eq-1Xz-p2} is estimated using the product estimates for either $\slinear[green][\leqM]$, $\SXXone$, or $\SXXtwo$ (Lemma \ref{analytic:lem-resonant-blue-green} or Proposition \ref{para:prop-resonant}). Depending on the choice of $\zeta^{(3)}$, the second factor in \eqref{pmain:eq-1Xz-p2} is estimated using one of our regularity estimates (Lemma \ref{analytic:lem-linear}, Proposition \ref{analytic:prop-quintic}, or Proposition \ref{ansatz:prop-Xj}).
The second summand  in \eqref{pmain:eq-1Xz-p1} can be bounded directly using Lemma \ref{analytic:lem-C15-C33} and, depending on the choice of $\zeta^{(2)}$, one of our regularity estimates (Lemma \ref{analytic:lem-linear} or Proposition \ref{ansatz:prop-Xj}). \\

\emph{Case 3.c: $\zeta^{(2)} = Y$ and $\zeta^{(3)} \in \Symb_{1/2}^p \backslash \{ Y \}$.} 
For any  $\zeta^{(3)} \in \Symb_{1/2}^p$, it holds that 
\begin{equation}\label{pmain:eq-1Yz-p1}
\Pi^\ast_{\leq N}  \Big( \slinear[blue][N_1], P_{N_2} Y_{\leq N}, P_{N_3} \zeta^{(3)}_{\leq N} \Big) 
= \slinear[blue][N_1] P_{N_2} Y_{\leq N} P_{N_3} \zeta^{(3)}_{\leq N} - \mathbf{1}\Big\{ \zeta^{(3)}=\squintic \Big\} \mathfrak{C}^{(1,5)}_{\leq N}[N_1,N_3] P_{N_2} Y_{\leq N}. 
\end{equation}
Since $Y_{\leq N}$ has the highest frequency in \eqref{pmain:eq-1Yz-p1}, the \oftt-cancellation is not essential in \eqref{pmain:eq-1Yz-p1} and we estimate both summands separately. For the first summand, we obtain from Strichartz estimates that 
\begin{align*}
 \Big\| \slinear[blue][N_1] P_{N_2} Y_{\leq N} P_{N_3} \zeta^{(3)}_{\leq N} \Big\|_{X^{-1/2+\delta_2,b_+-1}} 
 &\lesssim \Big\| \, \slinear[blue][N_1] P_{N_2} Y_{\leq N} P_{N_3} \zeta^{(3)}_{\leq N} \Big\|_{L_t^2 H_x^{-1/2+\delta_2}} \\
 &\lesssim \Big\| \, \slinear[blue][N_1] \Big\|_{L_t^\infty L_x^\infty} \| P_{N_2} Y_{\leq N} \|_{L_t^\infty L_x^2} \big\| P_{N_3} \zeta^{(3)}_{\leq N} \big\|_{L_t^\infty L_x^\infty} \\
 &\lesssim A N_2^{2\epsilon-\delta_2} \| Y_{\leq N} \|_{X^{1/2+\delta_2,b}} 
 \big\| \zeta^{(3)}_{\leq N} \big\|_{L_t^\infty L_x^\infty}. 
\end{align*}
Depending on $\zeta^{(3)}$, the desired conclusion now follows from the probabilistic Strichartz estimates for either $\slinear[green][\leqM]$, $\squintic$, $\SXXone$, or $\SXXtwo$ (Lemma \ref{analytic:lem-linear}, Proposition \ref{analytic:prop-quintic}, or Proposition \ref{ansatz:prop-Xj}). 
The second summand in \eqref{pmain:eq-1Yz-p1} can be estimated directly using Lemma \ref{analytic:lem-C15-C33} and the regularity information on $Y$. \\

\emph{Case 3.d: $\zeta^{(2)}=\zeta^{(3)}=Y$.} In this case, we prove that 
\begin{align*}
\Big\|  \Pi^\ast_{\leq N}  \Big( \slinear[blue][N_1], P_{N_2} Y_{\leq N}, P_{N_3} Y_{\leq N} \Big) \Big\|_{X^{-1/2+\delta_2,b_+-1}} \leq A N_{\textup{max}}^{-\epsilon} \| Y_{\leq N} \|_{X^{1/2+\delta_2,b}}^2. 
\end{align*}
To this end, we use the dyadic decomposition 
\begin{equation*}
 \Big( \Pi^\ast_{\leq N} - \HLL \Big) \Big( \slinear[blue][N_1], P_{N_2} Y_{\leq N}, P_{N_3} Y_{\leq N} \Big) 
 = \sum_{N_0} P_{N_0} \Big( \slinear[blue][N_1] \,  P_{N_2} Y_{\leq N} \,  P_{N_3} Y_{\leq N} \Big). 
\end{equation*}
Depending on the relative value of $N_0$ and $N_3$, we now utilize one of the following two deterministic Strichartz estimates. In the case $N_0 \lesssim N_3$, we estimate 
\begin{align*}
&\Big\| \Big( \Pi^\ast_{\leq N} - \HLL \Big) \Big( \slinear[blue][N_1], P_{N_2} Y_{\leq N}, P_{N_3} Y_{\leq N} \Big) \Big\|_{X^{-1/2+\delta_2,b_+-1}}\\
\lesssim& \, N_0^{1/2+\delta_2+\epsilon} \Big\| \slinear[blue][N_1] P_{N_2} Y_{\leq N} P_{N_3} Y_{\leq N} \Big\|_{L_t^2 L_x^1} \\ 
\lesssim& \,  N_0^{1/2+\delta_2+\epsilon} \Big\| \, \slinear[blue][N_1] \Big\|_{L_t^\infty L_x^\infty} \| P_{N_2} Y_{\leq N} \|_{L_t^\infty L_x^2} \| P_{N_3} Y_{\leq N} \|_{L_t^\infty L_x^2} \\
\lesssim& \,  A \Big(N_0/N_3\Big)^{1/2+\delta_2} N_{\textup{max}}^{2\epsilon-\delta_2} \| Y_{\leq N} \|_{X^{1/2+\delta_2,b}}^2. 
\end{align*}
In the case $N_0 \gtrsim N_3$, we estimate  
\begin{align*}
&\Big\| \Big( \Pi^\ast_{\leq N} - \HLL \Big) \Big( \slinear[blue][N_1], P_{N_2} Y_{\leq N}, P_{N_3} Y_{\leq N} \Big) \Big\|_{X^{-1/2+\delta_2,b_+-1}}\\
\lesssim& \, N_0^{-1/2+\delta_2} \Big\| \slinear[blue][N_1] P_{N_2} Y_{\leq N} P_{N_3} Y_{\leq N} \Big\|_{L_t^2 L_x^2} \\ 
\lesssim&\,  N_0^{-1/2+\delta_2}  \Big\| \, \slinear[blue][N_1] \Big\|_{L_t^\infty L_x^\infty} \| P_{N_2} Y_{\leq N} \|_{L_t^\infty L_x^2} \| P_{N_3} Y_{\leq N} \|_{L_t^2 L_x^\infty} \\
\lesssim& \, A \Big(N_0/N_3\Big)^{-1/2+\delta_2} N_{\textup{max}}^{2\epsilon-\delta_2} \| Y_{\leq N} \|_{X^{1/2+\delta_2,b}}^2. 
\end{align*}
In both cases, this yields an acceptable contribution. 
\end{proof}

We now turn to the proof of Proposition \ref{ansatz:prop-Y-no}, which is the last main estimate from Subsection \ref{section:main-estimates}. 

\begin{proof}[Proof of Proposition \ref{ansatz:prop-Y-no}:]
To simplify the notation, we focus only on $T=1$ and $\Jc=[-1,1]$, since the general case follows from minor modifications. Furthermore, we restrict all norms to $[-1,1]\times \T^3$. 
In the following, each estimate holds on $A$-certain events $E_A \in \mathcal{E}$, which will not be repeated below. As in the proof of Proposition \ref{ansatz:prop-Y-one}, we replace the $A$-factor in our desired estimate by $A^3$. In the statement of the proposition, we are concerned with terms of the form
\begin{equation}\label{pmain:eq-no-recalled}
    \Pi^\ast_{\leq N} \Big( P_{\leq N} \zeta^{(1)}, P_{\leq N} \zeta^{(2)}, P_{\leq N} \zeta^{(3)}  \Big), \qquad  \mathfrak{C}_{\leq N} \zeta_{\leq N}, \qquad \text{or} \qquad \big( \gamma_{\leq N} - \Gamma_{\leq N} \big) \zeta_{\leq N}, 
\end{equation}
where $\zeta^{(1)},\zeta^{(2)}, \zeta^{(3)}\in \Symb_0^{p}$ and $\zeta \in \Symb_0^b$. We first address the second and third terms in \eqref{pmain:eq-no-recalled}, whose estimates are rather simple. Using $X^{-1/2+\delta_2,0}\hookrightarrow X^{-1/2+\delta_2,b_+-1}$, Lemma \ref{diagrams:lem-mathfrakC}, and Lemma \ref{analytic:lem-difference-gamma}, we obtain for all $\zeta \in \Symb_0^b$ that 
\begin{align*}
&\big\| \mathfrak{C}_{\leq N} \zeta_{\leq N} \big\|_{X^{-1/2+\delta_2,b_+-1}} 
+ \big\| \big( \gamma_{\leq N} - \Gamma_{\leq N} \big) \zeta_{\leq N} \big\|_{X^{-1/2+\delta_2,b_+-1}} \\
\lesssim& \,\big\| \mathfrak{C}_{\leq N} \zeta_{\leq N} \big\|_{L_t^2 H_x^{-1/2+\delta_2}} 
+ \big\| \big( \gamma_{\leq N} - \Gamma_{\leq N} \big) \zeta_{\leq N} \big\|_{L_t^2 H_x^{-1/2+\delta_2}}  \\ 
\lesssim& \, \| \zeta_{\leq N} \|_{L_t^\infty H_x^{-1/2+\delta_2 + \epsilon}}.  
\end{align*}
The desired conclusion now follows from our regularity estimates for elements of $\Symb_0^b$. As a result, it remains to treat the first term in \eqref{pmain:eq-no-recalled}. To this end, we distinguish cases where the arguments $\zeta^{(1)},\zeta^{(2)}$, and $\zeta^{(3)}$ contain three, two, one, or zero cubic stochastic objects, respectively. \\

\emph{Case 1: $\zeta^{(1)}=\zeta^{(2)}=\zeta^{(3)}=\scubic$.} This case corresponds to the explicit nonic stochastic object, which has been estimated in Proposition \ref{analytic:prop-higher-order}. \\

\emph{Case 2: $\zeta^{(1)}=\zeta^{(2)}=\scubic$ and $\zeta^{(3)}\in \Symb_{1/2}^p$.} 
In this case, it holds that 
\begin{equation*}
\Pi^\ast_{\leq N} \Big( \scubic[\leqN], \scubic[\leqN], \zeta^{(3)}_{\leq N} \Big) 
= \Big( \Big( \scubic[\leqN] \Big)^2 - \mathfrak{C}^{(3,3)}_{\leq N} \Big) \zeta^{(3)}_{\leq N}. 
\end{equation*}
It follows that
\begin{align*}
\bigg\| \Big( \Big( \scubic[\leqN] \Big)^2 - \mathfrak{C}^{(3,3)}_{\leq N} \Big) \zeta^{(3)}_{\leq N} \bigg\|_{X^{-1/2+\delta_2,b_+-1}} 
&\lesssim \, \bigg\| \Big( \Big( \scubic[\leqN] \Big)^2 - \mathfrak{C}^{(3,3)}_{\leq N} \Big) \zeta^{(3)}_{\leq N} \bigg\|_{L_t^2 H_x^{-1/2+\delta_2}} \\
&\lesssim \,  \bigg\| \Big( \scubic[\leqN] \Big)^2 - \mathfrak{C}^{(3,3)}_{\leq N} \bigg\|_{L_t^\infty \C_x^{-1/2+\delta_1}} \Big\|  \zeta^{(3)}_{\leq N} \Big\|_{L_t^\infty H_x^{1/2-\delta_1}}.
\end{align*}
The first factor has been estimated in Proposition \ref{analytic:prop-sextic1}. The second factor can be bounded using our regularity estimates for elements of $\Symb_{1/2}^p$. \\

\emph{Case 3: $\zeta^{(1)}=\scubic$ and $\zeta^{(2)},\zeta^{(3)}\in \Symb_{1/2}^p$.}
In this case, we first decompose
\begin{equation}\label{pmain:eq-Y-no-p1}
\Pi^\ast_{\leq N}\Big( P_{N_1} \scubic[\leqN], \zeta^{(2)}_{\leq N}, \zeta^{(3)}_{\leq N}\Big) 
= \sum_{N_0} \sum_{N_1,N_2,N_3 \leq N} P_{N_0} \Big( 
 \scubic[\leqN] P_{N_2} \zeta^{(2)}_{\leq N} P_{N_3} \zeta^{(3)}_{\leq N}\Big).
\end{equation}
We now estimate the dyadic components in \eqref{pmain:eq-Y-no-p1}. Due to the symmetry in $\zeta^{(2)}$ and $\zeta^{(3)}$, we may assume that either $\zeta^{(3)} \in \Symb_{1/2}^p \backslash \{ Y \}$ or $\zeta^{(2)}=\zeta^{(3)}=Y$. In the second case, we can further assume that $N_2 \geq N_3$.

\emph{Case 3.a:  $\zeta^{(1)}=\scubic$, $\zeta^{(2)}\in \Symb_{1/2}^p$, and $\zeta^{(3)} \in \Symb_{1/2}^p \backslash \{ Y \}$.} 
Using Strichartz estimates, product estimates, and Lemma \ref{analytic:lem-cubic}, we first estimate 
\begin{align}
 &\Big\|  P_{N_0} \Big(  P_{N_1}
 \scubic[\leqN] P_{N_2} \zeta^{(2)}_{\leq N} P_{N_3} \zeta^{(3)}_{\leq N}\Big) \Big\|_{X^{-1/2+\delta_2,b_+-1}} \notag \\
 \lesssim& \, N_0^{-1/2+\delta_2+2\epsilon} \Big\|    P_{N_1}
 \scubic[\leqN] P_{N_2} \zeta^{(2)}_{\leq N} P_{N_3} \zeta^{(3)}_{\leq N} \Big\|_{L_t^2 H_x^{-2\epsilon}}\notag \\
 \lesssim& \, N_0^{-1/2+\delta_2+2\epsilon} N_1^{-\epsilon} \max(N_2,N_3)^{2\epsilon} \Big\| P_{N_1} \scubic[\leqN] \Big\|_{L_t^\infty \C_x^{-\epsilon}} 
 \| P_{N_2} \zeta^{(2)}_{\leq N} \|_{L_t^\infty L_x^2} \| P_{N_3} \zeta^{(3)}_{\leq N} \|_{L_t^2 L_x^\infty} \notag \\
 \lesssim& \,  A   N_0^{-1/2+\delta_2+2\epsilon} N_1^{-\epsilon} \max(N_2,N_3)^{2\epsilon} 
 \| P_{N_2} \zeta^{(2)}_{\leq N} \|_{L_t^\infty L_x^2} \| P_{N_3} \zeta^{(3)}_{\leq N} \|_{L_t^2 L_x^\infty}. \label{pmain:eq-Y-no-p2}
\end{align}
Since $\zeta\in \Symb_{1/2}^p$ and $\zeta^{(3)} \in \Symb_{1/2}^p \backslash \{ Y \}$, we obtain from regularity and probabilistic Strichartz estimates that 
\begin{equation}\label{pmain:eq-Y-no-p3}
  \| P_{N_2} \zeta^{(2)}_{\leq N} \|_{L_t^\infty L_x^2} \| P_{N_3} \zeta^{(3)}_{\leq N} \|_{L_t^2 L_x^\infty} \lesssim A^2 N_2^{-1/2+\delta_1} N_3^{-1/2+\delta_1} \Big(1 + \|v_{\leq N} \|_{X^{-1,b}}^4 + \| Y_{\leq N} \|_{X^{1/2+\delta_2,b}}^2 \Big).  
\end{equation}
After re-inserting \eqref{pmain:eq-Y-no-p3} into \eqref{pmain:eq-Y-no-p2}, we obtain an acceptable contribution. \\

\emph{Case 3.b: $\zeta^{(1)}=\scubic$, $\zeta^{(2)}=\zeta^{(3)}=Y$, and $N_2 \geq N_3$.} In this case, we further distinguish two-subcases, which depend on the relative size of $N_0$ and $N_3$. \\

\emph{Case 3.b.i: $\zeta^{(1)}=\scubic$, $\zeta^{(2)}=\zeta^{(3)}=Y$, and $N_0,N_2 \geq N_3$.} Arguing exactly as in \eqref{pmain:eq-Y-no-p2}, we obtain that 
\begin{equation}\label{pmain:eq-Y-no-p4}
\begin{aligned}
    &\Big\|  P_{N_0} \Big( 
 \scubic[\leqN] P_{N_2} \zeta^{(2)}_{\leq N} P_{N_3} \zeta^{(3)}_{\leq N}\Big) \Big\|_{X^{-1/2+\delta_2,b_+-1}} \\
 \lesssim&\, A N_0^{-1/2+\delta_2 + 2\epsilon} N_1^{-\epsilon} N_2^{2\epsilon} \big\| P_{N_2} Y_{\leq N} \big\|_{L_t^\infty L_x^2}
 \big\| P_{N_2} Y_{\leq N} \big\|_{L_t^2 L_x^\infty}.
 \end{aligned}
\end{equation}
Using deterministic Strichartz estimates and our assumption $N_0,N_2 \geq N_3$, we obtain that 
\begin{align*}
\eqref{pmain:eq-Y-no-p4} &\lesssim A N_0^{-1/2+\delta_2+2 \epsilon} N_1^{-\epsilon} N_2^{2\epsilon-1/2-\delta_2} N_3^{1/2+\epsilon-\delta_2} \\
&\lesssim A N_1^{-\epsilon} N_2^{-1/2-\delta+5\epsilon},
\end{align*}
which is acceptable.\\ 

\emph{Case 3.b.ii: $\zeta^{(1)}=\scubic$, $\zeta^{(2)}=\zeta^{(3)}=Y$, and $N_2 \geq N_3>N_0$.} Using Lemma \ref{analytic:lem-cubic}, we estimate 
\begin{align*}
 &   \Big\|  P_{N_0} \Big(  P_{N_1}
 \scubic[\leqN] P_{N_2} \zeta^{(2)}_{\leq N} P_{N_3} \zeta^{(3)}_{\leq N}\Big) \Big\|_{X^{-1/2+\delta_2,b_+-1}}  \\
 \lesssim& \,  N_0^{1/2+\delta_2+\epsilon} \Big\|    P_{N_1}
 \scubic[\leqN] P_{N_2} Y_{\leq N} P_{N_3} Y_{\leq N} \Big\|_{L_t^2 L_x^1} \\
 \lesssim& \,  N_0^{1/2+\delta_2+\epsilon} \Big\|    P_{N_1}
 \scubic[\leqN]  \Big\|_{L_t^\infty L_x^\infty} \big\| P_{N_2} Y_{\leq N} \big\|_{L_t^\infty L_x^2} \big\|  P_{N_3} Y_{\leq N} \big\|_{L_t^\infty L_x^2} \\
 \lesssim& \, A N_0^{1/2+\delta_2+\epsilon} N_1^\epsilon N_2^{-1/2-\delta_2} N_3^{-1/2-\delta_2} \| Y_{\leq N}\|_{X^{1/2+\delta_2,b}}^2. 
\end{align*}
Due to our frequency-restrictions, it holds that 
\begin{equation*}
N_0^{1/2+\delta_2+\epsilon} N_1^\epsilon N_2^{-1/2-\delta_2} N_3^{-1/2-\delta_2} 
\lesssim N_1^{\epsilon} N_2^{-1/2-\delta_2+\epsilon} \lesssim N_{\textup{max}}^{-1/2-\delta_2+2\epsilon}. 
\end{equation*}~

\emph{Case 4: $\zeta^{(1)}, \zeta^{(2)}, \zeta^{(3)} \in \Symb_{1/2}^p$.}
If $\zeta^{(1)} \in \Symb_{1/2}^p\backslash \{ Y \}$, we can use probabilistic Strichartz estimates and the same argument as in Case 3. By symmetry, it only remains to treat the case $\zeta^{(1)}=\zeta^{(2)}=\zeta^{(3)}=Y$. Since the cubic nonlinear wave equation is deterministically well-posed at regularities greater than or equal to $1/2$, this case can be treated using the standard well-posedness argument in $X^{1/2+\delta_2,b}$. 
\end{proof}

Finally, we prove Lemma \ref{global:lem-one-zero-revisited}, which is a variant of Proposition \ref{ansatz:prop-Y-one} and Proposition \ref{ansatz:prop-Y-no}. 

\begin{proof}[Proof of Lemma \ref{global:lem-one-zero-revisited}]
To simplify the notation, we focus only on $T=1$ and $\Jc=[-1,1]$, since the general case follows from minor modifications, and we restrict all norms to $[-1,1]\times \T^3$. 
In the following, each estimate holds on $A$-certain events $E_A \in \mathcal{E}$, which will not be repeated below. 

The argument is a minor modification of the proof of Proposition \ref{ansatz:prop-Y-one} and Proposition \ref{ansatz:prop-Y-no}. 
We first recall from \eqref{global:eq-one-zero-decomp} that 
\begin{equation}\label{pmain:eq-modification-q1}
v_{\leq N} = \XXone[v_{\leq N}, Y_{\leq N}] + \XXtwo[v_{\leq N}] + Y_{\leq N}. 
\end{equation}

If $\zeta^{(1)} \in \Symb^b$ and $\zeta^{(2)},\zeta^{(3)}\in \Symb_0^b$ contain at most one $v$-term, we directly insert the decomposition \eqref{pmain:eq-modification-q1} and argue as in the proof of Proposition \ref{ansatz:prop-Y-one} or Proposition \ref{ansatz:prop-Y-no}, which yields the desired estimate. Thus, we can assume that 
\begin{equation}\label{pmain:eq-modification-q2}
\zeta^{(2)}=\zeta^{(3)}= v. 
\end{equation}
We now distinguish the cases $\zeta^{(1)}\in \Symb^b \backslash\{v\}$ and $\zeta^{(1)}=v$. ~\\

\emph{Case 1: $\zeta^{(1)} \in \Symb^b \backslash\{v\}$.} 
If $\zeta^{(1)} = \slinear[blue]$, we have to estimate 
\begin{equation}\label{pmain:eq-modification-q3}
\Big( \Pi^\ast_{\leq N} - \HLL \Big) \Big( \slinear[blue][\leqN], v_{\leq N}, v_{\leq N} \Big) \\
= \sum_{\substack{N_1,N_2,N_3 \leq N \colon \\  \max(N_2,N_3) > N_1^\eta }} 
\slinear[blue][N_1] P_{N_2} v_{\leq N} P_{N_3} v_{\leq N}. 
\end{equation}
Alternatively, if $\zeta^{(1)} \in \Symb^b_0 \backslash\{v\}$, we have to estimate
\begin{equation}\label{pmain:eq-modification-q4}
\Pi^\ast_{\leq N}  \Big( \zeta_{\leq N}^{(1)}, v_{\leq N}, v_{\leq N} \Big) \\
= \sum_{\substack{N_1,N_2,N_3 \leq N  }} 
P_{N_1}\zeta_{\leq N}^{(1)}  P_{N_2} v_{\leq N} P_{N_3} v_{\leq N}. 
\end{equation}
We now treat the dyadically localized terms in  \eqref{pmain:eq-modification-q3} and \eqref{pmain:eq-modification-q4}.  Using symmetry in the second and third argument, we can reduce to the case $N_2 \geq N_3$. By using the decomposition \eqref{pmain:eq-modification-q1} and the definition of the nonlinear smoothing norm (Definition \ref{ansatz:def-nonlinear-smoothing}), we can write 
\begin{align}
P_{N_2} v_{\leq N} &=  P_{N_2} \XXone[v_{\leq N}, Y_{\leq N}] + P_{N_2} \XXtwo[v_{\leq N}] + P_{N_2} Y_{\leq N},  \label{pmain:eq-modification-p3} \\
P_{N_3} v_{\leq N} &=  P_{N_3} v_{\leq N}^{\str} + P_{N_3} v_{\leq N}^{\reg}, \label{pmain:eq-modification-p4} 
\end{align}
where 
\begin{equation*}
\big\| v_{\leq N}^{\str} \big\|_{L_t^\infty \C_x^{-\epsilon}\cap X^{-\epsilon,b}}
+ \big\| v_{\leq N}^{\reg} \big\|_{X^{1/2+\delta_2,b}} \leq 2 \big\| v_{\leq N} \big\|_{\NSN([0,T])}.
\end{equation*}
The superscripts in \eqref{pmain:eq-modification-p4} stand for ``Strichartz" and ``regularity", respectively. Using \eqref{pmain:eq-modification-p3} and \eqref{pmain:eq-modification-p4}, we decompose the dyadically localized terms in \eqref{pmain:eq-modification-q3} or \eqref{pmain:eq-modification-q4} as 
\begin{align}
&P_{N_1} \zeta_{\leq N}^{(1)} P_{N_2} v_{\leq N} P_{N_3} v_{\leq N} \notag \\
=& \, P_{N_1} \zeta_{\leq N}^{(1)} P_{N_2}\Big( \XXone + \XXtwo \Big) P_{N_3} v_{\leq N}^{\reg} 
\label{pmain:eq-modification-q5} \\
+& \, P_{N_1} \zeta_{\leq N}^{(1)} P_{N_2} Y_{\leq N} P_{N_3} v_{\leq N}^{\reg} 
\label{pmain:eq-modification-q6} \\
+& \, P_{N_1} \zeta_{\leq N}^{(1)} P_{N_2}\Big( \XXone + \XXtwo + Y_{\leq N} \Big) P_{N_3} v_{\leq N}^{\str}
\label{pmain:eq-modification-q7}.
\end{align}
The estimates of \eqref{pmain:eq-modification-q5} and \eqref{pmain:eq-modification-q6} are exactly as in the proof of Proposition \ref{ansatz:prop-Y-one} and Proposition \ref{ansatz:prop-Y-no}. Indeed, if $\zeta^{(1)}= \slinear[blue]$, the estimates of \eqref{pmain:eq-modification-q5} and \eqref{pmain:eq-modification-q6} are as in Case 3.b or Case 3.d of the proof of Proposition \ref{ansatz:prop-Y-one}, respectively. If $\zeta^{(1)}\in\Symb^b_0 \backslash\{v\}$, the estimates of \eqref{pmain:eq-modification-q5} and \eqref{pmain:eq-modification-q6} are as in Case 3 or Case 4 in the proof of Proposition \ref{ansatz:prop-Y-no}.

It remains to treat the term \eqref{pmain:eq-modification-q7}. Since the case $\zeta^{(1)} \in \Symb_0^b \backslash\{v\}$ is trivial, we focus on the case $\zeta^{(1)} =\slinear[blue]$.  The contributions from $\XXone$ and $\XXtwo$ can be bounded directly using Corollary \ref{para:cor-blueXj}. The remaining contribution of $Y_{\leq N}$ can be estimated by 
\begin{align*}
&\big\|  \, \slinear[blue][N_1] P_{N_2} Y_{\leq N}  P_{N_3} v_{\leq N}^{\str}\big\|_{X^{-1/2+\delta_2,b_+-1}} \\
\lesssim& \, \big\| \,  \slinear[blue][N_1] P_{N_2} Y_{\leq N}  P_{N_3} v_{\leq N}^{\str}\big\|_{L_t^2 H_x^{-1/2+\delta_2}}  \\
\lesssim& \, \Big( N_1^{-1/2+\delta_2} \mathbf{1}\big\{ N_1 \gg N_2 \big\}  + \mathbf{1}\big\{ N_1 \lesssim N_2 \big\} \Big) \big\|  \, \slinear[blue][N_1] \big\|_{L_t^\infty L_x^\infty} \big\| P_{N_2} Y_{\leq N} \big\|_{L_t^2 L_x^2} \big\|    P_{N_3} v_{\leq N}^{\str} \big\|_{L_t^\infty L_x^\infty} \\
\lesssim& \, A  \Big( N_1^{-1/2+\delta_2} \mathbf{1}\big\{ N_1 \gg N_2 \big\}  + \mathbf{1}\big\{ N_1 \lesssim N_2 \big\} \Big)  N_1^{1/2+\epsilon} N_2^{-1/2-\delta_2} N_3^\epsilon 
\|Y_{\leq N} \|_{X^{1/2+\delta_2,b}} \| v_{\leq N} \|_{\NSN}.  
\end{align*}
Since \eqref{pmain:eq-modification-q3} is restricted to frequency-scales satisfying $N_2 = \max(N_2,N_3) \geq N_1^\eta$, this yields an acceptable contribution. \\

\emph{Case 2: $\zeta^{(1)}=v$.} Due to \eqref{pmain:eq-modification-q2}, it holds that  $\zeta^{(1)}=\zeta^{(2)}=\zeta^{(3)}=v$.  As before, we utilize the dyadic decomposition 
\begin{equation*}
\Pi^\ast_{\leq N}  \Big( v_{\leq N}, v_{\leq N}, v_{\leq N} \Big) \\
= \sum_{\substack{N_1,N_2,N_3 \leq N }} 
P_{N_1}v_{\leq N} P_{N_2} v_{\leq N} P_{N_3} v_{\leq N}. 
\end{equation*}
By symmetry, we can assume that $N_1 \geq N_2 \geq N_3$. Using the decomposition \eqref{pmain:eq-modification-q1} and the nonlinear smoothing norm, it remains to estimate
\begin{equation}\label{pmain:eq-modification-p1}
    P_{N_1} \big( \XXone + \XXtwo + Y_{\leq N} \big) P_{N_2} \big( v_{\leq N}^{\str} + v_{\leq N}^{\reg} \big) P_{N_3} \big( v_{\leq N}^{\str} + v_{\leq N}^{\reg} \big). 
\end{equation}
All possible combinations in the product \eqref{pmain:eq-modification-p1} can be estimated easily using Strichartz estimates. For example, 
\begin{align*}
&\Big\|  P_{N_1} \big( \XXone + \XXtwo + Y_{\leq N} \big) P_{N_2}  v_{\leq N}^{\str}   
P_{N_3}  v_{\leq N}^{\str} \Big\|_{X^{-1/2+\delta_2,b_+-1}} \\
\lesssim&\, \Big\|  P_{N_1} \big( \XXone + \XXtwo + Y_{\leq N} \big) P_{N_2}  v_{\leq N}^{\str}   
P_{N_3}  v_{\leq N}^{\str} \Big\|_{L_t^2 L_x^2} \\
\lesssim&\, \Big\|  P_{N_1} \big( \XXone + \XXtwo + Y_{\leq N} \big) \Big\|_{L_t^\infty L_x^2} \Big\| P_{N_2}  v_{\leq N}^{\str}   \Big\|_{L_t^\infty L_x^\infty}
\Big\| P_{N_3}  v_{\leq N}^{\str} \Big\|_{L_t^\infty L_x^\infty} \\
\lesssim&\, N_1^{-1/2+\delta_1 + 2\epsilon}   \Big\|  P_{N_1} \big( \XXone + \XXtwo + Y_{\leq N} \big) \Big\|_{L_t^\infty H_x^{1/2-\delta_1}} \big\| v_{\leq N} \big\|_{\NSN}^2. 
\end{align*}
The desired conclusion now follows from the regularity estimate for $\XXone$ and $\XXtwo$ (Proposition \ref{ansatz:prop-Xj}). We omit the (standard) details for the remaining terms in the product \eqref{pmain:eq-modification-p1}. 
\end{proof}


\begin{appendix}

\section{The nonlinear stochastic heat equation with sharp frequency-cutoffs}\label{section:heat-appendix}

In Subsection \ref{section:diagram-parabolic}, we discussed the frequency-truncated stochastic nonlinear heat equation \eqref{diagram:eq-truncated-SNLH}, i.e.,
\begin{equation}\label{heat:eq-heat}
\begin{cases}
\big(\partial_s + 1 - \Delta \big) \Phi_{\leq N}^{\cos}  = -  P_{\leq N}  \Big( \lcol \big(P_{\leq N} \Phi_{\leq N}^{\cos} \big)^3 \rcol  + \gamma_{\leq N}  \Phi^{\cos}_{\leq N} \Big)+ \sqrt{2} \dW[\cos] \quad  (s,x) \in (s_0,\infty) \times \T^3,\\
\Phi^{\cos}_{\leq N} \big|_{s=s_0}= \phi^{\cos}. 
\end{cases}
\end{equation}
In the proof of Proposition \ref{ansatz:prop-caloric} at the end of Subsection \ref{section:diagram-parabolic}, we claimed that the local well-posedness of \eqref{heat:eq-heat} essentially follows from the previous literature \cite{CC18,GIP15,H14}. Due to the sharp frequency-cutoffs in \eqref{heat:eq-heat}, however, some modifications are necessary. \\

The sharp frequency-projections $(P_{\leq N})_N$, which are based on cubes, are uniformly bounded on $L_x^p$ for every $1<p<\infty$ but unbounded on $L_x^\infty$. However, for every $\delta>0$, $(P_{\leq N})_N$ is uniformly bounded as a map from $\C_x^\delta$ to $L_x^\infty$. More generally, for every $\alpha \in \R$ and $\delta >0$,  $(P_{\leq N})_N$ is uniformly bounded as a map from $\C_x^{\alpha+\delta}$ to $\C_x^\alpha$. Due to the smoothing properties of heat equations and the sub-criticality of \eqref{heat:eq-heat}, the resulting $\delta$-loss is acceptable, and most of the argument in \cite{CC18} applies verbatim to \eqref{heat:eq-heat}. \\

The main technical difficulty concerns commutator terms, which do not exhibit any gains under the sharp frequency-cutoffs. To illustrate this, let $N\geq 1$, let $n=(N,0,0) \in \Z^3$, and $m=(1,0,0)\in \Z^3$. Then, it holds that 
\begin{equation}\label{heat:eq-commutator-failure}
\big[ P_{\leq N}, e^{i\langle m,x\rangle}\big] e^{i\langle n, x \rangle}
= P_{\leq N}\big( e^{i \langle m, x \rangle}  e^{i\langle n, x \rangle} \big) - 
 e^{i \langle m, x \rangle} P_{\leq N}\big(  e^{i\langle n, x \rangle} \big) 
 = - e^{i \langle m+n , x \rangle}. 
\end{equation}
In particular, \eqref{heat:eq-commutator-failure} exhibits no gain in $N$. Due to the missing commutator estimates, a few steps in \cite{CC18} do not directly carry over to \eqref{heat:eq-heat}. For example, the decomposition in \cite[(3.3)]{CC18} has no direct counterpart in our setting. In the rest of this section, we show how the main estimate \cite[Proposition 3.8]{CC18} can be extended to our setting, but do not discuss any other (more minor) modifications. 

\subsection{Preparations}
Let $0<\epsilon \ll \delta \ll 1$ be parameters. 
In \eqref{diagram:eq-Hlin}, we defined the heat propagator $e^{-s (1-\Delta)}$ as the Fourier multiplier with symbol $n\mapsto e^{- s \langle n \rangle^2}$. For any initial time $s_0 \in \R$, we define the associated Duhamel integral as 
\begin{equation*}
\HeatDuh_{s_0} \big[ F \big](s) := \int_{s_0}^s \ds^\prime e^{-(s-s^\prime)(1-\Delta)} F(s^\prime).
\end{equation*}

We now recall the following Schauder-type estimate.

\begin{lemma}[Schauder-type estimate]\label{heat:lem-duhamel}
For any $\alpha \in \R$, $\theta \in (0,1)$, $s_0 \in [-1,0]$, and $F\colon [s_0,0]\times \T^3 \rightarrow \R$, it holds that 
\begin{equation}
\sup_{s\in [s_0,1]} |s-s_0|^\theta \, \big\| \HeatDuh_{s_0} \big[ F \big] \big\|_{\C_x^{\alpha+2-2\delta}} \lesssim \sup_{s\in [s_0,1]} |s-s_0|^\theta \,  \big\| F \big\|_{\C_x^{\alpha}}.
\end{equation}
Furthermore, if $F(s_0)=0$, we also have that 
\begin{equation}
\sup_{s\in [s_0,1]} |s-s_0|^\theta \, \big\| \HeatDuh_{s_0} \big[ \partial_s F \big] \big\|_{\C_x^{\alpha-2\delta}} \lesssim \sup_{s\in [s_0,1]} |s-s_0|^\theta \,  \big\| F \big\|_{\C_x^{\alpha}}.
\end{equation}
\end{lemma}

Except for a boundary term in $\Duh_{s_0}[\partial_s F]$, the $\delta$-loss in spatial derivatives can be used to obtain better weights in $|s-s_0|$. 

\begin{proof}
The estimates essentially follow from \cite[Proposition 2.7]{CC18}. For the second estimate, we also use the identity
\begin{equation*}
\HeatDuh_{s_0}\big[ F \big] = F - \HeatDuh_{s_0} \big[ (1-\Delta) F \big]. 
\end{equation*}
\end{proof}

Furthermore, we define the following bilinear para-products\footnote{In the definition of the para-products, we are not forced to choose the sharp frequency-projections and could have used a different (smooth) frequency-scale decomposition.}.

\begin{definition}[Bilinear para-products]
For any $f,g \colon \T^3 \rightarrow \R$, we define 
\begin{equation*}
f \parall g := \sum_{\substack{K,L\colon \\ K \ll L}} P_K f \, P_L g, \qquad 
f \parasim g := \sum_{\substack{K,L\colon \\ K \sim L}} P_K f \, P_L g, \qquad \text{and} \qquad 
f \paragg g := \sum_{\substack{K,L\colon \\ K \gg L}} P_K f \, P_L g. 
\end{equation*}
\end{definition}

\subsection{Main estimate}

We now prove the analogue of the main estimate \cite[Proposition 3.8]{CC18} with sharp frequency-cutoffs. In order to precisely state the estimate, we use the Wick-ordered square
\begin{equation*}
\shquadratic[\leqN][\cos] = \, \lcol \Big( \shlinear[\leqN][\cos] \Big)^2 \rcol . 
\end{equation*}

\begin{proposition}\label{heat:prop-main}
For all $A\geq 1$, there exists an $A$-certain event $E_A \in \mathcal{E}$ on which the following estimate holds: For all $\theta \in (0,1), N \geq 1$, $s_0\in[-1,0)$, and $\Psi \colon [s_0,0] \times \T^3 \rightarrow \R$, we have that 
\begin{equation}\label{heat:eq-main-estimates}
\begin{aligned}
   &\sup_{s\in [s_0,0]} |s-s_0|^{\theta} \, 
   \bigg\| \HeatDuh_{s_0} \bigg[ 9 \shquadratic[\leqN][\cos] \parasim \big( P_{\leq N} \HeatDuh_{s_0} \big) \Big[ \shquadratic[\leqN][\cos] \paragg P_{\leq N} \Psi \Big] - \gamma_{\leq N} P_{\leq N} \Psi \bigg]  \bigg\|_{\C_x^{3/2-6\delta}} \\
   \leq&  \,  A \sup_{s\in [s_0,0]} |s-s_0|^{\theta} \big\| \Psi (s) \big\|_{H_x^{1-\delta}}. 
   \end{aligned}
\end{equation}
\end{proposition}

\begin{remark} We make a few comments on Proposition \ref{heat:prop-main}. 
\begin{enumerate}[label=(\roman*)]
\item If instead of the nonlinear remainder $\Psi$ the low-frequency component is given by the cubic stochastic object 
\begin{equation*}
    \shcubic[\leqN][\cos],
\end{equation*}
the term can be treated as an explicit stochastic object and estimated more easily. 
\item  The regularity condition on $\Psi$ on the right-hand side of \eqref{heat:eq-main-estimates} is far from optimal, but sufficient for our applications to \eqref{heat:eq-heat}. It is the result of a simple (but crude) Sobolev-type embedding. 
\item The local well-posedness theory of \eqref{heat:eq-heat} also requires a minor variant of \eqref{heat:eq-main-estimates}, where the Duhamel integral is estimated at a lower spatial regularity but with better pre-factors of $|s-s_0|$. 
\end{enumerate}
\end{remark}

\begin{proof}
We will prove the estimate \eqref{heat:eq-main-estimates} with $\C_x^{1-\delta}$ on the right-hand side replaced by $H_x^{1-\delta}$, which is a stronger estimate. In the following, we denote the frequencies in the first and second factor of 
\begin{equation*}
    \shquadratic[\leqN][\cos]
\end{equation*}
by $(n_1,n_2)$ and $(n_3,n_4)$, respectively. We also denote the frequency of $\Psi$ by $n_5$. Inserting the definitions of the stochastic objects, Duhamel integral, and para-products, it follows that 
\begin{equation}
\begin{aligned}
9 \shquadratic[\leqN][\cos] \parasim P_{\leq N} \HeatDuh_{s_0} \Big[ \shquadratic[\leqN][\cos] \paragg P_{\leq N} \Psi \Big] 
= \sum_{N_5 \leq N} \sum_{n_5 \in \Z^3} 
\int_{s_0}^s \ds^\prime \widehat{P_{N_5} \Psi}(s^\prime,n_5) \, \Gc_{\leq N}(s,s^\prime,x;n_5,N_5), 
\end{aligned}
\end{equation}
where

\begin{equation}\label{heat:eq-G}
\begin{aligned}
&\Gc_{\leq N}(s,s^\prime,x;n_5,N_5) \\
=& 9 \cdot 2^{4/2}  
\sum_{n_1,n_2,n_3,n_4 \in \Z^3} 
\sum_{ \substack{N_{12},N_{345}\colon \\ N_{12} \sim N_{345}, \\ N_{345}\leq N }}
\sum_{\substack{ N_{34} \colon \\ N_{34} \gg N_5}} \bigg[ 
e^{i \langle n_{12345}, x \rangle} 
1_{N_{12}}(n_{12}) 1_{N_{345}}(n_{345}) 1_{N_{34}}(n_{34}) 1_{N_5}(n_5) \\
\times&  \Big( \prod_{1\leq j \leq 4} 1_{\leq N}(n_j) \Big) e^{-(s-s^\prime) \langle n_{345} \rangle^2} 
\Big( \int_{-\infty}^s \int_{-\infty}^s \dW[\cos][s_1][n_1] \dW[\cos][s_2][n_2] 
\prod_{j=1,2} e^{-(s-s_j) \langle n_j \rangle^2} \Big) \\
\times& \Big( \int_{-\infty}^{s^\prime} \int_{-\infty}^{s^\prime} \dW[\cos][s_3][n_3] \dW[\cos][s_4][n_4]
\prod_{j=3,4} e^{-(s^\prime-s_j) \langle n_j \rangle^2} \Big) \bigg].
\end{aligned}
\end{equation}

In order to proceed with our estimates, we decompose $\Gc_{\leq N}$ into terms with zero, one, and two resonances. To this end, we define the fourth-order chaos by 
\begin{equation}\label{heat:eq-G4}
\begin{aligned}
&\Gc^{(4)}_{\leq N}(s,s^\prime,x;n_5,N_5) \\
=& 9 \cdot 2^{4/2} \mathbf{1}\{ s^\prime \leq s\} 
\sum_{n_1,n_2,n_3,n_4 \in \Z^3} 
\sum_{ \substack{N_{12},N_{345}\colon \\ N_{12} \sim N_{345}, \\ N_{345}\leq N }}
\sum_{\substack{ N_{34} \colon \\ N_{34} \gg N_5}} \bigg[ 
e^{i \langle n_{12345}, x \rangle} \\
\times& 1_{N_{12}}(n_{12}) 1_{N_{345}}(n_{345}) 1_{N_{34}}(n_{34})  1_{N_5}(n_5)
  \Big( \prod_{1\leq j \leq 4} 1_{\leq N}(n_j) \Big) e^{-(s-s^\prime) \langle n_{345} \rangle^2} \\
\times& \Big( \int_{\R^4} \bigotimes_{1\leq j \leq 4} \dW[\cos][s_j][n_j]
\Big( \prod_{j=1,2} \mathbf{1}\{ s_j \leq s \} e^{-(s-s_j)\langle n_j \rangle^2} \Big) 
\Big( \prod_{j=3,4} \mathbf{1}\{ s_j \leq s^\prime \} e^{-(s^\prime-s_j) \langle n_j \rangle^2} \Big) \bigg) \bigg], 
\end{aligned}
\end{equation}
the second-order chaos by 
\begin{equation}\label{heat:eq-G2}
\begin{aligned}
&\Gc^{(2)}_{\leq N}(s,s^\prime,x;n_5,N_5) \\
=&18 \cdot 2^{4/2} \mathbf{1}\{ s^\prime \leq s\} 
\sum_{n_1,n_2,n_3,n_4 \in \Z^3} 
\sum_{ \substack{N_{12},N_{345}\colon \\ N_{12} \sim N_{345}, \\ N_{345}\leq N }}
\sum_{\substack{ N_{34} \colon \\ N_{34} \gg N_5}} \bigg[ \mathbf{1}\{n_{13}=0 \} 
e^{i \langle n_{245}, x \rangle} 
 \\
\times& 1_{N_{12}}(n_{12}) 1_{N_{345}}(n_{345}) 1_{N_{34}}(n_{34}) 1_{N_5}(n_5)  \Big( \prod_{2\leq j \leq 4} 1_{\leq N}(n_j) \Big) e^{-(s-s^\prime) \langle n_{345} \rangle^2} \\
\times& \bigg( \int_{-\infty}^s \int_{-\infty}^{s^\prime} \dW[\cos][s_2][n_2] \dW[\cos][s_4][n_4] \int_{-\infty}^{s^\prime} \ds_3 \Big( 
e^{-(s-s_2) \langle n_2 \rangle^2} e^{-(s+s^\prime-2s_3) \langle n_3 \rangle^2}
e^{-(s^\prime-s_4) \langle n_4 \rangle^2} \Big) \bigg) \bigg],
\end{aligned}
\end{equation}
and the zeroth-order chaos (or constant) by 
\begin{equation}\label{heat:eq-G0}
\begin{aligned}
&\Gc^{(0)}_{\leq N}(s,s^\prime,x;n_5,N_5) \\
=& 18 \cdot 2^{4/2} e^{i\langle n_5, x \rangle} 
\sum_{n_1,n_2,n_3,n_4 \in \Z^3} 
\sum_{ \substack{N_{12},N_{345}\colon \\ N_{12} \sim N_{345}, \\ N_{345}\leq N }}
\sum_{\substack{N_{34} \colon \\ N_{34} \gg N_5 }} \bigg[ \mathbf{1}\{ n_{13}=n_{24}=0 \}
1_{N_{12}}(n_{34}) 1_{N_{345}}(n_{345}) 1_{N_{34}}(n_{34}) \\
\times &  1_{N_5}(n_5)\Big( \prod_{3\leq j \leq 4} 1_{\leq N}(n_j) \Big) 
e^{-(s-s^\prime) \langle n_{345} \rangle^2} 
\Big( \int_{-\infty}^{s^\prime} \int_{-\infty}^{s^\prime} \ds_3 \ds_4 
\prod_{3\leq j \leq 4} e^{-(s+s^\prime-2s_j) \langle n_j \rangle^2} \Big) \bigg].
\end{aligned}
\end{equation}
We now treat the contribution of each chaos separately. \\

\emph{Step 1: Contribution of $\Gc^{(4)}$.} For expository purposes, we decompose the argument into three sub-steps. \\

\emph{Step 1.a:} In the first step, we define a dyadic decomposition of $\Gc^{(4)}$. We define
\begin{equation*}
\Gc^{(4)}(s,s^\prime,x;n_5,N_\ast) = \Gc^{(4)}(s,s^\prime,x;n_5,N_1,N_2,N_3,N_4,N_5,N_{12},N_{34},N_{345})
\end{equation*}
by inserting dyadic cut-offs in \eqref{heat:eq-G4}. More precisely, we remove the sum over (but keep the restrictions on) $N_{12}, N_{34}$, and $N_{345}$ and replace
\begin{equation*}
  \prod_{1\leq j \leq 4} 1_{\leq N}(n_j) 
    \qquad \text{with} \qquad 
 \prod_{1\leq j \leq 4} 1_{N_j}(n_j) . 
\end{equation*}
Due to our previous frequency-restrictions and symmetry, we can always restrict ourselves to the case  
\begin{equation}\label{heat:eq-G4-p1}
N_1 \geq N_2, \, N_3 \geq N_4, \, N_{12} \sim N_{34} \sim N_{345} \gg N_5. 
\end{equation}
In the following, we also write 
\begin{equation*}
N_{\max}:= \max\big( N_1,N_2,N_3,N_4,N_5 , N_{345} \big). 
\end{equation*}

\emph{Step 1.b:} In the second step, we prove for all $p\geq 2$  that
\begin{equation}\label{heat:eq-G4-p2}
\begin{aligned}
&\mathbb{E} \bigg[ \sup_{-1\leq s^\prime \leq s \leq 0} \Big\| P_{N_0} \Gc^{(4)}(s,s^\prime,x,n_5,N_\ast) \Big\|_{L_x^\infty}^p \bigg]^{1/p} \\
\lesssim&\,  p^2 \exp\Big(-1/8 |s-s^\prime| N_{345}^2\Big) N_0^{1/2} N_{345}^{2} \min\big( 1, N_0 N_5^{-1}\big) N_{\textup{max}}^{-1/2+\epsilon}.
\end{aligned}
\end{equation}
 Using the reduction arguments in Subsection \ref{prep:remark-reduction} (or Sobolev embedding in the $s$ and $s^\prime$-variables), it suffices to treat fixed $-1\leq s^\prime \leq s \leq 0$. For any $p\geq 2$, we obtain from Gaussian hypercontractivity that  
\begin{align}
&\E \Big[ \Big\| P_{N_0} \Gc^{(4)}(s,s^\prime,x,n_5,N_\ast) \Big\|_{L_x^\infty}^p \Big]^{2/p} \notag \\
\lesssim& \, p^2 N_{\textup{max}}^{\epsilon}   \E \Big[ \Big\| P_{N_0} \Gc^{(4)}(s,s^\prime,x,n_5,N_5) \Big\|_{L_x^2}^2 \Big] \notag \\
\lesssim&\, p^2 N_{\textup{max}}^{\epsilon} \sum_{\substack{n_0,n_1,n_2,n_3,n_4,n_5 \in \Z^3 \colon \\ n_0=n_{12345}}} \bigg[ \Big( \prod_{j=0}^5 1_{N_j}(n_j) \Big) 1_{N_{12}}(n_{12}) 1_{N_{34}}(n_{34}) 1_{N_{345}}(n_j) \exp\Big(-2(s-s^\prime) \langle n_{345} \rangle^2\Big) \notag \\
&\quad \times \int_{-\infty}^s \int_{-\infty}^s \ds_1 \ds_2 \int_{-\infty}^{s^\prime} \int_{-\infty}^{s^\prime} \ds_3 \ds_4 \Big( \prod_{j=1}^2 e^{-(s-s_j) \langle n_j \rangle^2} \Big) \Big( \prod_{j=3}^4 e^{-(s^\prime-s_j) \langle n_j \rangle^2} \Big) \bigg] \notag \\
\lesssim&\, p^2  N_{\textup{max}}^{\epsilon} \Big( \prod_{j=1}^4 N_j^{-2} \Big) 
\exp\Big(-1/2 (s-s^\prime) N_{345}^2\Big)  \notag  \\
&\quad \times 
\sum_{\substack{n_0,n_1,n_2,n_3,n_4,n_5 \in \Z^3 \colon \\ n_0=n_{12345}}} \bigg[ \Big( \prod_{j=0}^5 1_{N_j}(n_j) \Big) 1_{N_{12}}(n_{12}) 1_{N_{34}}(n_{34}) 1_{N_{345}}(n_j) \bigg].
\label{heat:eq-G4-p3}
\end{align}
By viewing $n_0,n_2,n_{345}$, and $n_4$ as the free variables and recalling \eqref{heat:eq-G4-p1}, we obtain that 
\begin{align}
\eqref{heat:eq-G4-p3}    &\lesssim p^2 N_{\textup{max}} \Big( \prod_{j=1}^4 N_j^{-2} \Big)  (N_0 N_2 N_{345} N_4 )^3 \exp\Big(-1/2 (s-s^\prime) N_{345}^2\Big)  \notag \\
&\lesssim N_{\textup{max}}^\epsilon N_0 N_{345}^4 \times \Big( N_0^2 N_1^{-1} N_3^{-1} N_{345}^{-1} \Big) \times \exp\Big(-1/2 (s-s^\prime) N_{345}^2\Big). \label{heat:eq-G4-p4}
\end{align}
In order to prove \eqref{heat:eq-G4-p2}, it only remains to bound the second factor in  \eqref{heat:eq-G4-p4}. Under our frequency-restrictions in \eqref{heat:eq-G4-p1}, it holds that $\max(N_0,N_5) \lesssim N_{345}$ and  $N_{345} \lesssim \min(N_1,N_3)$. As a result,
\begin{align*}
N_0^2 N_1^{-1} N_3^{-1}  N_{345}^{-1} 
&\lesssim \big( N_0^2 N_{345}^{-2} \big) \times \big( N_1^{-1} N_3^{-1} N_{345} \big)\\ &\lesssim  \min(1,N_0^2 N_5^{-2}) \max(N_1,N_3)^{-1} \\
&\sim  \min(1,N_0^2 N_5^{-2}) N_{\textup{max}}^{-1}. 
\end{align*}
This yields the desired estimate in \eqref{heat:eq-G4-p2}. \\

\emph{Step 1.c:}
In the final step, we control the contribution of $\Gc^{(4)}$ to the left-hand side of \eqref{heat:eq-main-estimates}. Due to Lemma \ref{heat:lem-duhamel}, it suffices to prove for all frequency-scales that 
\begin{equation}\label{heat:eq-G4-p5}
\begin{aligned}
&\sup_{s\in [s_0,0]} |s-s_0|^\theta \Big\| P_{N_0} \sum_{n_5 \in \Z^3} \int_{s_0}^s \ds^\prime 
\widehat{P_{N_5} \Psi}(s^\prime,n_5) \Gc^{(4)}(s,s^\prime,x;n_5,N_\ast) \Big\|_{\C_x^{-1/2-4\delta}} \\
\lesssim& \, A N_{\textup{max}}^{-\epsilon} \sup_{s\in [s_0,0]} |s-s_0|^\theta \big\| \Psi(s) \big\|_{H_x^{1-\delta}}.
\end{aligned}
\end{equation}
By using \eqref{heat:eq-G4-p2}, we obtain (on an $A$-certain event) that, for all $-1\leq s^\prime \leq s \leq 0$ and all frequency-scales,  
\begin{equation}
\Big\| P_{N_0} \Gc^{(4)}(s,s^\prime,x;n_5,N_\ast) \Big\|_{L_x^\infty} 
\leq A \exp\Big(-1/8 |s-s^\prime| N_{345}^2 \Big) N_0^{1/2} N_{345}^2 \min(1,N_0 N_5^{-1}) N_{\textup{max}}^{-1/2+\epsilon}. 
\end{equation}
As a result, 
\begin{align*}
&\Big\| P_{N_0} \sum_{n_5 \in \Z^3} \int_{s_0}^s \ds^\prime 
\widehat{P_{N_5} \Psi}(s^\prime,n_5) \Gc^{(4)}(s,s^\prime,x;n_5,N_\ast) \Big\|_{\C_x^{-1/2-4\delta}} \\
\lesssim& \, \int_{s_0}^s \ds^\prime \big\| P_{N_0}  \Gc^{(4)}(s,s^\prime,x;n_5,N_\ast) \big\|_{\C_x^{-1/2-4\delta}} \sum_{n_5} \big| \widehat{P_{N_5} \Psi}(s^\prime,n_5) \big| \\
\lesssim& \, A N_0^{-4\delta} \min\big(1,N_0 N_5^{-1}\big) N_{345}^2 N_5^{1/2+\delta} N_{\textup{max}}^{-1/2+\epsilon} \\
\times& \Big( \int_{s_0}^s \ds^\prime \exp\Big(-1/8 |s-s^\prime| N_{345}^2 \Big) |s^\prime-s_0|^{-\theta} \Big) \sup_{s^\prime \in [s_0,0]} |s^\prime -s_0|^\theta \big\| \Psi(s^\prime) \big\|_{H_x^{1-\delta}}. 
\end{align*}
From a direct calculation, we obtain that 
\begin{align*}
     &A N_0^{-4\delta} \min\big(1,N_0 N_5^{-1}\big) N_{345}^2 N_5^{1/2+\delta} N_{\textup{max}}^{-1/2+\epsilon} \Big( \int_{s_0}^s \ds^\prime \exp\Big(-1/8 |s-s^\prime| N_{345}^2 \Big) |s^\prime-s_0|^{-\theta} \Big) \\
     \lesssim& \, A N_{345}^{2\epsilon} N_5^{1/2-3\delta} N_{\textup{max}}^{-1/2+\epsilon}
      \Big( \int_{s_0}^s \ds^\prime  |s-s^\prime|^{-(1-\epsilon)} |s^\prime-s_0|^{-\theta} \Big) \\
      \lesssim& \,  A  N_{\textup{max}}^{3(\epsilon-\delta)} |s-s_0|^{-\theta}. 
\end{align*}
This yields the desired estimate \eqref{heat:eq-G4-p5}. \\

\emph{Step 2: Contribution of $\Gc^{(2)}$.} Since the argument is similar to the treatment of $\Gc^{(4)}$, we omit the details. \\

\emph{Step 3: Contribution of $\Gc^{(0)}$.} 
We emphasize that the contribution of $\Gc_{\leq N}^{(0)}$ has to be renormalized, i.e.,  has to cancel with $\gamma_{\leq N}$.
We first perform the sum over the  frequencies $n_1$ and $n_2$ and perform the integrals in $s_3$ and $s_4$, which yields
\begin{equation}\label{heat:eq-G0-simplified}
\begin{aligned}
&\Gc^{(0)}_{\leq N}(s,s^\prime,x;n_5,N_5) \\
=& 18 \cdot 1_{N_5}(n_5) \, e^{i\langle n_5, x \rangle} 
\sum_{n_3,n_4 \in \Z^3} 
\sum_{\substack{N_{12},N_{345} \colon \\ N_{12} \sim N_{345}, \\ N_{345} \leq N}} \sum_{\substack{N_{34} \colon N_{34} \gg N_5 }} \bigg[
1_{N_{12}}(n_{34}) 1_{N_{345}}(n_{345}) 1_{N_{34}}(n_{34}) \\
\times & 1_{\leq N}(n_{345}) \Big( \prod_{3\leq j \leq 4} 1_{\leq N}(n_j) \Big) 
e^{-(s-s^\prime) \langle n_{345} \rangle^2} \Big( \prod_{3\leq j \leq 4} \frac{e^{-(s-s^\prime) \langle n_j \rangle^2}}{\langle n_j \rangle^2} \Big) \bigg].
\end{aligned}
\end{equation}
We now isolate the main term in \eqref{heat:eq-G0-simplified}, which is given by 
\begin{equation}
\begin{aligned}
&\widetilde{\Gc}\vphantom{\Gc}^{(0)}_{\leq N}(s,s^\prime,x;n_5,N_5) \\ 
=& 18 \cdot 1_{N_5}(n_5)  e^{i\langle n_5, x \rangle} 
\sum_{n_3,n_4 \in \Z^3} \bigg[ 
 1_{\leq N}(n_{345}) \Big( \prod_{3\leq j \leq 5} 1_{\leq N}(n_j) \Big) 
e^{-(s-s^\prime) \langle n_{345} \rangle^2} \Big( \prod_{3\leq j \leq 4} \frac{e^{-(s-s^\prime) \langle n_j \rangle^2}}{\langle n_j \rangle^2} \Big) \bigg]. 
\end{aligned}
\end{equation}
In the error term $\Gc^{(0)}_{\leq N}-\widetilde{\Gc}\vphantom{\Gc}^{(0)}_{\leq N}$, one can gain a factor of $N_5 N_{345}^{-1}$, which makes the estimate rather easy.  
As a result, we focus only on the main term $\widetilde{\Gc}\vphantom{\Gc}^{(0)}_{\leq N}$. Using symmetry in $n_3,n_4$, and $n_{345}$, it follows that 
\begin{equation}\label{heat:eq-G0-symmetrized}
\begin{aligned}
&\widetilde{\Gc}\vphantom{\Gc}^{(0)}_{\leq N}(s,s^\prime,x;n_5,N_5) \\
=& - 6 \cdot 1_{N_5}(n_5) e^{i\langle n_5 ,x \rangle} \partial_s \bigg(
\sum_{n_3,n_4 \in \Z^3} \bigg[ 
 1_{\leq N}(n_{345})  \frac{e^{-(s-s^\prime) \langle n_{345} \rangle^2}}{\langle n_{345} \rangle^2} \Big( \prod_{3\leq j \leq 4} 1_{\leq N}(n_j) \frac{e^{-(s-s^\prime) \langle n_j \rangle^2}}{\langle n_j \rangle^2} \Big) \bigg] \bigg). 
\end{aligned}
\end{equation}
With a slight abuse of notation\footnote{In the main part of the paper, the second argument corresponds to the time-variable in the wave equation. In this appendix, the second argument corresponds to the time-variable in the heat equation. Furthermore, since we mostly work with $s-s^\prime \geq 0$, we changed the sign in the exponential.}, we define
\begin{equation}
\Gamma_{\leq N}(n_5,s;N_5) := 6 \cdot 1_{N_5}(n_5) \sum_{n_3,n_4 \in \Z^3} \bigg[ 
 1_{\leq N}(n_{345})  \frac{e^{-s \langle n_{345} \rangle^2}}{\langle n_{345} \rangle^2} \Big( \prod_{3\leq j \leq 4} 1_{\leq N}(n_j) \frac{e^{-s \langle n_j \rangle^2}}{\langle n_j \rangle^2} \Big) \bigg] \bigg).
\end{equation}
Equipped with this notation, \eqref{heat:eq-G0-symmetrized} reads
\begin{equation}\label{heat:eq-G0-Gamma}
\widetilde{\Gc}\vphantom{\Gc}^{(0)}_{\leq N}(s,s^\prime,x;n_5,N_5)= - \partial_s \Gamma_{\leq N}(n_5,s-s^\prime;N_5) \, e^{i\langle n_5, x \rangle}. 
\end{equation}
It follows that 
\begin{align}
&\sum_{n_5 \in \Z^3} \int_{s_0}^s \ds^\prime \widehat{P_{N_5} \Psi}(s^\prime,n_5) \widetilde{\Gc}\vphantom{\Gc}^{(0)}_{\leq N}(s,s^\prime,x;n_5,N_5) - \gamma_{\leq N} P_{N_5} \Psi(s) \notag  \\
=& \, - \partial_s \sum_{n_5 \in \Z^3} \int_{s_0}^s \ds^\prime \widehat{P_{N_5}\Psi}(s^\prime,n_5) \Gamma_{\leq N}(n_5,s-s^\prime;N_5) e^{i \langle n_5 , x \rangle} \label{heat:eq-G0-p1} \\
+& \, \sum_{n_5 \in \Z^3} \big( \Gamma_{\leq N}(n_5,0) - \gamma_{\leq N}\big) \widehat{P_{N_5} \Psi}(s,n_5) e^{i\langle n_5 ,x  \rangle} \label{heat:eq-G0-p2}. 
\end{align}
For the first term \eqref{heat:eq-G0-p1}, we use the second estimate in Lemma \ref{heat:lem-duhamel}, which reduces the desired estimate to
\begin{align*}
&\sup_{s\in [s_0,0]} |s-s_0|^\theta \, \Big\|   \sum_{n_5 \in \Z^3} \int_{s_0}^s \ds^\prime \widehat{P_{N_5}\Psi}(s^\prime,n_5) \Gamma_{\leq N}(n_5,s-s^\prime;N_5) e^{i \langle n_5 , x \rangle} \Big\|_{\C_x^{3/2-4\delta}} \\
\lesssim&\, \sup_{s\in [s_0,0]} |s-s_0|^\theta \, \big\| \Psi(s) \big\|_{H_x^{1-\delta}}. 
\end{align*}
This follows directly from the triangle inequality in $n_5$, inserting the definition of $\Gamma_{\leq N}$, and performing the $s^\prime$-integral. 

For the second term \eqref{heat:eq-G0-p2}, we use the first estimate in Lemma \ref{heat:lem-duhamel}, which reduces the desired estimate to 
\begin{align*}
\sup_{s\in [s_0,0]} |s-s_0|^\theta \, \Big\|  \sum_{n_5 \in \Z^3} \big( \Gamma_{\leq N}(n_5,0) - \gamma_{\leq N}\big) \widehat{P_{N_5} \Psi}(s,n_5) e^{i\langle n_5 ,x  \rangle} \Big\|_{\C_x^{-1/2-4\delta}} 
\lesssim\, \sup_{s\in [s_0,0]} |s-s_0|^\theta \, \big\| \Psi(s) \big\|_{H_x^{1-\delta}}. 
\end{align*}
This follows from the triangle inequality in $n_5$ and Lemma \ref{analytic:lem-difference-gamma}.\end{proof}

\section{Merging Estimates, Moment Method and Time Integrals}\label{section:mergetrimtime-appendix} We first recall the main technical tools needed from \cite{DNY20}, namely the merging estimate (Lemma \ref{counting:lem-merging}) and the moment method (Proposition \ref{counting:prop-moment}, called trimming estimate in \cite{DNY20}).
\begin{lemma}[Merging estimates, Proposition 4.11 in \cite{DNY20}]\label{counting:lem-merging}
Consider two tensors $h_{k_{A_1}}^{(1)}$ and $h_{k_{A_2}}^{(2)}$, where $A_1\cap A_2=C$. Let $A_1\Delta A_2=A$, define the semi-product
 \begin{equation}\label{combination}H_{k_A}=\sum_{k_C}h_{k_{A_1}}^{(1)}h_{k_{A_2}}^{(2)}.
 \end{equation} Then, for any partition $(X,Y)$ of $A$, let $X\cap A_1=X_1$, $Y\cap A_1 =Y_1$ etc., we have
 \begin{equation}\label{combinationbd}\|H\|_{k_X\to k_Y}\leq \|h^{(1)}\|_{k_{X_1\cup C}\to k_{Y_1}}\cdot\|h^{(2)}\|_{k_{X_2}\to k_{C\cup Y_2}}.
 \end{equation}
\end{lemma}
\begin{proposition}[Moment method, Proposition 4.14 in \cite{DNY20}]\label{counting:prop-moment}
Let $\mathcal{A}$, $\mathcal{X}$, and $\mathcal{Y}$ be disjoint finite index sets and let $h=h_{n_\mathcal{A} n_\mathcal{X} n_\mathcal{Y}}$ be a (deterministic) tensor. Let $N$ be a dyadic frequency-scale and assume that, on the support of $h$, $|n_j|\lesssim N$ for all $j\in \mathcal{A}\cup \mathcal{X} \cup \mathcal{Y}$. Furthermore, let $(\pm_j)_{j\in \mathcal{A}} \in \big\{ +, - \big\}^\mathcal{A}$ be a collection of signs. Finally, define the random tensor $H=H_{n_\mathcal{X} n_\mathcal{Y}}$ by 
\begin{equation}
H_{n_\mathcal{X} n_\mathcal{Y}} = \sum_{n_\mathcal{A}} h_{n_\mathcal{A} n_\mathcal{X} n_\mathcal{Y}} \SI[n_j,\pm_j \colon j\in \mathcal{A}],
\end{equation}
where the stochastic integrals $\SI$ are as in Subsection \ref{section:prelim-multiple-stochastic}. Then, we have for all $\delta>0$ and all $p\geq 1$ that 
\begin{equation}
\Big\| \big\| H_{n_\mathcal{X} n_\mathcal{Y}} \big\|_{n_\mathcal{X} \rightarrow n_\mathcal{Y}}\Big\|_{L^p_\omega}
\lesssim_\delta N^\delta p^{\# \mathcal{A}/2} \max_{\mathcal{B},\mathcal{C}} \big\| h_{n_\mathcal{A} n_\mathcal{X} n_\mathcal{Y}} \big\|_{n_\mathcal{B} n_\mathcal{X} \rightarrow n_\mathcal{C} n_\mathcal{Y}},
\end{equation}
where the maximum is taken over all partitions of $\mathcal{A}$. 
\end{proposition}
\begin{remark}\label{Bee3} Note that Proposition \ref{counting:prop-moment} is slightly different from Proposition 4.14 in \cite{DNY20}, due to the use of the renormalized product $\Sc\Ic$ instead of products of independent Gaussians. However, this difference does not affect the proof.
\end{remark}
In addition, we prove an almost $L^1$ estimate for iterated oscillatory time integrals (Proposition \ref{counting:timeint}). This is a special case of Lemma 10.2 in \cite{DH21}, and is used in the molecular analysis in Section \ref{section:analytic2}.
  \begin{proposition}[Some time integral estimates]\label{counting:timeint} Consider the following expressions:
\begin{align}\label{counting:timeint1}
\Hc_{1,3,3}(t,\Omega_0,\Omega_1,\Omega_2)&:=\int_0^t\int_0^t \chi(t)\chi(t_1)\chi(t_2)e^{i(\Omega_0t+\Omega_1t_1+\Omega_2t_2)}\,\mathrm{d}t_1\mathrm{d}t_2,\\
\label{counting:timeint2}
\Hc_{3,3,3}(t,\Omega_0,\Omega_1,\Omega_2,\Omega_3)&:=\int_0^t\int_0^t \int_0^t\chi(t)\chi(t_1)\chi(t_2)e^{i(\Omega_0t+\Omega_1t_1+\Omega_2t_2+\Omega_3t_3)}\,\mathrm{d}t_1\mathrm{d}t_2\mathrm{d}t_3,\\
\label{counting:timeint3}
\Hc_{1,1,5}(t,\Omega_0,\Omega_1,\Omega_2)&:=\int_0^t\int_0^{t} \int_0^{t_1}\chi(t)\chi(t_1)\chi(t_2)e^{i(\Omega_0t+\Omega_1t_1+\Omega_2t_2)}\,\mathrm{d}t_2\mathrm{d}t_1.
\end{align} Then we have the following estimates:
\begin{align}\label{counting:timeintest1}
\int_{\Rb^3}(\langle\Omega_0\rangle\langle \Omega_1\rangle \langle\Omega_2\rangle)^{8(1/2-b_+)}|(\Fc_t\Hc_{1,3,3})(\xi,\Omega_0,\Omega_1,\Omega_2)|\,\mathrm{d}\Omega_0\mathrm{d}\Omega_1\mathrm{d}\Omega_2&\lesssim \langle \xi\rangle^{4(1/2-b_+)},\\
\label{counting:timeintest2}
\int_{\Rb^4}(\langle\Omega_0\rangle\langle \Omega_1\rangle \langle\Omega_2\rangle\langle \Omega_3\rangle)^{8(1/2-b_+)}|(\Fc_t\Hc_{3,3,3})(\xi,\Omega_0,\Omega_1,\Omega_2,\Omega_3)|\,\mathrm{d}\Omega_0\mathrm{d}\Omega_1\mathrm{d}\Omega_2\mathrm{d}\Omega_3&\lesssim \langle \xi\rangle^{4(1/2-b_+)},\\
\label{counting:timeintest3}
\int_{\Rb^3}(\langle\Omega_0\rangle\langle \Omega_1\rangle \langle\Omega_2\rangle)^{8(1/2-b_+)}|(\Fc_t\Hc_{1,1,5})(\xi,\Omega_0,\Omega_1,\Omega_2)|\,\mathrm{d}\Omega_0\mathrm{d}\Omega_1\mathrm{d}\Omega_2&\lesssim \langle \xi\rangle^{4(1/2-b_+)}.
\end{align} The same estimates hold for all $\Omega_j$ derivatives of these functions.
\end{proposition}
\begin{proof} First, the $\Omega_j$ derivatives of the $\Hc$ functions satisfy the same estimates as themselves because taking one $\Omega_j$ derivative in (\ref{counting:timeint1})--(\ref{counting:timeint3}) just corresponds to multiplying by $t_j$ or $t$, which is bounded due to the cutoff $\chi$. Thus we will only prove the bounds for the original $\Hc$ functions. Moreover we will only prove (\ref{counting:timeintest3}) for $\Hc_{1,1,5}$ defined by (\ref{counting:timeint3}), since the proof for the other two will be similar (and simpler as they do not involve iterated time integrals).

Now consider $\Hc:=\Hc_{1,1,5}$. Let $1=\eta_0(\Omega)+\eta_\infty(\Omega)$ be a partition of unity supported on $|\Omega|\lesssim 1$ and $|\Omega|\gtrsim 1$ respectively, then integrating by parts we have
\begin{equation}\label{counting:ibp}\int_0^{t_1}\chi(t_2)e^{i\Omega_2t_2}\,\mathrm{d}t_2=\eta_0(\Omega_2)\psi(t_1,\Omega_2)+\eta_\infty(\Omega_2)\bigg(\frac{\chi(t_1)e^{i\Omega_2t_1}-\chi(0)}{i\Omega_2}-\frac{1}{i\Omega_2}\int_0^{t_1}\chi'(t_2)e^{i\Omega_2t_2}\,\mathrm{d}t_2\bigg),\end{equation} where $\psi$ is a fixed smooth function of $(t_1,\Omega_2)$. Note that in (\ref{counting:ibp}), the first term is a Schwartz function upon localizing in $t_1$ (which we always can do), and the last term can be integrated by parts again to get more decay in $\Omega_2$, so it is always manageable provided we can handle the other terms. Thus we will focus on the main term
\[\eta_\infty(\Omega_2)\frac{\chi(t_1)e^{i\Omega_2t_1}-\chi(0)}{i\Omega_2},\] and plug it into (\ref{counting:timeint3}) to get
\[\frac{\eta_\infty(\Omega_2)}{i\Omega_2}e^{i\Omega_0t}\int_0^t\chi(t_1)(\chi(t_1)e^{i\Omega_2t_1}-\chi(0))e^{i\Omega_1t_1}\,\mathrm{d}t_1.\] Integrating by parts again, we obtain some remainder terms (which are manageable provided we can handle the main term) plus the main term which is
\begin{equation}\label{counting:mainint}
e^{i\Omega_0t}\bigg[\frac{\eta_\infty(\Omega_2)}{i\Omega_2}\frac{\eta_\infty(\Omega_1+\Omega_2)}{i(\Omega_1+\Omega_2)}(\chi^2(t)e^{i(\Omega_1+\Omega_2)t}-\chi^2(0))-\frac{\eta_\infty(\Omega_2)}{i\Omega_2}\frac{\eta_\infty(\Omega_1)}{i\Omega_1}\chi(0)(\chi(t)e^{i\Omega_1t}-\chi(0))\bigg].
\end{equation}

For this main term we clearly have that\begin{equation}\label{counting:timeintex}|\Fc_x(\ref{counting:mainint})(\xi)|\lesssim\sum_{\Omega\in\{\Omega_0,\Omega_0+\Omega_1,\Omega_0+\Omega_1+\Omega_2\}}\bigg(\bigg|\frac{\eta_\infty(\Omega_2)}{i\Omega_2}\frac{\eta_\infty(\Omega_1+\Omega_2)}{i(\Omega_1+\Omega_2)}\bigg|+\bigg|\frac{\eta_\infty(\Omega_2)}{i\Omega_2}\frac{\eta_\infty(\Omega_1)}{i\Omega_1}\bigg|\bigg)\langle \xi-\Omega\rangle^{-10}.\end{equation} This easily implies (\ref{counting:timeintest3}), because both functions
\[\frac{\eta_\infty(\Omega_2)}{i\Omega_2}\frac{\eta_\infty(\Omega_1+\Omega_2)}{i(\Omega_1+\Omega_2)}\quad\textrm{and}\quad \frac{\eta_\infty(\Omega_2)}{i\Omega_2}\frac{\eta_\infty(\Omega_1)}{i\Omega_1}\] are almost $L^1$ in $(\Omega_1,\Omega_2)$, and becomes $L^1$ when multiplied by the negative weight $(\langle \Omega_1\rangle \langle\Omega_2\rangle)^{8(1/2-b_+)}$; the integrability in $\Omega_0$, as well as decay in $\xi$ in (\ref{counting:timeintest3}), follows because of the $\langle \xi-\Omega\rangle^{-10}$ factor in (\ref{counting:timeintex}) (so in particular we can assume $|\xi|\lesssim1+\max_j|\Omega_j|$). Now, once (\ref{counting:timeintest3}) is proved for the main term, the same argument can easily be applied to the remainder terms to get the same estimate. This finishes the proof of (\ref{counting:timeintest3}).
\end{proof}

\end{appendix}

\bibliographystyle{myalpha}
\bibliography{Cubic_Library}

\end{document}